\documentclass[a4paper,11pt,twoside]{report}
\usepackage[top=3cm,bottom=2cm,left=2cm,right=2cm,twoside,includefoot, heightrounded]{geometry}

\usepackage[final]{pdfpages} 
\usepackage{fancyhdr}
\usepackage{subfigure} 
\usepackage{pifont} 
\usepackage{color}
\usepackage{enumerate} 
\usepackage{amssymb,amsmath, mathtools} 
\usepackage[mathscr]{eucal}
\usepackage{multicol}
\usepackage[titles,subfigure]{tocloft} 
\usepackage{cite} 
\usepackage{url}
\usepackage{amsfonts} 
\usepackage{amsbsy} 
\usepackage{amscd}
\usepackage{amssymb}
\usepackage{algorithm2e} 
\usepackage{algorithmic}
\usepackage{morefloats}

\newcommand{\unnn}[1]{\underline{#1}}
\newcommand{\vveh}{\varkappa_{\text{v}}}
\newcommand{\vvehu}{\varkappa_{\text{v,u}}}
\newcommand{\spr}[2]{\left\langle #1; #2 \right\rangle}
\newcommand{\eer}[1]{\overline{#1}_{\text{est}}}

\newcommand{\ov}[1]{\overline{#1}}

\newcommand{\col}{\,\mbox{\bf col}\,}
\newcommand{\pf}{\paragraph{Proof}} 
\newcommand{\ovk}{\overline{\varkappa}}

\newcommand{\unk}{\underline{\varkappa}}
\newcommand{\ovkv}{\ovk_{\text{v}}} 

 \renewcommand{\Box}{\bullet}

\newcommand{\br}{\mathbb{R}}

\newcommand{\ve}{\varepsilon} 
  
 \newcommand{\epf}{$\quad \Box$}

\newcommand{\bt}{\boldsymbol{\mathscr{T}}}
\newcommand{\dtrig}{d_{\updownarrow}}

\newcommand{\sgn}{\text{\bf sgn}} \newcommand{\sat}{\text{\bf sat}}
\newcommand{\bldr}{\boldsymbol{r}} \newcommand{\bldv}{\boldsymbol{v}}
\newcommand{\blds}{\boldsymbol{s}} \newcommand{\bldu}{\boldsymbol{u}}
\newcommand{\bldw}{\boldsymbol{w}} \newcommand{\bldp}{\boldsymbol{p}}
\newcommand{\bldx}{\boldsymbol{x}}
\newcommand{\state}{\mathscr{S}} \newcommand{\contr}{\mathscr{U}}

\newcommand{\dist}{\text{\bf dist}}
\newcommand{\gv}[1]{\ensuremath{\mbox{\boldmath $ #1 $}}}
\newcommand{\grad}[1]{\gv{\nabla} #1} 

\newcommand{\dsafe}{d_{safe}}

\renewcommand{\phi}{\varphi}
\renewcommand{\epsilon}{\varepsilon}
\renewcommand{\kappa}{\varkappa}
\newcommand{\trs}{\top}

\newcommand{\blr}{r}
\newcommand{\interior}{\mathbf{int}}

\renewcommand{\phi}{\varphi}
\renewcommand{\epsilon}{\varepsilon}
\renewcommand{\kappa}{\varkappa}

\newtheorem{Remark}{Remark}[section]
\newtheorem{Definition}{Definition}[section]
\newtheorem{Assumption}{Assumption}[section]
\newtheorem{Theorem}{Theorem}[section]
\newtheorem{Lemma}{Lemma}[section]
\newtheorem{Proposition}{Proposition}[section]
\newtheorem{Algorithm}{Algorithm}[section]
\newtheorem{Requirement}{Requirement}[section]
\newtheorem{property}{Property}[section]

\newlength{\mylength} 

\numberwithin{equation}{chapter}

\providecommand{\abs}[1]{\lvert#1\rvert}
\providecommand{\norm}[1]{\lVert#1\rVert}

\begin{document}

\pagenumbering{roman}

\title {Methods for Collision--Free Navigation of Multiple Mobile Robots in Unknown Cluttered
Environments}
\author{\emph{Michael Colin Hoy}}

\maketitle






\begin{center}
\textbf{\large Abstract}
\end{center}

Navigation and guidance of autonomous vehicles is a fundamental problem in robotics, which has attracted intensive research in recent decades. This report is mainly concerned with provable collision avoidance of multiple autonomous vehicles operating in unknown cluttered environments, using reactive decentralized navigation laws, where obstacle information is supplied by some sensor system. 

Recently, robust and decentralized variants of model predictive control based navigation systems have been applied to vehicle navigation problems. Properties such as provable collision avoidance under disturbance and provable convergence to a target have been shown; however these often require significant computational and communicative capabilities, and don't consider sensor constraints, making real time use somewhat difficult. There also seems to be opportunity to develop a better trade-off between tractability, optimality, and robustness.

The main contributions of this work are as follows; firstly, the integration of the robust model predictive control concept with reactive navigation strategies based on local path planning, which is applied to both holonomic and unicycle vehicle models subjected to acceleration bounds and disturbance; secondly, the extension of model predictive control type methods to situations where the information about the obstacle is limited to a discrete ray-based sensor model, for which provably safe, convergent boundary following can be shown; and thirdly the development of novel constraints allowing decentralized coordination of multiple vehicles using a robust model predictive control type approach, where a single communication exchange is used per control update, vehicles are allowed to perform planning simultaneously, and coherency objectives are avoided.

Additionally, a thorough review of the literature relating to collision avoidance is performed; a simple method of preventing deadlocks between pairs of vehicles is proposed which avoids graph-based abstractions of the state space; and a discussion of possible extensions of the proposed methods to cases of moving obstacles is provided. Many computer simulations and real world tests with multiple wheeled mobile robots throughout this report confirm the viability of the proposed methods. Several other control systems for different navigation problems are also described, with simulations and testing demonstrating the feasibility of these methods.

\tableofcontents

\addcontentsline{toc}{chapter}{List of Tables}
\listoftables

\addcontentsline{toc}{chapter}{List of Figures}
\listoffigures



 

\pagenumbering{arabic}
\chapter{Introduction}
\label{chap:intro}

\section{Overview}

\label{sec:intro}

Navigation of autonomous vehicles is an important, classic research area in robotics, and many
approaches are well documented in the literature. However there are many aspects which remain an
increasingly active area of research. A review of recently proposed navigation methods
 applicable to collision avoidance is provided in Chapt.~\ref{chap:lit}.

Both single and coordinated groups of autonomous vehicles have many applications, such as industrial, office, 
and agricultural automation; search and rescue; and surveillance and inspection. All these problems
contain some similar elements, and collision avoidance in some form is almost universally needed.
Examples of compilations of potential applications may be readily found, see e.g. \cite{USDoD2005book1,
Shoenwald2000journ7}.

In contrast to initial approaches to vehicle navigation problems, the focus in the literature has
shifted to navigation laws which are capable of rigorous collision avoidance, such that
for some set of assumptions, it can be proven collisions will
  never occur. Overall, navigation systems with more general assumptions would be considered superior.
   Some examples of common assumptions are listed as follows:
  
  \begin{itemize}
    \item Vehicle models vary in complexity from velocity controlled linear models to realistic
    car-like models (see Sec.~\ref{sec:cluttered}). For example, collision avoidance for velocity
    controlled models is simpler since the vehicles can halt instantly if required; however this is physically
    unrealistic. This means a more complex model which better characterizes actual vehicles is desirable for use
    during analysis.
    \item Different levels of knowledge about the obstacles and other vehicles are required
     by different navigation strategies. This ranges from abstracted obstacle set information to allowance for the
     actual nature of realistic noisy sensor data obtained from range-finding sensors. In addition
     to the realism of the sensor model, the sensing requirements significantly varies
     between approaches -- in some vehicle implementations less powerful sensors are installed, 
     which would only provide limited 
     information to the navigation system (such as the minimum distance to the obstacle).
     \item Different assumptions about the shape of static obstacles have been proposed. To ensure
     correct behavior when operating near an obstacle when only limited information is available, it is
     often unavoidably necessary to presume smoothness properties relating to the obstacle boundary. However, approaches which have
     more flexible assumptions about obstacles would likely be more widely applicable to real world
     scenarios.
    \item Uncertainty is always present in real robotic systems, and proving a behavior occurs
    under an exact vehicle model does not always imply the same behavior will be exhibited when the system is
    implemented.
    To better reflect this, assumptions can be made describing bounded disturbance from the
  nominal model, bounded sensor errors, and the presence of communication errors.
 A review of the types of uncertainty present in vehicle systems is available 
 \cite{Dadkhah2012journ2}.

 \end{itemize}
 
While being able prove collision avoidance under broad circumstances is of utmost importance,
examples of other navigation law features which determine their effectiveness are listed as follows:
 
 \begin{itemize}
    \item Navigation laws which provably achieve the vehicle's goals are highly desirable.
  When possible, this can be generally shown by providing an upper bound on the time in which the
  vehicle will complete a finite task, however conservative this may be. However, proving goal
  satisfaction is possibly less critical than collision avoidance, so long as it can be experimentally
  demonstrated non-convergence is virtually non-existent.
  \item  In many applications, the computational ability of the vehicle is limited, and approaches
  with lower computational cost are favored. However, with ever increasing computer power, this concern
  is mainly focused towards small, fast vehicles such as miniature UAV's (for which the fast update rates required
  are not congruent with the limited computational faculties available).
  \item  Many navigation laws are constructed in continuous time. Virtually all digital control systems
  are updated in discrete time, thus navigation laws constructed in discrete time are more suitable
  for direct implementation.
  \end{itemize}

The key distinction between different approaches is the amount of information they have available
about the workspace. When full information is available about the obstacle set and a single vehicle
is present, \textit{global path planning} methods may be used to find the optimal path. When only 
local information is available, \textit{sensor based} methods are used. A subset of sensor based 
methods are \textit{reactive} methods, which may be expressed as a mapping between the sensor state and control input,
with no memory present.
\par
Recently, \textit{Model Predictive Control} (MPC) architectures have been applied to collision avoidance
problems, and this approach seems to show great potential in providing efficient
navigation, and easily extends to robust and nonlinear problems
(see Sec.~\ref{sec:mpc}). They have many favorable properties compared to the commonly used
\textit{Artificial Potential Field} (APF) methods and \textit{Velocity Obstacle} based methods, which could be generally more 
conservative when extended to higher order vehicle models. MPC continues to be developed 
and demonstrate many desirable properties for sensor based
navigation, including avoidance of moving
obstacles (see Sec.~\ref{sec:moving}), and coordination of multiple vehicles (see Sec.~\ref{sec:multiple}). Additionally,
 the use of MPC for sensor based boundary following problems has been proposed in this report
(see Chapt.~\ref{chap:convsingle}).

\section{Chapter Outline}

The problem statements and main contributions of each chapter are listed as follows:

 \begin{itemize}
 \item Chapt.~\ref{chap:lit} is based on \cite{hoy2013review}, and presents a review of the literature 
 relating to collision avoidance of static obstacles, moving
 obstacles, and other cooperating vehicles.
 
\item Chapt.~\ref{chap:singlevehicle} is based on \cite{Hoy2012journ5} and \cite{Hoy2012journ2}, and is concerned with the navigation of a single vehicle described by either holonomic or unicycle vehicle models with
bounded acceleration and disturbance. Collision avoidance is able to be proven, however convergence to the target is not analytically addressed (since assumptions regarding sensor data are not particularly realistic). The contributions of this chapter are the integration of the robust MPC concept
with reactive navigation strategies based on local path planning, applied to both holonomic and unicycle vehicle models with bounded acceleration and disturbance

\item Chapt.~\ref{chap:convsingle} is based on \cite{Hoy2011conf7}, and is concerned with navigation in unknown
environments where the information about the obstacle is limited to a discrete ray-based sensor
model. Collision avoidance, complete transversal of the obstacle and finite completion
time are able to be proven, however the vehicle model is assumed to be free from disturbance. To solve this problem, the MPC method from Chapt.~\ref{chap:singlevehicle} is extended to only consider the more limited obstacle information when performing the planning process.

\item Chapt.~\ref{chap:multiple} is based on \cite{Hoy2012journ2}, and is concerned with the navigation of multiple cooperating vehicles which are only allowed a limited number of communication exchanges per control update. Again, collision avoidance is able to be proven. To achieve this, a novel constraint is developed allowing decentralized robust coordination of multiple vehicles.  Unlike other approaches in this area, explicit
ordering of the vehicles is not required, each vehicle can plan trajectories simultaneously,
and no explicit coherency objectives are required.

\item Chapt.~\ref{chap:dead} is based on \cite{hoy2012dead}, and is concerned with preventing deadlocks between pairs of robots in cluttered environments. To achieve this, a simple yet original
method of resolving deadlocks between pairs of vehicles is proposed. The proposed method has the advantage of not requiring centralized computation or discrete abstraction of the state space.
\end{itemize}

In Chapts.~\ref{chap:singlevehicle} to \ref{chap:dead}, simulations and real 
world testing confirm the viability of the proposed methods. It should be emphasized they are all 
related since they all use the same basic planning and collision avoidance procedure proposed in 
Chapt.~\ref{chap:singlevehicle}.
\par
In Chapts.~\ref{chapt:rbf} to \ref{chapt:tcp}, the contribution of the author was the simulations 
and testing performed to validate some other varieties of navigation 
systems. In each of these chapters, some of the introductory
discussion, and all mathematical analysis, was contributed
by coauthors to the associated manuscript; however some details are included to assist 
understanding of the presented results. Each of these chapters are relatively independent (they solve quite distinct problems),
and some have introductions independent of Chapt.~\ref{chap:lit}.

\begin{itemize}
\item Chapt.~\ref{chapt:rbf} is based on \cite{Matveev2011conf9}, and is concerned with reactively avoiding obstacles 
and provably converging to a target using only scalar measurements about the minimum distance to 
obstacles and the angle to the target. The advantage of the proposed method is that it can be analytically shown to have the correct behavior, despite the extremely limited sensor information which was assumed to be available.

\item Chapt.~\ref{chapt:tf} is based on \cite{Savkin2012journ1}, and is again concerned with reactively avoiding obstacles 
and provably converging to a target. However in contrast to
Chapt.~\ref{chapt:rbf}, the sensor information available is the set of visible obstacle edges surrounding the vehicle. The advantage of the proposed method is that it explicitly allows for the kinematics of the vehicle, compared to equivalent control approaches proposed in the literature.

\item Chapt.~\ref{chapt:pf} is based on \cite{Matveev2010conf0}, and is concerned with outlining in more detail 
the path following approach employed in Chapt.~\ref{chap:singlevehicle}. In particular, extensive simulations and experiments
with an agricultural vehicle are documented. The advantage of the proposed method is that it explicitly allows for steering angle limits, and it has been shown to have good tracking performance in certain situations compared to other methods proposed in the literature.

\item Chapt.~\ref{chapt:fbf} is based on \cite{MSavkin2012journ1}, and is concerned with following the boundary of an obstacle using limited 
range information. In contrast to Chapt.~\ref{chap:convsingle}, only a single detection sector 
perpendicular to the vehicle is available. The advantage of the proposed method is that (compared to other equivalent methods proposed in the literature) it provides a single contiguous controller, and is analytically correct around transitions from concave to convex boundary segments.

\item Chapt.~\ref{chapt:ext} is based on \cite{MaHoSa13}, and is concerned with seeking the maximal point of a scalar environmental field. The advantage of the proposed method is that it des not require any type of derivative estimation, and it may be analytically proven to be correct in the case of time--varying environmental fields.

\item Chapt.~\ref{chapt:lst} is based on \cite{MaHoAnSa_sb}, and is concerned with
tracking level sets/isolines of an environmental field at a predefined set-point. The advantage of the proposed method is that it may be analytically proven to be correct in the case of time--varying environmental fields.

\item Chapt.~\ref{chapt:tcp} is based on \cite{matveev2013capt}, and is concerned with 
 decentralized formation control allowing a group of robots to form a 
circular `capturing' arrangement around a given target. The advantage of the proposed method is that it only requires local sensor information, allows for vehicle kinematics and does not require communication between vehicles.

\end{itemize}

Finally, Chapt.~\ref{chap:conclusions} presents a summary of the proposed methods, and outlines
several areas of possible future research. 

Appendix~\ref{apdx:helpexp} presents preliminary simulations
with a realistic helicopter model. Note the helicopter model and the text describing it was 
contributed by Dr. Matt Garratt. 


\chapter{Literature Review}
\label{chap:lit}

In this chapter both local and global approaches are reviewed, together with approaches applicable to
multiple vehicles and moving obstacles. Various types of vehicle and
sensor models are explored, and in the case of moving obstacles, various assumptions about their movement are discussed.
\par
This chapter is structured as follows. In Sec.~\ref{sec:cluttered} the
problem of navigating cluttered environments are described. In Sec.~\ref{sec:mpc} 
MPC-based navigation systems are outlined. In Sec.~\ref{sec:sensor2} methods of
sensor based navigation are introduced; in Sec.~\ref{sec:moving} methods of dealing with moving
obstacles are reviewed. Sec.~\ref{sec:multiple} deals with the case of multiple cooperating vehicles. 
Sec.~\ref{ch1:con} offers brief conclusions.

\section{Exclusions}
\label{sec:exclusion}

Because of the breadth of this research, the following areas are not reviewed, and only a
brief summary is provided where necessary:

\begin{itemize}
  \item \textit{Mapping algorithms}. Mapping is becoming very popular in real-world applications,
  where exploration of unknown environments is required (see e.g.
\cite{Durrant2006journ0}). While they are extremely useful, it seems unnecessary to build a map to
perform local collision avoidance, as this will only generate additional computational overhead.
One exception is the \textit{Bug} class of algorithms (see e.g. \cite{Ng2007journ8, Gabriely2008journ2}), 
which are possibly the simplest examples of convergent navigation.

\item \textit{Path tracking systems}. This continues to be an important, nontrivial problem in
the face of realistic assumptions, and several types of collision avoidance approaches assume the presence of
an accompanying path following navigation law. A review of methods applicable to agricultural vehicles
may be found in Chapt.~\ref{chapt:pf}.
	
  \item \textit{High level decision making}. The most common, classic approach to real world
implementations of autonomous vehicle systems seems to be a hierarchical structure, where a high level planner
provides general directions, and a low level navigation layer prevents collision and attempts to
follows the commands given by the higher layer. In pure reactive schemes, the high level is
effectively replaced with some heuristic. While this type of decision making is required in some
situations to show convergence, it becomes too abstracted from the basic goal of showing collision
avoidance.
Convergence tasks should only be delegated if they can be achieved within the same basic navigation
 framework (see e.g. \cite{Zhu2012journ}).

  \item \textit{Planning algorithm implementations}. Many of the approaches discussed may be used
  with several types of planning algorithms, thus the discussion may be separated. This review effectively
  focuses on the parameters and constraints given to path planning systems, and the subsequent use 
  of the output. Many other surveys
  have explored this topic, see e.g. \cite{Goerzen2009journ6, Lopez2008journ}. However, some local planning approaches
  are reviewed as they are directly relevant to this report.

  \item \textit{Specific tasks (including swarm robotics, formation control, target searching, area patrolling, and
 target visibility maintenance)}. In these cases the primary
  objective is not proving collision avoidance between agents (see e.g. \cite{COVER1}), so approaches to these problems are only
  included in cases where the underlying collision avoidance approach is not documented elsewhere. A review
  of some literature related to tracking environmental fields may be found in Chapts.~\ref{chapt:ext} and \ref{chapt:lst},
  and a review of literature related to formation control may be found on Chapt.~\ref{chapt:tcp}.

  \item \textit{Iterative Learning, Fuzzy Logic, and Neural Networks}. 

  While these are all important areas and are well suited to some applications, and also generate promising 
  experimental results,
it is generally more difficult to obtain guarantees of motion safety when applied directly to vehicle motion (see e.g.
\cite{Gomez2012journ2}). However, these may indirectly be used 
in the form of planning algorithms, which may be incorporated into some of the approaches discussed in this chapter.

\end{itemize}

 \section{Problem Considerations}
\label{sec:cluttered}

In this section, some of the factors which influence the design of vehicle navigation systems are outlined.

\subsection{Environment}
In this chapter, a cluttered environment consists of a $2$ or $3$ dimensional workspace, which
contains a set of simple, closed, untransversable obstacles which the vehicle is not allowed to
coincide with. The area outside the obstacle is considered homogeneous and equally easy to navigate. Examples of
cluttered environments may include offices, man made structures, and urban environments. An example
of classification of objects in an urban environment is available \cite{Douillard2011journ9}.
\par
The vehicle is spatially modeled as either a point, circle, or polygon in virtually
all approaches. Polygons can be conservatively bounded by a circle, so polygonal vehicle shapes are 
generally only required for tight maneuvering around closely packed obstacles, where an enclosing
circle would exclude marginally viable trajectories.

\subsection{Vehicle Kinematics}
\label{sec:vehdyn}

There are many types of vehicles which must operate in cluttered environments; such as ground
vehicles, unmanned air vehicles (UAV's), surface vessels and underwater vehicles. Most vehicles can be
generally categorized into three types of kinematic models -- \textit{holonomic}, \textit{unicycle} and \textit{bicycle} -- where the
 differences are characterized by different turning rate constraints. Reviews of different vehicle models
 are available, see e.g. \cite{Gracia2007journ2, Gracia2007journ2a, Kozlowski2009book0, Micaelli1993book3}.
In this chapter, the term \textit{dynamic} is used to describe models 
based on the resolution of physical forces, while the term \textit{kinematic} describes models based on  abstracted control inputs.

\begin{itemize}
  \item \textit{Holonomic kinematics}. In this report, the term holonomic is used to describe linear models which have equal control capability in any direction.
  Holonomic kinematics are encountered on helicopters, and certain types of wheeled robots equipped
  with omni-directional wheels.
  Holonomic motion models have no notion of body orientation 
for the purposes of path planning, and only the Cartesian coordinates are considered. However,
orientation may become a consideration when applying the resulting navigation law to real vehicles
(through this is decoupled from planning).
  
  \item \textit{Unicycle kinematics}. These describe vehicles which are associated with a particular angular
  orientation, which determines the direction of the velocity vector. Changes to the orientation
  are limited by a turning rate constraint.
  Unicycle models can be used to describe various types of vehicle, such as differential drive
  wheeled mobile robots and fixed wing aircraft, see e.g. \cite{Manchester2004journ9, Manchester2006journ2}.
  
  \item \textit{Bicycle kinematics}. These describe a car-like vehicle, which has a steerable front wheel separated
   from a fixed rear wheel. Kinematically this implies the maximum turning rate is
  proportional to
the vehicles speed. This places an absolute bound on the curvature of any path the vehicle may
follow regardless of speed. This constraint necessitates higher order planning to successfully navigate confined
environments.
\end{itemize}

It should be noted that nonholonomic constraints are in general a limiting factor only at low speed --
for example, realistic vehicles would likely be also subjected to absolute acceleration bounds.
More complex kinematics are also possible, but uncommon. 
In addition the these basic variants of kinematics, the associated linear and angular variables may
 be either \textit{velocity
controlled} or \textit{acceleration bounded}. Vehicles with acceleration bounded control inputs are
in general much harder to navigate; velocity controlled vehicles may stop instantly at any time if
required. 
\par
When predicting an vehicle's actual motion, nominal models are invariably subject to disturbance. The type of disturbance which may be modeled depends on
the kinematic model:

\begin{itemize}
  \item \textit{Holonomic models}. Disturbance models commonly consist of bounded additions 
  to the translational control inputs, see e.g. \cite{Richards2006journ1}.
  \item \textit{Unicycle models}. Bounded addends to the control inputs can be combined with a bounded difference between the 
  vehicle's orientation and actual velocity vector, see e.g. \cite{Lapierre2008journ2}.
More realistic models of differential drive mobile robots are also available, which are based on modeling 
wheel slip rates (see e.g. \cite{Albagul2004conf, Balakrishna1995journ}).
\item \textit{Bicycle models}. Disturbance can be modeled as slide slip angles on the front 
and rear wheels (see e.g. \cite{Matveev2010conf0}). Alternatively, more realistic disturbance models of car-like 
vehicles are available, which include factors such as suspension and type adhesion 
(see e.g. \cite{Bevan2007journ6, Yoon200journ}).
\end{itemize}

Vehicles with bicycle kinematics or vehicles with minimum speed constraints will be subject
to absolute bounds on their path curvature. This places some global limit on the types of
environments they can successfully navigate through, see e.g. \cite{Sergey2005conf1, Bicchi1996journ0}.
When lower bounds on allowable speed are present, the planning system is further complicated. 
For example, instead of halting, the vehicle must follow some holding pattern at the termination of a trajectory.

\subsection{Sensor Data}
\label{sec:sensor2data}

Most autonomous vehicles must base their navigation decisions on data reported by on-board sensors,
which provide information about the vehicles immediate environment. The main types of sensor model
are listed as follows:

\begin{itemize}
  \item \textit{Abstract sensor models}. This model informs the navigation law whether a given point lies within the
  obstacle set. Usually any occluded regions, without a line-of-sight to the vehicle, are considered to be
  part of the obstacle.
   Through this set membership property is impossible to determine precisely using a physical sensor, currently some 
   \textit{Light Detection and Ranging} (LiDAR) sensors have accuracy high enough for any sampling effects to be of minor concern.
   However, when lower resolutions are present, this model may be unsuitable for navigation law design. 
   \item \textit{Ray-based sensor models}. These models inform the navigation law of the distance to the obstacle in a finite number
   of directions around the vehicle, see e.g. \cite{Travis2005conf5, Moghadam2008conf9,
   Lapierre2007journ7}.
   This is a more physically realistic model of laser based sensors compared to the abstract sensor model,
   and may be suitable for determining the effect of low resolution sensors. A reduced version of
   this model is used in some boundary following applications, where only a single detection ray in a fixed
   direction (relative to the vehicle) is present.
  \item \textit{Minimum distance measurements}. This sensing model reports the distance to the nearest obstacle point.
   This may be realized by
  certain types of wide aperture acoustic or optic flow sensors. Using this type of measurement necessarily leads to less efficient movement
  patterns during obstacle avoidance, i.e. it is not immediately clear which side of the vehicle the obstacle is 
  on (see e.g. \cite{Matveev2011journ2}).
 \item \textit{Tangent sensors}. This sensor model reports the angles to visible edges of an obstacle as seen by the vehicle
  (see e.g. \cite{Shi2010journ9, Suri2008journ2}). This
  can be realized from a camera sensor, provided a method of detecting obstacle edges from a video stream is available 
  (see e.g. \cite{Huang2006journ7}).
  \item \textit{Optic Flow Sensors}. This model reports the average rate of pixel flow across a camera sensor 
  (see e.g. \cite{Bonin2008journ0, Green2008journ1}). While these types of sensor are very compact, unfortunately 
  rigorously provable collision avoidance does not seem to be currently possible.

\end{itemize}

There are a large number of ways in which noise and distortion may be compensated for in these 
models, and these tend to be quite specific to individual approaches.

\subsection{Optimality Criteria}
There are several different methods of preferentially choosing one possible trajectory over another.
Many path optimization algorithms may be implemented with various such measures or combinations of measures.
Common possibilities are listed below:

\begin{itemize}
  \item \textit{Minimum path length}. This is used in the majority of path planning schemes as it can be decoupled
  from the achievable velocity profile of the the vehicle.
  For moving between two configurations without obstacles, the classic result of Dubins describes the optimal motion of
  curvature bounded vehicles \cite{DUB57}. In this case, the optimal
  path consists of a sequence of no more than three maximal turns or straight segments. 
  Other similar results are available for vehicles with actuated speed \cite{Reeds1990journ}, 
  and for velocity controlled, omni-directional vehicles \cite{Balkcom2006journ5}. 
However these results are of little direct use in path planning, since obstacles 
have a complex effect on any optimal path. When acceleration constraints are absent, the minimum length
path may be constructed from the \textit{Tangent Graph} of an obstacle set, see e.g. \cite{Savkin2012journ1, LA92}.

  \item \textit{Minimum time}. Calculating the transversal time of a path depends on the velocity profile of the vehicle, 
  and thus includes kinematic (and possibly dynamic) constraints. In most situations it would be more appropriate
  than minimum length for selecting the trajectories that complete tasks in the most efficient manner. It is often used in MPC-based approaches, see e.g. \cite{Richards2006journ1}.
  
  \item \textit{Minimum control effort}. This may be more suitable for vehicles operating in limited
  energy environments, e.g. spacecraft or passive vehicles, however it is invariably combined with another measure for non-zero movement. Another formulation in the same vein, \textit{minimum wheel rotation}, applies to differential drive
  wheeled mobile robots. In most cases is only subtlety different from the minimum length
  formulation; however it may perform better in some situations, especially when fine movements are required \cite{Chitsaz2009journ8}.
  
  \item \textit{Optimal surveillance rate}. In unknown environments it may be better to select
  trajectories which minimize the occluded part of the environment (see e.g. \cite{Trevai2007journ6}). In cases where
  occluded parts of the environment must be treated as unknown dynamic obstacles, this could allow a more efficient transversal, 
  through it would unavoidably rely on stochastic
  inferences about the unknown portion of the workspace. This may be an interesting area of future research.
\end{itemize}

Other examples of requirements that can be applied to trajectories include higher order 
curvature rate limits, which may be useful to produce smoother trajectories (see e.g.
\cite{Armesto2012book}).

\subsection{Biological Inspiration}

Researchers in the area of robot navigation in complex environments find
much inspiration from biology, where the problem of
controlled animal motion is considered. 
This is prudent since biological systems are highly efficient and refined, 
while the equivalent robotic systems are in relative infancy.
Animals,
such as insects, birds, and mammals, are believed to use
simple, local motion control rules that result in remarkable
and complex intelligent behaviours. Therefore, biologically inspired or biomimetic algorithms
of collision free navigation play an important part in this research field.

In particular, ideas of the navigation along an equiangular spiral and the associated local
obstacle avoidance strategies have been proposed which 
are inspired by biological examples \cite{Teimoori2010journ1a,Teimoori2010journ1,Savkin2012conf0}).  It has been observed that peregrine falcons, which
are among the fastest birds on the earth, plummet toward their
targets at speeds of up to two hundred miles an hour along
an equiangular spiral \cite{Tucker2001journ}. Furthermore, in biology, a similar obstacle avoidance strategy 
is called `negotiating obstacles with constant curvatures' (see
e.g. \cite{Lee1998journ}). An example of such a movement is a squirrel running
around a tree. 
These ideas in reactive collision avoidance robotic systems  are further discussed in Sec.~\ref{sec:bf}. Furthermore, the sliding mode
control based methods of obstacle avoidance discussed in Sec.~\ref{sec:bf} are also inspired by biological examples such as the near-wall behaviour of a cockroach \cite{CaJh99}. Another example 
is the Bug family algorithms which are also inspired by bugs
behaviour on crawling along a wall.

Optical flow navigation is another important class of biologically inspired navigation methods.
The remarkable ability of honeybees and other insects like
them to navigate effectively using very little information is a
source of inspiration for the proposed control strategy. In particular,
the use of optical-flow in honeybee navigation has been explained, where
a honeybee makes a smooth landing on a surface without the
knowledge of its vertical height above the surface \cite{SRI00}. Analogous
to this, the control strategy we present is solely based on instantaneously
available visual information and requires no information on the
distance to the target. Thus, it is particularly suitable for robots
equipped with a video camera as their primary sensor (see e.g. \cite{Low2007journ}) 
As it is commonly observed in 
insect flight, the navigation command is derived from the average rate of pixel flow across a camera sensor (see e.g. \cite{Bonin2008journ0, Green2008journ1}). This is further discussed in Sec.~\ref{sec:opflow}.

Many ideas in multi-robot navigation are also inspired by biology,
where the problem of animal aggregation is central in both
ecological and evolutionary theory. Animal aggregations,
such as schools of fish, flocks of birds, groups of bees,
or swarms of social bacteria, are believed to use simple,
local motion coordination rules at the individual level that
result in remarkable and complex intelligent behaviour at the
group level (see e.g. \cite{FGLO99,BWF11}). Such intelligent behaviour
is expected from very large scale robotic systems.  Because of
decreasing costs of robots, interest in very-large-scale robotic
systems is growing rapidly. In such systems, robots should
exhibit some forms of cooperative behaviour.  We discuss it further in Sec.~\ref{sec:multiple}.

There is also some evidence that approaches resembling Model Predictive Control (MPC) are used
by higher animals to avoid obstacles \cite{Ahmadi2007journ1}. It seems natural
to achieve collision avoidance using some type of planning into the future, and
MPC-based navigation laws are discussed in Sec.~\ref{sec:mpc}.

Neural networks and fuzzy logics approaches are often identified as 
biologically inspired, see e.g. \cite{Mondada1993conf,Lewinger2006conf}. However,
implementations of these approaches are not generally concerned with rigorously showing the desired behaviour, and 
as such are not included in this survey (see Sec.~\ref{sec:exclusion}).

\subsection{Implementation Examples}

There are many review of current applications and implementations of real world vehicles, see e.g.
\cite{ROB:ROB20414}. An exhaustive list of 
reported applications would be excessive, however one particular application is highlighted.
\par
Semi-autonomous wheelchairs are a recent application in which a navigation law must be designed to
prevent collisions while taking high-level direction inputs from the user, see e.g. \cite{Boquete1999journ1, Wang2012journ}. 
In this case, a fundamental concern for these intelligent
wheelchairs is maintaining safety, thus the methods described in this review are highly relevant.
Several original collision avoidance
approaches were originally proposed for wheelchair applications, 
see e.g. \cite{Wang2012journ} (these are also discussed in Sec.~\ref{sec:movreact}).

\section{Summary of Methods}

\begin{figure}[ht]
  \centering 
  \includegraphics[width=\columnwidth]{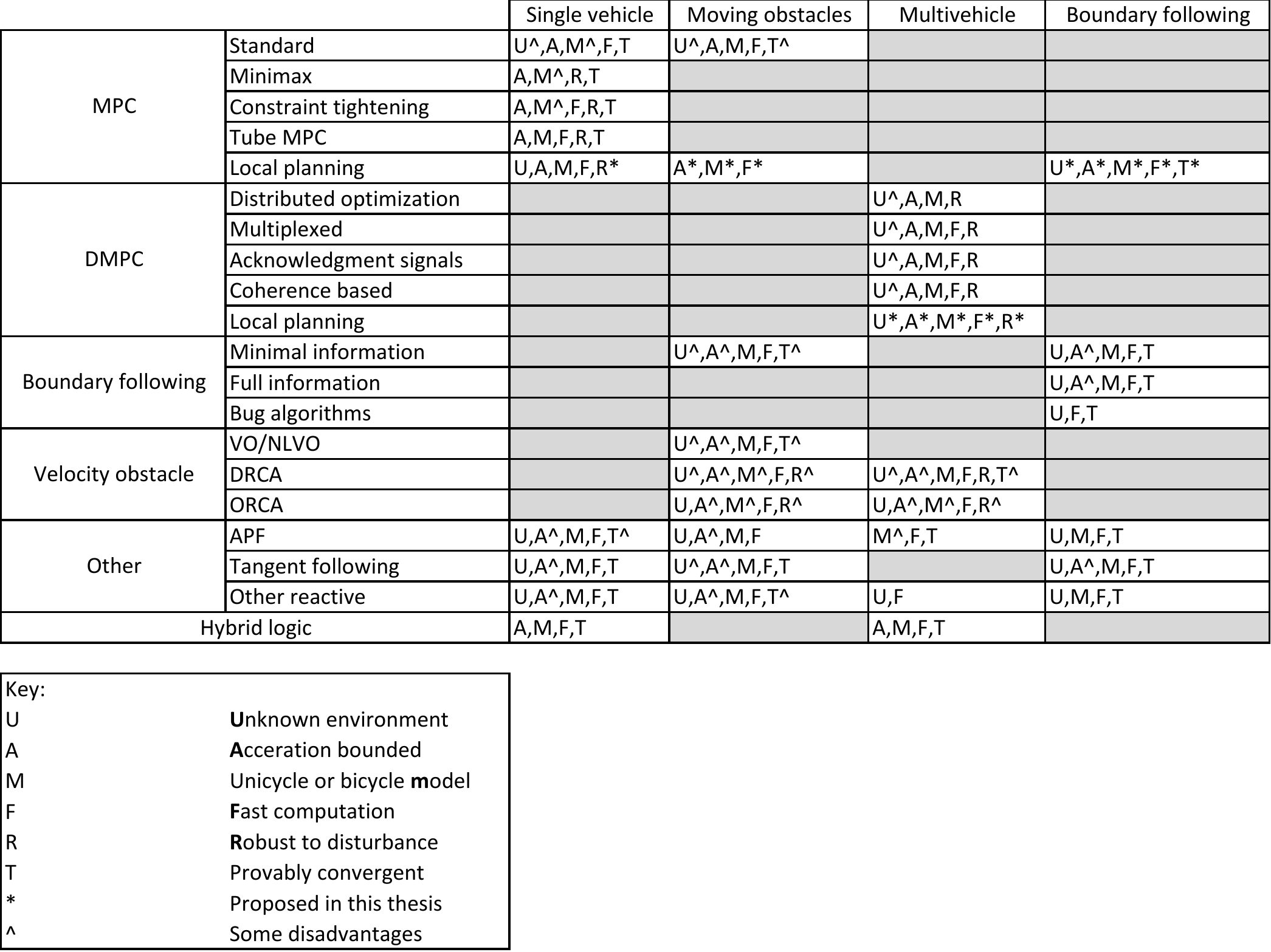}
  \caption{Summary of the methods reviewed.}
\label{fig:broad}
\end{figure}

A very broad summary of the methods considered in this review are listed in Fig.~\ref{fig:broad}, where the availability of certain traits is shown.

These methods will be discussed in the remainder of this chapter.

\clearpage \section{Model Predictive Control}

\label{sec:mpc}

If a obstacle-avoiding trajectory is planned off-line, there are many examples of path following
systems which are able to robustly follow it, even if subjected to bounded disturbance. However the lack of
flexibility means the environment would have to be perfectly known in advance, which is not conducive 
to on-line collision avoidance.
\par
\textit{Model Predictive Control} (MPC)\footnote{Equivalent to \textit{Receding Horizon Control} (RHC)}, is
increasingly being applied to vehicle navigation problems. It is useful as it combines path
planning with on-line stability and convergence guarantees, see e.g.
\cite{Richards2007journ, Mayne2003journ5}.
This is basically done by performing the path planning process at every time instant, then applying
the initial control related to the chosen trajectory to the vehicle. In most cases a partial path is planned (such that it would not usually arrive at the target),
which terminates with an invariant vehicle state. The navigation law would attempt to minimize some `cost-to-go' or navigation function corresponding to the target.
\par
In recent times, MPC has been increasingly being applied to vehicle navigation problems, and it seems
to be a natural method for vehicles to navigate. Note discussion of MPC-based approaches applicable
to unknown environments is reserved until Sec.~\ref{sec:sbpp}.

\subsection{Robust MPC}

The key advantage of MPC lies with its robust variants, which are able to account for set bounded disturbance (and are
the most useful for vehicle navigation). These 
can be categorized into three main categories:

\begin{itemize}
\item \textit{Min-max MPC}. In this formulation, the optimization is performed with respect
 to all possible evolutions of the
disturbance, see e.g. \cite{Scokaert1998journ}. While it is the optimal solution to linear robust control
problems, its high computational cost generally precludes it from being used for vehicle navigation.

\item \textit{Constraint Tightening MPC}. Here the state constraints are dilated by a given margin
so that a trajectory can guaranteed to be found, even when disturbance causes the state to evolve towards the
constraints imposed by obstacles (see e.g. \cite{Richards2006conf9, Richards2006journ1, Kuwata2007journ7}). The basic argument shows a future viable
trajectory exists using a feedback term, through a feedback input is not directly used for
updating the trajectory. This is commonly used for vehicle navigation problems -- for example a system has been described where an
obstacle avoiding trajectory is found based on a minimization of a cost functional compromising
the control effort and maneuver time \cite{Richards2006journ1}. In this case, convergence to the target and the
ability to overcome bounded disturbances can be shown. 

\item \textit{Tube MPC}. This uses an independent nominal model of the system, and employs a feedback
system to ensure the actual state converges to the nominal state (see e.g. \cite{Langson2004journ}). In
contrast, the constraint tightening system would essentially take the nominal state to be the actual state at
each time step. This formulation is more conservative than constraint tightening, since it wouldn't
take advantage of favorable disturbance. Thus it doesn't offer significant benefits for
vehicle navigation problems when a linear model is used.
However, it is useful for robust nonlinear MPC (see e.g. \cite{Mayne2011journ1}), and
problems where only partial state information is available (see e.g. \cite{Scholte2008journ8}). Also, the approach proposed
in this report, along with any approach which includes path
following with bounded deviation (see e.g. \cite{Defoort2009journ7}), is somewhat equivalent to tube MPC.
\end{itemize}

For robust MPC, the amount of separation required from the state constraints on an infinite horizon is determined
by the \textit{Robustly Positively Invariant} (RPI) set, which is the set of all possible state deviations 
that may be introduced by disturbance while a particular disturbance rejection law is operating. 
Techniques have been 
developed to efficiently calculate the smallest possible RPI set (the minimal RPI set) \cite{Rakovic2005journ0}.

If disturbance is Gaussian rather than set bounded, the MPC	
vehicle navigation problem may be reformulated stochastically so the overall risk of collision is bounded to
an arbitrary level, see e.g. \cite{Blackmore2011journ9, Du2012journ5}. This may be an interesting area of future research.

\subsection{Nonlinear MPC}

The current approaches to MPC-based vehicle navigation generally rely on linear kinematic models, usually with double
integrator dynamics. While many path planning approaches exist for vehicles with nonholonomic
kinematics, it is generally harder 
to show stability and robustness properties
\cite{Magni2009book0}. Approaches to robust nonlinear MPC are generally of the tube MPC type
\cite{Mayne2011journ1}.

In these cases, a nonlinear trajectory tracking system can be used to ensure the actual state 
converges to the nominal state. A proposition has been made to also use sliding mode control laws
for the auxiliary
system \cite{Rubagotti2010journ0}. Sliding mode control was employed in this report, and through such systems typically require 
continuous time analysis, disturbance rejection properties are typically easier to show.

In terms of vehicle navigation problems, examples of MPC which apply unicycle
kinematics while having disturbance present have been proposed, see e.g.
\cite{Defoort2009journ7, Defoort2009journ1}. However it seems more general
applications of nonlinear MPC to vehicle navigation problems should be possible; for example
 in this report a new control method employing tube MPC principles is proposed.

There are other methods in which MPC may be applied to vehicle navigation problems other than
performing rigorously safe path planning. In some cases the focus is shifted towards controlling
vehicle dynamics, see e.g. \cite{Park2009journ0, Yoon200journ, Gonzalez2011journ}. These use a realistic 
vehicle model during planning, and are able to give good practical
results, through guarantees of safety are currently easier with kinematic models. In other cases MPC may be used to
regulate the distance to obstacles, see e.g. \cite{Shim2006journ5}. However, this discussion of this type of method
 is reserved until Sec.~\ref{sec:fibf}.

\subsection{Planning Algorithms}

Global path planning is a relatively well studied
research area, and many thorough reviews are available see e.g. \cite{Goerzen2009journ6, Lopez2008journ}. MPC may be 
implemented with a number of different path planning algorithms. The main relevant measure of algorithm quality is \textit{completeness}, which indicates 
whether calculation of a valid path can be guaranteed whenever one exists. Some common
global path planning algorithms are summarized:

\begin{itemize}
  \item  \textit{Rapidly--Exploring Random Trees}.
  Creates a tree of possible actions to connect initial and goal configurations (see e.g. \cite{Diankov2007conf8, Karaman2011journ3}). Some variants are
  provably asymptotically optimal \cite{Karaman2011journ3}.
  
\item \textit{Graph search algorithms}. Examples include A* (see e.g. \cite{Sathyaraj2008journ}) and D* (see
e.g. \cite{Koenig2005journ7}). Most methods hybridize the environment into a square graph, with the
search calculating the optimal sequence of node transitions. However in many approaches the cells need not be
square and uniform, see e.g. \cite{Kallem2011journ9, Belta2005journ5}.

\item \textit{Optimization of parameterized paths}. Examples include Bezier curves
\cite{Skrjanc2010journ4}, splines \cite{Lau2009conf1} and polynomial basis functions \cite{Zh2004journ7}. While these are inherently smoother,
showing completeness may be more difficult in some situations.

\item \textit{Mathematical programming and optimization}. This usually
is achieved using \textit{Mixed Integer Linear Programming} (MILP) constraints to model obstacles as multiple convex polygons \cite{Abichandani2012conf}. 
Currently this is commonly used for MPC approaches.

\item \textit{Tangent Graph based planning}.
This limits the set of trajectories to cotangents between obstacles and obstacle boundary segments,
from which the minimum length path being found in general  \cite{Savkin2012journ1, Tovar2007journ3}.
The problem of shortest path planning in a known environment for unicycle-like vehicles with a hard constraint
on the robot's angular speed was solved in \cite{Savkin2012journ1}.
It is assumed 
that the environment consists a number of possibly non-convex obstacles
with a constraint on their boundaries curvature and a steady
target that should be reached by the robot. It has been proved
the shortest (minimal in length) path consists of edges
of the so-called tangent graph.
Therefore, the problem of the shortest path planning is
reduced to a finite search problem.

\item \textit{Artificial Potential Field Methods}. These methods are introduced in Sec.~\ref{sec:apfm}, and are ideally suited to online reactive navigation of vehicles. These can also be used as path planning approachs, essentially by solving the differential equations corresponding to the closed loop system (see e.g. \cite{MPCNF}). However, these trajectories would not be optimal and have the same drawbacks as the original method in general. 

\item \textit{Evolutionary Algorithms, Simulated Annealing, Particle Swarm Optimization}. These are
based on a population of possible trajectories, which follow some update rules until the optimal path is reached 
(see e.g. \cite{Zheng2005journ6, Besada2010journ3}). However these approaches seem to
be suited to complex constraints, and may have slower convergence for normal path planning problems. 

\item \textit{Partially Observable Markov Decision Processes}.
This calculates a type of decision tree for different realizations of uncertainty, and uses 
probabilistic sampling to generate plans that may be used for navigation 
over long time frames (see e.g. \cite{Kurniawati2011journ2}). However
this does not seem necessary for MPC-based navigation frameworks.
\end{itemize}

 \section{Sensor Based Techniques}

\label{sec:sensor2}

In comparison to path planning based approaches, sensor based navigation techniques only have limited 
local knowledge about the obstacle, similar to what would be obtained from range finding sensors, cameras,
or optic flow sensors. Reactive schemes are a subset of these which may be interpreted as a mapping between the current sensor
state and the actuator outputs; thus approaches employing even limited memory elements would not be considered reactive. 

A method for constraint-based task specification has been proposed for sensor-based vehicle systems \cite{De2007journ2}. 
This may provide a fixed design process for designing reactive navigation systems (an example 
of contour tracking is given), and this concept may be an interesting area of future work.

\subsection{Boundary Following}

\label{sec:bf}

Boundary following is a direct subproblem of obstacle avoidance, and in most cases a closed
loop trajectory bypassing an obstacle can be segmented into `boundary following' and `pursuit' actions, even if
this choice is not explicitly deliberated by the navigation law. Boundary following by itself also has
many direct uses such as border patrol, terrain tracking and structure monitoring; for application examples see e.g.
\cite{Girard2004conf0}.

\subsubsection*{Distance Based}

\label{sec:bfdis}

In many approaches, boundary following can be rigorously achieved by only measuring the minimum
distance to the obstacle, see e.g. \cite{Matveev2011journ2, Matveev2011conf9}. For example, a 
navigation strategy has been proposed using a feedback controller based on the
minimum obstacle distance, and is suitable for guiding nonholonomic vehicles traveling at
constant speed \cite{Matveev2011journ2}. In Chapt.~\ref{chapt:rbf} a similar feedback strategy is proposed which only
 requires the rate of change of the distance to the obstacle as input.

Other approaches have been proposed which use a single obstacle distance measurement at a specific
angle relative to the vehicle \cite{Toibero2009journ7, MSavkin2012journ1, 
Teimoori2010journ1, Teimoori2010journ1a, Kim2009journ7}. These can be classified based on the 
required measurement inputs; navigation can be based purely on the length of the detection ray 
\cite{Teimoori2010journ1, Teimoori2010journ1a}; or additionally based on the tangential angle of the
obstacle at the detected point \cite{Toibero2009journ7, MSavkin2012journ1}, or additionally based on
estimation of the boundary curvature \cite{Kim2009journ7}. Additional information would presumably 
result in improved behavior, through methods employing boundary curvature may be sensitive to noise, and performance
may degrade in such circumstances. Several of these methods are in the realm of switched controllers, for which
rigorous theoretical results are available \cite{Evans1,Evans2}.
Unfortunately, impartial comparisons of the closed loop performance of these approaches would be difficult.

Some other methods using similar assumptions are focused on following straight walls, see e.g.
\cite{Yata1998conf6, Bemporad2000journ0, Carelli2003journ0, Huang2009journ2}. However it seems,
at least theoretically, navigation laws capable of tracking contours are more general and therefore superior.

In most these examples the desired behavior can be rigorously shown. However, the common limitation
is that the vehicle must travel at constant speed, and this this speed must be set conservatively
according on the smallest feature of the obstacle. In some cases simple heuristic can partially
solve this problem; by instructing the vehicle to instantly stop and turn in place if the obstacle
distance becomes too small, collision may be averted \cite{Toibero2009journ7}.

\subsubsection*{Sliding Mode Control}
\label{sec:slidmode}

Special consideration should be given to sliding mode control based navigation approaches,
which are increasingly being applied to vehicle navigation problems where limited sensor information is available, 
see e.g. \cite{MSavkin2012journ1,Teimoori2010journ1, Teimoori2010journ1a,  Matveev2011conf9, Matveev2011journ4, Matveev2011journ9}. In this context,
sliding mode control consists of a discontinuous, switching navigation law which allows rigorous mathematical 
analysis, and has the additional benefit of having a high
resistance to noise, disturbance and model deviation, see e.g. \cite{Utkin1992book1}. In the context of collision avoidance, sliding mode control approaches
usually are designed as boundary following approaches, through they 
have also been applied to the avoidance of moving obstacles, and show promising results in that area (see Sec.~\ref{sec:movreact}).

\subsubsection*{Full Information Based}

\label{sec:fibf}
In situations where more information about the obstacle is available, a clearer view of
the immediate environment can be recreated. This means more informed navigation decisions may be able to be made. This
can lead to desirable behaviors, such as variable speed and offset distance from the obstacle. It
also allows us to loosen some of the assumptions on the obstacle shape and curvature. An example
of such behavior may be slowing down at concavities of a boundary and speeding up otherwise, or
completely skipping concavities of sufficiently small size that serve only to introduce
singularities into the motion. \cite{Matveev2011journ2, Kim2009journ7}.

 One such approach using abstract obstacle information is the VisBug class of algorithms, which navigates towards a visible edge of an
obstacle inside the detection range (see e.g. \cite{Lumelsky1990journ3, Ng2007journ8, Langer2007conf}).
However, these algorithms are concerned with the overall strategy, and are not concerned with
details relating to vehicle kinematics or the sensor model.
 Several approaches have been able to account for the vehicle dynamics, but still have inadequate
models of the vehicle sensor. This is similar to the \textit{joggers problem}, whose solution involves
ensuring safe navigation by ensuring the vehicle can stop in the currently sensed obstacle free set \cite{Shkel1997journ5}.
However, an abstract sensor model was used, which presumes the vehicle has continuous knowledge about
the obstacle set. A navigation approach which achieves boundary following by picking
\textit{instant goals} based on observable obstacles has been proposed \cite{Ge2005journ0}. A ray based sensor model is used, through 
a velocity controlled holonomic model is assumed.
Instant goals have also been used where allowance is made for the vehicle
kinematics, however in this case a ray-based obstacle sensor model was not used \cite{Shuzhi2007journ1}.

In Chapt.~\ref{chap:convsingle} a novel MPC-based approach to boundary following is proposed, which generates avoidance constraints
and suitable target points to achieve boundary following. This is an interesting, original application of MPC, and may be useful for
other types of sensor based navigation problems.

\subsubsection*{Bug Algorithms}
The boundary following navigation laws mentioned previously may perform target-convergent navigation
when coupled with high level behavior resembling Bug algorithms (see e.g. \cite{Matveev2011journ2,
Matveev2011conf9, Mastrogiovanni2009journ2}). Bug algorithms achieve global convergence by switching between
`boundary following' and `target pursuit' modes. By combining these systems, the main additional complexity 
involves finding and analyzing the conditions for switching between the two modes. 
While these can be proven to converge to the target, it is
important to note pure reactive navigation laws will fundamentally be subjected to local
minima problems and will not lead to provable target convergence -- this is
impossible in general with a reactive, deterministic methods \cite{Matveev2011conf9}.
A number of heuristics exist to prevent these, through they would not be classified as reactive (see e.g. \cite{Ordonez2008journ7}).

\subsection{Non--Trajectory Based Obstacle Avoidance}
\label{sec:apfm}

In this section, methods which neither explicitly generate a path nor explicitly perform boundary following are described. This 
includes for example potential field based methods.

Many approaches to this particular problem assume holonomic velocity controlled vehicles. However, this turns out not to be a severe limitation
since methods are available for extending these basic navigation laws to account for arbitrary
 dynamics (including acceleration constraints) are available, see e.g. \cite{Minguez2003conf, Minguez2009journ, Blanco2008journ}. 
This method is based on a coordinate transformation, which effectively provides a zone around the vehicle that contains all perturbations
introduced by the dynamics. This method may be applied to a range of navigation approaches, through it may be conservative in some situations. 
Alternatively a method has been 
proposed which guarantees collision avoidance
 simply by ensuring the distance to obstacles is always greater than the
stopping distance \cite{Manor2012conf}. Through more conservative, this approach may be useful in cases where little is known about the
vehicle model. 

\subsubsection*{Artificial Potential Field Methods}

A classic approach to reactive collision avoidance is to construct a virtual potential field that
causes repelling from obstacles and attraction to the target. These are termed \textit{Artificial Potential Field} (APF)
methods, and this continues to be an active area of research. Several improvements are listed as follows:

\begin{itemize}
  \item \textit{Unicycle kinematics}. Performance can be improved on vehicles with unicycle type kinematics. Specifically, this involves
   moving the vehicles reference point slightly away from the center of the vehicle, see e.g.
 \cite{Valbuena2012journ6, Ren2008journ5}.
 
  \item \textit{Local minima avoidance}. The shape of the potential field can be designed to flow around obstacle concavities; these
are termed harmonic potential fields and provide better performance with local minima, see e.g.
\cite{Masoud2012journ7, Masoud2010journ1}. However, it seems impossible to deterministically avoid local minima using reactive algorithms.

\item \textit{Closed loop performance}. Alteration to the shape of the
potential field leads to an improvement to the closed loop performance, see e.g. \cite{Kim2006journ8,
Cifuentes2012journ1}. However in general, the closed loop trajectories of APF based methods would not be optimal. 
Additionally, reductions of oscillation in narrow corridors may be achieved, 
 see e.g. \cite{Ren2006journ5,Ren2008journ5}. 
 
 \item \textit{Limited obstacle information}. Examples are available where only the
nearest obstacle point is available \cite{Chang2003conf8}. Several approaches assume global knowledge about the 
workspace, and thus would not suitable for sensor based navigation.

\item \textit{Actuator constraints}. Examples which focus on satisfying actuator constraints are
 also available, see e.g. \cite{Galicki2009journ2, Loizou2008journ}.
 
\end{itemize}

\subsubsection*{Tangent Based Methods}
\label{sec:tanmethod}

Many approaches can be classified as being tangent based, in the sense that they generally only consider
motions towards the tangents of obstacles. It has been shown the distance optimal transversal of a
cluttered environment can be taken from elements of the tangent graph, which is the set of all
tangents between objects, see e.g. \cite{Savkin2012journ1, LA92}. In these cases a method
of probabilistically convergent on-line
 navigation involves randomly choosing
tangents to travel down (see e.g. \cite{Savkin2012journ1}),
 or by use of the deterministic \textit{TangentBug} algorithm (see e.g. \cite{Kamon1998journ9}).

Tangent events can be detected from a ray based sensor model (see e.g.
 \cite{Shi2010journ9}) or by processing data from a camera sensor
 (see e.g. \cite{Huang2006journ7}). This results in an abstract tangent sensor which reports the angle to tangents around the
 vehicle.
A common method of achieving obstacle avoidance is to maintain a fixed angle
 between the tangent and the vehicles motion, see e.g. \cite{Sharma2012journ0, Huang2006journ7}.

\subsubsection*{Optic Flow Based Methods}. 
\label{sec:opflow}

This type of navigation is inspired by models of
insect flight, where the navigation command is derived from the average rate of pixel flow across a camera sensor 
(see e.g. \cite{Bonin2008journ0, Hrabar2005conf7, Green2008journ1, Muratet2005journ6}). 
From this rate of pixel flow, a navigation command may be reactively expressed, and good experimental results have been achieved. 
While this method has the advantage of using a compact sensor and requiring low computational overhead,
general mathematical analysis of such navigation laws for showing collision
avoidance seems more difficult than the equivalent analysis for range based sensors.

\subsubsection*{Other Reactive Methods}

There are some other variations of approaches which achieve collision avoidance:

\begin{itemize}
  \item The \textit{Safe Maneuvering Zone} is suited for kinematic unicycle model with
saturation constraints, when the nearest obstacle point is known \cite{Lapierre2012journ0}.
This is somewhat similar the \textit{Deformable Virtual Zone}, where the navigation
is based on a function of obstacle detection ray length \cite{Lapierre2007journ7}, through collision
avoidance is not explicitly proven.
\item  The \textit{Vector Field Histogram} directs the vehicle towards
sufficiently large gaps between detection rays \cite{ulrich_vfh*:_2000}.
The \textit{Nearness Diagram} is an improved version
which employs a number of behaviors for a number of different situations, providing good performance
even in particularly cluttered environments (see e.g. \cite{Minguez2004journ, Minguez2005journ}).
\item A collision avoidance system based on
 MPC has been proposed and shown to
successfully navigate real-world helicopters in unknown environments based on the nearest obstacle point
within the visibility radius \cite{Shim2006journ5}. However this is less concerned with a proof of collision avoidance,
and more with controlling vehicle dynamics.
\item A different class of navigation law is based on the \textit{Voronoi Diagram}, which essentially describes the set of
points equidistant from adjacent obstacles. In general it leads to longer paths than the tangent graph,
through it represents the smallest set of trajectories which span the free space in an environment. Navigation laws have
been developed to equalize the distance to obstacles, when a velocity controlled unicycle kinematic model is assumed (see e.g. \cite{Victorino2003journ}).
\end{itemize}

\subsection{Sensor Based Trajectory Planning}

\label{sec:sbpp}

Trajectory planning using only sensor information was originally termed the \textit{joggers problem}, since the
vehicle must always maintain a path which brings it to a halt within the currently sensor area, see e.g.
\cite{Shkel1997journ5, Alvarez1998conf4}.

The classic \textit{Dynamic Window} (see e.g. \cite{Fox1997journ4, Ogren2005journ3, Ogren2002journ5})
and \textit{Curvature Velocity Method}
(see e.g. \cite{Fernandez2004journ7, Shi2010journ9}) can be interpreted as a planning algorithm with a
prediction horizon of a single time step \cite{Ogren2005journ3}.
To this end, the range of considered control inputs is limited to those bringing the vehicle to a
halt within the sensor visibility area, using only circular paths. This can also be easily extended to other
vehicle shapes and models \cite{Schlegel1998conf}. Additionally, measures are available
which may reduce oscillatory behavior \cite{Stachniss2002conf}. A wider range of possible trajectory shapes
has also been considered, through it is unclear whether it significantly improves closed loop performance \cite{Blanco2008journ}. The \textit{Lane-Curvature Method} (see e.g.
\cite{Ko1998conf}), and the \textit{Beam-Curvature Method} (see e.g.
\cite{Fernandez2004journ7, Shi2010journ9}) are both variants based on a slightly different trajectory selection process. However, 
in all these cases a similar class of possible trajectories is employed.

In all these cases the justification for collision avoidance is essentially the same argument (the
vehicle can stop while moving along the chosen trajectory). The differences in performance are mainly
heuristic, and in particular they do not fully account for disturbance and noise. However, an approach similar 
to the dynamic window was extended to cases where safety constraints must be generated by
processing information from a ray-based sensor model \cite{Horn2010conf}. 

MPC-type approaches have previously been used to navigate vehicles in unknown
environments, see e.g. \cite{Krishnamurthy2007journ1, Brooks2009conf0, Yang2010journ7}. In most approaches the MPC 
navigation system is combined with some type of mapping algorithm; however, these often lack
the rigorous collision avoidance guarantees normally provided in full-information MPC approaches. 
In Chapt.~\ref{chap:singlevehicle}, a trajectory planning method is proposed, and while it is somewhat similar to the Dynamic Window class of approaches,
it implements a control framework somewhat similar to robust MPC. Accordingly, collisions avoidance may be shown even under disturbance.

 \section{Moving Obstacles}

\label{sec:moving}

Certain types of autonomous vehicle will unavoidably encounter moving obstacles, which are generally
more challenging to avoid than static equivalents. The main factors which affect the difficulty of
this problem are the characterization of the possible actions another object might take; the 
increased complexity of the search space and terminal constraints in the case of path planning;
and additional conservativeness in the case of sensor based systems. 

At one extreme, an obstacle translating at constant
speed and in a constant direction may be accounted for by merely considering the future position of
the obstacle. The other extreme is an
obstacle pursuing the vehicle, for which the set of potential locations grows polynomially along the
planning horizon. Several offerings also describe integrated approaches, including obstacle motion estimation 
from LiDAR sensors \cite{Montesano2008journ8}. However in this section discussion is focused on the avoidance behavior.

General planning algorithms suited for dynamic environments are also available, 
however in the absence of obstacle assumptions it is impossible to guarantee existence 
of a viable path, see e.g. \cite{Gecks2009book}. When planning in known environments, states which
  necessarily lead to collision -- the \textit{Inevitable Collision States} (ICS) -- 
  may also be abstracted and used to assist planning \cite{Petti2005conf4}. 
  If the motion of vehicles is known stochastically, the overall probability of collision for 
a probational trajectory may also be computed based on the expected behavior of other obstacles, see e.g.
\cite{Althoff2012journ2}. 
  
\subsection{Human--Like Obstacles}

Several works attempt to characterize the motion of moving obstacles. For avoiding humans,
several models of socially acceptable pedestrian behavior are available  
(see e.g. \cite{Sisbot2007journ7, Ohki2012journ3, Foka2010journ9, Ziebart2009conf6}). 
An approach which avoids obstacles based on the concept
of personal space has been proposed and works well in practice \cite{Ohki2012journ3}. Other
approaches can avoid human-like obstacles while also considering the reciprocal effect
of the vehicles motion have also been proposed \cite{Foka2010journ9, Ziebart2009conf6}.

\subsection{Known Obstacles}

Obstacles translating at constant speed and in a constant direction may be avoided using the concept
of a \textit{velocity obstacle},
 see e.g. \cite{Shiller2010book5, Fiorini1996conf5, Fiorini1998journ0}. This is essentially the set
 of vehicle velocities that will result in collision with the obstacle, and by avoiding these
 velocities, collisions may be avoided. This result may be extended to arbitrary (but known) obstacle
 paths and more complex vehicle kinematics using the \textit{nonlinear velocity obstacle}, see e.g.
 \cite{Large2005journ7}.
The velocity obstacle method also extends to 3D spaces, see e.g.
\cite{Shim2007conf1, Yang2012journ4}.

\subsection{Kinematically Constrained Obstacles}

When obstacles are only known to satisfy nominal kinematic constraints, the set of possible obstacle
positions grows over time. Avoidance may be ensured by either trajectory based methods or reactive methods.

\subsubsection*{Trajectory Based Methods}

There are three basic methods of planning trajectories which avoid such obstacles:

\begin{itemize}
  \item Ensuring that whenever a collision could possibly occur the vehicle is stationary -- this is
  referred to as \textit{passive motion safety} (see e.g. \cite{Bouraine2012journ2, Belkhouche2009journ1}). In some situations it is impossible
  to show any higher form of collision avoidance, through it ultimately relies on the behavior of
  obstacles to avoid collisions.
  \item Ensuring the vehicle can move arbitrarily far away from the obstacle over a infinite horizon. 
  In one such approach the the time-minimal paths to achieve a particular task were calculated \cite{van2008journ1}. 
  Similar examples of approaches include more allowance for other uncertainties, see e.g. \cite{Du2012journ5}.
  \item Ensuring the vehicle lies in a set of points that cannot be easily reached by the
  obstacle \cite{Wu2012journ4}.
  Under certain assumptions a non-empty set of points may be 
  found which lies just behind the obstacles velocity vector. This allows avoidance over a infinite horizon,
  while being possibly less conservative than the previous option.
\end{itemize} 

When performing path planing in a sensor based paradigm, the same types of approaches may be used, through
collision avoidance may be harder to show for general obstacle assumptions. The main additional assumption
is that any occluded part of the workspace must be considered as a potential dynamic obstacle
 \cite{Chung2009journ, Bouraine2012journ2}. Naturally this makes the motion of any vehicles even more 
 conservative.

\subsubsection*{Reactive Methods}
\label{sec:movreact}

When moving obstacles are present in the workspace, it is still possible to design reactive
navigation strategies which can provably prevent collisions, at least with some more restrictive assumptions about
the obstacles motion. These methods are outlined as follows:

\begin{itemize}
  \item  When obstacle are sufficiently spaced (so that multiple obstacles must not be simultaneously avoided), an extension of the velocity 
  obstacle method has been designed to prevent collisions \cite{Wang2012journ, MOVING1}.
  This effectively steers the vehicle towards the projected edge of the obstacle, and
   was applied to the semi-autonomous collision avoidance of robotic wheelchairs.
   
  \item Certain boundary following techniques proposed in Sec.~\ref{sec:bfdis} have been
successfully extended to moving obstacles and maintain provable collision avoidance, assuming 
the obstacles are sufficiently spaced and their motion and deformation is known to satisfy some 
smoothness constraints, see e.g. \cite{Savkin2012conf0, Matveev2012journ5}. These are all 
based on sliding mode control, which retains the advantages discussed previously in Sec.~\ref{sec:bf}. 

\item Some artificial potential field 
methods have been extended to moving obstacles, through rigorous 
justification is not provided (see e.g. \cite{Ren2008journ5, Ge2002journ1}).
\end{itemize}

 \section{Multiple Vehicle Navigation}

\label{sec:multiple}

Navigation of multiple vehicle systems has gained much interest in recent years. As autonomous vehicles
are used in greater concentrations, the probability of multiple vehicle encounters correspondingly increases,
and new methods are required to avoid collision.

The study of decentralized control laws for groups of
mobile autonomous robots has emerged as a challenging
new research area in recent years (see, e.g., \cite{SAV04,ST10a, COM1, COM2, COM4} and
references therein). Broadly speaking, this problem falls
within the domain of decentralized control, but the unique
aspect of it is that groups of mobile robots are dynamically
decoupled, meaning that the motion of one robot does not
directly affect that of the others.
This type of systems is viewed as a networked control system, which is an active field of
research. For examples of more generalized work in this area, see e.g. \cite{MS03,SAV06,Savkin2007journ6,
Matveev2009book7}. One of the important applications of navigation of multi-vehicle systems is  is in sensing
coverage. To improve coverage and
reduce the cost of deployment in a geographically vast area,
employing a network of mobile sensors for the coverage
is an attractive option. Three types of coverage problems
for robotic  were studied in recent years; blanket
coverage \cite{SJM12}, barrier coverage \cite{CS11,CS12}, and sweep
coverage \cite{CS11,CSJ11}. Combining existing coverage algorithms with effective collision avoidance methods in
an open practically important problem.

While there is an extensive literature on centralized navigation of multiple vehicles, it is only
briefly mentioned here, since it is generally not applicable to arbitrarily scalable 
on-line collision avoidance systems. 
Examples of off-line path planning systems which can find near optimal trajectories 
for a set of vehicles are available, see e.g. \cite{Skrjanc2010journ4}.
Another variation of this problem involves a precomputed prescription of the paths to be followed, 
where the navigation law must only find an appropriate velocity profile which avoids collisions 
(see e.g. \cite{Peng2005journ8, Cui2012journ6}).

\subsection{Communication Types}

There are three common modes of communication in multiple vehicle collision avoidance systems:

\begin{itemize}
  \item \textit{Direct state measurement.} This can be achieved using only sensor information to
  measure the state of the surrounding vehicles, and is used in many non-path based reactive approaches.
  \item \textit{Single direction broadcasting.} In addition to the physical state of the vehicle, additional variables 
  are also transmitted, usually relating to the current trajectory of the vehicle. As discussed in Chapt.~\ref{chap:multiple},
  this allows the projected states of other vehicles to be avoided during planning. 
  This type of communication is occasionally referred to as \textit{sign board} based.
  \item \textit{Two way communication.} This can range from simple acknowledgment signals to full decentralized optimization
  algorithms. These are commonly used for decentralized MPC, though some MPC variants have been proposed where 
  sign boards are sufficient.
\end{itemize}

A number of different models of communication delay and error are considered in networked navigation problems.
The ability to cope with unit communication delays, packet dropouts and finite communication ranges is 
definitely desirable in any navigation system.

\subsection{Reactive Methods}

The most basic from of this problem only considers a small number of vehicles. 
For example, navigation laws have been proposed to avoid collisions between two vehicles traveling at constant
speed with turning rate constraints, see e.g. \cite{Fujimori2000journ3, TFMasoud2012journ7}.
A common example of this type of system is an \textit{Air Traffic Controller} (ATC). However,
these types of navigation systems do not directly relate to avoiding collisions in cluttered environments
\cite{Kuchar2000journ2}.

\subsubsection*{Potential Field Methods}

Potential field methods may be constructed to mutually repel other vehicles. In some ways 
this approach is more satisfactory than the equivalent methods applies to static obstacles --
 for example local minima are less of an issue in the absence of contorted obstacle shapes.
Methods have been proposed which avoid collision between a unlimited number of velocity controlled unicycles
  or velocity controlled linear vehicles
 \cite{Mastellone2008journ8,Hernandez2011book2,Stipanovic2007journ5}. 
 One variant, termed the multi-vehicle navigation function, is able to show convergence to targets in 
 the absence of obstacles.
 However these still use similar types of repulsive and attractive fields, see e.g. \cite{Widyotriatmo2011journ5, 
 Dimarogonas2006journ7, Tanner2012journ}.  
 
 Other variants also
 include measures to provably maintainable cohesion between groups, see e.g. \cite{Dimarogonas2008journ9}, 
while others have also been applied to vehicles with limited sensing
 capabilities \cite{Dimarogonas2007journ4}. Many methods provide good results while neglecting mathematical analysis
 of collision avoidance, see e.g. \cite{Chang2003conf8, Samitha2009journ0, Fahimi2009journ3}.
 
 In cases where finite acceleration bounds are present (but still without any nonholonomic
constraints), an mutual repulsion based navigation system with a more sophisticated avoidance function
has been proven to avoid collisions for up to three vehicles \cite{Hoffmann2008conf0}. When more vehicles
are present, it is possible to back-step the additional dynamics into a velocity controlled model, through this 
does not lead to bounds on the control inputs of the dynamic model \cite{Loizou2008journ}. 

Many of these methods can be extended to static obstacles, and 
these combined systems are achieved
by the same avoidance functions as the single vehicle case, see e.g. \cite{Loizou2008journ}.
An interesting question may be whether transformation based approaches allowing arbitrary dynamics
(see e.g. \cite{Minguez2003conf, Minguez2009journ, Blanco2008journ}) may be extended to multiple vehicle cases.
As such, showing collision avoidance for an unlimited number of acceleration constrained vehicles using a
repulsive function seems to be an unsolved problem in robotics.

\subsubsection*{Reciprocal Collision Avoidance Methods}

Approaches termed \textit{Optimal Reciprocal Collision Avoidance} (ORCA) achieve collision avoidance by assuming each vehicle
takes half the responsibility for each pairwise conflict, with the resulting constraints forming a
set of viable velocities from which a selection can be made using linear programming, see e.g.
  \cite{van2009conf0, Snape2011journ2}. 
Some interesting extensions have been proposed to the RCA concept, for example it has been applied to both
nonholonomic vehicles and linear vehicles with acceleration constraints, while maintaining
collision avoidance \cite{Snape2009conf2, van2011conf2, Snape2011journ2}. 
 The method may be extended to arbitrary vehicle models, as rigorous avoidance is achieved though
 the addition of a generic bounded-deviation path tracking system  \cite{Rashid2012journ9, Snape2011journ2, Mora2012conf}. These methods
 are also able to integrate collision avoidance of static obstacles, which easily integrates into the navigation framework.
 
 RCA is an extension of an similar method based on collision cones called
  \textit{Implicit Cooperation} \cite{Abe2001conf5}.
  Another method has also been proposed which is based on collision cones, called \textit{Distributed Reactive
  Collision Avoidance} (DRCA).
This has the benefit of showing achievement of the vehicles objective in limited situations, ensuring minimum
speed constraints are met when global information is available, and
showing robustness to disturbance
 \cite{Lalish2008conf6, Lalish2012journ0}.

\subsubsection*{Hybrid Logic Methods}

For these approaches, discrete logic rules are used to coordinate vehicles. In most cases, this is achieved
through segregation of the workspace into cells, which can each only hold one vehicle, see e.g.
\cite{Belta2007journ, ghrist:1556, Reveliotis2011journ0, Nishi2005journ1}. In these cases,
collisions can be prevented by devising a scheme where two vehicles do not attempt to occupy the same
cell simultaneously.
Additionally, many methods of integrating this with control of the vehicle's dynamics has
been proposed, see e.g. \cite{Cowlagi2012journ1}. Hybrid control systems are becoming increasingly used 
to control real world systems (see e.g. \cite{HYBRID1}).

In some approaches the generation of cells may be on-line and ad-hoc. This is
useful when minimum speed constraints are present -- the vehicles may be
instructed
 to maintain a circular holding pattern, and then to shift their holding pattern
 appropriately when safe. In this case some different possibilities for the
 shifting logic have been proposed, for example based on vehicle priority
 \cite{Krontiris2011conf}, or traffic rules
\cite{Pallottino2007journ7}.

\subsection{Decentralized MPC}
\label{sec:dmpc}

While optimal centralized MPC is theoretically able to coordinate groups of vehicles, the
underlying optimization process is too complex for any scalable real time application. Examples of centralized MPC
for multiple vehicle systems are available (see e.g. \cite{Farrokhsiar2012conf}).

Decentralized variants
of MPC in general do not specifically address the problem of deadlock. For example in
\cite{Kuwata2011journ0} a distributed navigation system is proposed which is able to plan near optimal
solutions that robustly prevent collisions and allow altruistic behavior between the vehicles which
monotonically decrease the global cost function. However this does not equate to deadlock avoidance,
which is discussed in Sec.~\ref{sec:dead}.

A review of general decentralized MPC methods is available \cite{Bemporad2010book}, along with a review specific to 
vehicle navigation \cite{Shin2009journ}. 
There are currently four main methods of generating deconflicted trajectories which seem suitable for 
coordination of multiple vehicles:

\begin{itemize}
  \item Decentralized optimization can find the near-optimal
solution for a multi-agent system  using dual decomposition to find a set of trajectories for the
system of vehicles, see e.g. \cite{Raffard2004conf3, Wakasa2008conf5, Summers2012conf}. While this is more efficient than
centralized optimization, it requires many iterations of communication exchange between vehicles in
order to converge to a solution. Other types decentralized planning algorithms may also be effective,
for example a decentralized RRT approach has been proposed which implements the same type of processes
 \cite{Desaraju2012journ5}.
 
\item  Other approaches have been proposed using multiplexed MPC (see e.g. \cite{Kuwata2007journ7, Siva2011conf3}), 
and sequential decentralization (see e.g. \cite{Adinandra2012conf}). The robust control input for each vehicle may be computed by updating the
trajectory for each vehicle sequentially, at least when they are close. While multiplexed MPC is suited to real time
implementation, a possible disadvantage is path planning cannot occur simultaneously in two adjacent
vehicles. However, the same framework been extended
to provide collision avoidance in vehicle formation problems \cite{Weihua2011journ}.

\item Another possible solution is to require acknowledgment signals before implementing a
possible trajectory, and has the benefit of not requiring vehicles to be synchronized. This method seems an effective solution 
\cite {Bekris2009journ, Bekris2012journ6}, however interaction between vehicles may cause planning delays under certain 
conditions.

\item  Approaches also have been proposed which permit single communication exchanges per control
update \cite{Vaccarini2009conf6, Defoort2009journ1}. This is done by including a coherence objective
to prevent the vehicles from changing its planned trajectory significantly after transmitting it to
other vehicles. In Chapt.~\ref{chap:multiple}, an original set of trajectory constraints are proposed, which are possibly more 
general in that they do not explicitly enforce coherency objectives or limit the magnitude of trajectory alterations.

\end{itemize}

MPC may also easily include maintenance of objectives other than collision
avoidance. For example, radio propagation models have been included in the path evaluation function, so that communication between vehicles is maintained
 \cite{Gratli2012journ0, Grancharova2012conf}.

\subsection{Deadlock Avoidance}
 \label{sec:dead}

The collision avoidance techniques described in the Sec.~\ref{sec:dmpc} will not generally
guarantee that vehicles arrive at their required destination. System states which do not evolve to
their targets are referred to as \textit{deadlocks}\footnote{States where the vehicles are not
stationary indefinitely may also be further categorized as \textit{livelocks}}.
Naturally some control laws are more prone to deadlocks than others; this may be investigated using
probabilistic verification \cite{Pallotitino2007conf}.

It can be shown that when using suitable controllers, multiple vehicles may converge to their
targets in open areas, see e.g. \cite{Tomlin1998journ5, Tanner2012journ}. This means a more interesting question
relates to deadlock avoidance in unknown, cluttered environments.
It seems a generalized solution to the latter would be relatively sophisticated, and require
significant overlap with other areas such as distributed estimation, mapping, decision making and
control.

Currently, deadlock resolution systems are almost exclusively constructed based
on transitions between nodes on a graph or equivalent, and the solution can be
described as a resource allocation problem. Variations of this problem comprise
a well studied field, see e.g. \cite{Manca2011conf1, ghrist:1556, Jager2001conf5, Reveliotis2008conf2, Reveliotis2011journ0}. A
common simple example of such an algorithm is the bankers algorithm, which is
based on the concept of only allowing an action to be taken if it leaves
appropriate spatial resources so that every other agent may eventually run to
completion \cite{kim:849}.

However, these types of deadlock avoidance system are somewhat decoupled from the
actual sensor information available to the vehicle, and in some cases it would be advantageous to
eliminate the need for graph generation and use algorithms that correspond directly to the
continuous state space. These algorithms are also generally global solutions, requiring knowledge of
the parameters of all other vehicles. In some cases it could be more useful to use a more local
approach, especially when there is a very low density of vehicles operating in the environment.

In Chapt.~\ref{chap:dead} an initial solution to a simplified version of this problem is offered, where only two vehicles are present.
In future work it is hoped this can be extended to more general situations.

 \section{Summary}
\label{ch1:con}

This chapter provides a review of a range of techniques related to the navigation of
autonomous vehicles through cluttered environments, which can rigorously achieve collision avoidance
for some given assumptions about the system. This continues to be an active area of research, 
and a number of channels where current
approaches may be improved are highlighted. Approaches to avoiding collisions 
between multiple vehicles along with moving obstacles are considered. In
particular, approaches based on local sensor information are emphasized, which seems more
difficult and relevant than global approaches where full knowledge of the environment is assumed.
 Finally, the virtues of recently proposed MPC, decentralized MPC, and sliding mode
control based approaches are highlighted, when compared to existing methods.

\chapter{Collision Avoidance of a Single Vehicle}
\label{chap:singlevehicle}

In this chapter, the problem of preventing collisions between a vehicle and a set of
static obstacles while navigating towards a target position is considered. In this
chapter, the vehicle is presumed to have some information about the obstacle set;
 a full characterization is considered in Chapt.~\ref{chap:convsingle}. 
The aim of this chapter
 is to establish the navigation framework used in all subsequent chapters. Both
holonomic and unicycle kinematic motion models are considered (these were
informally described in Chapt.~\ref{chap:lit}).
\par 
The body of this chapter is organized as follows. In Sec.~\ref{ch2:ps}, the
problem statement is explicitly defined, and the vehicle model is given. In
Sec.~\ref{ch2:csa}, the navigation system structure is presented.
Sec.~\ref{ch2:sim} offers simulated results.
Finally, Sec.~\ref{ch2:con} offers brief conclusions.

\section{Problem Statement}
\label{ch2:ps}

A single autonomous vehicle traveling in a plane is considered, which is associated with a
steady point target $T$. The plane contains a set of unknown, untransversable,
static, and closed obstacles $D_j \not\ni T, j \in [1:n]$. The objective is to design a navigation
law that drives every vehicle towards the assigned target through the obstacle-free part of the plane
$F := \mathbb{R}^2 \setminus D$, where $D:= D_1\cup \ldots \cup D_n$.
Moreover, the distance from the vehicle to every obstacle and other vehicles should constantly exceed
the given safety margin $d_{sfe}$, which would naturally exceed the vehicle's physical radius. 

\subsection{Holonomic Motion Model}

Holonomic dynamics are generally encountered on helicopters and omni-directional
wheeled robots. A discrete-time point-mass model of vehicle is used, where for
simplicity the time-step $\Delta t$ is normalized to unity and the acceleration
capability of the vehicle is assumed to be identical in all
directions:\footnote{Note
Eq.\eqref{holo} may also be expressed in state-space notation, however in this work the notation used was found to be more convenient.}

\begin{subequations}
\label{holo}
\begin{gather} 
\blds(k+1) = \blds(k) + \bldv(k) + \frac{1}{2} \cdot \left( \bldu(k) + \bldw(k) \right) ,
\\ \bldv(k+1) = \bldv(k) + \bldu(k) + \bldw(k),
\label{d.2}
\\
v(k) := \norm{\bldv(k)} \leq v_{max}, \quad \norm{\bldu(k) } \leq u_{max}, \quad  \norm{\bldw(k)} \leq w_{max}.
\label{d.3}
\end{gather} 
\end{subequations}

Here $k \in \mathbb{N}$ is the
time index; $\blds=\col (s_x,s_y)$ is the vector of the vehicle's coordinates;
$\bldv = \col (v_x,v_y)$ is its velocity vector; $\bldu = \col(u_x,u_y)$ is the
control input; and the disturbance $\bldw = \col(w_x, w_y)$ accounts for any
kind of discrepancy between the real dynamics and their nominal model. In
particular, $\bldw$ may comprise the effects caused by nonlinear characteristics of
a real vehicle. Furthermore, $v_{max}$ is
the maximal achievable speed; $u_{max}$ is the maximal controllable
acceleration; and $w_{max}$ is an upper bound on the disturbance. Only the
trajectories satisfying all constraints from Eq.\eqref{holo} are feasible. The
state can be abstracted as $\state := \langle \blds, \bldv \rangle$, and the
control input as $\mathscr{U} \equiv \bldu$. The {\it nominal} trajectories,
 generated during planning, are created by setting $\bldw \equiv 0$ 
 (and thus may deviate from the actual ones). Here and throughout, $\norm{\cdot}$ 
is the standard Euclidean norm.

\subsection{Unicycle Motion Model}

Unicycle motion models describe vehicles which are associated with some heading
which determines the direction of movement, with changes to the heading limited
by a turning rate constraint.
A continuous-time point-mass model of the vehicle with discrete control
updates is considered. As before, the time step  $\Delta
t$ is normalized to unity.\footnote{For computation, Eq.\eqref{uni} may easily be analytically
converted into a fully discrete-time model. Also, through not done in this work,
an arbitrary acceleration bound $a_{max}$ such that $\ddot{\blds}(t) < a_{max}$ may be enforced. This was found to give 
better practical results as it leads to less conservative settings of the maximum rotation rate $u_{\theta, max}$. Unfortunately, it interferes with analysis 
of the auxiliary controller Eq.\eqref{cl.ma}; thus extending the analysis to cover this case remains an area of future research.}
\par
For $k < t \leq (k+1)$:
\begin{subequations}
\label{uni}
\begin{gather}
\blds(t) = \blds(k) + \int_k^t \left[ \begin{array}{c}\cos (\theta(\acute{t}) + \varphi(\acute{t}))\\ \sin (\theta(\acute{t}) + \varphi(\acute{t})) \\\end{array}\right]\cdot v(t) d\acute{t}\\
v(t) = v(k) + u_v(k)\cdot (t - k)+  \int_k^t w_{v}(\acute{t}) d\acute{t}, \quad |u_v(k)| \leq u_{v,max}\\
\theta(t) = \theta(k) + u_{\theta}(k)\cdot (t - k) + \int_k^t w_{\theta}(\acute{t})v(\acute{t}) d\acute{t}, \quad |u_{\theta}(k)| \leq u_{\theta, max}\\
0 \leq v(t) \leq v_{max}, \quad \abs{\varphi} \leq \varphi_{max}, \quad \abs{w_{\theta}} \leq w_{\theta, max}, \quad  \abs{w_{v}} \leq w_{v, max}
\end{gather}
\end{subequations}

Here $k \in \mathbb{N}$ is the time index; $\blds=\col (s_x,s_y)$ is the vector of the vehicle's coordinates; $v$ is the scalar
speed; $\theta$ is orientation angle; $v_{max}$ is the maximal achievable speed;
$\varphi$ is the angular difference between the vehicle orientation and the
velocity direction caused by the wheels slip; $w_{\theta}$ is the coefficient of
the rotation rate bias due to disturbance; and $w_v$ is the longitudinal
acceleration bias due to disturbance. This assumes that the maximal
feasible rotation rate bias due to disturbance is proportional to the vehicle
speed. The state can be abstracted as $\state :=
\langle \blds, v, \theta \rangle$, and the control input as $\mathscr{U} :=
\langle u_v, u_{\theta} \rangle$. The {\it nominal} trajectories
which are generated during planning are created by setting $\varphi \equiv
w_{\theta} \equiv w_{v} \equiv 0$, and naturally will deviate from the actual
ones.

\subsection{Sensor Requirements}
\label{sec:sensreq}

In this chapter, is is assumed that there is an arbitrary, time varying sensed area $F_{vis}(k)
\subset F$ corresponding to a subset of the obstacle free part of the plane. Normally, this would incorporate line-of-slight constraints together
with visibility range constraints, which are nominally taken to be $R_{nom} > 0$
for the purposes of designing the trajectory planning algorithm (this does not affect robustness properties). At time $k$, it is assumed
 that for any point $x$, the vehicle has knowledge of the distance to $\mathbb{R}^2 \setminus F_{vis}(k)$:

\begin{equation}
\dist_{F_{vis}(k)} [x] := \min_{y \in \mathbb{R}^2 \setminus F_{vis}} \norm{x - y}
\end{equation}

At time $k$, it is assumed that the vehicle has knowledge of its state $\state$. It is also
assumed that the vehicle has knowledge of the position of the target $T$.
Allowing for state measurement noise remains an area of future research, through
it seems likely it can be incorporated into the same robust navigation framework
proposed here.

 \section{Navigation System Architecture}
\label{ch2:csa}

\textit{Model Predictive Control} (MPC) is based on iterative, finite horizon
optimization of a trajectory corresponding to a given plant model. Chosen over a
relatively short time horizon $\tau$ into the future, a {\it probational}
trajectory, commencing at the current system state $\state$, is found over the
horizon $[k,k+\tau]$ at every time instant $k$. When disturbance is absent
(considered in Chapts.~\ref{chap:convsingle} and \ref{chap:dead}), the control
representing the first time-step is implemented on the vehicle.
When disturbance is present (considered in this chapter and
Chapt.~\ref{chap:multiple}), an auxiliary trajectory tracking controller updates the
control to correct deviations from the probational trajectory (this is a somewhat
similar concept to the tube MPC methodology). In both cases, at each subsequent
time-step, the trajectory planning process is repeated for the new vehicle state
and obstacle observations.
 \par
 Unlike general trajectory planning approaches usually employed by full-information MPC approaches, the planning system only considers a
 certain subset of possible trajectories.
 This is a similar concept to the Dynamic Window and Curvature-Velocity methods
 (see Chapt.~\ref{chap:lit}).
This means the simplified trajectory generation scheme proposed in Sec.~\ref{Sec.plan}
is not \textit{complete} -- to accommodate this, if a trajectory cannot be found, the trajectory from the
previous time-step is reused, with an appropriate time-shift.
\par
The upper index $(\cdot)^*$ is used to mark probational trajectories, and the associated variables depend
on two arguments $(j|k)$, where $k$ is the time instant when this trajectory is generated and $j
\geq 0$ is the number of time steps into the future (the related value concerns the state at time
$k+j$).

\subsection{Overview}

The navigation system consists of two modules; the \textit{Trajectory Planning Module} (TPM), and the 
\textit{Trajectory Tracking Module} (TTM). The TPM may be summarized by sequential execution of the following steps:

\begin{enumerate}[{\bf S.1}]
\item \label{step1} Generation of a finite set $\mathscr{P}$ of planned trajectories, each starting
at the current vehicle state $\state$ (see Sec.~\ref{Sec.plan}).

\item Refinement of $\mathscr{P}$ to only \textit{feasible} trajectories (see Sec.~\ref{subsec.ansm}).

\item \label{step3} Selection of a trajectory from $\mathscr{P}$:

\begin{itemize}
  \item If $\mathscr{P}$ is empty, the probational trajectory is \textit{inherited} from the previous time step
(with a proper time shift). 
 \item Otherwise, the probational trajectory is \textit{updated} by choosing an element of $\mathscr{P}$
   minimizing some cost function (see Sec.~\ref{sec.sel}).
\end{itemize}

\item Application of the first control $\mathscr{U}^*(0|k)$ related to the
selected probational trajectory to the vehicle.

\item $k:=k+1$ and go to { S.\ref{step1}}.

\end{enumerate}

The TTM runs concurrently to the TPM and
updates the vehicle's control $\mathscr{U}$ (in navigation problems where disturbance is present). The TTM architecture depends on the
vehicle model type:

\begin{itemize}
  \item \textbf{Holonomic:} In this case,
TTM updates occur at the same rate as TPM updates, and the TTM may be simply
added as an additional step after S.\ref{step3}. This means whenever the
probational trajectory is updated rather than inherited, the implemented control will be identical to
$\mathscr{U}^*(0|k)$.
\item \textbf{Unicycle:} In this case, TTM updates occur at a much faster rate than TPM
updates; this is because the TTM is implemented as a sliding mode control law (see Sec.~\ref{sec:uni}). 
Thus even when a trajectory is updated rather than inherited, the implemented control $\mathscr{U}(t)$ could only
be equal to $\mathscr{U}^*(0|k)$ only at exact, integral time instants.
\end{itemize}

\subsection{Safety Margins}
\label{subsec.ansm}
To ensure collision avoidance, TPM respects more conservative safety margins than $d_{sfe}$.
They take into account deviations from the probational trajectory caused by disturbances and are
based on estimation of these mismatches. 
\par
The first such estimate addresses the performance of TTM; let $d_{trk}$ be defined as an upper bound on
the translational deviation over the planning horizon between the probational trajectory and the
real motion of the vehicle driven by the TPM along this trajectory in the face of disturbances.
\par
 Computation of $d_{trk}$ takes into account the particular design of TPM and
 is discussed in Sec.~\ref{sec:ptm}. To decouple designs of TPM and TTM, the control capacity is
 a priory distributed between TPM and TTM.
Specifically, TPM must generate the controls only within the reduced bounds determined by:

\begin{subequations}
\label{v.nom}
\begin{gather}
\shortintertext{\bf Holonomic:}
u_{nom} := \mu \cdot u_{max}\\
\shortintertext{\bf Unicycle:}
u_{v,nom} := \mu_{u} \cdot u_{v,max},  \quad u_{\theta, max}:= \mu_{\theta} \cdot u_{\theta, max}, \quad v_{nom} :=  \mu_{v} \cdot v_{max}
 \end{gather}
\end{subequations}

The remaining control capacity is allotted to TTM. Here $\mu$, $\mu_{u}$,
$\mu_{\theta}$, $\mu_{v} \in (0,1)$ are tunable design parameters.

Additionally, for unicycle vehicles, one additional constraint on the absolute curvature of the path is present -- the curvature must fall
below a level accessible by the vehicle. Since this level depends on the actual speed $v(t)$, which
is influenced by unknown disturbances, it is required that the
curvature of the planned path must fit the worst case scenario. In other words, it should not exceed:

\begin{equation}
\label{kappa.nom}
\varkappa_{max} := \frac{u_{\theta,max}}{v_{max}} .
\end{equation}

In fact, a tighter bound is required -- the curvature should not exceed $\varkappa_{nom} := \mu_\varkappa \cdot \varkappa_{max}$.
Here $\mu_{\varkappa} \in (0,1)$ is a design parameter chosen similarly to those in Eq.\eqref{v.nom}.

Along with the system parameters, the control capability allotments $\mu$, $\mu_{u}$,
$\mu_{\theta}$, $\mu_{v}$, $\mu_{\varkappa}$ uniquely determine the estimate $d_{trk}$
(for the adopted design of TPM), as will be shown in Sec.~\ref{sec:ptm}. In
particular, it should be noted $d_{trk}$ does not depend on the choice of the probational
trajectory.
\par 
The following more conservative safety margins account for not only disturbances but also the
use of time sampling in measurements of relative distances:

\begin{equation}
\label{d.dist}
d_{tar}:= d_{sfe} + \frac{v_{max}}{2} + d_{trk}
\end{equation}

To illuminate their role, the following definition is introduced:

\begin{Definition}
\label{def:dmut}
A probational trajectory is said to be 
\textbf{feasible} if the distance from any way-point to any static obstacle exceeds $d_{tar}$ 
for any trajectory.
\end{Definition}

The motivation behind this definition is illuminated by the following:
 
\begin{Lemma}
\label{lem:feas0}
Let the probational trajectory adopted for use at time step $k-1$ be  feasible.
Then despite of the external disturbances, the vehicles do not collide with obstacles 
on time interval $[k-1,k]$, and also respect the required safety margin
$d_{sfe}$.
\end{Lemma}

\pf For any $t \in [k-1,k]$, it is known vehicle deviates from the related
position on the probational trajectory by no more than $d_{trk}$. Consider the
end $k_\ast$ of the interval $[k-1,k]$ nearest to $t$; the times $t$ and $k_\ast$ are separated by no
more than half time-step. So thanks to the upper bound on the speed, the nominal
positions at times $t$ and $k_\ast$ differ by no more than $\frac{v_{max}}{2}$. Overall, the
distance between the way-point on the probational trajectory that is related to time step $k_\ast$
and the actual position at time $t$ does not exceed $\frac{v_{max}}{2}+ d_{trk}$. Since the distance from this way-point to any static obstacle exceeds $d_{tar}$, the distance
from the real position at time $t$ is no less than $d_{tar} -  \frac{v_{max}}{2} -
d_{trk}  = d_{sfe}$. \epf

\begin{Remark} \rm
\label{rem.van}
\rm For arbitrary (not unit) sampling time $\Delta t$, the addend $\frac{v_{max}}{2}$ in Eq.\eqref{d.dist}
should be replaced by $\frac{v_{max}}{2}\Delta t$. This addend asymptotically vanishes as $\Delta t
\to 0$; so typically $d_{trk} $ would too provided that the { discrete} planning
horizon $\tau$ is upper bounded. Thus the conservatism imposed by the extra addends in
Eq.\eqref{d.dist} may be attenuated by reducing the time step employed during the planning process. 
As the time-step is conversely increased, $d_{trk} $ grows but remains bounded
even if $\Delta t \to \infty$, as will be discussed in Sec.~\ref{sec:ptm}. However the
other addend $\frac{v_{max}}{2}\Delta t$ grows without limits as $\Delta t \to \infty$.
\end{Remark}

 Overly enhanced margins Eq.\eqref{d.dist} may make the navigation objective unachievable; for example,
the corridors amidst the obstacles may be too narrow to accommodate a distance of no less than
$d_{tar}$ to the both sides of the corridor. In fact,
feasibility of the enhanced margins is the basic requirement regulating the practical choice of
$\Delta t$. In the case of heavy uncertainty about the scene, a reasonable option is to pick $\Delta
t$ so that the quantization error $\frac{v_{max}}{2}\Delta t$ is comparable with the limit bound on
the disturbance-induced error $d_{trk}$.

\subsection{Trajectory Planning}
\label{Sec.plan}
Trajectory generation could potentially use many forms, and some varieties of
commonly used approaches were outlined in Chapt.~\ref{chap:lit}.
However in the approach proposed here, a simplified version is designed for the following reasons:

\begin{itemize}
  \item Compared to commonly used optimization approaches, the proposed method
  is able to improve tractability while not significantly influencing closed
  loop performance. Similar types of simplified planning methods have been proposed previously, see
  e.g. \cite{Ogren2005journ3, Blanco2008journ}.
\item In Chapt.~\ref{chap:multiple}
it will be possible to make certain inferences about the possible trajectories
of other vehicles. 
\end{itemize}

Here two tunable parameters are used; $\Delta v$ and $\Delta \Lambda$,
which determine the `mesh' size for trajectory generation.

\subsubsection*{Longitudinal Control Pattern}
The generated set of trajectories $\mathscr{P}$ consists of two parts $\mathscr{P}_+$
({\it cruising} trajectories) and $\mathscr{P}_-$ ({\it slowing} trajectories), with each part being
composed of trajectories matching a given pattern of speed evolution over the planning horizon:

\begin{itemize}
\item
The pattern $p_+:=(+ - - \ldots)$ means that the initial speed $v^\ast(0|k) := v(k)$ is first increased by a given tunable speed increment $\Delta v$:
\begin{subequations}
\label{addit}
\begin{gather}
\shortintertext{\bf Holonomic:}
v_+^\ast(1|k)=\min\big\{ v^\ast(0|k)+\Delta v ; v_{max}\big\}, \quad \Delta v < u_{nom} \times (1 \,\text{time unit}) \equiv u_{nom}\Delta t\\
\shortintertext{\bf Unicycle:}
v_+^\ast(1|k)=\min\big\{ v^\ast(0|k)+\Delta v ; v_{nom}\big\}, \quad \Delta v < u_{v,nom} \times (1 \,\text{time unit}) \equiv u_{v,nom}\Delta t
 \end{gather}
\end{subequations}

The speed bounds follow from Eq.\eqref{v.nom}. The associated control inputs which achieve 
this pattern is given in Sec.~\ref{sec:lcp}. The planned speed is then constantly decreased by
 subtracting $\Delta v $ at all subsequent time steps $j \geq 1$ (while preventing the speed from taking a
 meaningless negative value):
 
\begin{equation}
\label{subtr}
v_+^\ast(j+1|k)= \max\big\{ v^\ast(j|k)-\Delta v, 0\big\}
\end{equation}

\item
The pattern $p_-=(- - - \ldots)$ means the speed is constantly decreased in accordance with Eq.\eqref{subtr}, creating a second speed profile $v_-^\ast(j|k)$.
\end{itemize}

The feasibility of this control pattern is illustrated in Fig.~\ref{fig.where}. The planning horizon $\tau$ is selected from the requirement that when following any generated
trajectory, the vehicle halts within this horizon:

\begin{equation}
\tau := \begin{cases}
\tau_-= \left\lceil \frac{{ v(k)}}{{\Delta v}}\right\rceil &
\text{for slowing trajectories}
\\[5.0pt]
\tau_+ =\left\lceil \frac{{ v(k)}}{{\Delta v}}\right\rceil+2 & \text{for cruising trajectories},
\end{cases}
\end{equation}

where $\left\lceil a\right\rceil$ is the integer ceiling of the real number $a$. This horizon is
determined by the current vehicle speed $v(k)$ and may vary over time, though it may also be fixed so it is 
equal to an upper bound of all possible values. To ensure obstacle avoidance,
attention is limited to trajectories for which the vehicle halts within the planning time horizon and
observed part of the environment.

\begin{Remark} \rm
\label{rem.stop}
\rm The length of any planned trajectory does not exceed
$ l = v_{max} + \frac{v_{max}^2}{2 \cdot u_{v,nom}}$. It seems prudent to adjust relevant parameters so the trajectory remains
within the sensed area, so it follows that $(l + d_{tar})$ is less than the nominal sensor range $R_{nom}$.
\end{Remark}

From this point, the subscript $(\cdot)_\pm$ is used to annotate variables which may correspond to either speed profile $p_+$ or $p_-$.

\subsubsection*{Lateral Control Pattern} 
\label{sec:lcp}

The lateral profile of the trajectories is constrained to invariably consist of a sharp turn, followed
by a reduced turn, followed by a straight section. 

\subsubsection*{Holonomic} 

To formulate this, first
consider how the angle of the vehicle's velocity vector changes over each time
step. Starting from an initial angle $\theta(0|k)$ given by the vehicle actual velocity vector $\bldv(k)$, 
the progression $\theta(j|k)$ is bounded by:

\begin{equation}
\abs{\underbrace{\theta(j+1|k) - \theta(j|k)}_{\Delta \theta(j|k)}} \leq \Delta \theta_{\max}(j|k).
\end{equation}

If $\bldv(k) \equiv 0$, $\theta(0|k)$ is arbitrary. If $\theta^*(j|k)$ is known, the velocity vector, and thus the vehicle
position vector, may be trivially found over the length of the trajectory by evaluating Eq.\eqref{holo}.

\begin{figure}[ht]
\centering
\subfigure[]{\includegraphics[width=0.35\columnwidth]{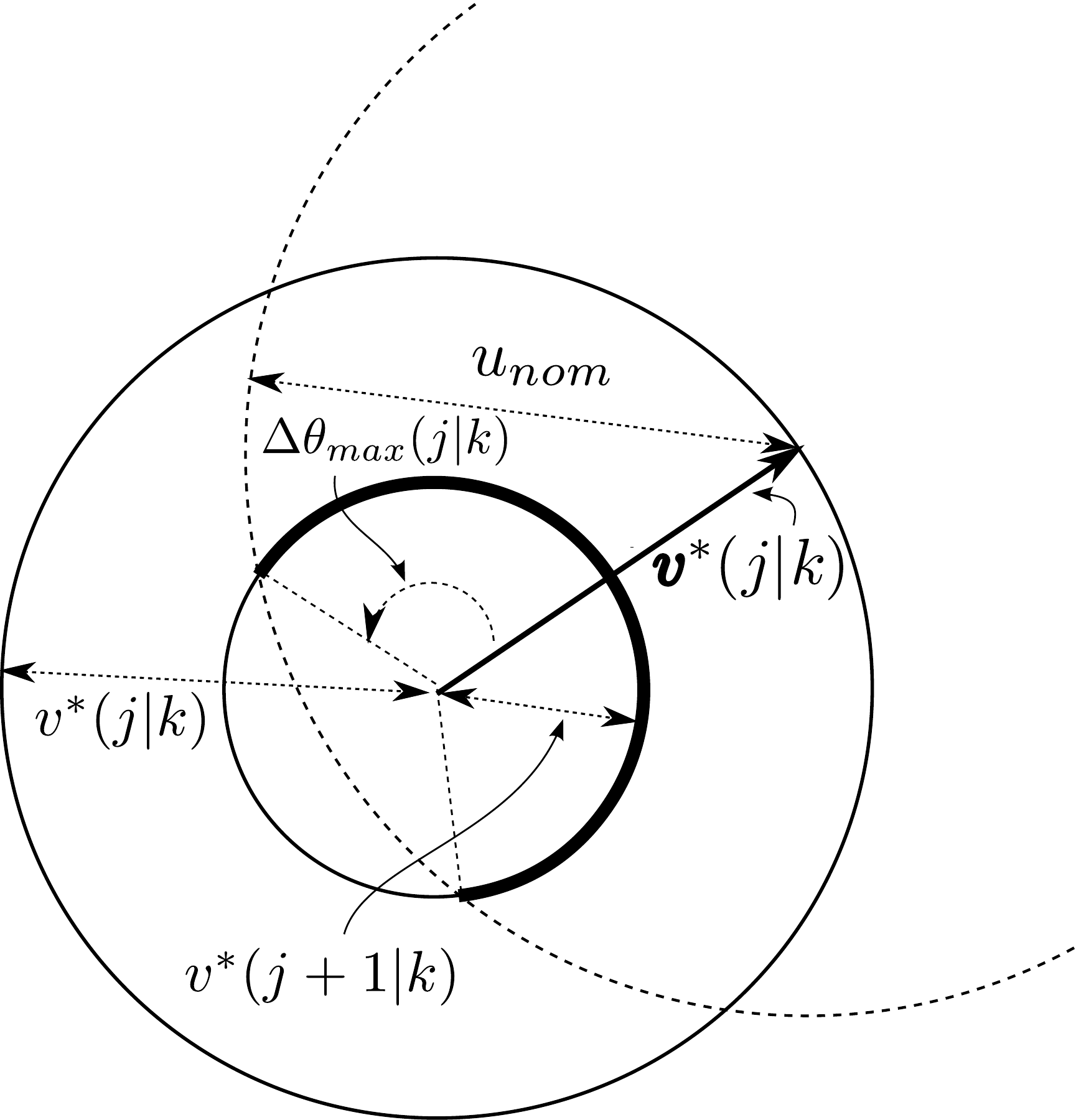}}
\subfigure[]{\includegraphics[width=0.35\columnwidth]{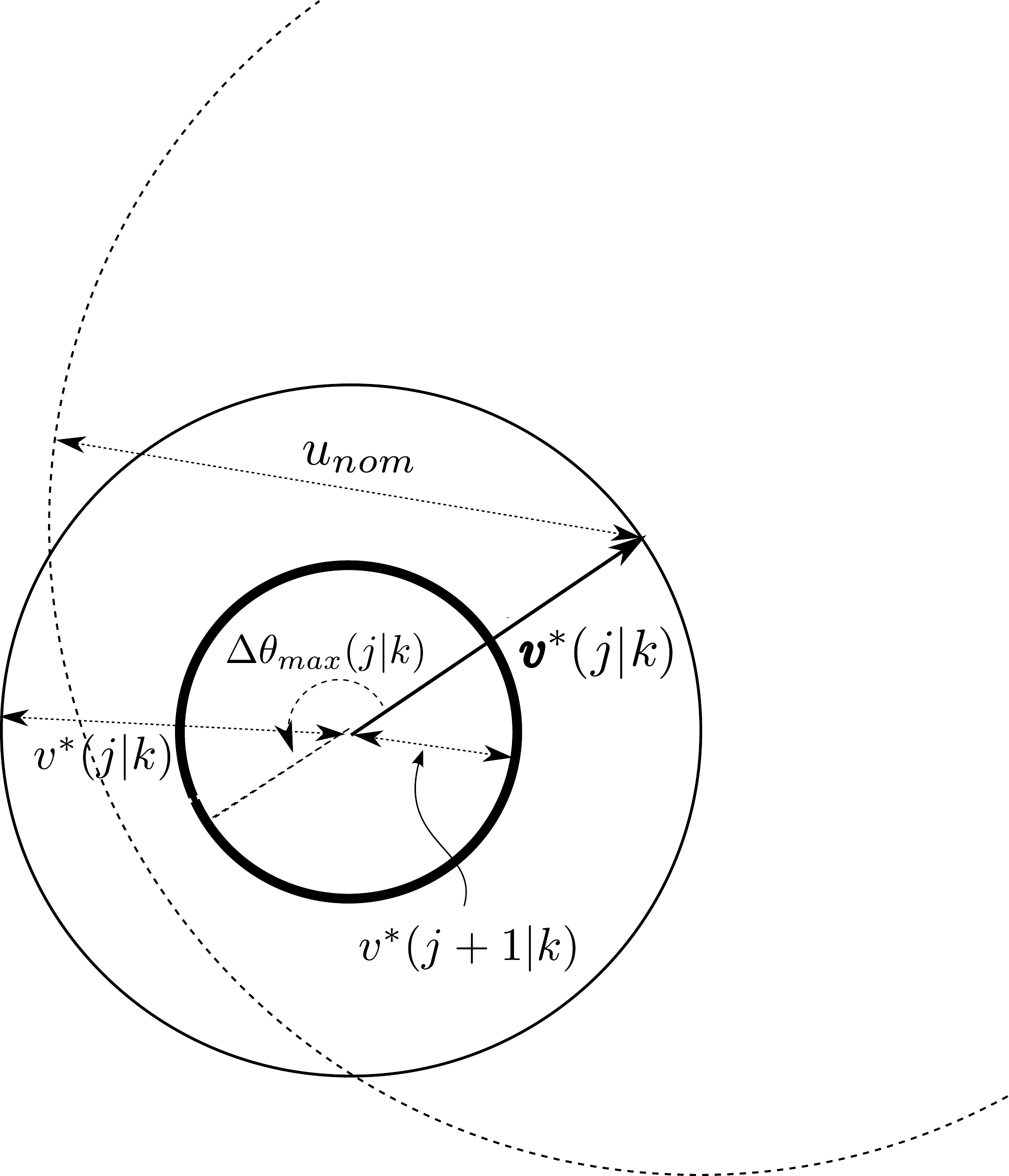}}
\subfigure[]{\includegraphics[width=0.35\columnwidth]{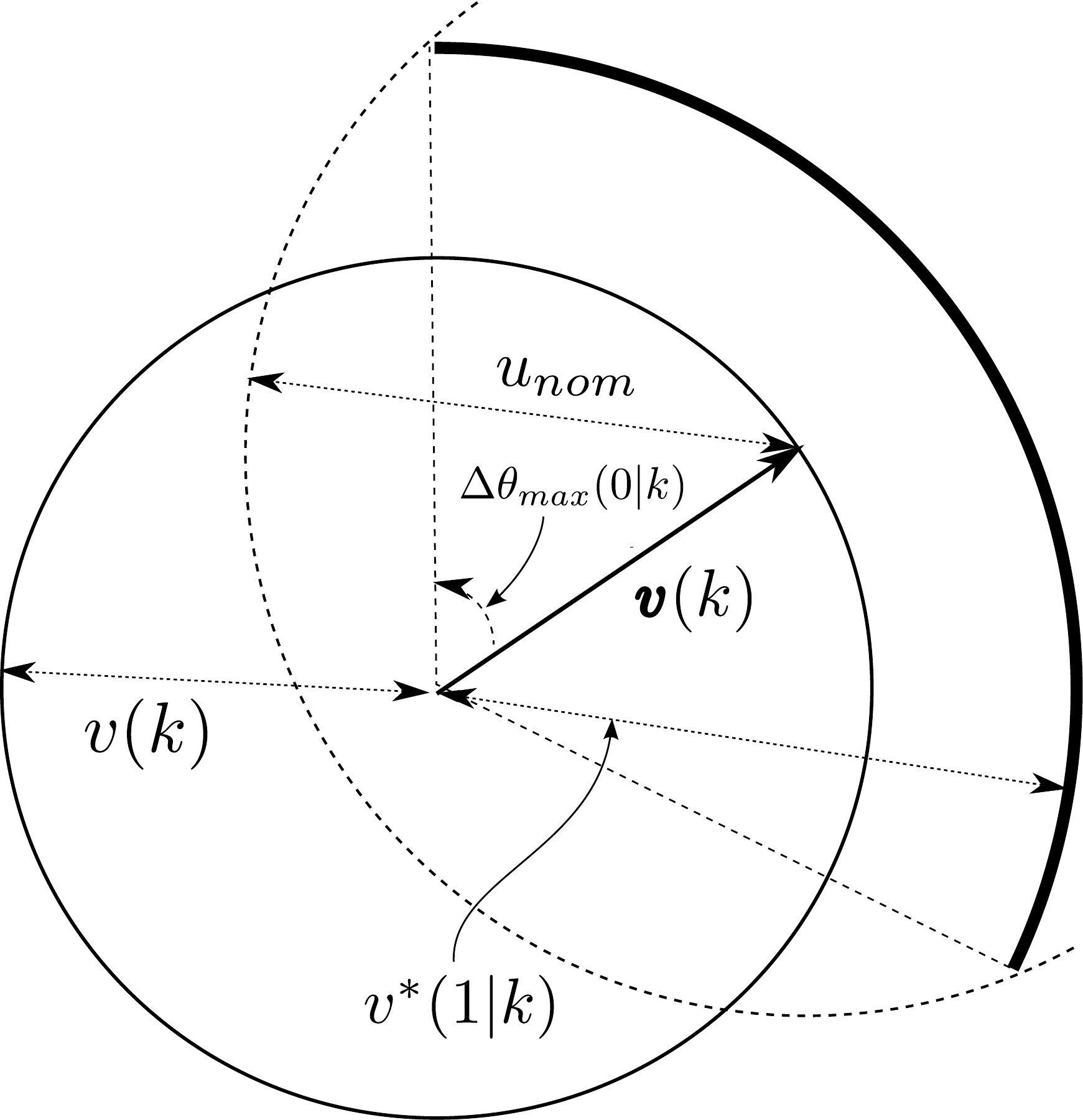}}
\caption{The set of possible values of the next probational velocity.}
\label{fig.where}

\end{figure}

The maximal angular deviation $\Delta \theta_{\max}(j|k)$ between the vectors $\bldv^\ast(j+1|k)$
and $\bldv^\ast(j|k)$ is determined by the speed profile $v_\pm^*(j|k)$ and
constant $u_{nom}$.
Specifically, $\Delta \theta_{\max}(j|k) = \pi$ if $v_\pm^*(j|k)+v_\pm^*(j+1|k)
\leq u_{nom}$; otherwise Fig.~\ref{fig.where} and the well-known formula giving the
angle of a triangle with three known sides imply that:

\begin{equation}
 \Delta \theta_{\max}(j|k) = \arccos \frac{v_\pm^*(j|k)^2+v_\pm^*(j+1|k)^2 - u_{nom}^2}{2 v_\pm^*(j+1|k) v_\pm^*(j|k)}.
\end{equation}

In this formulation, the overall turning profile can be uniquely described by a scalar $\Lambda(k)$. 
Small values correspond to short turns, while large values correspond to long turns:

\begin{equation}
\Delta \theta(j|k) =  \Delta \theta_{\max}(j|k) \cdot \sgn(\Lambda(k)) \cdot \sat_0^1(\abs{\Lambda(k)} - j).
\label{t.p}
\end{equation}

Here the function $\sat_b^a(x) := \min \{ \max \{x, b\}, a\}$, and the function $\sgn(\cdot)$ is the standard signum function. The possible values of
$\Lambda$ considered here are discretized by an integer $m$, and are given by:

\begin{equation}
\label{t.g}
\Lambda(k) = m \cdot \Delta \Lambda, \quad \left\{m \in \mathbb{I} : \abs{m} < \left\lceil\frac{\tau}{\Delta \Lambda}\right\rceil \right\}
\end{equation}

For each $m$, a planned trajectory may be uniquely determined using the nominal
equations Eq.\eqref{holo}, along the initial state $\state(k)$, and the speed profile $p \in \{p_+; p_-\}$.

\subsubsection*{Unicycle} 

The case of the unicycle vehicle model is similar to the holonomic model, 
and the expression for the turning control $u_\theta^*(j|k)$ is given directly by:

\begin{equation}
u_\theta^*(j|k) = \mu_{\varkappa}\varkappa_{max}\cdot \min  \big\{v^\ast(j|k), v^\ast(j+1|k)\big\} \cdot \sgn(\Lambda(k)) \cdot \sat_0^1(\abs{\Lambda(k)} - j).
\label{unilambda}
\end{equation}

\begin{figure}[ht]
	\centering
		\includegraphics[width=10cm]{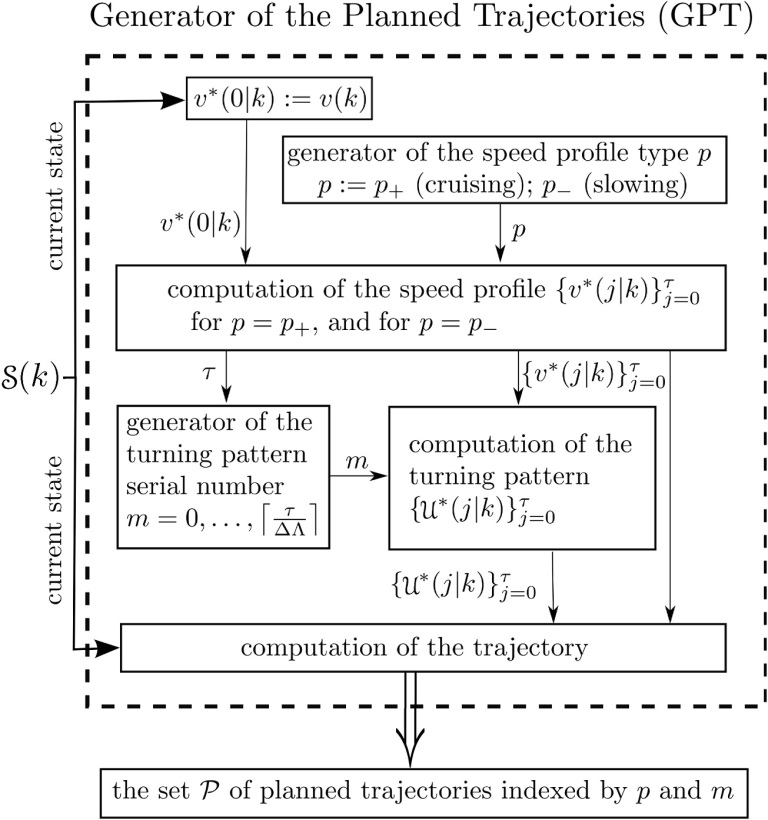}
	\caption{Generator of the planned trajectories.}
	\label{fig:gpt}
\end{figure}

This is similar to Eq.\eqref{t.p}, uses the same values of $\Lambda(k)$ from Eq.\eqref{t.g}, through it includes the absolute curvature constraints from Eq.\eqref{kappa.nom}. 
For each $m$, the planned trajectory may uniquely determined using the nominal
equations Eq.\eqref{uni}, along with the initial state $\state(k)$, and the speed profile $p$.
For unicycle vehicle models in particular, generation of the planned trajectories is summarized in Fig.~\ref{fig:gpt}.
 
\subsection{Trajectory Tracking}

\label{sec:ptm}

While the planning algorithm specified in the previous section is largely
similar for each vehicle model, the TTM is completely
different for holonomic and unicycle vehicle kinematic models.

\subsubsection*{Holonomic}
\label{sec:holotrac}

 In this case the TTM updates at the same rate as the TPM, as is only invoked
 after the failure of the trajectory planning algorithm to find a
feasible trajectory (which implies that the set $\mathscr{P}$ is empty, and thus a trajectory was inherited from the previous time-step).
 \par
Since this inherited
 trajectory might be traced for an extensive period of time, external
 disturbances may cause essential tracking errors. To compensate for this, the
 corrective control $\bldu_f(k)$ is added to the planned one, and so the overall
 control is as follows:
\begin{equation}
\bldu(k) = \bldu^*(0|k) + \bldu_f(k),
\end{equation}

The corrective control is formed by a linear feedback with saturation at the level $u_{exc} := (1-\mu)u_{max}$ from Eq.\eqref{v.nom}:

\begin{equation}
\bldu_f(k) = \sat^{u_{exc}}_{u_{exc}}\big[ k_0 \bar{\blds}(k) + k_1 \bar{\bldv}(k)\big].
\label{tfe}
\end{equation}

Here $\bar{\blds}$ and $\bar{\bldv}$ are attributed to the difference between the vehicle current state and
the predicted state on the trajectory being followed.
\par
For given $u_{exc}$ and $w_{max}$, the gain coefficients $a_{pos}$ and $a_{vel}$ from Eq.\eqref{tfe}
were optimized\footnote{For the research presented in this chapter, experimentally.} for minimum
distance deviation under disturbance. This leads to a value of the maximum positional error
$d_{trk}$, which can be computed in a similar manner as the robustly positively-invariant set
\cite{Rakovic2005journ0}.

There is some trade-off between the optimal distance deviation $d_{trk}$ and the tunable maximum
excess control input $u_{exc}$.

By following Lemma~\ref{lem:feas0}, it can be
guaranteed that whenever the last planned trajectory was feasible, the vehicle never moves closer than $d_{sfe}$
to the obstacle.

\begin{Remark} \rm
\rm In the model Eq.\eqref{holo}, the additive disturbance $\bldw(k)$
accounts for the un-modeled dynamics of the system, along with other factors. Whenever the overall
`disturbance' obeys the adopted upper bound, the above discussion argues the robustness of the
navigation law against the entire totality of uncertainties; the deviation from the planned trajectory does
not exceed $d_{trk}$. It however does not exclude that within this bound, behaviors like systematic
deviation from or oscillation about the planned trajectory may occur. If they yet occur and are undesirable,
further elaboration of the navigation law may be required.
\end{Remark}

An example of optimal feedback gains and positional errors for a particular set of parameters is
presented in Table~\ref{table1}.

\begin{table}[ht]
\centering

\begin{tabular}{| l | c |}
\hline
 $u_{exc}$ & $0.4$ \\
 \hline
 $w_{max}$ & $0.2$ \\
\hline

\end{tabular} \hspace{10pt} \begin{tabular}{| l | c |}
\hline
 $k_0$ & $1.33$ \\
\hline
 $k_1$ & $0.667$  \\
\hline
 $d_{max}$ & $0.30 $  \\
\hline
\end{tabular}
\caption{Trajectory tracking parameters for a holonomic vehicle.}
\label{table1}
\end{table}

It remains to discuss one special case that is not covered by the above instructions: the last
planned trajectory is exhausted, but the trajectory planning algorithm still fails. Since the planned probational velocity is zero, and the control Eq.\eqref{tfe} is employed, the
actual vehicle speed is typically small in this case. This leads to the possibility of a deadlock-type 
behavior. While it may be possible rely on a more complete trajectory planning method in these circumstances,
it cannot be generally guaranteed that a feasible trajectory exists, since the evolution of the sensed free space $F_{vis}(k)$ was 
not characterized. A solution is presented in Chapt.~\ref{chap:convsingle} which uses a better
characterization of $F_{vis}(k)$, and an alternative solution is explored in Chapt.~\ref{chap:dead},
which uses a navigation function to prevent deadlock.

\subsubsection*{Unicycle}
\label{sec:uni}
A more sophisticated trajectory tracking controller is required for unicycle vehicles, due to the reduced
maneuverability and nonlinear constraints. In this section a sliding mode control law is used
to perform the compensation.
\par
To compensate for space errors, the sliding mode path tracking approach
presented in \cite{Matveev2010conf0} is adopted (this is included in Chapt.~\ref{chapt:pf}). 
This navigation law is designed for pure steering control of a car-like vehicle governed by bicycle kinematics,
traveling at the constant nominal longitudinal speed. Chosen for its provable bounds on tracking
error in the face of a bounded wheels slip, this navigation law is adapted for trajectory tracking
purposes in this chapter. In particular, a Longitudinal Tracking Module (LTM) must be added to
the Path Tracking Module (PTM) to maintain the temporal difference between the planned and actual
positions. Other examples of navigation approaches suitable for this task include,  MPC based
trajectory tracking (see e.g. \cite{Gonzalez2009conf4}).

\subsubsection*{Path Tracking}
Three types of disturbance are considered in Eq.\eqref{uni} -- bias in the longitudinal
acceleration, bias in the turning rate, and a side slip causing an angular
difference between the vehicle orientation and the velocity direction.
\par
The proposed trajectory tracking strategy belongs to the class of sliding mode control laws
\cite{Utkin1992book1}. Due to the well-known benefits such as high insensitivity to noises, robustness against
uncertainties, and good dynamic response \cite{Utkin1992book1}, the sliding mode approach attracts a growing
interest in the area of motion control. The major obstacle to implementation of sliding mode
controllers is a harmful phenomenon called "chattering", which is undesirable finite frequency
oscillations around the ideal trajectory due to un-modeled system dynamics and constraints. The
problem of chattering elimination and reduction has an extensive literature (see e.g.
\cite{Edwards2006book3,LeUtMa09}). It offers a variety of effective approaches, including
smooth approximation of the discontinuity, insertion of low-pass observers into the control
loop, the combination of sliding mode and adaptive control or fuzzy logic techniques, and higher order sliding
modes. Whether chattering be encountered in applications of the proposed controller, it can be
addressed using these techniques.
\par 
In these simulations and experiments, the discontinuous signum function employed by the
navigation law was replaced with the saturated linear function $\sat^{+1}_{-1}(\cdot)$. This was
found to give more than satisfactory performance during trajectory tracking.

\begin{figure}[ht]
	\centering
		\includegraphics[width=0.6\mylength]{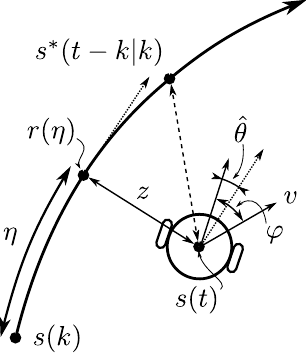}
	\caption{The trajectory tracking model of the unicycle.}
	\label{fig:track}
\end{figure}

Since PTM updates controls at much higher rate than TPM, the
continuous-time model of the vehicle motion is now employed. Furthermore, since disturbances are
the main concern now, they are included in Eq.\eqref{uni}. To come into details, the following variables are introduced (see
Fig.~\ref{fig:track}):

 \begin{itemize}
 \item $\eta$ -- the curvilinear abscissa of the point on the
 reference path closest to the robot;
  \item $\bldr(\eta) \in \mathbb{R}^2$ -- a regular parametric representation of the reference path in the world frame;
 \item $\varkappa(\eta)$ -- the signed curvature of the reference path at the point $r(\eta)$;
 \item $z$ -- the distance from the robot to the reference path;
\item $\varphi$ -- the angular difference between the robot orientation and the velocity direction caused by the wheels slip;
\item $\hat{\theta}$ -- the angle between the vehicle orientation and the orientation of a virtual vehicle perfectly tracking the reference path, positioned on $r(\eta)$;
\item $w_{\theta}$ -- the coefficient of the rotation rate bias due to disturbance;
\item $w_v$ -- the longitudinal acceleration bias due to disturbance.
 \end{itemize}

\begin{Lemma} 

\label{lem.kineq} The trajectory tracking kinematic model is as follows:
\begin{subequations}
\begin{gather}
\label{eq1} \dot{\eta} = \frac{v \cos (\hat{\theta} - \varphi) }{1-
\varkappa(\eta) z }, \; \dot{z} =  - v \sin (\hat{\theta} - \varphi),
\\
\label{eq21}
{\dot{v} = \begin{cases}u_v + w_v & \text{\rm if}\; 0 < v < v_{max} \;\;\;\text{\rm or}\; \left\{
\begin{array}{l}
v=0
\\
u_v + w_v >0
\end{array}
\right.
\;\text{\rm or}\; \left\{
\begin{array}{l}
v=v_{max}
\\
u_v + w_v <0
\end{array}
\right.
\\
0 & \text{\rm otherwise}
\end{cases}}
\\
\label{eq2} \dot{\hat{\theta}} = \varkappa(\eta)  \frac{v \cos (\hat{\theta} -
\varphi) }{1- \varkappa(\eta) z } +
u_{\theta} + {v} w_{\theta}.
\end{gather}
\end{subequations}
\end{Lemma}
\par
The  proof of this lemma is similar to that of Lemma~3.1 from \cite{Micaelli1993book3}, where the bicycle error tracking model is presented. To adapt this model to the case at hand, it suffices to put $u_{\theta} + {v}w_{\theta}$ in place of $\frac{v}{L} \left[
\tan(\delta+\beta) - \tan \varphi \right]$ and take into account the differences between a unicycle and bicycle, as well as that the velocity is now actuated.
\par
To generate the steering control signal, the navigation law from Chapt.~\ref{chapt:pf} is employed, slightly adopted to the current context:

\begin{gather}
\label{cl.ma}
u_{\theta}(t) = \boldsymbol{sat}_{-u_{\theta, max}}^{\phantom{-}u_{\theta, max}} \left[ v\frac{\vveh - p(\abs{S},z)}{\vveh - \unk(z)} \Xi  +  v p(\abs{S},z) \sgn (S) \right],\\
\nonumber
 \text{where}\quad
 \Xi := \frac{\varkappa(\eta) \cos(\hat{\theta} )}{1-\varkappa(\eta)z} + \frac{d \chi}{dz}(z) \sin (\hat{\theta}  ), \quad
 S := \hat{\theta} - \chi(z).
\end{gather}

Here the lower bound $\varkappa_{max}$ on the reference path curvature is given by Eq.\eqref{kappa.nom}.
The functions $\chi(z)$, $p(\abs{S},z)$ and $\unk(z)$ are user selectable, subject to some restrictions described in Chapt.~\ref{chapt:pf}.

\begin{Assumption}
\label{ass,ref}
There exists a choice of the coefficients $\mu$ (with indices) from 
Eq.\eqref{v.nom}, Eq.\eqref{t.p} and functions $\chi(z)$, $p(\abs{S},z)$, and
$\unk(z)$, such that all assumptions from Chapt.~\ref{chapt:pf} hold.
\end{Assumption}

These assumptions employ some constants, which are as follows in the context of
this chapter: $\overline{\varphi} := \varphi_{max}, \varphi_{est}:=\varphi_{max},
\varkappa_{v} := \varkappa_{max}, \overline{\varkappa}:= \mu_{\varkappa}
\varkappa_{max}, \varkappa_{v,\bldu}:= \varkappa_{max} - w_{\theta,max}$.

\par
Constructive sufficient conditions for Assumption~\ref{ass,ref} to hold, as well
as details of the required functions may be found in Chapt.~\ref{chapt:pf}.
At the same time,
it may be noted that these conditions can be satisfied by proper choice of the
coefficients $\mu_u, \mu_{\theta}, \mu_{v}$ in
Eq.\eqref{v.nom}, provided that minor and partly unavoidable requirements are met.
According to Chapt.~\ref{chapt:pf}, the resultant tracking errors
$\overline{z}_{err} \geq |z|$ ($\overline{z}_{est}$ in notations of
Chapt.~\ref{chapt:pf}) and $\overline{\theta}_{err} \geq |\hat{\theta}|$ are
explicitly determined by these coefficients (for given parameters of the
vehicles) and are not influenced by both the path selected by TPM and the speed
profile $v(\cdot)$\footnote{To reduce the case at hand to the constant-speed
vehicle considered in Chapt.~\ref{chapt:pf}, a change of the
independent variable is employed; the time $t$ is replaced by the covered path $\int v \;
dt$. Due to the multiplier $v$ in the right-hand sides of Eq.\eqref{eq1},
Eq.\eqref{eq2}, Eq.\eqref{cl.ma}, this hint sweeps $v$ away and reduces the equations
to the form considered in Chapt.~\ref{chapt:pf}, with the only exception in
the saturation thresholds in Eq.\eqref{cl.ma}. They become $\pm
\frac{u_{\theta,max}}{v}$ and may vary over time. It is easy to see by
inspection that all arguments from Chapt.~\ref{chapt:pf} concerning the
vehicle with constant thresholds $\pm \varkappa_{max}$ remain valid for the case
at hand since $\frac{u_{\theta,max}}{v} \geq \varkappa_{max}$ ($\Leftarrow v
\leq v_{max}$).}. It is worth noting that typically $\overline{z}_{err}
\mu_{\varkappa}\varkappa_{max}\ll 1$ and $\overline{\theta}_{err} +
\varphi_{max} \ll \pi/2$ (see Chapt.~\ref{chapt:pf}).

\subsubsection*{Longitudinal Tracking}

Longitudinal tracking adjusts the vehicles speed so the correct position along the trajectory is
maintained throughout path tracking. Though the trajectory to be tracked is given by the sequence of way-points, it
can be easily interpolated on any interval between sampling times to give rise
to its continuous-time counterpart. The related variables will be marked by
$^\ast$.
Whereas PTM makes the curvilinear ordinate $z$ close to that $z^\ast =0$ of the
traced trajectory, the goal of LTM is to equalize the abscissas $\eta$ and
$\eta^\ast$. To this end, the following longitudinal tracking
controller is employed:

\begin{gather}
\label{longtrack}
u_v(t) = u_v^\ast(t) - u_{v,exc} \;\sgn \Big\{ k_0 \big[\eta-\eta^{\ast}(t)\big]+ k_1 \big[ v - v^\ast(t) \big] \Big\},
\\
\nonumber
\text{where} \quad u_{v,exc} :=u_{v, max} - u_{v, nom} \overset{\text{Eq.\eqref{v.nom}}}{=} \big( 1- \mu_u\big) u_{v, max} >0
\end{gather}

Here $k_i>0$ are tunable controller parameters. In order that the longitudinal control be realistic,
the vehicle should be controllable in the longitudinal direction, i.e., the available control range exceeds that of disturbances: $u_{v,max}> w_{v,max}$.

\begin{Proposition}
Suppose that the controller parameters $\mu_u  \in (0,1), k_0>0, k_1>0$  are chosen so that :
\begin{equation}
\label{c.par}
u_{v, exc}> w_{v,max}, \quad u_{v, exc} - w_{v,max} > \frac{k_0}{k_1}\left[ 2 v_{max}+w_{\eta,max}\right],
\end{equation}
where
\begin{equation}
\label{w.eta}
w_{\eta,max}:= \frac{\overline{z}_{err} \mu_{\varkappa}\varkappa_{max}+1/2 \left(\overline{\theta}_{err} + \varphi_{max} \right)^2}{1-\overline{z}_{err} \mu_{\varkappa}\varkappa_{max} }.
\end{equation}

 Then the maximum longitudinal deviation $|\eta-\eta^\ast|$ along the trajectory is bounded by $\frac{k_1 w_{\eta,max}}{k_0}$.
\end{Proposition}

\pf First put $\eta_{\Delta}(t):= \eta(t) - \eta^\ast(t), v_{\Delta}(t) := v(t) -
v^\ast(t), S:= k_0 \eta_\Delta + k_1 v_\Delta$ and note that due to Eq.\eqref{eq1}
and Eq.\eqref{eq21},

\begin{gather}
\label{eq4}
\dot{\eta}_{\Delta} = \frac{v \cos (\hat{\theta} - \varphi) }{1-
\varkappa(\eta) z } - v^\ast = v_{\Delta}+w_\eta,
\\
{\dot{v}_{\Delta} = \begin{cases} \sigma:= - u_{v,exc}\,\sgn (S) + w_v & \text{\rm if}\; 0 < v < v_{max} \;\text{\rm or}\; \left|
\begin{array}{l}
v=0
\\
\sigma >0
\end{array}
\right.
\;\text{\rm or}\; \left|
\begin{array}{l}
v=v_{max}
\\
\sigma <0
\end{array}
\right.
\\
-\bldu^\ast_v(t) & \text{\rm otherwise}
\end{cases}}
\nonumber
\end{gather}
In Eq.\eqref{eq4},

 \begin{equation*}
  w_\eta:= \frac{v \cos (\hat{\theta} - \varphi) }{1- \varkappa(\eta) z } - v
\Rightarrow |w_\eta| \overset{\text{Eq.\eqref{w.eta}}}{\leq} w_{\eta,max}.
\end{equation*} 

First it is shown in the domain $0<v<v_{max}$, the discontinuity surface $S=0$ is
sliding. Indeed, in a close vicinity of any concerned point, it follows:

\begin{multline*}
\dot{S}\,\sgn(S) = \left[k_0 \dot{\eta}_\Delta + k_1 \dot{v}_\Delta \right] \,\sgn(S) = k_1 ( - u_{v,exc} + w_v\,\sgn(S)) + k_0(v_{\Delta}+w_\eta)\,\sgn(S)
\\
\leq - k_1 (u_{v,exc} - w_{v,max}) + k_0(|v_{\Delta}|+|w_\eta|) \leq  - k_1 (u_{v,exc} - w_{v,max}) + k_0(2 v_{max}+w_{\eta,max})
\\
 \leq  - k_1 \left[\big( 1- \mu_u \big) u_{v, max} - w_{v,max} \right] + k_0(2 v_{max}+w_{\eta,max}) \overset{\text{Eq.\eqref{c.par}}}{<} 0.
\end{multline*}

According to the architecture of the proposed navigation system, the reference trajectory always starts at the real state. So tracking is commenced with zero error $\eta_\Delta=0, v_\Delta=0$, i.e., on the sliding surface. Hence tracking proceeds as sliding motion over the surface $S=0$ while $0 < v < v_{\max}$. During this motion,

\begin{multline*}
S=k_0 \eta_\Delta + k_1 v_\Delta=0 \Rightarrow v_\Delta = - \frac{k_0}{k_1} \eta_\Delta \overset{\text{Eq.\eqref{eq4}}}{\Rightarrow} \dot{\eta}_{\Delta} = v_{\Delta}+w_\eta = - \frac{k_0}{k_1} \eta_\Delta +w_\eta
\\
\Rightarrow \eta_\Delta(t) = \int_0^t e^{- \frac{k_0}{k_1}(t-\tau)} w_\eta(\tau)\; d \tau \Rightarrow |\eta_\Delta(t)| \leq \int_0^t e^{- \frac{k_0}{k_1}(t-\tau)} w_{\eta,max}\; d \tau \leq \frac{k_1 w_{\eta,max}}{k_0}.
\end{multline*}

Now the marginal cases where $v=0$ or $v=v_{max}$ are analyzed.
If a state where $v=0$ is reached with $\eta<\eta^\ast$, it follows $\sigma = -
u_{v,exc} \, \sgn  \,\left[ k_0(\eta-\eta^\ast)- k_1 v_\ast\right] +w_v >0$.
Then the formula for $\dot{v}_{\Delta}$ is identical to that in the domain
$0<v<v_{max}$, and the above analysis remains valid. If conversely, $\eta \geq
\eta^\ast$, the vehicle is steady $\dot{\eta}=0$, whereas $\dot{\eta}_\ast \geq
0$, so the error $\eta - \eta^\ast$ does not increase. If a state where
$v=v_{max}$ is reached with $\eta>\eta^\ast$, it follows $\sigma = - u_{v,exc}
\,\sgn  \,\left[ k_0(\eta-\eta^\ast)+ k_1 (v_{max} - v_\ast)\right] +w_v < 0$,
and the arguments concerning the domain $0<v<v_{max}$ remain valid. If
conversely, $\eta \leq \eta^\ast$, the vehicle moves at the maximal speed
$\dot{\eta}=v_{max}$, whereas $\dot{\eta}_\ast \leq \mu_v v_{max} < v_{max}$, so
the error $\eta - \eta^\ast$ decrease.
\par
Thus the longitudinal deviation does not exceed $\frac{k_1 w_{\eta,max}}{k_0}$ in any case. \epf

It should be noted that conditions Eq.\eqref{c.par} can always be satisfied by proper choice of the controller parameters $\mu_v,k_0,k_1$.
 In conclusion, it is noted that the overall positional deviation
 $\norm{\blds(t) - \blds^\ast(t)}$ between the actual and reference trajectories
 does not exceed $|z|+|\eta-\eta^\ast|$.
 So the estimate introduced in Sec.~\ref{subsec.ansm} can be
 taken in the form:

\begin{equation}
d_{trk}  = \overline{z}_{est} + \frac{k_1 w_{\eta,max}}{k_0}.
\end{equation}

\subsection{Implementation Details}

The selection of the trajectory from $\mathscr{P}$ is arbitrary -- any possible choice of cost functional will not alter the collision avoidance properties. In addition this cost functional could also be used with other optimization routines, which would likely result in similar results. 
For this implementation, the trajectory from $\mathscr{P}$ was chosen minimizing of the cost functional:

    \begin{equation}
    J := \norm{\blds^*(\tau|k) - T} - \gamma_0\abs{v^*(1|k)}.
    \end{equation}
    
Here $\gamma_0 > 0$ is a given tunable parameter. This minimization aims to achieve the best progression to the target
and to simultaneously increase the speed at the first time step, and was found to give excellent experimental results for both holonomic and unicycle vehicles.

\section{Simulations} 
\label{ch2:sim}

\begin{figure}[ht]
    \centering
        \includegraphics[width=0.8\mylength]{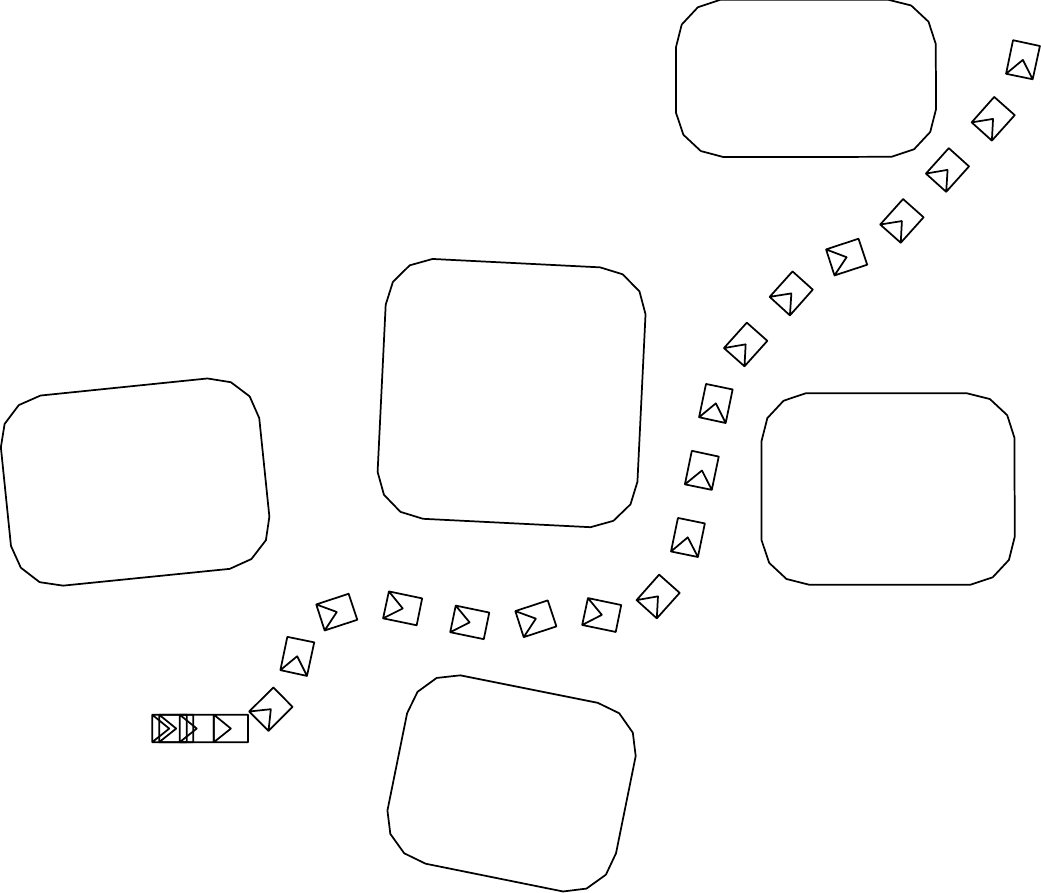}
    \caption{Trajectory obtained using the basic controller with a holonomic vehicle model.}
    \label{fig:holotraj}
\end{figure}

\begin{figure}[ht]
    \centering
        \includegraphics[width=1.2\mylength]{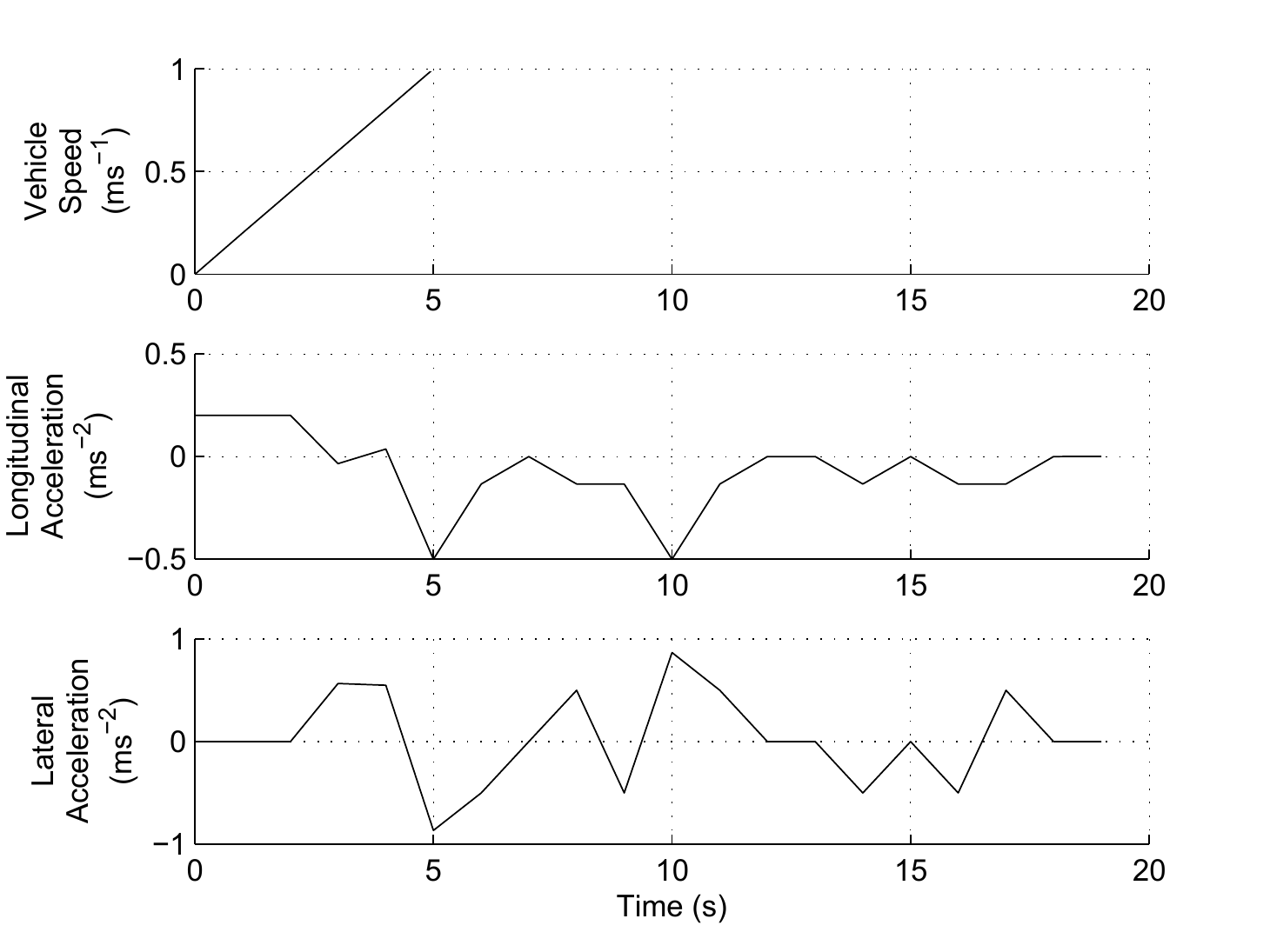}
    \caption{Control time history of the basic controller with a holonomic vehicle model.}
    \label{fig:holotimehist}
\end{figure}

\begin{figure}[ht]
    \centering
        \includegraphics[width=0.8\mylength]{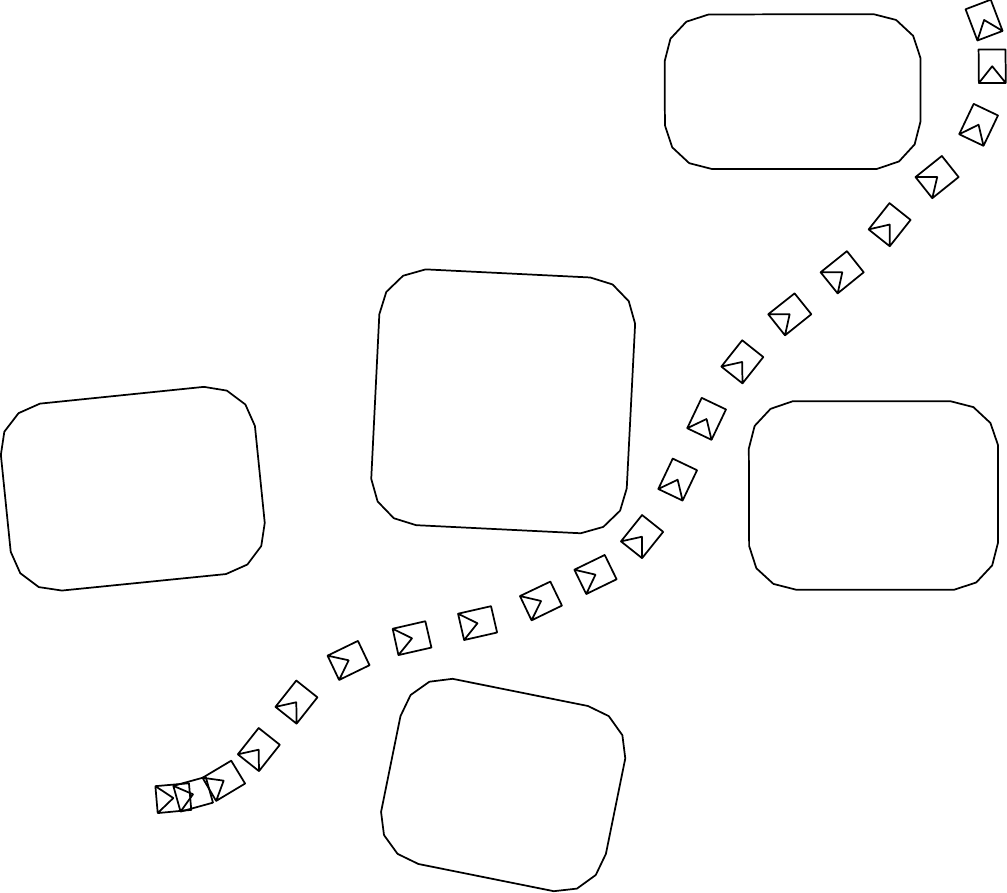}
    \caption{Trajectory obtained using the basic controller with a unicycle vehicle model.}
    \label{fig:unitraj}
\end{figure}

\begin{figure}[ht]
    \centering
        \includegraphics[width=1.2\mylength]{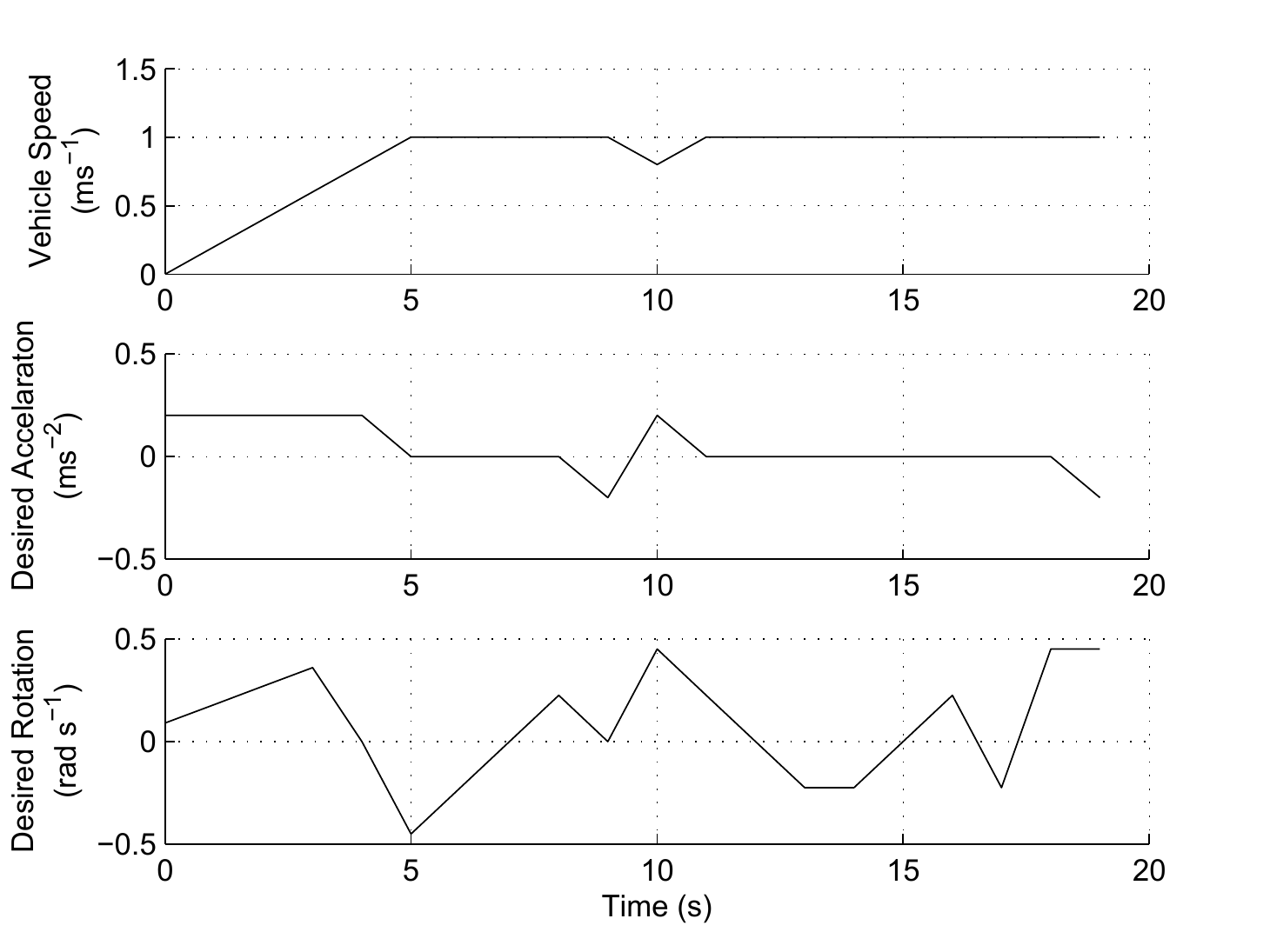}
    \caption{Control time history of the basic controller with a unicycle vehicle model.}
    \label{fig:unitimehist}
\end{figure}

For these simulations, both holonomic and unicycle vehicle models were employed, and the controller sampling time was set to unity. 
The parameters used for
navigation can be found in Tables~\ref{param} and \ref{paramuni}. The visible region $F_{vis}(k)$ was taken to be the 
points from $F$ within $4 m$ of the vehicle, which also had line-of-sight visibility to the vehicle through $F$. 
As may be seen in Figs.~\ref{fig:holotraj} to \ref{fig:unitimehist}, the proposed method was able to guide the vehicle satisfactorily without collision for both the vehicle models.

\begin{table}[ht]
\centering
 \begin{tabular}{| l | c |}
 \hline
 $u_{max}$ & $1.0 ms^{-2}$ \\
\hline
 $u_{nom}$ & $0.5 ms^{-2}$  \\
\hline
 $u_{sp}$ & $0.2 ms^{-1}$  \\
\hline

 \end{tabular} \hspace{10pt} \begin{tabular}{| l | c |}
\hline
 $d_{tar}$ & $0.5 m$  \\
\hline
 $v_{max}$ & $1.0 ms^{-1}$  \\
\hline
 $\gamma_0$ & $10$  \\
\hline
\end{tabular}

\caption{Simulation parameters for holonomic controller.}
\label{param}
\end{table}

Figs.~\ref{fig:holotimehist} and \ref{fig:unitimehist} show that in this case, the vehicle accelarated to the maximum speed and mostly stayed around this speed. In later chapters (Chapts.~\ref{chap:convsingle} to \ref{chap:dead})
where more complicated scenarios are considered, more speed fluctations will be required to give the vehicles more maneuverability when avoiding obstacles.

\begin{table}[ht]
\centering

\begin{tabular}{| l | c |}
\hline
 $d_{tar}$ & $0.5 m$  \\
 \hline
 $v_{max}$ & $1.1 ms^{-1}$ \\
  \hline
  $u_{\theta, max}$ & $0.6 rads^{-1}$ \\
    \hline
  $u_{v,max}$ & $0.3 ms^{-2}$ \\
  \hline
 $v_{nom}$ & $1.0 ms^{-1}$ \\
 \hline
  $u_{\theta, nom}$ & $0.5 rads^{-1}$ \\
 \hline
 \end{tabular} \hspace{10pt} \begin{tabular}{| l | c |}
 \hline
 $u_{v,nom}$ & $0.2 ms^{-2}$ \\
\hline
 $\Delta \Lambda$ & $0.5$  \\
 \hline
 $k_0$ & $1$  \\
 \hline
 $k_1$ & $1$  \\
\hline
$\gamma_0$ & $10$ \\
\hline
$\mu_{\varkappa}$ & $0.9$ \\
\hline

\end{tabular}
\caption{Simulation parameters for unicycle controller.}
\label{paramuni}
\end{table}

\clearpage \section{Summary}
\label{ch2:con}

In this section, the basic navigation framework which achieves collision avoidance is presented.
By combining a simplified trajectory planning scheme and a robust MPC type navigation approach, 
collision avoidance can be proven under bounded disturbance, for both holonomic and unicycle 
vehicle motion models, when static obstacles are present.

\chapter{Collision Avoidance with Limited Sensor Information}
\label{chap:convsingle}

In this chapter, the previous navigation system is extended to cases where only
sampled information about the obstacle boundary is available, similar to what
would be obtained from a LiDAR device.
Here the problem of following an obstacle boundary in a convergent
manner is considered, however this approach could easily be extended to general collision avoidance, and a
possible method of achieving target convergent navigation is outlined.
Additionally, a possible extension allowing avoidance of moving obstacles is
provided.
 \par
 By forming some assumptions about both the shape of the obstacle and the resolution
 of the obstacle sensor, constraints suitable for navigation using a MPC type approach are generated. 
An algorithm for selecting a target point to enable navigation in-line with
the goal of boundary following is also proposed. Overall this leads to a navigation system
more suitable for real world implementation, compared with the system from Chapt.~\ref{chap:singlevehicle}.
\par 
The body of this chapter is organized as follows. In Sec.~\ref{ch3:ps} the problem is formally defined, and in Sec.~\ref{ch3:ca} the navigation
approach is described. A method of extending the controller to target convergence problems is
described in Sec.~\ref{ch3:tar}, and to moving obstacle in Sec.~\ref{ch3:mov}. Simulations and
experiments are presented in Secs.~\ref{ch3:sim} and \ref{ch3:exp}. Finally, brief conclusions are given
in Sec.~\ref{ch3:con}.

\section{Problem Statement}
\label{ch3:ps}

A single vehicle with unicycle kinematics traveling in a plane is
considered, which contains a single static obstacle $D$. It is assumed
$D$ is simply connected, and thus does not enclose the vehicle. The
boundary of $D$ is denoted $\partial D$, which is assumed to have
continuous curvature. The desired \textit{transversal direction} is
given, which is encoded by the variable $\Gamma \in \{-1, +1\}$, where
clockwise corresponds to $-1$ and counterclockwise corresponds to
$+1$. The opposite direction $\bar{\Gamma}$ is referred to as the
\textit{anti transversal direction}.
\par
The distance to be maintained from any possible obstacle point during
planning is $d_{rad}$.  Due to the quantized nature of the range data
available, the navigation approach cannot
plan vehicle positions arbitrarily close to the obstacle.  Thus a user
selectable threshold $d_{tar} > d_{rad}$ is introduced, and navigation
between two adjacent obstacle segments is guaranteed if the minimum
distance between them is greater than $2d_{tar}$.  The value of
$d_{tar}$ partially determines the sensing and obstacle requirements,
and the margin between $d_{tar}$ and $d_{rad}$ is denoted $d_{ob} :=
d_{tar} - d_{rad}$.
\par
In this problem, the obstacle boundary $\partial D$ should be tracked
so that any point on $\partial D$ is within a distance $2d_{tar}$ of
at least one detected point on the obstacle at some time
during the transversal. In a obstacle avoidance problem this could be
interpreted as assurance that any corridors, possibly leading to the
target, will be transversed.
\par
In this chapter, only the unicycle model is considered for brevity, through these results could easily be extended to the holonomic model.
Also for brevity, the system is assumed free from disturbance -- extending these results to cases where disturbance 
is present remains a topic of future research. As the open loop control inputs can be applied directly without employing the TTM, there is no need
for measurement of the vehicles position $\blds(k)$ and heading $\theta(k)$. However it is assumed
 at every time $k$, every vehicle has knowledge of its speed $v(k)$ (as required by the TPM).
 
\subsection{Obstacle Requirements}

It is necessary to characterize the obstacle so that the resolution
requirements of the range sensor can be formulated. More contorted
obstacles would have greater sensing requirements than smooth
obstacles in order to navigate successfully.  The first requirement is
concerned with preventing collisions:

\begin{Requirement}
\label{as:shape}
The distance $d_{ob}$ is required to be sufficiently large such that for any point $\bldp_1 \in \partial
D$, there exists another point $\bldp_2$ separated by a distance $d$, $0 < d \leq d_{ob}$, such that
exactly one of the two boundary segments connecting them lies completely within $d_{ob}$ of both
points for any $\bldp_2$.
\end{Requirement}

\begin{figure}[ht]
\label{fig:shape}
\centering
\includegraphics[width=10cm]{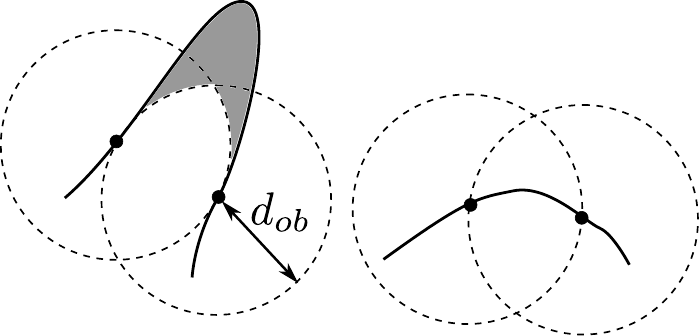}
\caption{Obstacle boundary assumption related to $d_{ob}$.}
\end{figure}

Requirement~\ref{as:shape} is demonstrated in Fig.~\ref{fig:shape}. This excludes narrow obstacle
protrusions, which may not be detected in some circumstances and cause a collision.

The next definition is similar to Requirement~\ref{as:shape} in that it provides a bounded area which
the obstacle boundary must lie within between two points. While the last definition was focused on
safety, the following ensures the overall navigation objective can be achieved:

\begin{Requirement}
\label{as:all}
The distance $d_{tar}$ is required to be sufficiently large such that for any point $\bldp_1 \in \partial
D$, there exists another point $\bldp_2$ separated by a distance $d$, $0 < d \leq 2d_{tar}$, such that
at least one of the two boundary segments connecting the points must lie completely within
$2d_{tar}$ of either point for any $\bldp_2$.
\end{Requirement}

\begin{figure}[ht]
\label{fig:asall}
\centering
\includegraphics[width=7cm]{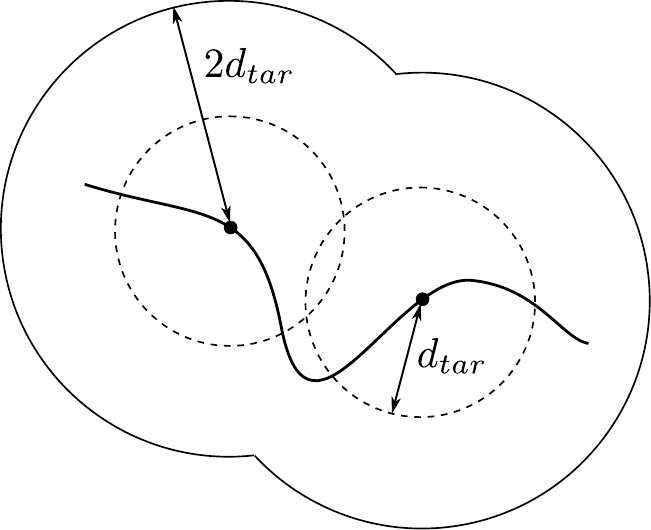}
\caption{Obstacle boundary assumption related to $d_{tar}$.}
\end{figure}

This requirement is demonstrated in Fig.~\ref{fig:asall} and avoids
long narrow corridors which may be skipped by the vehicle resulting in
incomplete navigation.

\subsection{Sensor Model}

\label{sec:sensor}
At every time step $k$, the vehicle emanates a finite sequence of $N_r$ rays radially around it, with
a constant angular spacing $\Delta \phi = \frac{2\pi}{N_r}$ between them. The vehicle has knowledge of the point
$R_i$ for which each ray first intersects an obstacle, $i = 1:N_r$. These points are taken from the
vehicles reference frame, so that the distance from $R_i$ to the vehicle may be given by $\lVert R_i\rVert $.
Also, for convenience the notation is wrapping, i.e. $R_i \equiv R_{i \pm N_r}$.

\begin{Assumption}
\label{as:range}
It is assumed that the range sensor has a maximum detection radius of $R_{sen}$, which means $R_i$ is
unknown if $\lVert R_i\rVert  > R_{sen}$. This maximum range is tightened based on the obstacle feature
size, and the effective maximum detection radius $R_{max}$ is defined as follows:

\begin{equation}
\label{eq:gamma}
R_{max} := \min \left\{ R_{sen}, \frac{d_{ob}}{\sqrt{\frac{8}{3}(1 - \cos(\Delta \phi))}}\right \}
\end{equation}

It is also assumed $R_{max}$ is greater than $4d_{tar}$.
\end{Assumption}

The denominator in Eq.\eqref{eq:gamma} arises during analysis; its origin
may be found in the proof of Lemma~\ref{prop:ok}.

\subsection{Concluding Remarks}

\begin{figure}[ht]
  \label{fig:explain}
  \centering
  \includegraphics[width=9cm]{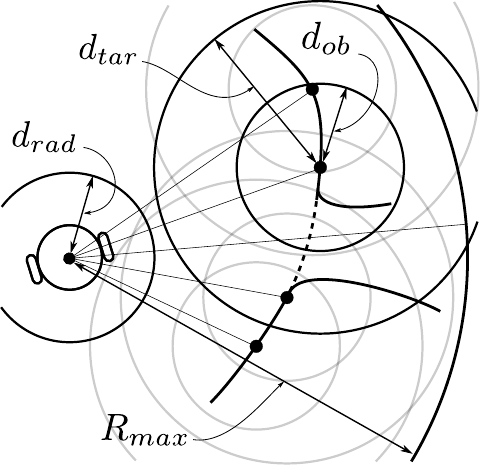}
  \caption{Relationship between the distances employed by the
    navigation law.}
\end{figure}

The distances employed by the navigation law are summarized in
Fig.~\ref{fig:explain}.  In this figure, since all circles of radius
$d_{tar}$ overlap, Requirement~\ref{as:all} applies to every pair of
points. This means the protruding part of the obstacle must lie within
$2d_{tar}$ of every pair of points.  Also, the circles of radius
$d_{ob}$ cover subsequent detection points, except for the center
detection ray. This means between these points a bounded area which
contains the obstacle boundary may be expressed, according to
Requirement~\ref{as:shape}. A circle of radius $d_{rad}$ is also shown
around the vehicle, which must exclude any point which is potentially
part of the obstacle (see Algorithm~\ref{algorithm}).

\begin{Remark}
  Note Requirement~\ref{as:shape} is automatically fulfilled if the
  boundary curvature radius is greater than $\frac{1}{2}d_{ob}$. This
  limit would also eliminate small protrusions which may not be
  detected by the obstacle sensor, through it should be emphasized the
  curvature does not need not be bounded for
  Requirement~\ref{as:shape} to hold.  Also note
  Requirement~\ref{as:all} may be omitted if the entire boundary is
  not required to be covered by the vehicle.
\end{Remark}

\clearpage \section{Navigation System Architecture}
\label{ch3:ca}

The idea behind this approach is to make the vehicle move towards the edge of the points that are
known to be part of a solid obstacle boundary, while adjusting the speed to avoid collisions. As in Chapt.~\ref{chap:singlevehicle}
a type of receding horizon control is employed, but in this case the probational trajectories are constrained to achieve
the navigation objective.
\par
As in Chapt.~\ref{chap:singlevehicle}, a probational trajectory commencing at the current system state is found over a variable finite horizon at every
time instant $k$. As before, this probational trajectory is specified by a finite sequence of way-points, which
satisfy the constraints of the vehicle model. Probational trajectories
necessarily have zero velocity at the terminal planning time, which ensures availability of at least
one safe probational trajectory at subsequent time steps. 
\par
Trajectory generation may potentially use many different forms, through
in the implementation presented in this chapter, the planning
algorithm proposed in Chapt.~\ref{chap:singlevehicle} is used.
\par
There are two caveats for the trajectory planning system:

\begin{itemize}
  \item Firstly, regardless of the path planning method employed,assume the probational trajectory may be \textit{inherited} from the previous 
time-step. The reasoning is the constraints are time-altering, this it is unknown whether 
a feasible trajectory exists at any particular time;  the inherited trajectory was subject 
to constraints from a previous time-step and thus is independent of the current ones. 

\item Secondly, the algorithm may alternatively select specifically prescribed trajectory,
which is formulated in Lemma~\ref{lem:trajexist}.
This may be regarded as a contingent `recovery scheme', which can guarantee a
trajectory is always available when the vehicle is stationary.\footnote{This
is allowed to facilitate analysis, but
  it should be emphasized it was never required by the navigation 
  system while testing the implementation presented in
  Sec.~\ref{sec:implem}.}
\end{itemize}

The final algorithm consists of sequential execution
of the following steps:

\begin{enumerate}[{\bf S.1}]
\item Construction of the constraints based on the range detections. These include \textit{avoidance constraints} (see Sec.~\ref{sec:obcon}) and
\textit{convergence constraints} (see Sec.~\ref{sec:tarpoint}).

\item Arbitrary selection of a probational trajectory satisfying constraints.
  
\item Application of the first control $\mathscr{U}^*(0|k)$ related to the
selected probational trajectory to the vehicle.

\item $k:=k+1$ and go to S.1.

\end{enumerate} 

Note in this formulation, the best choice of the probational trajectory from those
which satisfy constraints is purely arbitrary.
The optimal trajectory based on the given sensor information would likely depend on guessing the shape of the unknown
part of the obstacle, and this characterization would likely be non-trivial (and may be an area of future research).
The simplified heuristic used in for testing in this chapter may be found in
Sec.~\ref{sec:implem}.
\par
It is assumed the vehicle initially can plan a feasible probational trajectory:

\begin{Assumption}
Initially at $k=0$, at least one detection ray intersects the obstacle, and the distance from the obstacle exceeds $d_{tar}$.
\end{Assumption}

Trajectory related variables depend on two arguments $(j|k)$, where $k$ is the time instant when the
probational trajectory is generated and $j \geq 0$ is the number of time steps into the future; the related
value concerns the state at time $k+j$. Note all points are given from the vehicles reference frame,
and the origin is denoted $O$.

\subsection{Avoidance Constraints}
\label{sec:obcon}
The first step in determining trajectory feasibility is finding the regions of the surroundings guaranteed to be obstacle free.
This deliberation is based on Requirement~\ref{as:shape}; the algorithm for determining this is
described as follows:

\begin{Algorithm}
\label{algorithm}
Given a point $\boldsymbol{p}$ which bisects the two adjacent obstacle detection rays associated with $R_i$ and
$R_{i+1}$ for some $i$. There are several cases:

\textbf{$\boldsymbol{R}_i$ and $\boldsymbol{R}_{i+1}$ defined:}
Assume $\norm{R_i} \leq \norm{R_{i+1}}$ (A similar argument follows if
$R_{i+1}$ is nearer the vehicle than $R_i$), and define $d_i := \norm{R_i - R_{i+1}} $;

{\bf 1)} If $\lVert \boldsymbol{p} \rVert  > (R_{max} - d_{ob})$, $\boldsymbol{p}$ is labeled a \textbf{potential obstacle}.

{\bf 2)} Else if $d_i \leq d_{ob}$, and the minimum distance from the line segment $\overline{O\boldsymbol{p}}$ to $R_i$ and the
minimum distance from $\overline{O\boldsymbol{p}}$ to $R_{i+1}$ are both less than $d_i$, $\boldsymbol{p}$ is labeled a potential
obstacle.

{\bf 3)} Else, construct a circle of radius $d_{ob}$ centered on $R_i$. Intersect this with the
detection ray associated with $R_{i+1}$ and choose the intersection point furthest from the vehicle,
labeled $\boldsymbol{q}$. If the minimum distance from $\overline{O\boldsymbol{p}}$ to $\boldsymbol{q}$ and the minimum distance from the line segment
$\overline{O\boldsymbol{p}}$ to $R_i$ are both less than $d_{ob}$, $\boldsymbol{p}$ is labeled a \textbf{potential obstacle}.

{\bf 4)} Otherwise $\boldsymbol{p}$ is not a potential obstacle.

\textbf{$\boldsymbol{R}_i$ defined and $\boldsymbol{R}_{i+1}$ not defined:}
Employ steps 1, 3 and 4 only.

\textbf{$\boldsymbol{R}_i$ and $\boldsymbol{R}_{i+1}$ not defined:}
Employ steps 1 and 4 only.

\end{Algorithm}

\begin{figure}[ht]
\label{fig:areas}
\centering
\includegraphics[width=10cm]{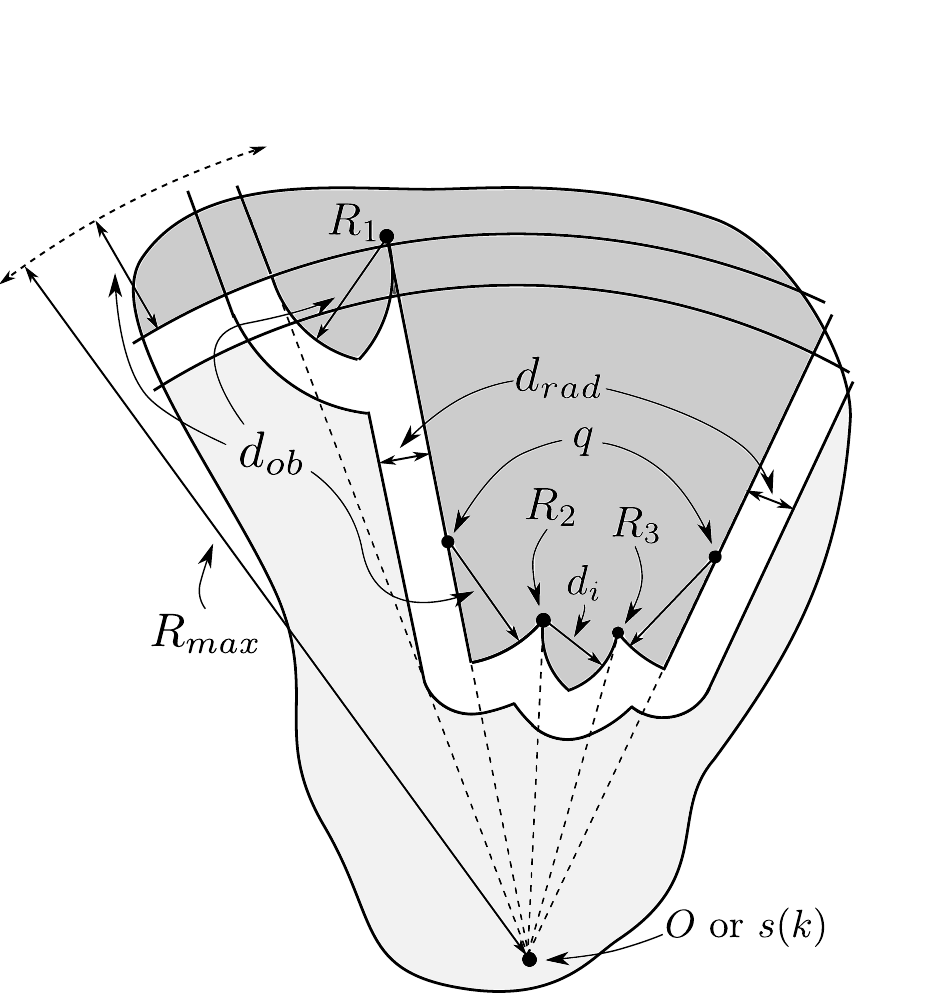}
\caption{The regions which are assumed to be part of the obstacle.}
\end{figure}

 Conversely the point is \textit{safe} if it
can be assumed to clear the set of potential obstacles by the minimum distance $d_{rad}$. 
An example of the results obtained using this method is illustrated in Fig.~\ref{fig:areas}.
In this figure the regions which must be assumed to be part of the obstacle are highlighted 
in dark gray, while the areas which can be assumed to be obstacle free are highlighted in light gray.

\begin{Proposition}
\label{prop:safe}
If point $\boldsymbol{p}$ is deemed not a potential obstacle by Algorithm~\ref{algorithm}, it is in fact not part of the obstacle.
\end{Proposition}

\pf First, note adjacent range readings are sufficient to determine if
a point is a potential obstacle as $\partial D$ cannot intersect the
line segment $\overline{OR_{i+1}}$. When $d_i < d_{ob}$ the result
follows directly from Assumption~\ref{as:shape}. When $d_i > d_{ob}$,
a section of $\partial
D$ 
must pass through a arc of radius $d_{ob}$ centered on $R_i$, between
the adjacent detection rays, on the far side from the vehicle.  In the
worst case it will pass through the intersection of the arc and
$\overline{OR_{i+1}}$, which is how the result follows from
Assumption~\ref{as:shape}. If the distance $\lVert
\boldsymbol{p}\rVert $ is greater than $R_{max} - d_{ob}$, and the
first two rules have not eliminated it, it indicates $\lVert R_i\rVert
$ and $\lVert R_{i+1}\rVert $ are both larger than $R_{max}$. The
sensor may not be able to detect obstacles further than $R_{max}$ as
they either are outside the range of the sensor or $\lVert R_{i} -
R_{i+1}\rVert > d_{ob}$, so that Assumption~\ref{as:shape} no longer
applies. \epf

Note that any potential probational trajectory must lie
completely with the set of currently sensed safe points.

\begin{Remark} \rm
Note a computational issue may be encountered when determining whether a probational trajectory meets avoidance constraints.
The easiest method of computing whether a probational trajectory is safe involves sampling it at a sufficiently high resolution, 
then applying the test from Proposition~\ref{prop:safe} on each point. By doing this, there are no sampling effects related to
the maximum velocity, as there were in Chapt.~\ref{chap:singlevehicle} 
\end{Remark}

\subsection{Convergence Constraints}
\label{sec:tarpoint}
The constraints in this section ensure the vehicle make at least a non-zero movement along the obstacle boundary.
To do this two definitions are introduced; the \textit{contiguous points}, which is a set of adjacent range detections near the vehicle,
and the \textit{target point}, which is calculated to ensure progression along the boundary. Firstly, the 
definition of the contiguous point set is introduced:

\begin{Definition}
The point $C_{end}(0)$ is a arbitrary valid obstacle detection taken before the navigation program commences.
The set of \textbf{contiguous points} $C(k)$ is defined by first selecting the obstacle detection
nearest $C_{end}(k-1)$. Once this point is found, adjacent
obstacle detections that are within an Euclidean distance of $2d_{tar}$ with respect to any element from the
current set of points can then be sequentially included. The final point to be included in $C(k)$ is labeled $C_{end}(k)$.

\end{Definition}

\begin{figure}[ht]
\label{fig:control}
\centering
\includegraphics[width=5cm]{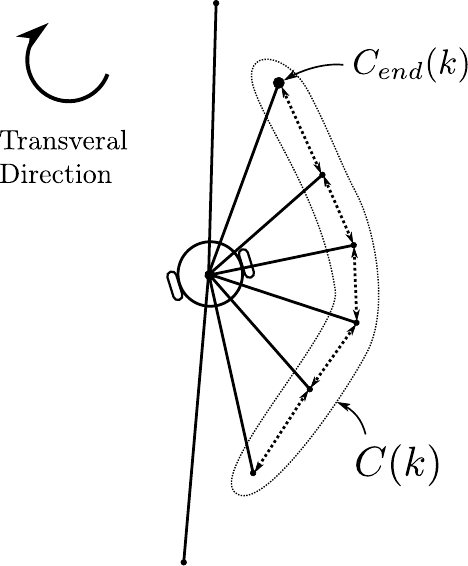}
\caption{Description of the contiguous set $C(k)$.}
\end{figure}

The purpose of finding $C(k)$ is that it can be guaranteed to be part of a obstacle
boundary -- it can contain no hidden corridors or gaps. Thus the vehicle can confidently navigate towards the
edge of this region, which may result in a shorter trajectory and faster allowable speeds. The set $C(k)$
and its endpoint $C_{end}(k)$ are illustrated in Fig.~\ref{fig:control}.

The \textit{target point} is taken from a circle of radius $d_{tar}$ around $C_{end}(k)$, intersected with the next detection ray after $C_{end}(k)$. This represents the
target offset distance which the vehicle may need to assume in order to track the obstacle boundary appropriately. 
The calculation of the target point is given as follows:

\begin{Definition}
\label{defak}
Let $R_n$ be the detected point corresponding to the next detection ray after  $C_{end}(k)$ 
in the anti transversal direction. The \textbf{target point} $A(k)$ is given by:

\begin{equation}
A(k) := \begin{cases} \hat{R}_n \left(\begin{array}{l}\norm{C_{end}(k)}\cos\Delta\phi\\ + \sqrt{d_{tar}^2 - (\norm{C_{end}(k)}\sin\Delta\phi)^2} \end{array}\right) & \norm{R_n} > \norm{C_{end}(k)}\\
C_{end}(k) & \norm{R_n} < \norm{C_{end}(k)}
\end{cases}
\end{equation}

Here $\hat{R}_n$ is the unit vector corresponding to $R_n$. 
\end{Definition}

\begin{figure}[ht]
\label{fig:safe}
\centering
\includegraphics[width=9cm]{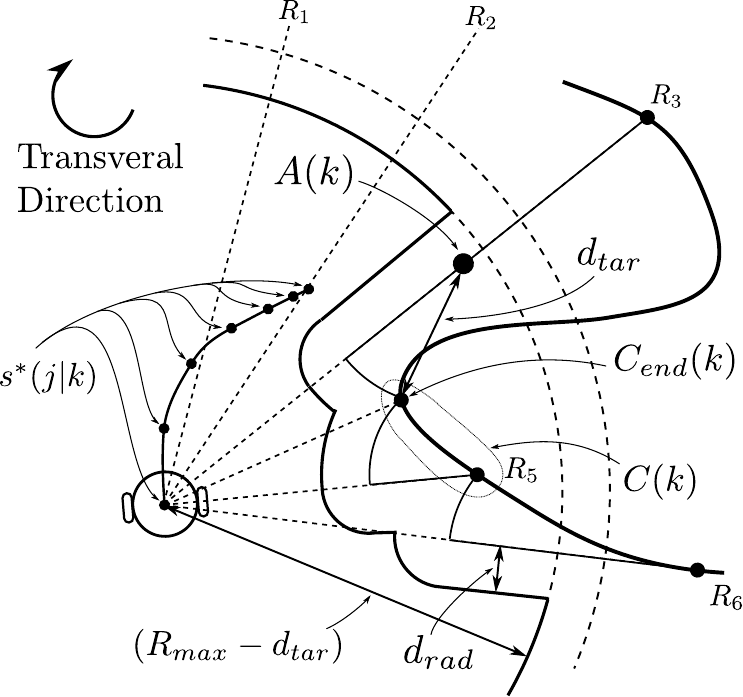}
\caption{Description of the target point $A(k)$.}
\end{figure}

The formula follows from elementary application of the cosine rule. The target point $A(k)$ is illustrated in Fig.~\ref{fig:safe} (note $A(k)$ is invariably an unsafe point).
To constrain the trajectory, the following convergence constraints are enforced:

\begin{enumerate}[{\bf r.1)}]
\item As the vehicle moves towards $A(k)$, there exists some $b > 0$, so that for all $j$
bounded by $1\leq j \leq \tau$:

\begin{equation}
\label{as:gotoa}
\lVert \blds^*(j|k) - A(k)\rVert  < \lVert A(k)\rVert  - b
\end{equation}

\item For each planned point on the probational trajectory, at least one pair of adjacent detection rays must enclose $C_{end}(k)$,
without occlusion from other potential obstacles.
\end{enumerate}

The second constraint would not normally significantly affect the operation
of the path planning system, however it is required to prove correct operation of the 
navigation law.

\subsection{Analysis}
\label{sec:anal}
In this section analysis of the navigation law in achieving the stated
goals is documented. The first step is to show existence and
smoothness of $C(k)$:

\begin{Lemma}
The set $C(k)$ is nonempty, the minimum pairwise distance to $C(k-1)$ is less than $2d_{tar}$, and
the point $C_{end}(k)$ exists.

\label{prop:ok}
\end{Lemma}

\begin{figure}[ht]
\centering
\includegraphics[width=7cm]{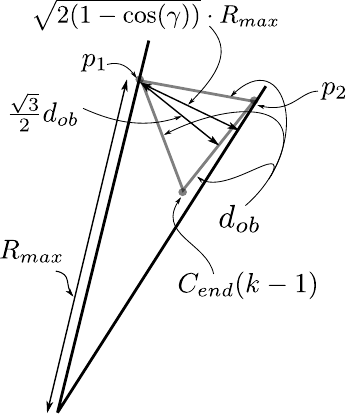}

\caption{Accompaniment to Lemma~\ref{prop:ok}.}
\label{fig:aendthere}
\end{figure}

\pf To show $C(k)$ is nonempty, consider a triangle of points in
$\partial D$, where one point is $C_{end}(k-1)$. The other two points
$\boldsymbol{p}_1, \boldsymbol{p}_2$ are chosen to bisect
$C_{end}(k-1)$ on $\partial D$, and also chosen so that
$\norm{C_{end}(k-1) - \boldsymbol{p}_1} = \norm{C_{end}(k-1) -
  \boldsymbol{p}_2} = d_{ob}$.  It is known $\norm{C_{end}(k-1)} <
(R_{max} - d_{tar})$, as either $\norm{C_{end}(k-1)} < 2d_{tar}$, or
$\norm{C_{end}(k-1)}$ decreased over the last time step. If there are
multiple choices for these two points $\boldsymbol{p}_1$ and
$\boldsymbol{p}_2$, then the points furthest away from $C_{end}(k)$
along $\partial D$ are chosen. It is clear that
$\norm{\boldsymbol{p}_1 - \boldsymbol{p}_2} \geq d_{ob}$, since if
this was not the case it would possible to select
$\acute{\boldsymbol{p}}_1$ and $\acute{\boldsymbol{p}}_2$, separated
by $d_{ob}$ and enclosing $\boldsymbol{p}_1$ and $\boldsymbol{p}_2$,
and for which both must be further than $d_{ob}$ from $C_{end}(k)$,
violating Assumption~\ref{as:shape}.  It is sufficient to show
intersection of this triangle corresponding to the points $\bldp_1$,
$C_{end}(k-1)$, $\bldp_2$ with at least one range ray. A trigonometric
argument using Assumption~\ref{as:range} is able to show intersection
for distances up to $R_{max}$. The relation is shown in
Fig.~\ref{fig:aendthere}; if one of the undetected points is nearer
than $R_{max}$ to the vehicle, it is evident $\sqrt{2(1-\cos(\Delta
  \phi))} \cdot R_{max} > \frac{\sqrt{3}}{2}d_{ob}$ which is
ascertained by Eq.\eqref{eq:gamma}, where it was assumed $R_{max} := \min
\left\{ R_{sen}, \frac{d_{ob}}{\sqrt{\frac{8}{3}(1 - \cos(\Delta
      \phi))}}\right \}$.

To show the minimum distance between $C(k)$ and $C(k-1)$ is less than
$2d_{tar}$, first note the vehicle will be able to detect the triangle
corresponding to the points $\bldp_1$, $C_{end}(k-1)$,
$\bldp_2$. Since at least one detected point will be between the
triangle and the vehicle, it will be within $2d_{tar}$ of the previous
$C_{end}(k)$.

Existence of $C_{end}(k)$ follows from Requirement~\ref{as:all}, as if
the distance between every pair of adjacent range points is less than
$2d_{tar}$ it indicates the obstacle is enclosing the vehicle, which
is assumed to not be the case. \epf

Next a guarantee that $A(k)$ is sufficiently spaced from the vehicle is provided. The argument
is based on the geometry of the detection rays:

\begin{Lemma}
\label{lem:speed}

The distance $\norm{A(k)}$ is lower bounded by $d_{tar}$.
\end{Lemma}

 \pf Consider a point $\bldp$ generated by moving $d_{tar}$ towards the origin from $A(k)$.
 If $\bldp$ is not the origin, the lemma follows.
If $\lVert R_n\rVert > \lVert C_{end}(k)\rVert$, using Definition~\ref{defak}
and the fact $d_{tar} > d_{ob}$, it can be inferred $\bldp$ is a potential obstacle and does not 
coincide with the origin. When $\lVert R_n\rVert \leq \lVert
 C_{end}(k) \rVert$, note $\norm{C_{end}(k)} > 2d_{tar}$ and thus
 $\norm{A(k)} > d_{tar}$. \epf

 Next there is shown to be at least one feasible trajectory that may
 be produced by the algorithm. Note this trajectory is not a
 particularly efficient method of transversing the obstacle, it is
 merely to ensure the analytical claims are valid.

\begin{Lemma}
\label{lem:trajexist}
There always exists at least one valid trajectory $\blds^*(j|k)$
\end{Lemma}

 \pf Since the trajectory from the previous time-step can always be
 used, the case where it is exhausted and the vehicle is stationary is
 considered. The recovery trajectory consists of a turn-in-place so
 the vehicle is aligned towards $A(k)$, followed by straight movement
 towards $A(k)$ of length $b$. The straight movement would be achieved
 by a given longitudinal acceleration at one time-step followed by an
 opposite deceleration at the proceeding time-step. Although the
 individual segments of this trajectory would not meet convergence
 constraints, the trajectory as a whole does. Thus the entire recovery
 trajectory would be executed until completion. \epf

Next it is shown that following feasible trajectories results in the point $C_{end}(k)$ progressing along the obstacle boundary:

\begin{Lemma}
\label{lem:advance}
There exists a finite number of time steps $c_t$ over which time $C_{end}(k)$ will advance by at least
$d_{ob}$ along the obstacle boundary. 
\end{Lemma}

 \pf Let $k_c$ refer to the time-step at the beginning of the interval
 for which advancement of $d_{ob}$ is being shown. Let the previous
 $C_{end}(k_c)$ and $A(k_c)$ be expressed in the vehicles current
 reference frame. Absurdly suppose $C_{end}(k)$ never advances by more
 than $d_{ob}$ along the boundary $\partial D$. It can be deduced that
 $C_{end}(k)$ will never regress by more than $d_{ob}$ along $\partial
 D$ since the triangle argument from Lemma~\ref{lem:speed} can be
 applied to $C_{end}(k_c)$. Since $\norm{A(k)} > d_{tar}$ it can be
 inferred that $\norm{C_{end}(k_c)}$ will eventually decrease below
 $2d_{tar}$ in finite time due to Eq.\eqref{as:gotoa} and
 Lemma~\ref{lem:trajexist}. It is known $C_{end}(k)$ lies in a disk of
 radius $d_{ob}$ around $C_{end}(k_c)$.  Next, define $\nu(k)$ to be
 the heading of $\blds(k)$ from a reference frame centered at
 $C_{end}(k_c)$, measured in the transversal direction from $k_c$ such
 that $\nu(k_c) := 0$.  Decreasing the distance to $A(k)$ will
 increase $\nu$ since all possible values of $C_{end}(k)$ are
 contained in a disk of radius $d_{ob}$ around $C_{end}(k_c)$. This
 movement will occur in finite time due to Eq.\eqref{as:gotoa} and
 Lemma~\ref{lem:trajexist}. Once $\nu > 2\pi$, the vehicle is
 guaranteed to sense a different side of the triangle corresponding to
 $C_{end}(k_c)$. Due to the triangle side length, this advancement of
 $C_{end}(k)$ along the boundary is at least $d_{ob}$. The expression
 for $c_t$ would depend on $b$ from Eq.\eqref{as:gotoa}.\epf

 Note that showing collision avoidance is elementary given the
 trajectories are always feasible. The final statement combines the
 above Lemma's to prove the final result.

\begin{Theorem}
\label{lem:beh}
The vehicle exhibits boundary following behavior, avoiding collision with the boundary while
making a range reading within $2d_{tar}$ of every
point on a finite obstacle boundary in finite time.
\end{Theorem}

\pf From the existence of $C_{end}(k)$ it follows the accumulated $C(k)$, given by:

\begin{equation}
\mathcal{C}(k) := \bigcup_{\acute{k} = 0}^{k} C(\acute{k})
\end{equation}

contains points within $2d_{tar}$ of every point on $\partial D$ between some start and end point.
From Lemma~\ref{lem:advance}, it follows that the span of $\mathcal{C}$ is monotonically
increasing.
From Lemma~\ref{lem:speed} and Lemma~\ref{lem:advance}, the position of $C_{end}(k)$ along the obstacle
increases by at least some increment in a finite time. Thus it may be inferred the span of
$\mathcal{C}(k)$ will increase arbitrarily in a finite time, and thus the obstacle will be transversed
in a finite time.\epf

\subsection{Implementation Details}
\label{sec:implem}

To implement the planning system, the same simplified planning system from
Chapt.~\ref{chap:singlevehicle} is employed, and a finite set of heuristically given
trajectories $\mathscr{P}$ is found.
As this planning algorithm is known not to be complete, whenever $\mathscr{P}$ 
is empty and the inherited probational trajectory is exhausted, the trajectory used in
Lemma~\ref{lem:trajexist} may be used in extraordinary circumstances (see Remark~\ref{tracks}).
\par
Only feasible trajectories may be chosen so any choice of cost functional will not alter the collision avoidance properties.  
For this implementation, the trajectory from $\mathscr{P}$ minimizing of the following cost functional was selected:

    \begin{equation}
    \label{sdsd}
  J := \lVert A(k) - s^*(\tau|k) \rVert
    \end{equation}
 
 Due to convergence constraints, any feasible trajectory will exhibit boundary following properties, 
thus this selection is heuristic. The formulation Eq.\eqref{sdsd} was found to give good performance when implemented.

\begin{Remark} \rm
\label{tracks}
Simulations have demonstrated the approach described is sufficiently complete for the proposed navigation law.
During the simulations and experiments in this chapter, a trajectory inheritance event was unable to be induced, so it seems somewhat unlikely
that trajectories from previous time steps will be routinely employed.
\end{Remark}

 \section{Extension to Target Convergence}
\label{ch3:tar}
It is reasonably straightforward to extend the proposed navigation
method to cases where a known target $T \notin D$ must be converged to
in finite time. The main modification is the addition of a trajectory
leading the target, which is followed whenever feasible. The reactive
navigation system is somewhat similar that described in Chapt.~\ref{chapt:rbf}.
However, the ability to
imitate the exact same switching conditions from `boundary following'
to `target pursuit' would be required to show the bi-similarity of the
algorithms, and an area of future research.

 \subsection{Angular Progression}
 Since the heading $\theta$ of the vehicle is not necessarily
 correlated with any parametric representation of the obstacle, other
 means are required to track the \textit{angular progression} along
 the obstacle boundary. Define $\alpha$ to be the relative heading
 from the vehicle between $C_{end}(k)$ and the target $T$:

\begin{equation*}
\alpha(k) := \angle(C_{end}(k) - \blds^*(0|k)) - \angle(C_{end}(k) - T) + \pi
\end{equation*}

Here $\angle(\cdot)$ represents the angle of a vector. The angle
$\alpha(k)$ evolves continuously over time, such that it may become
arbitrary large in certain types of obstacles with local minima having
complicated shapes. Thus a \textit{projected} value $\alpha_\tau(k)$
is defined, where $\blds^*(0|k)$ is replaced with the projected
position at the end of the trajectory $\blds^*(\tau|k)$. When
transitioning into boundary following mode, $\alpha$ may be reset to
lie within the interval $[-\pi, \pi]$.

\begin{figure}[ht]
  \centering \subfigure[]{\label{fig:beta}\includegraphics[width=4cm]{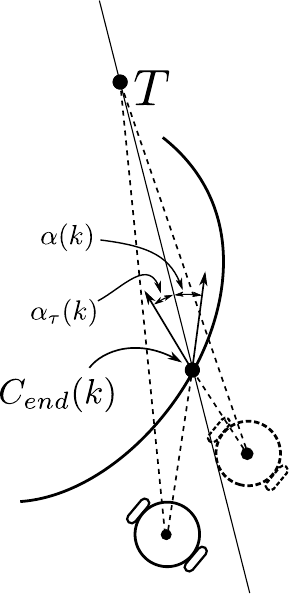}}
  \hspace{20pt} \subfigure[]{\label{fig:beta2}\includegraphics[width=4cm]{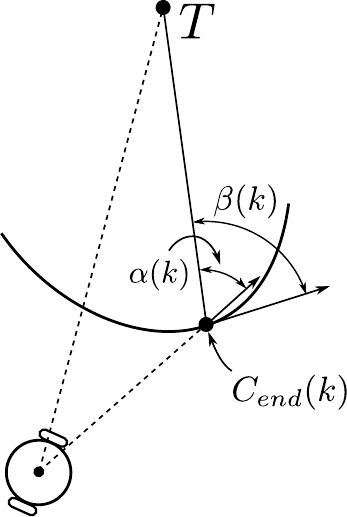}}
\caption{Relationship between $\alpha$ and $\beta$.}
\end{figure}

Define $\beta(k)$ to be the angle between a tangential ray from the
obstacle at $C_{end}(k)$ in the transversal direction, and the line
segment $\overline{C_{end}(k)T}$. Similarly to $\alpha(k)$, $\beta(k)$
evolves continuously over time, however it is impossible for the
vehicle to calculate the exact vale of $\beta(k)$ is in
general. However it is somewhat equivalent to the measurement used in
Chapt.~\ref{chapt:rbf} to calculate mode switching conditions, so by
showing $\alpha(k)$ and $\beta(k)$ are similar, it may mean some of
those results are relevant.

\begin{Proposition}
\label{prop:alpbet}
  The following inequalities hold: $|\alpha(k) - \beta(k)| < \pi$ and
  $|\beta(k)| \geq |\alpha(k)|$
\end{Proposition}

\pf The angle of the actual tangent ray of the obstacle at
$C_{end}(k)$ subtended the line segment $\overline{OC_{end}(k)}$,
labeled $\nu$, is contained in the interval $[0, \Gamma\pi]$, as any
other direction would require $C_{end}(k)$ to be a corner point, which
is precluded by Assumption~\ref{as:shape}. By summing the angles
around $C_{end}(k)$, $(\pi - \alpha(k)) + \beta(k) + \nu \equiv 2\pi$,
from which the proposition follows. \epf

In Fig.~\ref{fig:beta}, the switching condition based on $\Gamma
\cdot \alpha_\tau(k)$ is illustrated, while in Fig.~\ref{fig:beta2}
the correlation between $\alpha(k)$ and $\beta(k)$ is illustrated.

\subsection{Reference Trajectory}

Next, the \textit{reference trajectory} is defined by the choice of
navigation parameters which satisfy the following criteria:

\begin{itemize}
\item Monotonically decreases the distance to target relative to the
  current vehicle position.
\item The projected $\Gamma \cdot \alpha_\tau(k) < 0$; this is
  illustrated in Fig.~\ref{fig:beta2}.
\item Does not violate collision avoidance constraints.
\end{itemize}

If multiple valid choices are available, a selection is made to
minimize the distance to the target for the final trajectory
point. Whenever the reference trajectory exists, it is used in
preference to the trajectory obtained from the boundary following
method. However, the contiguous set would be updated as normal during
this mode (in some cases it may become empty, such as when departing
an obstacle).

 \subsection{Initial Selection of Contiguous Set}
 Whenever boundary following is engaged, an \textit{alternative
   contiguous set} is generated by first finding the set of obstacle
 detections within $d_{tar}$ of any point on the reference trajectory,
 then selecting the point nearest the vehicle as the seed point.

\begin{figure}[ht]

   \centering \subfigure[]{\label{fig:initsel}\includegraphics[width=3cm]{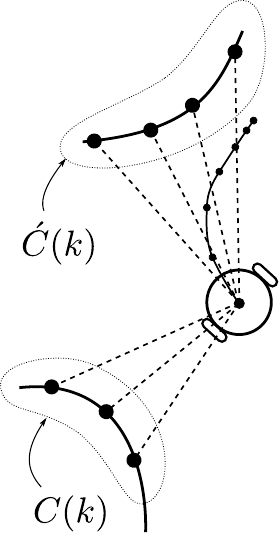}}
   \hspace{20pt} \subfigure[]{\label{fig:moving}\includegraphics[width=3cm]{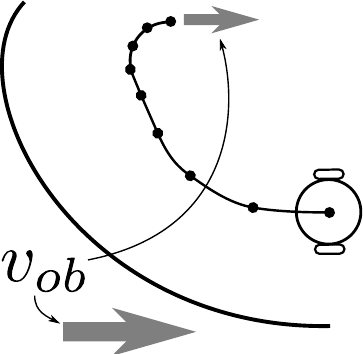}}
   \caption{(a) The alternate contiguous set; (b) Constraints for
     moving obstacles.}
\end{figure}

 If the alternative $\acute{C}(k)$ matches the current $C(k)$, or if
 the alternative contiguous set is empty, it can be presumed that the
 vehicle is still following the same boundary segment, and no changes
 to the transversal direction are required. However if they are
 different, it must assumed that the vehicle has encountered a new
 boundary segment; this is demonstrated in
 Fig.~\ref{fig:initsel}. In this case, the transversal direction
 can be randomly chosen in-line with Chapt.~\ref{chapt:rbf}.

 \begin{Remark} \rm It is not claimed that the proposed navigation law
   converges to the target in general. However due to the properties
   of $\alpha$, the switching condition is similar
   \cite{Matveev2011conf9}, so convergence to target seems likely, at
   least in conservative scenarios.
\end{Remark}

 \section{Extension to Moving Obstacles}
\label{ch3:mov}

 Allowance can be made for the obstacle to move with a constant
 translational velocity. This may have applications such as following
 large moving objects such as ships, or movement in environments with
persistent disturbances, such as water currents and wind.

\begin{Assumption}
  The obstacle is traveling at a constant speed $v_{ob}$ in constant
  direction $\theta_{ob}$, which is known to the vehicle. The obstacle
  velocity does not exceed that of the vehicles such that $v_{ob} <
v_{max}$.
\end{Assumption}

While the bounds on the obstacle speed are not restrictive, lower
obstacle velocities are likely to result in better performance of the
algorithm. The velocity of the obstacle could be easily estimated from
the range measurements. Many algorithms have been developed for this
task, an example being Iterative Closest Point matching
\cite{Zhang1994journ2}. However, during the simulations presented
here, the obstacle velocity was assumed to be known by the vehicle,
thus a obstacle tracking algorithm was not implemented.

Collision avoidance can be predicted by subtracting the future
movement of the obstacle from the future vehicle position, and
applying Proposition~\ref{prop:safe} to determine interference with
the obstacle as normal.  There is an additional modification to the
terminal constraint -- the vehicle must be traveling in the same speed
and direction as the obstacle at the end of the trajectory. To do
this, the speed along the trajectory is not reduced beyond
$v_{ob}$. The trajectory is then appended with a turn at maximal
actuation aligning the vehicle heading $\theta(j|k)$ with the
obstacle's direct $\theta_{ob}$, while maintaining the same
speed. This is demonstrated in Fig.~\ref{fig:moving}.

\begin{Remark} \rm It is impossible to design a trivial recovery
  scheme as was the case for the stationary obstacle.  Because of
  this, correct boundary following behavior of the vehicle cannot be
  shown. However collision avoidance with the obstacle can still be
  shown, as there is zero relative movement between the vehicle and
  obstacle at the end of the trajectory.
\end{Remark}

 \section{Simulations}
\label{ch3:sim}
The control law was simulated using the perfect discrete time model,
updated at a sample time of $1s$. The control parameters may be found
in Fig.~\ref{fig:param}. Simulations were carried out on a 3.0 GHz
Pentium processor, using MATLAB interfaced with C++ Mex files.

\begin{table}[ht]
\label{fig:param}
\centering
\begin{tabular}{| l | c |}
\hline
    $u_{\theta,max}$ & $0.5 rads^{-1}$ \\
\hline
  $u_{v,nom}$ & $0.2 ms^{-2}$ \\
\hline
 $v_{max}$ & $1.0 ms^{-1}$ \\
\hline
\end{tabular} \hspace{10pt} \begin{tabular}{| l | c |}
\hline
    $d_{ob}$ & $1.0 m$  \\
\hline
 $d_{rad}$ & $0.5 m$ \\
\hline
\end{tabular} \hspace{10pt} \begin{tabular}{| l | c |}
\hline
    $R_{max}$ & $5 m$  \\
    \hline
    $N_r$ & $40$  \\
\hline
\end{tabular} 
  \caption{Simulation parameters for boundary-following controller.}
\end{table}

\begin{figure}[ht]
\centering
  \includegraphics[width=6cm]{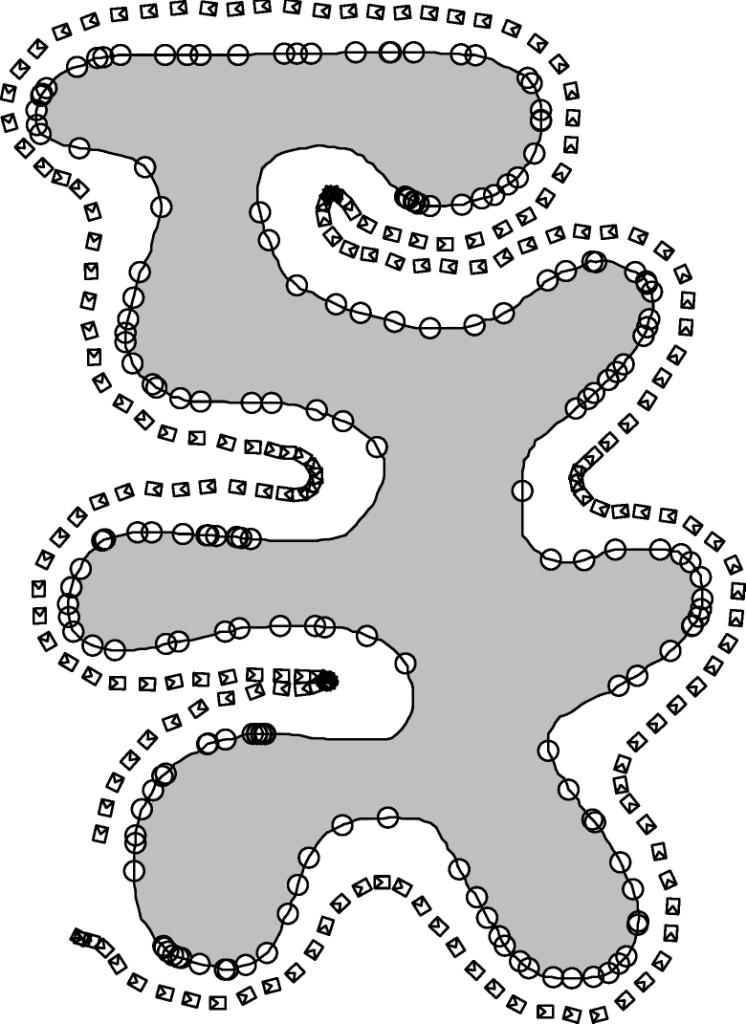}
  \caption{Trajectory for Simulations with a relatively simple obstacle.}
\label{fig:simple}
\end{figure}

\begin{figure}[ht]
\centering
  \includegraphics[width=1.2\mylength]{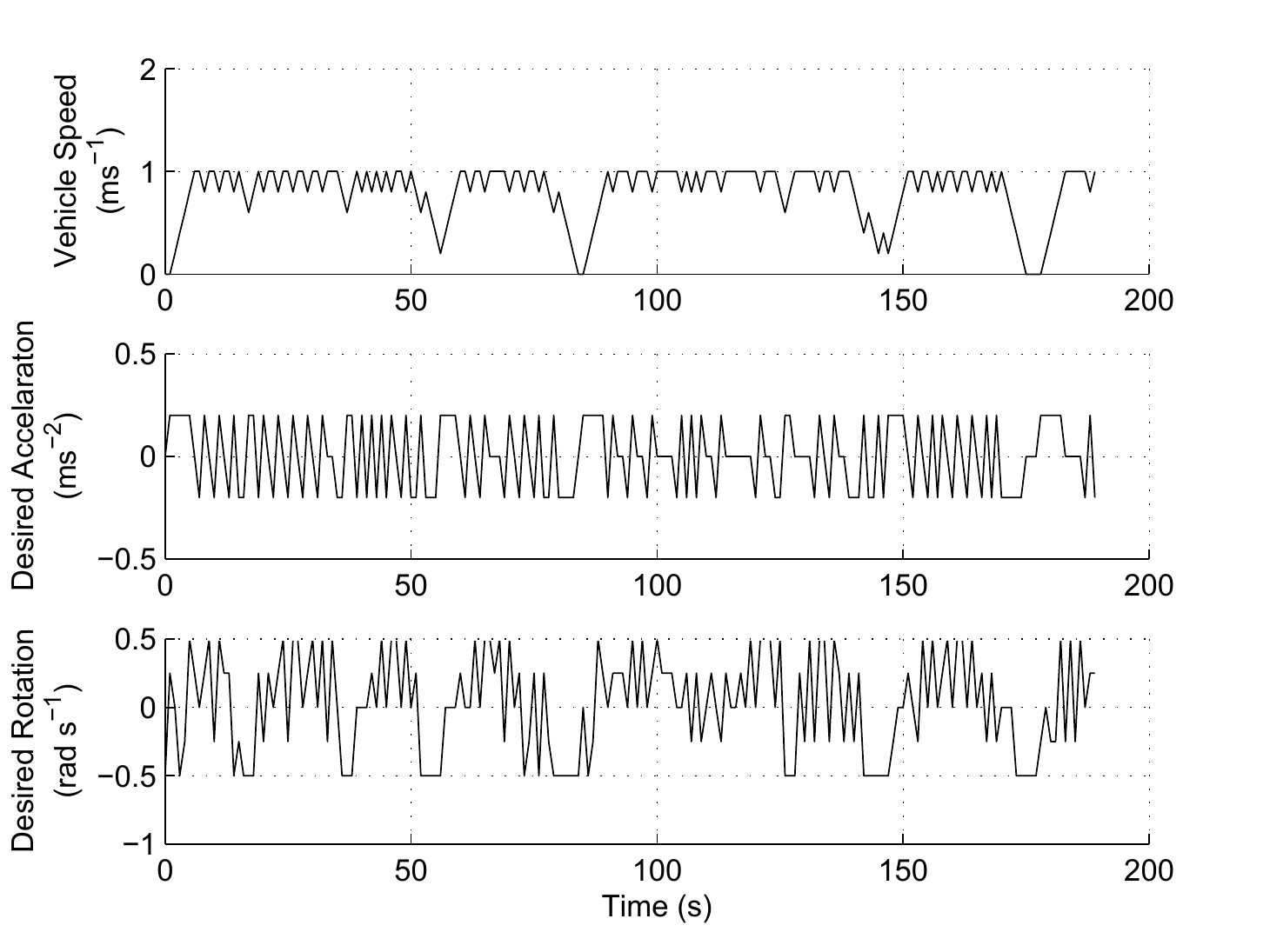}
  \caption{Control time history for simulations with a relatively simple obstacle.}
\label{fig:bfsimth}
\end{figure}

\begin{figure}[ht]
\centering
  \includegraphics[width=5cm]{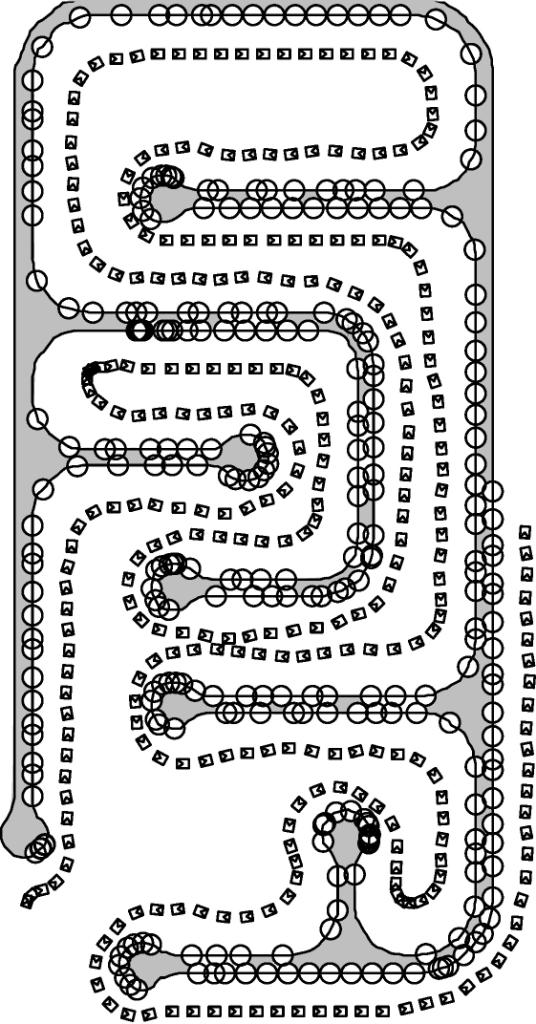}
\caption{Trajectory for simulations with a more complex obstacle.}
\label{fig:complex}
\end{figure}

\begin{figure}[ht]
\centering
  \includegraphics[width=1.2\mylength]{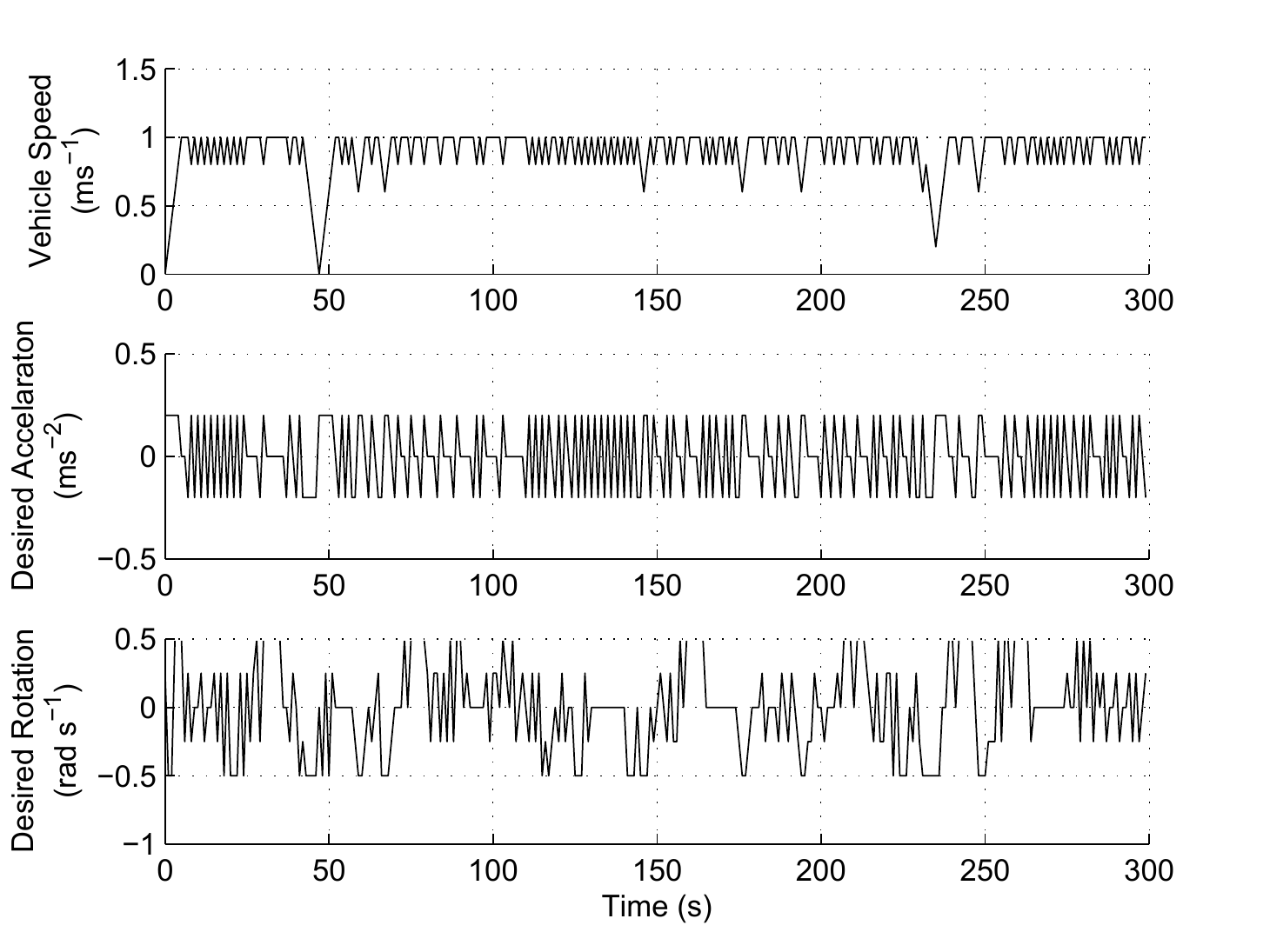}
\caption{Control time history for simulations with a more complex obstacle.}
\label{fig:bfcompth}
\end{figure}

It can be seen in Figs.~\ref{fig:simple} to \ref{fig:bfcompth} that
relatively straight segments were transversed efficiently and quickly.
The vehicle slowed down appropriately to transverse concave
corners -- the nominal turning radius of the vehicle decreases as the
vehicle slows down and allows tighter paths to be planned. 
In some situations there was some `overshoot' which is a
nuance of the path planning approach, and does not affect the claims
of navigation characteristics in any way. The average time taken to
compute the control law was 1.8~ms with a standard deviation of
0.7~ms.

\begin{figure}[ht]
\centering
  \includegraphics[width=6cm]{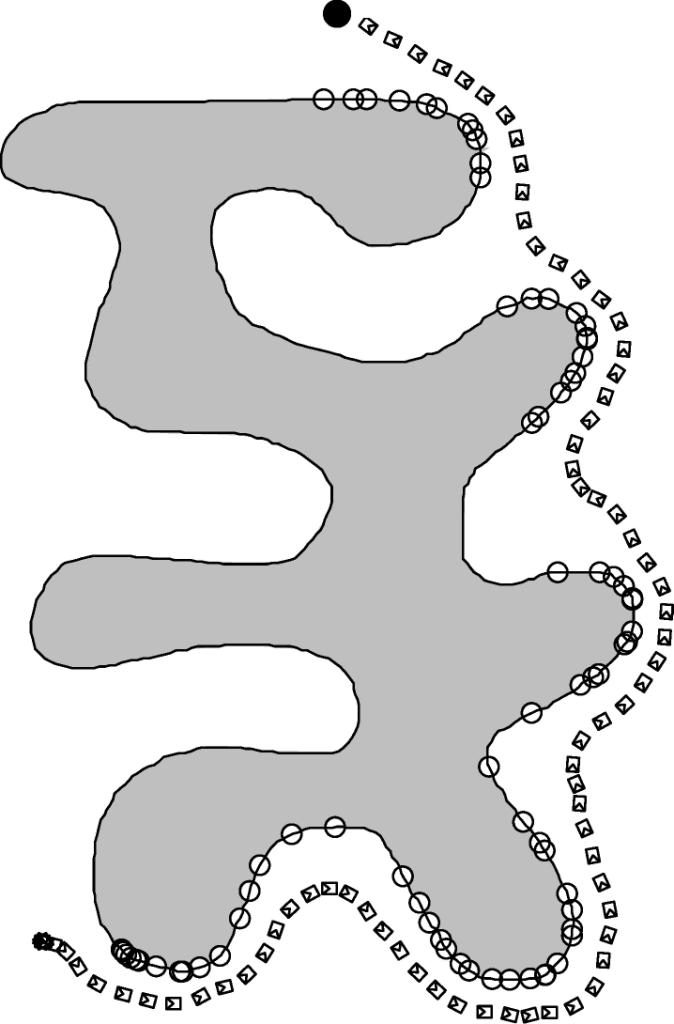}
  \caption{Trajectory for simulations with target convergence.}
\label{fig:target}
\end{figure}

\begin{figure}[ht]
\centering
  \includegraphics[width=1.2\mylength]{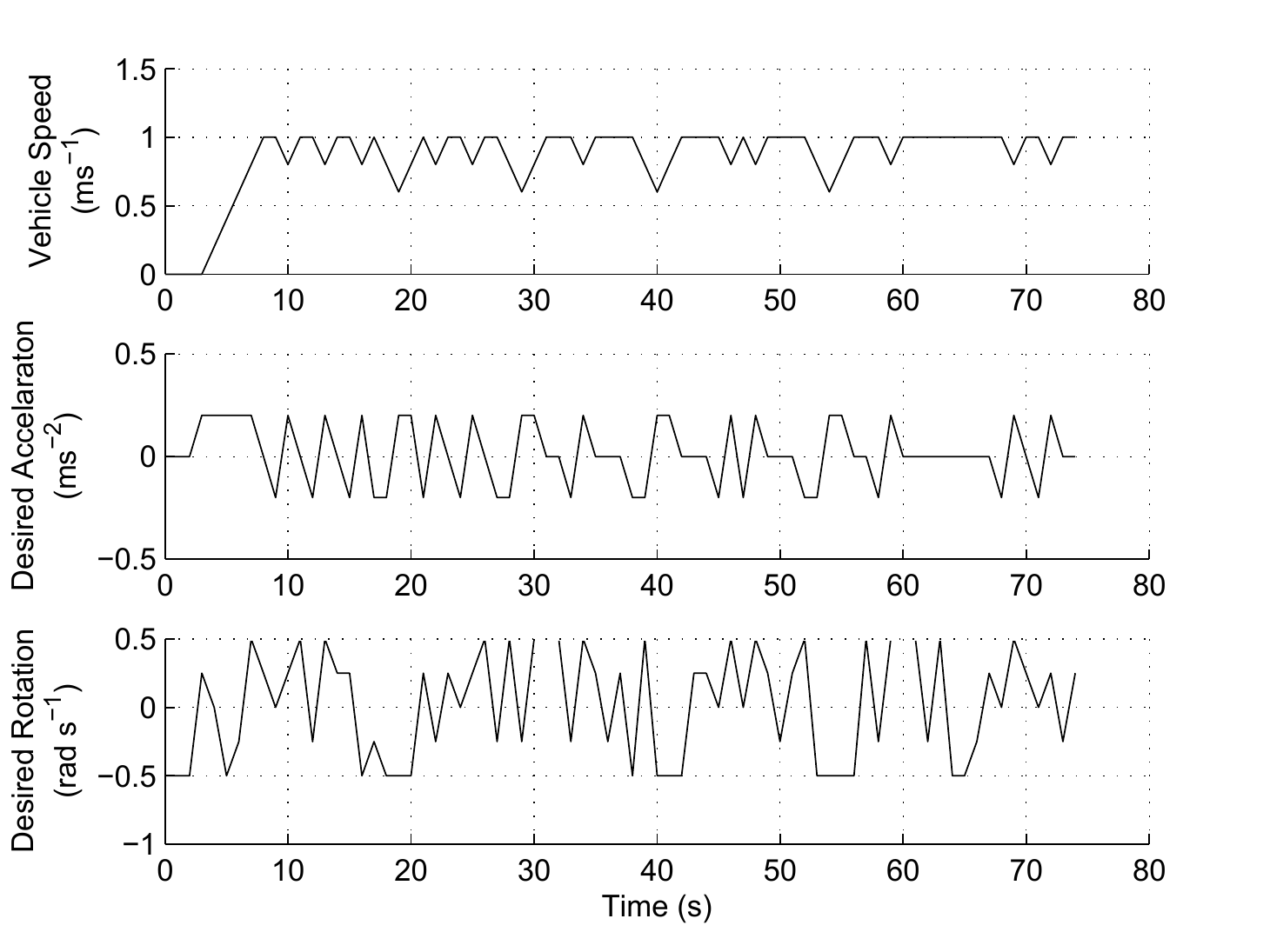}
  \caption{Control time history for simulations with target convergence.}
\label{fig:bfconvth}
\end{figure}

When a target point was defined (for the simulations shown in 
Figs.~\ref{fig:target} and \ref{fig:bfconvth}), the vehicle was
able to efficiently converge to the target. A random decision was made
after the two straight segments present, as there were insufficient
obstacle detections to guarantee the position of the obstacle boundary. 
This means a different value would lead to a different realization
of the trajectory, however it would still converge to the target. As
stated previously, high level decision making would provide a better
alternative to this simplified method.

\begin{figure}[ht]
  \centering \subfigure[]{\includegraphics[width=6cm]{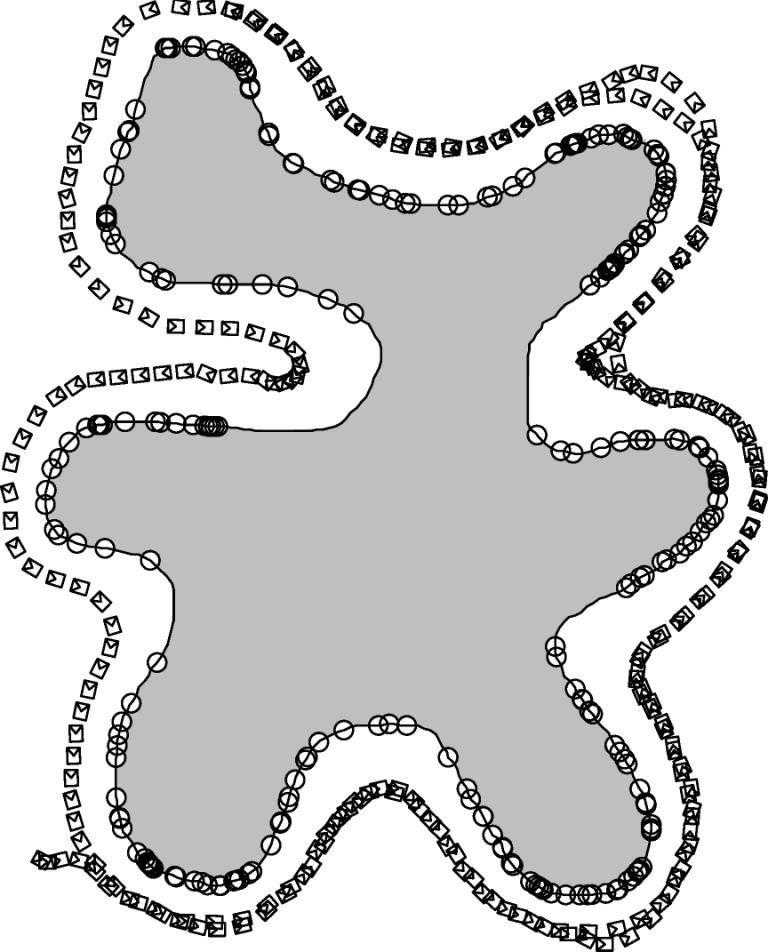}}
  \subfigure[]{\includegraphics[width=10cm]{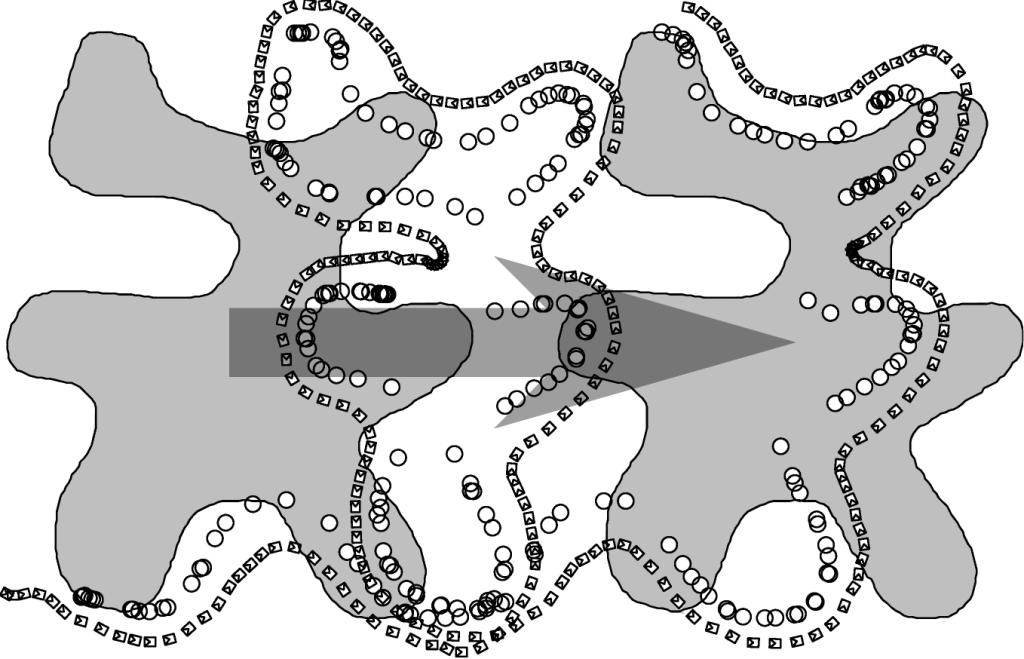}}
  \caption{Simulations where the obstacle is moving; (a) Obstacle's reference frame; (b) Vehicle's reference frame.}
\label{fig:move}
\end{figure}

\begin{figure}[ht]
\centering
  \includegraphics[width=1.2\mylength]{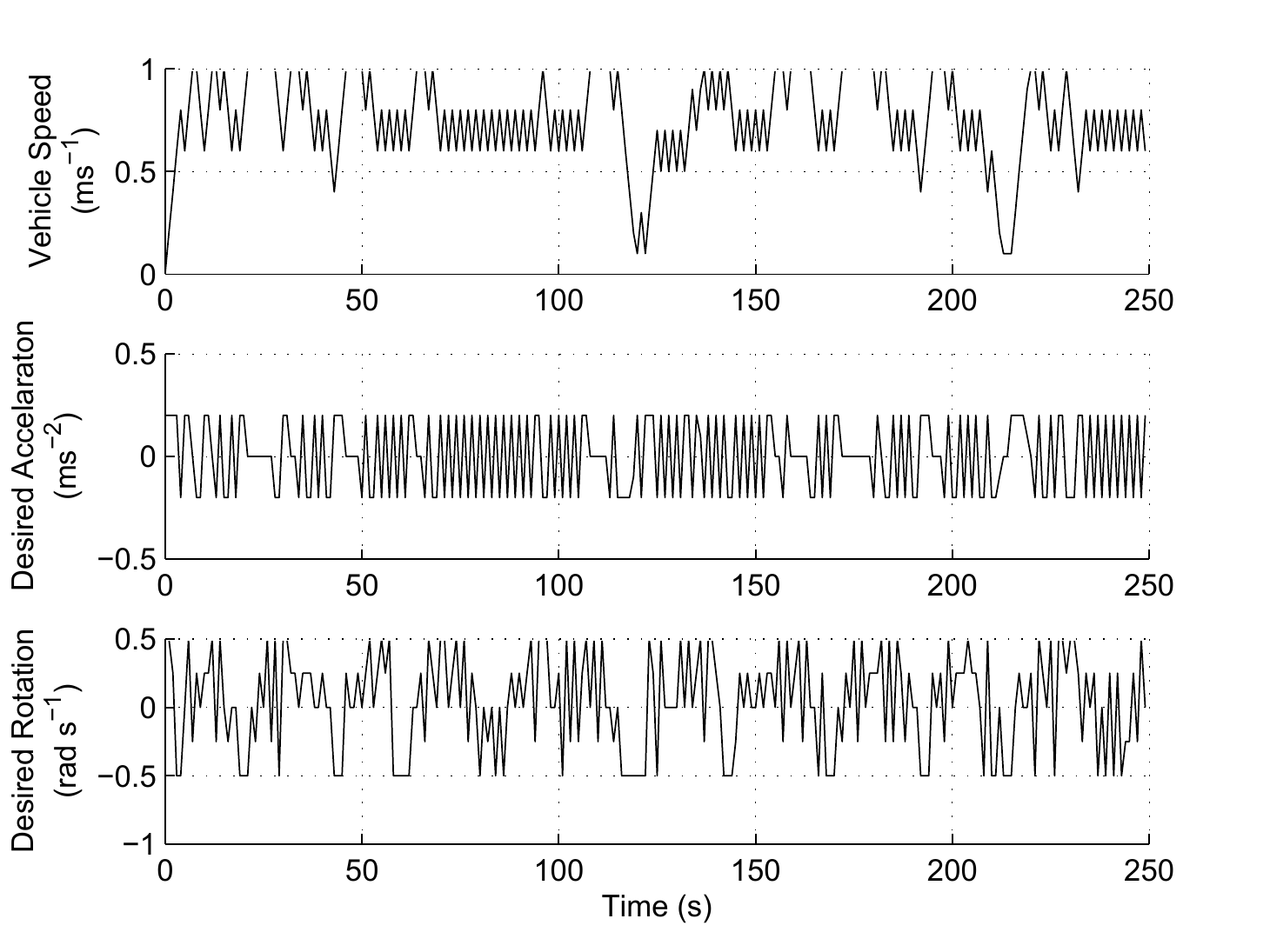}
  \caption{Control time history for simulations where the obstacle is moving.}
\label{fig:bfmovth}
\end{figure}

For Figs.~\ref{fig:move} and \ref{fig:bfmovth}, 
during simulation the obstacle moves to the right at a speed
of $0.1 ms^{-1}$ (the trajectory is drawn from both the obstacle's and the vehicle's reference
frame). The vehicle was given the precise velocity vector of the
obstacle, which negated the need to implement an obstacle tracking
algorithm.  The effect of the moving obstacle can be observed especially around the
points where the vehicle is moving slowly around corners - the motion
appears to be sideways in the moving frame. While the only assumption
on obstacle velocity was that it was below the maximum vehicle
velocity, higher obstacles velocities within this limit were observed
to reduce the likelihood of successful navigation (through safety was
not affected).

The proposed control law (PCL) was also compared with an existing
boundary following control law from the literature
\cite{Ge2005journ0}, which is referred to as the alternate control law
(ACL).  The method of generating target points for the ACL (called
`instant goals') is similar, but not identical, to the method employed
by the PCL. However, the key difference is that the ACL uses a
potential field type term to maintain the required distance from the
boundary, whereas the PCL employs path planning. The ACL also assumes
a velocity controlled holonomic vehicle model, and to highlight the
advantage of the PCL, the ACL was applied to the unicycle model. To do
this, the maximal control inputs tending the vehicle's heading and speed to the
output of the ACL were implemented. The tunable parameter $\xi$ was
taken to be 20, in-line with the recommendations of the method, and
$V_{opt}$ was taken to be the same as $v_{max}$ from the PCL.

\begin{figure}[ht]
  \centering \subfigure[]{\includegraphics[width=5cm]{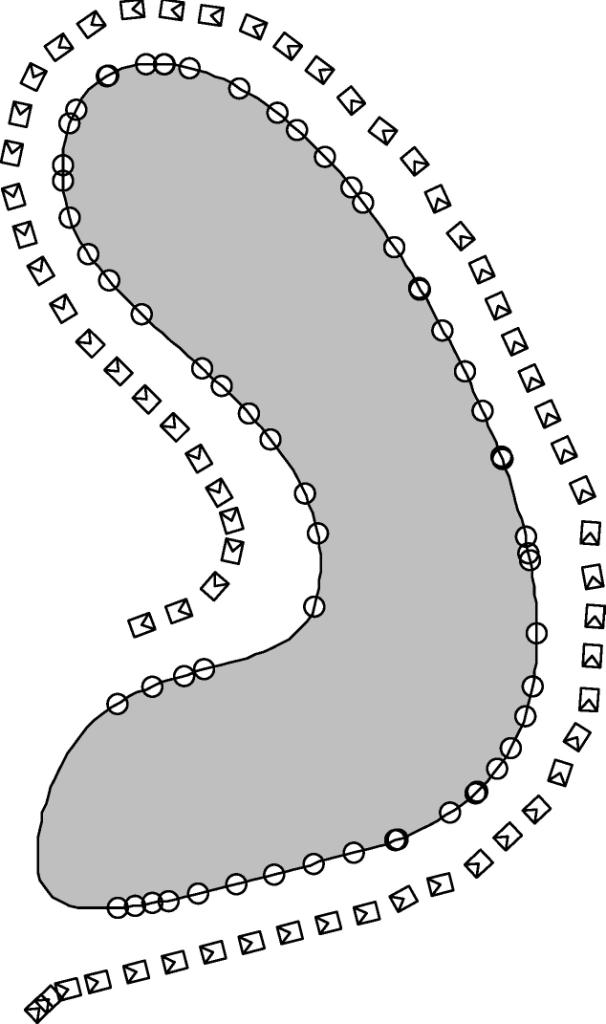}}
  \subfigure[]{\label{fig:move2}\includegraphics[width=5cm]{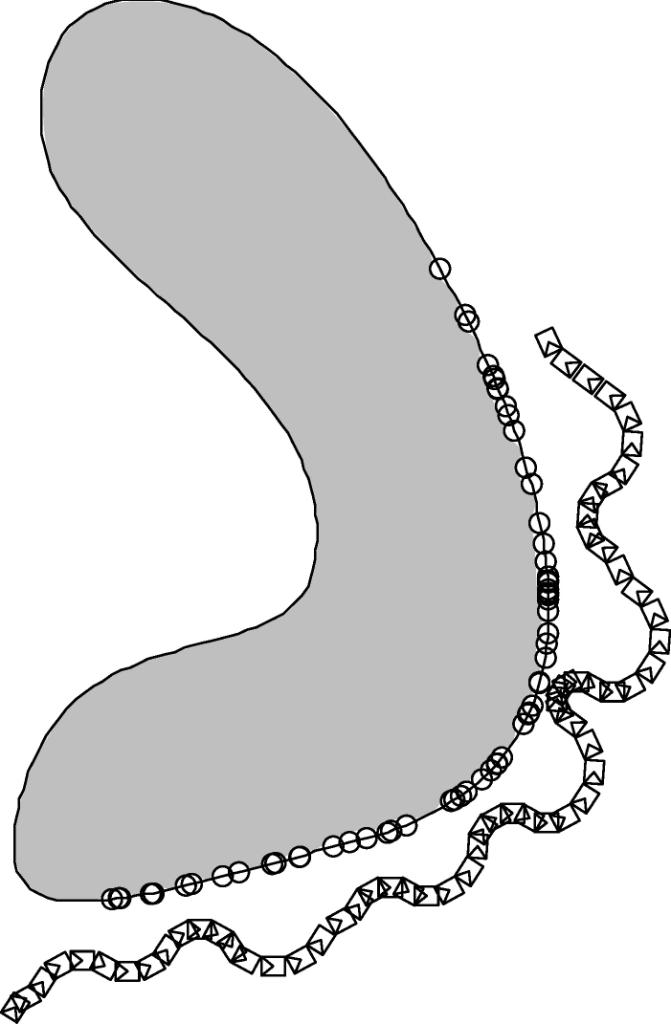}}
  \caption{Comparison which an alternate control law from the
    literature; (a) PCL; (b) ACL.}
    \label{fig:move1}
\end{figure}

\begin{figure}[ht]
\centering
  \includegraphics[width=1.2\mylength]{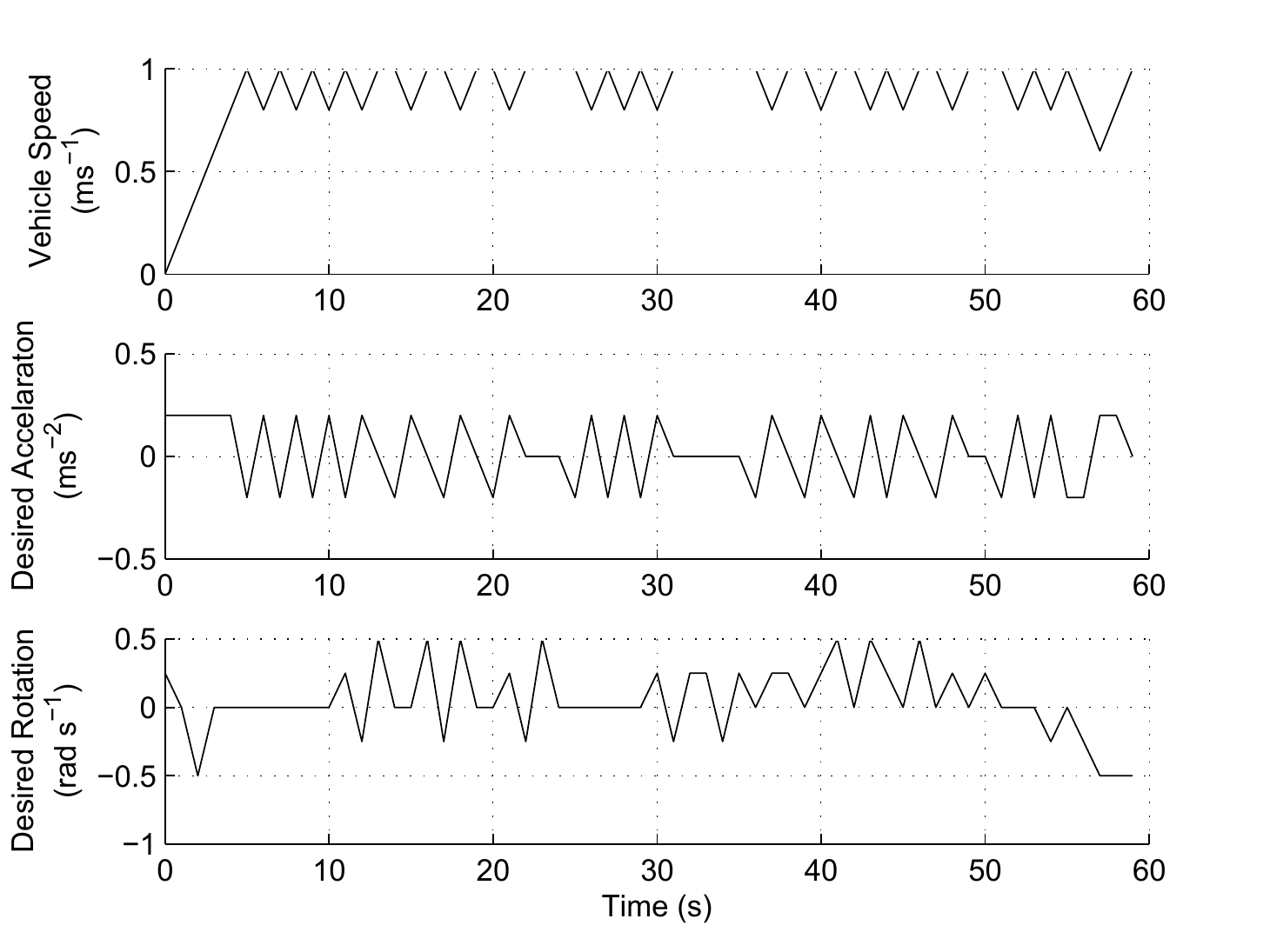}
  \caption{Control time history for PCL.}
\label{fig:bfourth}
\end{figure}

\begin{figure}[ht]
\centering
  \includegraphics[width=1.2\mylength]{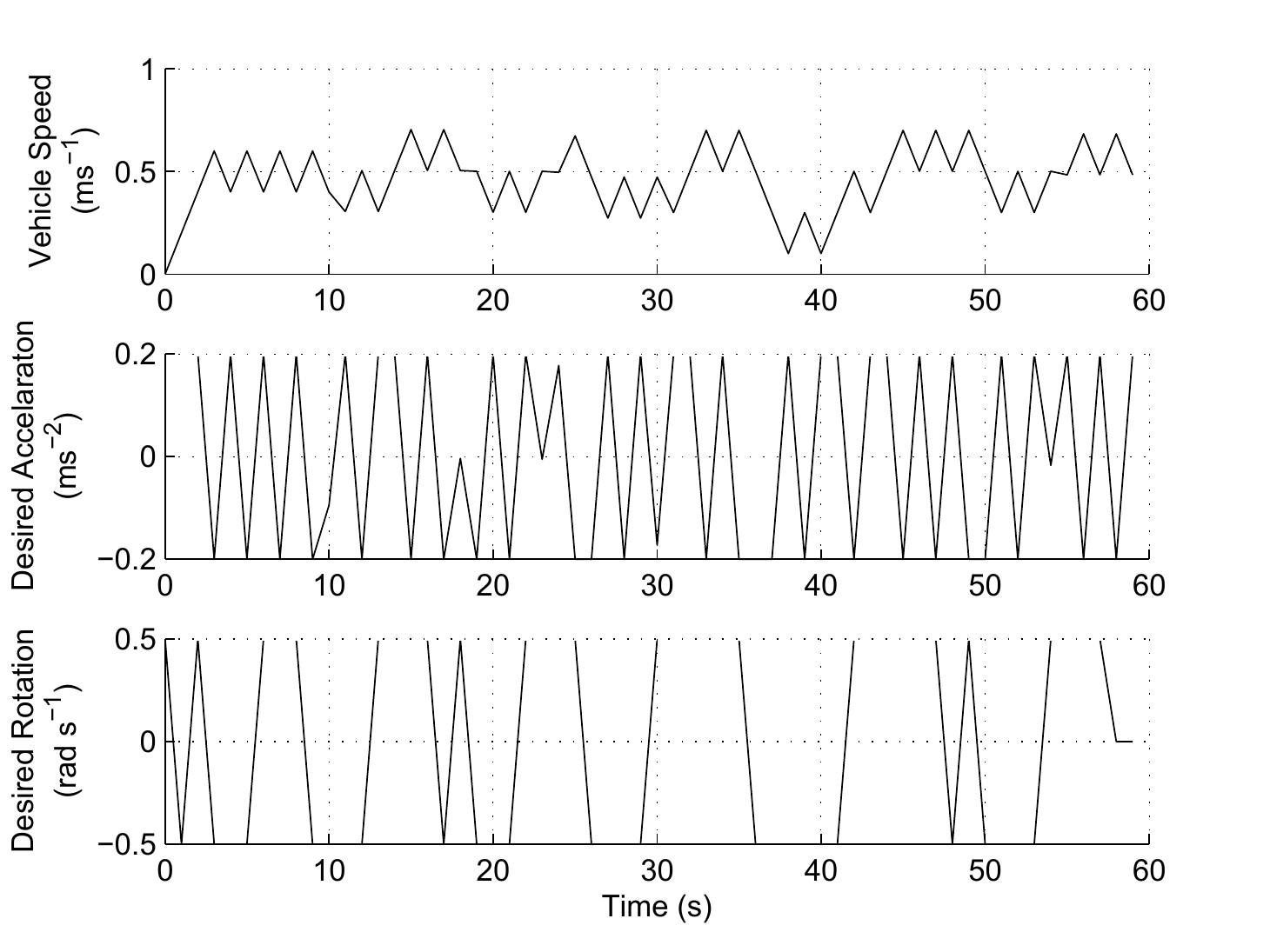}
  \caption{Control time history for ACL.}
\label{fig:bfgeth}
\end{figure}

The results can be seen in Fig.~\ref{fig:move1}, with the control time histories shown
in Figs.~\ref{fig:bfourth} and \ref{fig:bfgeth}. 
It can be seen the ACL is slower, and is affected by oscillation. The
speed is set by a fixed computation and is not directly adjustable. The
oscillation is presumably caused by the extra kinematics, which are
not allowed for by the ACL.

\clearpage \section{Experiments}
\label{ch3:exp}

\begin{figure}[ht]
\centering
\includegraphics[width=0.6\mylength]{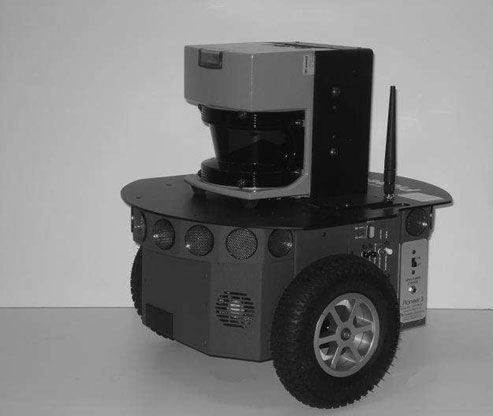}
\caption{Pioneer P3-DX mobile vehicle used for testing.}
\label{fig:photo}
\end{figure}

Experiments were performed with a Pioneer P3-DX mobile robot to show
real-time applicability of this system (see Fig.~\ref{fig:photo}). A
SICK LMS-200 laser range-finding device was used to detect obstacles
in a vicinity of the vehicle. This device has a nominal accuracy of
$15mm$ along each detection ray (however this may be significantly
degraded in real world circumstances \cite{Sanz2011journ7}). Each
measurement was used directly as a separate detection ray in the
sensor model described by Assumption~\ref{as:range}, after being
cropped to be less than $R_{max}$. The standard on-board PC for a P3-DX
was used, which is equipped with a 1.6 GHz Pentium processor.
\par
The robot is equipped with the ARIA library (version 2.7.4) which
provides commands to control the motor drives.  At each control
update, the desired translational acceleration and turning rate of the
ARIA interface were set to match the outputs generated by the
navigation algorithm.  Note the recovery scheme was not implemented as it was
never requested by the navigation algorithm (this may require
additional feedback control to stabilize the position of the vehicle).
The navigation calculation and robot control update was carried out at
$0.1 s$ sampling periods, however the sampling period for generating
trajectories was increased to $1.0 s$ to reduce computational load
(this does not affect the properties of the navigation law). The
initial $C_{end}$ was taken to be the nearest detected point. The
values of the parameters used in the experiments are listed in
Table~\ref{paramexp}. The physical parameters of the vehicles result
from conservative estimates, whereas exact identification was not
carried out. The obstacle was assumed to be stationary.

\begin{table}[ht]
\label{tabular}
\centering
\begin{tabular}{| l | c |}
\hline
 $u_{\theta,max}$ & $0.8 rads^{-1}$ \\
\hline
 $u_{v,nom}$ & $0.1 ms^{-2}$ \\
\hline
 $v_{max}$ & $0.4 ms^{-1}$ \\
\hline
\end{tabular} \hspace{10pt} \begin{tabular}{| l | c |}
\hline
$d_{ob}$ & $0.15 m$  \\
\hline
    $d_{rad}$ & $0.3 m$ \\
\hline
    $R_{max}$ & $4 m$  \\
\hline
$N_r$ & $90$ (semicircular) \\ 
\hline
\end{tabular} 
  \caption{Experimental parameters for boundary-following controller.}
\end{table}

Note that while the algorithm specification assumes a full circular
field of view, the sensor used can only sense the half plane in front
of the vehicle. Through not covered by the analysis, it seems logical
to continue to enforce the same constraints. Conveniently, the
constraint that requires the visibility of the obstacle to be
maintained may compensate for the deficiency; however, because of
disturbance, there is a possibility (however unlikely) for the vehicle
to lose sight of the obstacle while navigating certain types of
corner. However, an experimental design that produced these conditions
was unable to be found.
\par

\begin{figure}[ht]
\centering
  \subfigure[]{\scalebox{0.45}{\includegraphics[width=\columnwidth]{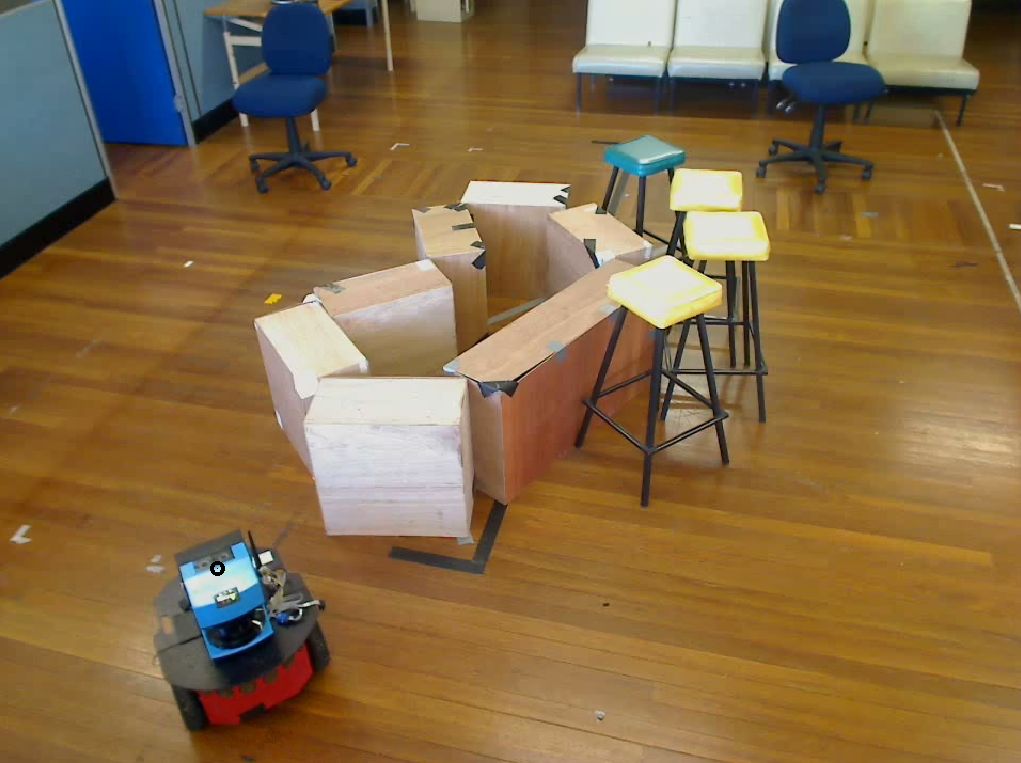}}}
  \subfigure[]{\scalebox{0.45}{\includegraphics[width=\columnwidth]{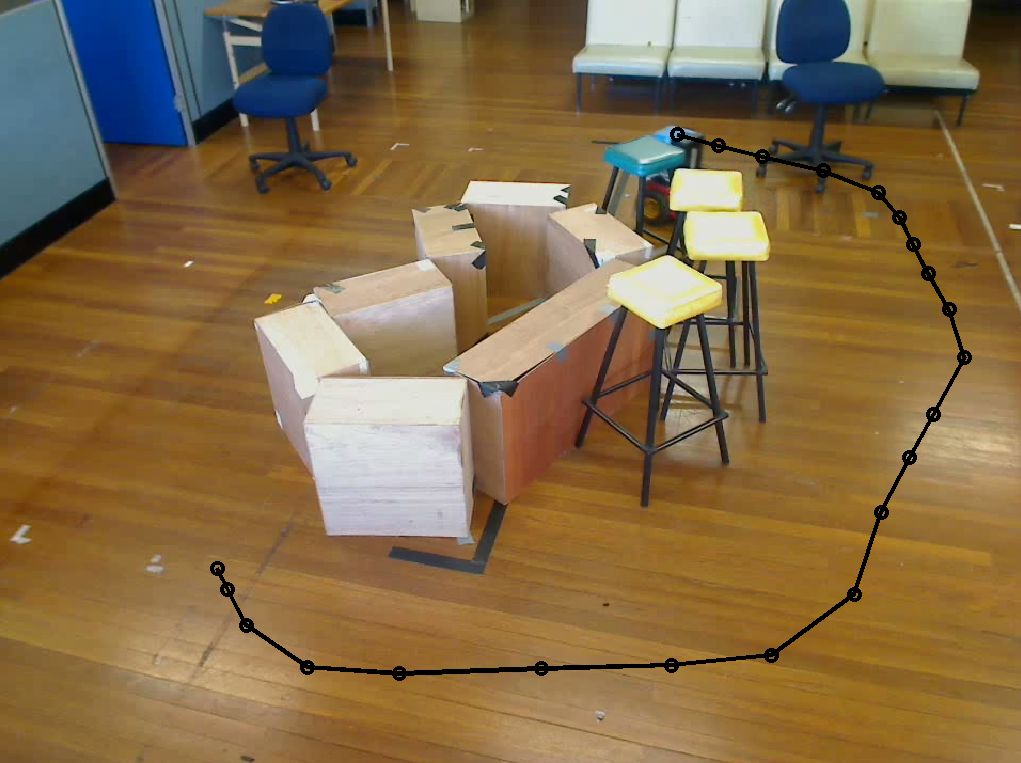}}}
  \subfigure[]{\scalebox{0.45}{\includegraphics[width=\columnwidth]{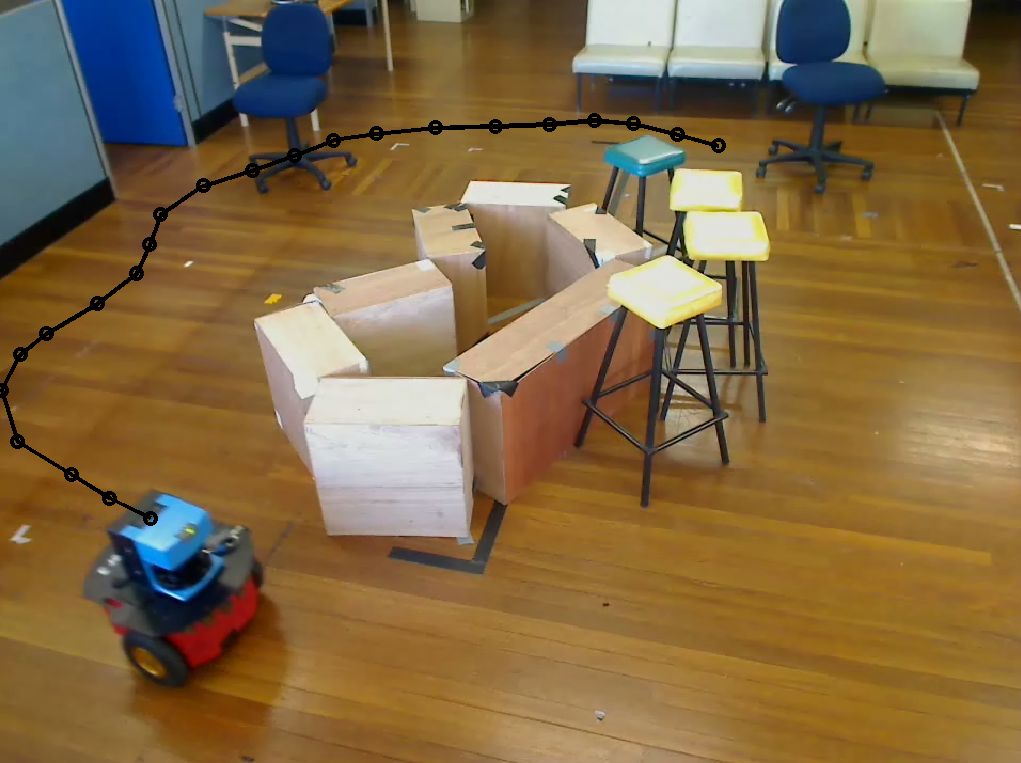}}}
  \subfigure[]{\scalebox{0.45}{\includegraphics[width=\columnwidth]{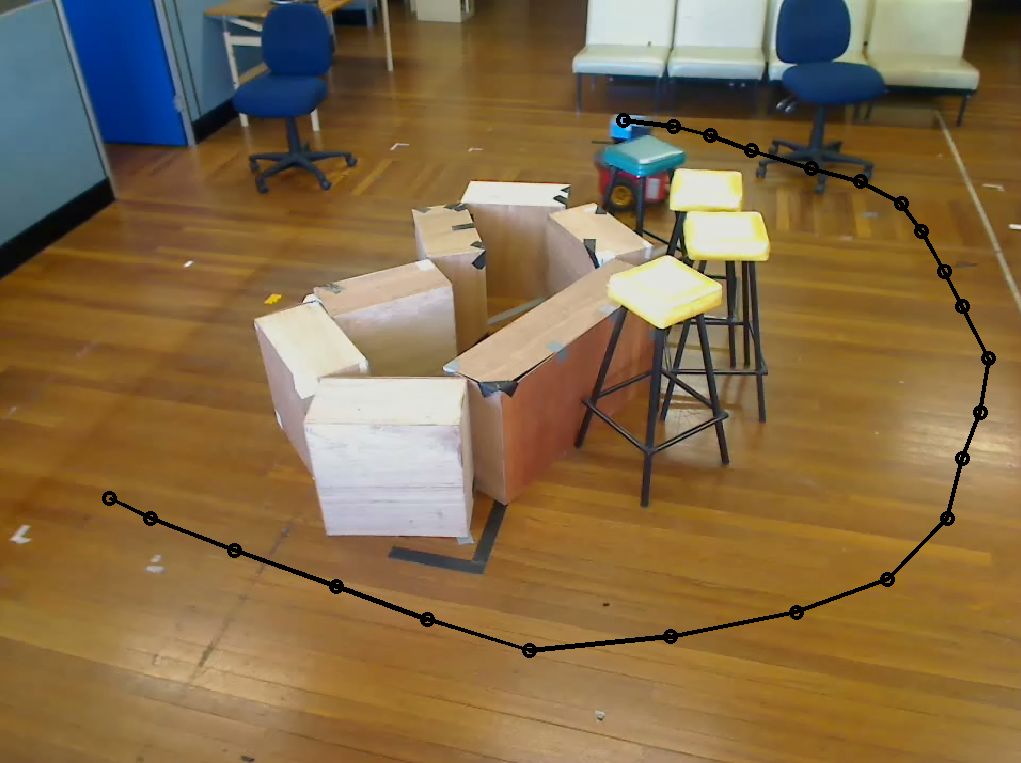}}}
  \subfigure[]{\scalebox{0.45}{\includegraphics[width=\columnwidth]{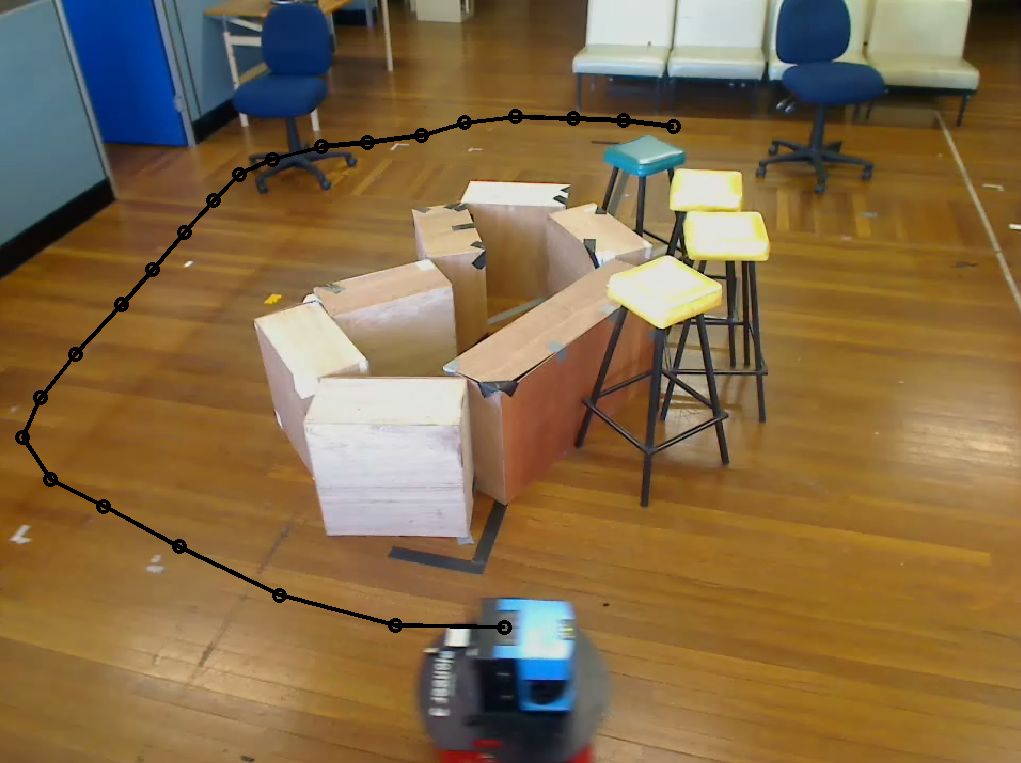}}}

  \caption{Sequence of images showing the experiment.}
\label{fig:exp}
\end{figure}

In Fig.~\ref{fig:exp}, the obtained closed loop trajectory recorded
using a camera is presented, where the trails are plotted by manually
marking the video sequence. It can be seen that the vehicle performed
as expected, making reasonably consistent loops along the boundary.

\begin{figure}[ht]
  \centering
  \includegraphics[width=8cm]{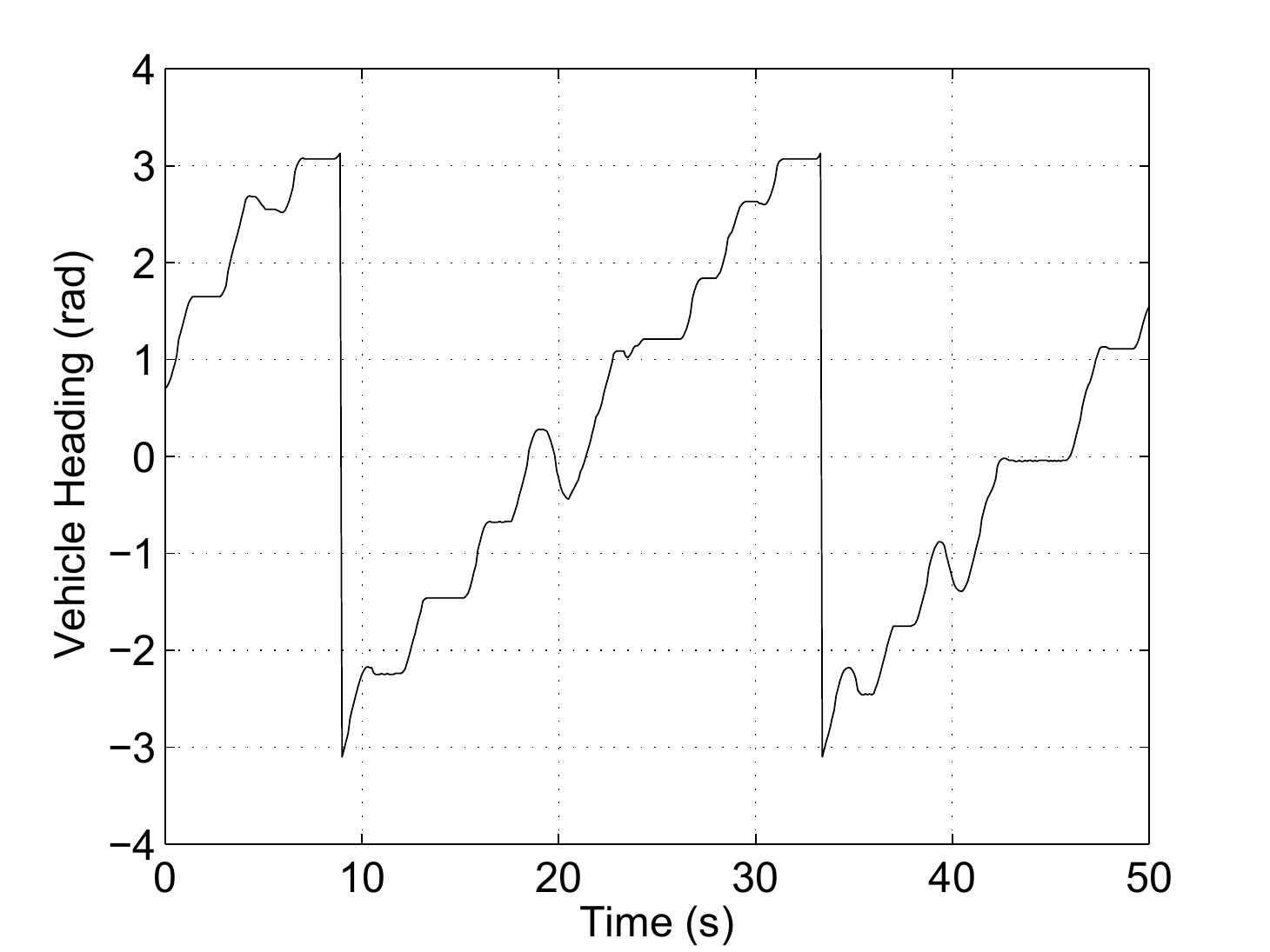}
  \caption{Evolution of the heading of the robot over the
    experiment.}
  \label{fig:expdir}
\end{figure}

\begin{figure}[ht]
  \centering
  \includegraphics[width=8cm]{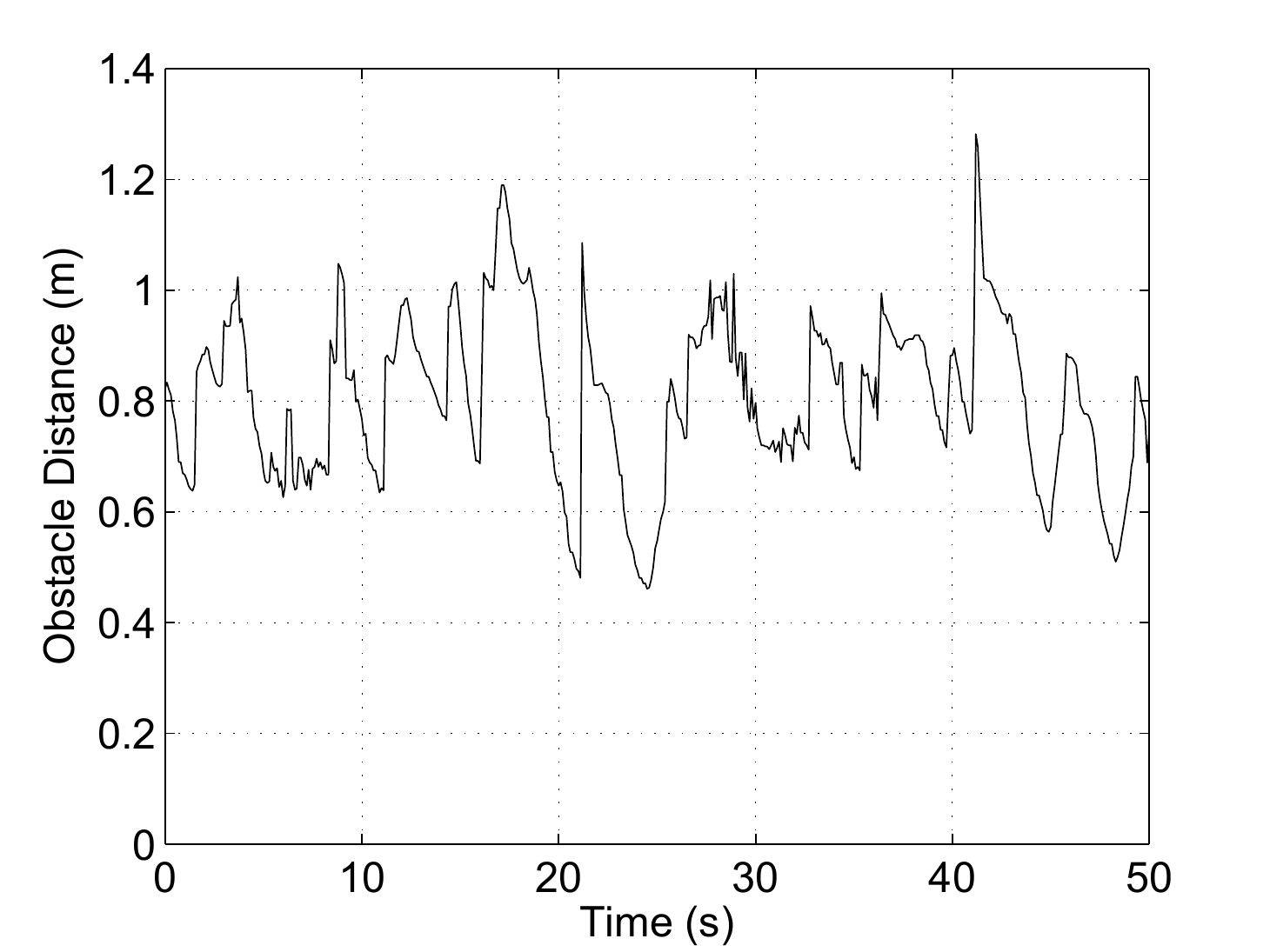}
  \caption{Minimum distance measured by the LiDAR sensor over the
    course of the experiment.}
  \label{fig:expdis}
\end{figure}

\begin{figure}[ht]
  \centering
  \includegraphics[width=8cm]{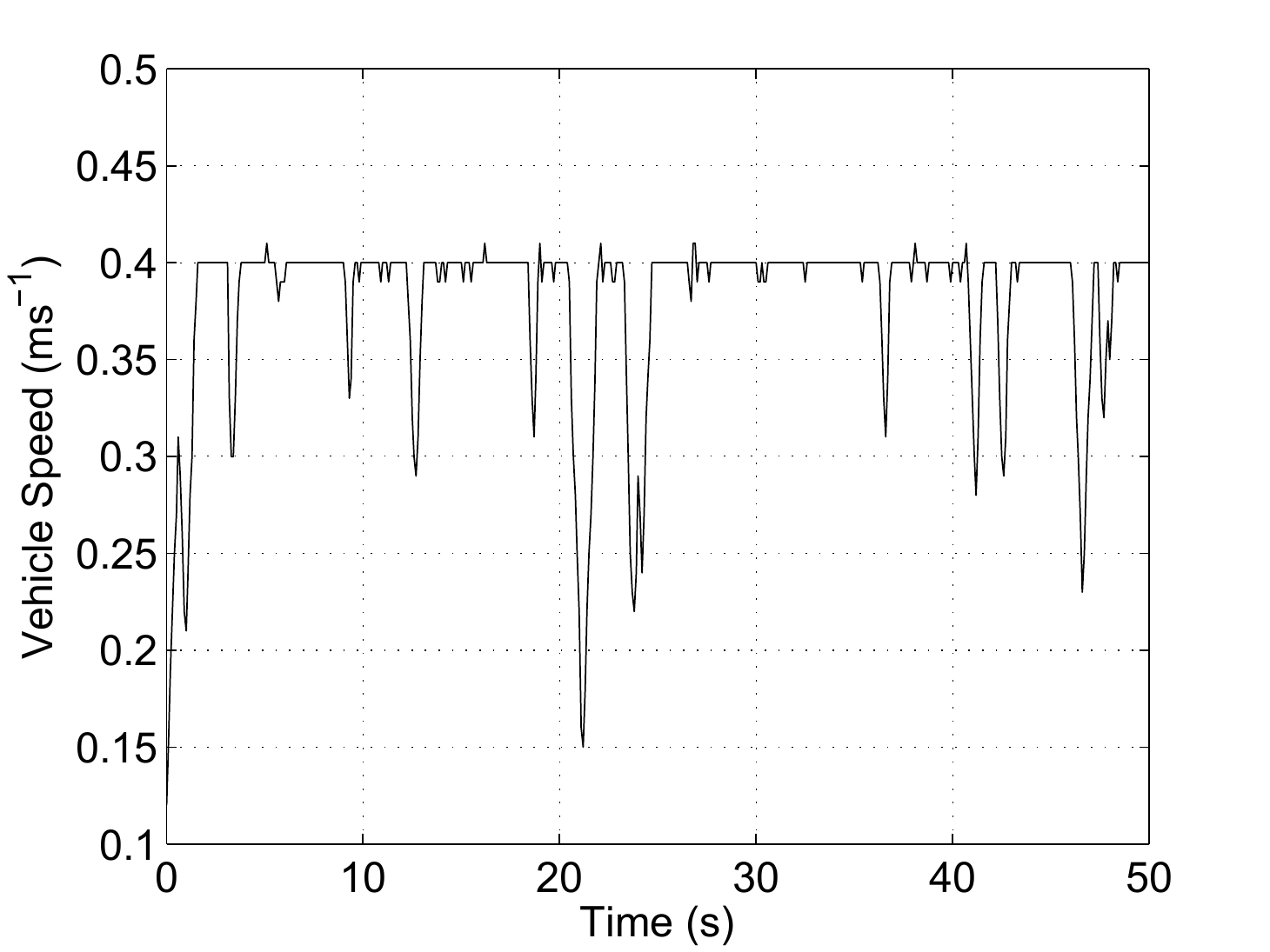}
  \caption{Speed of the robot over the course of the experiment.}
  \label{fig:expvel}
\end{figure}

\begin{figure}[ht]
  \centering
  \includegraphics[width=8cm]{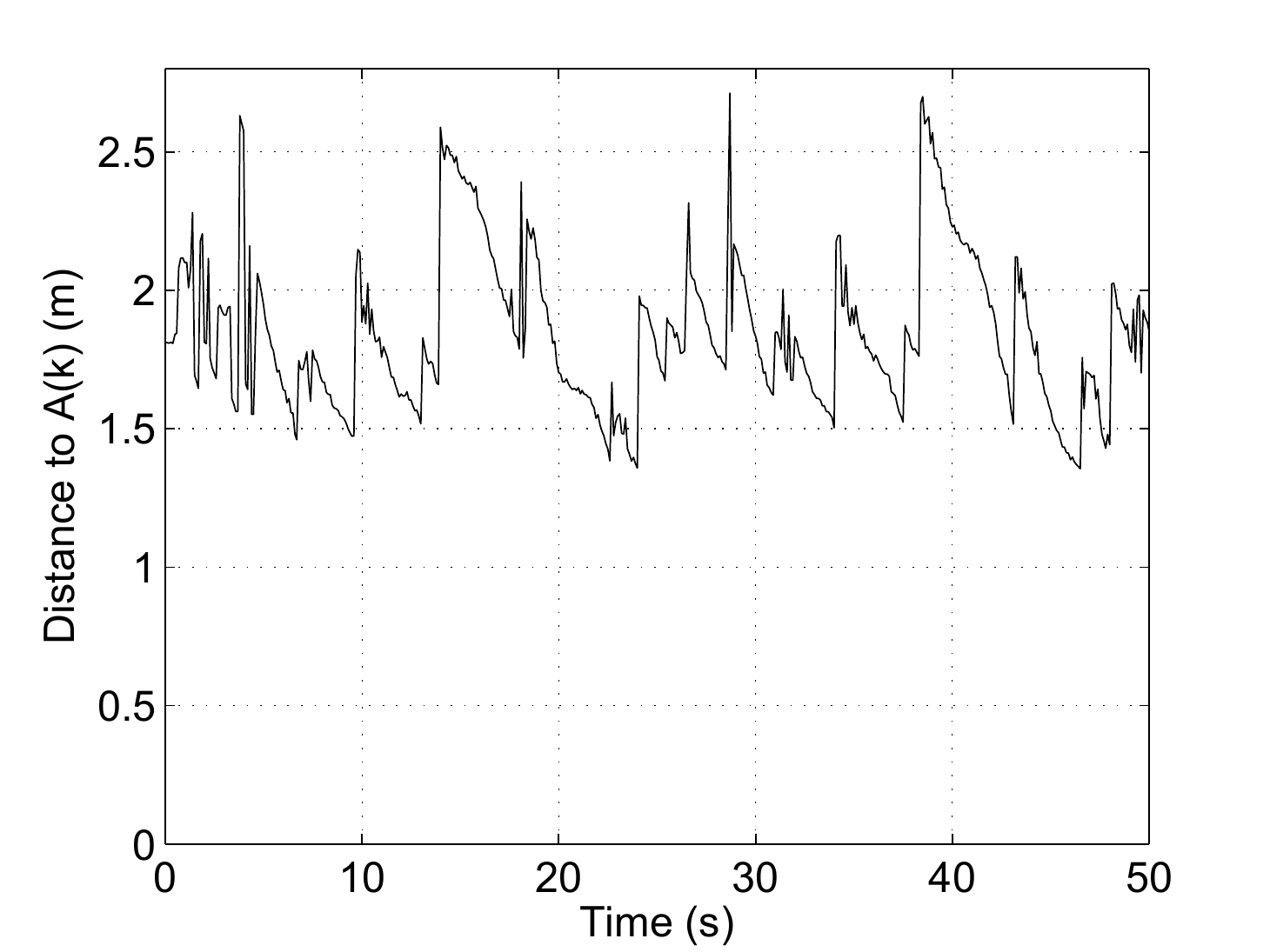}
  \caption{Distance to the target point $A(k)$ over the course of
    the experiment.}
  \label{fig:exptar}
\end{figure}

In Fig.~\ref{fig:expdir} the estimated global heading recorded by
odometry is shown, and it may be seen that the vehicle performs just
over two complete loops around the obstacle. In Fig.~\ref{fig:expdis}
the minimum distance to the obstacle was always greater than 0.43 $m$,
which is congruent with the chosen $d_{tar}$ (0.45 $m$). In
Fig.~\ref{fig:expvel} the vehicle velocity mostly stayed near the
maximum (nominal) velocity (0.4 $ms^{-1}$), except for some brief reductions in
speed. These reductions in speed are most likely caused by the maximum (nominal)
speed being too fast to successfully navigate certain segments of the obstacle.
In Fig.~\ref{fig:exptar} the distance to the target point was
always well over 1.3 $m$, which is expected.

Future work may look at performing experiments with target convergence
and with moving obstacles; however these tests indicate proof of
concept for the basic control law. 

\clearpage \section{Summary}
\label{ch3:con}
In this chapter a method for navigating a vehicle along an obstacle boundary was proposed when local
information from a range finding device is available. The MPC-type approach proposed here is able to
plan a trajectory towards the edge of the points known to be a solid obstacle boundary. This allows
the obstacle offset and vehicle speed to change appropriately to the obstacle. It was shown that
every point on a finite obstacle boundary will have a nearby detection within a finite time.
Possible extensions to target convergence and moving obstacles were presented. Computer simulations
and real world testing confirm the methods validity.

\chapter{Collision Avoidance with Multiple Vehicles}
\label{chap:multiple}

In this chapter, the problem of decentralized collision avoidance of
multiple vehicles operating in a common workspace is considered. Briefly, this is
accomplished by introducing constraints to the trajectory selection process which are
able to maintain non-interfering trajectories despite communication delay. This
is offered as an alternative to the methods of multiplexed MPC and coherency
objectives, which were discussed in Chapt.~\ref{chap:lit}.
\par
Here, a MPC type method is proposed that requires only a single communication
exchange per control update and addresses the issue of communication delay. The
proposed approach does not employ imposing an artificial and auxiliary coherence
objective, and may be suited to real-time implementations, while retaining
robustness properties. 
\par
The body of this chapter is organized as follows. In Sec.~\ref{ch4:ps}, the problem statement is
explicitly defined and the vehicle model is given; in Sec.~\ref{csa}, the structure of
the navigation system is presented. Simulation results with a perfect unicycle model are in Sec.~\ref{ch4:sim}; Experiments 
are in Sec.~\ref{ch4:exp}. Finally, Sec.~\ref{ch4:con} offers brief conclusions.

\section{Problem Statement}
\label{ch4:ps}

In this chapter, $N_h$ autonomous vehicles traveling in a plane are considered, each of which is
associated with a steady point target $T_i, i \in [1:N_h]$. As before, the plane
contains a set of unknown, untransversable, static, and closed obstacles $D_j
\not\ni T_i, j \in [1:n]$. As in Chapt.~\ref{chap:singlevehicle}, the objective is to design a navigation law that
drives every vehicle towards the assigned target through the obstacle-free part of
the plane $F := \mathbb{R}^2 \setminus D$, where $D:= D_1\cup \ldots \cup D_n$.
As before the distance from the vehicle to every obstacle and other vehicles should
constantly exceed the given safety margin $d_{sfe}$. Thus it follows the vehicle's physical radius would be 
less than $\frac{1}{2}d_{sfe}$.
\par
In this chapter, only the unicycle model is considered for brevity, through these
results could easily be extended to the holonomic model.
As in Chapt.~\ref{chap:singlevehicle}, it is assumed that the vehicle has knowledge of an
arbitrary subset of the obstacle free part of the plane, $F_{vis, i}(k) \subset F$;
the vehicle has knowledge of its current state $\state$; and the vehicle has access
to the position of its target $T_i$.

\begin{figure}[ht]
	\centering
		\includegraphics[width=0.75\mylength]{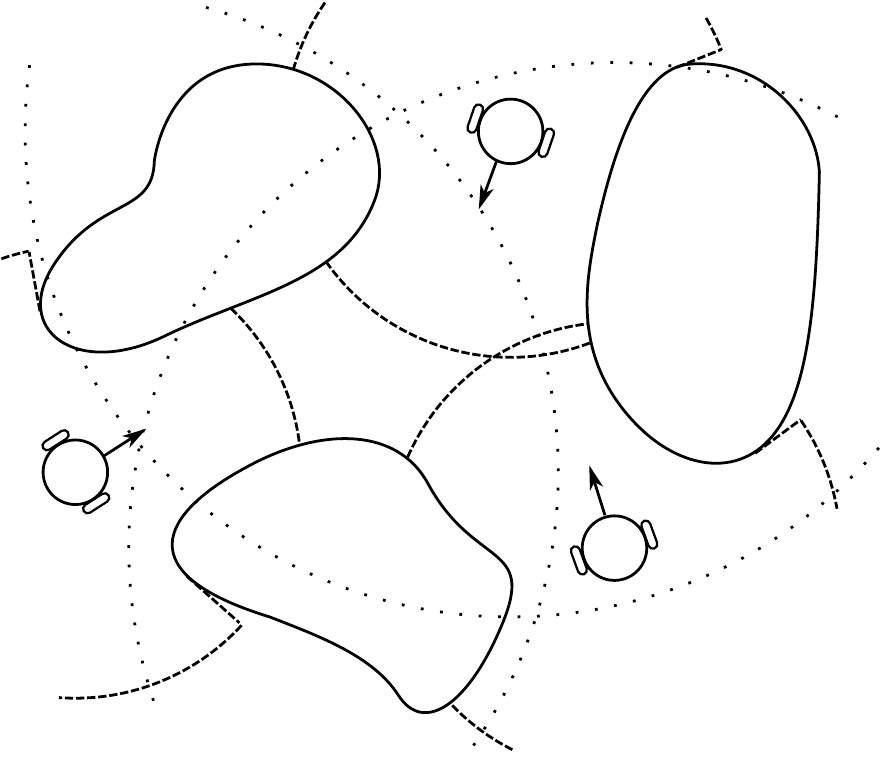}
	\caption{Scenario containing several vehicles and static obstacles.}
	\label{fig:scen}
\end{figure}

It is assumed the vehicles have capability of communication with nearby companions --
every vehicle is able to broadcast to the other vehicles within a given radius $C$.
Any communication delay is assumed to be less than the time step, so that any
data transmitted at time step $k$ from one vehicle is available at time $k+1$ at
all vehicles within the communication radius.
 \par
 The vehicles need not be identical -- the model parameters
and the characteristics of the sensed area $F_{vis, i}(k)$ may depend on the vehicle index $i$.
 However it is assumed that the communication radius $C$ is the same for all
 vehicles (in other words, the communication graph is undirected). Furthermore,
 it is assumed that each vehicle has knowledge of all the fixed parameters of other
 vehicles inside the communication radius, by broadcast or otherwise.
 \par
 It is required that the sampling times of each vehicle are synchronized, however
communication based algorithms are available for the task of synchronization
\cite{Giridhar2006conf4, Sun2006journ, Zhou2007journ}.

\begin{Remark} \rm
There is a limit to the number of agents that lie within the communication range of a particular
vehicle, owing to the constraint on the minimum distance between the vehicles. This entails a bound
on the maximum communication burden on each vehicle caused by scaling the number of vehicles.
\end{Remark}

Finally, as in Chapt.~\ref{chap:singlevehicle}, it is assumed that the vehicles are
subjected to external disturbances. 
To highlight the major points by dropping secondary and standard technical
details, perfect measurement of the vehicle state is assumed. At the same time,
bounded state estimation and obstacle measurement errors can be taken into account in a
similar way, which point of view transforms estimation error into
disturbance.

 \section{Navigation System Architecture}
\label{csa}

As  in Chapt.~\ref{chap:singlevehicle}, navigation is based on generating probational trajectories over the relatively short planning horizon
$[k,(k+\tau)]$ at every time step $k$. Probational trajectory is nominal and is specified by a
finite sequence of way-points, which satisfy the constraints of the nominal model. 
As before, probational trajectories necessarily halt at the terminal planning time. When necessary, the trajectory can be
being prolonged arbitrarily by the `stay still' maneuver.
\par
At any sampling time $k$, every vehicle broadcasts its current state $\state$
 and probational trajectory, which data arrive at the vehicles in the communication range
one time step later. To avoid collisions with the companions, the vehicle also attempts to
reconstruct their current probational trajectories, starting with reconstruction of the planned
ones. However, the latter start at the current states of the companions, which are unknown to the
vehicle at hand. So the vehicle first estimates these states on the basis of the received states 
at the previous time step (the estimate may differ from the
actual state because of the disturbances). After this the vehicle copycats generation of the planned
trajectories for every vehicle in the communication range, with substituting the estimate in place
of the true state in doing so. This gives rise to the set of {\it presumable planned trajectories}.

\subsection{Overview}

The \textit{Trajectory Tracking Module} (TTM) is identical to the presentation in Chapt.~\ref{chap:singlevehicle}. 
However, the \textit{Trajectory Planning Module} (TPM) is similar to the process employed 
in Chapt.~\ref{chap:singlevehicle}. Every vehicle iteratively executes the following steps; in this enumeration, only Steps~S.2 and S.5 are different:

\begin{enumerate}[{\bf S.1}]
\item \label{ch4:step1} Generation of a finite set $\mathscr{P}$ of planned trajectories, each starting
at the current vehicle state $\state$.

\item  {\bf Refinement of $\boldsymbol{\mathscr{P}}$ to only \textit{mutually feasible} trajectories (see Sec.~\ref{sec:refine})}.

\item Selection of a trajectory from $\mathscr{P}$:

\begin{itemize}
  \item If $\mathscr{P}$ is empty, the probational trajectory is \textit{inherited} from the previous time step
(with a proper time shift). 
 \item Otherwise, the probational trajectory is \textit{updated} by choosing an element of $\mathscr{P}$
   minimizing some cost function (see Sec.~\ref{sec.sel}).
\end{itemize}

\item {\bf Transmission of the chosen probational trajectory $\blds^*(j|k)$ and the current state $\boldsymbol{\state}$
to all vehicles within the communication radius $\boldsymbol{C}$}.

\item $k:=k+1$ and go to { S.\ref{ch4:step1}}.
\end{enumerate}

\begin{figure}[ht]
	\centering
		\includegraphics[width=0.75\mylength]{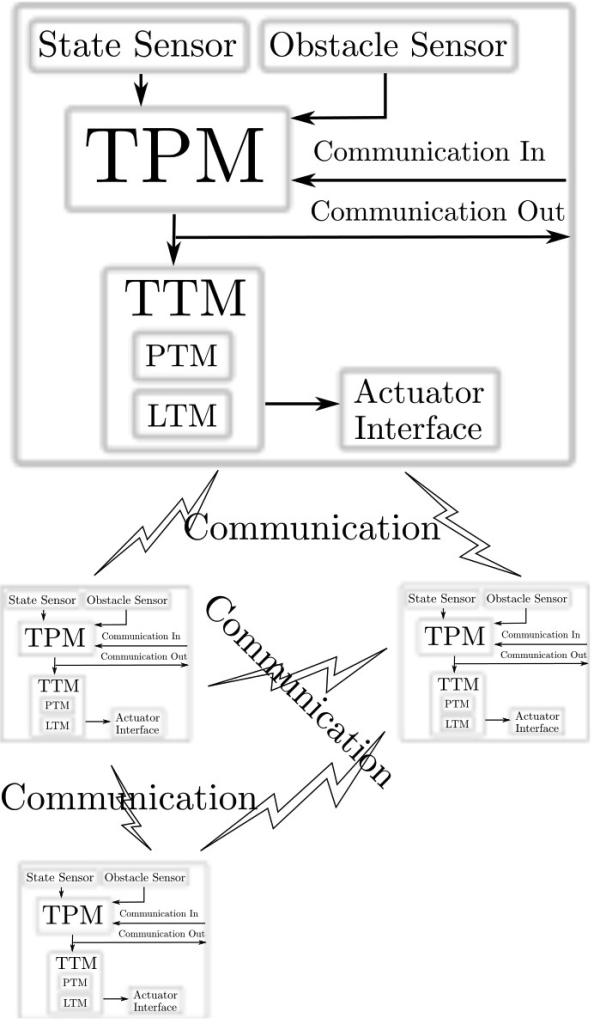}
	\caption{The overall architecture of the navigation system.}
	\label{fig:ovearll}
\end{figure}

This navigation process is executed co-currently in each vehicle. The overall architecture of 
the navigation system is illustrated in Fig.~\ref{fig:ovearll}. Step~S.2 is solely based on the
currently observed part of the environment and data received from and sent to other vehicles in
communication range since time step $k-1$. This step involves construction of the set of presumable
planned trajectories for all vehicles in the communication range.

\par
The upper index $^*$ is used to mark variables associated with the planned trajectories; this
index is discarded to emphasize that the concerned trajectory is in fact probational; the hat
$\widehat{\{\}}\vphantom{\}}^{\ast}$ is added to signal that a presumable planned trajectory is
concerned. These variables depend on two arguments $(j|k)$, where $k$ is the time instant when the
trajectory is generated, and $j \geq 0$ is the number of time steps into the future: the related
value concerns the state at time $k+j$. Whenever a certain vehicle is considered, the lower index
$\{\}_i$ refers to other vehicles within the communication range $C$, whereas its absence indicates
reference to the vehicle at hand.

\subsection{Safety Margins}
\label{ch4.ansm}
To ensure collision avoidance, as before TPM respects more conservative safety margins than $d_{sfe}$.
They take into account deviations from the probational trajectory caused by disturbances and are
based on estimation of these mismatches. 
\par
Recall $d_{trk}$ is the maximum deviation between the vehicles actual state and the nominal
state predicted by the probational trajectory. To simplify notations, it is assumed that $d_{trk}$ is the same for all
vehicles; however the navigation scheme is still feasible if this quantity varies between them.
\par 
The following more conservative safety margins account for not only disturbances but also the
use of time sampling in measurements of relative distances:

\begin{equation}
\label{ch4.dist}
d_{tar}:= d_{sfe}+ \frac{v_{max}}{2} + d_{trk}, \qquad d_{mut} :=  d_{sfe} + 2 \left[ \frac{v_{max}}{2} +  d_{trk} \right].
\end{equation}

To illuminate their role, the following definition is introduced:

\begin{Definition}
\label{ch4:dmut}
An ensemble of probational trajectories generated by all vehicles at a given time is said to be \textbf{
mutually feasible} if the following claims hold:
\begin{enumerate}[i)]
\item The trajectory is feasible, i.e. the distance from any way-point to any static obstacle exceeds $d_{tar}$ for any trajectory;
\item For any two trajectories, the distance between any two matching way-points exceeds $d_{mut}$.
\end{enumerate}
\end{Definition}

The motivation behind this definition is illuminated by the following:

\begin{Lemma}
\label{ch4:lem:feas0}
Let an ensemble of probational trajectories adopted for use at time step $k-1$ be mutually feasible.
Then despite of the external disturbances, the vehicles do not collide with the static obstacles and
each other  on time interval $[k-1,k]$ and, moreover, respect the required safety margin
$d_{sfe}$.
\end{Lemma}
 
\pf The static obstacle case follows from Lemma~\ref{lem:feas0}. For any pair of vehicles, the distance between the
concerned way-points exceeds $d_{mut}$. So the  the distance between the real positions at time $t$
is no less than $d_{mut} -  2[\frac{v_{max}}{2} + d_{trk}]  = d_{sfe}$. \epf

\subsection{Trajectory Constraints}
\label{sec:refine}

Refinement is prefaced by the following computations, which are carried out
 at the vehicle at hand for
every vehicle $i$ that was within the communication range at the previous time step $k-1$:

\begin{itemize}
\item The nominal motion equations are integrated for vehicle $i$ from $k-1$ to $k$ with
the initial data $\state_i(k-1)$ and controls $\contr^\ast(0|k-1)$ to acquire the estimated current state $\widehat{\state}_{i}(k)$;\footnote{Recall that these data and controls were
broad-casted by the $i$th vehicle at time $k-1$. So they are received by the vehicle at hand at time
$k$.} 
\item The procedure from Sec.~\ref{Sec.plan} is carried out for the $i$th vehicle,
with the true current state $\state_i(k)$ replaced by the estimated one
$\widehat{\state_i}(k)$.
\end{itemize}

\begin{figure}[ht]
	\centering
		\scalebox{0.2}{\includegraphics{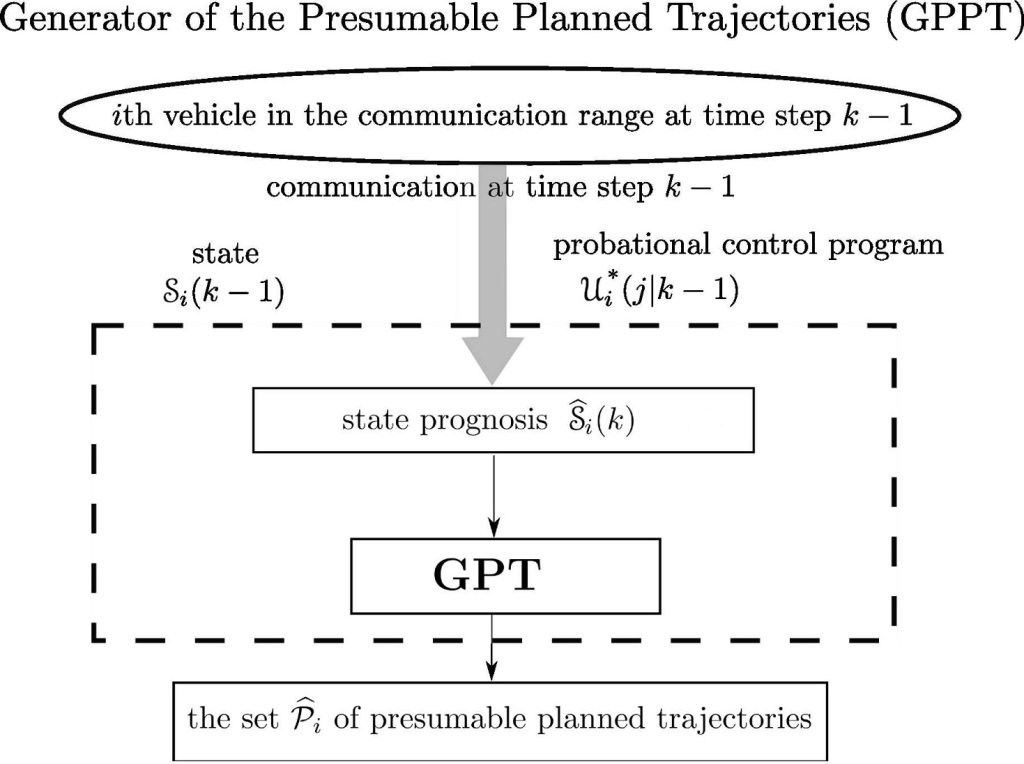}}
	\caption{Generator of the planned trajectories.}
	\label{fig:gppt}
\end{figure}

This procedure is summarized in Fig.~\ref{fig:gppt} and results in the set
$\widehat{\mathscr{P}}_i$ of presumable planned trajectories of the $i$th vehicle.

To proceed, the following estimates are introduced:

\begin{itemize}

\item $\Delta_{osd}$ -- an upper bound on the deviation over the single time step between the
nominal and real state $\state$ of the vehicle driven by the TTM along the probational trajectory in
the face of disturbances;
 \item  $d_{\tau}$ -- an upper bound on the translational
deviation between two nominal trajectories over the entire planning horizon of duration $\tau$,
provided that the difference between their initial states does not exceed $\Delta_{osd}$ and the both 
follow common speed and turning patterns selected by the trajectory generation rule from
Sec.~\ref{Sec.plan}.
\end{itemize}

\begin{figure}[ht]
	\centering
		\includegraphics[width=0.8\mylength]{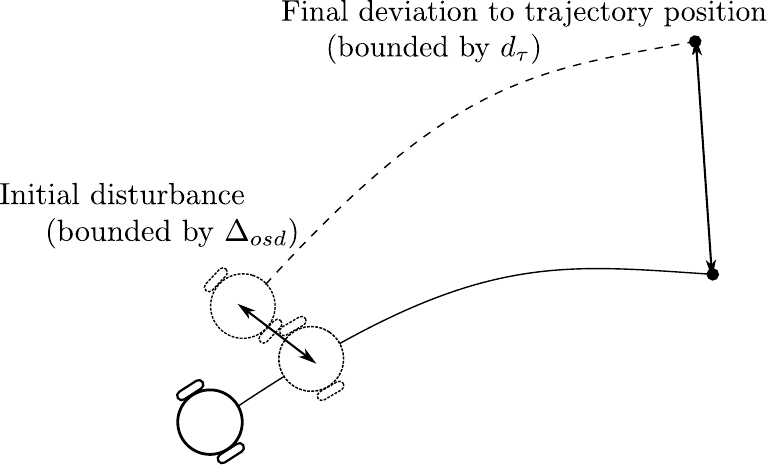}
	\caption{Effect of disturbance over a single time-step.}
	\label{fig:dist}
\end{figure}

 The constant $\Delta_{osd}$ is similar to $d_{trk}$; its computation is discussed
 in Sec.~\ref{sec:ptm}. $\Delta_{osd}$ and $d_{\tau}$ are illustrated in Fig.~\ref{fig:dist}.
 
 Given $\Delta_{osd}$, the bound $ d_{\tau}$ can be easily computed based on the vehicle model.
  More or less conservative bounds can be used, up to the tight bounds based on
 computer-aided calculation of the reachable sets. However the experimental results show that even
 conservative estimates are able to entail good performance. Partly in view of this and partly due
 to the fact that computation of conservative bounds is elementary, the related
 details are omitted.
\par
Note that Remark~\ref{rem.van} is evidently extended on $\Delta_{osd}$ and $ d_{\tau}$. The
following remark is immediate from the foregoing and partly elucidates the role of $ d_{\tau}$:

\begin{Remark} \rm
\label{rem.mismatch}
The mismatch $\max_{j=0,\ldots,\tau}$ $\norm{\blds_i(j|k)$ $-$ $\widehat{\blds}\vphantom{\blds}_i(j|k)}$ between the
presumable planned trajectory and the true planned trajectory does not exceed $ d_{\tau}$.
\end{Remark}

In order to ensure that the probational trajectories are mutually feasible, the set of planned
trajectories is subjected to a series of refinements.
First, the vehicle at hand refines the generated set of presumable planned trajectories
$\{\widehat{\blds}\vphantom{\blds}_i^{\;\ast}(j|k) ,\ldots\}_j$ of any other vehicle $i$ that was in the
communication range at the previous time step $k-1$. Specifically, it discounts any presumable
trajectory that interferes with the true probational trajectory of the vehicle at hand
$\{\blds(j|k-1),\ldots\}_j$ selected at the previous time step $k-1$, i.e., such that:

\begin{equation}
\norm{\widehat{\blds}\vphantom{\blds}_i^{\;\ast}(j|k) - \blds(j+1|k-1)} \leq d_{mut} -  d_{\tau} \quad \text{for some} \quad j.
\end{equation}

\begin{figure}[ht]
	\centering
		\scalebox{0.35}{\includegraphics{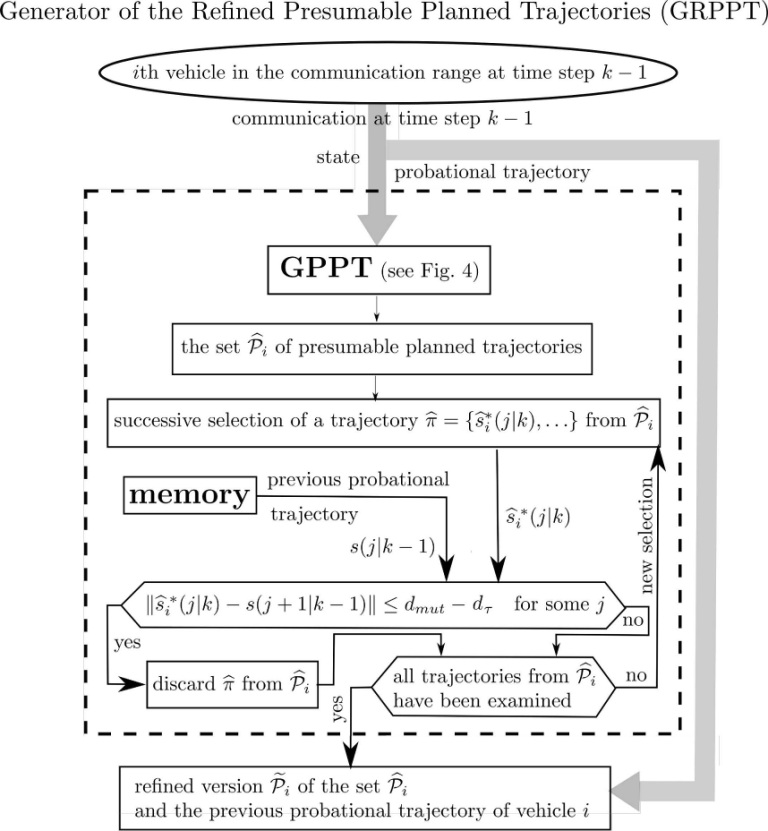}}
	\caption{Generator of the presumable planned trajectories.}
	\label{fig:grppt}
\end{figure}

This operation is summarized in Fig.~\ref{fig:grppt} and gives rise to the set of the refined
presumable planned trajectories $\widetilde{\mathscr{P}}_i$ of the $i$th vehicle. 

After this, the vehicle at hand refines its own set of planned trajectories by successively
discounting those interfering with:

\begin{enumerate}[{\bf r.1)}]
\item Static obstacles, such that:
\begin{equation}
 \norm{\blds^*(j|k) - x } < d_{tar} \quad \text{for some} \quad  x \in D ,\quad \text{ and} \quad j
 \end{equation}
\item { Probational trajectories of other vehicles received
 from them within the time interval $[k-1,k]$}, such that:
\begin{equation}
\label{kminus}
 \norm{\blds^*(j|k) - \blds_i(j+1|k-1) } < d_{mut} \quad \text{for some}\; j
  \end{equation}
  \item  Any refined presumable planned trajectory of any other vehicle $i$ that was in the
  communication range of the vehicle at hand at the previous time step $k-1$:
\begin{equation}
\norm{\blds^*(j|k) - \widehat{\blds}\vphantom{\blds}_i^{\;\ast}(j|k)} < d_{mut}  +  d_{\tau}  \; \text{for some} \; j \; \text{and} \; \{\widehat{\blds}\vphantom{\blds}_i^{\;\ast}(j|k), \ldots \} \in \widetilde{\mathscr{P}}_i.
\label{excl}
  \end{equation}
\end{enumerate}

As stated previously, if this refinement sweeps away all trajectories, the probational
trajectory from the previous time step can be used. This leads to the following:

\begin{Lemma}
\label{lem:feas}
The proposed refinement procedure does ensure that the probational trajectories generated at time
step $k$ are mutually feasible provided that the probational trajectories were mutually feasible at
the previous time step $k-1$.
\end{Lemma}

\pf Let us interpret vehicles $j$ and $i$ as the `vehicle at hand' and the `other
 vehicle', respectively; according to the adopted notation, this means the
 index $j$ may be dropped everywhere. Also note the probational trajectories were selected from the set
 of planned trajectories, i.e. $\blds(\cdot|k)=\blds^\ast(\cdot|k), \blds_i(\cdot|k)=\blds_i^\ast(\cdot|k)$.

The trajectory $\blds_i^\ast(\cdot|k)$ has passed the test from r.2) at the vehicle
$i$ and so does not interfere with the trajectory transmitted from the vehicle
at hand at time step $k-1$, such that $\norm{\blds_i^{\ast}(r|k) - \blds(r+1|k-1) } >
d_{mut}\; \forall r$. For the related presumable planned trajectory, it follows:

\begin{multline}
\label{udbrc}
\norm{\widehat{\blds}\vphantom{\blds}_i^{\ast}(r|k) - \blds(r+1|k-1) }  \geq \norm{\blds_i^{\ast}(r|k) - \blds(r+1|k-1) } \\ - \underbrace{\norm{ \widehat{\blds}\vphantom{\blds}_i^{\ast}(r|k) - \blds_i^{\ast}(r|k) } }_{\leq  d_{\tau} \;\text{by Remark~\ref{rem.mismatch}}} \geq d_{mut} -  d_{\tau} .
\end{multline}

So this trajectory was among those against which $\blds^\ast(\cdot|k)$ was
tested in correspondence with Eq.\eqref{excl}. Since $\blds^\ast(\cdot|k)$ has
survived the refinement procedure, this test was passed:

\begin{equation*}
\norm{ \blds^*(r|k) - \widehat{\blds}\vphantom{\blds}_i^{\;\ast}(r|k)} > d_{mut} +  d_{\tau} \qquad \forall r.
\end{equation*}

The proof is completed by invoking the under-braced inequality from
Eq.\eqref{udbrc}:

\begin{multline*}
\norm{ \blds(r|k) - \blds_i(r|k)} = \norm{ \blds^\ast(r|k) - \blds_i^\ast(r|k) } \geq \norm{ \blds^\ast(r|k) - \widehat{\blds}\vphantom{\blds}_i^\ast(r|k)}
\\
- \norm{ \blds_i^*(r|k) - \widehat{\blds}\vphantom{\blds}_i^{\;\ast}(r|k)} \geq  d_{mut} +  d_{\tau} -  d_{\tau} = d_{mut}.
\end{multline*} 
This completes the proof. \epf

 By retracing the arguments of the first part of the proof and invoking Remark~\ref{remiitt}, the following statement can be made:
 
\begin{Remark} \rm
\label{remiitt}
It is assumed that initially there are no vehicles in communication range of each
other and every vehicle is far enough from the static obstacles so that at least
one planned trajectory survives r.1). Then the generated set of probational
trajectories is well-defined and mutually feasible, so it also follows that this
set remains mutually feasible in subsequent time steps.
\end{Remark}

The overall performance of the entire navigation system is addressed in the following:

\begin{Proposition}
Let every vehicle be driven by the proposed navigation law. Then the vehicles do not collide with the
static obstacles and each other and, moreover, respect the required safety margin $d_{sfe}$.
\end{Proposition}

\pf By Lemma~\ref{lem:feas} and Remark~\ref{remiitt}, the ensemble of probational trajectories adopted
for use at any step is mutually feasible. Lemma~\ref{ch4:lem:feas0} completes the proof. \epf

\begin{figure}[ht]
	\centering
		\scalebox{0.2}{\includegraphics{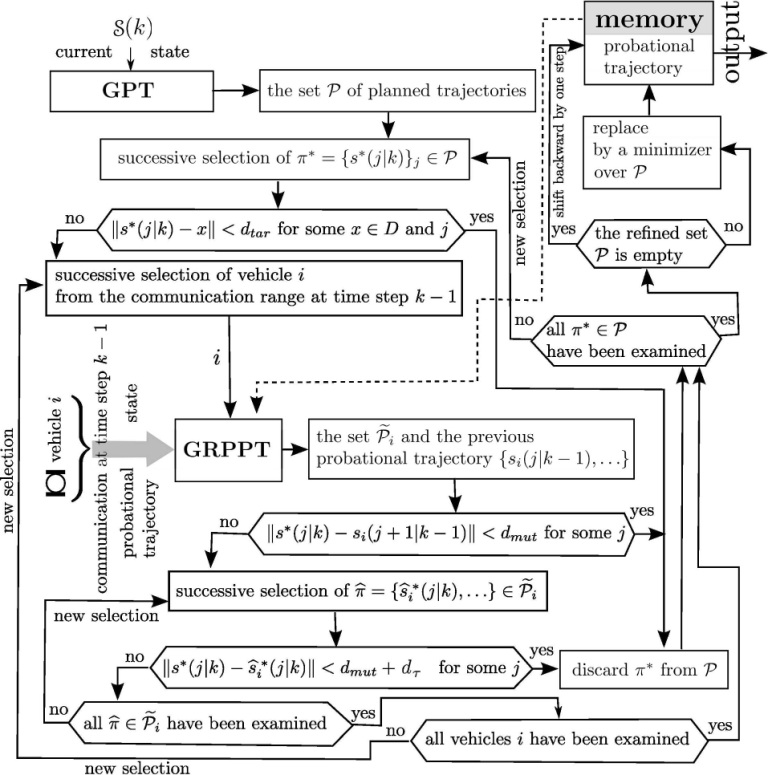}}
	\caption{Trajectory Planning Module.}
	\label{fig:tpm}
\end{figure}

\begin{algorithm}

\caption{Operation of TPM.}
\KwIn{The current state of the vehicle: $\state(k) \equiv \langle \blds(k), \theta(k), v(k)  \rangle$\\
The probational trajectory that is the previous output of the algorithm: $\{ \blds^*(j|k-1), \ldots \}$\\
The trajectories received from other vehicles: $\{ \blds_i^*(j|k-1), \ldots \}$\\
The currently sensed obstacle set $F_{vis,i}(k)$}
\KwOut{Current probational trajectory: $\{ \blds^*(j|k), \ldots \}$}
\textbf{Process:}\\
Generate the set of planned trajectories  $\mathscr{P}$\\
Remove from $\mathscr{P}$ the trajectories that are closer than $d_{tar}$ to the currently sensed obstacle set\\

\ForEach{Received trajectory $\{ \blds_i^*(j|k-1) , \ldots \}$}{
Remove from $\mathscr{P}$ the trajectories that are closer than $d_{mut}$ to $\blds_i^*(j|k-1)$\\
Generate the set of presumable planned trajectories $\widehat{\mathscr{P}}_i$ based on $\state^\ast(1|k-1)$\\
Remove from $\widehat{\mathscr{P}}_i$ the trajectories that are closer than $d_{mut}- d_{\tau}$ to $\blds^*(j|k-1)$\\
Remove from $\mathscr{P}$ the trajectories that are closer than $d_{mut}+ d_{\tau}$ to some trajectory in the remaining set $\widehat{\mathscr{P}}_i$\\
}

\If{$\mathscr{P} = \emptyset$}{
Reuse the trajectory from the previous time-step
}
\Else{
Select a minimizer of Eq.\eqref{minimizer} over $\mathscr{P}$
}
\label{alg1}
\end{algorithm}

This operation is summarized in Fig.~\ref{fig:tpm} and by the pseudo-code Algorithm~\ref{alg1}. Note an upper bound on the total computation load is proportional to:

\begin{equation}
O(N_v, N_t, N_p) = N_v \cdot N_t^2 \cdot N_p, \quad N_t = 4 \cdot \left\lceil\frac{\tau}{\Delta \Lambda}\right\rceil + 2, \quad N_p \equiv  \tau
\end{equation}

where $N_v$ is the number of the vehicles in the communication range and $N_t$ is the number of
planned trajectories.
As for the communication load, only four reals are required to encode the state for communication,
with extra two integers being needed to encode the index of the trajectory inside $\mathscr{P}$ and
the progression along the trajectory.
\par 
Modulo slight extension, the proposed algorithm also displayed good resistance to packet
dropouts, which are common in wireless communications \cite{Wijesinha2005conf0}.
The only trouble caused by a packet dropout is that the vehicle loses access to the state and
probational trajectory of some `other' vehicle at the previous time step and so is unable to 
carry out the proposed refinement procedure; see Sec.~\ref{sec:refine}. This trouble can be
easily overcome by employing the latest known data about the `other' vehicle, which is brought by
the last successful transmission. This however gives more conservative approximations of the true
planned trajectories If this conservatism appears to be too large, trajectory tracking would be engaged.
Thus it can be reasonably inferred that packet loss will reduce the system performance, though can
hardly violate robustness.
\par 
The main advantage of the proposed approach is its `reactiveness': collision avoidance only
requires a single time step of latency for communication between vehicles before mutual action is
taken.

\subsection{Implementation Details}
\label{sec.sel}
From the set of the planned trajectories surviving the refinement procedure, the
final trajectory is selected to furnish the minimum of the cost functional:

    \begin{equation}
    \label{minimizer}
    J = \norm{\blds^*(\tau|k) - T} - \gamma_0 \cdot v^*(1|k),
    \end{equation}
    
where $T$ is the target for this vehicle and $\gamma_0 > 0$ is a given weighting coefficient, which may
depend on the vehicle. This minimization aims to align the vehicle with the target while preferring
trajectories with faster initial planned speeds, and gave good results when implemented.
\par
Simulations and experiments (see Secs.~\ref{ch4:sim} and \ref{ch4:exp}, respectively) have
shown that the proposed navigation strategy does allow every vehicle to efficiently converge to its
assigned target in most circumstances. Occasionally in scenarios with high vehicle densities, some
deadlocks may be avoided by instructing the vehicles to rotate in place to present different initial
states to the planner.

 \section{Simulations}
\label{ch4:sim}

Simulations were carried out in MATLAB with a perfect unicycle model, with parameters given in
Table~\ref{paramsim}. TPM and TTM updates occurred at 1 Hz and 10 Hz, respectively. To examine the
sensor noise implications, vehicle state measurements were corrupted by bounded, random, uniformly
distributed zero mean errors. The magnitude of these was $0.1 m$ for the translational position
measurement and $0.02 rad$ for the heading measurement. Acquisition of data about obstacles was
modeled as rotation of the detection ray with the step of $0.157 rad$, finding the distance to the
obstacle reflection point, and interpolation of the reflection points by straight line segments.
Once the vehicles were within $0.5 m$ of their targets, the navigation task was assumed to be completed
and so the navigation signals were set to zero.

\begin{table}[ht]
\centering

\begin{tabular}{| l | c |}
\hline
 $R_{sen}$ & $4.0m$  \\
\hline
 $C$ & $8.5m$  \\
 \hline
 $d_{tar}$ & $0.50 m$  \\
 \hline
 $v_{max}$ & $1.1 ms^{-1}$ \\
  \hline
  $u_{\theta, max}$ & $3.0 rads^{-1}$ \\
    \hline
  $u_{v,max}$ & $0.30 ms^{-2}$ \\
  \hline
   $v_{nom}$ & $1.0 ms^{-1}$ \\
  \hline
  \end{tabular} \hspace{10pt} \begin{tabular}{| l | c |}
  \hline
  $u_{\theta, nom}$ & $1.0 rads^{-1}$ \\
 \hline
 $u_{v,nom}$ & $0.20 ms^{-2}$ \\
\hline
 $\Delta \Lambda$ & $0.25$  \\
 \hline
 $k_0$ & $1$  \\
 \hline
 $k_1$ & $1$  \\
 \hline
 $\gamma_0$ & $1$\\
 \hline
 $\mu_{\varkappa}$ & $0.9$\\
\hline

\end{tabular}

\caption{Simulation parameters for multiple-vehicle controller.}
\label{paramsim}
\end{table}

The path tracking surfaces were taken to be:

\begin{align*}
\chi(z) &= \sgn(z) \cdot \min\{ 0.50 \cdot |z|, 0.50 \} \\
\unk(z) &= 0.05 + \min\{ 1.0 \cdot |z|, 0.70\vveh\}\\
p(S,z) &= \min\{\unk(z) + 1.0S, \vveh\} \\
\end{align*}

\begin{figure}[ht]
	\centering
		\includegraphics[width=\hsize]{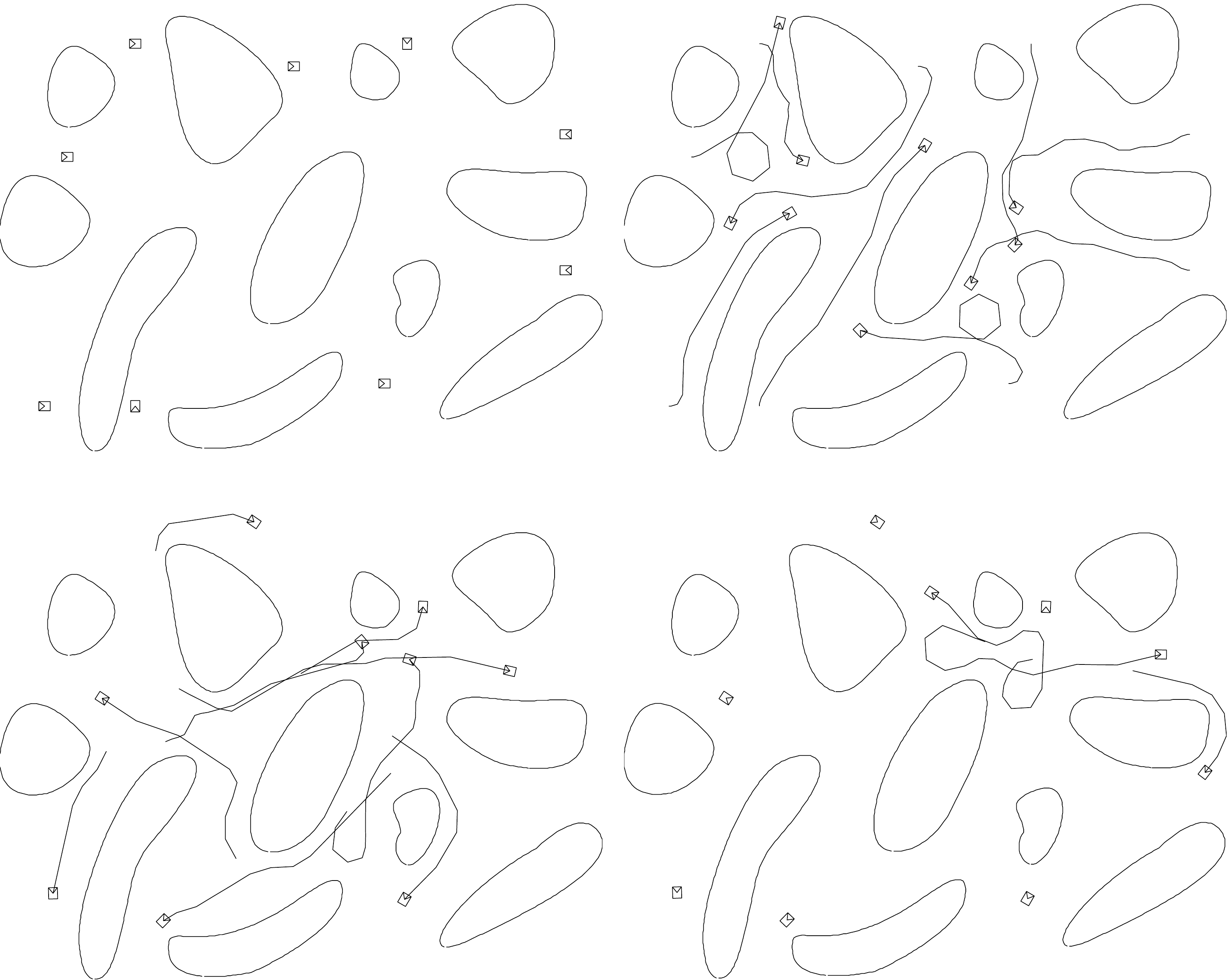}
	\caption{Simulations with nine vehicles in a complex scene.}
	\label{fig:four}
\end{figure}

\begin{figure}[ht]
	\centering
		\includegraphics[width=0.75\hsize]{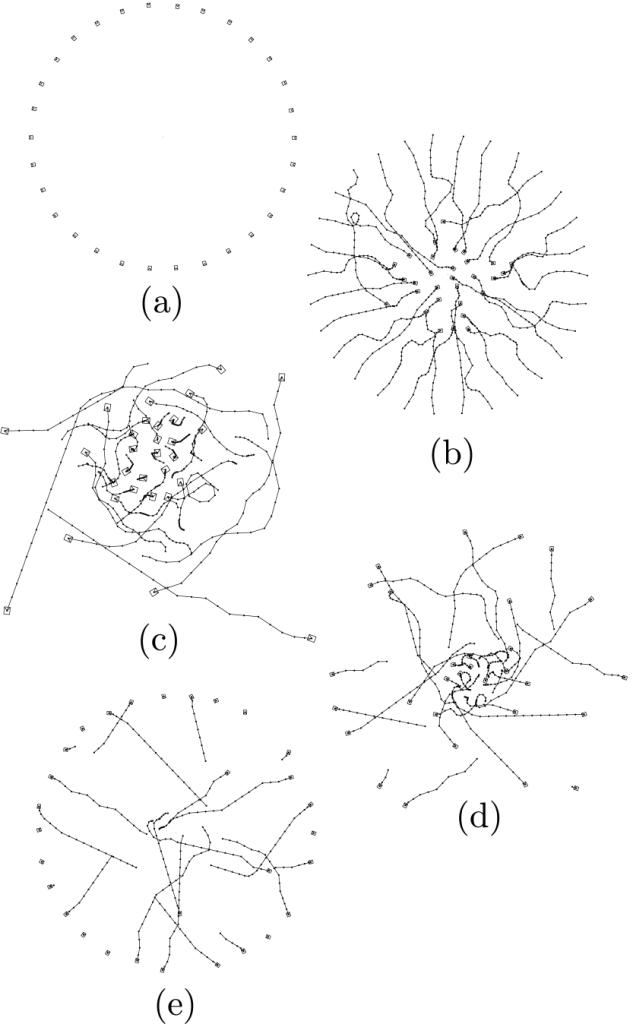}
	\caption{Simulations with thirty vehicles in an open scene.}
	\label{fig:thirty}
\end{figure}

\begin{figure}[ht]
	\centering
		\includegraphics[width=0.7\columnwidth]{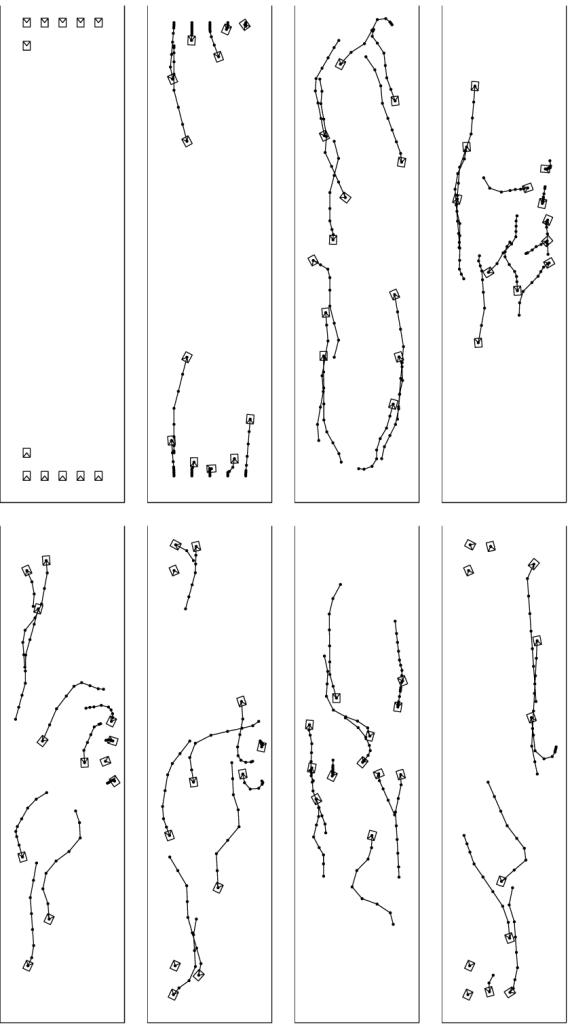}
	\caption{Simulations with twelve vehicles in a constrained environment.}
	\label{fig:sm}
\end{figure}

The obtained closed loop trajectories are displayed in Figs.~\ref{fig:four}, \ref{fig:thirty}, and
\ref{fig:sm}. In all cases when the noise was absent from the system, the inter vehicle distance was constantly 
monitored by the simulation, and was observed to never
drop below $0.995 m$. This confirms the theoretical results -- any slight deacreases below $d_{tar}$ can be explained by effects such as limited computational accuracy of the computer. The vehicle-obstacle distance never dropped below $0.495 m$, in harmony
with the chosen value of $d_{tar}$. When the sensor noise was added, as in the simulations shown
in Figs.~\ref{fig:four} to \ref{fig:sm}, the minimum observed inter-vehicle distance was $0.6547 m$. This can be treated as disturbance, and the desired inter-vehicle distance can be increased accordingly if this separation is too small. In every case, the vehicles successfully converged to the assigned targets, for example in Fig.~\ref{fig:sm} the time taken for all vehicles to reach their targets was $70 s$.

Fig.~\ref{fig:four} shows a collision avoidance maneuver for nine vehicles. This scenario may be
encountered in an office, warehouse, factory or urban environment. The vehicles smoothly move around
each other and adjust their speeds to prevent collisions. When a vehicle slows down its turning radius decreases, 
thus allowing sharper turns to be made (any of the points on the trajectory where the curvature is small would have been carried out at less than full speed). In addition, an interesting phenonemon was `wheeling', where vehicles would move in a full circle, presumably to dealy its entry into an adjacent area.

Fig.~\ref{fig:thirty} represents an obstacle free environment where thirty vehicles are initially
aligned around the edge of a circle; every vehicle must exchange positions with the vehicle on the
opposite side of the circle. Similar scenarios may be found in the literature
\cite{van2009conf0}. While the motion looks like relatively chaotic, all vehicles
still converged to the desired locations. At some points, some vehicles were stationary, and
performing the rotation maneuver allows to find a
trajectory after a short delay.

Fig.~\ref{fig:sm} shows the results with twelve vehicles in a constrained environment. This
experiment highlights that in complex situations higher order planning is useful to generate
efficient solutions. However even in this situation, the vehicles were all able to progress towards
their respective targets. Most importantly, collisions were always avoided.

\begin{figure}[ht]
	\centering
		\includegraphics[width=0.8\mylength]{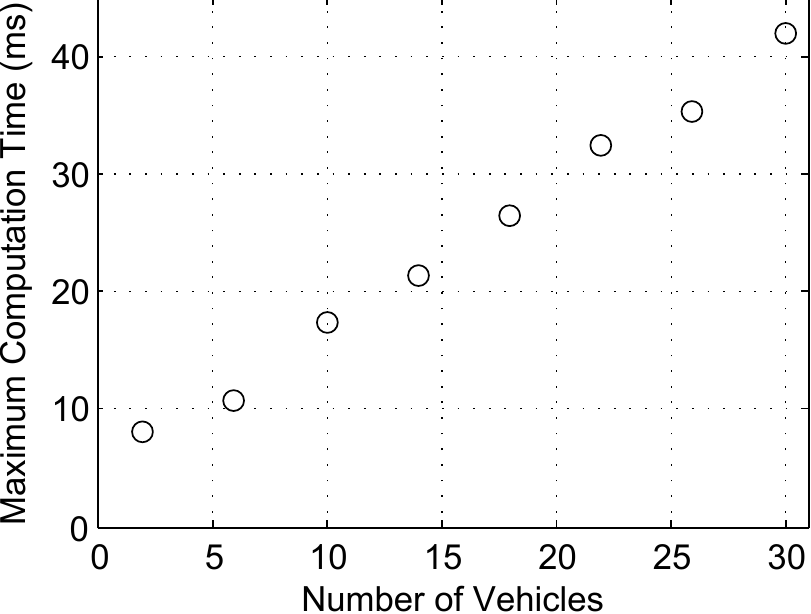}
	\caption{Computation time relative to the number of vehicles present.}
	\label{fig:comptime}
\end{figure}

Fig.~\ref{fig:comptime} shows the computation time of each execution of the algorithm on each
vehicle on a 3 GHz Desktop PC. The maximum CPU time was under $50ms$ per vehicle per control update, even for thirty vehicles,
which implies that the proposed approach is suitable for real time implementation.

\subsection*{Comparisons}

In order to compare the performance of the proposed algorithm with other types of collision
avoidance controllers, simulations were also carried out with a potential field type method
\cite{Chang2003conf8, Hernandez2011book2}. Various elements from these papers were combined to produce a
controller that solves a problem similar to ours. In doing so, the nearest point on the static
obstacle was treated as a virtual vehicle, as in \cite{Chang2003conf8}. The overall navigation
law is given by the following expression, loosely based on \cite{Hernandez2011book2}:

\begin{equation*}
\left[ \begin{array}{c}v_t(t) \\w_t(t) \end{array} \right] = A^{-1}(\theta) \cdot \left(
-\frac{1}{2}\grad{\gamma} -\sum_i \grad{V_i} \right)
\end{equation*} 

Here $\gamma$ is a navigation function
accounting for the distance to the target, $V_i$ is a repulsion potential with respect to the $i$th
companion vehicle, and $A(\theta)$ translates the control signals from the vehicle's reference frame
into the global one:

\begin{equation*}
A(\theta) = \left[ \begin{array}{c c}\cos(\theta) & - l \sin(\theta)\\ \sin(\theta) & l
\cos(\theta) \end{array} \right] 
\end{equation*}

\begin{equation*}
V_i = \left\{ \begin{array}{c c} \left(\frac{1}{\beta_i^2} - \frac{1}{d_{sep}^2} \right)^2 &
\beta_i \geq d_{sep}\\0 & \beta_i < d_{sep} \end{array} \right\}, 
\end{equation*} 

where $\beta_i$ is the distance
to the $i$th vehicle.
The control signals produced by the above formulas were subjected to saturation with regard to their
maximally feasible levels $u_{\theta, nom}$ and $u_{v,nom}$.

The controls were updated with a sampling time of $0.1 s$, and $d_{sep} := 4 m, l := 1 m$. The
vehicle completed the task when entering the disk of the radius $5 m$ centered at its target.

\begin{figure}[ht]
	\centering
		\includegraphics[width=0.7\columnwidth]{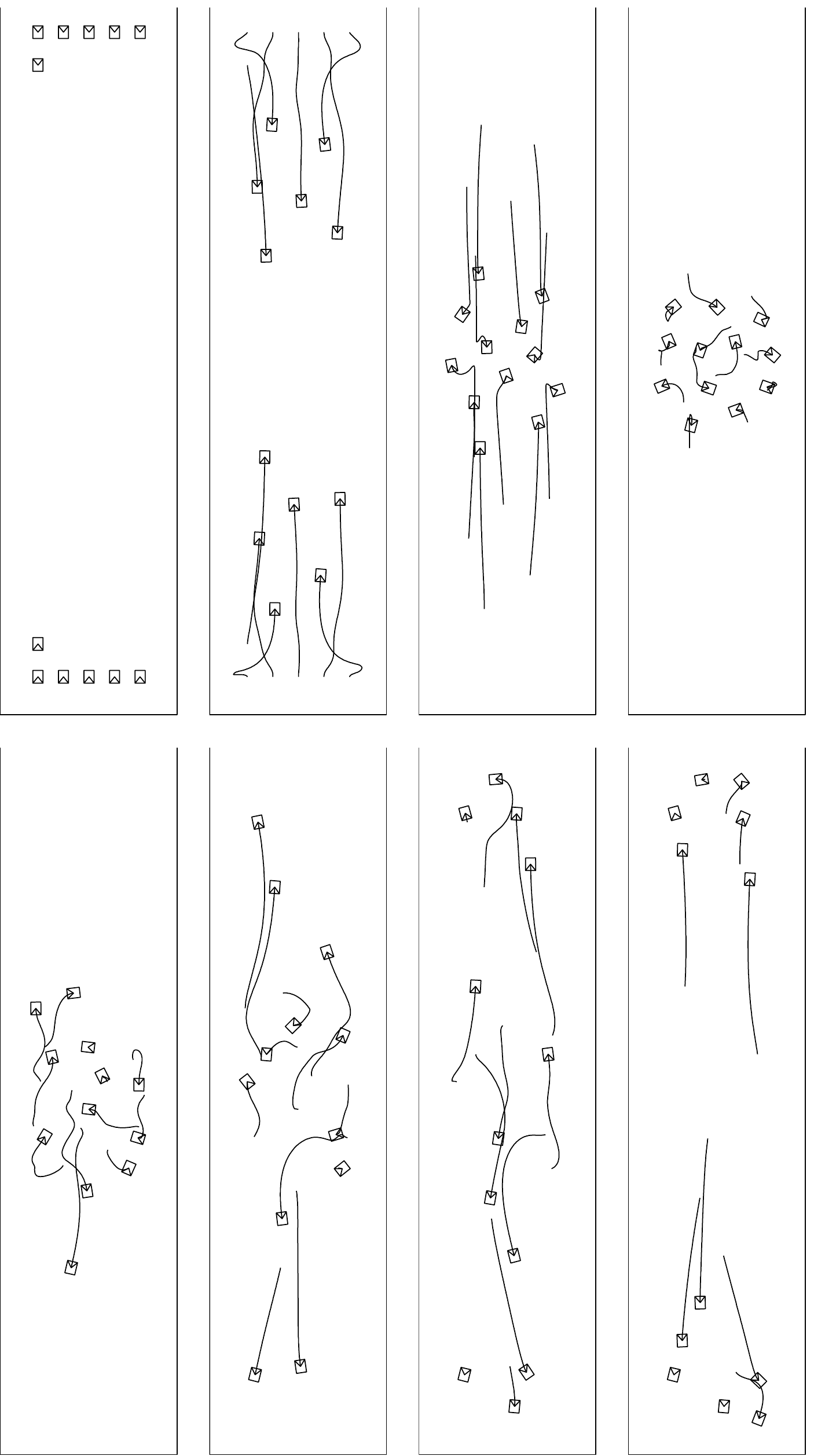}
	\caption{Simulations using a artificial potential field approach.}
	\label{fig:sp}
\end{figure}

The potential field method (PFM) was tested in favorable circumstances where its well known
disadvantages, like deadlocks caused by local minima due to, e.g., obstacle concavities, are
unlikely to manifest themselves. The corresponding scene and results are shown in Fig.~\ref{fig:sp}.
Simulations showed that PFM took  $96.7 s$ to perform the maneuver, whereas the proposed navigation law
(abbreviated to PCL in this discussion for clarity) took only $70 s$, which is essentially quicker. At the same time, the total distance traveled
by all vehicles is approximately the same for both methods: $280.52 m$ for PFM and $282.40 m$ for
PCL. The minimum distance between the vehicles was $0.951 m$ for PFM, whereas it was $1.001 m$ for
PCL. It is also worth noting that it was not difficult to find a simple alternative scenario where
the above PFM failed to drive the vehicles to their assigned targets due to local minima, whereas
PCL still succeeded.

\clearpage \section{Experiments}
\label{ch4:exp}
Experiments were carried out with two Pioneer P3-DX mobile robots to show real-time applicability of
this system. As in Chapt.~\ref{chap:convsingle}, a SICK LMS-200 laser range-finding device was used to detect obstacles in a vicinity of
the vehicle. In order to localize the
vehicles, they were launched from known starting positions and headings, and odometry feedback was
used during the maneuver.
\par
In the experiments presented here, this method resulted in a good enough performance, which is in harmony with
\cite{Martinelli2002journ1}, where the systematic error was reported to be under $1.5 \%$ of the distance
traveled.\footnote{For longer experiments, accumulation of odometry errors over time is probable.
To compensate for this, absolute positioning systems such as Indoor GPS or camera tracking may be
employed.}
 At each control update, as before the target translational acceleration and turning rate of the low level
 wheel controllers were set to match the results of the navigation algorithm. The laser detections
 in the circle of radius $d_{tar}$ around the estimated position of the other vehicle were excluded
 when computing the navigation approach to prevent detections of the other vehicle affecting the path
 computations. The values of the parameters used in the experiments are listed in
 Table~\ref{ch4:paramexp}.
 
\begin{table}[ht]
\centering

\begin{tabular}{| l | c |}
\hline
 $R_{sen}$ & $4.0m$  \\
\hline
 $d_{tar}$ & $0.30 m$  \\
 \hline
 $v_{max}$ & $0.50 ms^{-1}$ \\
  \hline
  $u_{\theta, max}$ & $0.80 rads^{-1}$ \\
    \hline
  $u_{v,max}$ & $0.30 ms^{-2}$ \\
  \hline
 $v_{nom}$ & $0.40 ms^{-1}$ \\
 \hline
  $u_{\theta, nom}$ & $0.60 rads^{-1}$ \\
 \hline
 \end{tabular} \hspace{10pt} \begin{tabular}{| l | c |}
 \hline
 $u_{v,nom}$ & $0.10 ms^{-2}$ \\
\hline
 $\Delta \Lambda$ & $0.50$  \\
 \hline
 $k_0$ & $1$  \\
 \hline
 $k_1$ & $1$  \\
\hline
$\gamma_{vel}$ & $1$ \\
\hline
$\mu_{\varkappa}$ & $0.9$ \\
\hline

\end{tabular}

\caption{Experimental parameters for multiple-vehicle controller.}
\label{ch4:paramexp}
\end{table}

The navigation scheme proposed in this chapter assumes decentralized communication. However for purely
technical reasons, this decentralization was emulated on a common base station in the experiments so
that the vehicles communicated with each other not directly but through this station. Communication
was established via  sending UDP packets over a 802.11g WLAN, with a checksum error detection scheme
being employed \cite{Wijesinha2005conf0}. A substantial rate of packet dropouts was observed, which according to
Sec.~\ref{sec.sel} does not affect robustness but does affect performance. The event of packet
dropout was handled similarly to the event where refinement of the set of planned trajectories results
in the empty set -- the previous probational trajectory remains in use.
TTM and TPM were run at a common rate, the communication range covered the entire operational zone
of the vehicles.

 \begin{figure}[ht]
\centering
\subfigure[]{\includegraphics[width=0.7\mylength]{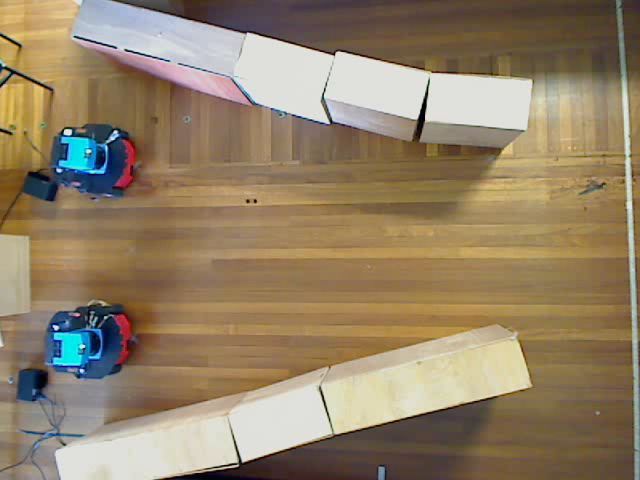}}
\subfigure[]{\includegraphics[width=0.7\mylength]{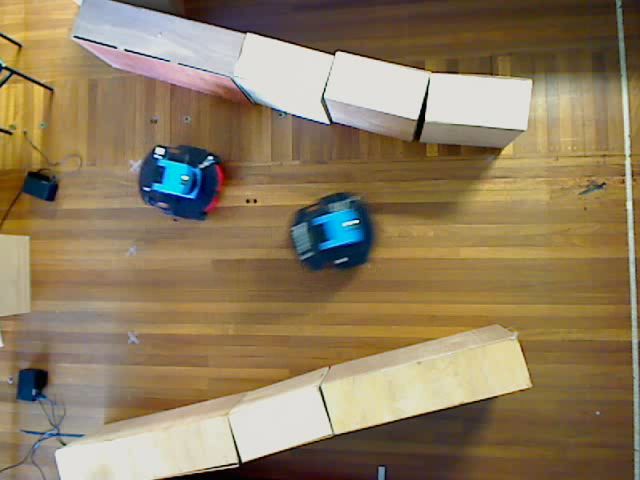}}
\subfigure[]{\includegraphics[width=0.7\mylength]{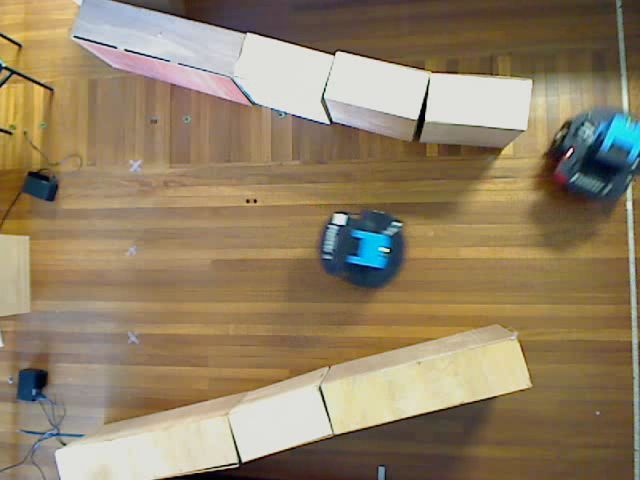}}
\subfigure[]{\includegraphics[width=0.7\mylength]{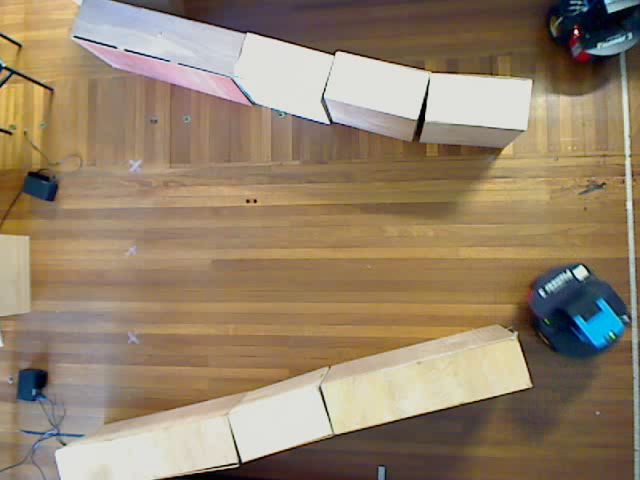}}
\caption{Sequence of images showing the experiment (parallel encounter).}
\label{ch4:para}
\end{figure}

 \begin{figure}[ht]
\centering
\includegraphics[width=1.2\mylength]{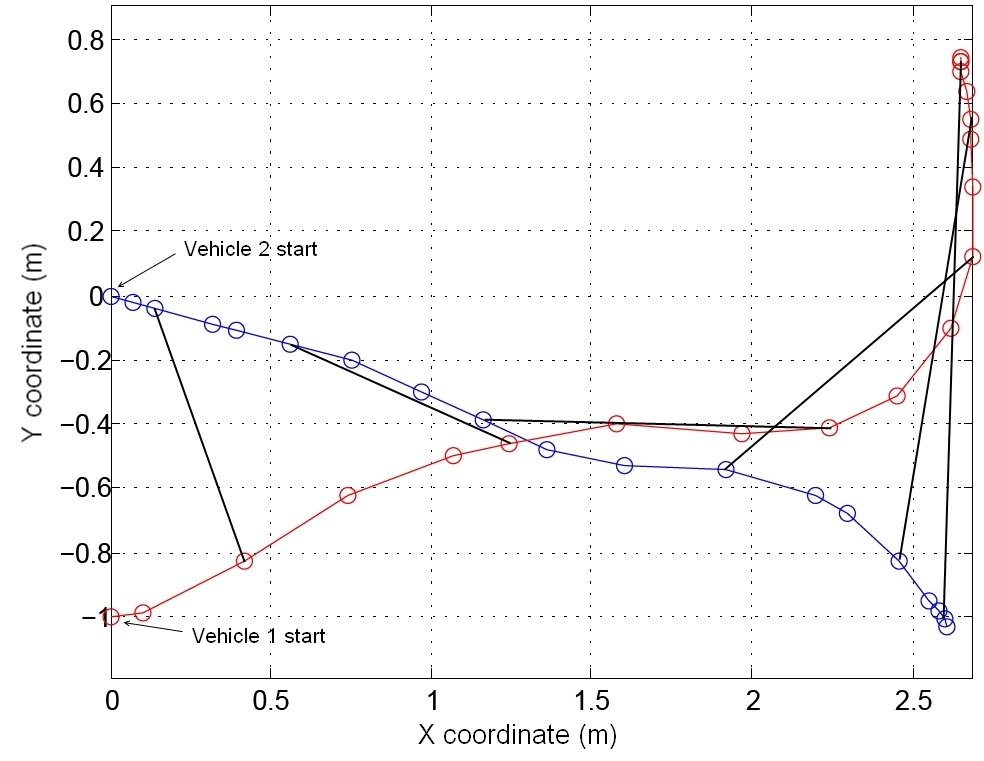}
\caption{Trajectory for the experiment (parallel encounter).}
\label{ch4:paraplot}
\end{figure}

In the first scenario demonstrated by Figs.~\ref{ch4:para} and \ref{ch4:paraplot}, the vehicles
are both trying to move through a gap which is only wide enough for one vehicle.
The top vehicle allows the bottom vehicle to pass through a narrow region before moving into the gap itself.

 \begin{figure}[ht]
\centering
\subfigure[]{\includegraphics[width=0.7\mylength]{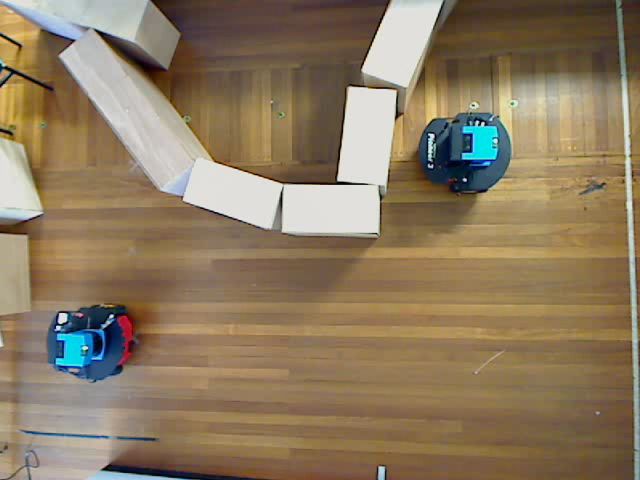}}
\subfigure[]{\includegraphics[width=0.7\mylength]{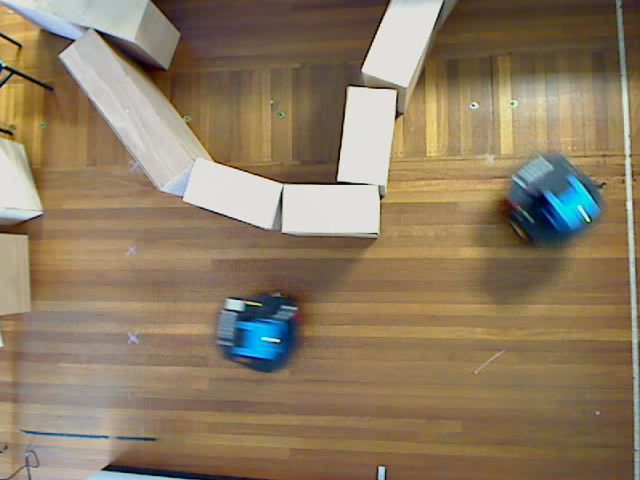}}
\subfigure[]{\includegraphics[width=0.7\mylength]{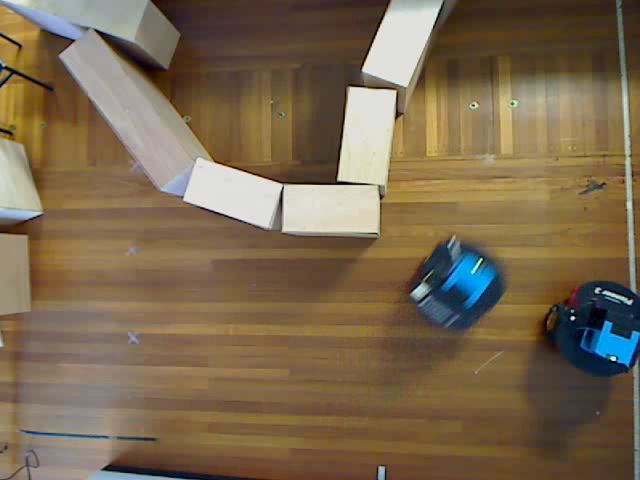}}
\subfigure[]{\includegraphics[width=0.7\mylength]{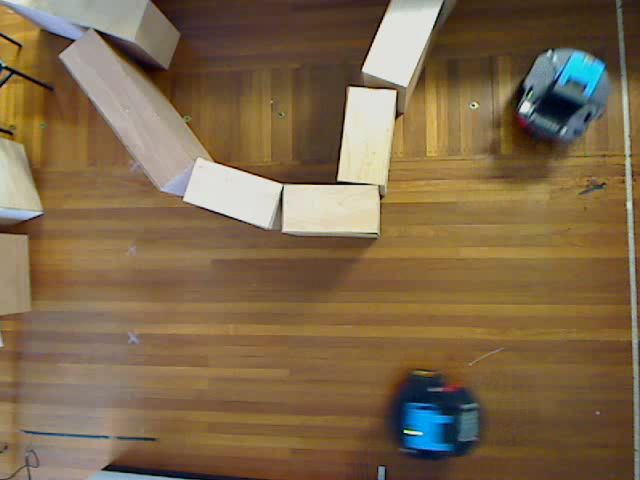}}
\caption{Sequence of images showing the experiment (head-on encounter).}
\label{fig:headon}
\end{figure}

 \begin{figure}[ht]
\centering \includegraphics[width=1.2\mylength]{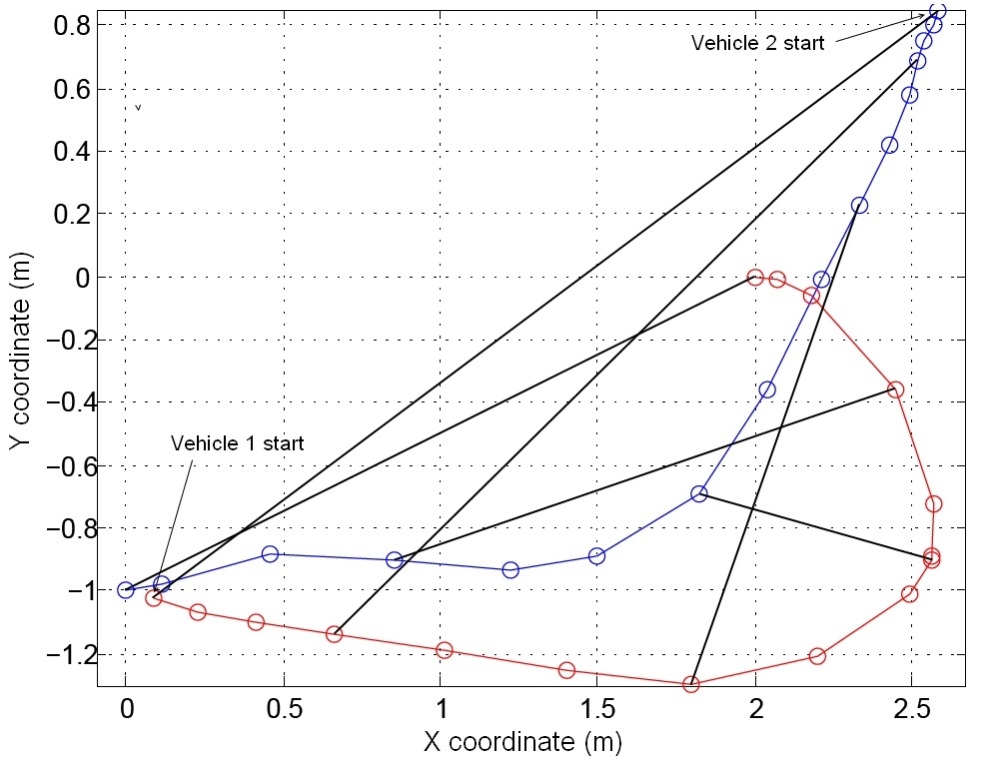}
\caption{Trajectory for the experiment (head-on encounter).}
\label{fig:headonplot}
\end{figure}

In the second scenario demonstrated by Figs.~\ref{fig:headon} and \ref{fig:headonplot}, the
vehicles must pass each other while a static obstacle is present. The vehicle which begins on the
left takes a longer path to prevent collision occurring.

\clearpage \section{Summary}
\label{ch4:con}

In this chapter, a robust cooperative collision avoidance system is proposed, which is suitable for maintaining safety within a
group of vehicles operating in unknown cluttered environments. By calculating a set of presumable planned trajectories for each
vehicle within communication radius, it is possible to generate constraints which maintain mutually feasible trajectories, from which
collision avoidance behaviour may be shown. The system was verified with
simulations and tests on a pair of real wheeled mobile robots.

\chapter{Deadlock Avoidance}
\label{chap:dead}

In recent times, many types of decentralized navigation strategies have been proposed to coordinate
multiple vehicles and prevent collisions between them. In these situations, showing robust
convergence to target is a equally important issue, however deadlock avoidance\footnote{In this
chapter no differentiation is made between the terms deadlock and livelock} is generally a more
complex and non-trivial subproblem which must be solved.
\par
The approach presented in this chapter is motivated by the fact that, at least anecdotally, 
deadlocks most commonly occur in narrow corridors, i.e. where there is not enough
room for vehicles to move past each other. Here, a hierarchy is enforced between the vehicles during
encounters, which means the highest ranked vehicle will always have priority for transversing the
corridor environment. This work is most applicable to unknown environments, where a map suitable for
segmenting into a graph is not available. Also, no assumptions are made about the shape or geometry of
the obstacle as was done in Chapt.~\ref{chap:convsingle}. However to substitute this information,
 it will be assumed that a suitable navigation function exists.
 \par
\begin{figure}[ht]
\centering
\subfigure[]{\includegraphics[width=0.49\mylength]{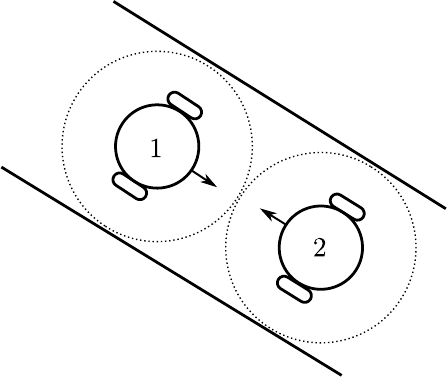}}
\subfigure[]{\includegraphics[width=0.49\mylength]{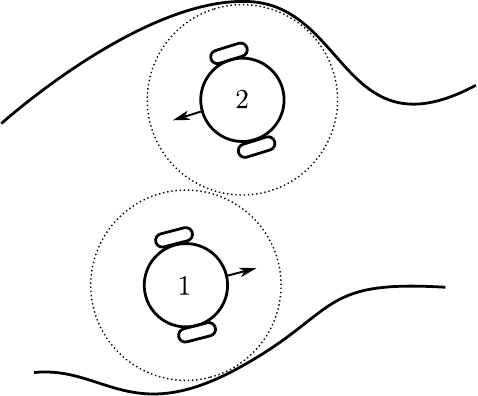}}
\subfigure[]{\includegraphics[width=0.49\mylength]{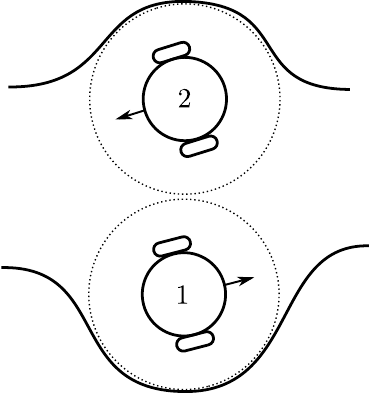}}

\caption{Different categories of deadlock; (a) Type I; (b) Type II; (c) Type III.}
\label{fig:para}
\end{figure}

 As mentioned in the introduction, most research into deadlock free motion planning divides the
workspace into some type of discrete graph, where the vehicle travels between nodes. This chapter is
concerned with preventing deadlocks from a continuous state space perspective, where twin vehicle deadlocks
may be classified into three categories; Type I deadlock may be resolved without the primary vehicle
moving contrary to its objective; Type II deadlock requires the primary vehicle to temporarily
move contrary to its objective; and finally, Type III deadlocks are impossible to solve. These categories are
illustrated in Fig.~\ref{fig:para}.

In this chapter, almost all deadlocks are of the first kind. The second type are extremely rare, and
would require use of the recovery scheme proposed in Sec.~\ref{sec:recovery} to solve. The third type may only exist if the vehicles
starting positions are inappropriate, and this is avoided by Assumption~\ref{as:space}.
 
\par
The body of this chapter is organized as follows. In Sec.~\ref{ch5:ps} the
problem statement is given; in Sec.~\ref{ch5:csa} a new navigation system is
proposed to avoid deadlocks. Simulations in a corridor environment are in
Sec.~\ref{ch5:sim}; experiments are in Sec.~\ref{ch5:emp}. Finally, Sec.~\ref{ch5:con}
offers brief conclusions.

\section{Problem Statement}
\label{ch5:ps}

Two autonomous vehicles $\mathcal{R}_1, \mathcal{R}_2$ are considered, which travel in a plane with maximum
speed $v_{max}$, each of which is associated with a steady point target $T_1, T_2$ respectively.
$\mathcal{R}_1$ is designated the primary vehicle, which is given priority for navigation. As before the plane
contains a set of untransversable, static, and closed obstacles $D_j \not\ni T_i, j \in [1:n]$. The
objective is to design a navigation law that drives every vehicle to their assigned targets in finite
time through the obstacle-free part of the plane $\mathbb{R}^2 \setminus D, D:= D_1\cup \ldots \cup
D_n$. Moreover, the distance from the vehicle to every obstacle and other vehicle should constantly
exceed the given safety margin $d_{sfe}$.  As in previous chapters, a unit time-step is used for
simplicity.
\par
In this chapter, only the unicycle model is considered for brevity, through these
results could easily be extended to the holonomic model.
As in Chapt.~\ref{chap:convsingle}, it is assumed that the system model is free from
disturbance.
\par
The trajectory structure and feasibility criteria employed in this chapter are identical
to those in Chapt.~\ref{chap:multiple}.
The main difference are minimum requirements on the sensed ares $F_{vis}(k)$.
This generalized obstacle set characterization is described in the following assumption:

\begin{Assumption}
The vehicle has perfect knowledge of the obstacle set within some distance $R_{con} < R_{sen}$ of the vehicle, such that:

\begin{equation}
F_{vis}(k) \supset \left\{\bldx \in F: \norm{\bldx - \blds(k)} \leq R_{con}\right\}
\end{equation}
\end{Assumption}

Communication assumptions are mostly identical to Chapt.~\ref{chap:multiple}, where vehicles
have their control updates synchronized and experience unit communication delay.
However in contrast to Chapt.~\ref{chap:multiple}, it is assumed that
communication is always possible regardless of the distance between the
vehicles.
\par
In this chapter, the trajectory planning approach presented
in Chapt.~\ref{chap:singlevehicle} must be slightly modified.
This is because minimum requirements of the trajectory planner must be formulated in
order to ensure convergent navigation.
This involves allowing for sampling effects through a finite quantity $d_{pln} >
0$ as follows:

\begin{Definition}
A point is \textbf{plannable} if the minimum distance to the obstacle set $D$ exceeds $d_{tar} +
d_{pln}$. The set of plannable points is denoted $G \subset F$.
\end{Definition}

The corresponding assumptions on the path planning system are located in Assumption~\ref{ass:planreq}.
Any point not within $d_{pln}$ of the set of plannable points would be
effectively be considered as non-feasible, through this is achieved without the consideration of the navigation system.
\par
It is also assumed that a global navigation
function (see e.g. \cite{Ogren2005journ3}) is available and known to the vehicle:

\begin{Assumption}
For any point $\bldx$, the vehicle has access to a navigation function $\mathscr{N}_i(\bldx)$. If $\bldx \in G$, $\mathscr{N}_i(\bldx)$ is defined as the minimum
length of any curve intersecting both $\bldx$ and $T_i$, which is contained within the set of plannable points $G$. 
If such a curve does not exist, $\mathscr{N}_i(\bldx) := \infty$.
If $\bldx \not\in G$, let $\acute{\bldx} := \arg\min_{\acute{\bldx} \in G}\norm{\bldx - \acute{\bldx}}$.
Then $\mathscr{N}(\bldx) := \mathscr{N}(\acute{\bldx}) + \norm{\bldx - \acute{\bldx}}$.
\end{Assumption}

It is assumed that the target $T_i$ is connected to the vehicle, such that the navigation function $\mathscr{N}$ is defined at both
vehicles' starting positions.
\par
Based on the navigation function, the definition of a convergent trajectory can be introduced as follows:

\begin{Definition}
A probational trajectory $s^*(j|k)$ is \textbf{convergent} if the navigation function $\mathscr{N}$  decreases by at least a decrement
$d_{dec}$ from the initial point to any other point on the trajectory:

\begin{equation}
\max_{1\leq j \leq \tau} \mathscr{N}(\blds^*(j|k)) \leq \mathscr{N}(\blds(k)) - d_{dec}
\end{equation}
 
It is permissible for only $\blds^*(\tau|k)$ to be considered so long that the trajectory $\blds^*(j|k)$ is followed to
completion. In this case $\mathscr{N}(\blds^*(\tau|k)) \leq \mathscr{N}(\blds(k)) - d_{dec}$.
\end{Definition}

It is required that the given trajectory planning system is sufficiently complete, so that if there exists a convergent trajectory
consisting of plannable points, the trajectory generation scheme will not fail to find at least one convergent,
feasible trajectory.
\par
Finally it is assumed that the targets and initial positions are sufficiently spaced from the obstacle, so
that the other vehicle may successfully navigate around them:

\begin{Assumption}
\label{as:space}
The points $T_1$, $T_2$,$\blds_1(0)$, $\blds_2(0) \in G$, and are spaced by at least $(2d_{tar} + 2d_{pln})$ from $D$. 
\end{Assumption}

\begin{Remark} \rm
The use of a navigation function assumes the environment is previously known, and
prevents the method from being classified as sensor based. It may be possible to replace the navigation function with
a heuristic similar to the target convergence scheme in Chapts.~\ref{chap:convsingle} and \ref{chapt:rbf}.
\end{Remark}

 \section{Navigation System Architecture}
\label{ch5:csa}

The basic premise for deadlock avoidance is for the secondary vehicle to attempt to ensure the
primary vehicle always has a valid trajectory to follow in the future. Additionally, the planning
algorithm is extended with an additional trajectory shape. Also, a prescribed
maneuver (the recovery scheme) is described for the vehicles to execute when all other possibilities are exhausted. This maneuver is
guaranteed to result in at least a given progression along the navigation function of the primary
vehicle.
\par
Recall in Chapt.~\ref{chap:singlevehicle} only two possible longitudinal control 
patterns were considered -- either the pattern $p_+:=(+ - - \ldots)$ or the pattern
$p_-:=(- - - \ldots)$. In this chapter a third pattern $p_s$ is introduced, which is only considered when the vehicle
is stationary ($v(k) = 0$), where $p_s:=(0 \thickspace 0 \thickspace 0 \thickspace \ldots \thickspace 0 \thickspace 0 \thickspace + -)$. 
The overall effect of executing this trajectory is to introduce a small step into the vehicle's
position (see Lemma~\ref{lem:stepp}). The remainder of the planning algorithm is identical to the system
 documented in Chapt.~\ref{chap:singlevehicle}.

\subsection{Primary Vehicle}
Probational trajectories for the primary vehicle must be convergent. In additional to transmitting
the chosen probational trajectory $\blds^*(j|k)$, the primary vehicle also
considers its projected state in the next time step $\blds^\ast(1|k)$, and generates a set of
trajectories which are also convergent. This set of trajectories is labeled $\mathscr{P}_{DF}$, and
is transmitted to the secondary vehicle.

\subsection{Secondary Vehicle}
For the secondary vehicle, in addition to the collision avoidance constraints, any chosen probational trajectory of the
secondary vehicle must not interfere with at least one trajectory from $\mathscr{P}_{DF}$. It must
also attempt to maximize some cost function $J$, which reflects some trade-off between the maximum distance the primary vehicle progresses
in the non-intersected subset of $\mathscr{P}_{DF}$, the distance between the vehicles if below some
threshold, and the distance the secondary vehicle progresses to its target:

\begin{multline}
J = \sup_{\blds^*_{DF,i}(j|k) \in \mathscr{P}_{DF}} \Biggl[\mathscr{N}_i(\blds^*_{DF,i}(\tau|k))\\
 - \gamma_0 \cdot \max \left \{ (2 + c) \cdot d_{tar}, \inf_{1 \leq j \leq \tau} \left \Vert \blds^*_{DF,i}(j|k) - \blds^*(j|k) \right \Vert \right \} \Biggr] \\
+ \gamma_1 \cdot \mathscr{N}(\blds^*(\tau|k ))
\end{multline}

Here $\gamma_0$, $\gamma_1$ and $c$ are tunable parameters. Once the primary vehicle is positioned 
at its target the secondary vehicle assumes it is stationary and modifies its navigation function to
treat it as a static point obstacle.

\begin{Assumption}
\label{ass:planreq}
There are several relations between variables which must be maintained to show the proposed trajectory planning system
is sufficiently complete:

\begin{subequations}
\begin{align}
d_{sns} &:= u_{v,nom} < \min\left\{d_{pln}; R_{con} - (d_{tar} + d_{pln}) \right\}\\
 d_{spc} &:= \Delta \Lambda \cdot u_{\theta,nom} \cdot d_{sens} < d_{pln}\\
 \theta_{rng} &:= (\tau-2) \cdot u_{\theta,nom} > \pi\\
 d_{dec} &< d_{sens} - \frac{d_{spac}}{2}\\
\end{align}
\end{subequations}
\end{Assumption}

\begin{Lemma}
\label{lem:stepp}
The proposed trajectory generation system always generates a convergent trajectory when both
vehicles are stationary and spaced by at least $d_{mut} + 2d_{sens}$
\end{Lemma}

\begin{figure}[ht]
\centering
\includegraphics[width=0.8\mylength]{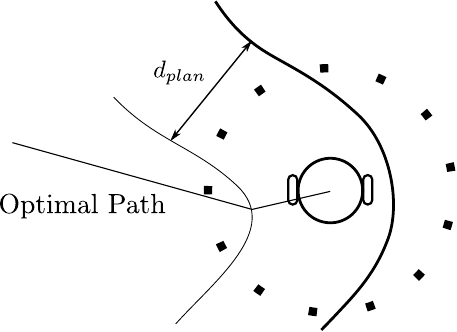}
\caption{Contingency trajectories for the planning system.}
\label{fig:plan}
\end{figure}

\pf The set of the possible final vehicle positions $\blds^*(\tau|k)$ corresponding to the current set of probational
trajectories is shown in Fig.~\ref{fig:plan}, labeled $\mathscr{F}$, and under the specified conditions they
are independent for each vehicle. These consist of a rotation with zero longitudinal velocity, followed by a
longitudinal control of $u_{v,nom}$, followed by a control of $-u_{v,nom}$. The figure also shows
the line corresponding to the steepest descent of the navigation function, which consists of
possibly a straight line to the set of plannable points if outside this set, followed by a curve
directed along the best route to the target. Intersect this curve with a circle of radius $d_{sns}$
around the vehicle, and choose the first intersection point $\mathscr{G}$. It may be seen the
minimum distance from $\mathscr{G}$ to $\mathscr{F}$ is less than $\frac{d_{spc}}{2}$. This closest
point $\mathscr{F}_c$ will always be valid since $\mathscr{G}$ is plannable. Thus the total
decrement to the navigation function by navigating to $\mathscr{F}_c$ will be at least $d_{sns} -
\frac{d_{spc}}{2}$. \epf

\subsection{Recovery Scheme}
\label{sec:recovery}

The recovery scheme is a fixed maneuver which may be undertaken to ensure deadlocks are avoided even
when the main scheme fails. Note this scheme is applicable to the worst case scenario and shortcuts
may be feasible that would improve efficiency without affecting overall motion. Assume at time-step
$\kappa$ the primary vehicle is stationary and cannot plan a convergent trajectory:

\begin{itemize}
\item The vehicles backtrack along their trajectories using opposite control inputs until the
secondary vehicle is able to plan and execute a trajectory resulting in a final stationary position
that satisfies Assumption~\ref{as:space}. In the absolute worst case scenario this will be the
secondary vehicle's starting position.

\item The secondary vehicle then remains stationary while the primary vehicle moves until its
navigation function has decreased by at least $d_{dec} + (\pi - 2) \cdot (d_{tar} + d_{pln})$
relative to its previous value at time-step $\kappa$. There is a slight modification to the planning
structure, as the primary vehicle assumes the secondary vehicle remains stationary, and modifies its
navigation function to treat it as a static point obstacle.

\item Normal operation ensues.
\end{itemize}

\begin{Remark} \rm
The recovery scheme at first appearances does not seem to be the most general solution, however it
has some justification. In this chapter, for simplicity, discussion the obstacle
boundary are avoided -- for example the properties of its shape and the available sensing are not discussed
like in Chapt~\ref{chap:convsingle}. For example, counterexamples
possibly could be constructed where the vehicle must `slide' along the other vehicle and the boundary,
only possible in a continuous time paradigm. The only information accessible to the navigation law is the validity
of a particular trajectory, and in such a case previous trajectories are the only given
means of guaranteeing resolution of a conflict.
\end{Remark}

\begin{Lemma}
The recovery scheme is always feasible and the primary vehicles navigation function always
decrements by at least $d_{dec}$ during each application.
\end{Lemma}

\pf No assumption about the availability of any other trajectory other than those already used has been
made. When the secondary vehicle is stationary and satisfies Assumption~\ref{as:space}, the maximum
absolute change in the navigation function is bounded by $2(d_{pln} + d_{tar})$, which bounds
the difference between going around and through the vehicle when it is treated as a static point. \epf

\begin{Remark} \rm
Note the recovery scheme may be completed in a finite number of time-steps. 
\end{Remark}

\begin{Proposition}
Both vehicles converge to their targets in finite time.
\end{Proposition}

\pf The primary vehicle only chooses trajectories that decrement the navigation function by $d_{dec}$
over the trajectory, and engages the recovery scheme if this cannot be done. Thus the navigation
function will always decrease by $d_{dec}$ in a finite number of time-steps. \epf

 \section{Simulations}
\label{ch5:sim}

\begin{figure}[ht]
\centering
\includegraphics[width=\columnwidth]{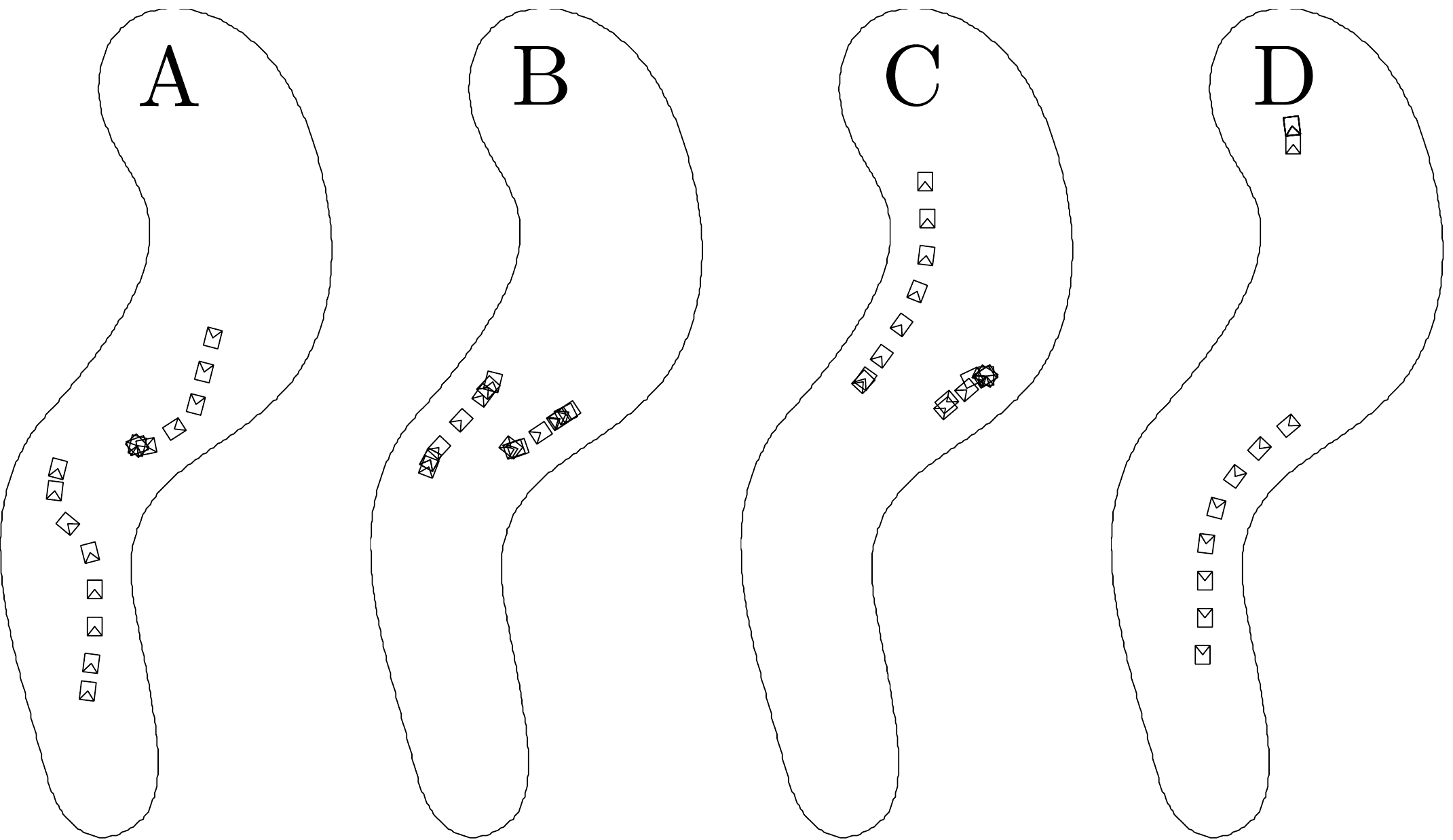}

\caption{Simulations with deadlock avoidance (head-on encounter).}
\label{fig:simon}
\end{figure}

\begin{figure}[ht]
 
\centering
\subfigure[]{\includegraphics[width=0.3\mylength]{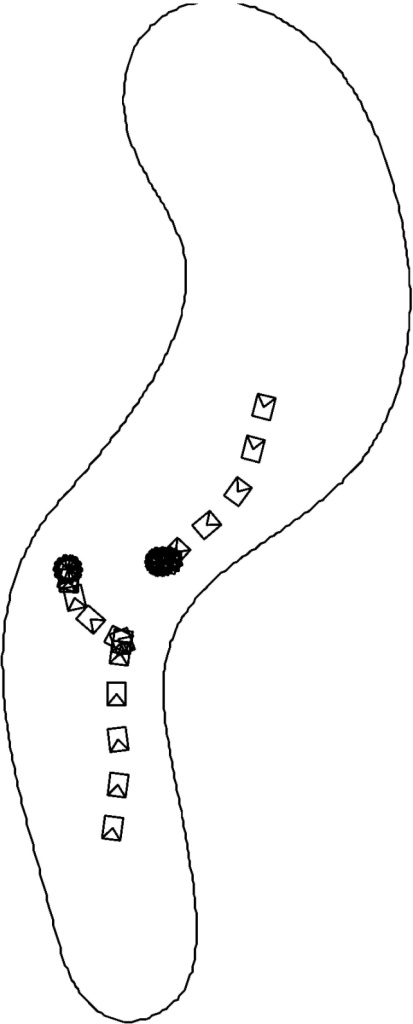}}
\hspace{20pt}
\subfigure[]{\includegraphics[width=0.5\mylength]{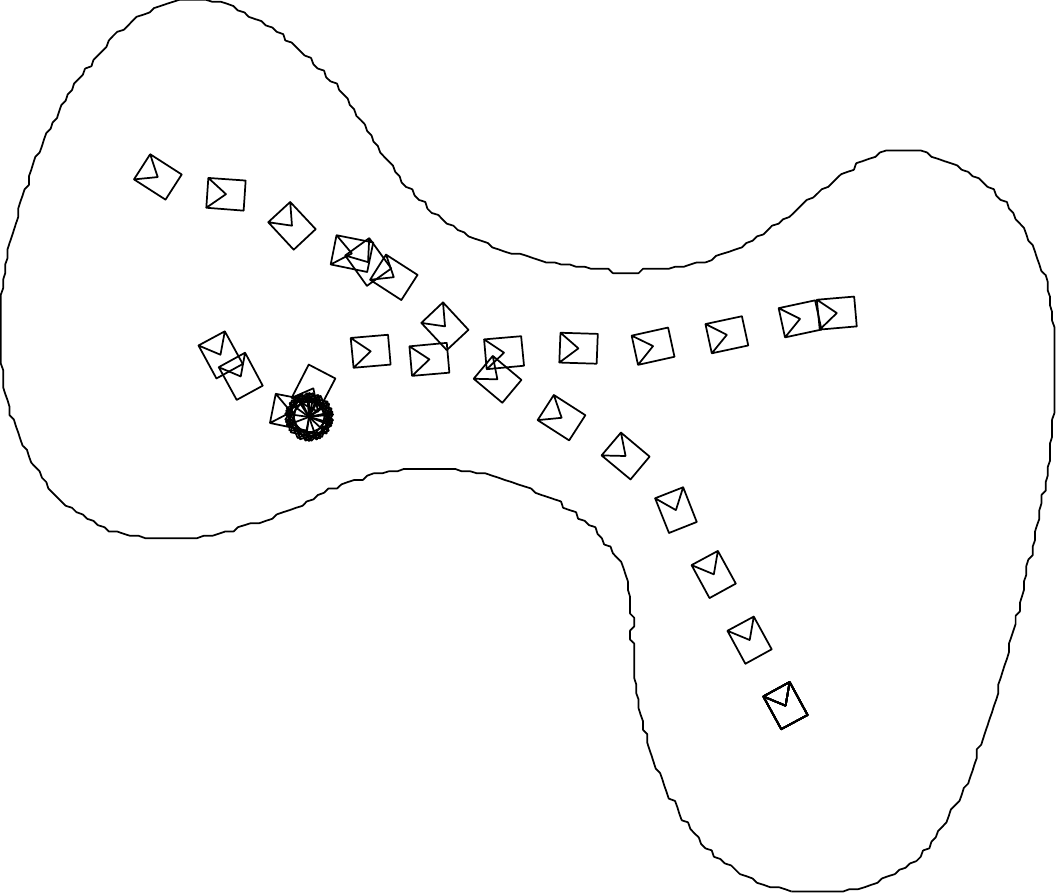}}
\caption{(a) Simulations without deadlock avoidance; (b) Simulations with deadlock avoidance (parallel encounter).}
\label{fig:sim2}

\end{figure}

The parameters corresponding to that navigation law are shown in the Table~\ref{ch5:paramsim}.
The navigation system was tested with both the deadlock avoidance constraints turned on and off. The
trajectory with the deadlock avoidance turned on is shown in Fig.~\ref{fig:simon}, and the
vehicles took $69$ time steps to both reach their respective targets. It can be seen the secondary
vehicle promptly moved backwards to allow the primary vehicle to pass. When the deadlock avoidance
constraints were turned off (see Fig.~\ref{fig:sim2}(b)) the vehicles failed to converge to their
targets. Note the latter scenario is not identical to Chapt.~\ref{chap:multiple} as in
that chapter the trajectory selection included a velocity maximization term which attempts 
to mitigate simple deadlocks.

\begin{table}[ht]
\centering

\begin{tabular}{| l | c |}
\hline
 $R_{sensor}$ & $4.0m$  \\
 \hline
 $d_{tar}$ & $1.0 m$  \\
\hline
 $v_{max}$ & $0.5 ms^{-1}$ \\
  \hline
  $u_{\theta, max}$ & $0.5 rads^{-1}$ \\
    \hline
 $u_{v,nom}$ & $0.20 ms^{-2}$ \\
\hline
\end{tabular} \hspace{10pt} \begin{tabular}{| l | c |}
\hline
 $\Delta \Lambda$ & $0.25$  \\
 \hline
 $c$ & $1$\\
 \hline
 $\gamma_1$ & $0.1$\\
 \hline
 $\gamma_2$ & $0.001$\\
 \hline
\end{tabular}
\label{ch5:paramsim}
\caption{Simulation parameters for deadlock-avoiding controller.}
\end{table}

In Fig.~\ref{fig:sim2}(b), the approach is used a different topology of encounter, showings its
versatility. The secondary vehicle loitered before the corridor before passing through.

Note that conditions that induce the recovery scheme were unable to be devised, and thus simulations
of this part of the controller are not included here.

\clearpage \section{Experiments}
\label{ch5:emp}
Experiments were carried out with the identical experimental setup as Chapt.~\ref{chap:multiple}.
 The updated values of the parameters used in the experiments are listed in Table~\ref{paramexp}.
 
\begin{table}[ht]
\centering

\begin{tabular}{| l | c |}
\hline
 $R_{sen}$ & $4.0m$  \\
\hline
 $d_{tar}$ & $0.3 m$  \\
 \hline
 $v_{max}$ & $0.5 ms^{-1}$ \\
  \hline
  $u_{\theta, max}$ & $0.8 rads^{-1}$ \\
 \hline
 $u_{sp}$ & $0.1 ms^{-2}$ \\
 \hline
 \end{tabular} \hspace{10pt} \begin{tabular}{| l | c |}
\hline
 $\Delta \Lambda$ & $0.50$  \\
\hline
 $c$ & $1$\\
 \hline
 $\gamma_1$ & $0.1$\\
 \hline
 $\gamma_2$ & $0.001$\\
 \hline
\end{tabular}

\caption{Experimental parameters for deadlock-avoiding controller.}
\label{paramexp}
\end{table}

\begin{figure}[ht]
\centering
\includegraphics[width=1.2\mylength]{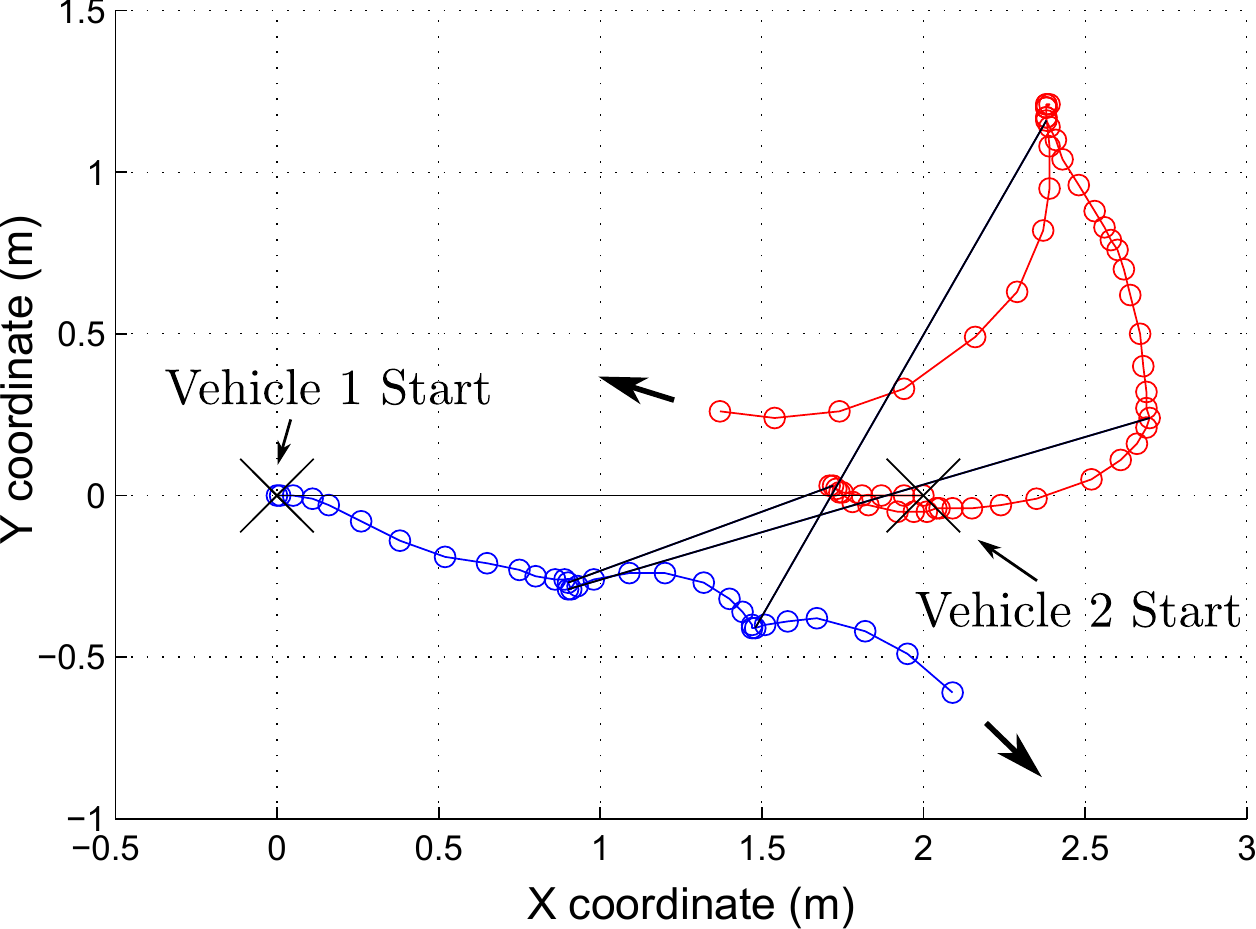}
\caption{Trajectory for the experiment.}
\label{fig:paraplottt}
\end{figure}

 \begin{figure}[ht]
\centering
\subfigure[]{\includegraphics[width=0.7\mylength]{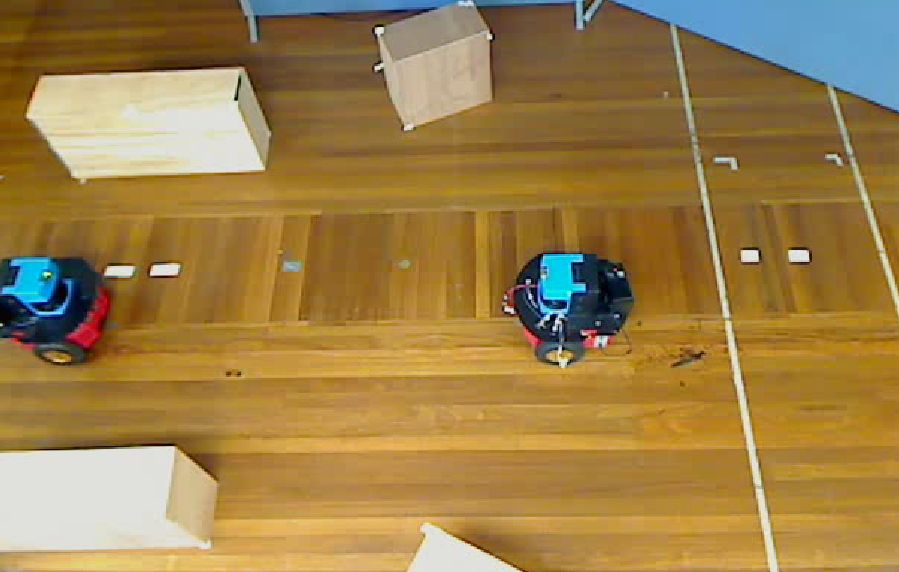}}
\subfigure[]{\includegraphics[width=0.7\mylength]{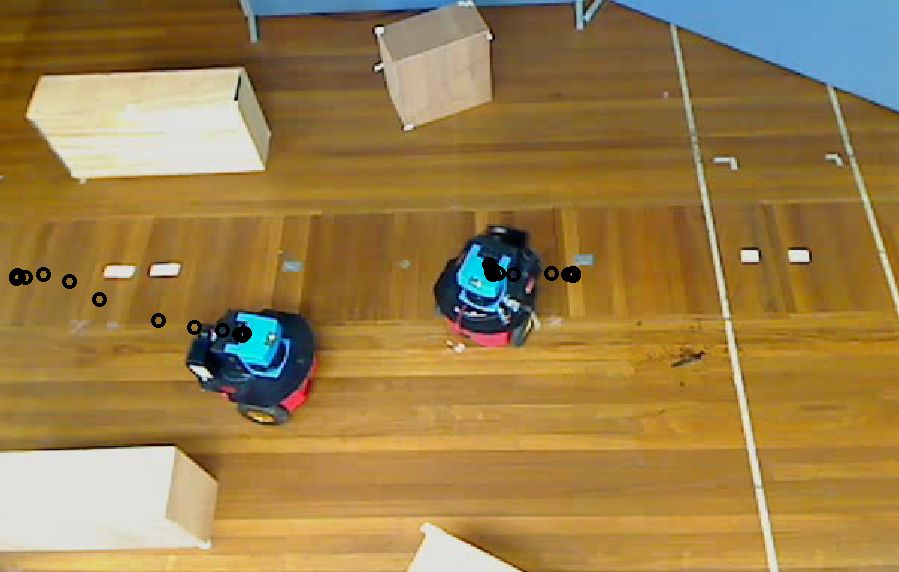}}
\subfigure[]{\includegraphics[width=0.7\mylength]{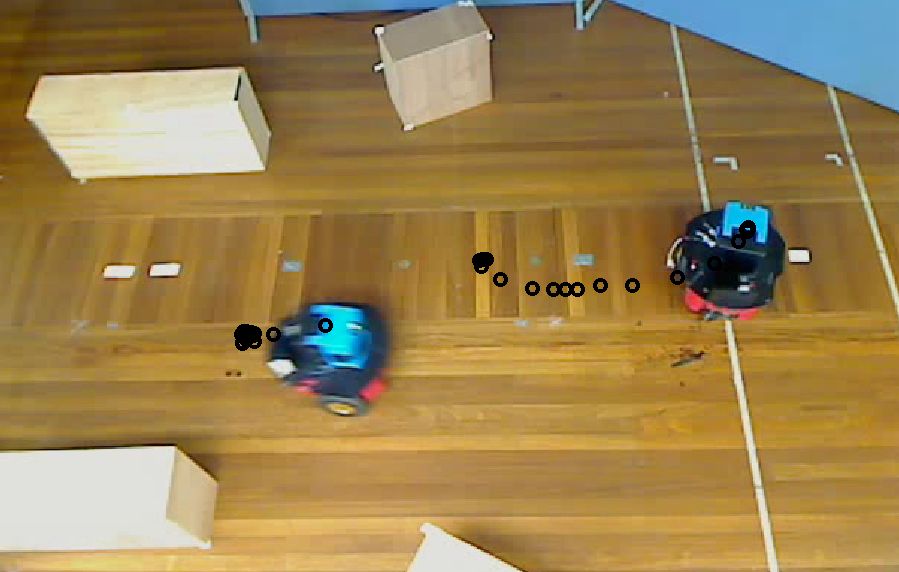}}
\subfigure[]{\includegraphics[width=0.7\mylength]{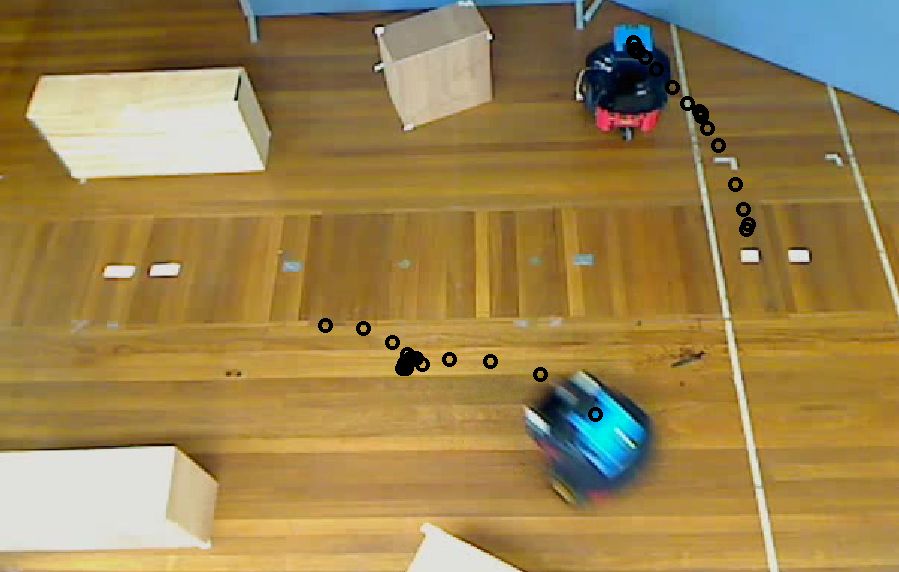}}
\subfigure[]{\includegraphics[width=0.7\mylength]{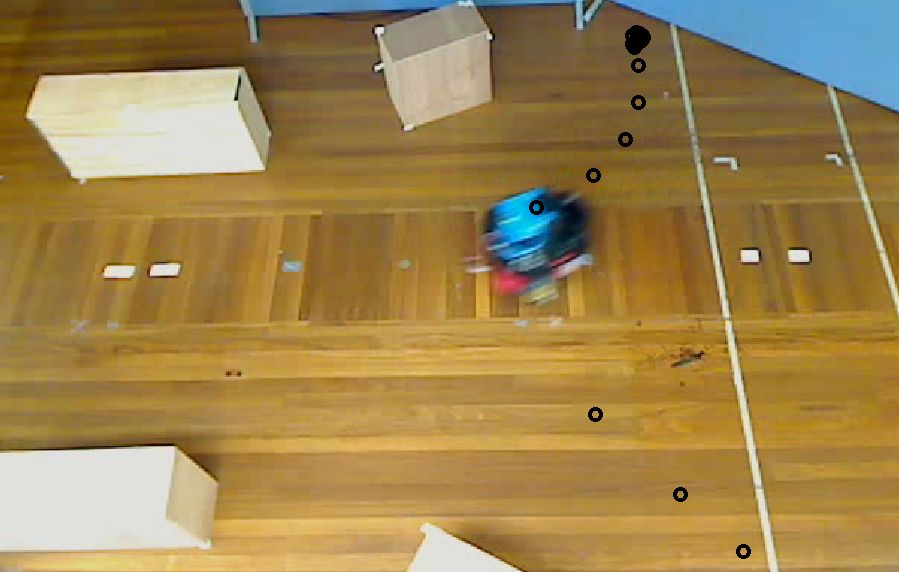}}
\caption{Sequence of images showing the experiment.}
\label{ch5:para}
\end{figure}

In the first scenario demonstrated by Figs.~\ref{ch5:para} and \ref{fig:paraplottt}, the vehicles
are both trying to move through a gap which is only wide enough for one vehicle.
The secondary vehicle to the right allows the primary vehicle to the left to pass through before
moving through itself.

\clearpage \section{Summary}
\label{ch5:con}

In this chapter a method for reactively preventing deadlocks between two mobile vehicles is proposed. A
hierarchy between the vehicles is enforced and the secondary vehicle always attempts to ensure the
primary vehicle has a convergent trajectory to follow. While simulations show this approach is
effective, it does not precipitate provable deadlock avoidance in all cases. As
such, a less efficient recovery scheme which provably works in all conditions is provided. Future
work will include more assumptions about obstacle and sensing requirements, which would
be required to prove the main navigation system alone causes the vehicles to converge to their
respective targets.

\chapter{Convergent Reactive Navigation using Minimal Information}
\label{chapt:rbf}

In this chapter, information about the obstacle is limited to just the rate 
of change in the distance to the nearest point on the obstacle. This
leads to a completely new control approach compared to the navigation
 laws proposed in Chapts.~\ref{chap:singlevehicle} to \ref{chap:dead}, where more
information was assumed about the surrounding environment.  Instead 
of path planning, a sliding mode control law is proposed which
directly controls the control inputs of the vehicle. This means
the navigation law has very low computational complexity, making it
suitable for vehicles which have low computational capability -- these are well
 exemplified by micro and miniature aerial vehicles.
\par
Another advantage of the method proposed in this chapter is that
global convergence to a target can be proven, where information about the
target is limited to its bearing from the vehicles. Many reactive
navigation methods such as artificial potential fields are not able to
show this property. It can be achieved
without randomization by using the Bug family of algorithms. Also, 
 the Pledge \cite{AbSe80} and Angulus
\cite{LuTi91} algorithms are in a similar vein. 
\par
These algorithms typically drive the robot directly towards the target when
possible; otherwise, the obstacle boundary is followed in a close
range. However, while very simple, they are not classed as reactive
since a limited amount of information must be stored between control 
updates \cite{LuTi91,ChLHKBKT05}.  In addition, a common problem with 
these strategies is that they
typically take as granted the capability of carrying out specific maneuvers
such as wall following, are disregard the kinematic constraints of the
vehicle. Extending these methods to include such considerations 
basically lies in an uncharted territory.
\par
Inspired by behaviors of animals, which are believed to use simple,
local motion control rules that result in remarkable and complex
intelligent behaviors \cite{Thomp66,CaJh99,Faj01}, a navigation strategy is proposed 
that is aimed at reaching a steady
target in a steady arbitrarily shaped maze-like environment and is
composed of the following reflex-like rules:

\begin{enumerate}[{\bf s.1)}]
\item At considerable distances from the obstacle,
  \begin{enumerate}
  \item
    \label{rbf:rule.a}
    turn towards the target as fast as possible;
  \item move towards the target when headed for it;
  \end{enumerate}
\item In close proximity of the obstacle,
  \begin{enumerate}
    \setcounter{enumii}{2}
  \item Follow (a,b) when moving away from the obstacle;
  \item Otherwise, quickly avert the collision threat by making a
    sharp turn.
  \end{enumerate}
\end{enumerate}

Studies of target pursuit in animals, ranging from dragonflies to
humans, have suggested that they often use pure pursuit method s.1) to
catch both steady and moving targets. The obstacle avoidance rule s.2)
is also inspired by biological examples such as the near-wall behavior
of a cockroach \cite{CaJh99}.
\par
To address the issues of nonholonomic constraints, control saturation,
and under-actuation, a vehicle of the Dubins car type is considered. It
travels forward with a constant speed along planar paths of bounded
curvatures and is controlled by the upper limited angular
velocity. To implement s.1), s.2), only
perceptual capabilities are needed that are enough to judge
whether the distance to the obstacle is small or not, to estimate the
sign of its time derivative, and to determine the polar angle of the
target line-of-sight in the vehicle reference frame.
\par
The convergence and performance of the proposed navigation and guidance
law are confirmed by computer simulations and real world tests with a
Pioneer P3-DX robot, equipped with a SICK LMS-200 LiDAR sensor.
\par
All proofs of mathematical statements are omitted here; they are
available in the original manuscript \cite{Matveev2011conf9, MaHoSa11ar}.
\par
The body of this chapter is organized as follows. In Sec.~\ref{rbf:S2} the problem is formally defined, and in Sec.~\ref{rbf:S3} the
main assumptions are described. The main results are outlined in Sec.~\ref{rbf:S5}. Simulations and
experiments are presented in Secs.~\ref{rbf:sec.sim} and \ref{rbf:sec.exp}. Finally, brief conclusions are given
in Sec.~\ref{rbf:S7}.

\section{Problem Statement}
\label{rbf:S2}

A planar vehicle is considered that travels forward with a constant
speed $v$, and is controlled by the angular velocity $u$ limited by a
given constant $\ov{u}$. There also is a steady point target $\bt$ and
a single steady obstacle $D \not\ni \bt$, which is an arbitrarily
shaped compact domain. Its boundary $\partial D$ is Jordan piece-wise
analytical curve without inner corners (see
Fig.~\ref{rbf:fig1}(a)). Modulo smoothened approximation of such
corners, this requirement is typically met by all obstacles considered
in robotics. The objective is to drive the vehicle to the target with
respecting a given safety margin $d(t) \geq d_{\text{safe}} >0 \;
\forall t$. The minimum distance to the obstacle is given as:

 \begin{equation*}
d(t) := \text{\bf dist}_{D}[\bldr(t)], \qquad
\dist{\bldr}:= \min_{\bldr_\ast \in D} \|\bldr_\ast - \bldr\|,
\end{equation*}

Here $\bldr(t)$ is the robot location
given by its abscissa $x(t)$ and ordinate $y(t)$ in the world
frame. The orientation of the vehicle is described by the angle
$\theta$ introduced in Fig.~\ref{rbf:fig1}(b). The kinematics of the
considered vehicles are classically described by the following
equations:

\begin{equation}
  \label{rbf:1}
  \begin{array}{l}
    \dot{x} = v \cos \theta,
    \\
    \dot{y} = v \sin \theta,
  \end{array}, \quad\dot{\theta} = u \in [-\overline{u},
  \overline{u}], \quad
  \begin{array}{l}
    \bldr(0) = \bldr_0 \not\in D
    \\
    \theta(0) = \theta_0
  \end{array} .
\end{equation}

Thus the minimal turning radius of the vehicle is equal to:

\begin{equation}
  \label{rbf:Rmin} R= v/\overline{u}.
\end{equation}

The vehicle has access to $d(t)$ and the signum $\sgn \, \dot{d}(t)$
if $d(t)\leq d_{\text{range}}$, where $ d_{\text{range}} >
d_{\text{safe}}$ is a given sensor range, and to the relative target
bearing $\beta$ (see Fig.~\ref{rbf:fig1}(b)). Mathematically, the examined strategy is described by the following equations:

\begin{figure}[ht]
  \subfigure[]{\scalebox{2.0}{\includegraphics{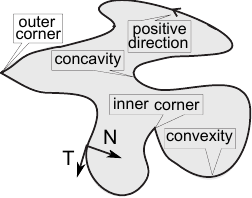}}}
  \hspace{20pt}
  \subfigure[]{\scalebox{2.0}{\includegraphics{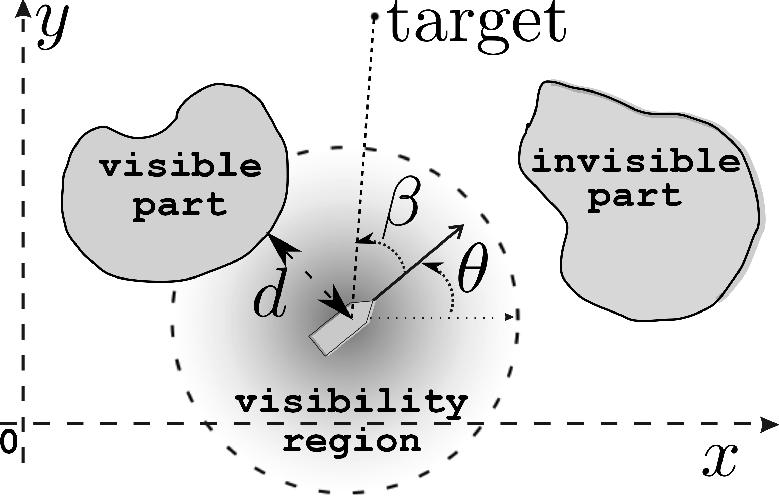}}}
  \subfigure[]{\scalebox{1.5}{\includegraphics{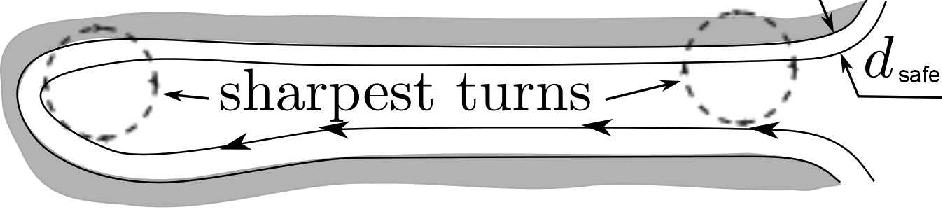}}}
  \hspace{20pt}
  \subfigure[]{\scalebox{1.5}{\includegraphics{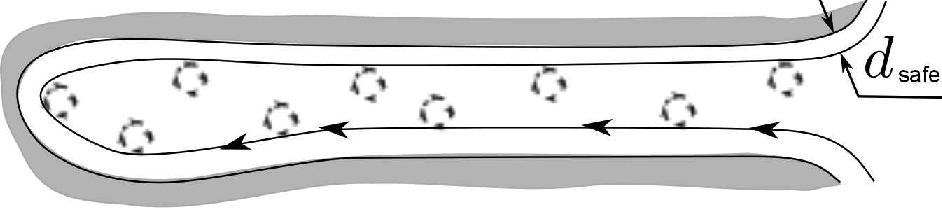}}}
  \caption{(a) Obstacle; (b) Planar vehicle; (c) Unavoidable
    collision; (d) Maneuverable enough vehicle.} \label{rbf:fig1}
\end{figure}

\begin{gather}
  \label{rbf:c.a}
  u = \overline{u} \times
  \begin{cases}
    \sgn \beta \hfill \mid & \!\!\! \text{if}\; d >\dtrig \;
    (\text{mode}\,\mathfrak{A})
    \\
    \left\{
      \begin{array}{l l}
        \sgn \beta & \text{if} \; \dot{d} > 0
        \\
        - \sigma & \text{if} \; \dot{d} \leq 0
      \end{array} \right| & \!\!\! \text{if}\; d \leq \dtrig \; (\text{mode}\,\mathfrak{B})
  \end{cases}.
\end{gather}

Here $\sigma$ and $\dtrig \in (0, d_{\text{range}})$ are controller
parameters, and $\sigma$ can assume the values ``+1'' or ``-1''. The
former gives the turn direction in (d) and the latter regulates mode
switching: $\mathfrak{A} \mapsto \mathfrak{B}$ when $d$ reduces to
$\dtrig$; $\mathfrak{B} \mapsto \mathfrak{A}$ when $d$ increases to
$\dtrig$.  When mode $\mathfrak{B}$ is activated, $\dot{d} \leq 0$; if
$\dot{d}=0$, the `turn' submode $u:=-\sigma \ov{u}$ is set up. 

In the {\it basic} version of the algorithm, $\sigma = \pm 1$ is
fixed.  As will be shown, the algorithms with $\sigma=+1$ and
$\sigma=-1$ have basically identical properties except for the
direction of bypassing the obstacle, which is counter clockwise and
clockwise, respectively.  To find a target hidden deep inside the
maze, a {\it randomized} version can be used: whenever $\mathfrak{A}
\mapsto \mathfrak{B}$, a new value of $\sigma$ is randomly generated
from a fixed Bernoulli distribution over $+1,-1$.

 \section{Main Assumptions}
\label{rbf:S3}

It is assumed that the vehicle is maneuverable enough to avoid trapping in
narrows of the maze (see Fig.~\ref{rbf:fig1}(c,d)).
Specifically, when following $\partial D$ with $ d(t) =
d_{\text{safe}}$, the full turn can always be made without violation
of the safety margin. Furthermore, it is assumed this can be done without crossing a
center of curvature of a concavity of $\partial D$, the normal radius
of the osculating circle at a distance $\leq R$ from this center, and
a location whose distance $d$ from $D$ is furnished by multiple points
in $D$. This is required since at the respective points, the
distance $d$ is uncontrollable -- even if $\dot{d}(0) = 0$, no control
can prevent convergence to $D$: $\dot{d}(t) <0\; \forall t>0, t
\approx 0$. For safety reasons, it is also assumed that $\dsafe$ exceeds
the unavoidable forward advancement $R$ during the sharpest turn.
\par
Absence of inner corners implies that for any point from $\{\bldr
\not\in D : \text{\bf dist}_D[\bldr] < d_\star \}$ with sufficiently
small $d_\star>0$, the distance $\text{\bf dist}_D[\bldr]$ is attained
at only one point $\bldr_\star \in \partial D$ and the curvature
center does not lie on the straight line segment $[\bldr_\star,\bldr]$
directed from $\bldr_\star$ to $\bldr$ \cite{Krey91}.  The {\em
  regular margin} $d_\star(D)>0$ of $D$ is the supremum of such
$d_\star$'s. Thus it follows:

\begin{equation}
  \label{rbf:uniq.rad}
  d_\star (D) \leq R_D := \inf_{\bldr \in \partial D: \varkappa(\bldr) <0} R_\kappa(\bldr)
\end{equation}

Here the variable $d_\star(D) = \infty$ if $D$ is convex.  The variable
$\varkappa(\bldr_\ast )$ is the signed curvature
($\varkappa(\bldr_\ast) < 0$ on concavities),
$R_{\varkappa}(\bldr_\ast) := |\varkappa(\bldr_\ast)|^{-1}$, and
$\inf$ over the empty set is set to be $+\infty$.  Since the distance
$d$ to $D$ may be increased by $2R$ during the full turn, the above
assumptions can be summarized by the following equation:

\begin{equation}
  \label{rbf:ass.d}
  d_\star(D) > \dsafe +2R, \quad R_D > \dsafe +3R, \quad \dsafe > R.
\end{equation}

It is also assumed that the sensor range is large enough to avoid
violation of the safety margin after detection of $D$:

\begin{equation}
  \label{rbf:ass.r}
  d_{\text{range}} > 2R + \dsafe.
\end{equation}

This takes into account that even the sharpest turn may decrease $d$
by $2R$. As $\dsafe \to R$, Eq.\eqref{rbf:ass.d} and Eq.\eqref{rbf:ass.r} shape into
$\frac{v}{\ov{u}} =R < \min \left\{ d_{\star}(D)/3;
  R_D/4;d_{\text{range}}/3 \right\}$ and mean that the robot speed $v$
should not be large to cope with the maze.
The following choice of $\dtrig$ is feasible thanks to Eq.\eqref{rbf:ass.d},
Eq.\eqref{rbf:ass.r}:

\begin{equation}
  \label{rbf:eq.trig} \dsafe + 2R< \dtrig < d_\star(D) , d_{\text{range}},
  R_D-R,
\end{equation}

 \section{Summary of Main Results}
\label{rbf:S5}

\begin{figure}[ht]
  \begin{center}
    \scalebox{2.0}{\includegraphics{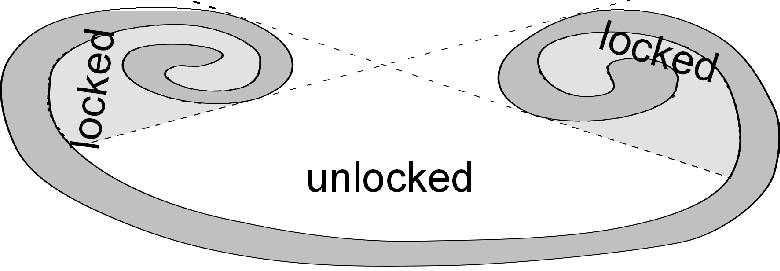}}
    \caption{Locked locations} \label{rbf:fig3}
  \end{center}
\end{figure}

The $d$-{\it neighborhood}
$\mathcal{N}(d)$ of $D$ is defined by $\mathcal{N}(d) := \{\bldr : \text{\bf
  dist}_D[\bldr] \leq d \}$.  Let $\mathscr{D} \subset \br^2$ be a
domain, $\bldr_\lozenge, \bldr_\ast \in \partial \mathscr{D}$ lie on a
common ray emitted from $\bt \not\in \mathscr{D}$, and
$(\bldr_\lozenge, \bldr_\ast)\cap \mathscr{D} = \emptyset$.  The
points $\bldr_\lozenge, \bldr_\ast$ divide $\partial \mathscr{D}$ into
two arcs. Being concatenated with $[\bldr_\lozenge, \bldr_\ast]$, each
of them gives rise to a Jordan loop encircling a bounded domain, one
of which is the other united with $\mathscr{D}$. The smaller domain is
called the {\it cave of $\mathscr{D}$ with corners} $\bldr_\lozenge,
\bldr_\ast$. The location $\bldr$ is said to be {\it locked in
  $\mathscr{D}$} if it lies in some cave of $\mathscr{D}$ (see
Fig.~\ref{rbf:fig3}).

\begin{Theorem}
  \label{rbf:th.main}
  Suppose that Eq.\eqref{rbf:ass.d}--Eq.\eqref{rbf:eq.trig} hold, both the vehicle
  initial location $\bldr_0$ and the target $\bt$ are unlocked in
  $\mathcal{N}(\dtrig)$ and

  \begin{equation}
    \label{rbf:far.enough} \dist{\bldr_0}> \dtrig + 2R,  \|\bldr_0-\bt\|>2R, \dist{\bt} >
    \dtrig.
  \end{equation}

  Then the basic control law brings the vehicle to the target in
  finite time without violation of the safety margin.
\end{Theorem}

\begin{figure}[ht]
  \begin{center}
    \subfigure[]{\scalebox{2.5}{\includegraphics{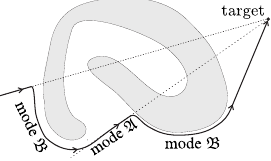}}}\hspace{20pt}
    \subfigure[]{\scalebox{2.5}{\includegraphics{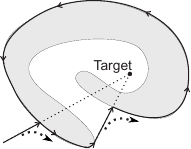}}}
    \subfigure[]{\scalebox{3.0}{\includegraphics{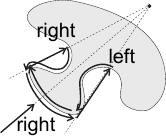}}}
    \caption{(a) Obstacle avoidance with two AM's; (b,c) Insufficiency
      of (b) only-right-turns and (c) cycle-left-and-right-turns
      options.} \label{rbf:fig.insuf}
  \end{center}
\end{figure}

Dealing with a given obstacle may comprise several AM's (i.e., motions
in mode $\mathfrak{B}$); see Fig.~\ref{rbf:fig.insuf}(a). By ii) and
iii), at most one both AM and SMEC is performed if $D$ is convex.
\par
Fig.~\ref{rbf:fig.insuf}(b) offers an example where the basic control
law fails to find a locked target. Fig.~\ref{rbf:fig.insuf}(c) shows
that repeatedly interchanging left and right turns is not enough
either. It can be shown that periodic repetition of any finite
deterministic sequence of left and right turns is not enough to find
the target in an arbitrary maze.
Randomization overcomes the insufficiency of deterministic algorithms
and aids to cope with uncertainty about the global geometry of complex
scenes.

 \section{Simulations}
\label{rbf:sec.sim}

\begin{figure}[ht]
\centering
  \subfigure[]{\scalebox{3.2}{\includegraphics{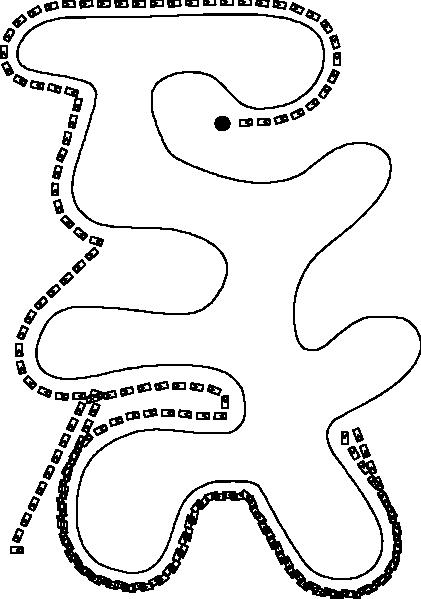}}}
  \subfigure[]{\scalebox{3.2}{\includegraphics{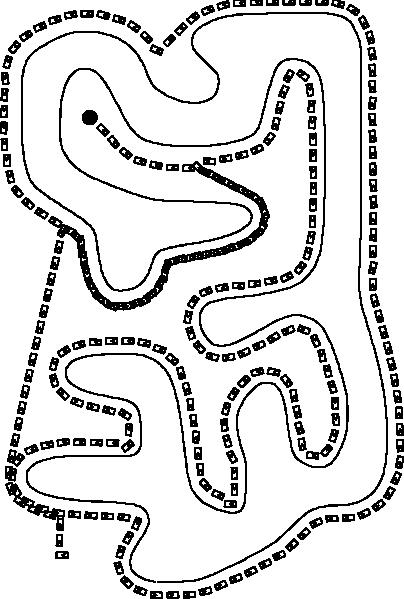}}}  
  \subfigure[]{\scalebox{1.0}{\includegraphics{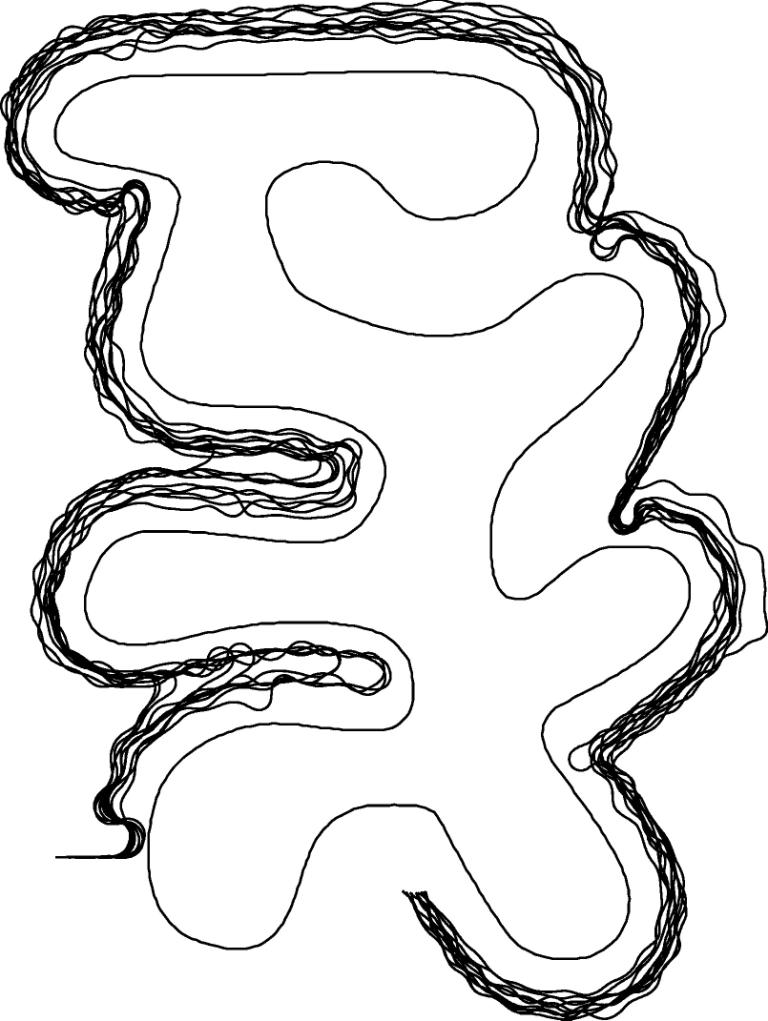}}}
  
  \caption{(a,b) Traversal to the target inside highly concave
    obstacles; (c) Performance under random noises.} \label{rbf:fig.sim1}
\end{figure}

In simulations, the control was updated every $0.02$ seconds, $\dtrig
=8m$, $\ov{u}=2.5rad/s$, $v=3 m/s$. Figs.~\ref{rbf:fig.sim1}(a,b) present typical results for the
randomized controller, where the realizations $\sigma=+,-,-,-$ and
$\sigma= -,+,+,+,-$ of the array of random turns were observed in (a)
and (b), respectively.

To address performance under noises and un-modeled dynamics, the
sensor readings and system equations were corrupted by independent
Gaussian white noises with the standard deviation $0.1m$ for $d$,
$0.1rad$ for $\beta$, $0.5 rad/s$ for $\dot{\theta}$.  To access
$\dot{d}$, the difference quotient of the noisy data $d$ was computed
with the time-step $0.1s$. The upper bound $\ov{u}$ on the turning
rate $u$ was decreased to $0.7 rad/s$, the threshold $\dtrig$ was
increased to $16m$, and extra dynamics were added via upper limiting
$\dot{u}$ by $0.7 rad/s^2$.  Fig.~\ref{rbf:fig.sim1}(c) shows that the
basic control law satisfactory guides the vehicle to the target for
various realizations of the noises and disturbances.

\clearpage \section{Experiments}
\label{rbf:sec.exp} 

The use of switching regulation often gives rise to concerns about its
practical implementation and implications of the noises and un-modeled
dynamics in practical setting, including possible chattering at
worst. To address these issues, experiments were carried out with an
Activ-Media Pioneer 3-DX mobile robot using its on-board PC and the
Advanced Robot Interface for Applications (ARIA 2.7.0), which is a C++
library providing an interface to the robot's angular and
translational velocity set-points.

\begin{figure}[ht]
\centering
  \subfigure[]{\scalebox{0.5}{\includegraphics{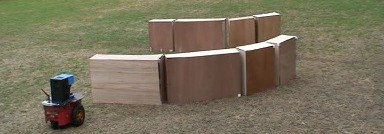}}}
  \subfigure[]{\scalebox{0.5}{\includegraphics{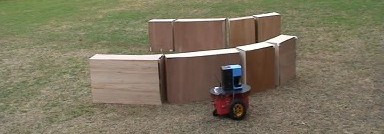}}}
  \subfigure[]{\scalebox{0.5}{\includegraphics{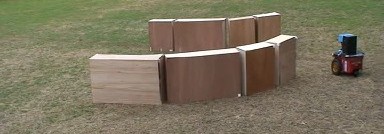}}}
  \subfigure[]{\scalebox{0.5}{\includegraphics{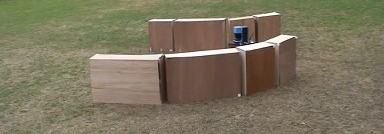}}}  
  \subfigure[]{\scalebox{0.5}{\includegraphics{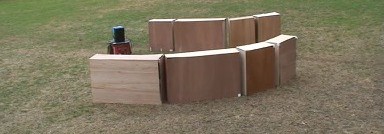}}} 
  \subfigure[]{\scalebox{0.5}{\includegraphics{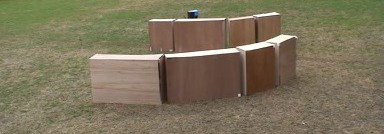}}}
  \subfigure[]{\scalebox{0.7}{\includegraphics{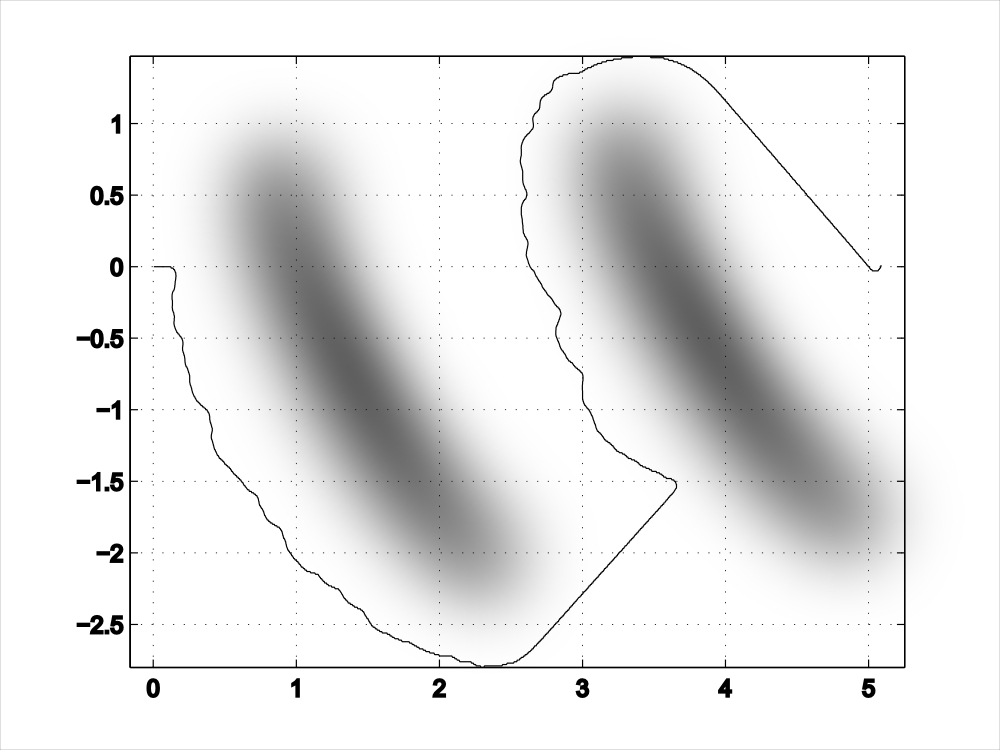}}}
  \caption{(a,b,c) Sequence of images obtained from real world
    experiments; (d) The trajectory of the robot, obtained through
    odometry.} \label{rbf:fig.test1}
\end{figure}

The position relative to the target was obtained through odometry, and
$\dtrig$ was taken to be 0.6 m and 0.9 m. The distance to the
obstacles was accessed using both LIDAR and sonar sensors.  A typical
experimental result is presented on Fig.~\ref{rbf:fig.test1}. In this
experiment, like in the others, the robot reaches the target via safe
navigation among the obstacles, with no visible mechanical chattering
of any parts being observed.

\clearpage \section{Summary}
\label{rbf:S7}

The problem of safe navigation of a Dubins-car like robot to a target
through a maze-like environment has been considered and a new
bio-inspired control law has been proposed.
 The algorithm alternates between the pure pursuit navigation to
 the target and the sharpest turns away from the obstacle, with the
latter being activated only at a short distance from obstacle and
only if this distance decreases.
Via computer simulations and
experiments with real robots, it has been shown that this law
constitutes an effective target-seeking strategy that can be realized
even in complex maze-like environments.

\chapter{Convergent Reactive Navigation using Tangent Tracking}
\label{chapt:tf}

In previous chapters it was assumed that a finite sensor range was
available to the vehicle. In this chapter it is assumed that instead of
range information, the vehicle has knowledge of the relative heading of
visible obstacle edges surrounding the vehicle. This is more
similar to what would be obtained from a vision sensor.
\par
The advantage of using a tangent sensor is that under certain
conditions, the vehicle is able to follow the minimum length path.
In this chapter, firstly the problem of global shortest path
planning to a steady target in a known environment is considered. It is assumed that 
the environment contains a number of
potentially non-convex obstacles, which have constraints imposed on their maximum
curvature. 
\par
Secondly, a randomized navigation
algorithm is propoposed, for which it can be shown that the robot will reach the target with
probability $1$, while avoiding collision with the obstacles.  The 
robot studied in this chapter is a unicycle like vehicle, described
by the standard nonholonomic model with a hard constraint on the
angular velocity. 
\par
Unlike many other papers on this area of robotics which present
heuristic based navigation strategies, mathematically rigorous
analysis is available for the navigation algorithm proposed in this chapter. 
Moreover, it should be pointed out that many other
papers on this topic (see e.g. \cite{S00,
  Kamon1997journ6, Kamon1998journ9, LA92}) do not assume non holonomic
constraints on robot's motion, which is a severe limitation in
practice.  The performance of this real-time navigation
strategy is confirmed with extensive computer simulations and outdoor experiments
with a Pioneer P3-DX mobile wheeled robot.  
\par
All proofs of mathematical statements are omitted here; they are
available in the original manuscript \cite{Savkin2012journ1}.
\par
The body of this chapter is organized as follows. In Sec.~\ref{tf:S2} the problem is formally defined, and in Sec.~\ref{tf:S3} the optimal
off-line paths are described. A on-line navigation method is presented in Sec.~\ref{tf:S4}. Simulations and
experiments are presented in Secs.~\ref{tf:S5} and \ref{tf:S6}. Finally, brief conclusions are given
in Sec.~\ref{tf:S7}.

\section{Problem Statement}
\label{tf:S2}

A Dubins-type vehicle travelling at constant speed in the
plane is considered, which has its maximum angular velocity limited by
a given constant. The model of the vehicle is given as follows:

\begin{equation}
  \label{tf:1}
  \begin{array}{l}
    \dot{x} = v \cos \theta
    \\
    \dot{y} = v \sin \theta
    \\
    \dot{\theta} = u \in [-u_M, u_M]
  \end{array}, \qquad
  \begin{array}{l}
    x(0) = x_0
    \\
    y(0) = y_0
    \\
    \theta(0) = \theta_0.
  \end{array} 
\end{equation}

Here $(x,y)$ is the vector of the vehicle's Cartesian coordinates,
$\theta$ gives its orientation, $v$ and $u$ are the speed and angular
velocity, respectively. The maximal angular velocity $u_M$ is given.
The robot satisfies the standard non-holonomic constraint:

  \begin{equation}
    \label{tf:nonh}
    |u(t)|\leq u_M
  \end{equation} 
  
Also, the minimum turning radius of
the robot is given by:

\begin{equation}
  \label{tf:Rmin} R_{min}= \frac{v}{u_M}.
\end{equation}

Any path $(x(t),y(t))$ of the robot (\ref{tf:1}) is a plane curve
satisfying the following constraint on its so-called average curvature
(see \cite{DUB57}). Let $P(s)$ be this path parametrized by arc
length; it is assumed the following contraint applies:

\begin{equation}
  \label{tf:ac}
  \|P'(s_1)-P'(s_2)\|\leq \frac{1}{R_{min}}|s_1-s_2|.
\end{equation}

Notice
that the constraint (\ref{tf:ac}) on average curvature is used because
the standard definition of curvature from differential
geometry is not suitable (see e.g. \cite{POG}), since the curvature may not exist at
some points of the robot's trajectory.

There is a steady point-wise target ${T}$ and several disjoint
obstacles $D_1,\ldots,D_k$ in the plane. Let the safety margin $d_0>0$
be given.  The objective is to drive the vehicle to the target through
the obstacle-free part of the plane while keeping the safety margin.

Let $D$ be a closed set, $p$ be a point in the plain.  Introduce the
distance $\rho(D,p)$ given by:

\[
\rho(D,p):=\min_{q\in D}\|p-q\|.
\]

is set to zero if $p\in D$.
\par
The following definition defines the neighborhood set of a point:

\begin{Definition}
  For $d_0>0$, the $d_0-${\em neighborhood} of the domain $D \subset
  {\bf R}^2$ is the set formed by all points at the distance $\leq
  d_0$ from $D$, i.e., $ {\cal N}[D,d_0]$ $:=$ $\{p \in {\bf R}^2:
    \rho(D,p)\leq$ $d_0 \}$
\end{Definition}

The correct behavior of the system is defined as follows:

\begin{Definition}
  A path $p(t)=(x(t),y(t))$ of the robot (\ref{tf:1}) is said to be
  {\em target reaching with obstacle avoidance} if there exists a time
  $t_f>0$ such that $p(t_f) ={T}$ and $p(t)$ does not belong to
  ${\cal N}[D_i,d_0]$ for all $i , t\in [0,t_f]$.
\end{Definition}

There are several assumptions which can now be made about the obstacle set:

\begin{Assumption}
  \label{tf:Aso} \label{tf:Aso1} \label{tf:Aso2a} \label{tf:Aso2} For all $i$ the set ${\cal N}[D_i,d_0]$ is a closed,
  bounded, connected and linearly connected set. The sets ${\cal N}[D_i,d_0]$ and ${\cal N}[D_j,d_0]$
  do not overlap for any $i\neq j$. For all $i$ the boundary $\partial D_i$ of the
  obstacle $D_i$ is a closed, non-self-intersecting analytic curve.  For all $i$ the boundary $\partial D_i(d_0)$ of the
  set ${\cal N}[D_i,d_0]$ is a closed, non-self-intersecting analytic
  curve with curvature $k_i(p)$ at any point $p$ satisfying
  $|k_i(p)|\leq \frac{1}{R_{min}}$.
\end{Assumption}

\begin{figure}[ht]
  \centering
  \includegraphics[width=0.5\columnwidth]{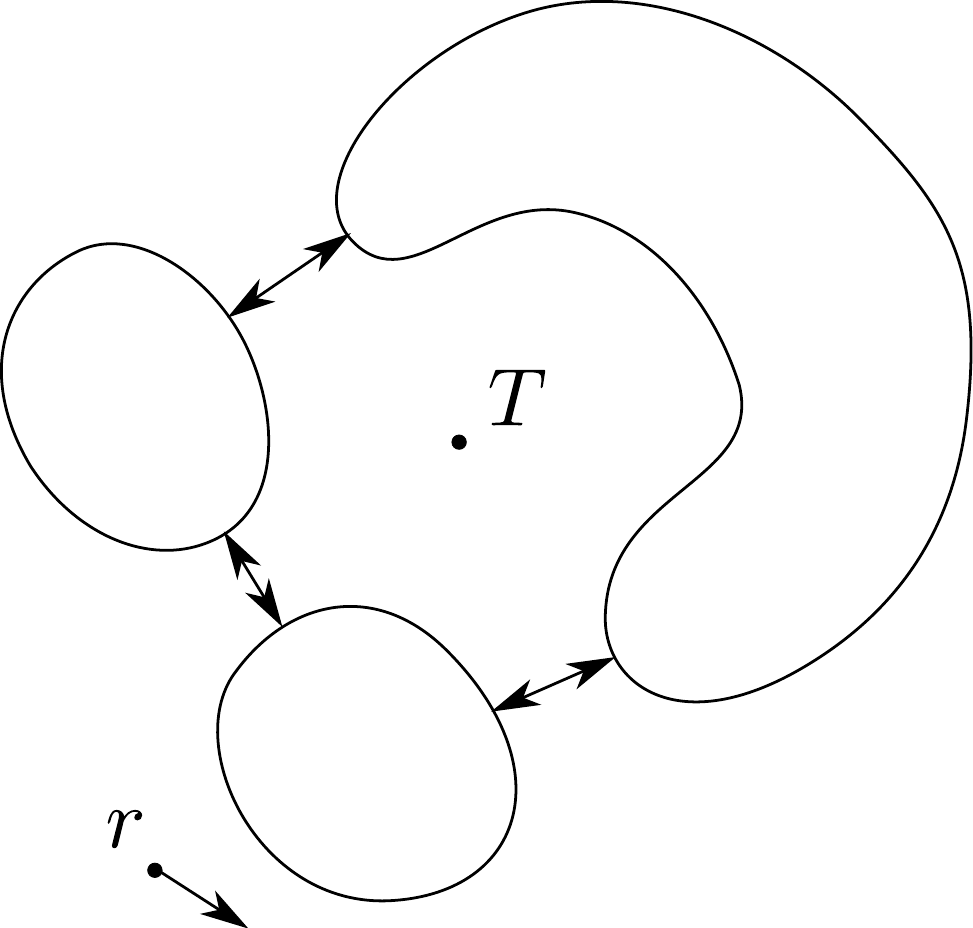}
  \caption{Conditions for reaching the target.}
  \label{tf:fig2}
\end{figure}

  It is obvious that if Assumption \ref{tf:Aso1} does not hold, target
  reaching with obstacle avoidance may be impossible (see
  e.g. Fig. \ref{tf:fig2}).

It is also assumes that the initial position $p(0)=(x(0),y(0))$ of the
robot is far enough from the obstacles and the target:

\begin{Assumption}
  \label{tf:Aso3} \label{tf:Aso4} The following inequalities hold: $\rho({\cal
    T},p(0))$ $\geq$ $8R_{min}$ and $\rho({\cal N}[D_i,$ $d_0],p(0))$ $\geq$ 
  $8R_{min}$ for all $i$. The robot initial heading $\theta(0)$ is not tangent
  to any boundary $\partial D_i(d_0)$.
\end{Assumption}

  There are two circles with the radius $R_{min}$ that cross the
  initial robot position $p(0)$ and tangent to the robot initial
  heading $\theta(0)$. These are labeled initial circles.

 \section{Off-Line Shortest Path Planning}
\label{tf:S3}
In this section, the shortest or minimum length target
reaching paths with obstacle avoidance is described.

\begin{Definition}
  A straight line $L$ is said to be a tangent line if one of the
  following conditions holds:
  \begin{enumerate}
  \item The line $L$ is simultaneously tangent to two boundaries
    $\partial D_i(d_0)$ and $\partial D_j(d_0)$ where $i\neq j$.
  \item The line $L$ is tangent to a boundary $\partial D_i(d_0)$ on
    two different points.
  \item The line $L$ is simultaneously tangent to a boundary $\partial
    D_i(d_0)$ and an initial circle.
  \item The line $L$ is tangent to a boundary $\partial D_i(d_0)$ and
    crosses the target ${T}$.
  \item The line $L$ is tangent to an initial circle and crosses the
    target ${T}$.
  \end{enumerate}
  Points of boundaries $\partial D_i(d_0)$ and initial circles
  belonging to tangent lines are called tangent points. 
  Consider only finite segments of tangent lines such that their
  interiors do not overlap with any boundary $\partial D_i(d_0)$. The
  segments of tangent lines of types $1)$ and $2)$ are called
  $(OO)-$segments, the segments of tangent lines of type $3)$ are
  called $(CO)-$segments, the segments of tangent lines of type $4)$
  are called $(OT)-$segments, and the segments of tangent lines of
  type $5)$ are called $(CT)-$segments. Furthermore, only
  segments that do not intersect the interiors of the initial circles are considered.
\end{Definition}

\begin{figure}[ht]
  \centering
  \includegraphics[width=0.4\columnwidth]{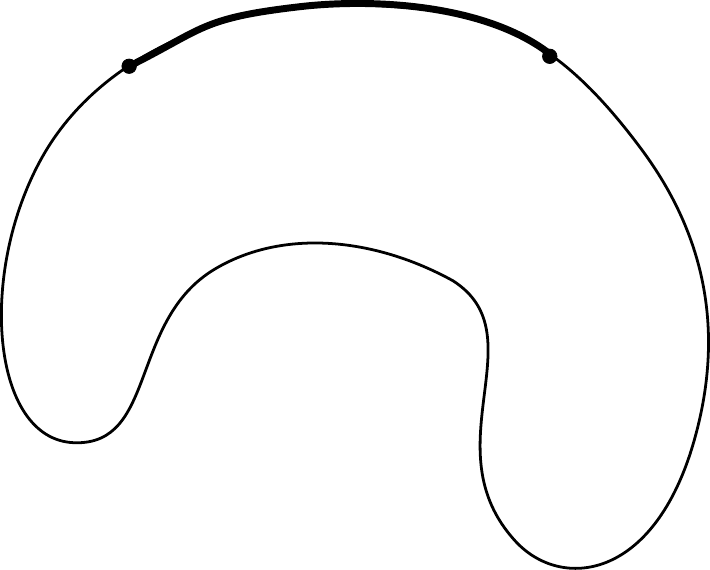}
  \caption{A $(B)-$segment.}
  \label{tf:fig3}
\end{figure}

\begin{figure}[ht]
  \centering
  \includegraphics[width=0.3\columnwidth]{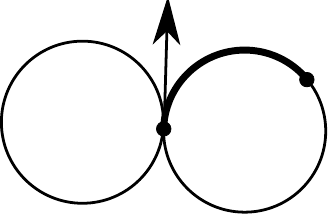}
  \caption{A $(C)-$segment.}
  \label{tf:fig4}
\end{figure}

\begin{Definition}
  A segment of a boundary $\partial D_i(d_0)$ between two tangent
  points is called $(B)-$segment if curvature is non-negative at any
  point of this segment (see Fig. \ref{tf:fig3}).  A segment of an
  initial circle between the initial robot position $p(0)$ and a
  tangent point is called $(C)-$segment (see Fig. \ref{tf:fig4}).
\end{Definition}

Now the main result of this section may be presented:

\begin{Theorem}
  \label{tf:T1} Suppose that Assumptions~{\rm \ref{tf:Aso} --
    \ref{tf:Aso4}} hold.  Then there exists a shortest (minimum
  length) target reaching path with obstacle avoidance.  Furthermore,
  a shortest target reaching path consists of $n\geq 2$ segments
  $S_1,S_2,\ldots,S_n$ such that if $n=2$ then $S_1$ is a
  $(C)-$segment and $S_2$ is a $(CT)-$segment.  If $n\geq 3$ then
  $S_1$ is a $(C)-$segment, $S_2$ is a $(CO)-$segment and $S_n$ is a
  $(OT)-$segment. If $n>3$ and $3\leq k\leq n-1$ then any $S_k$ is
  either $(OO)-$segment or $(B)-$segment.
\end{Theorem}

\begin{figure}[ht]
  \centering
  \includegraphics[width=0.5\columnwidth]{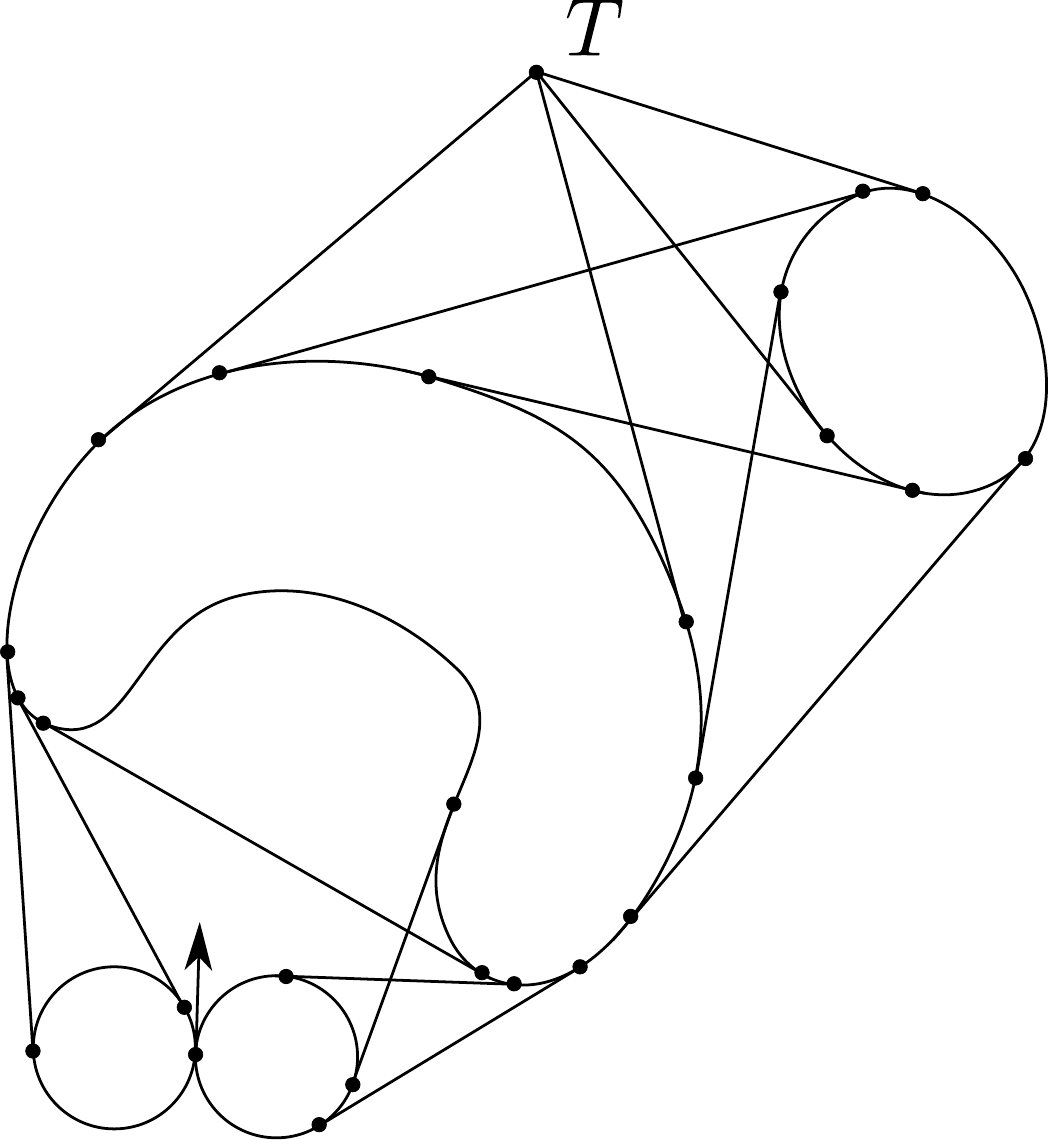}
  \caption{An example of the tangent graph, with a target point $T$.}
  \label{tf:fig5}
\end{figure}

 \section{On-Line Navigation}
\label{tf:S4}

\begin{figure}[ht]
  \centering
  \includegraphics[width=0.5\columnwidth]{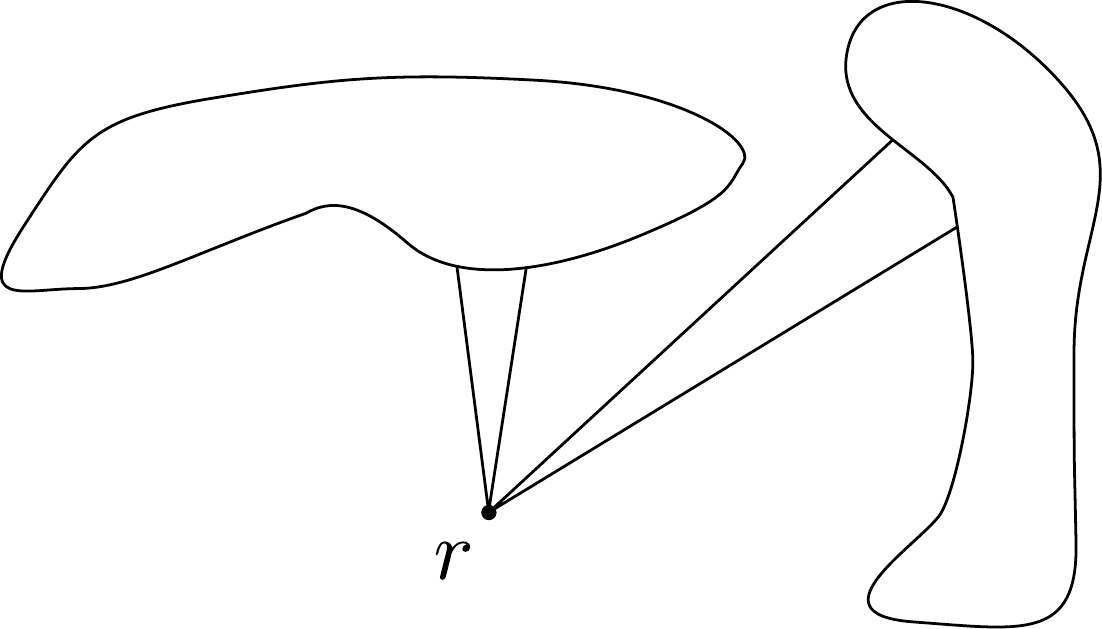}
  \caption{Vision sensor model of the robot.}
  \label{tf:fig6}
\end{figure}

In this section, the case of sensor based on-line navigation is considered, where
the robot does not know the
location of the target and the obstacles a priori. The robot is
equipped with a vision type sensor which is able to determine
coordinates of the target and points of the boundaries $\partial
D_i(d_0)$ if the straight line segment connecting the robot current
coordinates and the point of interest does not intersect any obstacle
$D_i$ (see Fig. \ref{tf:fig6}).

\begin{figure}[ht]
  \centering
  \subfigure[]{\includegraphics[width=0.3\columnwidth]{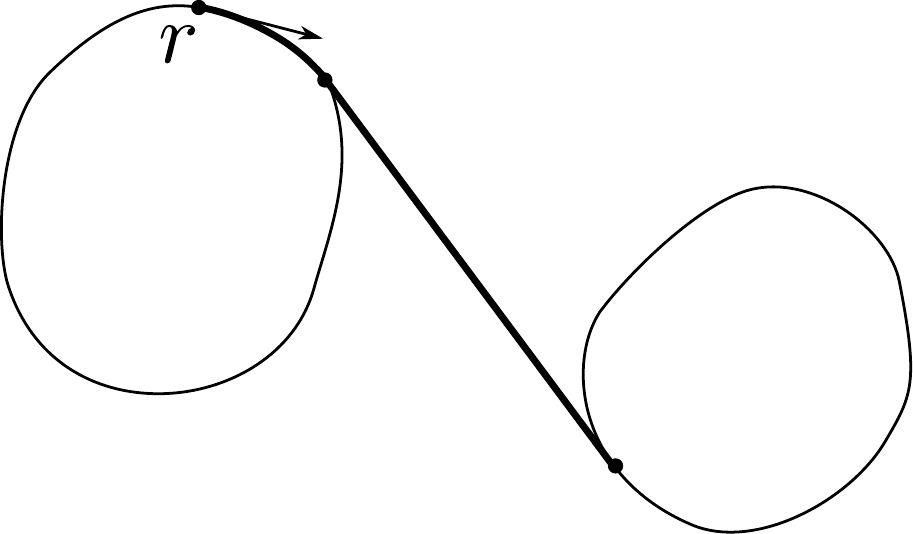}}
  \hspace{20pt}
  \subfigure[]{\includegraphics[width=0.25\columnwidth]{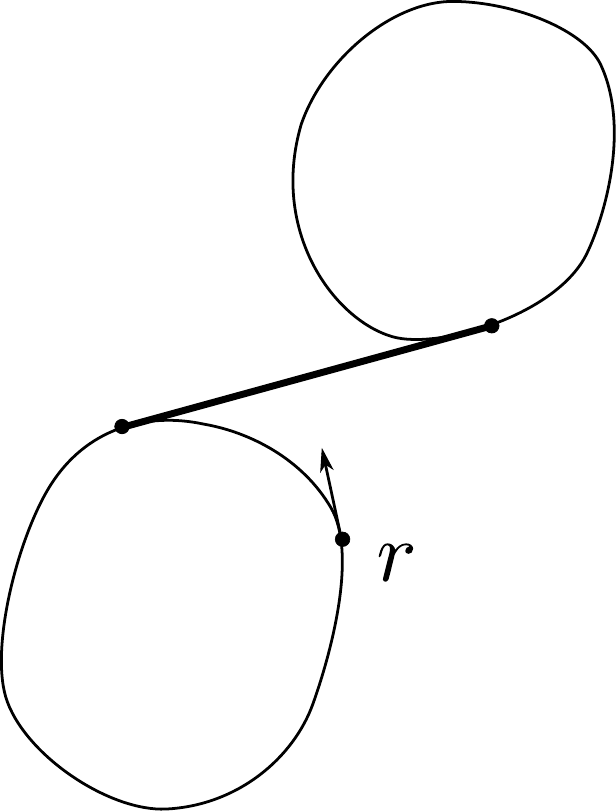}}
  \caption{Conditions for a tangent segment to be transversed.}
  \label{tf:fig12}
\end{figure}

\begin{Definition}
  When the robot moves along the boundary of an obstacle or an initial
  circle and reaches a tangent point, it can leave the boundary and
  move along the corresponding straight line edge of the tangent graph
  if its heading is equal to direction of this edge (see
  Fig. \ref{tf:fig12}a). In this case, the robot is said to have reached
  an exit tangent point.  (A case when the robot cannot leave the
  boundary at a tangent point is shown on Fig. \ref{tf:fig12}b).
  Furthermore, if the straight line edge corresponding to this exit
  tangent point is an $(OT)$ or $(CT)$ segment, then this
  point is called an exit tangent point of $T-$type. Otherwise, if the straight
  line edge corresponding to this exit tangent point is an $(OO)$ or
  $(CO)$ segment, this point is called an exit tangent point of
  $O-$type.
\end{Definition}

Let $0<p<1$ be a given number. To solve this problem, the
following probabilistic navigation algorithm is proposed:

{\bf A1:} The robot starts to move along any of two initial circles.

{\bf A2:} When the robot moving along an obstacle boundary or an
initial circles reaches an exit tangent point of $T-$type, it starts
to move along the corresponding $(OT)$ or $(CT)$ segment.

{\bf A3:} When the robot moves along an obstacle boundary or an
initial circle and reaches an exit tangent point of $O-$type, with
probability $p$ it starts to move along the corresponding $(CO)$ or
$(OO)$ segment, and with probability $(1-p)$ it continues to move
along the boundary.

{\bf A4:} When the robot moves along a $(CO)$ or $(OO)$ segment and
reaches a tangent point on an obstacle boundary, it starts to move
along the obstacle boundary.

Now the main result of this section may be presented:

\begin{Theorem}
  \label{tf:T2}
  Suppose that Assumptions~{\rm \ref{tf:Aso} -- \ref{tf:Aso4}} hold.  Then for any $0<p<1$, the algorithm {\bf
    A1--A4} with probability $1$ defines a target reaching path with
  obstacle avoidance.
\end{Theorem}

 \section{Simulations}
\label{tf:S5}

In this section, computer simulations of the reactive
navigation algorithm {\bf A1--A4} are presented.  This navigation strategy was
realized as a sliding mode control law, by switching between a
boundary following approach proposed in \cite{Matveev2011journ2}, and
the pure pursuit navigation approach (see e.g. \cite{ST10}). The
navigation law can be expressed as follows:

\begin{equation}
  \label{tf:eq:control}
  u(t) = \left \{  \begin{array}[h]{cc}\pm u_M & R1 \\\Gamma sgn\left[ \phi _{tan}(t) \right]u_M & R2 \\ \Gamma sgn\left[\dot{d}_{min}(t) + X(d_{min}(t)-d_{0} )\right]u_M & R3 \\ \end{array} \right.\\
\end{equation}

\begin{eqnarray}
  \label{tf:eqn:trans}
  R1 \rightarrow R2 &:& CO \text{ or } CT \text{ detected} \nonumber \\
  R2 \rightarrow R3 &:& d_{min}(t) < d_{trig}, \dot{d}_{min}(t) < 0 \nonumber\\
  R3 \rightarrow R2 &:& \left \{ \begin{array}{l}OT \text{ detected} \\OO \text{ detected, probability } p\end{array}\right.  \nonumber\\
\end{eqnarray}

 This navigation law is a rule for switching between three separate
  modes $R1-R3$.  Initially mode $R1$ is active, and transitions to
  other modes are determined by equation (\ref{tf:eqn:trans}).  Mode
  $R1$ describes motion along the initial circle with maximal
  actuation. Mode $R2$ describes pursuit navigation, where
  $\phi_{tan}(t)$ is defined as the angle between the vehicle's
  heading and a line segment connecting the vehicle and currently
  tracked tangent edge (see equation (\ref{tf:eq:tint})). Mode $R3$
  describes boundary following behaviour, where the control
  calculation is based on the minimum distance to the nearest
  obstacle, defined as $d_{min}(t)$. This control law is subject to
  some restrictions which are inherited from \cite{Matveev2011journ2},
  however these are satisfied due to the assumptions in Sec.~\ref{tf:S2}.  
  The variable $\Gamma$ is defined as $+1$ if the
  obstacle targeted is on the left of the tangent being tracked, $-1$
  if it is on the right. A constant $d_{trig} > d_0$ is also
  introduced which determines when the control system transitions to
  boundary following mode \cite{Matveev2011journ2}. The saturation
function $X$ is defined as follows:

\begin{equation}
  \label{tf:eq:sat}
  X(r) = \left \{  \begin{array}[h]{cc} lr & |r| < k\\ lk\ sgn(r) & otherwise\\\end{array} \right.
\end{equation}

Here $l$ and $k$ are tunable constants. Because of any potential
chattering in the output, once it is decided to not pursue a tangent,
there is a short pause until tangent following can potentially be
engaged again.

  \begin{Remark}
    Notice that it follows from (\ref{tf:1}) that the robot's
    trajectory is differentiable. This and Assumption \ref{tf:Aso2a}
    imply that the function $d_{min}(t)$ is continuous, however, in
    the case on a non-convex obstacle, the function $d_{min}(t)$ may
    be non-differentiable for some $t$. It does not really matter in
    practice, in these computer simulations and experiments with a real
    robot, the first order difference approximation of
    $\dot{d}_{min}(t)$ in the equation (\ref{tf:eq:control}) is used.
  \end{Remark}

  In these simulations and experiments, the robot is assumed to posess a LiDAR-type
  device, which informs the vehicle of the distance from the vehicle to
  obstacle in a finite number of directions around the obstacle. The
  length of these detection rays are defined by $\eta_i$, where
  $\eta_0$ refers to a detection ray directly in front of the
  vehicle. The angular spacing between successive rays is defined
  $\Delta \theta$. In this case a tangent in front of the vehicle is
  detected by monitoring $|\Delta \eta_0|$ for any changes beyond some
  threshold $d_{thresh}$. The transversal direction can be determined
  by comparing the immediately adjacent obstacle detections; $\Gamma
  := \sgn(\eta_1 - \eta_{-1})$.

  To calculate the error parameter $\phi_{tan}$, a point $T_{int}$ is
  defined as an intermediate target. This is calculated by iteration
  of the following steps:

  {\bf 1)} Initially when $R2$ is engaged, the $T_{int}$ is set to be
  a constant offset from the detected tangent point:

  \begin{equation}
    \label{tf:eq:tint}
    T_{int} = \left [ \begin{array}{c} x(t) + cos(\theta(t) + i\Delta\theta + \Gamma \tan^{-1} (d_{0}/d_{cen}(t)))\\y(t) + sin(\theta(t) + i\Delta\theta + \Gamma \tan^{-1} (d_{tar}/d_{cen}(t))) \end{array} \right ]
  \end{equation}
  with $i := 0$.

  {\bf 2)} $\phi_{tan}$ is calculated by finding the angle between the
  vehicle heading and a line connecting the vehicle and $T_{int}$.

  {\bf 3)} In subsequent time steps a successive tangent
  point is found; a search occurs for the appropriate $i$ so that $\eta_i -
  \eta_{i + \Gamma} > d_{thresh}$ (so the same transversal
  direction as the stored intermediate target is maintained), and also has the
  smallest Euclidean distance to $T_{int}$. The point is calculated
  using equation (\ref{tf:eq:tint}).

\begin{Remark}
  While this calculation calls for an estimate of position to be
  available to the robot, this estimate only needs to be accurate for
  a relatively short time between control updates. Thus in the studied
  case, robot odometry is sufficient since robot odometry gives an
  accurate estimate over short time intervals.
\end{Remark}

\begin{table}[ht]
  \centering {
    \begin{tabular}{| l | c |}
      \hline
      $u_{max}$ & $1.3 rads^{-1}$ \\
      \hline
      $v$ & $1.5 ms^{-1}$ \\
      \hline
      $d_{tar}$ & $5 m$  \\
      \hline
   \end{tabular} 
\hspace{10pt}
    \begin{tabular}{| l | c |}
      \hline
      $d_{trig}$ & $10 m$ \\
      \hline
      $l$ & $0.33$  \\
      \hline
      $k$ & $15m$  \\
      \hline
   \end{tabular} 
\hspace{10pt}
    \begin{tabular}{| l | c |}
      \hline
      $p$ & $0.7$  \\
      \hline
      $d_{thresh}$ & $10 m$  \\
      \hline
    \end{tabular} 
  }
  \caption{Simulation parameters for tangent-following controller.}
  \label{tf:simpar}
\end{table}

\begin{figure}[ht]
  \centering
  \includegraphics[width=0.5\columnwidth]{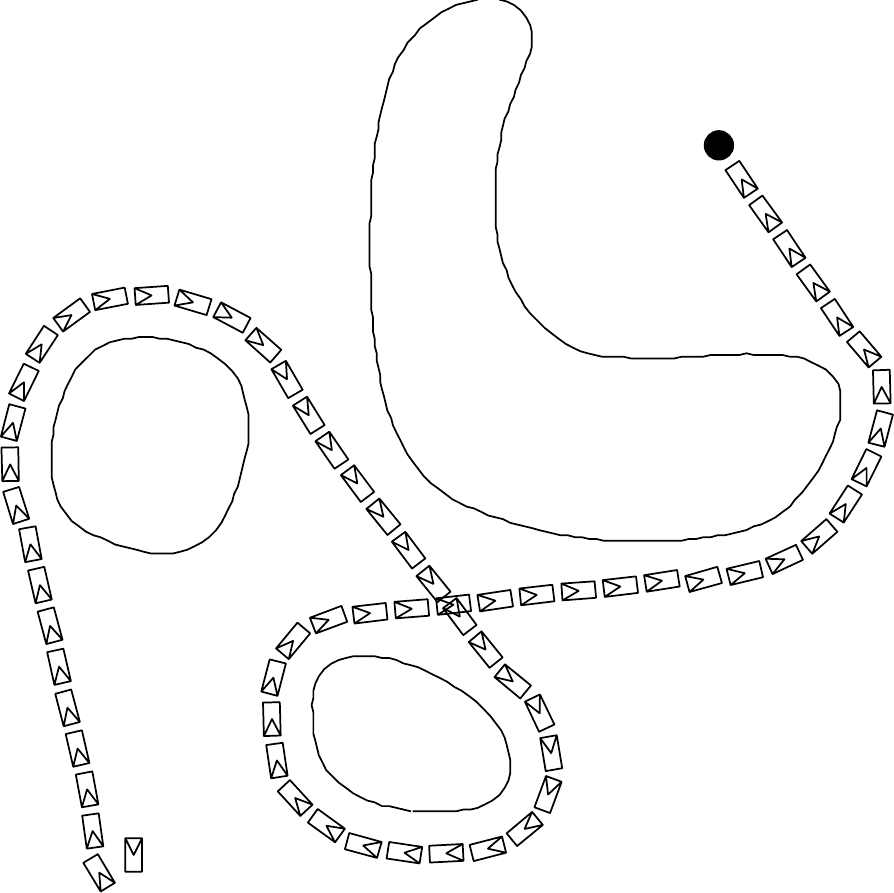}
  \caption{Simulation with a simple environment.}
  \label{tf:fig:sim1}
\end{figure}

\begin{figure}[ht]
  \centering
  \includegraphics[width=0.5\columnwidth]{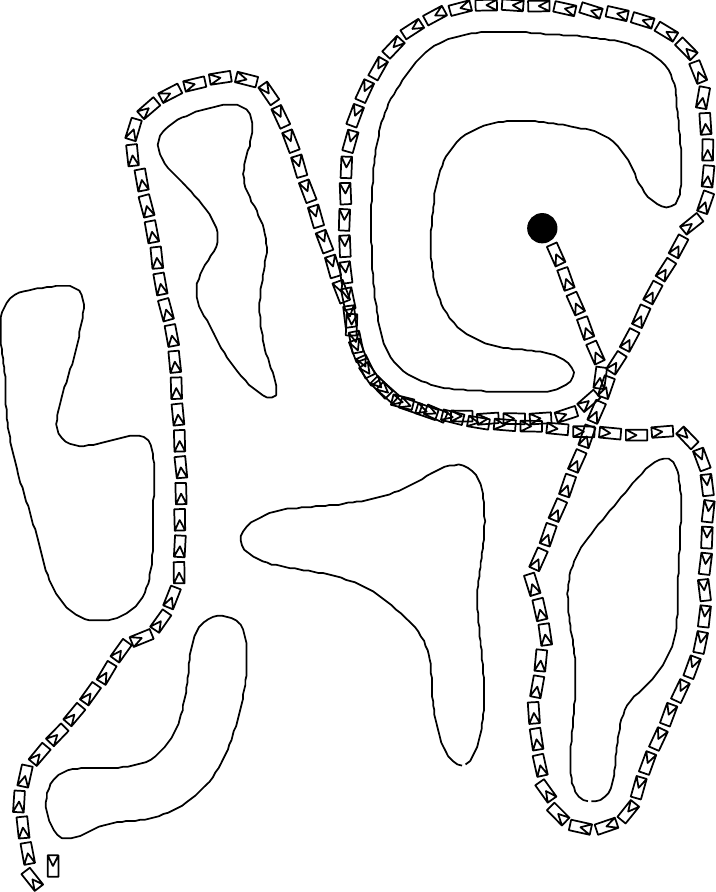}
  \caption{Simulation with a challenging environment.}
  \label{tf:fig:sim2}
\end{figure}

\begin{figure}[ht]
  \centering
  \includegraphics[width=0.5\columnwidth]{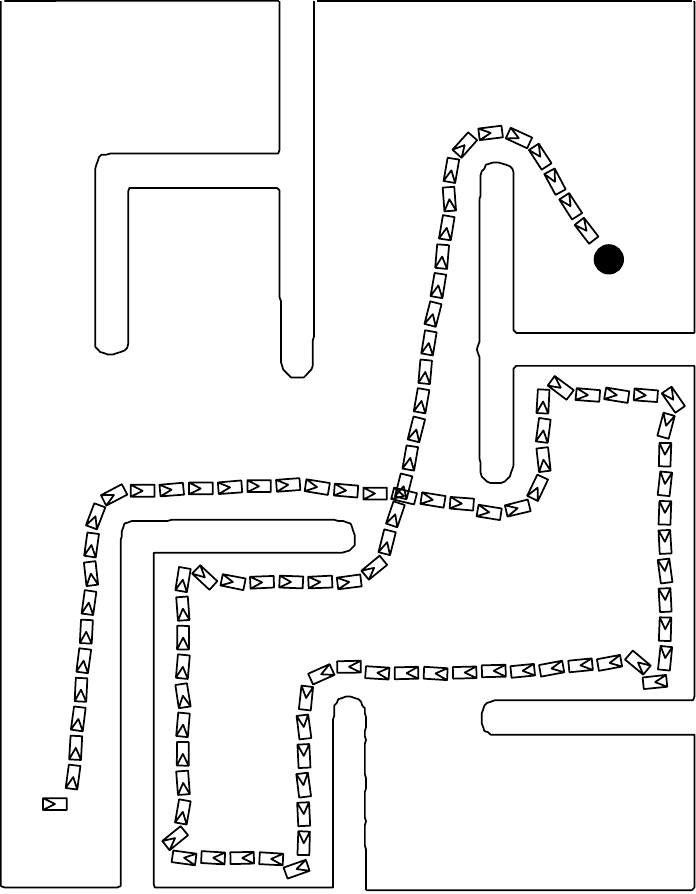}
  \caption{Simulation with a more challenging environment.}
  \label{tf:fig:sim3}
\end{figure}

The following simulation was carried out with the unicycle model
(\ref{tf:1}) and the navigation law updated at 10 Hz. Parameters 
used for simulations may be found in  Table~\ref{tf:simpar}. In Figs
\ref{tf:fig:sim1}, \ref{tf:fig:sim2} and \ref{tf:fig:sim3}, it can be
seen that the robot converges to the target without any problems, as
expected. Different sequences of random numbers would of course lead
to different paths around the obstacles.

\clearpage \section{Experiments}
\label{tf:S6}

In this section, an implementation of the reactive
navigation algorithm {\bf A1--A4} is presented.  Experiments were carried out with
a Pioneer P3-DX mobile robot. A LiDAR device with an angular
resolution of $1^o$ was used to detect tangents directly ahead of
the robot as well as obstacles in the vicinity of the robot. In the
scenario tested the vehicle was not provided with a target, rather
it was allowed to patrol the area indefinitely. Odometry information
available from the robot was used over a single time step to
compensate for the movement of the previously calculated tangent point
relative to the vehicle, as explained in the previous section (see
Remark 5.1). Range readings over a maximum threshold $R_{max}$
were truncated, in order to prevent any object outside the test area
from influencing the results. Parameters used for experiments are shown in Table~\ref{tf:exppar}.

Measures were also included to reduce control chattering, which can
cause detrimental effects when using real robots
\cite{Teimoori2010journ1a}. The signum function of equation
(\ref{tf:eq:control}) was replaced by a saturation function - equation
(\ref{tf:eq:sat}) with $l = k = 1$), and the standard low level
controller on the robot was modified to filter high frequency control
inputs.

\begin{table}[ht]
  \centering {
    \begin{tabular}{| l | c |}
      \hline
      $v$ & $0.25 ms^{-1}$ \\
      \hline
      $d_{tar}$ & $1 m$  \\
      \hline
      $d_{trig}$ & $1.5 m$ \\
      \hline
   \end{tabular} 
\hspace{10pt}
    \begin{tabular}{| l | c |}
      \hline
      $l$ & $0.1$  \\
      \hline
      $k$ & $1.0m$  \\
      \hline
      $p$ & $0.7$  \\
      \hline
   \end{tabular} 
\hspace{10pt}
    \begin{tabular}{| l | c |}
      \hline
      $R_{max}$ & $6 m$  \\
      \hline
      $d_{thresh}$ & $1 m$  \\
      \hline
    \end{tabular} 
  }
  \caption{Experimental parameters for tangent-following controller.}
  \label{tf:exppar}
\end{table}

\begin{figure}[ht]
  \centering
    \subfigure[]{\includegraphics[width=0.45\columnwidth]{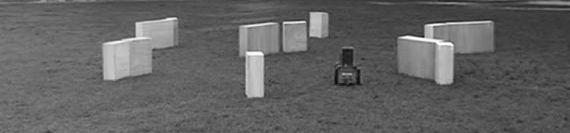}}
    \subfigure[]{\includegraphics[width=0.45\columnwidth]{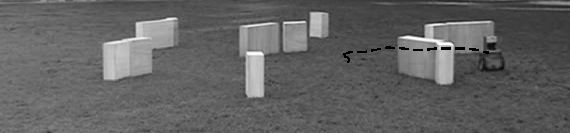}}
    \subfigure[]{\includegraphics[width=0.45\columnwidth]{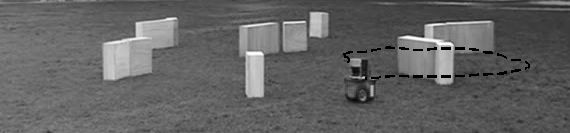}} 
    \subfigure[]{\includegraphics[width=0.45\columnwidth]{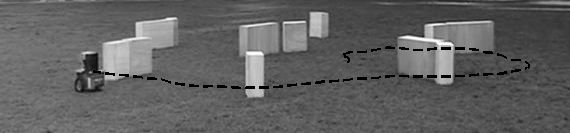}}
    \subfigure[]{\includegraphics[width=0.45\columnwidth]{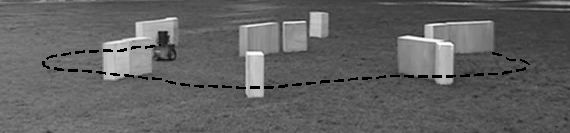}} 
    \subfigure[]{\includegraphics[width=0.45\columnwidth]{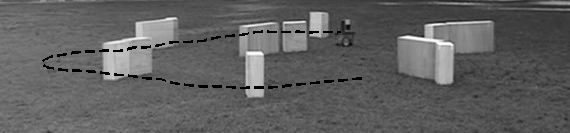}}
    \subfigure[]{\includegraphics[width=0.45\columnwidth]{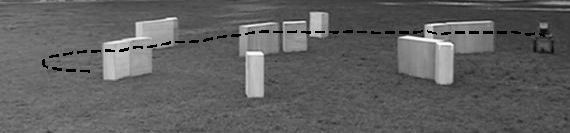}}  
  \caption{Sequence of images showing the experiment.}
  \label{tf:fig:real}
\end{figure}

The vehicle successfully navigates around the obstacles without
  collision, as indicated in Fig. \ref{tf:fig:real}.

\clearpage \section{Summary}
\label{tf:S7}

The shortest (minimal in length) path for a
unicycle-like mobile robot in a known environment is described, which contains smooth
(possibly non-convex) obstacles. Furthermore, a reactive randomized
algorithm for robot navigation in unknown environments is proposed. The performance of this algorithm has been confirmed
by computer simulations and outdoor experiments with a Pioneer P3-DX
mobile wheeled robot.

\chapter{Nonlinear Sliding Mode Control of an
Agricultural Tractor}
\label{chapt:pf}

This chapter considers the problem of automatic path tracking by
autonomous farming vehicles subject to wheel slips, which are
characteristic for agricultural applications. The control law used
is similar to the one employed in Chapt.~\ref{chap:singlevehicle} 
to regulate the position of a unicycle-type vehicle. Simulations compare 
the performance of this mixed nonlinear-sliding mode control law with both 
a pure sliding mode control law and another control law proposed in the literature.
Real world experiments confirm the applicability and performance of the
proposed guidance approach.

\section{Introduction}
 \label{pf:sec.intr}

Automatic guidance of agricultural vehicles is a widely studied problem in
robotics. Such automatic systems are
aimed at replacement of human beings in tedious or
hazardous operations, and improve productivity while reducing work hardness. 
Driver assistance devices and automatic steering
systems have attracted considerable interest (see e.g.
\cite{KeSE00,ReZhNoDi00}). The reported systems provide
satisfactory performance if the motion is nearly pure rolling without slip -- wheel slips usually
causing visible performance degradation. Unfortunately, slipping
effects commonly occur during agricultural tasks since the
farming terrain is often non-smooth and undulating, 
and the vehicle is affected by disturbances from the towed
implements interacting with the ground.
\par
There is an extensive literature on path following control of
wheeled vehicles in the no-slip circumstances, see e.g.
\cite{Bloch03,BuLe05,Raj06}. Even in this
case, the problem is rather challenging. Due to the
nonholonomic nature, the system cannot be asymptotically stabilized
by a smooth time-invariant feedback \cite{Brock83}. Slip effects can be incorporated into this analysis
if exact knowledge of the slipping
parameters is available \cite{MoCa00}. Slipping has also been modeled as fast
dynamics, with a singular perturbation approach being used to achieve
robust tracking under sufficiently small sliding \cite{AnCaBa95,LeDan97}. Asymptotic tracking has been achieved for
slightly curved reference paths using a time varying
controller \cite{LeDan96}. Adaptive control has been proposed in the form of adaptive back-stepping
and adaptive re-scheming of the desired path \cite{LeThCaMa03,FaLeThMa05,
FaRuThMa06, LRTBC06}. A differential flatness approach has been applied to a
car-like robot controlled by both steering and front-wheel drive
speed, however rear wheel drive vehicles are mostly employed in
agricultural applications \cite{AgRy08}. A MPC element has 
been proposed to account for actuator delay when added to a control law \cite{LRTBC06}.
Kinematic control law in combination with dynamic
observers estimating slip parameters have been used to 
achieve improved path tracking accuracy \cite{LTCM10,LTHM11,FDCLTM11}.
\par
Due to the well-known benefits, such as stability under large
disturbances and robustness against system uncertainties, the
sliding mode approach attracts an increasing interest in the area of
agricultural vehicles control. A discrete-time sliding mode control
has been proposed for trajectory tracking in the
presence of skidding effects \cite{CoOr02}. The sliding mode
approach has also been used to achieve global lateral stabilization 
via steering actuation only \cite{HaLeThMa04}. Here the drive speed was assumed
constant and the slipping effects were represented by additive
disturbances in the ideal kinematic model, both in the original and
in the transformed chained forms. The ideal
kinematic equations can also be modified to include slip-induced biases in the
steering angle and drive speed  \cite{EaKaPoSi09}. Assuming full
 actuation, this chapter reports successful use of
the combined sliding mode and back-stepping techniques to achieve
both lateral and longitudinal stabilization. A similar problem has been treated
 along the lines of the singular
perturbation approach, and the problem was
solved using a second order sliding mode control 
\cite{HaAcFlPe05,AnCaBa95}. However, for some
agricultural applications, the precise longitudinal stabilization is
of minor importance, in which case the only-steering control is a
reasonable option.
\par
The characteristic feature of the above works is the lack of concern
for the actual steering capabilities of the vehicle. The stability 
of path tracking is guaranteed only under the implicit and unjustified assumption
that the outputs of the proposed sophisticated nonlinear controller
always obey the mechanical steering angle limit. This is definitely 
plausible in favorable circumstances. Unfortunately, wheel
slippage, large lateral disturbances and sharp contortions in the 
reference path may require the vehicle to operate near the 
steering angle limit, and this may lead to excessively
large controller outputs and thus violation of the above assumption.
This is a likely cause for systematic occasional degradation of the 
tracking accuracy observed for existing laws.
\par
From the more general perspective, there is only limited literature 
available about stabilization and tracking of nonholonomic wheeled 
vehicles with input saturations. For nonholonomic systems in chained form, 
discontinuous controllers obeying a given bound and stabilizing an 
equilibrium have been obtained \cite{LuTs00,AlMa03}.
A passivity-based approach has been applied to stabilize equilibria 
of special controllable drift-less systems using a time-varying smooth 
state feedback controller \cite{Lin96}.
Semi-global practical stabilizing control schemes for wheeled robots 
with two steering wheels and one castor wheel have also been proposed
 \cite{Wang08}. However, all these works did not treat path tracking. This issue has been addressed using passivity-based, saturated, Lipschitz continuous, time-varying feedback
laws for unicycles with forward and turning velocity bounds \cite{JLN01}.
Trajectory tracking controllers for nonholonomic
robots with fixed-wing UAV-like input constraints and experimental results on their performance are also available \cite{RSBML05}.
The characteristic feature of these works is the
use of unnecessary and non-conservative bounds on controller outputs, which in general 
does not permit the controller to employ the full range of the vehicle's steering capabilities 
and may cause performance degradation in challenging circumstances. Another restriction is that 
these works assume pure wheels rolling without slipping. Off-road path tracking with 
both wheel slip and explicit concern for steering saturation has been addressed, 
using a model predictive algorithm for the longitudinal speed adjustment.
 This enables the main steering controller to be decoupled and be borrowed 
 from \cite{LthCM07,LTCM10,LTHM11}, which can adaptively cope with contortions of 
 the reference path during cornering in the face of wheel slip and steering angle 
 limit. However even in this work, the design of the steering controller still neglects this limit.
\par
To fill this gap, in this chapter the problem of globally stable tracking of 
an arbitrarily curved path via bounded steering
only is considered -- the nominal drive speed is constant and the controller 
must ensure that the steering angle bound is always satisfied. Following
\cite{FaRuThMa06,LeThCaMa06,Micaelli1993book3,FFTM06}, the slipping effects are
modeled via biases of the steering angle and vehicle velocity in the
ideal kinematic equations, where they are treated as bounded
uncertainties. To design a controller, a sliding mode approach is employed. 
Unlike the previous research, the requirements to the sliding surface 
that are imposed by the steering angle limit are explicitly disclosed. 
To fully use the steering capabilities of the vehicle, the optimal 
surface among those satisfying these requirements is found -- this surface minimizes 
the maximal steady-state error.
Furthermore, the proposed sliding mode control laws do not artificially 
impede the steering angle within the given limits.
Two control laws are examined and compared.
The first of them is a pure sliding-mode
controller that formally requests only limit values of the steering angle.
As was discovered via simulation tests and experiments, this controller may 
cause an oscillatory behavior in some cases. So it was used as a theoretical 
platform for design of the second law that combines a smooth nonlinear 
control with  switching in the sliding-mode fashion between reduced limits 
determined only by the need to
reject the slipping effects. The smooth control is in fact the equivalent 
control for the first law \cite{Utkin1992book1}. More precisely, since 
computation of the equivalent control requires unknown slip parameters, a reasonable
approximation is used.
\par
The applicability and performance of the proposed control law is validated via computer
simulations and real world testing using an agricultural vehicle.
\par
All proofs of mathematical statements are omitted here; they are
available in the original manuscript \cite{Matveev2010conf0}.
\par
The body of this chapter is organized as follows. Secs.~\ref{pf:sec1} and
\ref{pf:sec.model} introduce the problem setup and the employed
kinematic model. Sec.~\ref{pf:sec.dd} offers the
study of sliding surfaces; the proposed control laws are discussed in
Secs.~\ref{pf:sec.ma} and \ref{pf:sec.ra}. Sec.~\ref{pf:sec.sr}
presents the simulation results, and experimental results are 
presented in Sec.~\ref{pf:sec:expresult}. Finally, brief conclusions are given
in Sec. \ref{pf:con}.

 \section{Problem Statement} 
\label{pf:sec1} 

\begin{figure}[ht]
  \centerline{\includegraphics[width=10cm]{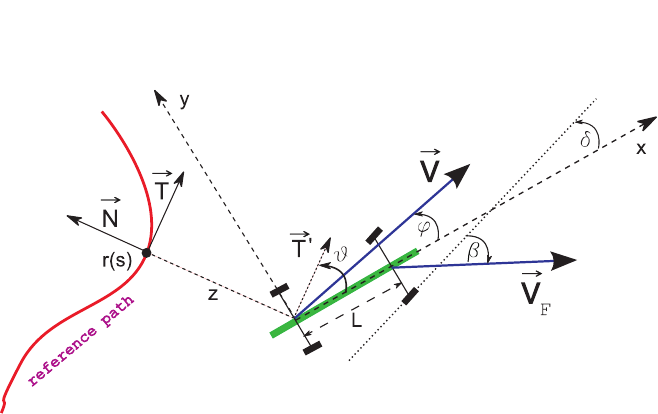}}
  \caption{Basic variables and constants.} \label{pf:f.1}
\end{figure}

A planar wheeled mobile robot modeled as a bicycle with front
steering is considered. It is controlled by the front wheels steering angle and
travels with the constant rotating angular velocity of the rear
wheels. The objective is to follow a given reference path as close as
possible. To describe the robot motion, the {\it
  relative} vehicle-fixed normally oriented Cartesian coordinate
system with the $x$-axis directed along the vehicle centerline from
rear to front is introduced, along with the following variables (see
Fig.~\ref{pf:f.1}):

\begin{itemize}
\item $\mathfrak{F}$ -- the center of the front wheels axle;
\item $\mathfrak{R}$ -- the center of the rear wheels axle;
\item $s$ -- the curvilinear abscissa of the point on the reference
  path that is closest to $\mathfrak{R}$, the path is oriented so that
  $s$ ascends in the required motion direction;
\item $r(s) \in \mathbb{R}^2$ -- a regular parametric representation
  of the reference path in the world frame;
\item $[T(s), N(s)]$ -- the Frenet frame of the path at $r(s)$:
  $T(s)$ is the unit positively oriented tangential vector, $N(s)$ is
  the unit normal, the frame is normally oriented;
\item $\varkappa(s)$ -- the signed curvature of the reference path:
  $\varkappa = \spr{T^\prime}{N}$, where $^\prime$ is differentiation
  with respect to $s$ and $\spr{\cdot}{\cdot}$ is the standard inner
  product in the plane $\br^2$;
\item $z$ -- the $N$-coordinate of $\mathfrak{R}$ in the Frenet
  frame;
\item $\vec{V}$ -- the velocity of the point $\mathfrak{R}$ in the
  world frame;
\item $\vec{V}_{\text{F}}$ -- the velocity of the point
  $\mathfrak{F}$ in the world frame;
\item $v$ -- the speed of the point $\mathfrak{R}$, i.e.,
  $v=\|\vec{V}\|$;

\item $\varphi$ -- the rear wheel slip angle, i.e., the angular polar
  coordinate of $\vec{V}$ in the relative coordinate system;
\item $\theta$ -- the angular polar coordinate of the tangential
  vector $T(s)$ in the relative coordinate system;
\item $\delta$ -- the steering angle of the front wheels;
\item $\delta_{\max}$ -- the maximal steering angle: $\delta \in
  [-\delta_{\max}, \delta_{\max}]$;
\item $\beta$ -- the steering angle bias due to sliding;
\item $L$ -- the vehicle wheelbase;

\item $v_{\text{rw}}$ -- the driving speed: $v_{\text{rw}} = \rho
  \cdot \omega_{\text{r}}$, where $\rho$ is the radius of the rear
  wheel and $\omega_{\text{r}}$ is its angular velocity;
\item $v_{\text{rw}}^{\text{s}}$ -- the bias of the driving speed due
  to wheels slip, i.e., the projection of the velocity of the rear
  wheel point of contact with the ground on the vehicle centerline.
\end{itemize}

The objective is to design a controller such that

$$
z \to 0 \qquad \text{as} \quad t \to \infty.
$$

To this end, the vehicle is equipped with sensors that ensure access
to the positional $z$ and angular $\theta$ errors of the path
following. There also is an access to the path curvature
$\varkappa(s)$. 
\par
Direct measurement of the slip parameters $\varphi, \beta$, and
$v_{\text{rw}}^{\text{s}}$ is hardly possible at a reasonable cost
(see e.g., \cite{StChMe04}). So they are treated as bounded system
uncertainties:

\begin{equation}
  \label{pf:up.bound} |\varphi| \leq \overline{\varphi}, \quad |\beta|
  \leq \overline{\beta}, \quad \left| v_{\text{rw}}^{\text{s}} \right|
  \leq \overline{v}.
\end{equation}

It is assumed that the wheels slip can reverse neither the translational
nor, under the maximum steering angle, angular directions of the
vehicle motion:

\begin{equation}
  \label{pf:up.bound2} \overline{v} < v_{\text{rw}}, \qquad
  \overline{\beta} < \delta_{\max}.
\end{equation}

In agricultural applications, the biases $\varphi, \beta$ are usually
small ($\approx$ several degrees) \cite{FFTM06}.  For technical
reasons, it is assumed:

\begin{equation}
  \label{pf:up.bound1} \ov{\varphi} < \pi/2, \qquad \ov{\beta} < \pi/2.
\end{equation}

Through not required, indirect on-line estimates $\widehat{\varphi}$ and
$\widehat{\beta}$ of $\varphi$ and $\beta$, may be
available, with the errors:

\begin{equation}
  \label{pf:up.boun} \left| \varphi - \widehat{\varphi} \right| \leq
  \eer{\varphi} , \qquad \left| \beta - \widehat{\beta} \right| \leq
  \eer{\beta}.
\end{equation}

If estimation is not carried out,
$\widehat{\varphi}:= \widehat{\beta}:= 0$ and so $\eer{\varphi} =
\ov{\varphi}, \eer{\beta} = \ov{\beta}$. The estimates are supposed to
be improvable via neither annihilation nor saturation at the bound
from Eq.\eqref{pf:up.bound}:

\begin{equation}
  \label{pf:non.impr} \eer{\varphi} \leq \overline{\varphi}, \quad
  \eer{\beta} \leq \overline{\beta}, \quad |\widehat{\varphi}| \leq
  \ov{\varphi}, \quad |\widehat{\beta}| \leq \ov{\beta}.
\end{equation}

It is also assumed that $\varphi, \widehat{\varphi}, \beta,
\widehat{\beta},v_{\text{rw}}^{\text{s}}$ continuously depend on time
$t$.

 \section{Kinematic Model}
\label{pf:sec.model}

Basic assumptions are given as follows:

\begin{Assumption}
  The reference path is $C^1$-smooth and $C^2$-piece-wise smooth
  regular curve. The curvature of the reference path is upper bounded:
  \begin{equation}
    \label{pf:b.cur} \overline{\varkappa} := \sup_s |\varkappa(s)|  <
    \infty.
  \end{equation}
\end{Assumption}

\begin{Assumption}
  \label{pf:can.track} The vehicle is capable of tracking the
  reference path: the path curvature radius at any point is no less
  than the minimal turning radius $R= \frac{L}{\tan \delta_{\max}}$ of
  the vehicle. Moreover, the associated strict inequality holds:
  \begin{equation}
    \label{pf:encan} \ov{\varkappa} < \vveh:=R^{-1}= L^{-1}\tan
    \delta_{\max}.
  \end{equation}
\end{Assumption}

Vehicle maneuvers are examined at distances not exceeding the
curvature radius of the reference path:

\begin{equation}
  \label{pf:b.z} |z | \leq \overline{\varkappa}^{-1}\lambda
\end{equation}

Here $0 < \lambda <1$ is a given parameter.\footnote{This parameter is used to underscore
  that the distance $|z|$ should be uniformly less than the path
  curvature radius; typically, $\lambda \approx 1$. Furthermore, an
  implicit requirement to the controller is in fact imposed here: the
  vehicle should not leave the domain Eq.\eqref{pf:b.z}. This requirement is
  enhanced as $\lambda \downarrow$. If the path is straight,
  $\overline{\varkappa}^{-1} = 0^{-1}:=\infty$, and Eq.\eqref{pf:b.z}
  imposes no bounds on $z$.} Hence the point on the path that is
closest to $\mathfrak{R}$ is well defined and smoothly depends on the
vehicle position. The kinematic equations given by the following lemma
are borrowed from \cite{Micaelli1993book3,FFTM06,LeThCaMa06,
  20110613643991}:

\begin{Lemma}
  \label{pf:lem.kineq} The vehicle kinematic model is as follows:
  \begin{equation}
    \label{pf:eq1} \dot{s} = \frac{v \cos (\theta - \varphi) }{1-
      \varkappa(s) z }, \; \dot{z} =  - v \sin (\theta - \varphi), \; v =
    \frac{v_{\text{rw}} + v_{\text{rw}}^{\text{s}}}{\cos \varphi} ,
  \end{equation}
  \begin{equation}
    \label{pf:eq2} \dot{\theta} = \varkappa(s)  \frac{v \cos (\theta -
      \varphi) }{1- \varkappa(s) z } - \frac{v}{L} \left[
      \tan(\delta+\beta) - \tan \varphi \right] .
  \end{equation}
\end{Lemma}

 \section{Desired Dynamics}
\label{pf:sec.dd}

This chapter is concerned with relations of the form:

\begin{equation}
  \label{pf:c.a} \theta  = \chi(z) +\widehat{\varphi}.
\end{equation}

\begin{figure}[ht]
  \centerline{\includegraphics[width=8cm]{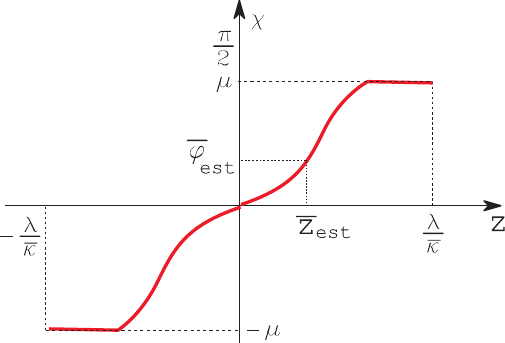}} \caption{The
    desired relation $\theta=\chi(z)$.} \label{pf:f3}
\end{figure}

It is assumed that $\chi(\cdot)$ is an odd ascending
continuous piece-wise smooth function such that $\frac{d
  \chi}{dz}(z\pm 0 )>0$ whenever $-\mu < \chi(z) < \mu:=
\chi(\lambda/\ov{\varkappa})$ (see
Fig.~\ref{pf:f3}). A simple sample employed here is the linear function with
saturation: $\chi(z) := \gamma z$ if $|z| \leq \Delta$ and $\chi(z) :=
\mu \sgn (z)$ otherwise.  Here $\mu := \gamma\Delta$ and $\gamma >0,
0<\Delta \leq \lambda/\ov{\varkappa}$ are design parameters.

\subsection{Conditions for Robustly Stable Path Tracking}

The first lemma characterizes the positional error caused by the
errors in estimation of the wheels slip parameters.
\begin{Lemma}
  \label{pf:lem.conv} Let the motion starts in the domain Eq.\eqref{pf:b.z},
  \begin{equation}
    \label{pf:mu.b} \eer{\varphi} < \mu := \chi\left(\lambda/\ov{\varkappa}
    \right) < \pi/2 - \eer{\varphi}, \footnote{The last inequality is
      not essential and is imposed only to simplify computations.}
  \end{equation}
  and let $\eer{z}$ denote the unique root of the equation
  $\chi(\eer{z}) = \eer{\varphi}$.  Whenever Eq.\eqref{pf:c.a} is maintained
  during the motion, the vehicle remains in the domain Eq.\eqref{pf:b.z} and
  the positional error $|z|$ of the path following monotonically
  decays to $\eer{z}$.  This decay either is asymptotical $|z(t)|
  \xrightarrow{t \to \infty} \eer{z}, |z(t)| > \eer{z}\; \forall t$ or
  is completed for a finite time $\exists t_\ast: |z(t)| \leq \eer{z}
  \; \forall t \geq t_\ast$. In any case, the decay is at no less than
  the exponential rate:
  \begin{equation}
    \label{pf:exp.rate} |z(t)| \leq \eer{z} + c e^{- \gamma \left(
        v_{\text{rw}} - \ov{v}\right) t} , \qquad \gamma := \frac{d
      \chi}{dz}(\eer{z}+0).
  \end{equation}
\end{Lemma}

\subsection{Feasibility of the Desired Dynamics}

Restrictions on the vehicle steering angle may make maintenance of
relation Eq.\eqref{pf:c.a} unrealistic. To disclose the conditions under
which this is not the case, the following is introduced:

\begin{Definition}
  Relation Eq.\eqref{pf:c.a} is said to be {\it nominally maintainable in
    the domain} Eq.\eqref{pf:b.z} if whenever it is achieved $S:= \theta -
  \chi(z) -\widehat{\varphi} =0$ at a vehicle position from this
  domain, it can be maintained by a proper choice of the front wheels
  steering angle: $\exists \delta \in [- \delta_{\max},
  \delta_{\max}]: \dot{S}(\delta)=0$.
\end{Definition}

It is assumed here that the choice of $\delta$ is based on not only
the measured data $z,\theta, \varkappa(s)$ and available estimates
$\widehat{\varphi}, \widehat{\beta}$ but also on access to the slip
uncertainties $\varphi, \beta, v^s_{\text{rw}}$. It will be shown that
conditions necessary for maintainability of Eq.\eqref{pf:c.a} in these
idealized circumstances are `almost sufficient' for its
maintainability in the realistic case where the above uncertainties
are unknown.

\begin{Lemma}
  \label{pf:lem.fies} Relation Eq.\eqref{pf:c.a}
  is nominally maintainable in the region Eq.\eqref{pf:b.z} if
  \begin{multline}
    \label{pf:cond.feaslide} \vvehu:=
    \frac{\tan(\delta_{\max}-\ov{\beta})}{L} - \frac{\tan
      \ov{\varphi}}{L} -
    \frac{\ov{d\widehat{\varphi}}}{v_{\text{wr}}-\ov{v} } - \frac{2
      \ov{\varkappa}}{1-\lambda} \eer{\varphi}
    \\
    \geq \frac{ \ov{\varkappa} \cos [\chi(z)+\eer{\varphi}] }{1-
      \ov{\varkappa} z } + \frac{d \chi}{dz}(z\pm 0) \sin \left[
      \chi(z) + \eer{\varphi}\right]
  \end{multline}
  for all $z \in [0, \lambda/\ov{\varkappa}]$, where
  \begin{equation}
    \label{pf:der.upb} \ov{d\widehat{\varphi}} \geq \sup_t |\frac{d
      \widehat{\varphi}}{dt}(t)|
  \end{equation}
  is an available upper bound of the estimate derivative.
\end{Lemma}

 \section{Sliding Mode Control with Maximal Actuation}
\label{pf:sec.ma}

The simplest sliding mode control law can be expressed as follows:

\begin{equation}
  \label{pf:c.law0s} \delta= \delta_{\max} \sgn(S), \qquad S:= \theta -
  \chi(z) - \widehat{\varphi}.
\end{equation}

Let $\mathfrak{D}$ denote the set of points $\zeta \in \br^2$ for
which Eq.\eqref{pf:b.z} holds and the distance from $\zeta$ to the reference
path $\mathfrak{P}$ is not furnished by an end-point of $\mathfrak{P}$
if the path is not a closed curve.
\par
To state the results in the case where $S(0) \neq 0$, consider the
motions of the vehicle driven by $\delta \equiv \delta_{\max} \sgn
S(0) $ from the initial state in the absence and presence of the
wheels slip, respectively; the values attributed to the second motion
are marked by $\widetilde{\phantom{x}}$. The vector $\zeta(t)$ of the
vehicle absolute Cartesian coordinates runs over the {\it initial
  circle} $C^{\text{in}}$ (of the radius $R$ given by Eq.\eqref{pf:b.z}) for
$T:=2 \pi R/v_{\text{rw}}$ time units. Let $d_\text{dev}$ denote the
maximal positional deviation $\max_{t \in [0,2T]}\|\zeta(t) -
\widetilde{\zeta}(t)\|$ that may be caused by the wheels slip under
the given constraints Eq.\eqref{pf:up.bound}. The {\it initial disc}
$D^{\text{in}}$ is that encircled by the initial circle
$C^{\text{in}}$.

\begin{Assumption}
  \label{pf:sm.sl1} The initial position of the vehicle lies in the
  domain $\mathfrak{D}$. Moreover, if $S(0) \neq 0$, this domain
  covers the $d_{\text{dev}}$-neighborhood $D_\ast$ the initial disc
  $D^{\text{in}}$.
\end{Assumption}

Let $\alpha$ be the maximal angular discrepancy
between them that may be caused for $t \in [0,2T]$ by slipping limited
according to Eq.\eqref{pf:up.bound}:

\begin{Assumption}
  \label{pf:sm.sl2} The angular deviation $\alpha < \pi$ if $S(0) \neq
  0 $
\end{Assumption}

There exists a function $\chi(\cdot)$ for
which Eq.\eqref{pf:cond.feaslide} holds with the strict inequality sign as
well -- for $z \in \left[0, \lambda/\ovk \right]$:

\begin{equation}
\label{pf:unif.disc} \vvehu > \frac{ \ov{\varkappa} \cos
[\chi(z)+\eer{\varphi}] }{1- \ov{\varkappa} z }
 + \frac{d \chi}{dz}(z\pm 0)
 \sin \left[ \chi(z) + \eer{\varphi}\right].
\end{equation}

\begin{Proposition}
  \label{pf:prop.ssl} Let Assumptions~\ref{pf:sm.sl1} and
  \ref{pf:sm.sl2} hold and the controller Eq.\eqref{pf:c.law0s} employ the
  function $\chi(\cdot)$ satisfying Eq.\eqref{pf:unif.disc}. Then the
  vehicle driven by this controller converges to the reference path at
  no less than the exponential rate and then tracks this path with the
  precision $\eer{z}$ -- Eq.\eqref{pf:exp.rate} holds, where $\eer{z}$ is the
  root of the equation $\chi(\eer{z})=\eer{\varphi}$.
\end{Proposition}

 \section{Sliding Mode Control with Reduced
  Actuation} \label{pf:sec.ra}

  The sliding mode control law
Eq.\eqref{pf:c.law0s} requests highly frequent oscillations of the steering
angle between the extreme values $-\delta_{\max}$ and
$\delta_{\max}$. In practice, such a steering pattern may not only
seem strange but also be unrealistic since the turn of the front
wheels between the extreme angles may require excessive effort and
time. In fact, this pattern merely serves as a machinery to generate
the equivalent control \cite{Utkin1992book1}, i.e., that driving the
system over the sliding surface $S:=\theta - \chi(z) -
\widehat{\varphi}=0$. In the ideal situation of pure rolling without
wheels slipping, this machinery is not required since the equivalent
control can be directly computed from the measurements. In the face of
the wheel slips, its use is motivated by the lack of access to the
required slip parameters $\varphi, \beta$. However, since these
uncertainties are usually small (at most several degrees according to
\cite{FFTM06}) in practical applications, beating them with the
maximal steering angles looks superfluous.
\par
In this section, a control algorithm that combines the
sliding mode control at the reduced amplitude
$\underline{\delta}<\delta_{\max}$ with direct computation of the
reasonable approximation of the equivalent control is examined. This approximation
results from substitution of $\widehat{\varphi}, \widehat{\beta},
v_{\text{rw}}$ in place of $\varphi, \beta, v$, respectively:

\begin{multline}
  \label{pf:xi.def} \delta = - \widehat{\beta} + \arctan \left[ \tan
    \widehat{\varphi} + L \Xi\right], \quad \text{where} \quad \Xi :=
  \\
  \frac{\varkappa(s) \cos(\theta-\widehat{\varphi})}{1-\varkappa(s)z}
  + \frac{d \chi}{dz}(z) \sin (\theta-\widehat{\varphi}) -
  \frac{1}{v_{\text{rw}}} \frac{d \widehat{\varphi}}{dt} .
\end{multline}

So far as the equivalent control is reasonable only on the sliding
surface $S=0$, these controls are mixed so that the participation of
the equivalent control decays as $|S|\uparrow$.
\par
Though the idea is to make the amplitude of control oscillations
small, it should be large enough to cope with slipping. To specify its
choice, start with the following:

\begin{Lemma}
  \label{pf:lem.delicate} Let Eq.\eqref{pf:unif.disc} hold. Then for $z \in
  Z:= \left[-\frac{\lambda}{\ovk}, \frac{\lambda}{\ovk} \right]$,
  \begin{multline}
    \label{pf:mueta} \lim_{\sigma \to 0} \left[
      \mu_\sigma(z)+\eta_\sigma(z) \right] + \frac{\tan
      \ov{\varphi}}{L} < \ovk_{\text{v}}:=
    \frac{\tan(\delta_{\max}-\ov{\beta})}{L} ,
    \\
    \text{\rm where} \quad \eta_\sigma(z):= \sup_{|S|\leq \sigma}
    |\Xi|,
    \\
    \mu_{\sigma}(z) := \sup_{|S|\leq \sigma} \left|v^{-1}\dot{S} +
      \frac{\tan[\delta+\beta]-\tan \varphi}{L} - \Xi \right|,
  \end{multline}
  the limit exists, is uniform over $z\in Z$, and $\mu_\sigma(\cdot),
  \eta_\sigma(\cdot)$ are continuous functions. Here $S= \theta -
  \chi(z) - \widehat{\varphi}$ and $\sup$ is over $\theta, \varphi,
  \widehat{\varphi}, d \widehat{\varphi}/dt, \beta,
  \widehat{\beta},\delta, v$ that along with the condition $|S|\leq
  \sigma$, satisfy the uncertainty constraints Eq.\eqref{pf:up.bound},
  Eq.\eqref{pf:up.boun}, Eq.\eqref{pf:non.impr}, Eq.\eqref{pf:der.upb}.
\end{Lemma}

Lemma~\ref{pf:lem.delicate} permits us to pick $\sigma>0, \varepsilon
> 0$, and a continuous function $\underline{\varkappa}(z)$ such that
whenever $|z| \leq \lambda/\ovk$,
\begin{multline}
  \label{pf:bb.ii} \eer{\varphi} + \sigma < \chi(\lambda/\ovk);
  \\
  \mu_\sigma(z)+\eta_\sigma(z) + \frac{\tan \ov{\varphi}}{L} \leq
  \ovkv - \varepsilon, \quad \vveh> \unk(z) \geq \mu_\sigma(z)
  \\
  + \frac{\sin \eer{\varphi}}{L\cos^2\ov{\varphi}} + \tan
  \delta_{\max} - \tan(\delta_{\max} - \eer{\beta}) + \varepsilon .
\end{multline}

A continuous function $p(c,z), c \geq 0$ is picked such that:

\begin{multline}
  \label{pf:unka} p(0,z) = \unk(z), \; p(c,z) = \vveh \quad \forall c
  \geq \sigma, \quad p(c,z) \in [\unk(z), \vveh]\quad \forall c \geq
  0.
\end{multline}

Finally, the following control law is introduced:

\begin{multline}
  \label{pf:cl.ma} \delta := \boldsymbol{SAT} \left[ \arctan \Upsilon
    - \widehat{\beta} \right], \quad \text{where}
  \\
  \Upsilon:= \frac{\vveh -p(|S|,z)}{\vveh - \unk(z)} \big( \tan
  \widehat{\varphi} + L \Xi\big) + L p(|S|,z) \sgn S,
\end{multline}

$S:=\theta - \chi(z) - \widehat{\varphi}$, $\Xi$ is given by
Eq.\eqref{pf:xi.def}, and $\boldsymbol{SAT}$ is saturation at the maximal
steering angle: $\boldsymbol{SAT}(\delta) := \delta$ if $|\delta| \leq
\delta_{\max}$ and $\boldsymbol{SAT}(\delta) := \delta_{\max}$
otherwise.

\begin{Proposition}
  \label{pf:prop.sslr} Let Assumptions~\ref{pf:sm.sl1},
  \ref{pf:sm.sl2} hold and the controller Eq.\eqref{pf:cl.ma} employs the
  functions $\chi(\cdot)$ and $\unk(\cdot)$ satisfying
  Eq.\eqref{pf:unif.disc} and Eq.\eqref{pf:bb.ii}.  Then the vehicle driven by
  this controller converges to the reference path at no less than the
  exponential rate and then tracks this path with the precision
  $\eer{z}$: Eq.\eqref{pf:exp.rate} holds, where $\eer{z}$ is the root of
  the equation $\chi(\eer{z})=\eer{\varphi}$.
\end{Proposition}

 \section{Simulations}
\label{pf:sec.sr}

Simulations were carried out to evaluate the two control laws (with
the maximal and reduced actuation, respectively) on the kinematic
model of the vehicle specified in Lemma \ref{pf:lem.kineq}. To emulate
difficult off-road conditions, the model was subjected to randomly
varying side-slip parameters $\varphi$ and $\beta$.  These were
independently drawn from the uniform Bernoulli distribution over
$\{-0.05, +0.05\}$ with a sampling period of $10 s$.  In addition to
the upper limit on the steering angle, steering dynamics were also
added to the system - the steering angle was not allowed to change
faster than a rate of $1.0 rad s^{-1}$.  This was done to better model
the actual steering capabilities of autonomous tractors.  A time step
of $0.05 s$ was used to update the control law. A reference trajectory
consisting of line and circle segments was generated as a sequence of
way-points, with the values of $z$, $\theta$, and $\chi$ being
calculated by tracking the perpendicular projection of the vehicle
onto the trajectory. The circular segments had a radius of $5 m$.  The
following user-selectable functions $\chi(z)$ and $p(|S|,z)$ were
used:

$$\chi(z) := sign(z)\cdot \min \left\{\frac{\mu |z|}{\Delta}, \mu \right\}, \quad p(|S|,z) := \min \{\unk + \frac{|z| \cdot (\vveh - \unk)}{\sigma}, \vveh \}. $$

\begin{table}[ht]
  \centering
  \begin{tabular}{| l | c |}
\hline
    $v_{rw}$ & $1 ms^{-1}$ \\
    \hline
    $\delta_{max}$ & $0.57 rad$ \\
    \hline
    $L$ & $1.69 m$  \\
    \hline
\end{tabular}
\hspace{10pt}
 \begin{tabular}{| l | c |}
\hline
    $\mu$ & $1.3 rad$\\
    \hline
    $\Delta$ & $1.2 m $\\
    \hline
\end{tabular}
\hspace{10pt}
 \begin{tabular}{| l | c |}
\hline
    $\sigma$& $0.2 m$\\
    \hline
    $\unk$ & $0.05m^{-1}$\\
    \hline
  \end{tabular}
  \caption{Control parameters used during simulation.}
  \label{pf:fig:paramsim}
\end{table}

The remaining relevant parameters are indicated in
Table~\ref{pf:fig:paramsim}.
\par
\begin{figure}[ht]
  \centering \includegraphics[width=10cm]{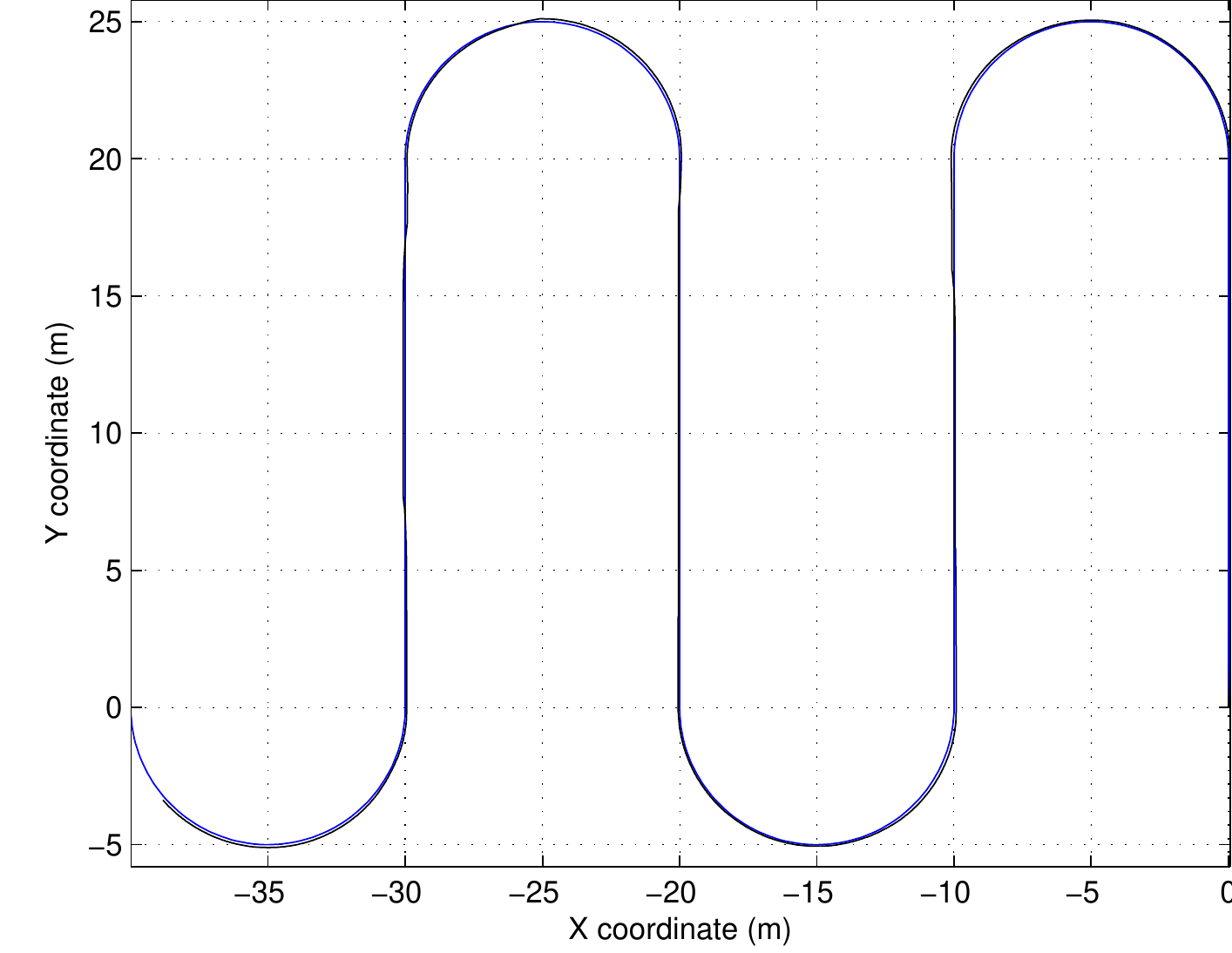}
  \caption{Simulations using the reduced actuation
    controller.}
  \label{pf:fig:trajsmin}
\end{figure}
\begin{figure}[ht]
  \centering \includegraphics[width=10cm]{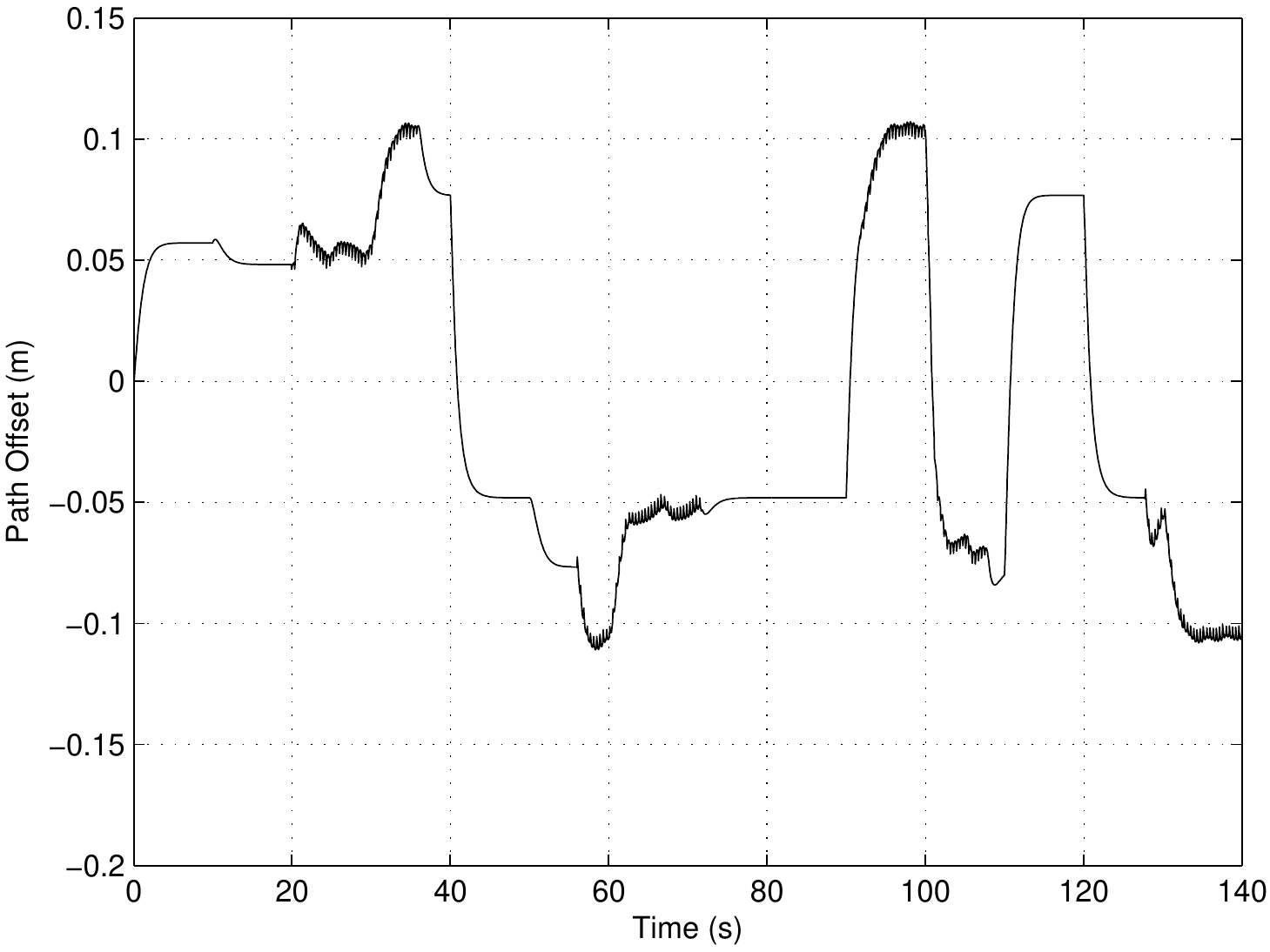}
  \caption{Offset error from the trajectory when using the reduced
    actuation controller.}
  \label{pf:fig:zerrorsmin}
\end{figure}
\begin{figure}[ht]
  \centering \includegraphics[width=10cm]{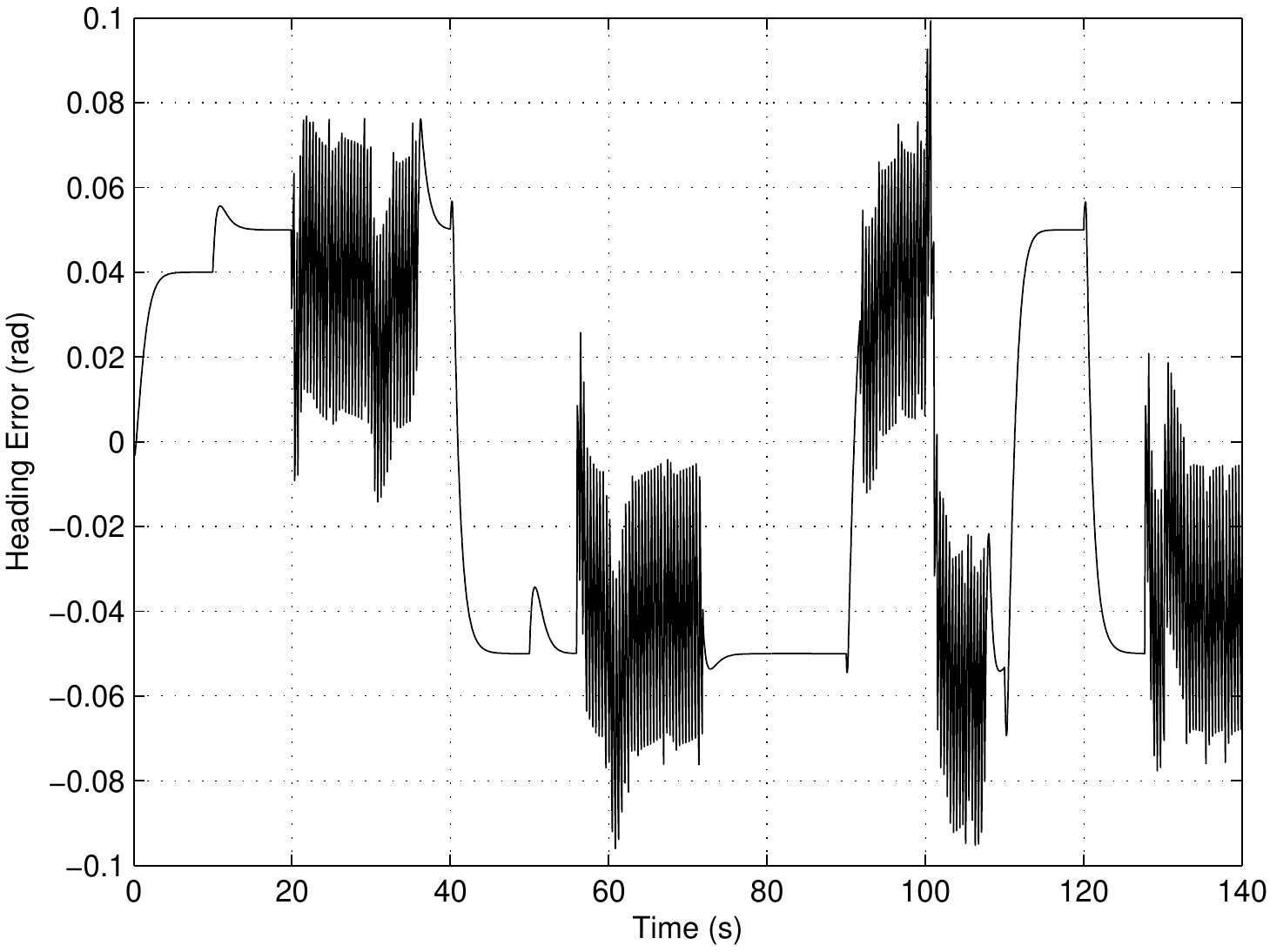}
  \caption{Heading error obtained when using the reduced actuation
    controller.}
  \label{pf:fig:therrorsmin}
\end{figure}
\begin{figure}[ht]
  \centering \includegraphics[width=10cm]{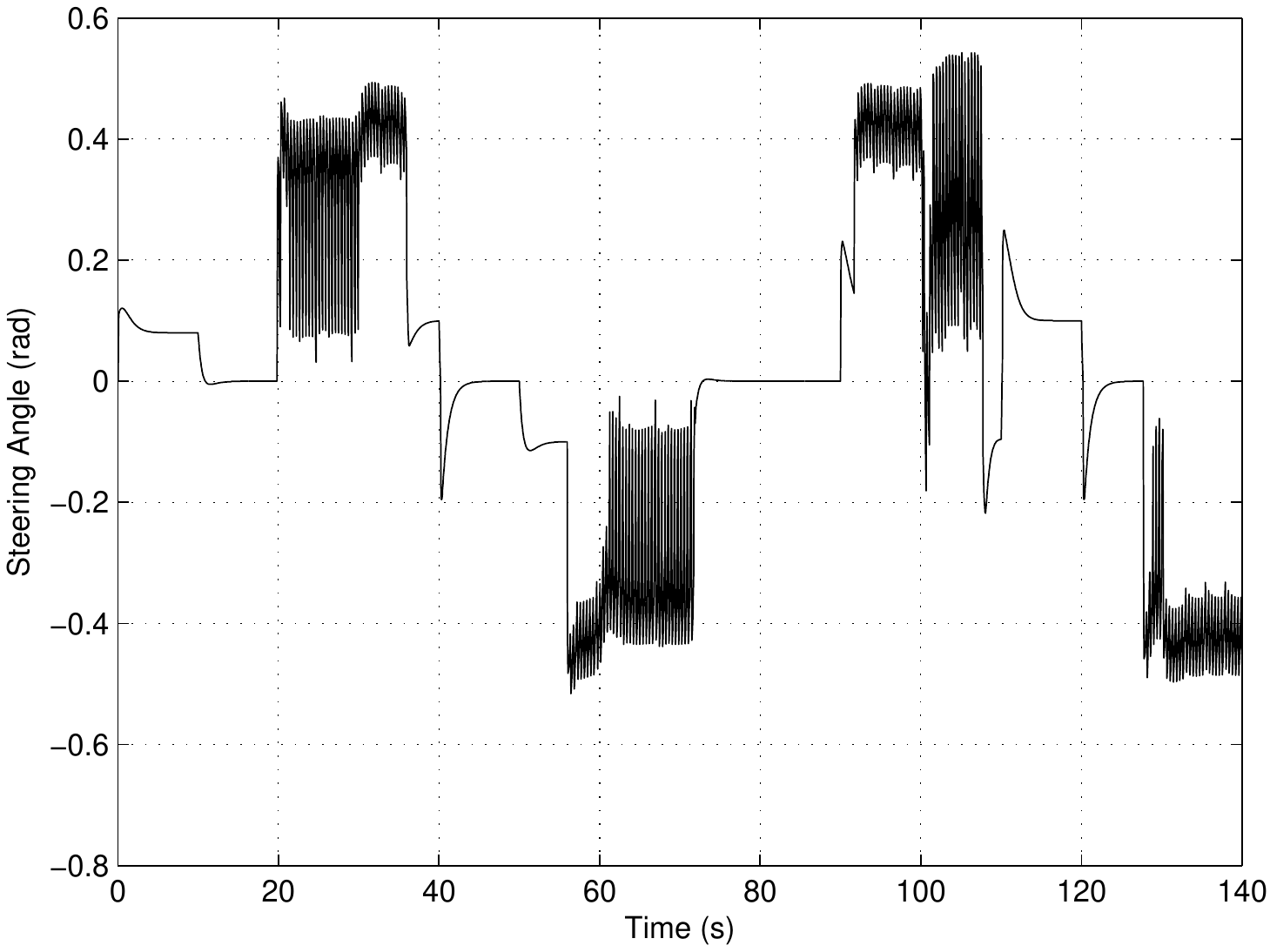}
  \caption{Steering angle requested when using the reduced
    actuation controller.}
  \label{pf:fig:steerangsmin}
\end{figure}
\begin{figure}[ht]
  \centering \includegraphics[width=10cm]{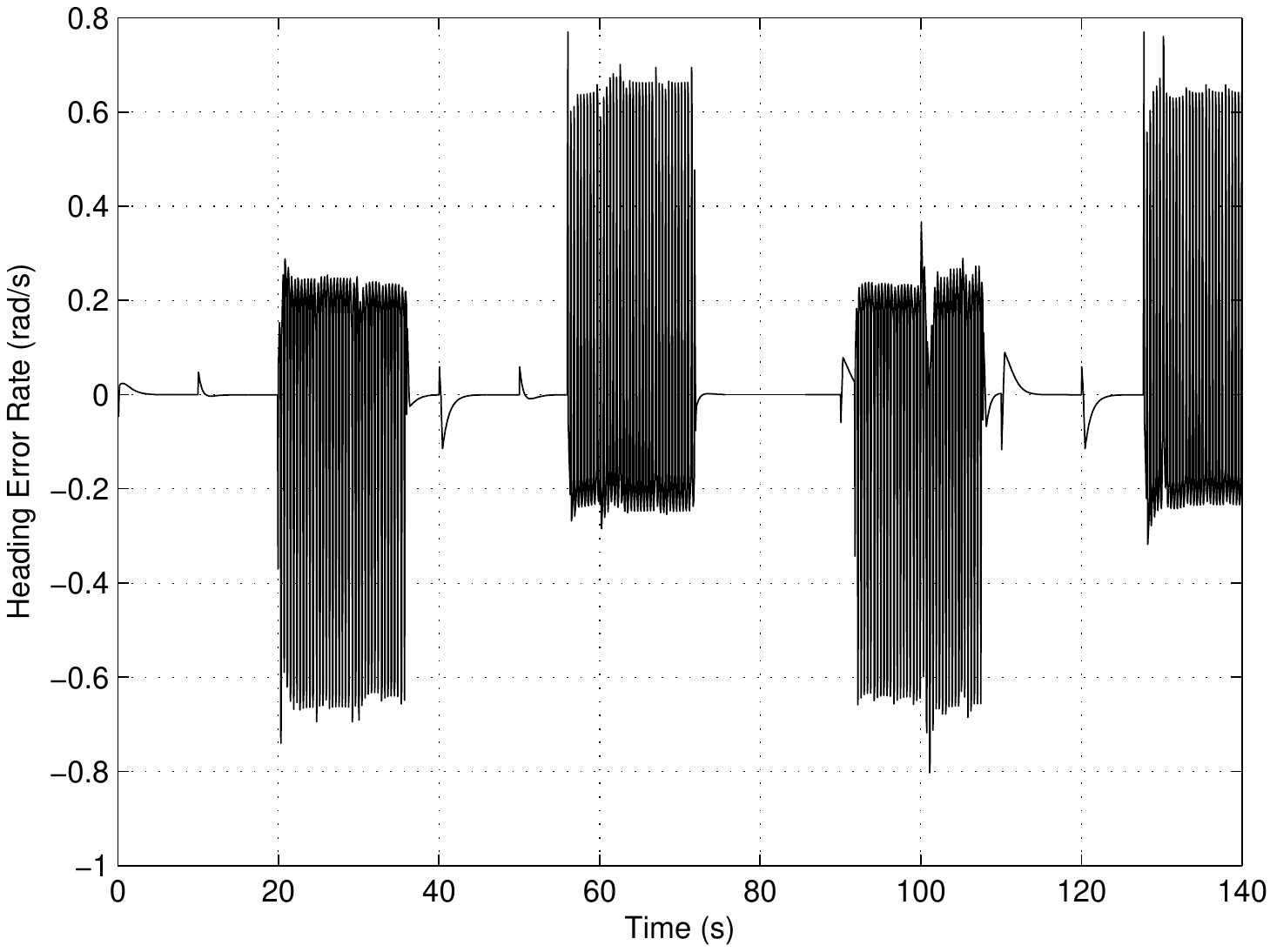}
  \caption{Rate of change of the heading error obtained when using the reduced actuation
    controller.}
  \label{pf:fig:therratesmin}
\end{figure}

Typical results obtained for the non-linear controller with reduced
actuation (see Sec.~\ref{pf:sec.ra}) are presented in
Figs.~\ref{pf:fig:trajsmin}, \ref{pf:fig:zerrorsmin},
\ref{pf:fig:therrorsmin}, \ref{pf:fig:steerangsmin} and 
\ref{pf:fig:therratesmin}.  The spatial
tracking error does not exceed $12cm$, whereas the orientation error
almost always lies within $[-0.1rad, 0.1rad]$, which features good
enough performance in the face of two $0.05rad$ actuation errors
($\varphi$ and $\beta$) due to wheels slipping.
\par
\begin{figure}[ht]
  \centering \includegraphics[width=10cm]{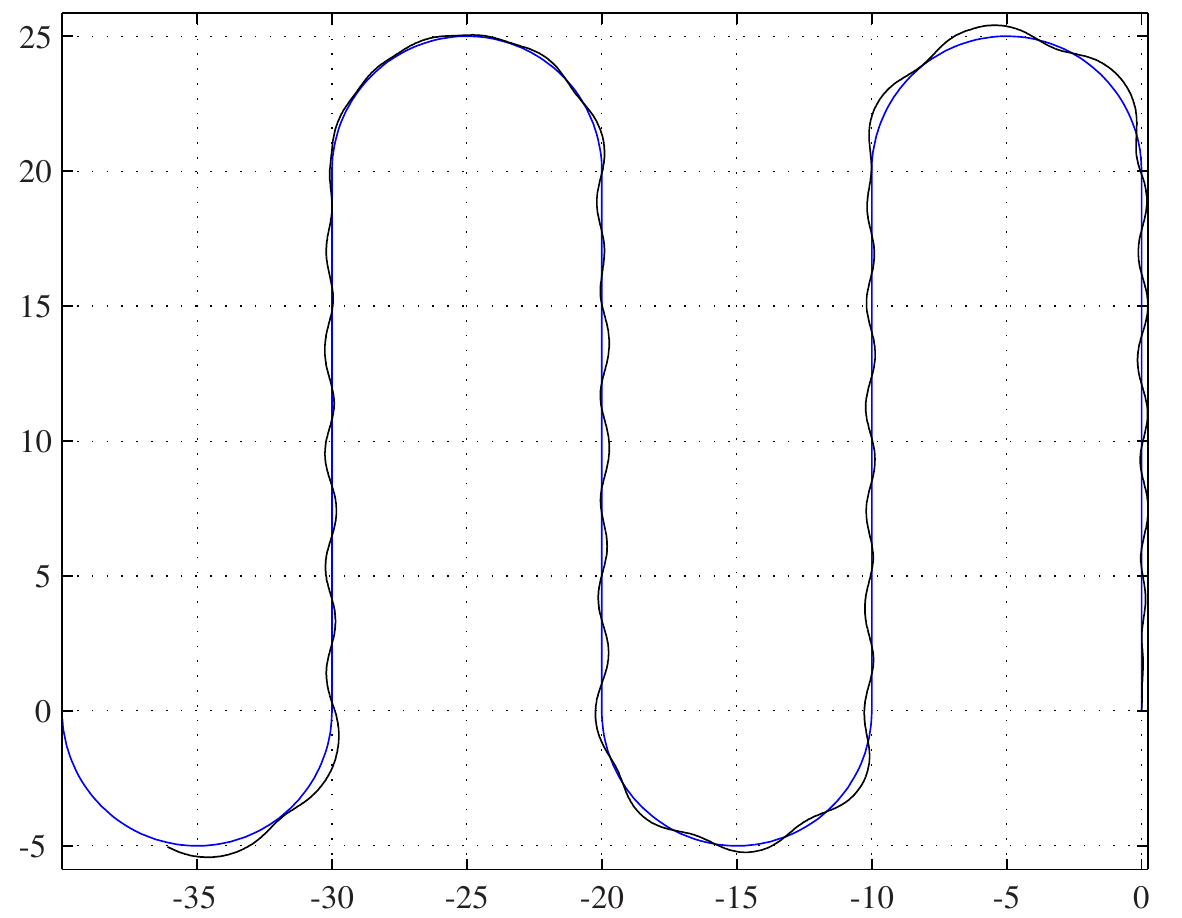}
  \caption{Simulations using the maximum actuation
    controller.}
  \label{pf:fig:trajsmax}
\end{figure}
\begin{figure}[ht]
  \centering \includegraphics[width=10cm]{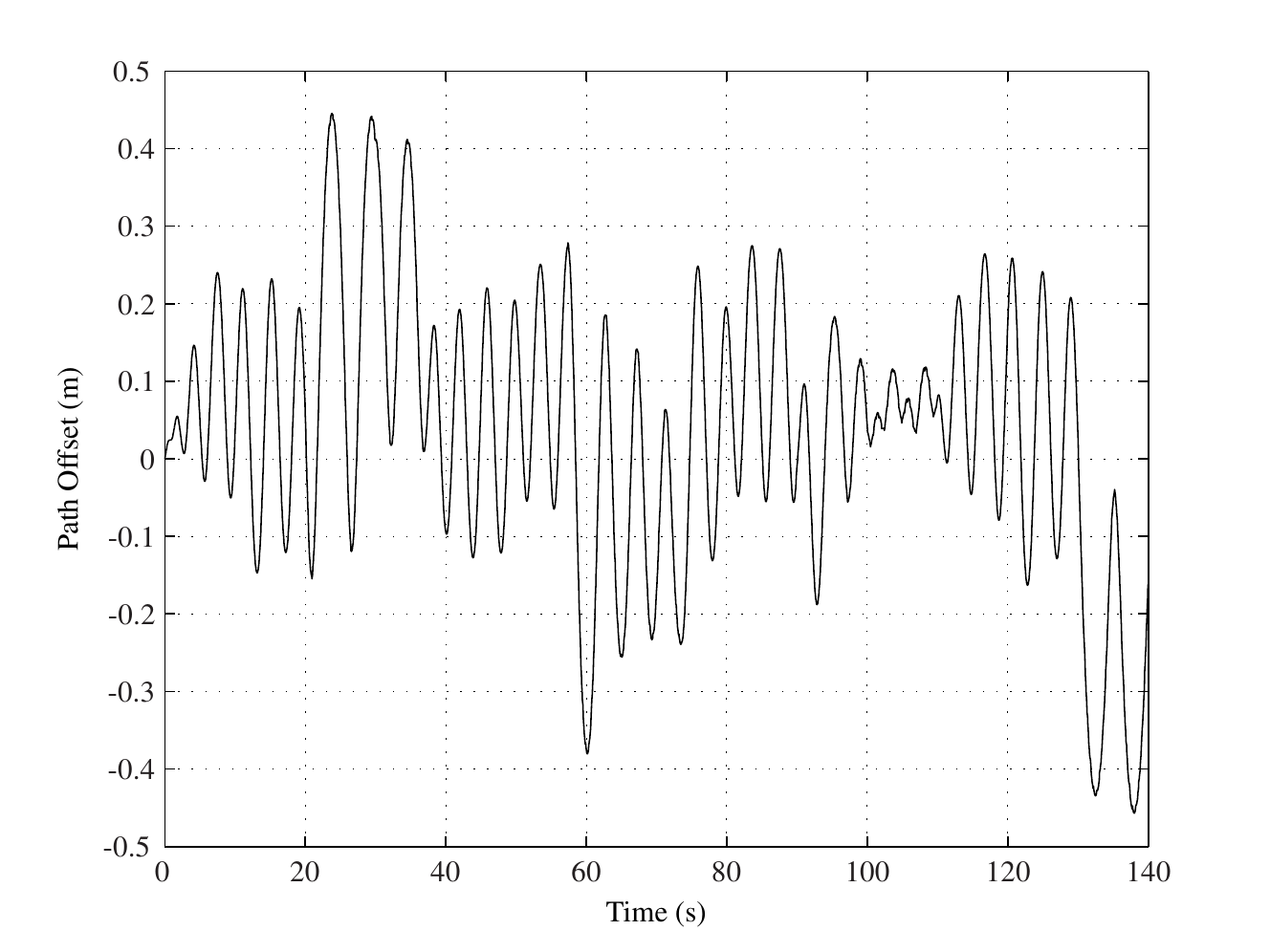}
  \caption{Offset error from the trajectory when using the maximum
    actuation controller.}
  \label{pf:fig:zerrorsmax}
\end{figure}
\begin{figure}[ht]
  \centering \includegraphics[width=10cm]{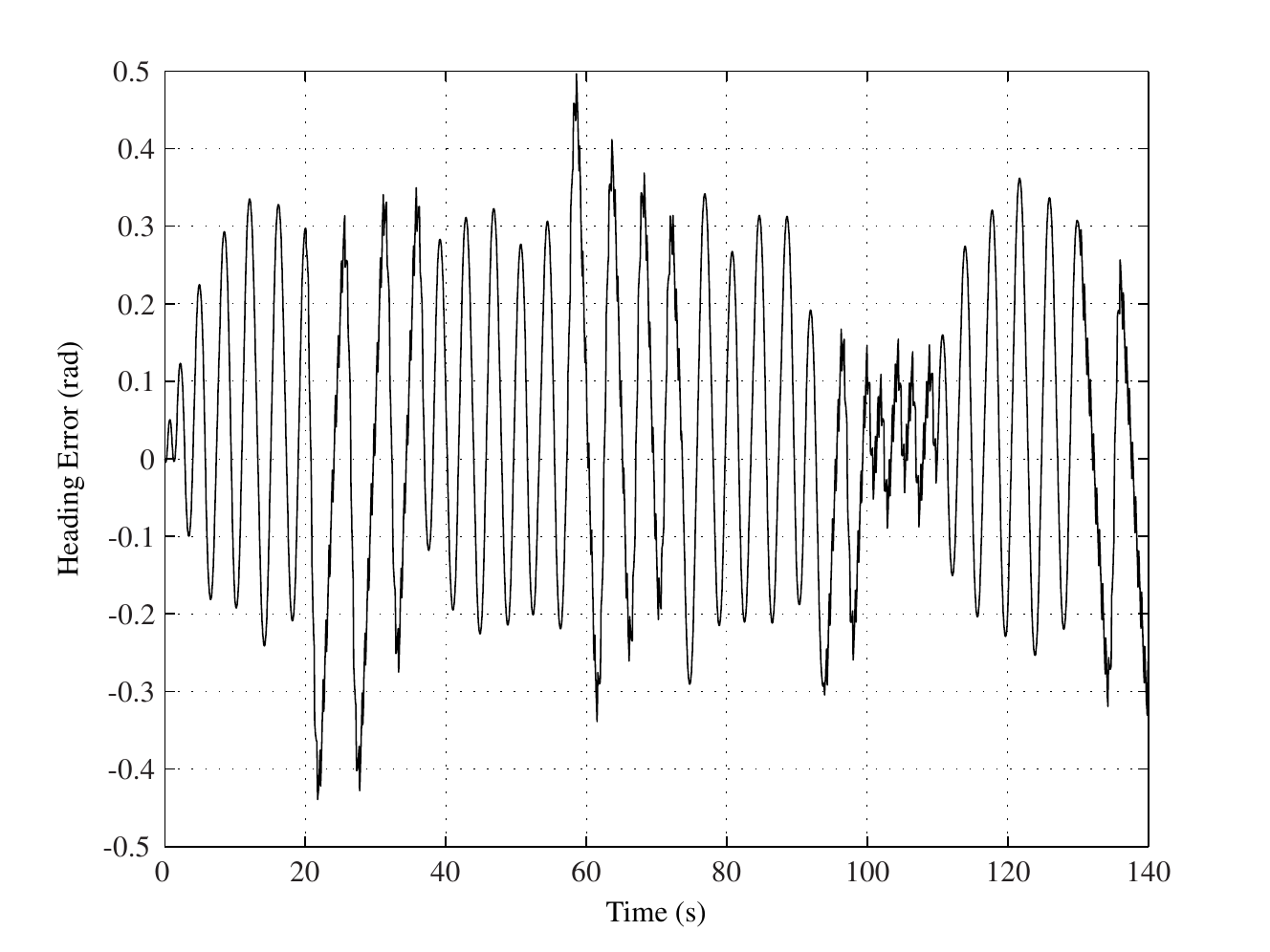}
  \caption{Heading error obtained when using the maximum actuation
    controller.}
  \label{pf:fig:therrorsmax}
\end{figure}

Typical results obtained for the simpler controller introduced in
Sec.~\ref{pf:sec.ma} are shown in Figs.~\ref{pf:fig:trajsmax},
\ref{pf:fig:zerrorsmax} and \ref{pf:fig:therrorsmax}. Even in
Fig.~\ref{pf:fig:trajsmax}, it is visible that the performance is
worse than for the controller from Sec.~\ref{pf:sec.ra}. This is
fleshed out by Figs.~\ref{pf:fig:zerrorsmax} and
\ref{pf:fig:therrorsmax}, which display systematic maximal spatial
error $\approx 45cm$ and angular error $\approx 0.4 rad$. As was
discovered via extended simulation tests, the degradation of
performance of the controller from Sec.~\ref{pf:sec.ma} is mostly
due to the extra un-modeled steering dynamics. In particular,
increasing the upper limit on the rate of steering $\dot{\delta}$
typically caused decay of the above errors. In the case where this
rate is not limited, the spatial tracking error typically did not
exceed $7cm$, which is even better than for the controller with the
reduced actuation. This provides an evidence that the controller from
Sec.~\ref{pf:sec.ma} can be viewed as a worthwhile option if the
vehicle is equipped with high-speed high-torque steering
servo. However, this is not the case for most autonomous agricultural
tractors nowadays due to a variety of reasons, including the hardware
cost.
\par
\begin{figure}[ht]
  \centering \includegraphics[width=10cm]{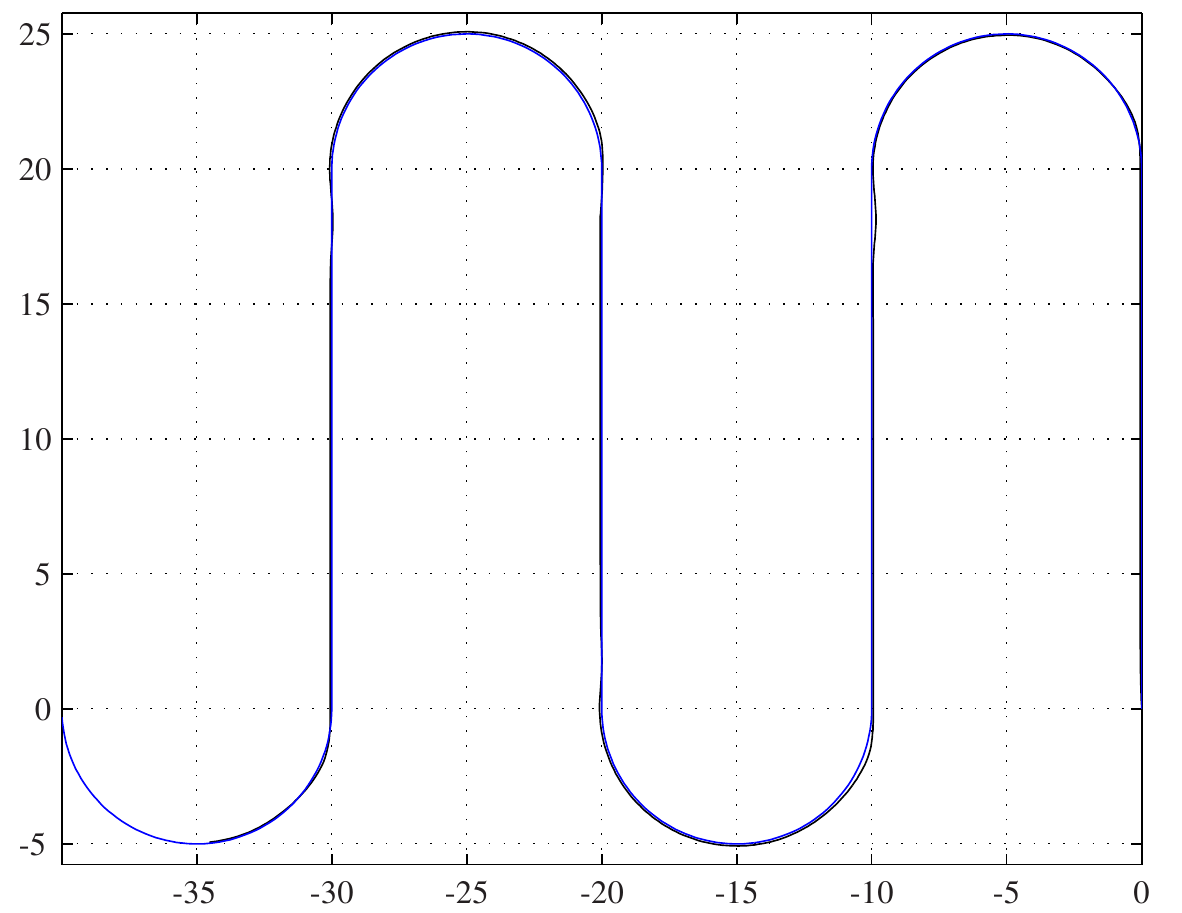}
  \caption{Simulations using the reduced actuation controller at
    high speed.}
  \label{pf:fig:trajhs}
\end{figure}
\begin{figure}[ht]
  \centering \includegraphics[width=10cm]{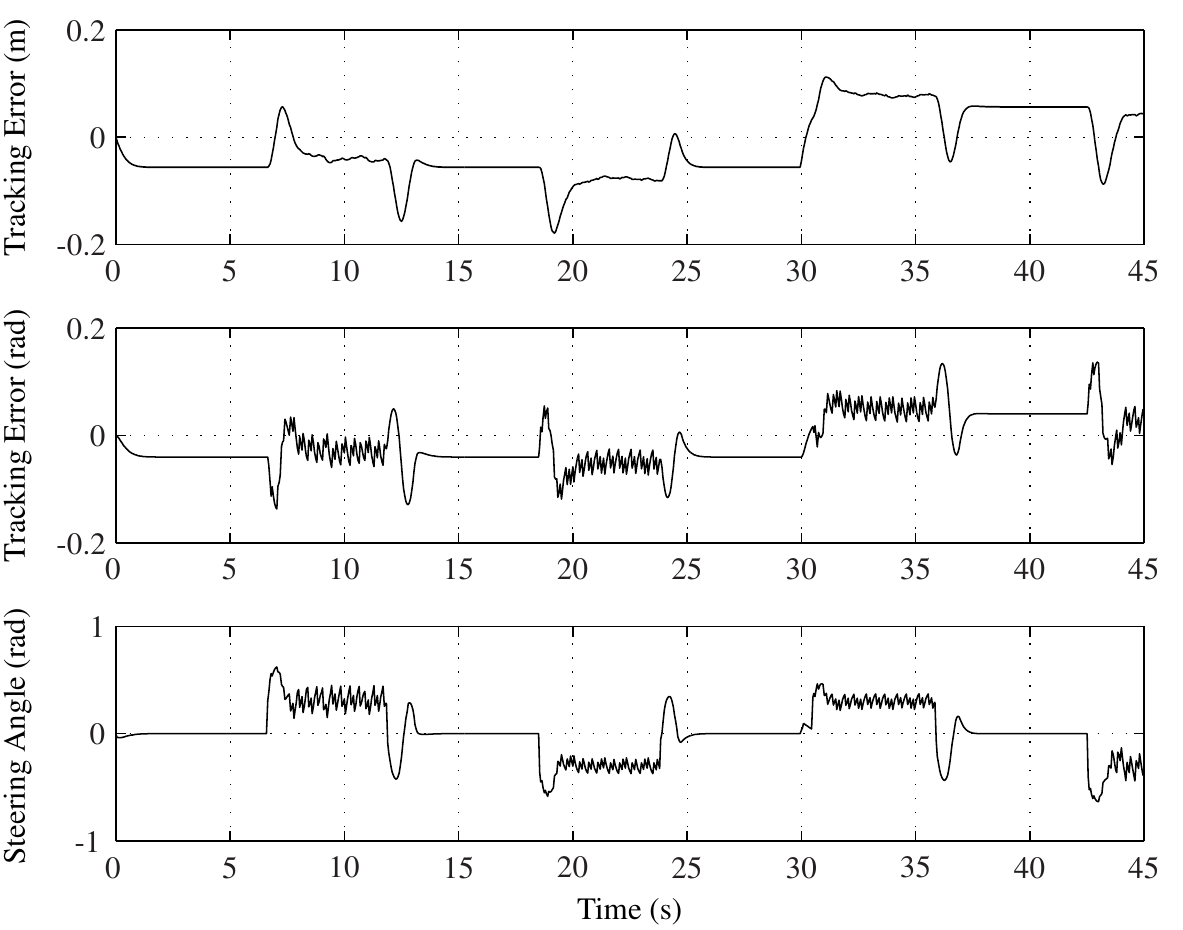}
  \caption{Performance of the reduced actuation controller at high
    speed.}
  \label{pf:fig:simhs}
\end{figure}

In the above simulation tests, the tractor speed was moderate. In the
next test, performance of the closed-loop system was examined in the
case of rather high (for agricultural vehicles) speed $\approx 10.8
km/h$. The parameters used for this simulation are given in
Table~\ref{pf:fig:paramsimhs1}. Typical results obtained for the
non-linear controller with reduced actuation (see
Sec.~\ref{pf:sec.ra}) are presented in Figs.~\ref{pf:fig:trajhs}
and
\ref{pf:fig:simhs}. 
The spatial tracking error does not exceed $18cm$, whereas the
orientation error almost always lies within $[-0.15rad,
0.15rad]$. This exhibits only slight performance degradation as
compared with operation at a moderate speed (see
Figs.~\ref{pf:fig:trajsmin}, \ref{pf:fig:zerrorsmin},
\ref{pf:fig:therrorsmin}, and \ref{pf:fig:steerangsmin}) and
demonstrates a good enough overall performance at both moderate and
high speeds.

\begin{table}[ht]
  \centering
  \begin{tabular}{| l | c |}
   \hline
    $v_{rw}$ & $3 ms^{-1}$ \\
    \hline
    $\delta_{max}$ & $0.7 rad$ \\
    \hline
    $L$ & $1.69 m$  \\
    \hline
\end{tabular}
\hspace{10pt}
 \begin{tabular}{| l | c |}
    \hline
    $\mu$ & $2.7 rad$\\
    \hline
    $\Delta$ & $2.2 m $\\
    \hline
\end{tabular}
\hspace{10pt}
 \begin{tabular}{| l | c |}
   \hline
    $\sigma$& $0.4 m$\\
    \hline
    $\unk$ & $0.01m^{-1}$\\
    \hline
  \end{tabular}
  \caption{Control parameters used during high speed simulation test.}
  \label{pf:fig:paramsimhs1}
\end{table}

In the discussed experiments, saturation of the steering angle rarely
occurred.  The next simulation test illuminates the controller
performance in the case of systematic visible control saturation. In
this experiment, the tractor was slightly slowed down and traveled at
the speed $v=0.7 ms^{-1}$, and the steering angle limit was reduced to
$\delta_{\max} = 0.45 rad \approx 30^{\circ}$, which is only $\approx
7^{\circ}$ grater than the value necessary to perfectly track the
path. The other parameters are shown in Table~\ref{pf:fig:paramsimls}.

\begin{table}[ht]
  \centering
  \begin{tabular}{| l | c |}
    \hline
    $v_{rw}$ & $0.7 ms^{-1}$ \\
    \hline
    $\delta_{max}$ & $0.45 rad$ \\
    \hline
    $L$ & $1.69 m$  \\
    \hline
\end{tabular}
\hspace{10pt}
 \begin{tabular}{| l | c |}
    \hline
    $\mu$ & $1.0 rad$\\
    \hline
    $\Delta$ & $0.4 m $\\
    \hline
\end{tabular}
\hspace{10pt}
 \begin{tabular}{| l | c |}
    \hline
    $\sigma$& $0.5 m$\\
    \hline
    $\unk$ & $0.01m^{-1}$\\
    \hline
  \end{tabular}
  \caption{Control parameters used during low speed simulation.}
  \label{pf:fig:paramsimls}
\end{table}

\begin{figure}[ht]
  \centering \includegraphics[width=10cm]{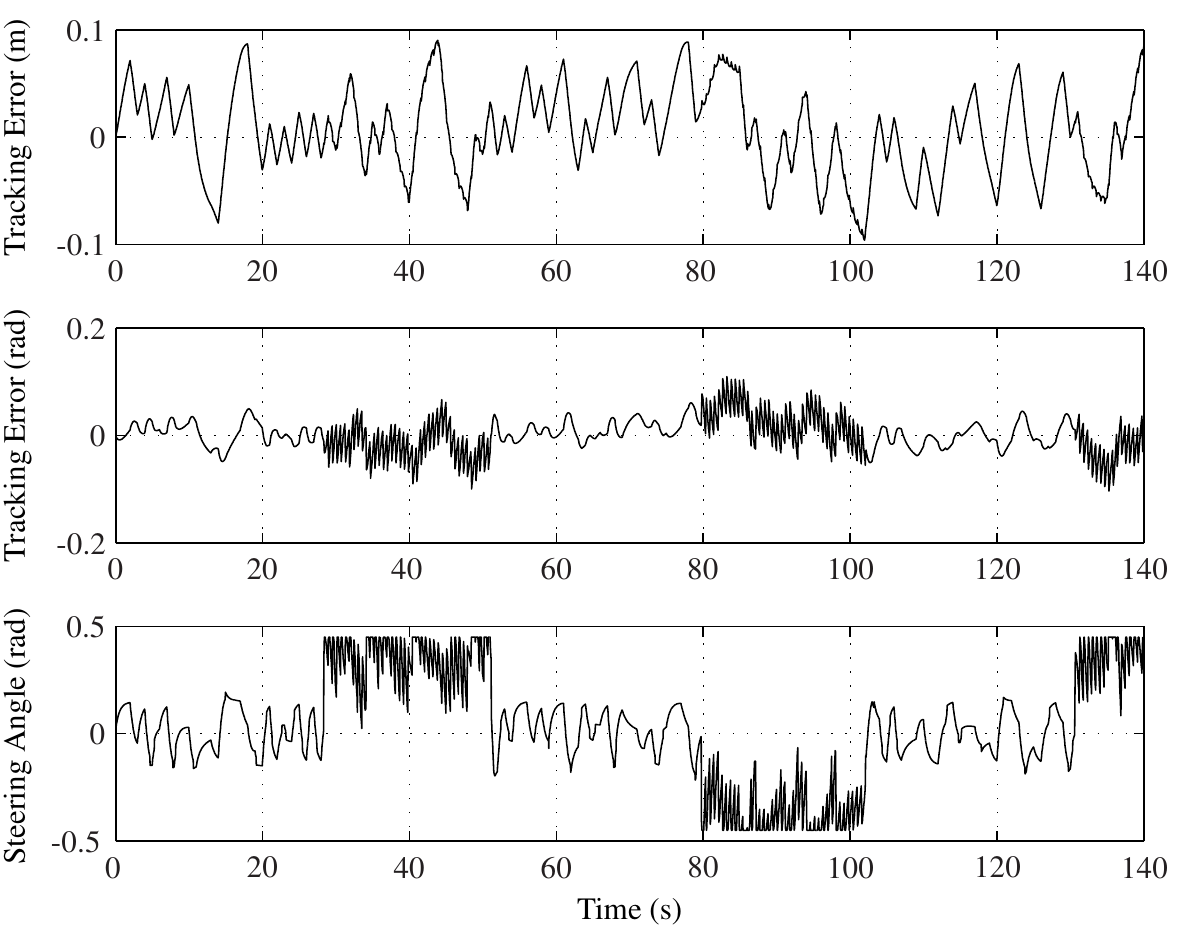}
  \caption{Offset errors and steering angle for the reduced
    actuation controller.}
  \label{pf:fig:zerrorsmin1}
\end{figure}

The trajectory did not visibly alter as compared with
Fig.~\ref{pf:fig:trajsmin}, whereas the offset errors and steering
angle are displayed in Fig.~\ref{pf:fig:zerrorsmin1}. In spite of
systematic saturation of the control signal, the controller
demonstrated good performance, with the spatial and angular errors
being no more than $9 cm$ and $0.1 rad$, respectively. As was found
out via extra simulation tests, slowing down the tractor makes the
closed-loop system more amenable to saturation. This is worthwhile of
attention since some agricultural operations are performed at very low
speeds.
\par
Thus the simulation tests have shown that both controllers behaved as
expected. However, the controller with reduced actuation displays a
significant advantage in terms of stability and performance in the
face of steering dynamics, which was the rationale for its
introduction in Sec.~\ref{pf:sec.ra}.

For a comparative study, the adaptive and predictive controller
from \cite{EJC07} was chosen, which is among those with the best reported
performances in off-road conditions. The simulation setup was
basically borrowed from the previous subsection. The parameters of the
adaptive controller were taken from \cite{EJC07}; the sideslip
observer was not implemented to provide a more distilled comparison of
the control laws themselves.\footnote{Such an observer can be
  commissioned to aid both compared control laws.} The tractor
wheelbase and steering dynamics were taken from \cite{EJC07}, through
the tractor speed was different; the speed is given in
Table~\ref{pf:fig:comparesim}.

\begin{table}[ht]
  \centering
  \begin{tabular}{| l | c |}
    \hline
    $v_{rw}$ & $0.7 ms^{-1}$ \\
    \hline
    $\delta_{max}$ & $0.7 rad$ \\
    \hline
    $L$ & $2.75 m$  \\
    \hline
\end{tabular}
\hspace{10pt}
 \begin{tabular}{| l | c |}
    \hline
    $\mu$ & $1.1 rad$\\
    \hline
    $\Delta$ & $0.2 m $\\
    \hline
\end{tabular}
\hspace{10pt}
 \begin{tabular}{| l | c |}
    \hline
    $\sigma$& $0.5 m$\\
    \hline
    $\unk$ & $0.01m^{-1}$\\
    \hline
  \end{tabular}
  \caption{Control parameters used during comparison test.}
  \label{pf:fig:comparesim}
\end{table}

\begin{figure}[ht]
  \centering \includegraphics[width=10cm]{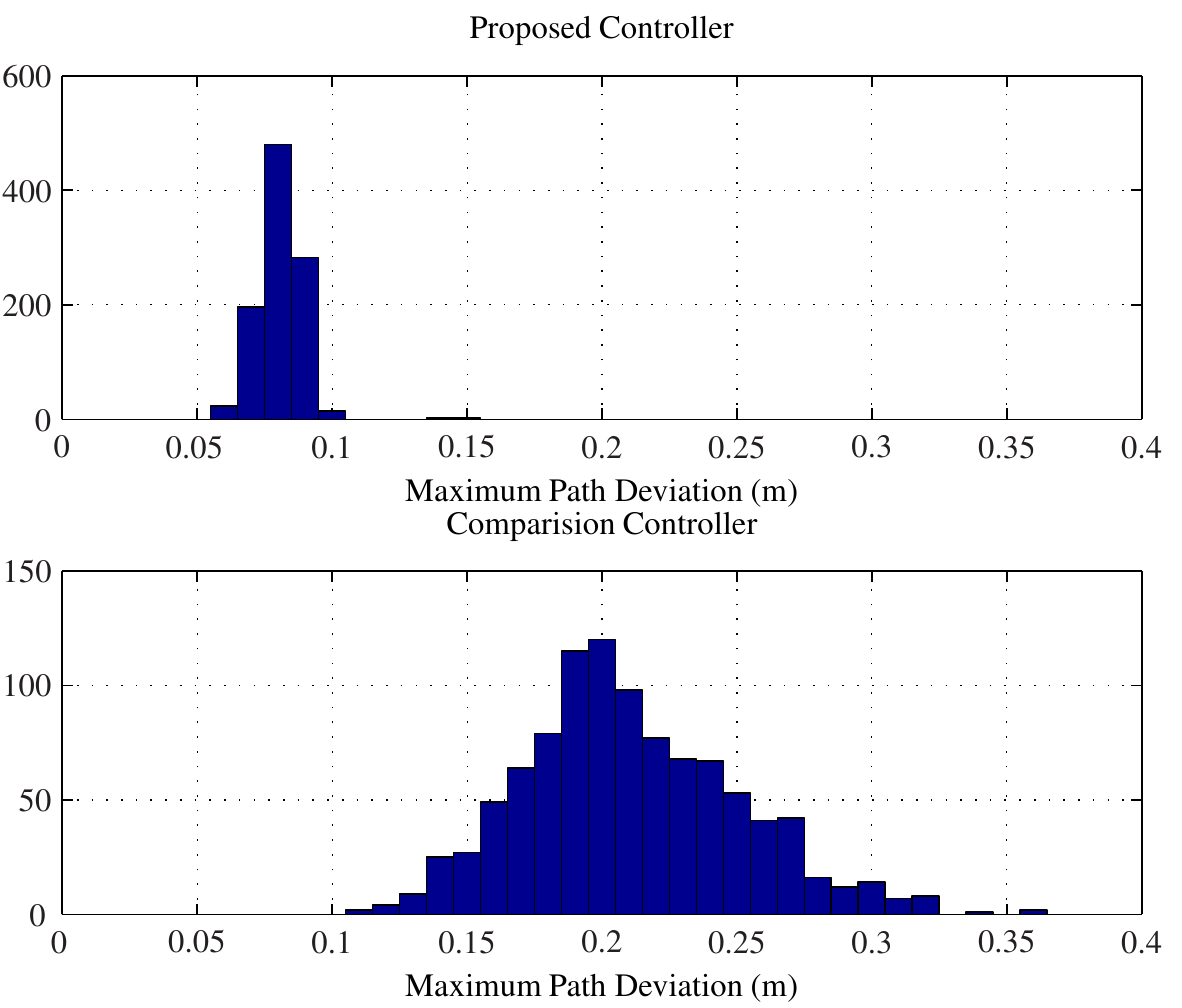}
  \caption{Comparative distribution of the maximal path deviation.}
  \label{pf:fig:hist}
\end{figure}

\begin{figure}[ht]
  \centering \includegraphics[width=10cm]{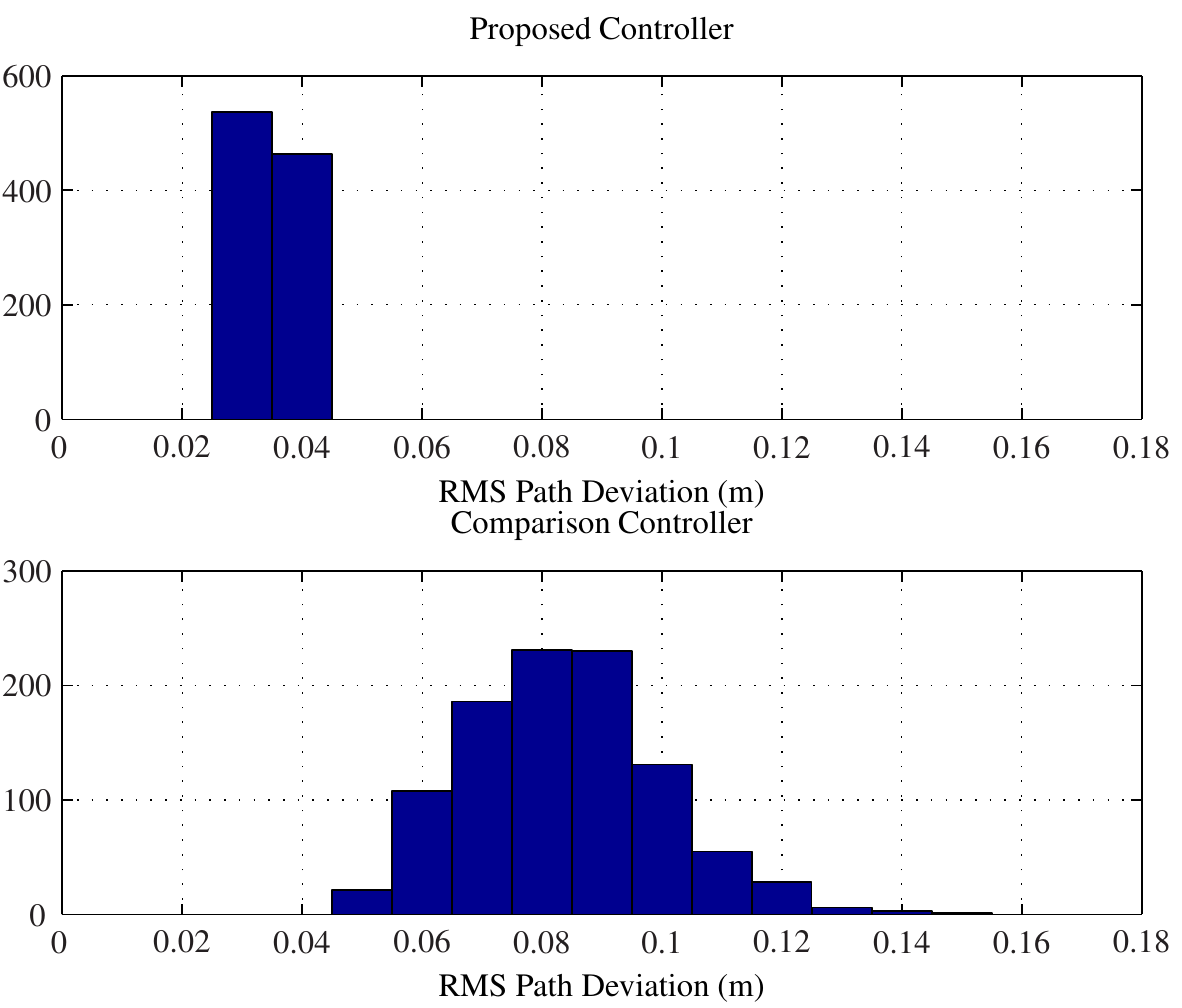}
  \caption{Comparative distribution of RMS path deviation.}
  \label{pf:fig:histrms}
\end{figure}

The two compared laws were tested against identical sequences of
random slip parameters and sensor noises; $1000$ sequences were
randomly generated and used in respective tests. Since the main
concern of this chapter is about systematic occasional degradation of
the tracking accuracy, the maximal spatial deviation from the
reference path was chosen as the basic performance index for a given
test.  However, the root mean square (RMS) deviation was also taken
into account.  The diagrams in Figs.~\ref{pf:fig:hist} and
\ref{pf:fig:histrms} show the number of tests (the vertical axis) with
a given value of the performance index (the horizontal axis). It
follows that the proposed controller provides $\approx 50\%$ of
performance improvement on average, at least in the considered
specific circumstances.

 \clearpage \section{Specifications of the Agricultural Vehicle}

\begin{figure}[ht]
   \centering
   \includegraphics[width=0.5\columnwidth]{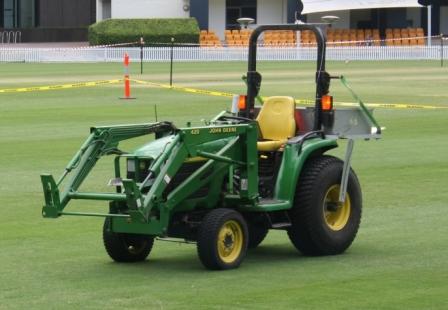}
   \caption{Autonomous tractor used for testing.}
   \label{pf:fig:photo}
 \end{figure}

 The controller was implemented on a fully autonomous compact
agricultural tractor (see Fig.~\ref{pf:fig:photo}) developed at the
University of New South Wales, Sydney, Australia. The tractor is fully
custom instrumented and automated with the integration of a complete
sensor suite and the accompanying software \cite{20113014186103}. The base tractor was a
John Deere 4210 Compact Utility Tractor \cite{tractor}.

 The control inputs to the tractor are the steering control signal and
 the propulsion control signal. The tractor's propulsion system is
 driven by a hydrostatic transmission system, which allows the control
 of the speed through the control of its swash plate. The swash plate
 can be controlled through electronic means by directly applying the
 required voltage by interfacing to the built-in computer of the
 tractor. This facilitated a non-invasive means of controlling the
 vehicle speed. Two separate analog voltages, generated by the
 on-board computer (note that the built-in and on-board are two
 different computers), one for forward motion and another for reverse
 motion are used to control the forward and reverse speeds of the
 tractor. These two voltage channels correspond to the forward pedal
 and reverse pedal on the tractor. All other logic to take care of
 unacceptable scenarios such as both pedals being pressed at the same
 time or a pedal voltage exceeding its range are provided by the
 built-in computer of the tractor.

 The automation of the steering however did not have the same ease of
 interfacing. Hence it was necessary to mount an additional actuator
 on the steering system.  There was no positive coupling between the
 steered position of the wheels and the steering column angular
 position. Hence the mounting of a separate steering sensor was
 necessary. A non-linear relationship exists between the steering
 sensor position and the actual steered angle of the front wheels. A
 standard well tuned PI controller was implemented to control the
 steering. The desired steering is specified as the desired steering
 sensor position. Such a command can easily be generated to correspond
 a desired steered angle through a simple table look-up. The steering
 control system operates very well with a speed of response an order
 higher than the speed of response of the tractor motion.

 Within limits (as determined by the manually controlled gear
 position), the speed of the tractor can be automatically
 controlled. The speed is obtained by averaging the rate of change of
 the encoder counts of the left and right hand rear wheel
 encoders. The rear wheels have an external diameter of 0.9 m. and the
 encoders generate 40000 counts per revolution. The on-board
 computer's time stamp is used to calculate the actual speed in m/s.

 \subsection{Safety Subsystem}

 The safety subsystem is centered around the built-in computer's halt
 mode that is wired to the seat switch. A simple re-triggerable one
 shot timer circuit mimics the seat switch signal, also known as the
 watch-dog signal, and when active selects the propulsion signals
 generated by the on-board computer instead of those generated by the
 pedals. The re-triggering of the timer must take place every 128
 ms. or less, or else the watchdog signal will become inactive thereby
 switching all propulsion signals to the manual pedals. This requires
 the pedals to be pressed afresh to effect motion. If a pedal is not
 pressed or a driver is not present the tractor will be halted.

 The control software is written in such a manner that if any of the
 critical software modules fail then the watchdog will be disabled. In
 addition an operator at the remote console may also choose to disable
 the watchdog signal.

 \subsection{Sensors}

 The tractor is equipped with a rich set of sensors. To measure its
 global position and orientation two GPS systems, both with 2 cm. RTK
 accuracy are
 available. 
 Other sensors include an encoder each for the two rear wheels and a
 steering position sensor sensing the displacement of the steering
 system's power cylinder. For communication between individual
 on-board components and the remote station outside the tractor, an
 Ethernet network is established and connected to an access point with
 a high gain antenna capable of 300m line-of-sight Wi-Fi
 communication.

 \subsection{Software}

 \begin{figure}
   \centering
   \includegraphics[width=4in]{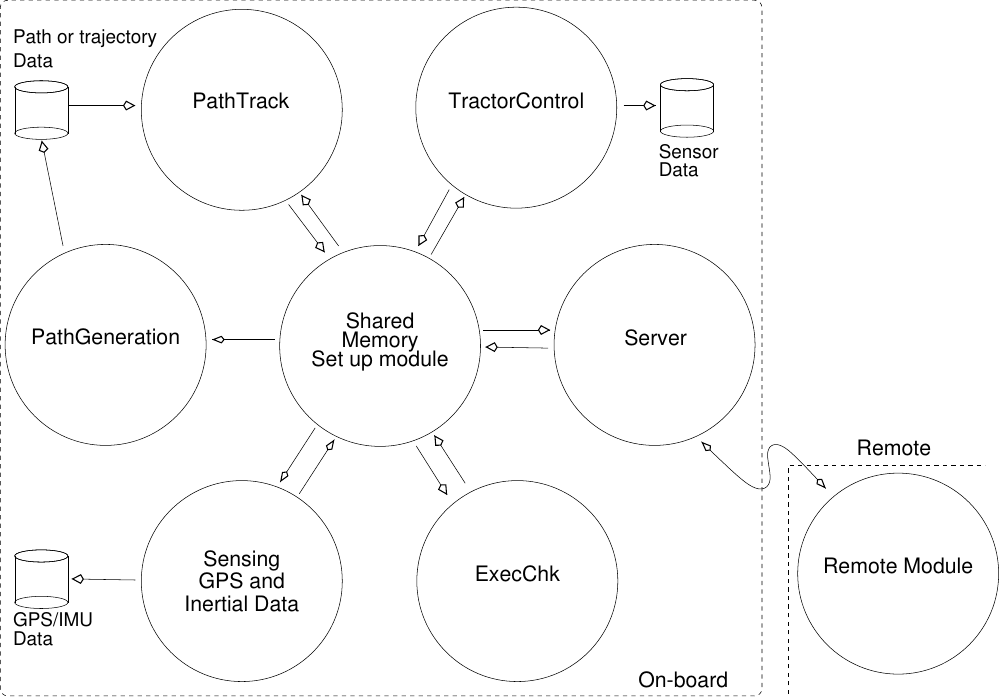}
   \caption{Software architecture.}\label{pf:fig:software}
 \end{figure}

 A complete suite of software is put in place to take care of a
 variety of tasks. In addition to the on-board software modules, an
 additional module must run on a remote client. The data communication
 between the external client and the tractor's on-board computer takes
 place via the earlier mentioned Wi-Fi network. The remote module
 predominantly control the watchdog signal in addition to the mode
 selection and remote control.

 The software modules are shown in Fig.~\ref{pf:fig:software}. Each
 circle represents a software module. Each link shows the data flow
 and the cylinders represent disk storage. The mode of operation of
 the tractor is chosen by the remote module through buttons on a
 joystick attached to the remote computer. The tractor can operate in
 one of three modes, (i) Path generation, (ii) Remote control and
 (iii) Autonomous control. Depending on the mode of operation various
 different software modules move onto the critical path. Based on the
 mode chosen the \texttt{ExecChk} module sequentially executes the
 required modules and monitors the live execution of all required
 modules adhering to the stipulated precedence. If any of the critical
 modules fail, the \texttt{ExecChk} module will initiate a
 shutdown. If the \texttt{ExecChk} itself fails, it will be detected
 by all other modules. In that situation, the \texttt{TractorControl}
 module will disable the watchdog signal. As an additional safety
 measure, the human operator stands guard of the watchdog signal at
 the remote station.

 The \texttt{PathGen} module can generate a path of a predefined shape
 to align with the tractor's current position and
 orientation. Obviously, for this to function the GPS readings must be
 available hence the GPS module must run before \texttt{PathGen} can
 be run. The \texttt{SharedMem} module facilitates data exchange
 between modules. The \texttt{Server} module maintain contact with the
 remote module and receives the mode select, watchdog signal and
 remote control commands. \texttt{Sensing} module acquires GPS and
 inertial data. Two separate sub-modules, one for GPS and one for the
 IMU are used within it due to their different data rates. The
 \texttt{PathTrack} module executes high level control algorithms. The
 \texttt{TractorControl} module carries out the low level control of
 the steering and propulsion while recording all low level sensor
 data.

  \section{Experiments}
 \label{pf:sec:expresult}

 The dimensions of the experimental test trajectory were identical to
 the simulated scenario. The controller parameters for the experiment
 are shown in Table~\ref{pf:fig:paramexp}. The average speed of the
 tractor during the test was approximately $1.9 km/h$. No filtering
 was done on the raw GPS measurements, which is acceptable due the
 high noise resistance of sliding mode control laws (the main type of
 error observed in the results is not attributable to random noise,
 but rather systematic offsets). The GPS receivers were evenly and
 equally spaced $0.230 m$ behind the center of the rear axle, and the
 position and orientation of the rear axle was calculated directly
 based on this assumption.
 \par
 The steering angle requested by the controller was used as the
 set-point for the steering controller specified in the previous
 section. The steering dynamics can be approximated by a rate limit,
 and while it was not accurately characterized a full transition of
 the steering between the steering limits was observed to take up to
 one second (dependent on vehicle speed). The mismatch between the
 actual steering dynamics and those assumed during the the theoretical
 development were found to not cause malfunction of the control system
 provided the tractor speed was not too large -- for lower speeds, the
 delay in changing the steering angle has a lesser effect on the
 tractors state. More information about possible solutions to this
 problem are discussed in \cite{20101612850657}.

\begin{figure}[ht]
  \centering \includegraphics[width=10cm]{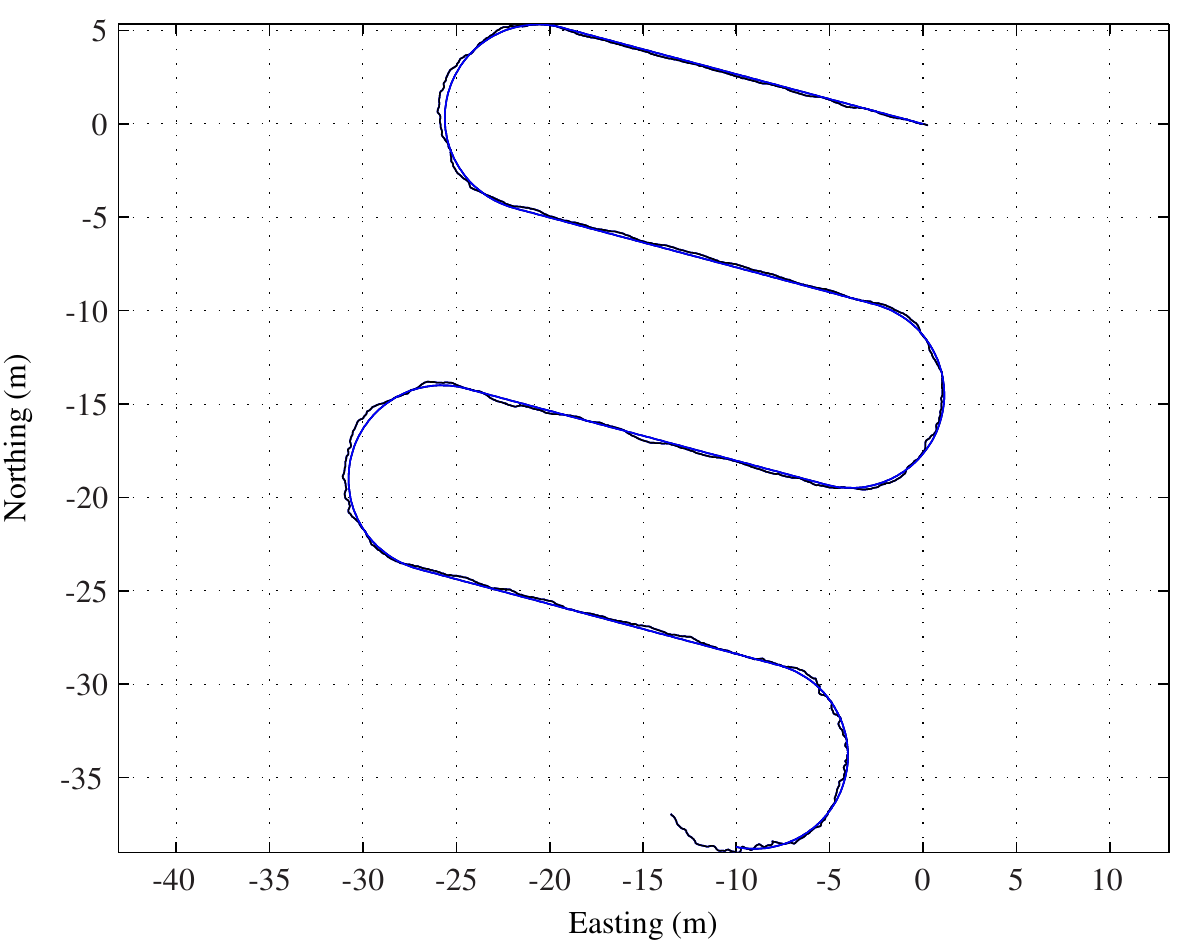}
  \caption{Trajectory obtained during the experiment.}
  \label{pf:fig:traj}
\end{figure}

\begin{figure}[ht]
  \centering \includegraphics[width=10cm]{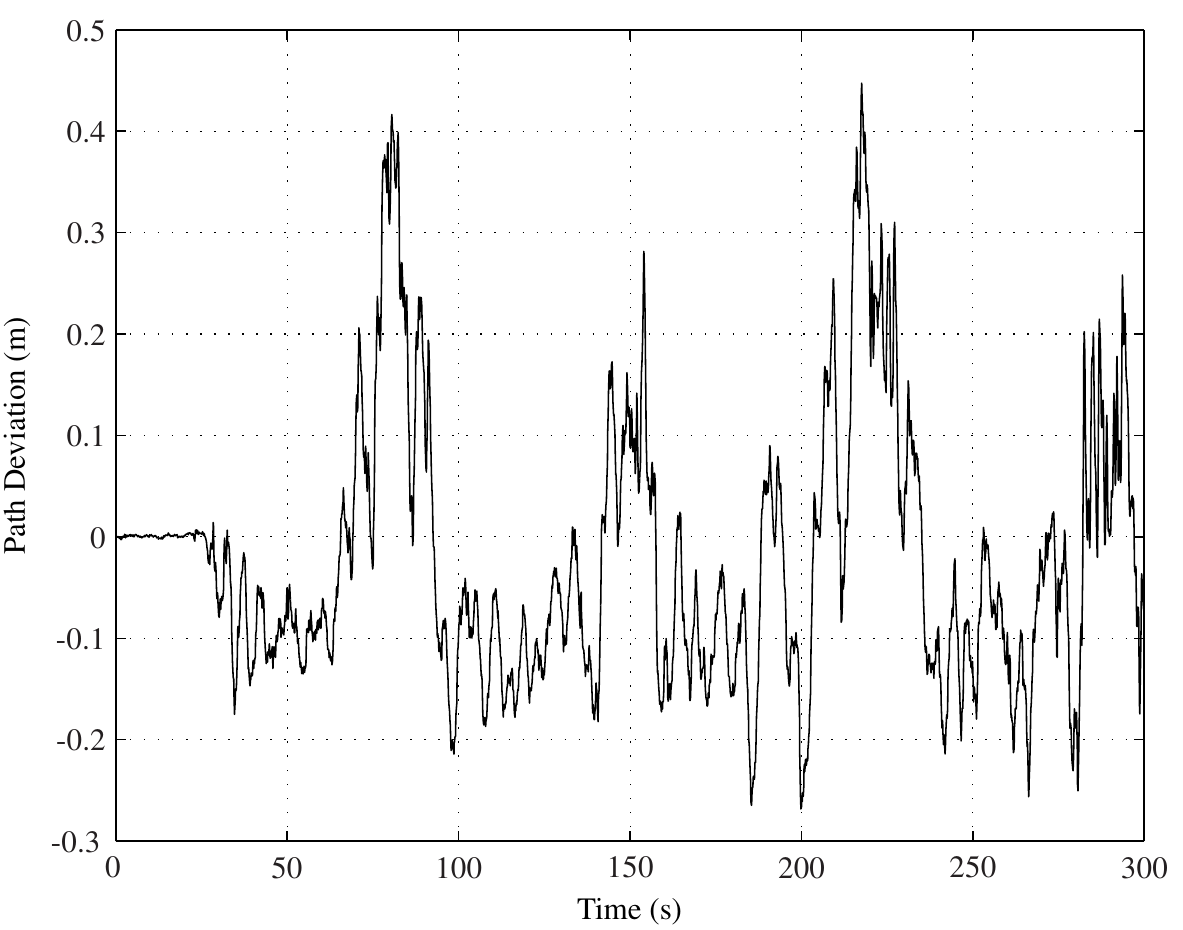}
  \caption{Path offset error obtained during the experiment.}
  \label{pf:fig:zerror}
\end{figure}

\begin{figure}[ht]
  \centering \includegraphics[width=10cm]{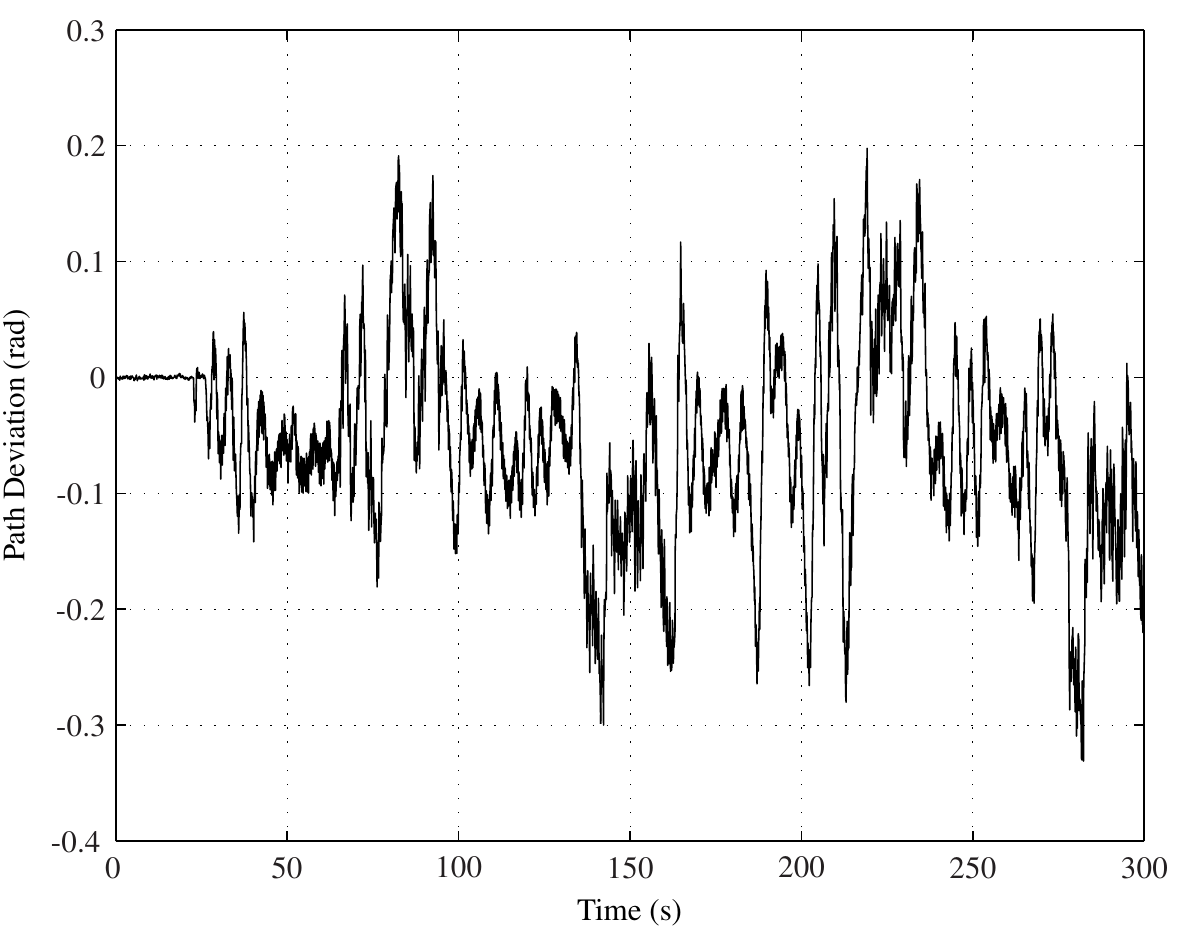}
  \caption{Heading error obtained during the experiment.}
  \label{pf:fig:therror}
\end{figure}

\begin{figure}[ht]
  \centering \includegraphics[width=10cm]{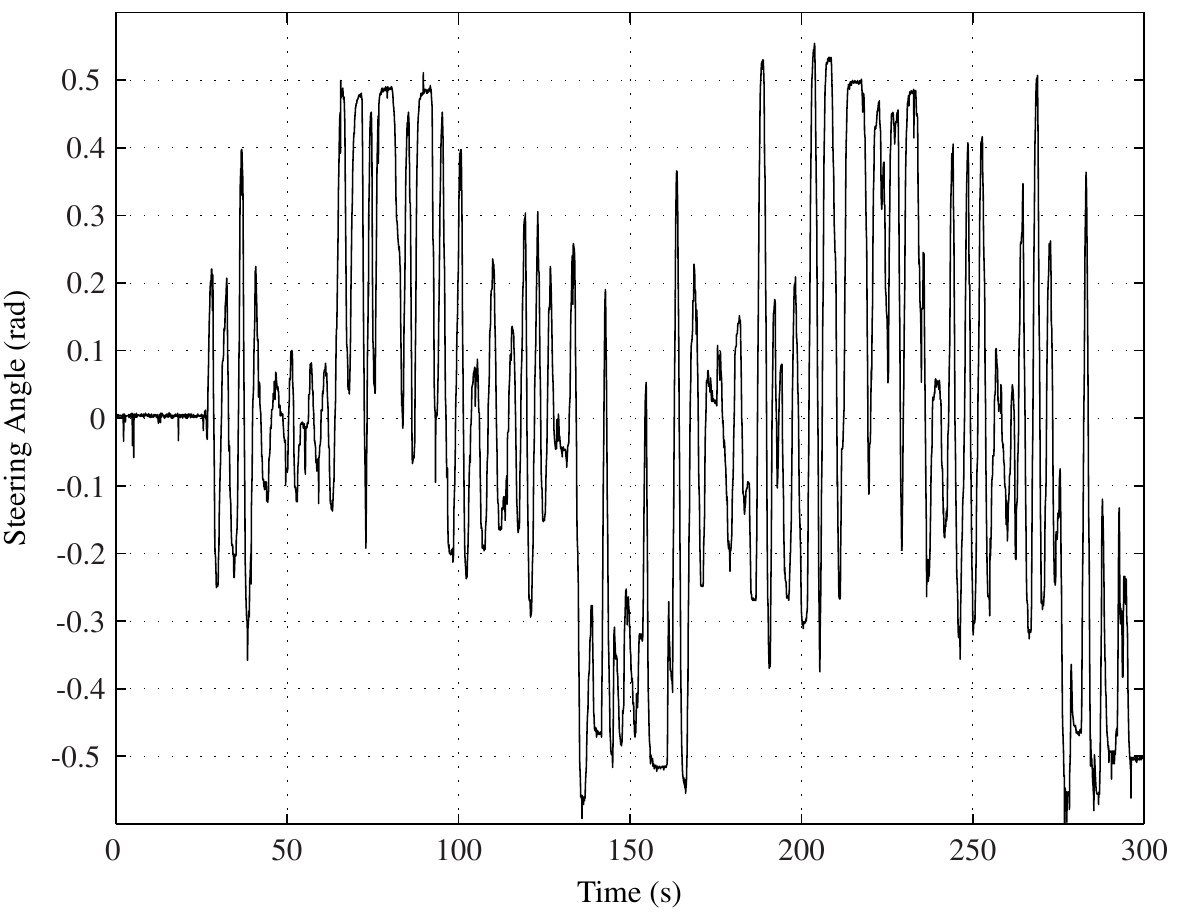}
  \caption{Actual steering angle of the tractor obtained during the
    experiment.}
  \label{pf:fig:steerang}
\end{figure}

 The results of the tractor test are shown in Figures
 \ref{pf:fig:traj}, \ref{pf:fig:zerror}, \ref{pf:fig:therror} and
 \ref{pf:fig:steerang}, where it may be seen the tractor behaved as
 expected. In Fig.~\ref{pf:fig:steerang}, the actual steering angle
 of the tractor is shown. It follows that the system operates near the
 mechanical steering angle limit, which is approximately equal to $0.6
 rad$.

 The highest measured path offset error obtained was approximately
 $0.4 m$, whereas the RMS error was $0.3199m$ and the RMS heading
 error was $0.2918rad$. There are some systematic errors when tracking
 straight path segments, and these are most likely caused by
 systematic errors in the measurement system. They could likely be
 eliminated by introducing some type of adaptive observer, however
 that is outside the scope of this chapter.

\begin{table}[ht]
  \centering
  \begin{tabular}{| l | c |}
    \hline
    $\delta_{max}$ & $0.5 rad$ \\
    \hline
    $L$ & $1.69 m$  \\
    \hline
    $\mu$ & $0.3 rad$\\
    \hline
  \end{tabular}
\hspace{10pt}
  \begin{tabular}{| l | c |}
    \hline
    $\Delta$ & $0.15 m$\\
    \hline
    $\sigma$& $0.5 m$\\
    \hline
    $\unk$ & $0.05 m^{-1}$\\
    \hline
  \end{tabular}
  \caption{Control parameters used for experiments.}
  \label{pf:fig:paramexp}
\end{table}

The experiment was undertaken on a very rough field, with the
amplitude of undulation up to $15.0cm$. Though direct comparison with
other controllers under similar circumstances is highly troublesome
since every "rough" field has its own individual features (which can
be hardly described and reproduced in details), the observed tracking
error is comparable with the best results reported in the literature.

\begin{figure}[ht]
  \centering \includegraphics[width=10cm]{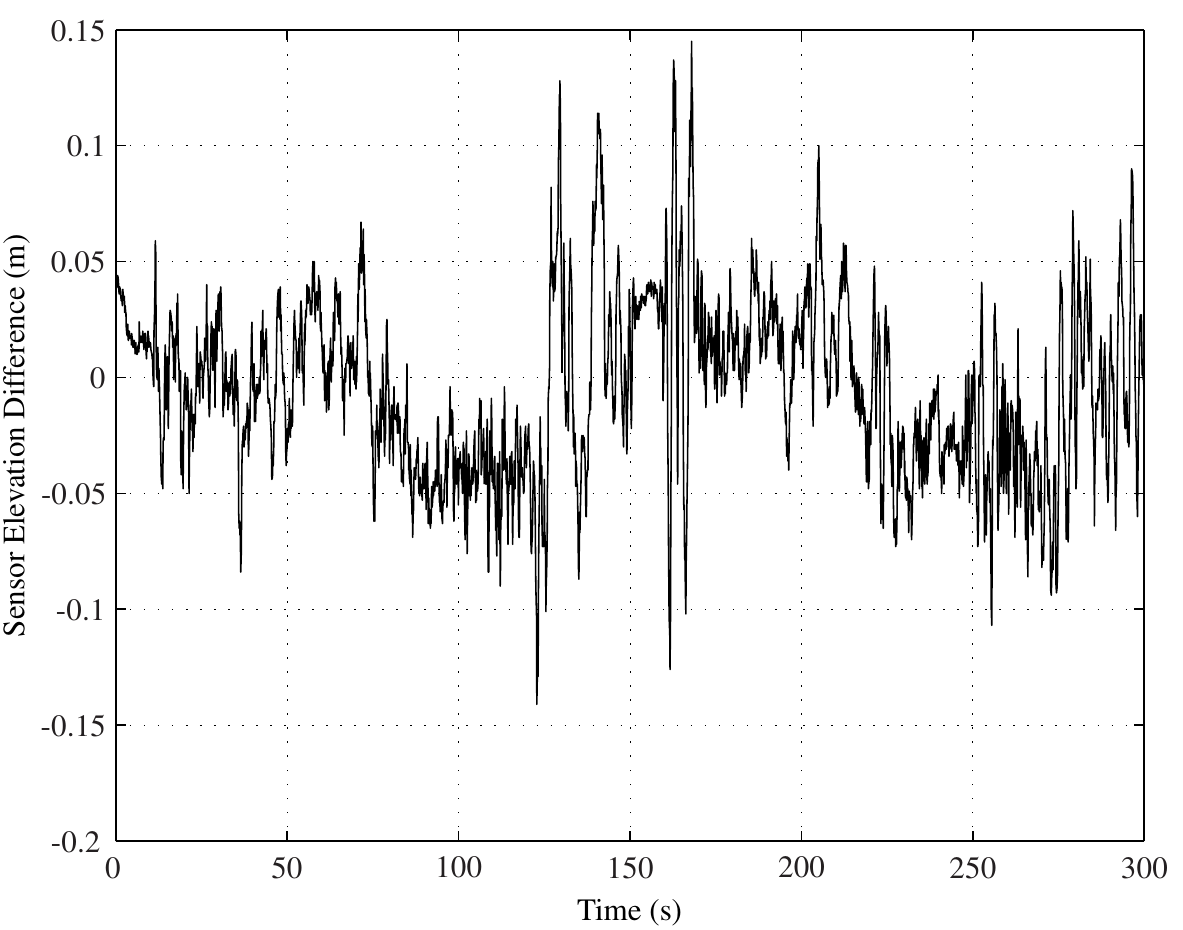}
  \caption{Difference in measured elevation between the two GPS
    sensors.}
  \label{pf:fig:hdif}
\end{figure}

\clearpage \section{Summary}
\label{pf:con}

Two approaches for control of farm tractors in the presence of
slippage are presented. One of them uses maximal actuation of the
vehicle, whereas fluctuation of the control input is reduced with the
aid of a smooth nonlinear control law for the other. Simulation
studies showed that the nonlinear controller produced more stable
trajectories with lower control efforts as compared to the pure
sliding mode controller. Experiments
were also carried out on a real agricultural vehicle, which confirmed
the real world viability of the proposed control system.

\chapter{Boundary Following using Minimal Information}
\label{chapt:fbf} 

The boundary following method proposed in Chapt.~\ref{chap:convsingle} relies on
numerous obstacle detections in order to successfully track the boundary
of an obstacle. In this chapter, an alternative method is proposed, which
is based on only
the distance along and the reflection angle of the ray perpendicular
to the vehicle centerline. Such a situation holds if the
measurements are supplied by several range sensors rigidly mounted to
the vehicle body at nearly right angles from its centerline, or by a
single sensor scanning a similarly narrow sector.  This
perception scheme is used in some applications to reduce the
complexity, cost, weight, and energy consumption of the sensor system
and to minimize detrimental effects of mechanical external disturbances
on the measurements.
\par
This approach retains many of the advantages of the method proposed in
Chapt.~\ref{chap:convsingle}, such as provably correct behaviour, and no requirement for
measurement of the boundary curvature. The main difference is that a
equidistant curve from the boundary is tracked at constant speed, whereas the 
method in Chapt.~\ref{chap:convsingle} exhibits no fixed speed and offset from the boundary
(which may make it unsuitable
for certain applications). In addition, the method proposed in this chapter
has very low computational requirements, and can be expressed as a
simple sliding mode control law. The main disadvantage is that global
assumptions are present which place bounds on the obstacle curvature
(however these are unavoidable given the reduced amount of information
available). Also, deficit of sensor data makes
most of known navigation solutions 
inapplicable. This gives rise to special
challenges, like inability to detect a threat of head-on collision in
certain situations (see Fig.~\ref{fbf:fig.srct}).
\par
This chapter proposes a novel sliding-mode navigation strategy that
does not employ curvature estimates and homogeneously handles both
concavities and convexities of the followed boundary, as well as
transitions between them. This strategy asymptotically steers the
vehicle to the pre-specified distance to the boundary and afterwards
ensures stable maintenance of this distance. In addition, mathematically rigorous
justification results for non-local convergence of the proposed strategy are
available. In doing so, possible abrupt jumps of the sensor readings are
taken into account. Furthermore, much attention is given to revealing
requirements to the global geometry of the boundary that make it
possible to avoid front-end collisions with it based on only side view
sensors, thus making extra front-view sensors superfluous. The
convergence and performance of the proposed navigation and guidance
law are confirmed by computer simulations and real world tests with a
Pioneer P3-DX robot, equipped with a SICK LMS-200 LiDAR sensor.
\par
All proofs of mathematical statements are omitted here; they are
available in the original manuscript \cite{MaHoSa11ar}.
\par
The body of this chapter is organized as follows. In Sec.~\ref{fbf:sec.pst} the problem is formally defined, and in Sec.~\ref{fbf:sec.ass} the 
main assumptions are described. A summary of the main results are in Sec.~\ref{fbf:sec.mr}. Simulations and
experiments are presented in Secs.~\ref{fbf:sec.sim} and \ref{fbf:sec.expr}. Finally, brief conclusions are given
in Sec.~\ref{fbf:sec.concl}.

 \section{Problem Statement}
\label{fbf:sec.pst}
\begin{figure}[ht]
  \centering
  \subfigure[]{\scalebox{0.45}{\includegraphics{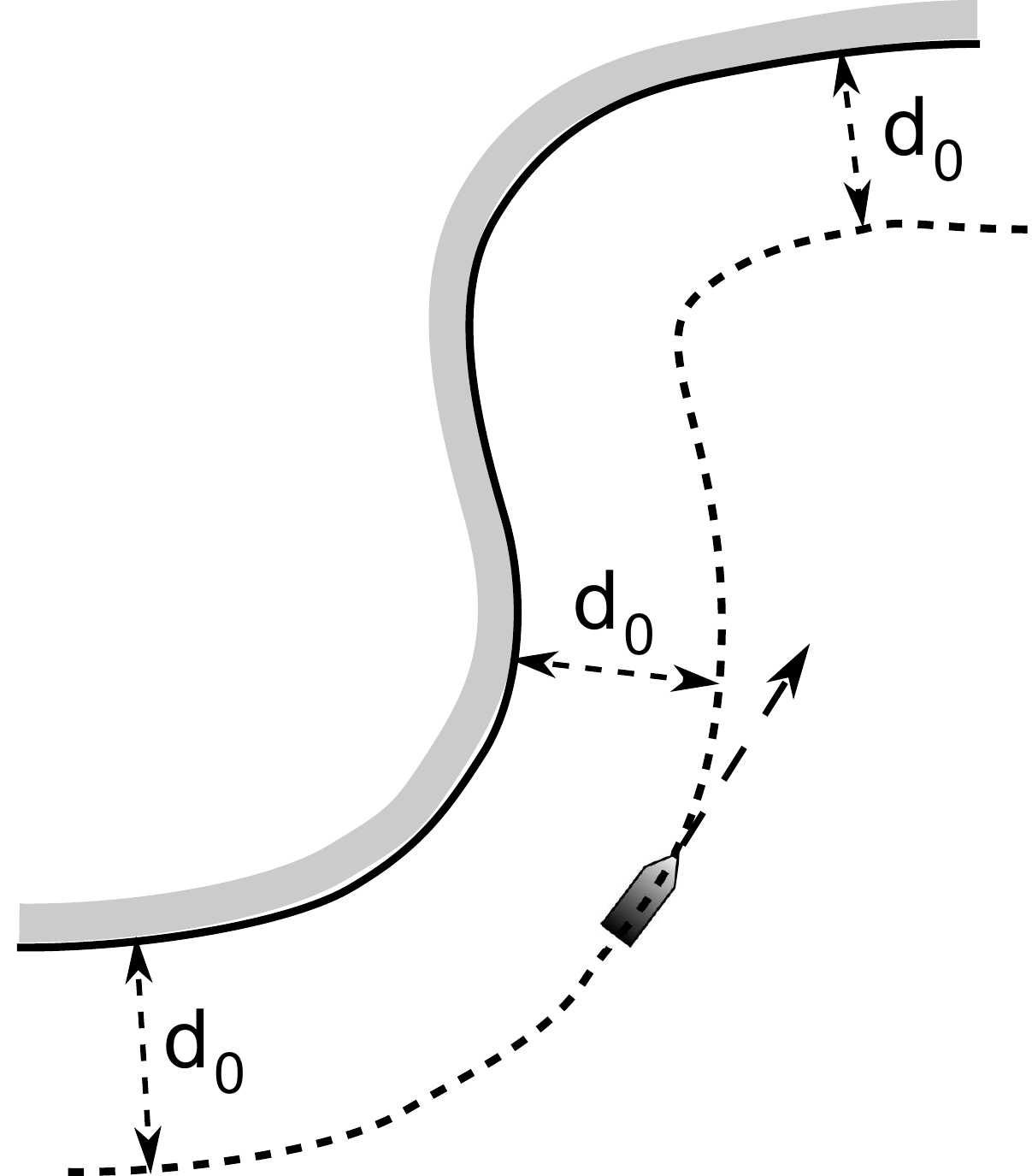}}}
\hspace{20pt}
  \subfigure[]{\scalebox{0.45}{\includegraphics{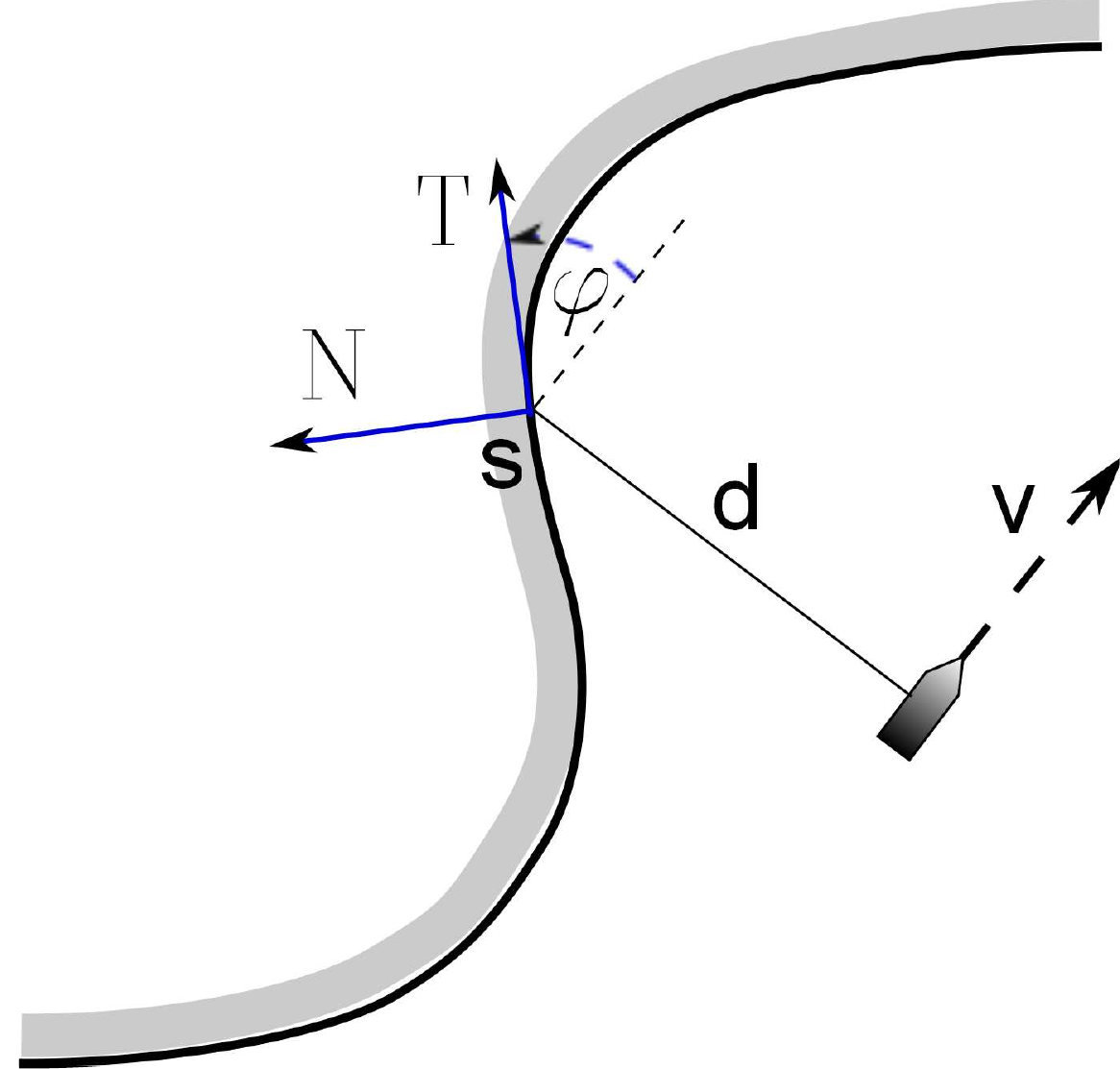}}} \caption{(a)
    Motion over the equidistant curve; (b) Vehicle with a rigidly
    mounted range sensor.} \label{fbf:fig.1}
\end{figure}

A Dubins-type vehicle travels in the plane with the constant speed
$v$. It is controlled by the angular velocity $u$ limited by a given
constant $\ov{u}$. There also is a domain $D$ with a smooth boundary
$\partial D$ in the plane. The objective is to drive the vehicle over
the equidistant curve of the domain $D$ separated from it by the
pre-specified distance $d_0$ (see Fig.~\ref{fbf:fig.1}(a)). 
The vehicle is equipped with a narrow-aperture range sensor
directed perpendicularly to the vehicle centerline and to the
left. This sensor provides the distance $d$ from the vehicle to the
nearest point of $D$ in the sensed direction (see
Fig.~\ref{fbf:fig.1}(b)).

The scan within the aperture and processing of the collected data
provides the vehicle with access to the angle $\varphi$ from its
forward centerline to the tangential direction of the boundary at the
reflection point. Whenever the sensor does not detect the obstacle,
$d:= \infty, \varphi := 0$.
\par
Apart from stable maintenance of the motion over the equidistant
curve, it is required to ensure transition to this motion from a given
initial state. In doing so, the vehicle must not collide with the
boundary $\partial D$. The distance from the boundary if defined 
as follows:

\begin{equation}
  \label{fbf:true.dist}
  \text{\bf dist}\,[\blr,D]:= \min_{\blr^\prime \in D} \|\blr-\blr^\prime\|
\end{equation}

The distance from the vehicle location $\blr$ to the boundary should constantly
exceed the given safety margin $d_- < d_0$.
\par
The kinematics of the considered vehicles are classically described by
the following equations:

\begin{equation}
  \label{fbf:1}
  \begin{array}{l}
    \dot{x} = v \cos \theta,
    \\
    \dot{y} = v \sin \theta,

  \end{array}, \quad\dot{\theta} = u \in [-\overline{u},
  \overline{u}], \quad
  \begin{array}{l}
    \bldr(0) = \bldr_0 \not\in D
    \\
    \theta(0) = \theta_0
  \end{array} .
\end{equation}

\begin{figure}[ht]
  \centering
  \subfigure[]{\scalebox{0.3}{\includegraphics{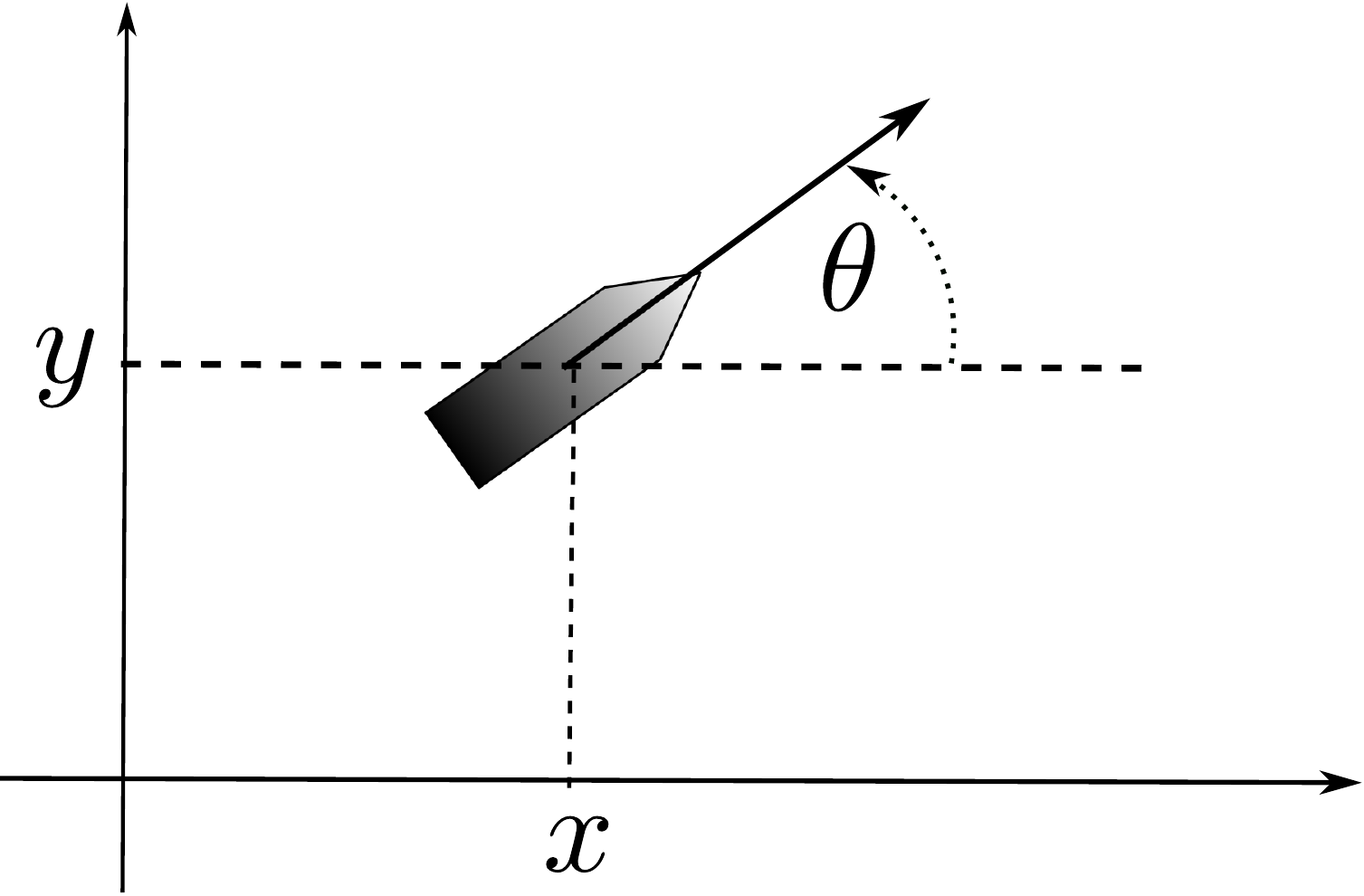}}}
\hspace{20pt}
  \subfigure[]{\scalebox{0.3}{\includegraphics{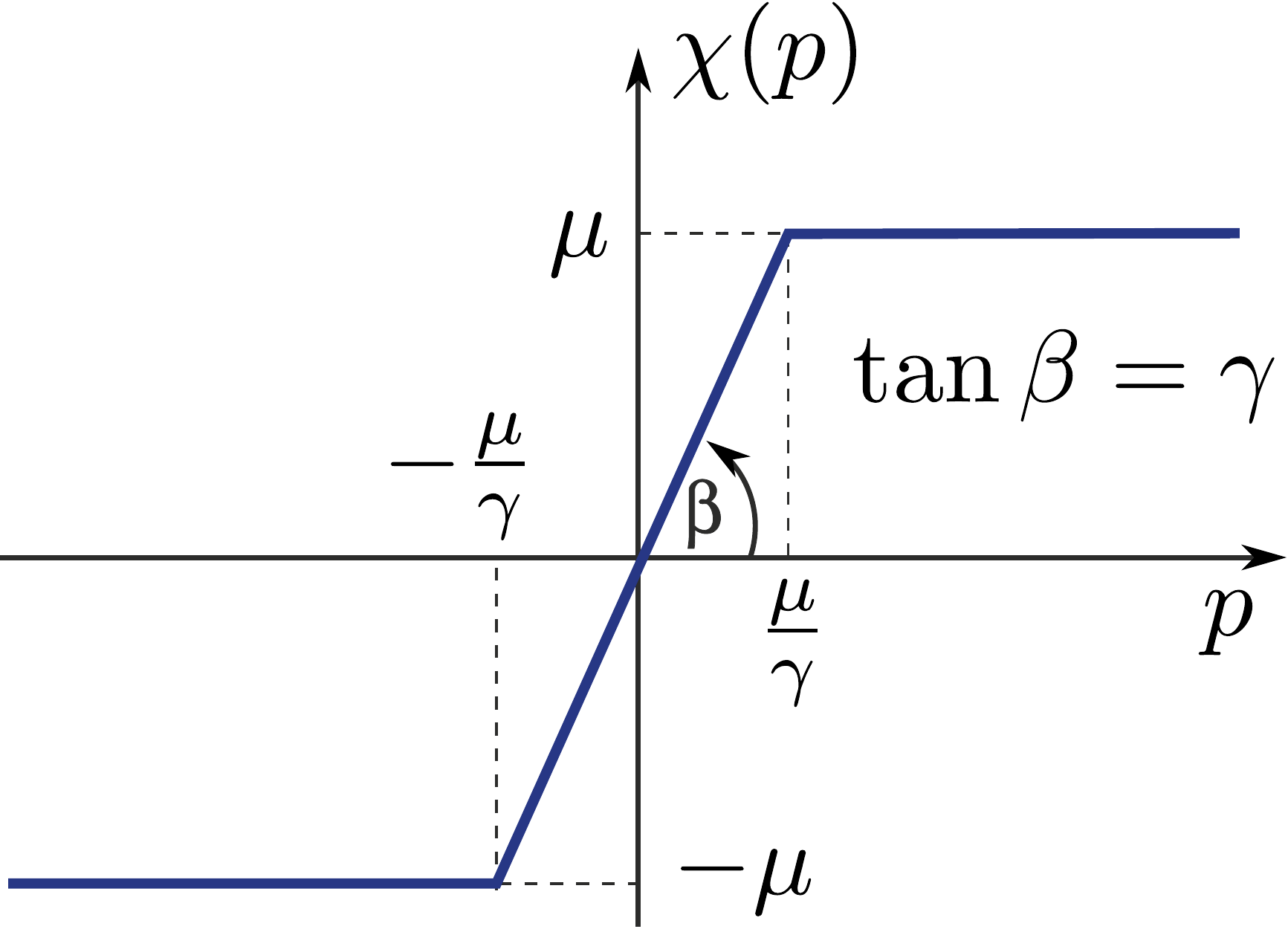}}} \caption{(a)
    Planar vehicle; (b) Linear function with
    saturation.} \label{fbf:fig.str}
\end{figure}

Here $x,y$ are the coordinates of the vehicle in the Cartesian world
frame, whereas $\theta$ gives its orientation, the angular velocity
$u$ is the control parameter (see Fig.~\ref{fbf:fig.str}(a)). (The set
of initial states $\bldr_0 \not\in D$ for which the problem has a
solution will be specified later on in Theorem~\ref{fbf:th.m}.)

 In the case at hand, the minimal turning radius is given by:
 
\begin{equation}
  \label{fbf:Rmin} R= v/\overline{u}.
\end{equation}

In this chapter, the following navigation law is examined:

\begin{equation}\label{fbf:c.a}
  u = \left\{
    \begin{array}{l}
      - \ov{u} \; \text{if} \; S:= \varphi + \chi[d-d_0] \leq 0
      \\
      u_+:= \min \left\{\ov{u}; vd^{-1}\right\} \; \text{if} \; S >0
    \end{array} \right. .
\end{equation}

Here $\chi(\cdot)$ is a linear function with saturation:

\begin{equation}
  \label{fbf:chi} \chi(p):= \begin{cases} \gamma p & \text{if} \; |p| \leq
    \mu/\gamma
    \\
    \sgn (p) \mu   & \text{otherwise}
  \end{cases}.
\end{equation}

The gain coefficient $\gamma >0$ and the saturation level $\mu \in
\left(0, \frac{\pi}{2}\right)$ are design parameters (see
Fig.~\ref{fbf:fig.str}(b)).

 \section{Main Assumptions}
\label{fbf:sec.ass}

For the control objective to be achievable, the vehicle should be
capable of tracking the $d_0$-equidistant curve of the boundary
$\partial D$. However this is impossible if this curve contains cusp
singularities, so far as any path of the unicycle Eq.\eqref{fbf:1} is
everywhere smooth. Such singularities are typically born whenever the
boundary contains concavities and the required distance $d_0$ exceeds
the critical value, which is equal to the minimal curvature radius of
the concavity parts of the boundary \cite{Arnold93}. Moreover, even if
there are no singularities, the equidistant curve should not be much
contorted since the robot is able to trace only curves whose curvature
radius exceeds Eq.\eqref{fbf:Rmin}. These observations are detailed in the
conditions necessary for the $d_0$-equidistant curve to be trackable
by the robot that are given by Lemmas~3.1 and 3.2 from
\cite{Matveev2011journ2}. Being slightly enhanced by putting the
uniformly strict inequality sign in place of the non-strict one, they
come to the following nearly unavoidable assumption:

\begin{Assumption}
  \label{fbf:ass.prelim}
  The following inequalities hold:
  \begin{equation}
    \label{fbf:mecess}
    R_\varkappa^+(D):= \inf_{\blr \in \partial D : \varkappa(\blr)> 0}
    R_\kappa > R - d_0,
\quad
    R_\varkappa^-(D):= \inf_{\blr \in \partial D : \varkappa(\blr) <
      0} R_\kappa > d_0 + R
  \end{equation}
  and the curvature is bounded $K:=\displaystyle{\sup_{\blr
      \in \partial D}} |\varkappa(\blr)| < \infty$.
\end{Assumption}

Here $\varkappa=\varkappa(\blr)$ is the signed curvature of $\partial
D$ at $\blr$, $R_\varkappa:= |\varkappa|^{-1}$ is the curvature
radius, and $\inf$ over the empty set is defined to be $+\infty$.  The
signed curvature is non-negative on convexities of the boundary and
negative on concavities.
\par
\begin{figure}[ht]
  \centerline{\scalebox{0.5}{\includegraphics{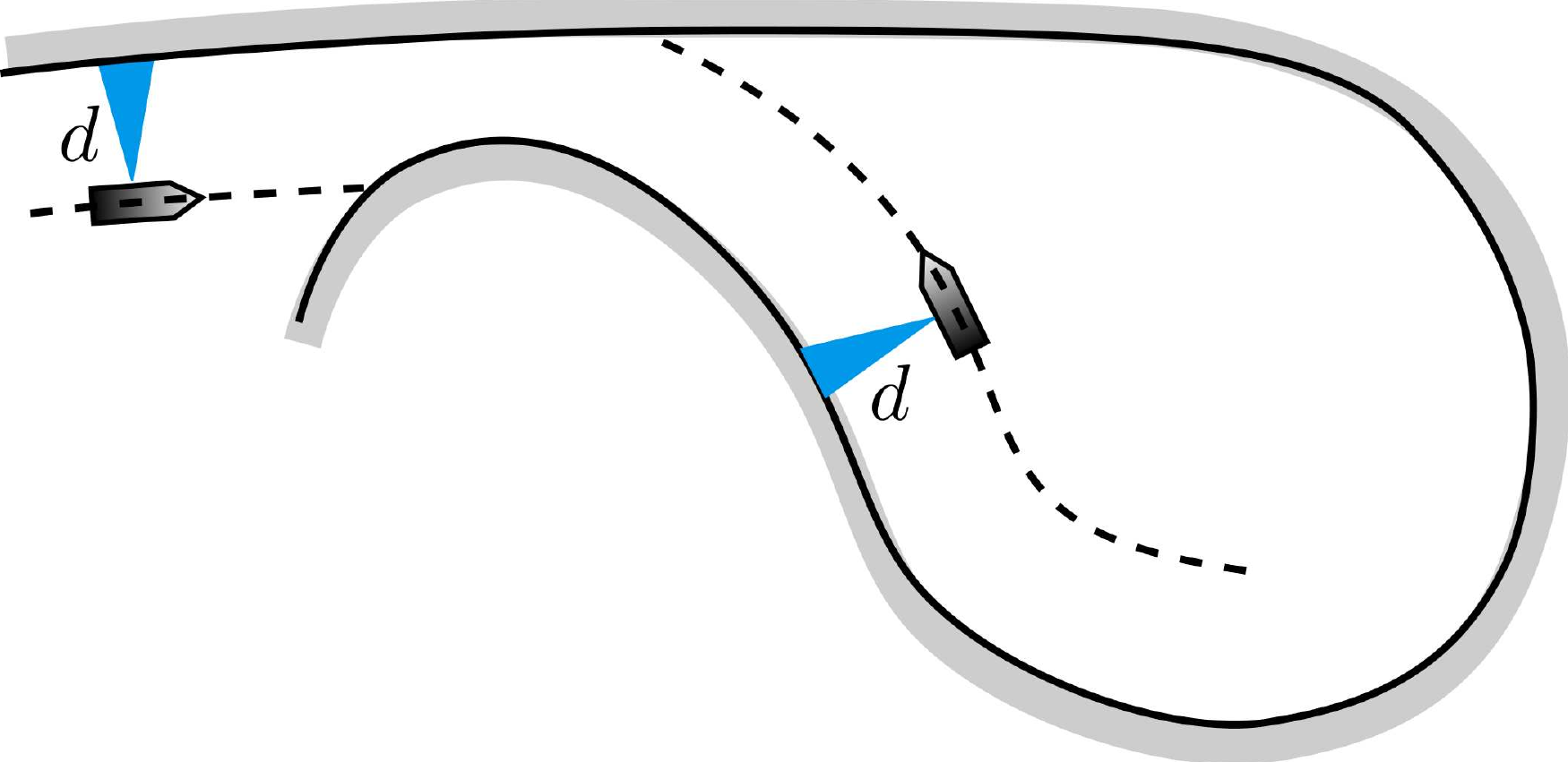}}}
  \caption{Insufficiency of the side sensor to ensure safety.}
  \label{fbf:fig.srct}
\end{figure}

As is illustrated in Fig.~\ref{fbf:fig.srct}, the sensor system of the
vehicle is deficient in capability of on-line detection of head-on
collisions with the domain $D$. In the absence of extra forward-view
sensors, a partial remedy may be systematic full turns to accomplish
environment mapping within the entire vicinity of the vehicle given by
the sensor range, which however consumes extra resources and may be
unacceptable. If no special measures are taken to explore the forward
direction, collisions can be excluded only due to special geometric
properties of the obstacle $D$. They should guarantee a certain amount
of free forward space on the basis of circumstances sensible in the
side direction. These guarantees should cover the entire operational
zone including the transient.

Local guarantees of such a kind are given by Eq.\eqref{fbf:mecess}.  To
highlight this, the outer normal ray
to $\partial D$ (rooted at $\blr \in \partial D$) is denote by $\mathscr{N}(\blr)$. Also, the
point of $\mathscr{N}(\blr)$ at the distance $L$ from $\blr$ is denoted by $\blr(L)$.  Then
the second inequality from Eq.\eqref{fbf:mecess} implies that some piece
$\partial_s D$ of $\partial D$ surrounding $\blr$ does not intersect
the open disk of the radius $d_0+R$ centered at $\blr(d_0+R)$
\cite{Arnold93}, as is illustrated in Fig.~\ref{fbf:fig.locguar}. In
turn, this implies that first, the smaller disk of the radius $R$
centered at $\blr(d_0+R)$ is separated from $\partial_s D$ by a
distance of not less than $d_0 > d_-$ and second, the smaller open
disk $\mathfrak{D}(\blr,d_0)$ of the radius $d_0$ centered at
$\blr(d_0)$ does not intersect $\partial_s D$.  Now these
guarantees are extended on the entire boundary and operational zone, which is
assumed to be upper limited $d\leq d_\ast$ by a constant $d_\ast$:

\begin{figure}[ht]
  \centering
  \subfigure[]{\scalebox{0.4}{\includegraphics{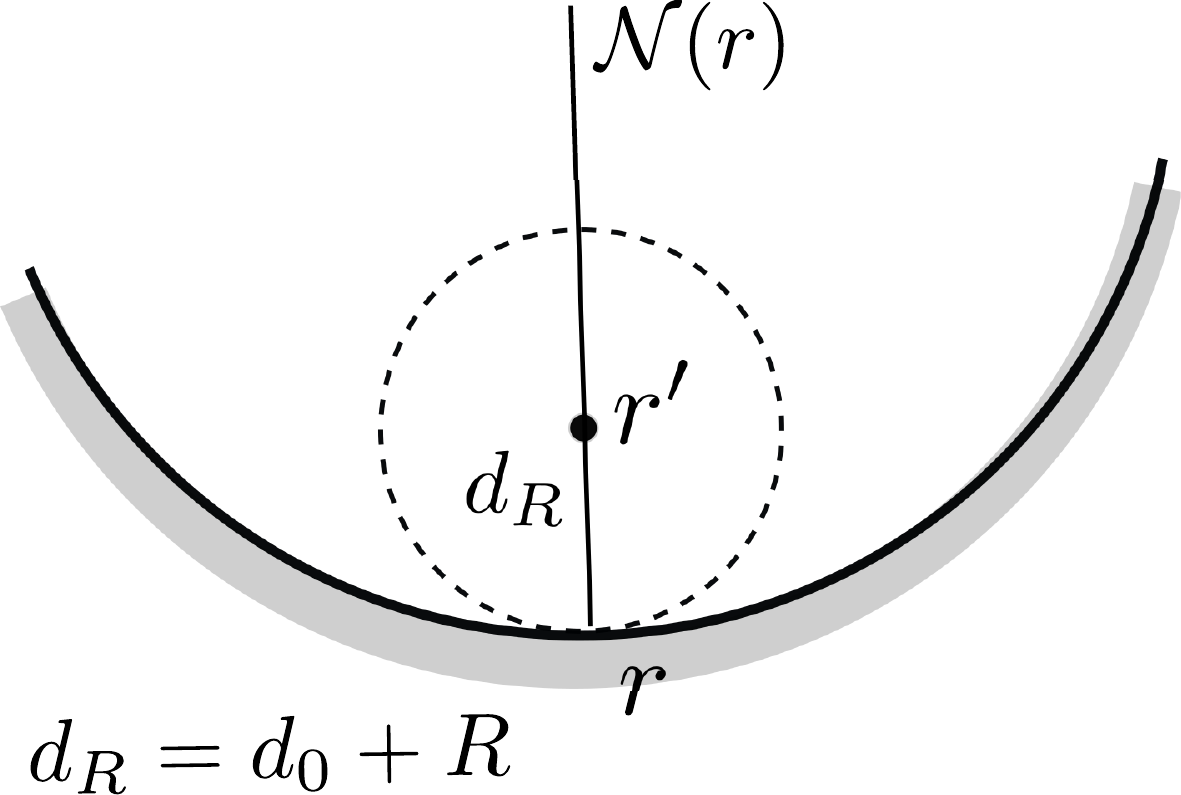}}}
\hspace{20pt}
  \subfigure[]{\scalebox{0.4}{\includegraphics{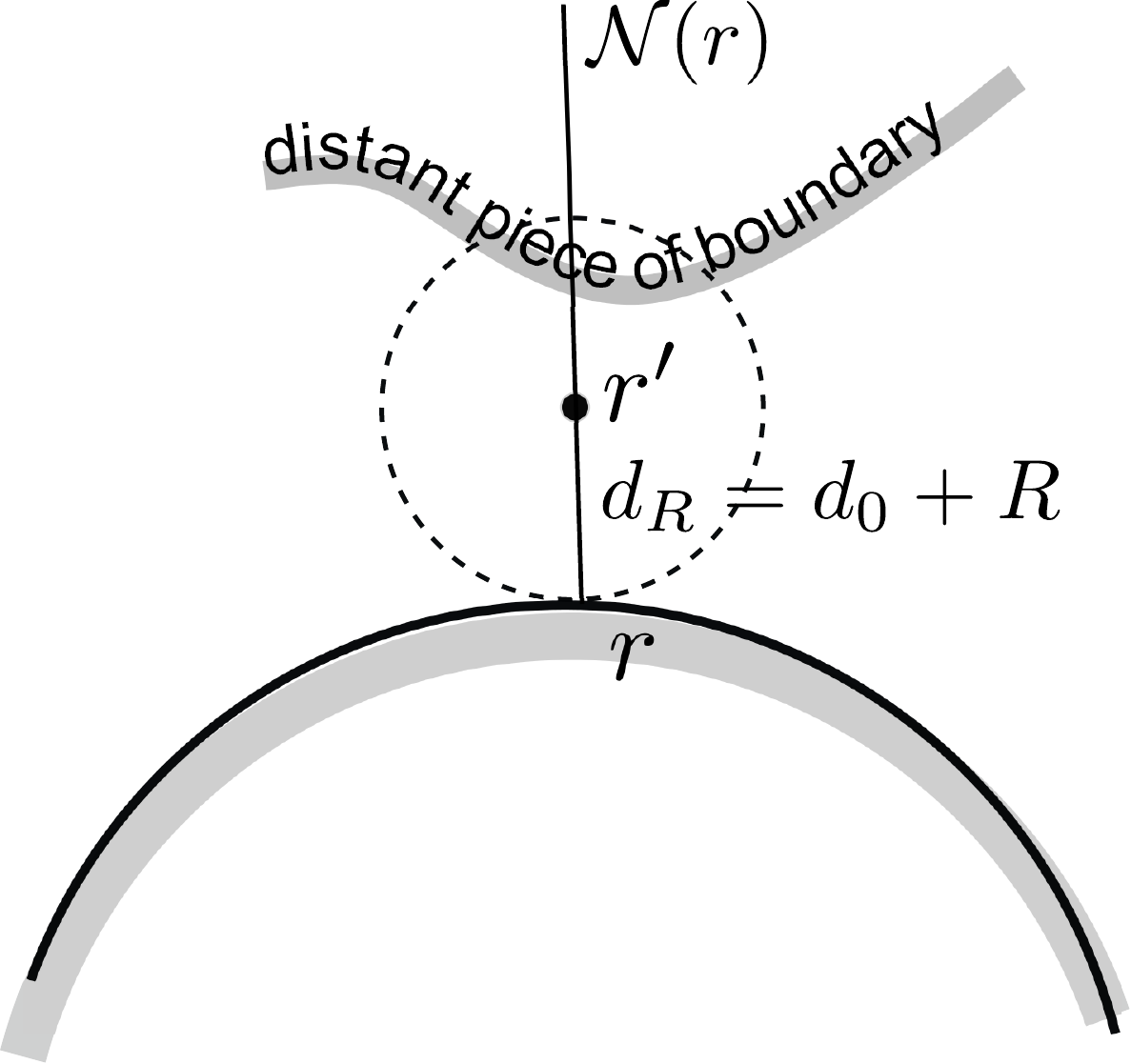}}}
  \caption{A disk free of collision with the local part of the
    boundary.}
  \label{fbf:fig.locguar}
\end{figure}

\begin{figure}[ht]
  \centering
  \subfigure[]{\scalebox{0.4}{\includegraphics{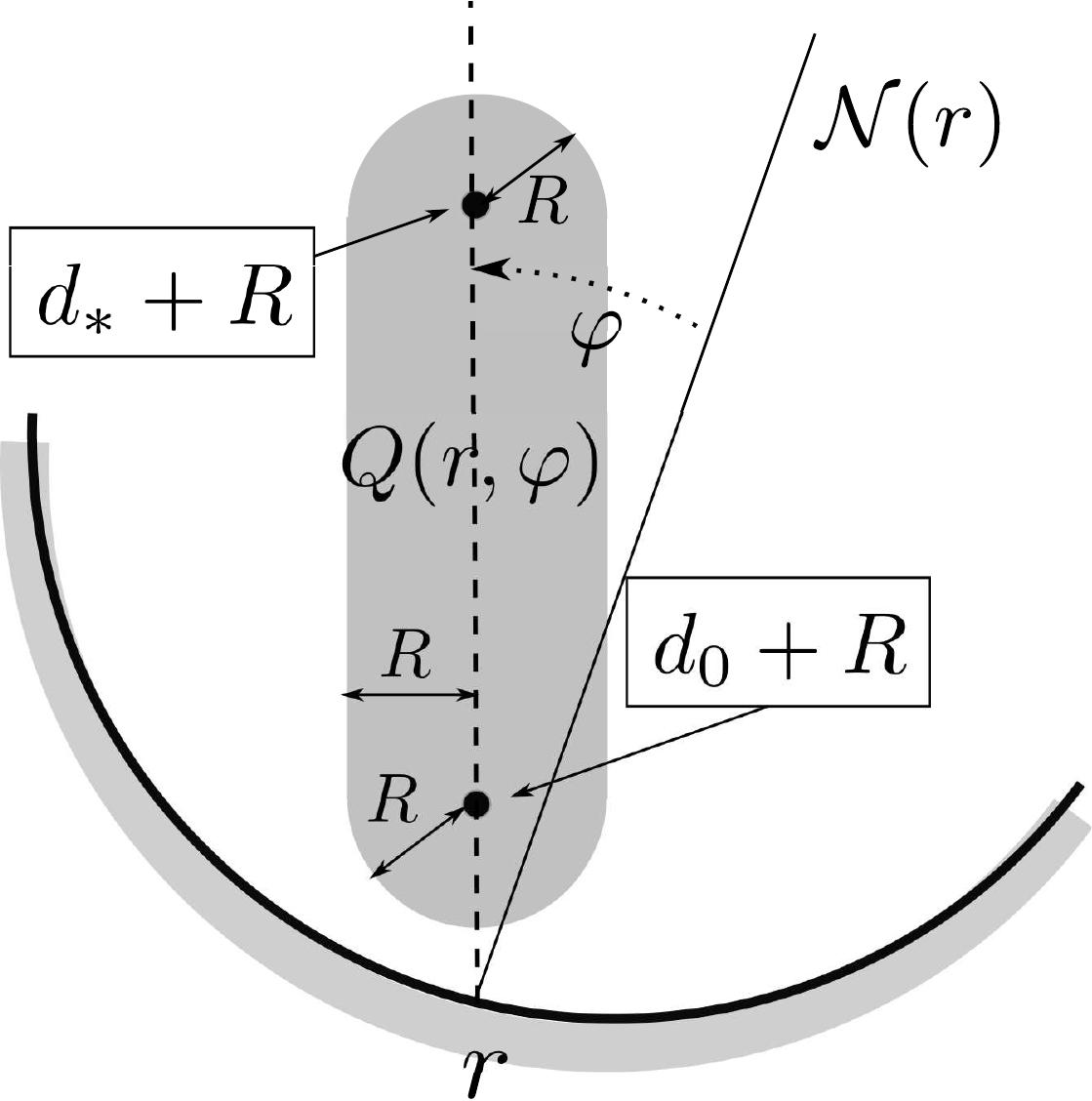}}}
\hspace{20pt}
  \subfigure[]{\scalebox{0.4}{\includegraphics{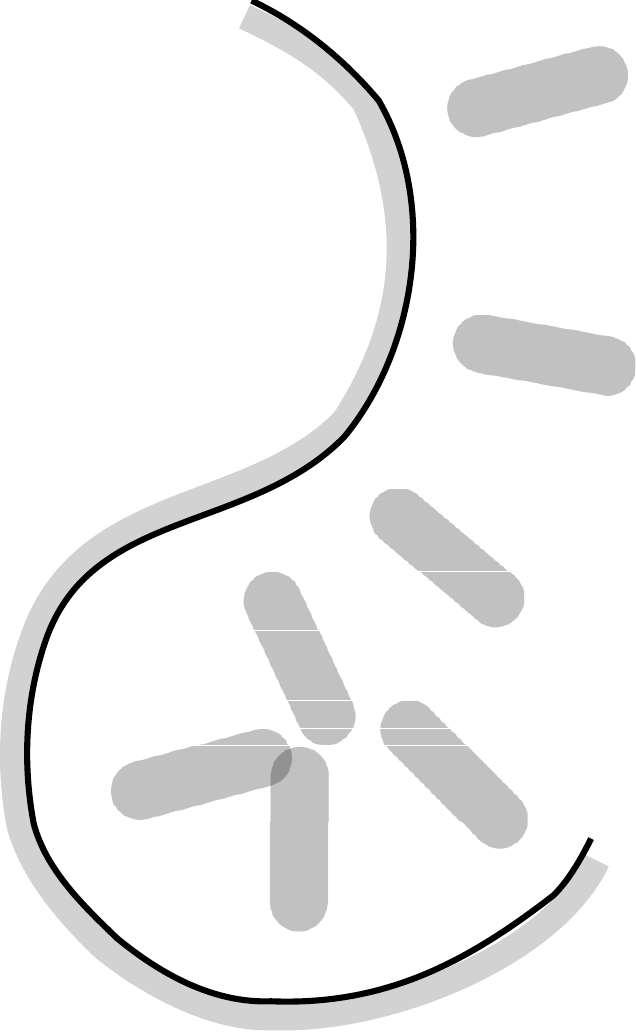}}}
  \caption{(a) The set $Q(\blr,\varphi)$; (b)
    Assumption~\ref{fbf:ass.globsl}.}
  \label{fbf:fig.q}
\end{figure}

\begin{Assumption}
  \label{fbf:ass.globsl}
  The following two claims hold:
  \begin{enumerate}[{\bf (i)}]
  \item For any $\blr \in \partial D$, the open disk
    $\mathfrak{D}(\blr,d_0)$ is disjoint with the entire boundary
    $\partial D$;
  \item There exist $d_\ast > d_0$ and $\eta >0$ such that for any
    $\blr \in \partial D$, the set $Q(\blr,0)$ introduced by Fig.~{\rm
      \ref{fbf:fig.q}(a)} is separated from the boundary $\partial D$
    by a distance of not less than $d_-+ \eta$; see Fig.~{\rm
      \ref{fbf:fig.q}(b)}.
  \end{enumerate}
\end{Assumption}

It follows from Fig.~{\rm \ref{fbf:fig.q}(a) that $d_-+ \eta < d_0$.
  \par
  Finally it is assumed that the distance to $D$ is locally controllable
  when operating at the safety margin $d_-$: it can be maintained
  constant, increased, and decreased by selecting respective
  controls. As was shown in \cite{Matveev2011journ2}, this is
  equivalent to the following enhancement of the first inequality from
  Eq.\eqref{fbf:mecess}:

  \begin{Assumption}
    \label{fbf:ass.control}
    The following inequality holds:
    \begin{equation}
      \label{fbf:mecess1}
      R_\kappa^+(D) + d_-  > R .
    \end{equation}
  \end{Assumption}

  The similar condition $R_\kappa^-(D) > d_- + R$ for the concavity
  parts $\varkappa<0$ of the boundary follows from the second
  inequality in Eq.\eqref{fbf:mecess} since $d_0>d_-$.

  Due to Eq.\eqref{fbf:mecess}, Eq.\eqref{fbf:mecess1}, there exists $\mu \in \left(
    0, \frac{\pi}{2}\right)$ such that:

  \begin{gather}
    \label{fbf:contr.tuning1a}
    R_\varkappa^+(D) \cos \mu + d_- > R ,
    \\
    \label{fbf:contr.tuning2a}
    R_\varkappa^-(D) \cos \mu > R + d_0.
  \end{gather}

  The parameter $\eta_\ast \in (0,\eta)$ may be picked, where $\eta$ is taken from (ii) of
  Assumption~\ref{fbf:ass.globsl}, and by decreasing $\mu$ if
  necessary, ensure the following property, which is possible thanks
  to Assumption~\ref{fbf:ass.globsl}:

  \begin{property}
    \label{fbf:prop.1}
    For any $\blr \in \partial D$ and $\varphi \in [-\mu,\mu]$, the
    following two claims hold:
    \begin{enumerate}[{\bf (i)}]
    \item For any point $\blr^\prime$ such that $\|\blr^\prime -
      \blr\| \leq d_0$ and the angle subtended by $\blr^\prime - \blr$
      and the normal $\mathscr{N}(\blr)$ equals $\varphi$,
      \begin{itemize}
      \item[{\bf i.1)}] the straight line segment with the end-points
        $\blr$ and $\blr^\prime$ has only one point $\blr$ in common
        with $\partial D$;
      \item[{\bf i.2)}] the distance from $\blr^\prime$ to $\partial
        D$ is no less than $d_-$ provided that $\|\blr^\prime - \blr\|
        \geq d_- +\eta_\ast$;
      \end{itemize}
    \item the set $Q(\blr,\varphi)$ is separated from the boundary
      $\partial D$ by a distance of not less than $d_- + \eta_\ast$.
    \end{enumerate}
  \end{property}

  Finally, the variable $[s]_+$ is set to $\max\{s,0\}$, and $\gamma$ is picked so that:

  \begin{subequations}
  \begin{gather}
    \label{fbf:gammab}
    \gamma < \frac{R_\varkappa^+(D)\cos \mu -
      [R-d_-]_+}{R_\varkappa^+(D) [R-d_-]_+\sin \mu} \; \text{if} \;
    R>d_-,
    \\
    \label{fbf:gamma1ab}
    \gamma < \frac{R_\varkappa^-(D) \cos \mu - (R + d_0)}{(R + d_0)
      R_\varkappa^-(D) \sin \mu},
    \\
    \label{fbf:freeab}
    \gamma < \frac{\cos \mu}{(R+ d_0) \sin \mu}.
  \end{gather}
  \end{subequations}

  This choice is possible, since the right-hand sides of all
  inequalities are positive due to Eq.\eqref{fbf:contr.tuning1a} and
  Eq.\eqref{fbf:contr.tuning2a}.
  \par
  If the boundary $\partial D$ is compact, inequalities
  Eq.\eqref{fbf:mecess}--Eq.\eqref{fbf:gamma1ab} can be checked in the point-wise
  fashion. In doing so, any inequality involving $R_\varkappa^+(D)$
  should be checked at any point $\blr \in \partial D$ of convexity
  $\varkappa (\blr) >0$ with substituting $R_\varkappa(\blr)$ in place
  of $R_\varkappa^+(D)$. Similarly, any inequality involving
  $R_\varkappa^-(D)$ should be checked at any point $\blr \in \partial
  D$ of concavity $\varkappa (\blr) < 0$ with substituting
  $R_\varkappa(\blr)$ in place of $R_\varkappa^-(D)$.

   \section{Summary of Main Results}
  \label{fbf:sec.mr}

 When examining
  convergence of the control law, it is assumed that the initial state is
  in the set $\mathscr{V}$ of all states $\blr \not\in D,\theta$ for
  which the domain is visible. Let $d(\blr,\theta)$ denote the
  corresponding measurement $d$. \par
  The set $\mathfrak{C}$ of initial states $\blr,\theta$ from which
  convergence to the required equidistant curve can be theoretically
  guaranteed is composed of three parts $\mathfrak{C}_0,
  \mathfrak{C}_-$, and $\mathfrak{C}_+$. They contain initial states
  with $S=0, S<0$, and $S>0$, respectively, where $S$ is defined in
  Eq.\eqref{fbf:c.a}:

  \begin{equation}
    \label{fbf:c.0}
    \mathfrak{C}_0 := \Big\{ (\blr, \theta) \in \mathscr{V} : S = 0 \; \text{and}\; d_-+\eta_\ast \leq d \leq d_\ast \Big\}.
  \end{equation}

  To introduce $\mathfrak{C}_-$, attention is first paid to the
  following:

  \begin{Lemma}
    \label{fbf:lem.termin1}
    Under the control law Eq.\eqref{fbf:c.a}, motion with $S<0$ necessarily
    terminates with arrival at $S=0$, provided that the vehicle does
    not collide with the domain $D$.
  \end{Lemma}

  Let $C^-_{\blr,\theta}$ be the circle of radius Eq.\eqref{fbf:Rmin} traced
  clockwise from the initial state $\blr,\theta$ and $(\blr_\ast,
  \theta_\ast)$ be the first position on this circle that either
  belongs to $D$ or is such that $S=0$. By
  Lemma~\ref{fbf:lem.termin1}, this position does exist. Let also
  $\widehat{C}^-_{\blr,\theta}$ denote the arc of $C^-_{\blr,\theta}$
  between these two positions, and $\textbf{dist}(A,B) := \inf_{\blr
    \in A, \blr^\prime \in B} \|\blr-\blr^\prime\|$ denote the
  distance between the sets $A$ and $B$. The second part of the set
  $\mathfrak{C}$ is given by:

  \begin{multline}
    \label{fbf:c.1}
    \mathfrak{C}_- := \Big\{ (\blr, \theta) \in \mathscr{V} : S < 0,
    \quad \textbf{dist}[\widehat{C}^-_{\blr,\theta},D] \geq d_-,
    \\
    \text{and} \quad d (\blr_\ast,\theta_\ast) \geq d_-+\eta_\ast
    \Big\}.
  \end{multline}

  Introduction of the last part $\mathfrak{C}_+$ is prefaced by the
  following:

  \begin{Lemma}
    \label{fbf:lem.termin2}
    Under the control law Eq.\eqref{fbf:c.a} and for initial states with
    $S>0$, there may be the following three scenarios:
    \begin{enumerate}[{\rm i)}]
    \item With maintaining $S > 0$, the safety margin is violated;
    \item With respecting the safety margin and maintaining $S>0$, the
      vehicle arrives at a position where the view of $D$ becomes
      obstructed by another part of $D$; at this moment, $S$ abruptly
      jumps down at a negative value;
    \item With respecting the safety margin and maintaining both the
      view of $D$ unobstructed and $S>0$, the vehicle arrives at
      $S=0$.
    \end{enumerate}
  \end{Lemma}

  While $S>0$, the vehicle moves counter-clockwise over the circle of
  the radius $d$ centered at the reflection point (which does not
  move) whenever $d>R$. Otherwise it moves with the maximal turning
  rate $\ov{u}$ over a circle of radius Eq.\eqref{fbf:Rmin}.
  \par
  Finally the set $\mathfrak{C}_+$ is introduced, which contains all initial states
  $(\blr, \theta) \in \mathscr{V}$ for which the following claims
  hold:

  \begin{itemize}
  \item Scenario i) from Lemma~\ref{fbf:lem.termin2} does not hold and
    $S(0)>0$;
  \item In the case ii) from Lemma~\ref{fbf:lem.termin2}, the vehicle
    arrives at a state from the set $\mathfrak{C}_-$ when $S$ becomes
    negatives
  \item In the case iii) from Lemma~\ref{fbf:lem.termin2}, the vehicle
    arrives at a state from the set $\mathfrak{C}_0$ when $S$ becomes
    zero.
  \end{itemize}

  Now the main result of the chapter may be stated:
  
  \begin{Theorem}
    \label{fbf:th.m}
    Let Assumptions~{\rm \ref{fbf:ass.prelim}--\ref{fbf:ass.control}}
    be true and in Eq.\eqref{fbf:c.a}, the parameters be chosen so that
    Eq.\eqref{fbf:contr.tuning1a}--Eq.\eqref{fbf:freeab} and Property~{\rm
      \ref{fbf:prop.1}} hold. If the initial state lies in the set
    $\mathfrak{C}:= \mathfrak{C}_0 \cup \mathfrak{C}_- \cup
    \mathfrak{C}_+$, the vehicle driven by the navigation and guidance
    law Eq.\eqref{fbf:c.a} does not lose track of the domain $D$, respects
    the safety margin, and asymptotically follows the boundary of $D$
    at the required distance: $d(t) \to d_0, \varphi(t) \to 0$ as $t
    \to \infty$.
  \end{Theorem}

 \section{Simulations}
\label{fbf:sec.sim}

Simulations were performed using the perfect kinematic model of the
vehicle Eq.\eqref{fbf:1}. To estimate the angle $\varphi$, the tangent at the
reflection point was approximated by the secant between this point and
another point slightly in front; the angular separation between these
points was $9$ deg. The control law was updated with the sampling
period of $0.1 s$. Other parameters used for simulation are shown in
Table~\ref{fbf:fig:paramsim}.

\begin{table}[ht]
  \centering
  \begin{tabular}{| l | c |}
    \hline
    $\overline{u}$ & $45.8 deg/s$ \\
    \hline
    $v$ & $0.3 m/s$ \\
    \hline
    $\mu$ & $57.3 deg$  \\
    \hline
  \end{tabular}
\hspace{10pt}  
\begin{tabular}{| l | c |}
    \hline
    $\gamma$ & $171.9 deg/m$  \\
    \hline
    $d_{0}$ & $1.0 m$ \\
    \hline
  \end{tabular}
  \caption{Simulation parameters for fixed-sensor boundary following controller.}
  \label{fbf:fig:paramsim}
\end{table}

In the first simulation test, the domain $D$ fits the maneuverability
of the robot: the minimal turning radius of the vehicle exceeds the
radius required for perfectly tracking the boundary of $D$ with the
requested margin $d_0$. Fig.~\ref{fbf:fig:simsimp} shows that after a
short transient, the proposed control law provides a visibly perfect
motion over the desired equidistant curve and successfully copes with
both convexities and concavities of the obstacle, as well as with
transitions from convexities to concavities and vice versa.

\begin{figure}[ht]
  \centering
  \includegraphics[width=0.6\columnwidth]{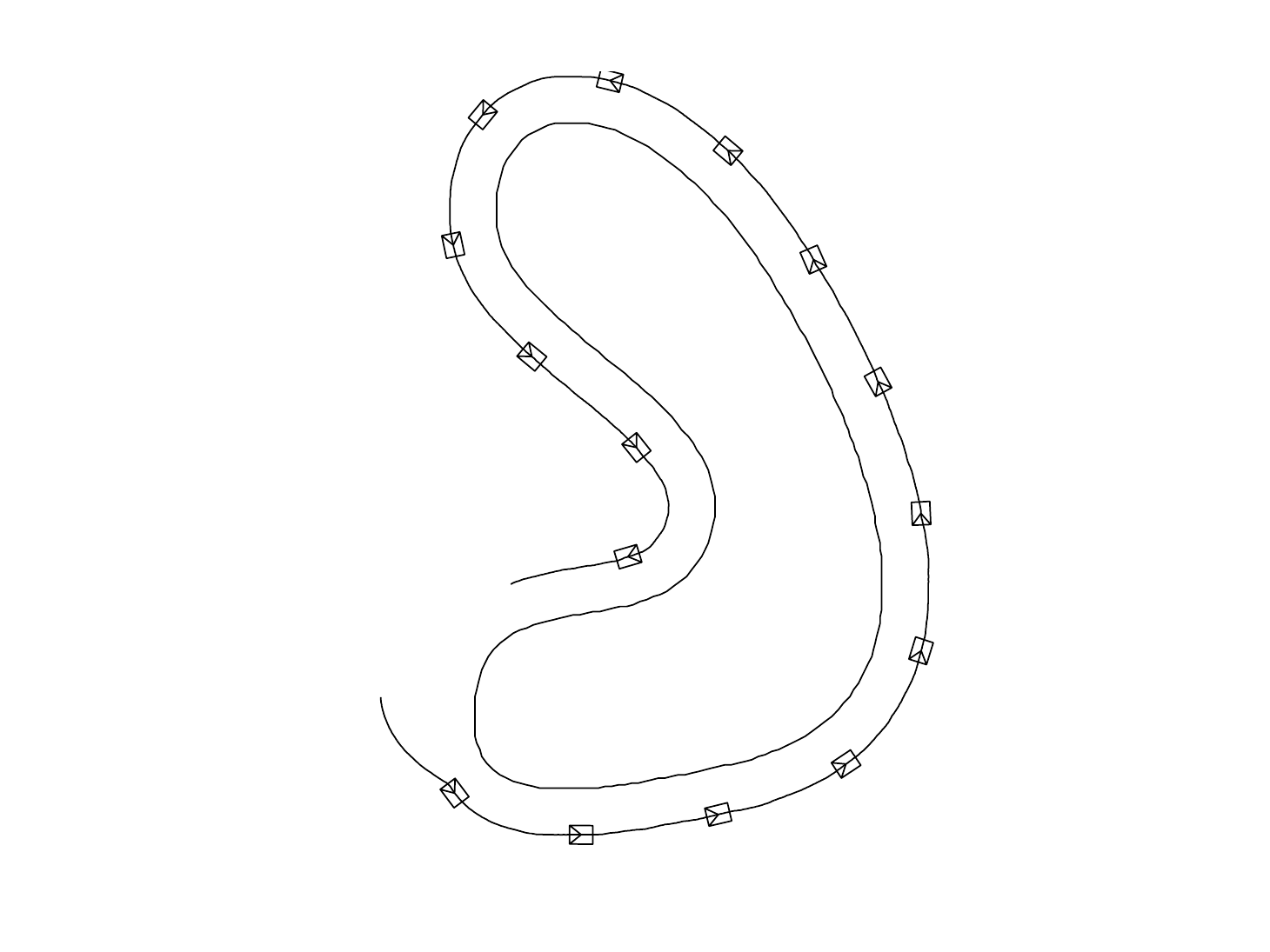}
  \caption{Simulations with a simple boundary.}
  \label{fbf:fig:simsimp}
\end{figure}

Figs.~\ref{fbf:fig:dissimsimp} and \ref{fbf:fig:angsimsimp} provide a
closer look at the boundary following errors. After the transient is
completed ($t \gtrapprox 18 sec$), the error in the true distance to
the obstacle Eq.\eqref{fbf:true.dist} does not exceed $1 cm$, whereas the
error in the estimated angular discrepancy $\varphi$ between the
tangent at the reflection point and the vehicle centerline does not
exceed $12.6^\circ$. However, this good exactness proceeds from taking
into account only the non-idealities that are due to control sampling
and numerical evaluation of the relative tangent angle.

\begin{figure}[ht]
  \centering
  \includegraphics[width=10cm]{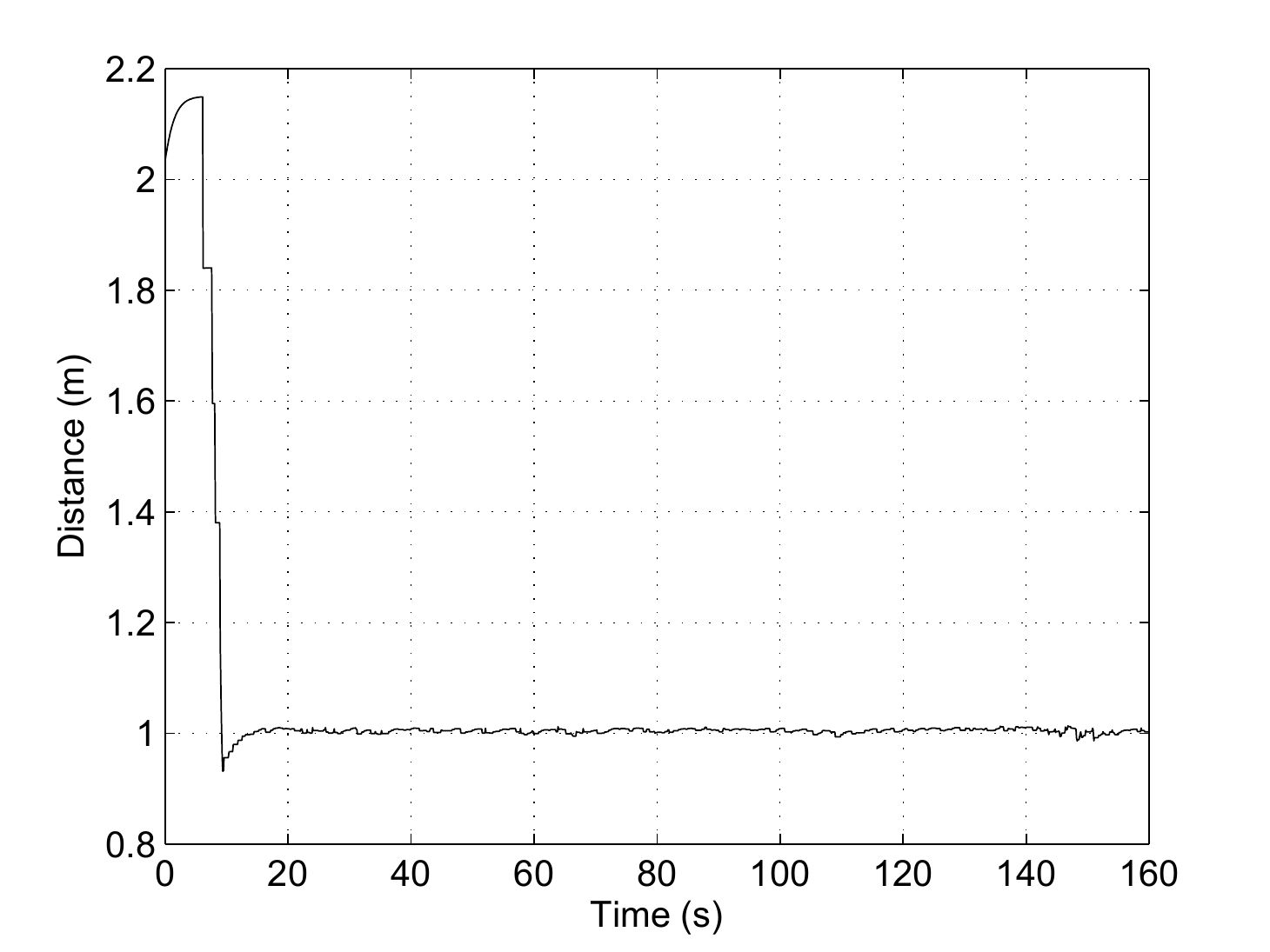}
  \caption{Distance to the obstacle during
    simulation.}
  \label{fbf:fig:dissimsimp}
\end{figure}

\begin{figure}[ht]
  \centering
  \includegraphics[width=10cm]{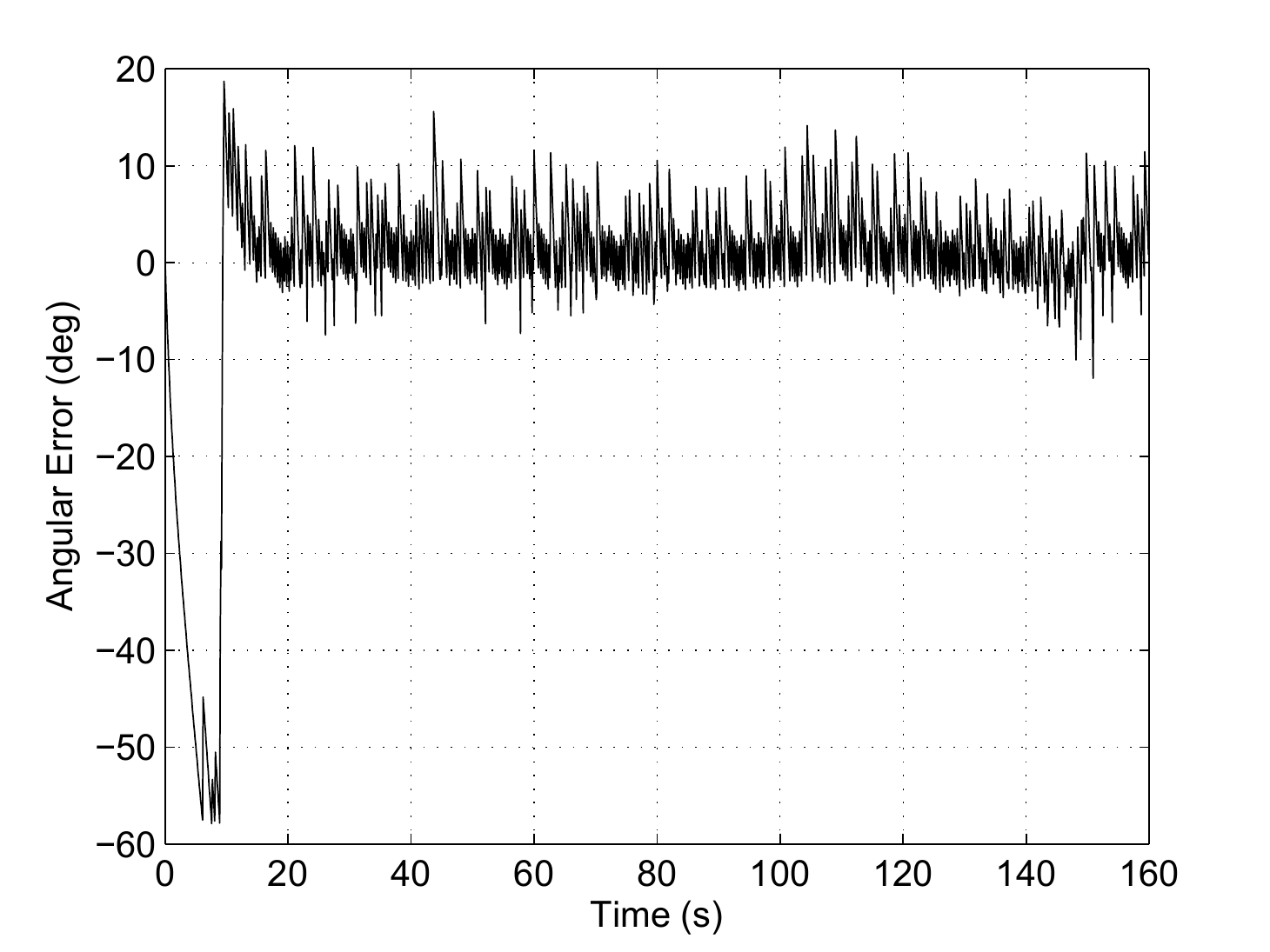}
  \caption{Estimated relative tangent angle
    $\varphi$ during simulation.}
  \label{fbf:fig:angsimsimp}
\end{figure}

The second group of simulation tests provides deeper insights into the
effects of real-live non-idealities on the performance of the
closed-loop system. These tests were carried out in the previous scene
with additionally taking into account sensor and actuator noises and
un-modeled dynamics. To this end, a bounded random and uniformly
distributed offset was added to every relevant quantity at each
control update. Specifically, the noises added to $d$ and $\varphi$
were $0.3m$ and $17.2^\circ$, respectively; the noise added to the
control signal was $11.5deg/s$, and the control signal was not allowed
to change faster than $4 rad s^{-2}$. These are relatively large
noises that would be unlikely met for typical modern sensors and
actuators. The test was repeated several times, each with its own
realization of the random noises. Ten typical results are depicted in
Fig.~\ref{fbf:fig:simnoisesimp}. They show that the control objective
is still achieved with the distance error $\leq 0.3 m$. Since this is
the accuracy of the distance sensor, the result seems to be more than
satisfactory.

\begin{figure}[ht]
  \centering
  \includegraphics[width=10cm]{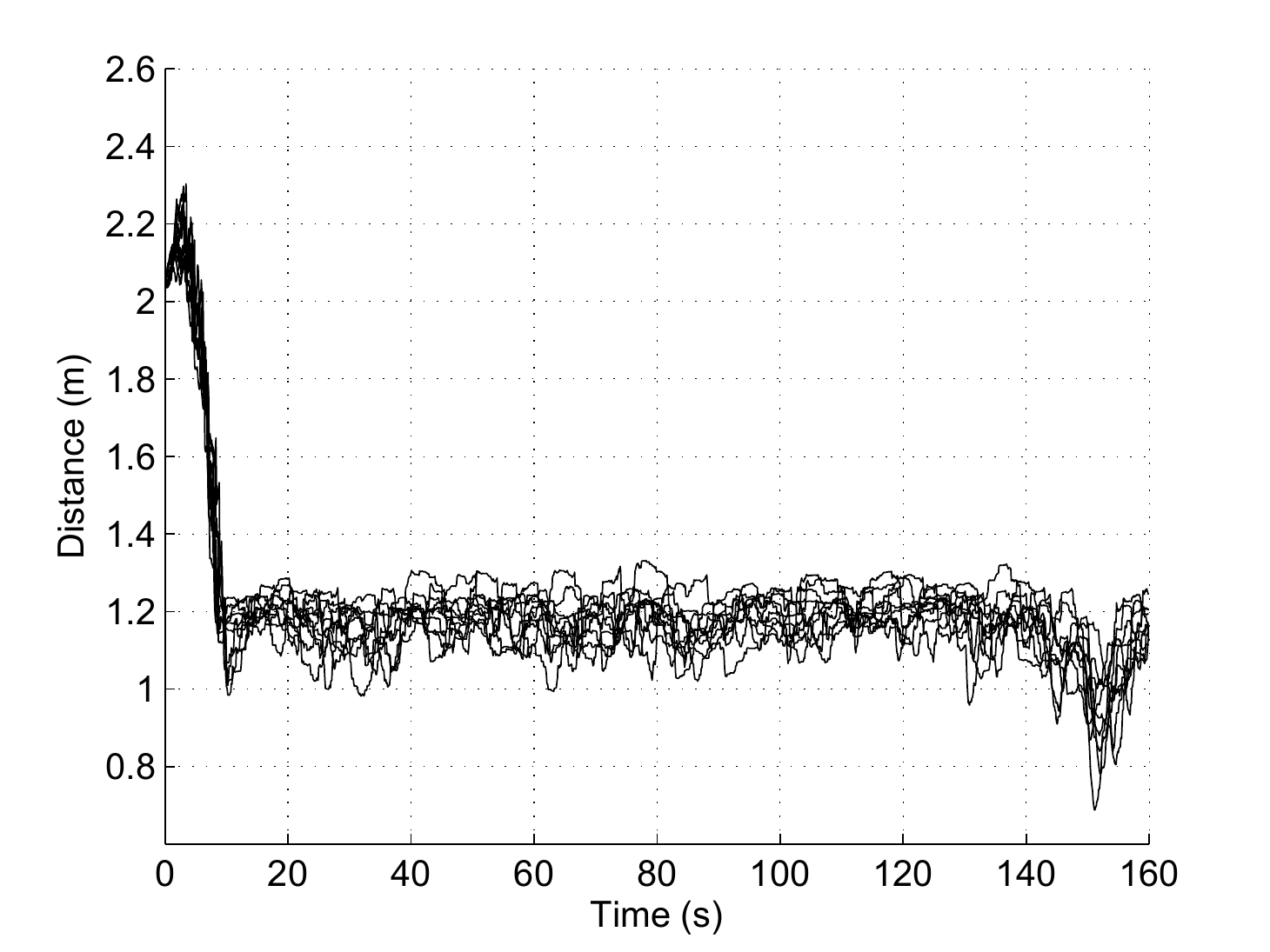}
  \caption{Distance to obstacle during
    simulations with actuator dynamics and sensor noises.}
  \label{fbf:fig:simnoisesimp}
\end{figure}

The purpose of the next test is to examine the performance of the
algorithm in the case when the obstacle boundary $\partial D$ contains
points where the vehicle is absolutely incapable of maintaining the
required distance $d_0$ to $D$ because of the limited turning radius
(see Fig.~\ref{fbf:fig:imturn}) -- to move over the equidistant curve,
the vehicle should turn sharper than feasible.  Though the main
theoretical results of the work are not concerned with this case, it
may be hypothesized that if these points constitute only a small piece
of the boundary, the overall behavior of the closed-loop system
remains satisfactory\footnote{since this behavior is expected to be
  close to that in the absence of such points}. To verify this
hypothesis and reveal details, the obstacle depicted in
Fig.~\ref{fbf:fig:sim} is considered, for which concavities do contain the
afore-mentioned points.

\begin{figure}[ht]
  \centering
  \includegraphics[width=0.6\columnwidth]{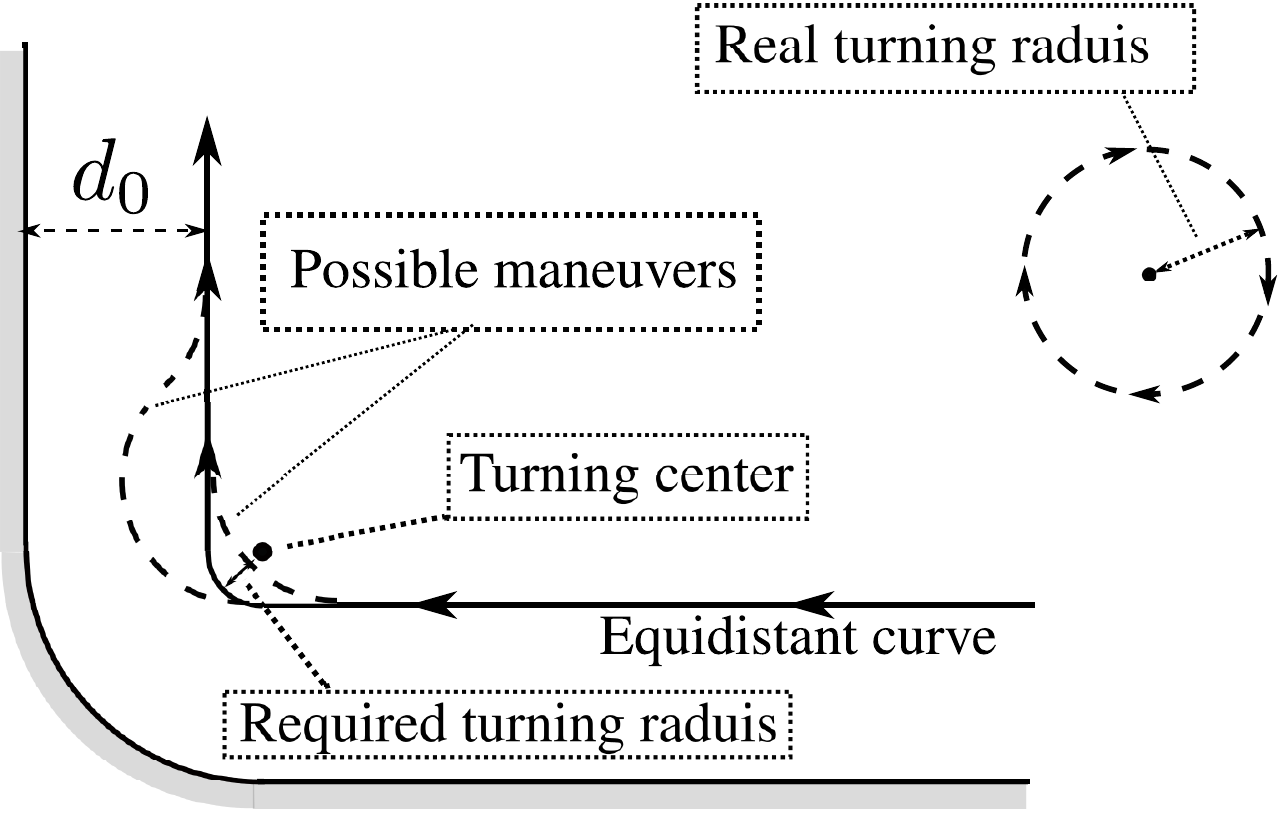}
  \caption{Unavoidable disturbance of the distance to $D$ because of
    the limited turning radius.}
  \label{fbf:fig:imturn}
\end{figure}

The related results are shown in Figs.~\ref{fbf:fig:sim},
\ref{fbf:fig:dissim} and \ref{fbf:fig:angsim}. All abrupt local both
falls of the distance $d$ in Fig.~\ref{fbf:fig:dissim} and deviations
of the angle $\varphi$ from $0$ in Fig.~\ref{fbf:fig:angsim} hold at
$t \approx 90, 130, 190, 300, 340, 400, 490 sec$ when the vehicle
passes points where it is absolutely incapable to maintain the
distance to $D$ at the level $d_0$. However, the algorithm
demonstrates a good capability of quickly recovering after these
unavoidable distance errors, and except for the related short periods
of time and the initial transient, keeps the distance error within a
very small bound of $1 cm$ and always maintains the margin of safety
at the level $\geq 0.3 m$. This good overall performance is
illustrated by Fig.~\ref{fbf:fig:sim}, where the vehicle's path looks
like nearly perfect except for few very local violations of the
required distance, with the most of them being hardly visible. The
worst distance error is observed when passing the left upper concavity
in Fig.~\ref{fbf:fig:sim}. This concavity is similar to that from
Fig.~\ref{fbf:fig:imturn} -- it is of the right-angle type. Such
concavities are challenging for vehicles equipped with only side-view
sensors.  For example, when perfectly following the horizontal part of
the equidistant curve from Fig.~\ref{fbf:fig:imturn}, the vehicle with
side-view sensor is incapable to detect the need for turn until the
reflection point leaves the flat part of the boundary. Even if the
vehicle performs the sharpest turn after this, its forward advancement
towards the vertical part of the boundary may be nearly equal to the
minimal turning radius Eq.\eqref{fbf:Rmin} (which holds if in
Fig.~\ref{fbf:fig:imturn}, the `requested turning radius' is close to
$0$). For the vehicle examined in the simulation tests, the minimal
turning radius amounts to $0.375 m$. So theoretically the distance to
the obstacle may be reduced to $1m-0.375 m = 0.625 m$ when passing the
above concavity. Practically it reduces to $\approx 0.3 m$, where the
difference is basically caused by errors in numerical evaluation of
the tangent angle $\varphi$. Except for this troublesome concavity,
the control law ensures the safety margin of $\geq 0.6m$ according to
Fig.~\ref{fbf:fig:dissim}.

\begin{figure}[ht]
  \centering
  \includegraphics[width=0.4\columnwidth]{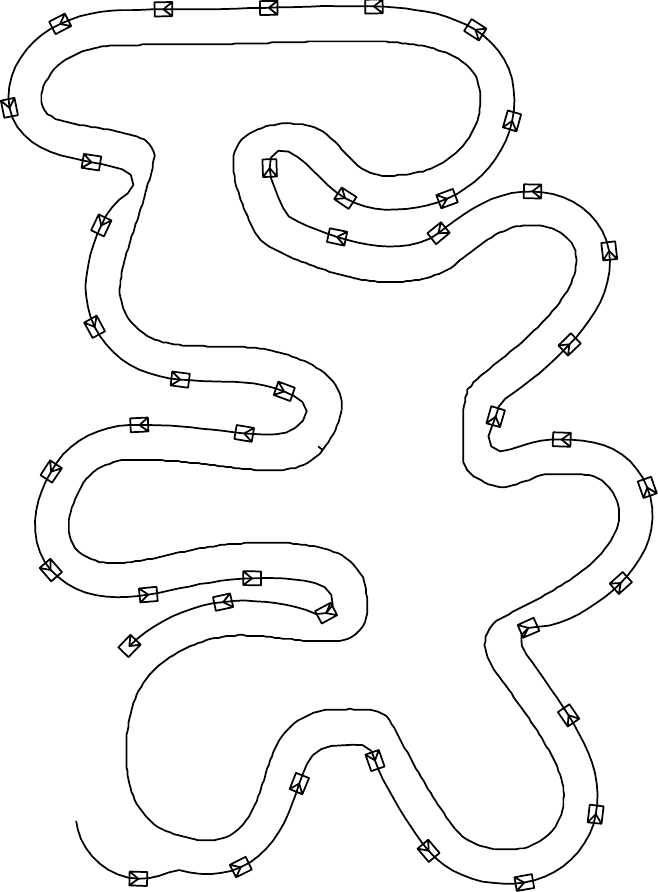}
  \caption{Simulation with a tight boundary.}
  \label{fbf:fig:sim}
\end{figure}

\begin{figure}[ht]
  \centering
  \includegraphics[width=10cm]{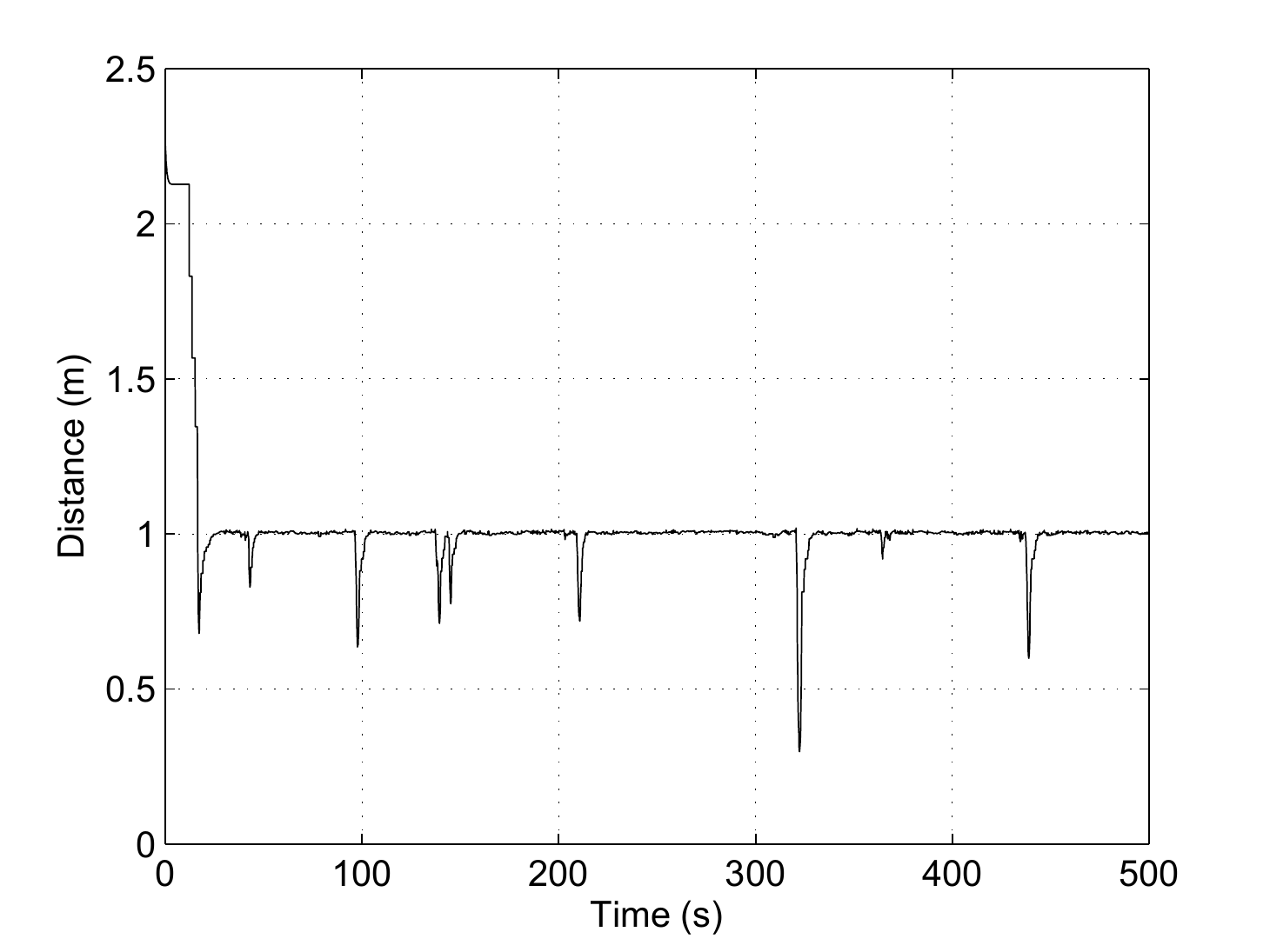}
  \caption{Distance to the obstacle
    during simulation.}
  \label{fbf:fig:dissim}
\end{figure}
\begin{figure}[ht]
  \centering
  \includegraphics[width=10cm]{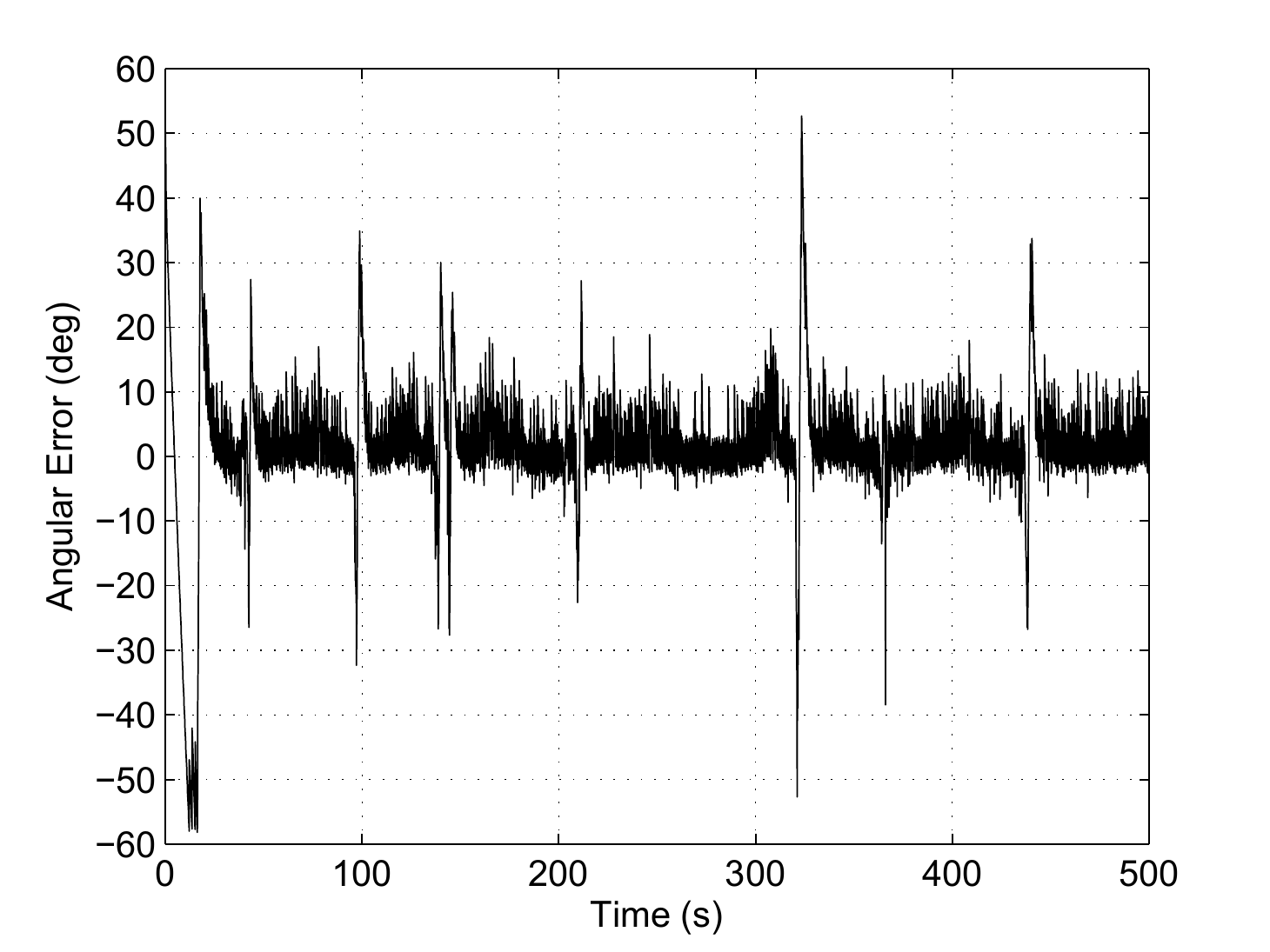}
  \caption{Estimated relative tangent angle $\varphi$
    during simulation.}
  \label{fbf:fig:angsim}
\end{figure}

In this test like in the first one, only non-idealities related to
control sampling and numerical evaluation of the relative tangent
angle were taken into account. Simulations were also carried out to
test the additional effect of un-modeled dynamics, sensor noise, and
actuator noise on the performance of the closed-loop system in the
environment from Fig.~\ref{fbf:fig:sim}.  All these phenomena were
modeled like in the second group of simulations.  Ten typical results
corresponding to various realizations of the random noises are
depicted in Fig.~\ref{fbf:fig:simnoise}. Similarly to the second group
of simulations, the distance accuracy degradation is approximately
equal to the error of the distance sensor, which seems to be the fair
price for using imperfect sensors. The vehicle still successfully
follows the boundary without collision with the obstacle.

\begin{figure}[ht]
  \centering
  \includegraphics[width=10cm]{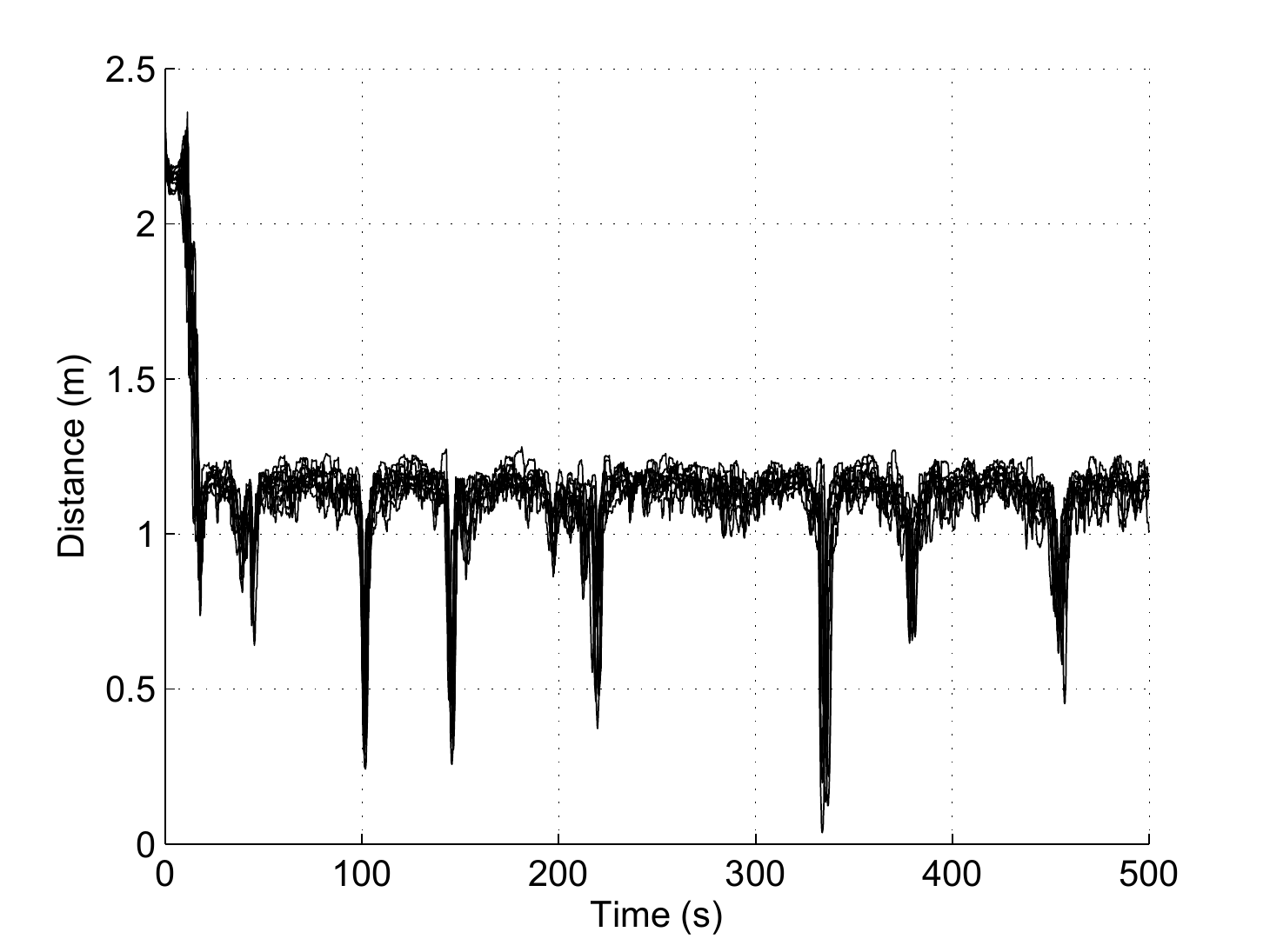}
  \caption{Distance to the obstacle during
    simulations with actuator dynamics and sensor noises.}
  \label{fbf:fig:simnoise}
\end{figure}

The next series of simulation tests were aimed at illustration of the
algorithm performance in environments with many obstacles. The
objective of the algorithm is to follow the boundary of the selected
obstacle $D_0$ despite presence of the others. According to the above
discussion, it does follow the boundary provided that firstly, the
view of $D_0$ is not obstructed and secondly, the assumptions of
Theorem~\ref{fbf:th.m} are satisfied, in particular, the boundary of
$D_0$ is everywhere smooth. The first requirement can be typically met
by picking the desired distance of boundary following $d_0$ small
enough as compared with spacing between obstacles. So the focus in the
tests was on following boundaries with fractures (see
Fig.~\ref{fbf:fig:fracture}).

\begin{figure}[ht]
  \centering
  \subfigure[]{\scalebox{0.4}{\includegraphics{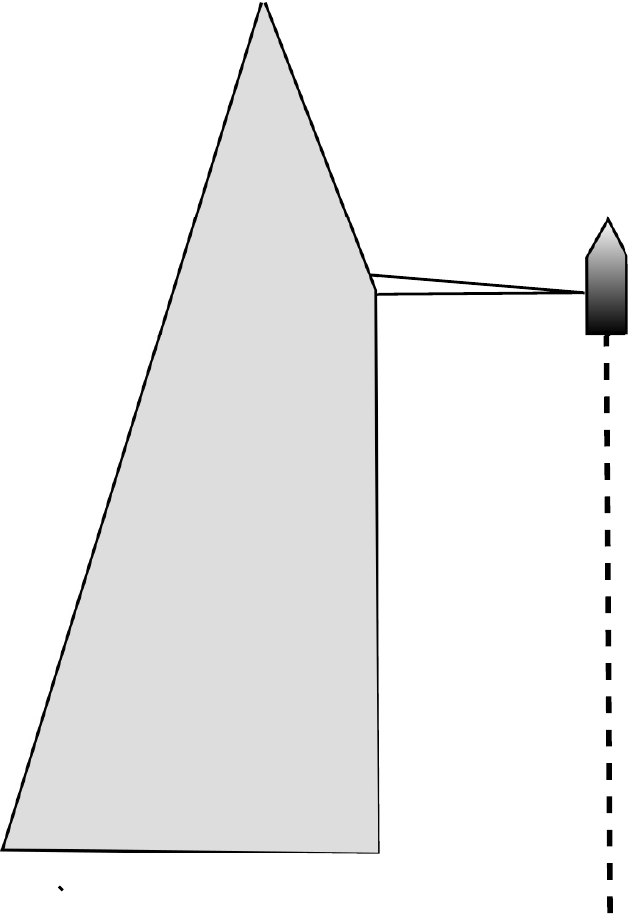}}}
  \subfigure[]{\scalebox{0.4}{\includegraphics{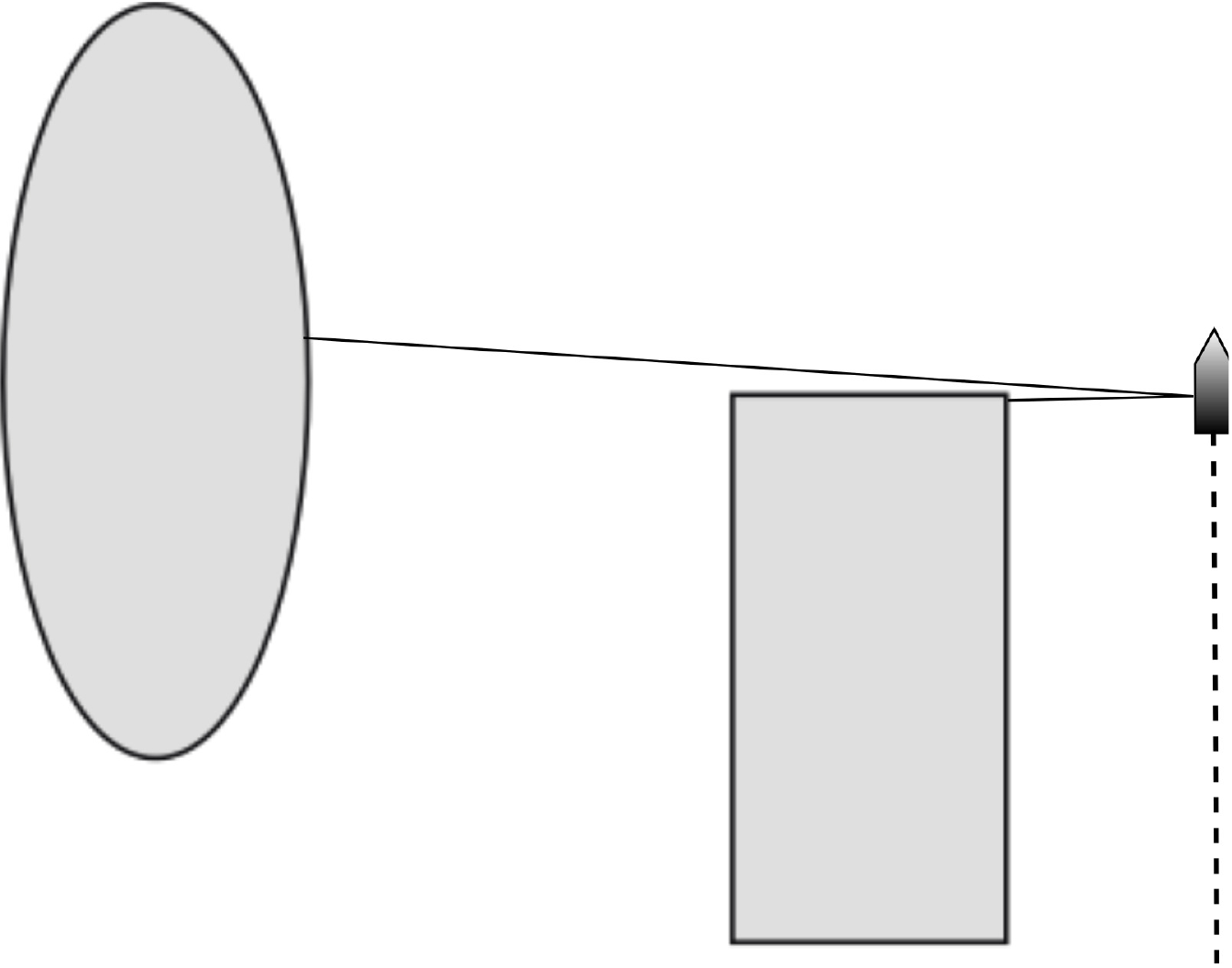}}}  
  \subfigure[]{\scalebox{0.4}{\includegraphics{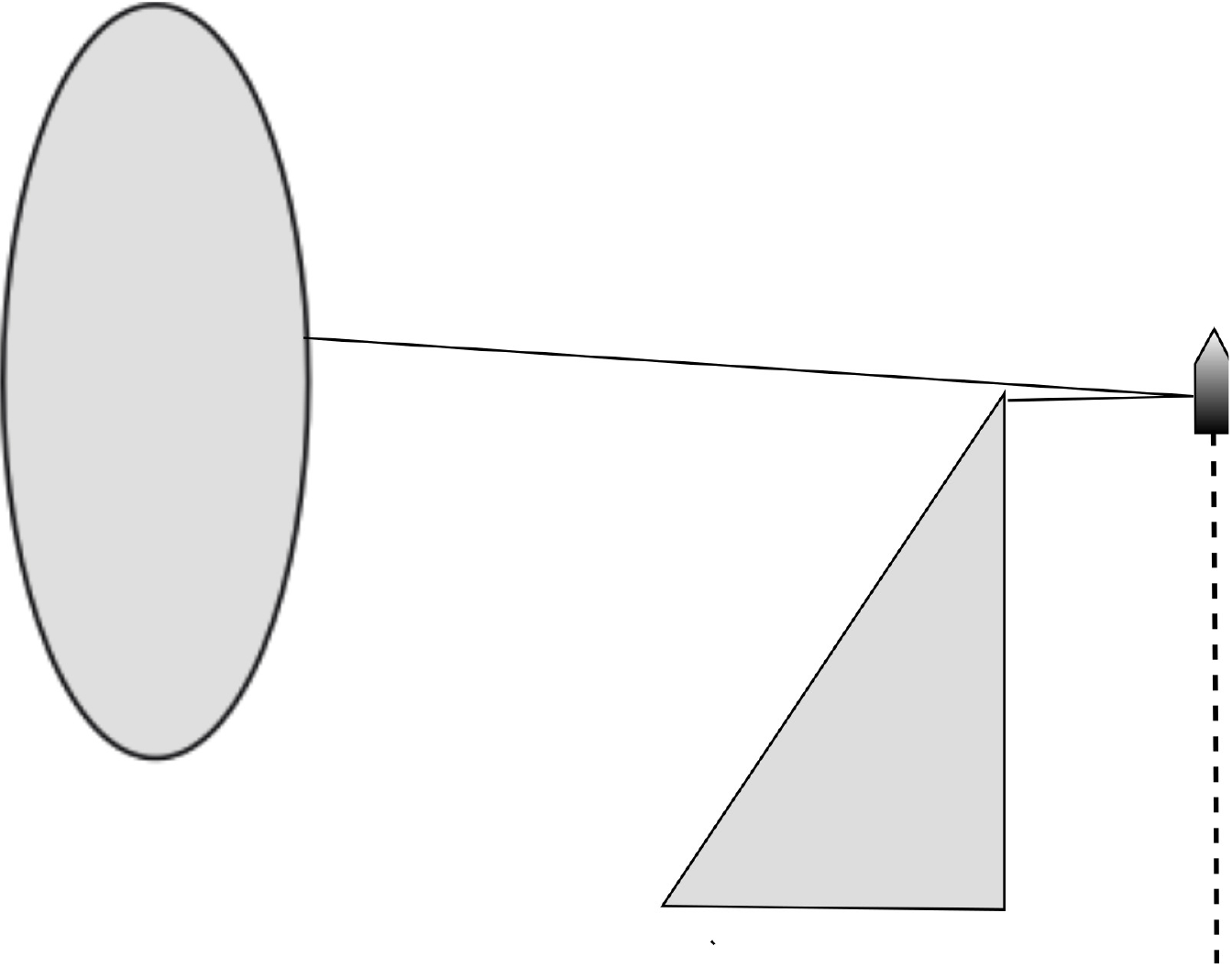}}}
  \caption{Following boundaries with fracture points.}
  \label{fbf:fig:fracture}
\end{figure}

At the fracture point, the tangent $T$ to the boundary and the angle
$\varphi$ between $T$ and the vehicle centerline, which is used in the
control law Eq.\eqref{fbf:c.a}, are strictly speaking undefined. So the
employed method to access $\varphi$ may be puzzled at fractures, with
the outcome being dependent on the method. In these experiments, the
tangent $T$ was approximated by the secant between the points
$\bldr_\ast$ and $\bldr_+$ of incidence of the perpendicular ray
$R_\ast$ and a ray $R_+$ slightly in front, respectively. When
arriving at a fracture point, this method implies an abrupt increase
of $\varphi$. Let for simplicity the vehicle perfectly follow the
boundary $d \equiv d_0, \varphi \equiv 0$ prior to this event. Then
the angle becomes positive $\varphi >0$. By Eq.\eqref{fbf:c.a}, this causes
vehicle rotation about $\bldr_\ast$ at the distance $d_0$ from
$\bldr_\ast$; see Fig.~\ref{fbf:fig:rotation}.

\begin{figure}[ht]
  \centering
  \subfigure[]{\scalebox{0.8}{\includegraphics{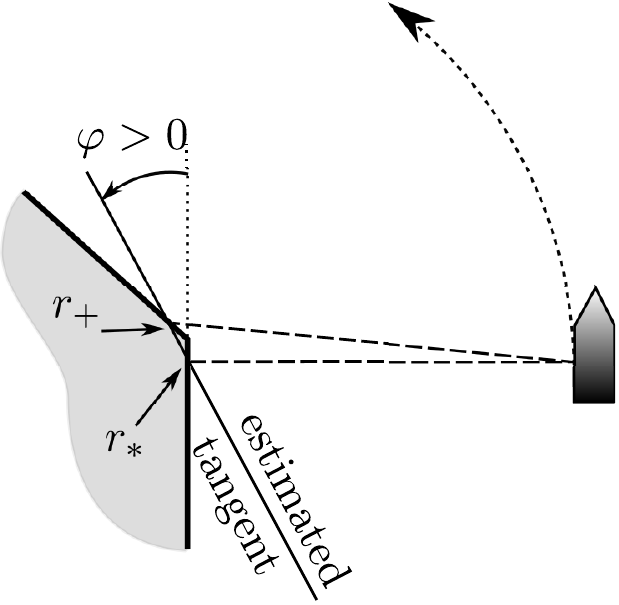}}}
  \subfigure[]{\scalebox{0.4}{\includegraphics{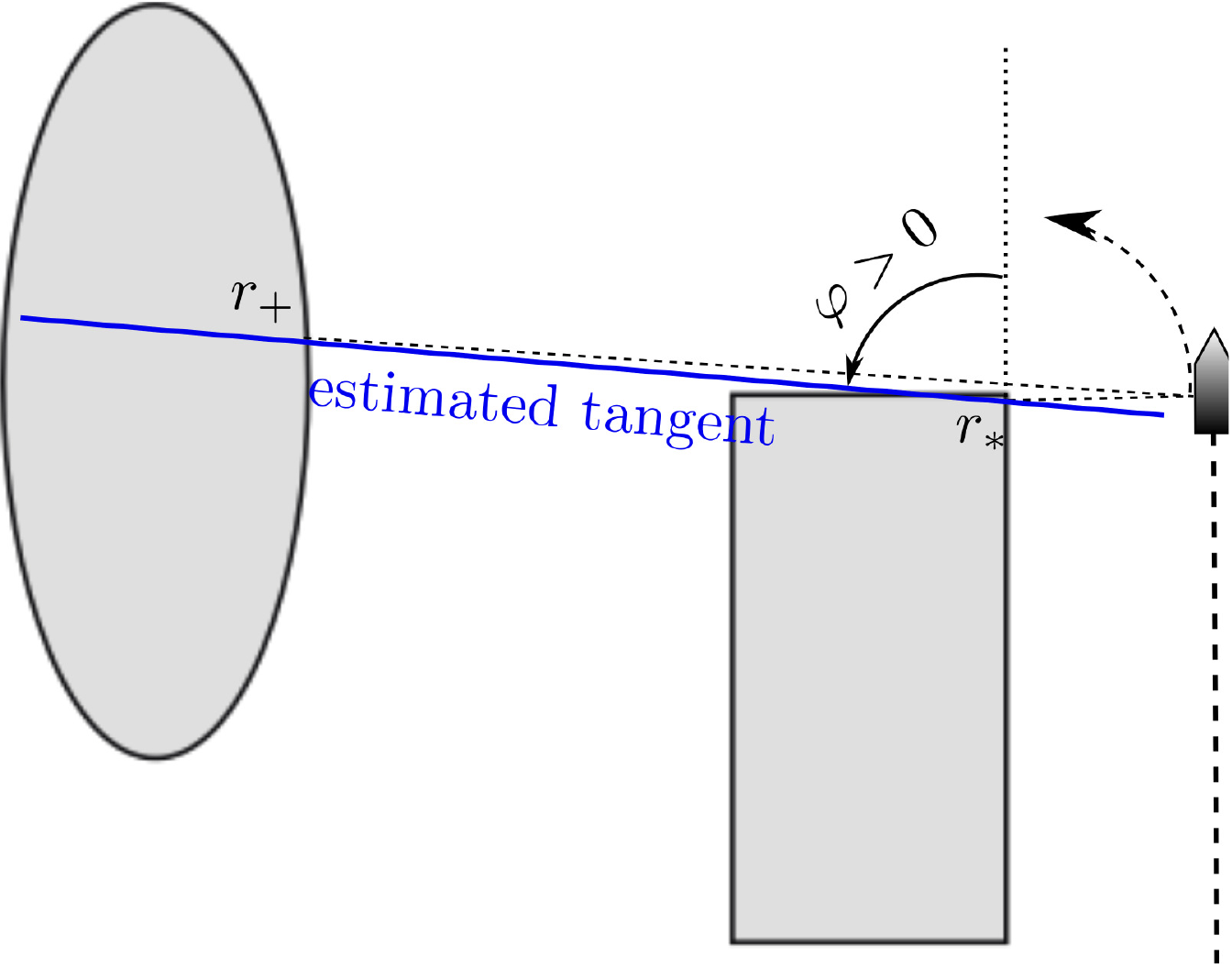}}}
  \subfigure[]{\scalebox{0.4}{\includegraphics{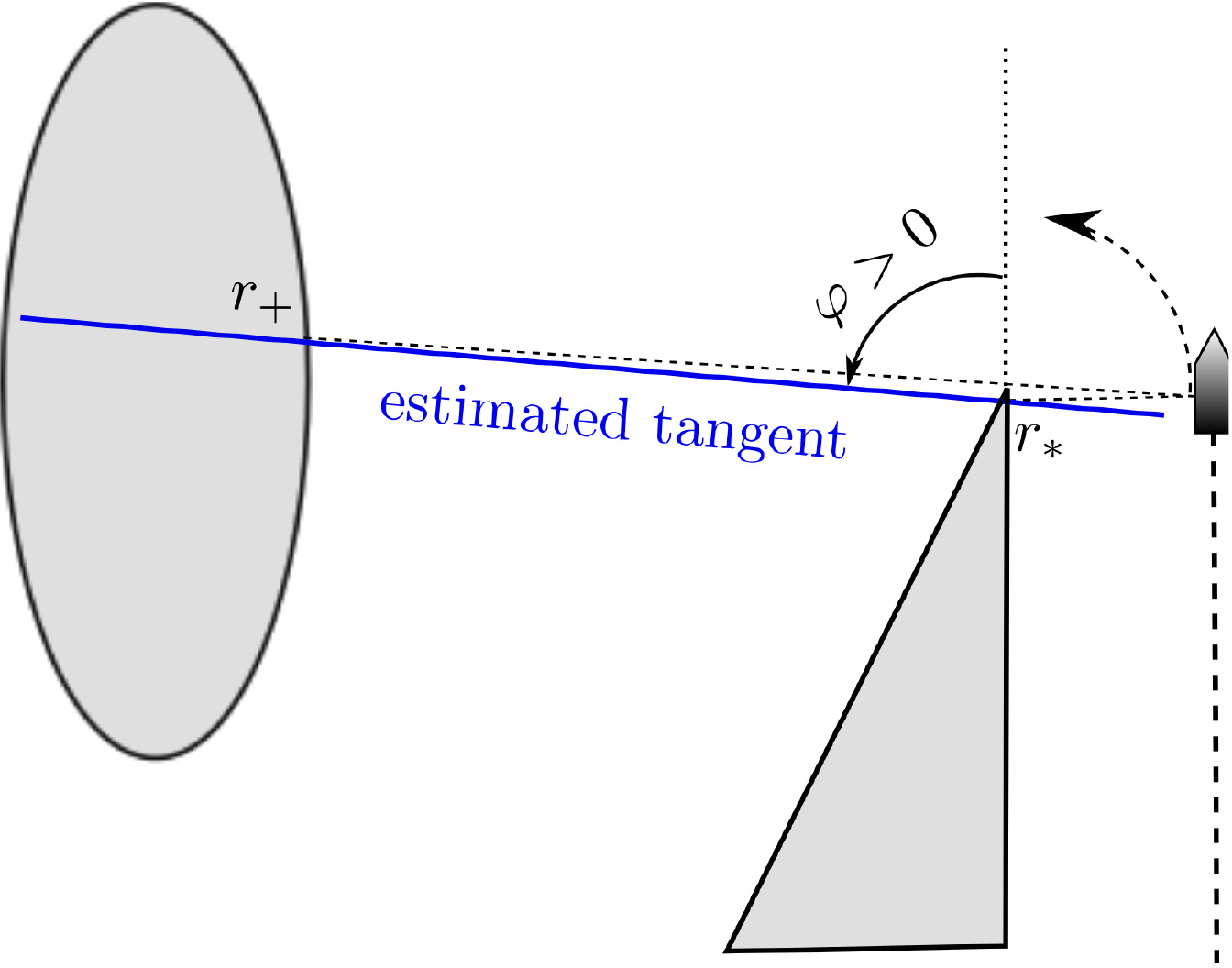}}}
  \caption{Rotation about $\bldr_\ast$.}
  \label{fbf:fig:rotation}
\end{figure}

This rough analysis gives a first evidence that the vehicle maintains
following the boundary of the obstacle at hand, as is desired.

\begin{figure}[ht]
  \centering
  \subfigure[]{\scalebox{0.4}{\includegraphics{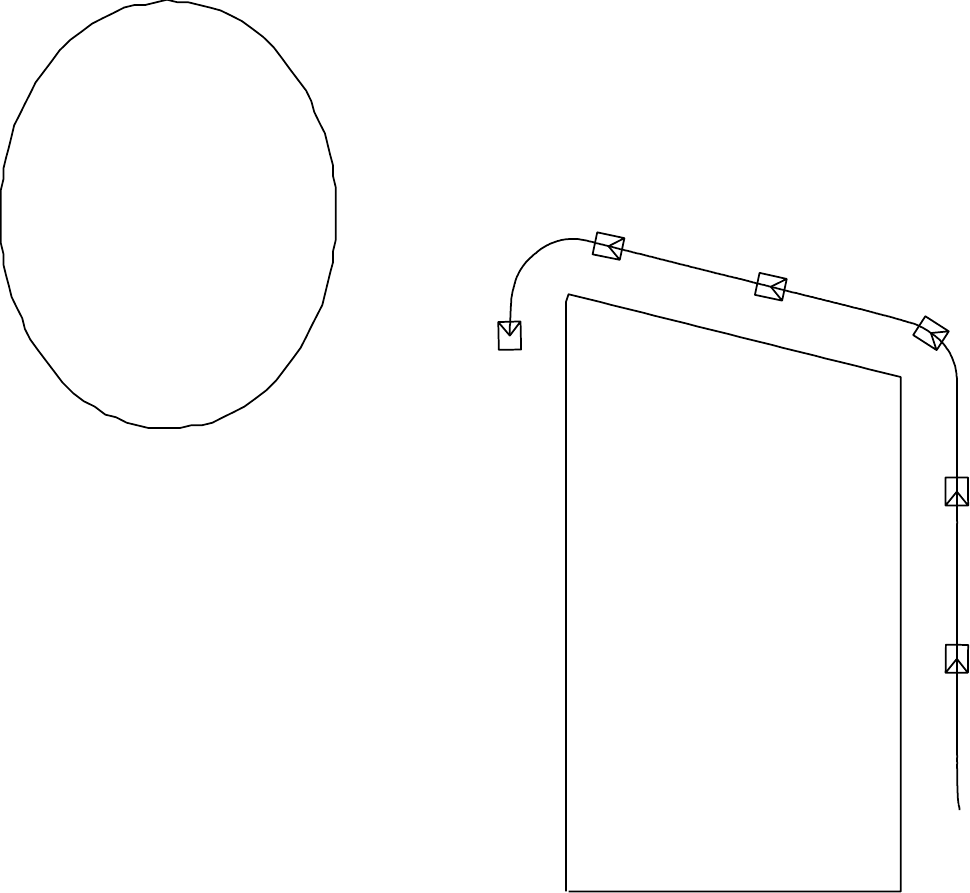}}}
  \subfigure[]{\scalebox{0.4}{\includegraphics{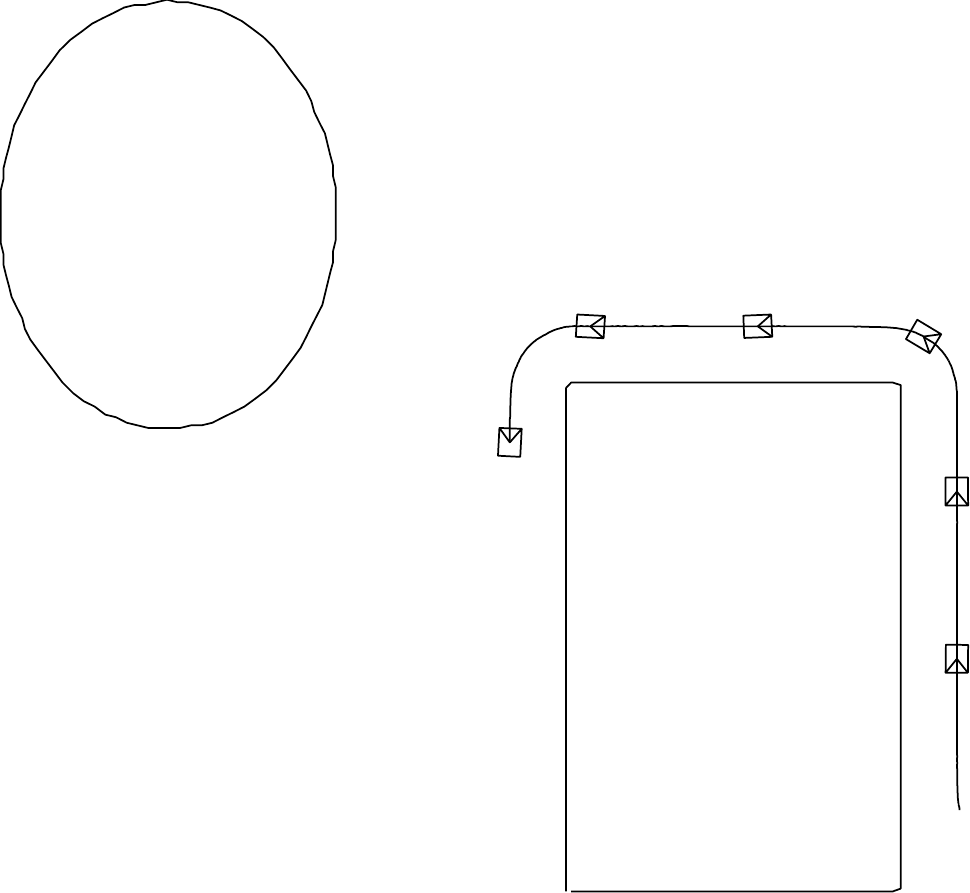}}}
  \subfigure[]{\scalebox{0.4}{\includegraphics{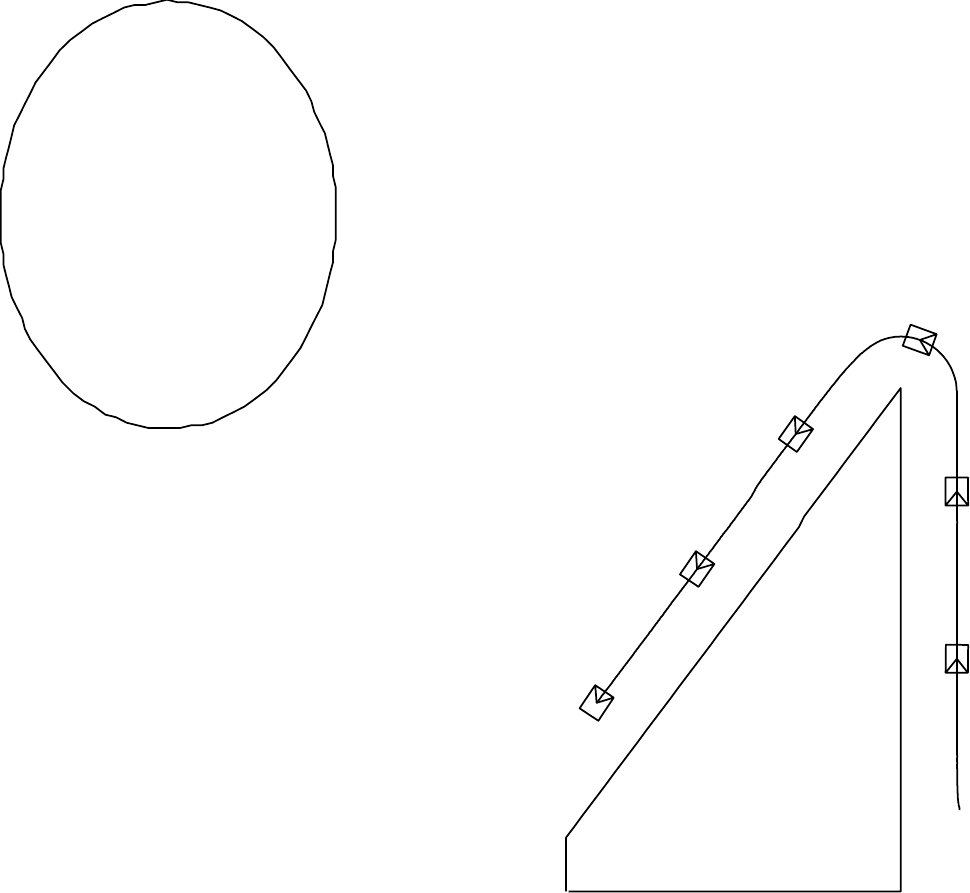}}}  
  \subfigure[]{\scalebox{0.4}{\includegraphics{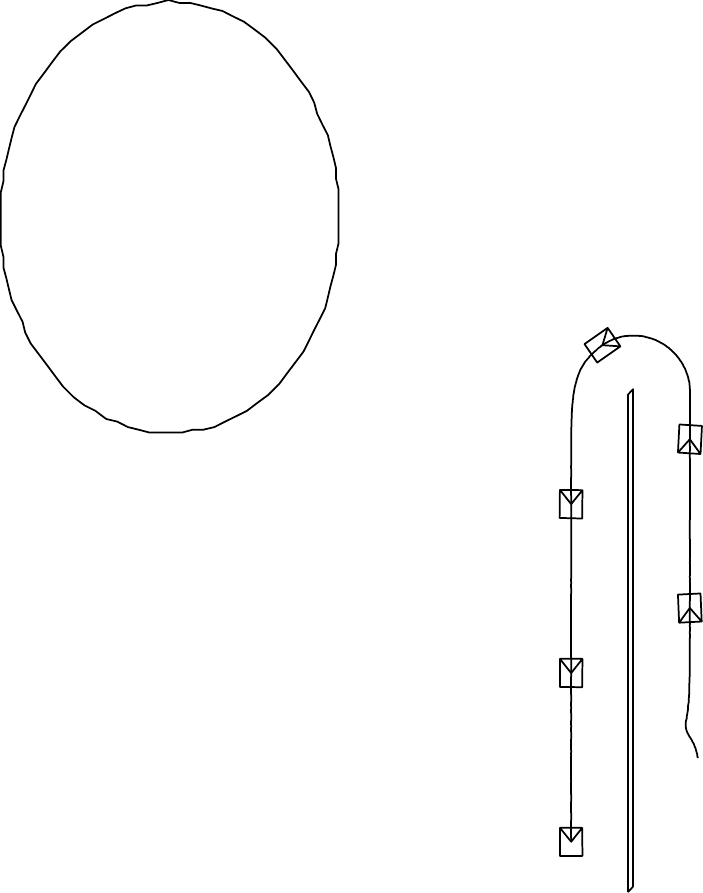}}}
  \caption{Following a fractured boundary in a cluttered environment.}
  \label{fbf:fig:bypassfrac}
\end{figure}

Similar results with a slightly worse performance are obtained in the
case where the `preceding' ray $R_+$ is replaced by a ray $R_-$
slightly behind $R_\ast$, and the distance $d$ in the perpendicular
direction is computed as $d:=\min \{d_\ast, d_-\}$, where $d_-$ and
$d_\ast$ are the distances along the respective rays; see
Fig.~\ref{fbf:fig:bypassfrac1}. The idea to compute the distance to
the nearest obstacle as the minimum of available distances along two
close rays conforms to common sense. If it is illogically set to the
maximum of them, the vehicle correctly passes obtuse fractures but for
right and acute ones, starts circular motion about the farthest
incidence point in accordance with Eq.\eqref{fbf:c.a}, thus switching to a
bypass of the competing obstacle; see Fig.~\ref{fbf:fig:bypassfrac2}.

\begin{figure}[ht]
  \centering
  \subfigure[]{\scalebox{0.4}{\includegraphics{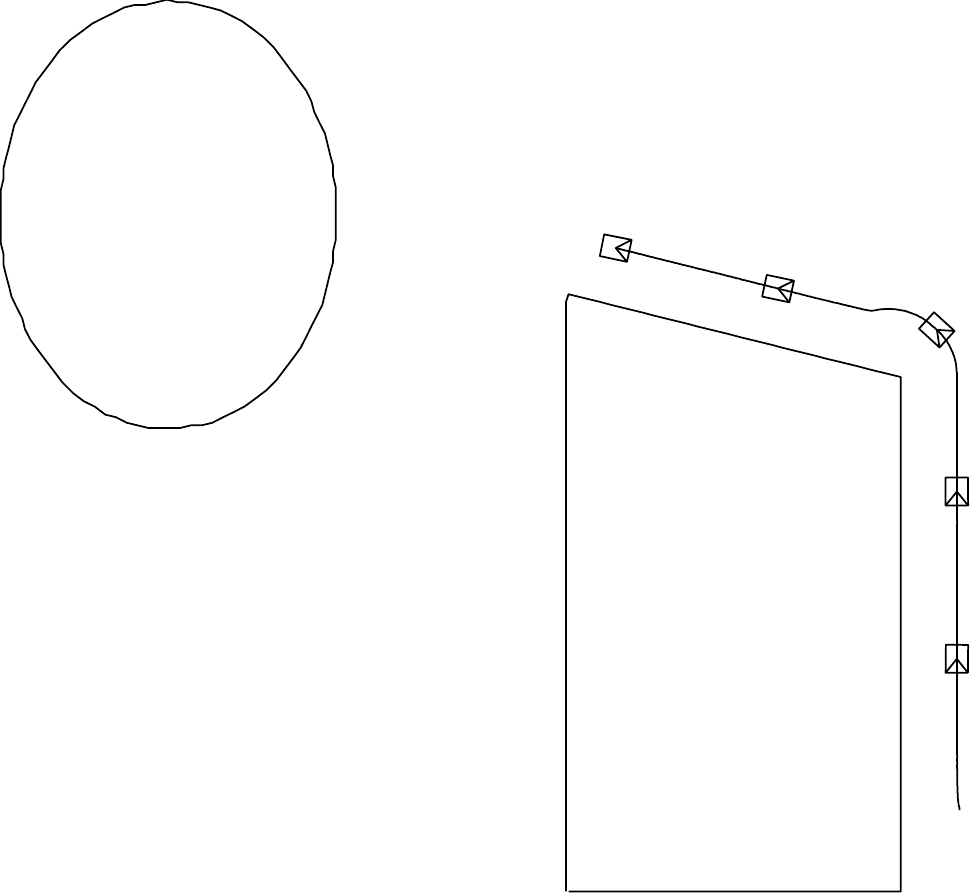}}}
  \subfigure[]{\scalebox{0.4}{\includegraphics{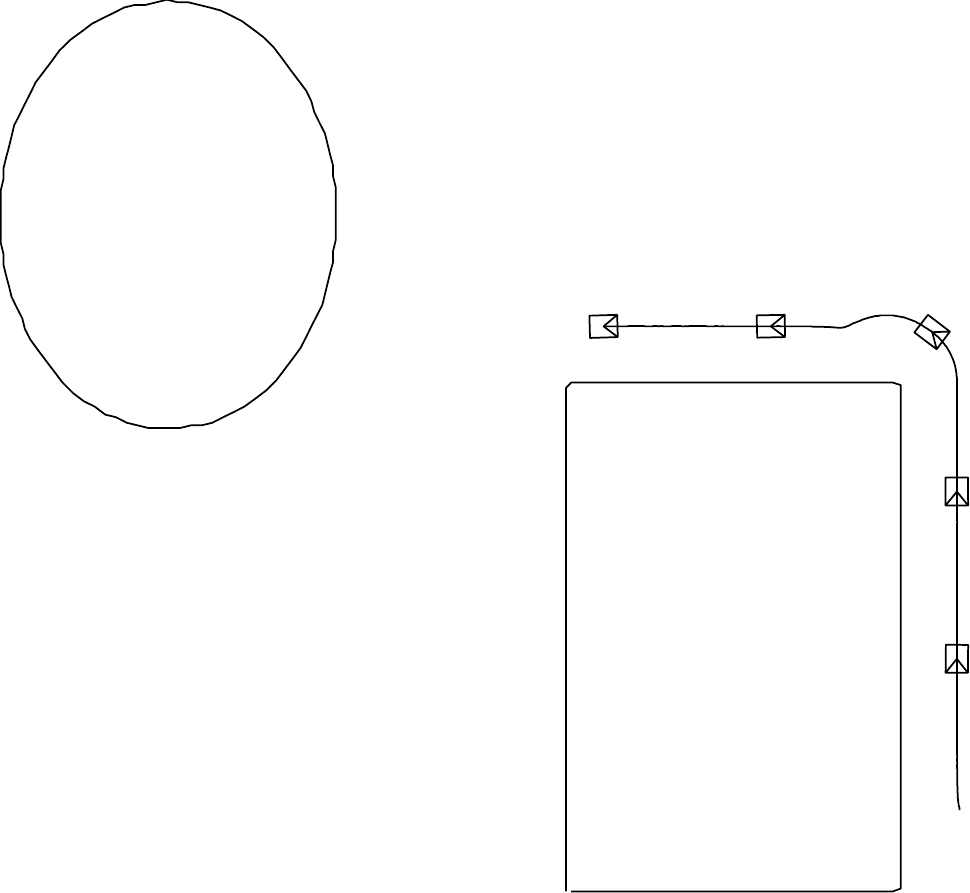}}}
  \subfigure[]{\scalebox{0.4}{\includegraphics{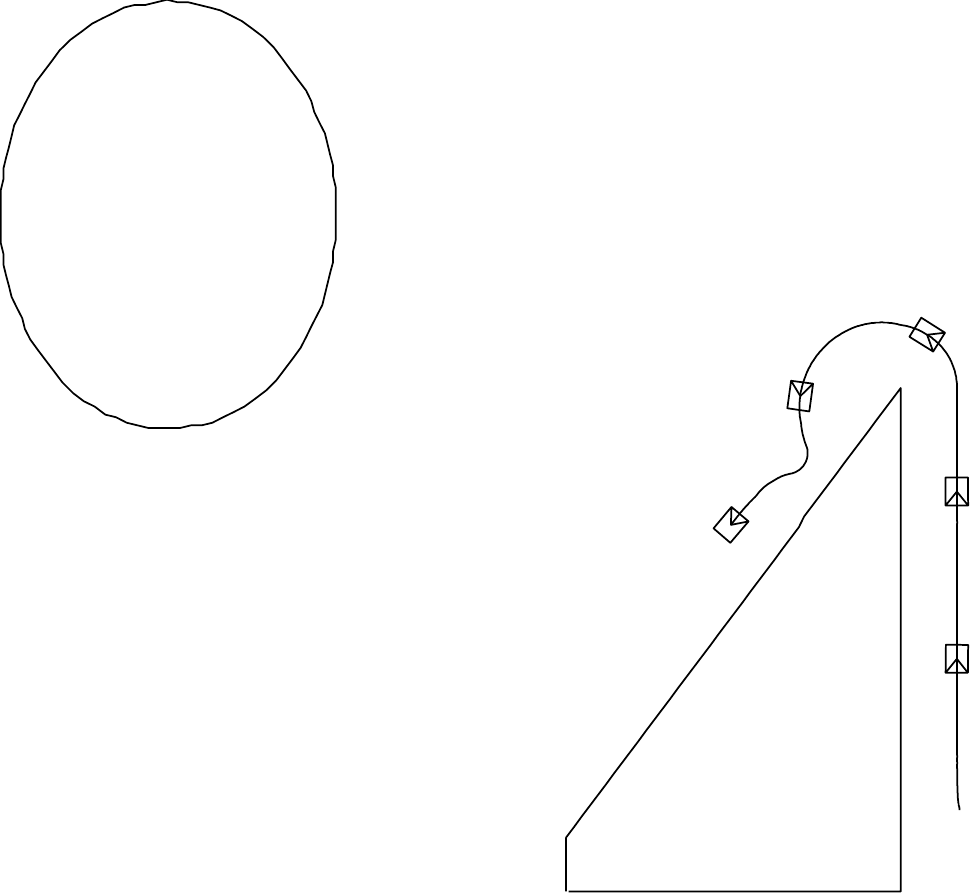}}}
  \caption{Following a fractured boundary with the use of the
    `following' ray $R_-$.}
  \label{fbf:fig:bypassfrac1}
\end{figure}

\begin{figure}[ht]
  \centering
  \subfigure[]{\scalebox{0.4}{\includegraphics{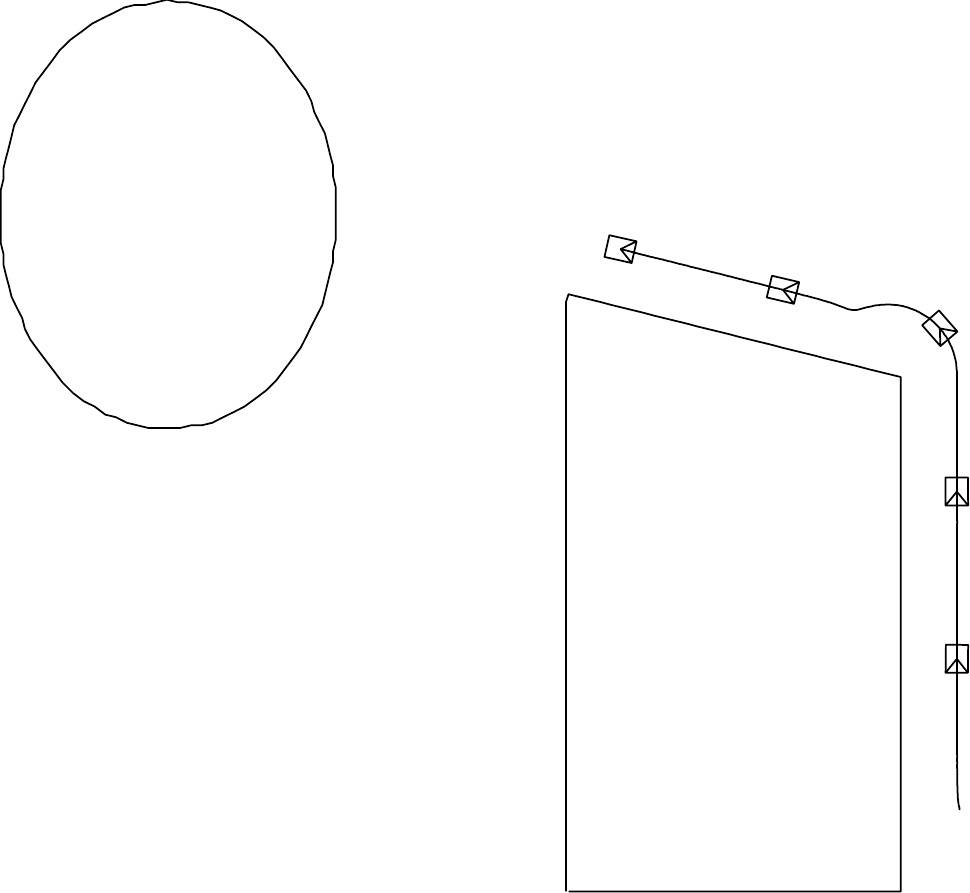}}}  
  \subfigure[]{\scalebox{0.4}{\includegraphics{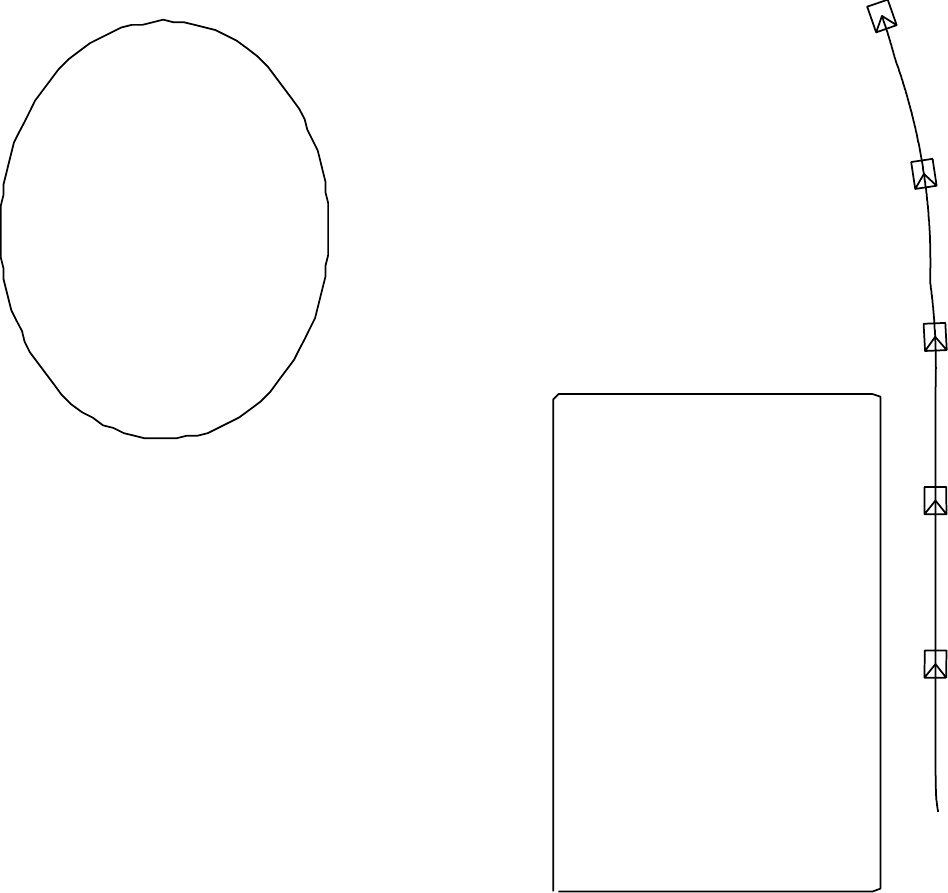}}}
  \subfigure[]{\scalebox{0.4}{\includegraphics{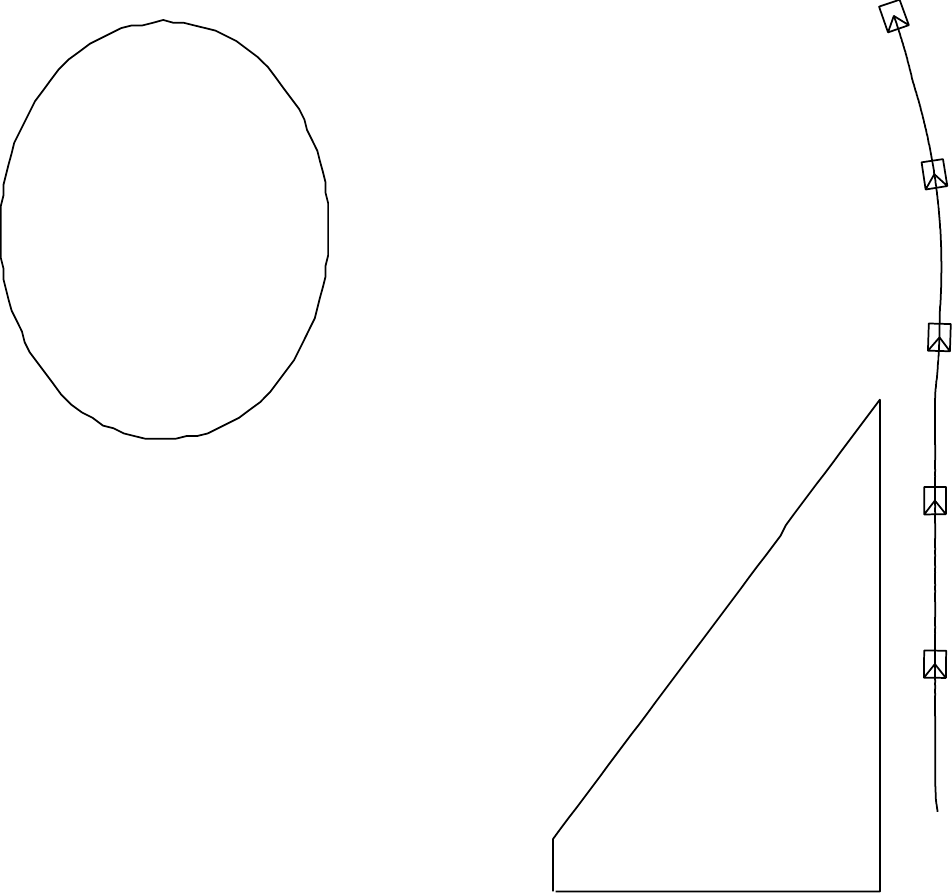}}}
  \caption{Following a fractured boundary in the case where
    illogically $d:= \max\{d_\ast,d_-\}$.}
  \label{fbf:fig:bypassfrac2}
\end{figure}

Simulations were also carried out to compare the performance of the
proposed navigation method with that from \cite{Kim2009journ7}. The
latter algorithm employs the boundary curvature. By following
\cite{Kim2009journ7}, a geometric estimate of this curvature was used
in these simulations. It results from drawing a circle through three
boundary points detected by three close rays; the main
perpendicular ray and two auxiliary rays on either side of the main
one. The lengths to the boundary along every detection ray were
corrupted by iid noises uniformly distributed over the very small
interval $[-5 mm, 5 mm]$. The same noises affected the two-point
approximation of the tangent to the boundary for the control law
proposed in this chapter. To enhance the difference, the both methods
were equally challenged by not pre-filtering the noisy measurements.
The simulation scenario was the obtuse corner from
Fig.~\ref{fbf:fig:bypassfrac}.  The parameters of the controller from
\cite{Kim2009journ7} were taken to be $\kappa_M := 0.6; \epsilon :=
0.4; \epsilon_2 = 0.2; \mu = 1; \mu_2 = 5; \mu_3 = 5$, which meets the
guidelines given in \cite{Kim2009journ7}. All simulation tests started
at a common position, which corresponds to perfectly following the
boundary at the required distance $d_0$. The tests were repeated $500$
times with individual realizations of the noises; for each test, the
maximal (over the experiment) deviation from the desired distance to
the boundary $\max_t |d(t) - d_0|$ was recorded.
\par
The overall results are demonstrated by the histogram in
Fig.~\ref{fbf:fig:hist}, which shows the number of experiments with a
given maximal deviation. The displayed better performance of the
proposed method is presumably due to getting rid of the second
derivative property (the boundary curvature) in the control law, which
is particularly sensitive to both the distance measurement noises and
violations of the boundary smoothness.

\begin{figure}[ht]
  \centering
  \includegraphics[width=10cm]{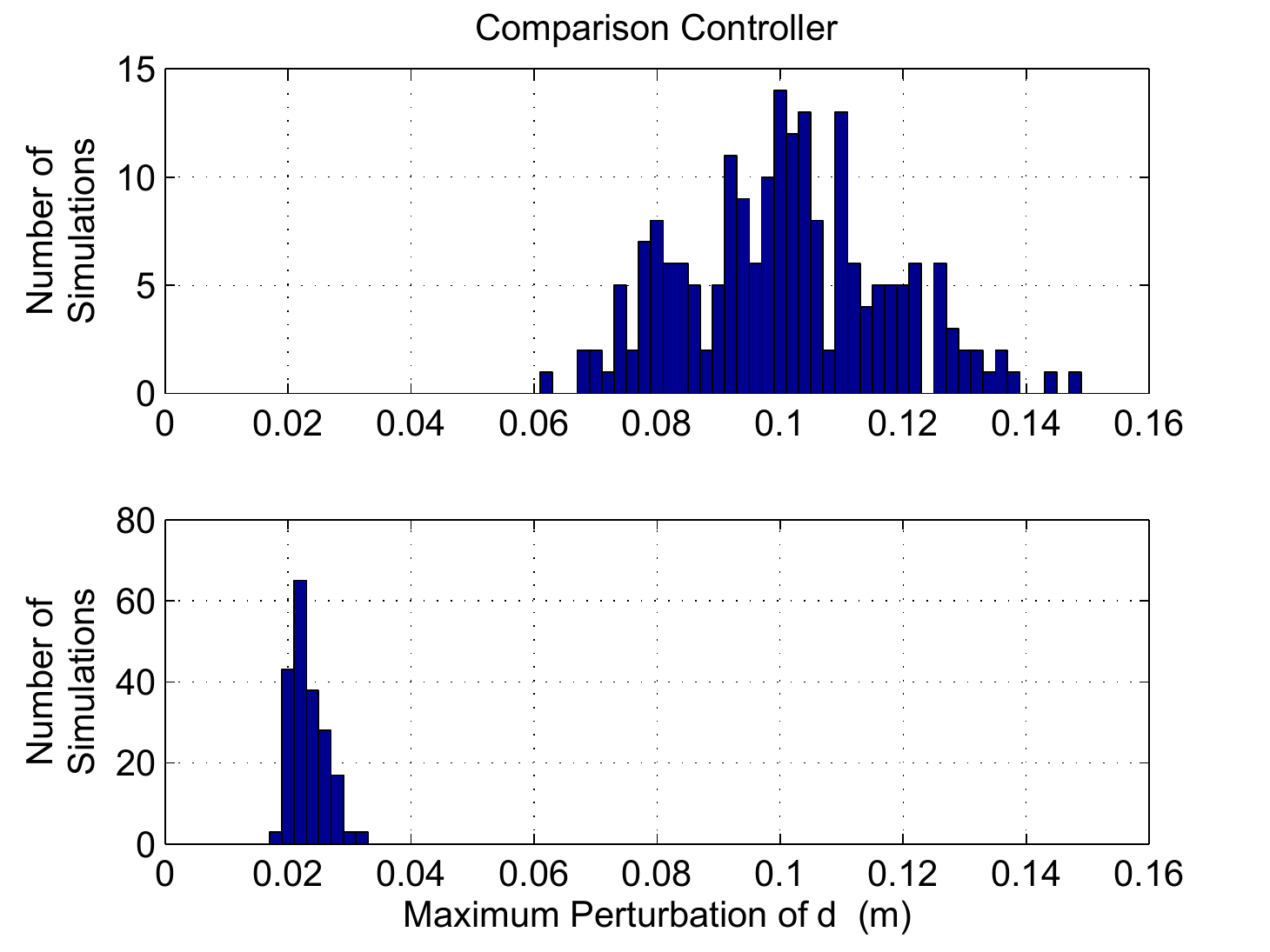}
  \caption{Distribution of the maximum distance error for two boundary
    following methods.}
  \label{fbf:fig:hist}
\end{figure}

\clearpage \section{Experiments}
\label{fbf:sec.expr}

Experiments were performed with a Pioneer P3-DX mobile robot, equipped
with a SICK LMS-200 LiDAR device. The
controller was slightly modified by continuous approximation of the
nominal law Eq.\eqref{fbf:c.a}:

\begin{equation*}
  u_{\text{act}} := \sgn (u) \cdot \min \{ |u|, \lambda |S| \}
\end{equation*}

Here $u$ is given by Eq.\eqref{fbf:c.a} and $\lambda$ is a tunable parameter;
$\lambda = 2$ in the experiments.  Continuous approximation in a
boundary layer is a common approach in practical implementation of
discontinuous control laws
\cite{Utkin1992book1,Edwards2006book3,Teimoori2010journ1a}, basically
aimed at chattering elimination.
\par
The control requested by the navigation law was forwarded as the
desired steering input to the ARIA library associated with the robot
(version 2.7.4). 
The control was updated at the rate of $0.1 s$, and the angle between
the two rays used to determine the secant that approximates the
tangent to the obstacle boundary was set to be $10$ deg. The
parameters used in testing are shown in Table~\ref{fbf:fig:paramexp}.

\begin{table}[ht]
  \centering
  \begin{tabular}{| l | c |}
    \hline
    $\ov{u}$ & $28.6 deg/s$ \\
    \hline
    $v$ & $0.2 ms^{-1}$ \\
    \hline
    $\mu$ & $57.3 deg$  \\
    \hline
  \end{tabular}
\hspace{10pt}
\begin{tabular}{| l | c |}
    \hline
    $\gamma$ & $171.9 deg/m$  \\
    \hline
    $d_{0}$ & $0.5 m$ \\
    \hline
  \end{tabular}
  \caption{Experimental parameters for fixed-sensor boundary following controller.}
  \label{fbf:fig:paramexp}
\end{table}

 The path obtained is shown in
Fig.~\ref{fbf:fig:exp} and is captured by manually marking the
position of the robot on the video frames. Both the video and
Fig.~\ref{fbf:fig:exp} demonstrate that the robot behaves as expected,
successfully circling the obstacle with both concavities and
convexities without collision and with a safety margin $\geq
0.46cm$. The related tracking measurements are displayed in
Figs.~\ref{fbf:fig:disexp} and \ref{fbf:fig:angexp}.  As is seen in
Fig.~\ref{fbf:fig:disexp}, the deviation of the distance measurement
$d$ from the desired value $0.5 m$ is always within the interval
$[-4.0cm, +6.0cm]$ and moreover, is typically in $[-1.0cm, +3.0cm]$,
except for several very short jumps out of this smaller range.
However after these jumps, the vehicle quickly recovers. Except for
the respective short periods of time, the distance error typically
does not exceed $3.0 cm$, with the dangerous deviation from the
required distance towards the obstacle being no more than $1cm$. As is
seen in Fig.~\ref{fbf:fig:angexp}, the estimate of the boundary
relative tangential angle $\varphi$ is subject to a significant amount
of noise. This is due to amplifying the noise in the distance
measurement during numerical evaluation of $\varphi$ via approximation
of the tangent by the secant between two close points of the
boundary. This approximation also appeared to be sensitive to small
angular perturbations accompanying the motion of the vehicle. Even
with these detrimental effects, the controller has demonstrated the
good ability to achieve the desired outcome.

\begin{figure}[ht]
  \centering
  \subfigure[]{\scalebox{0.35}{\includegraphics[width=\columnwidth]{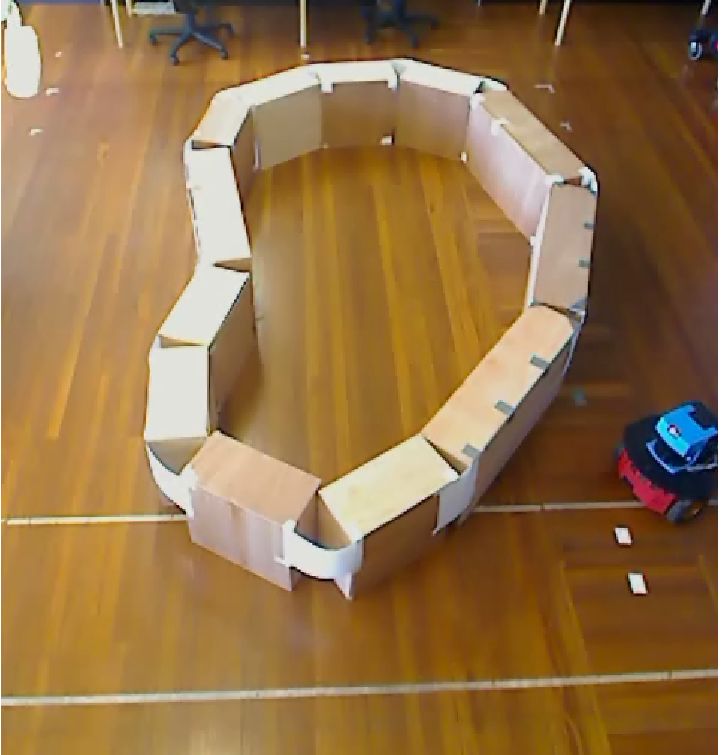}}}
  \subfigure[]{\scalebox{0.35}{\includegraphics[width=\columnwidth]{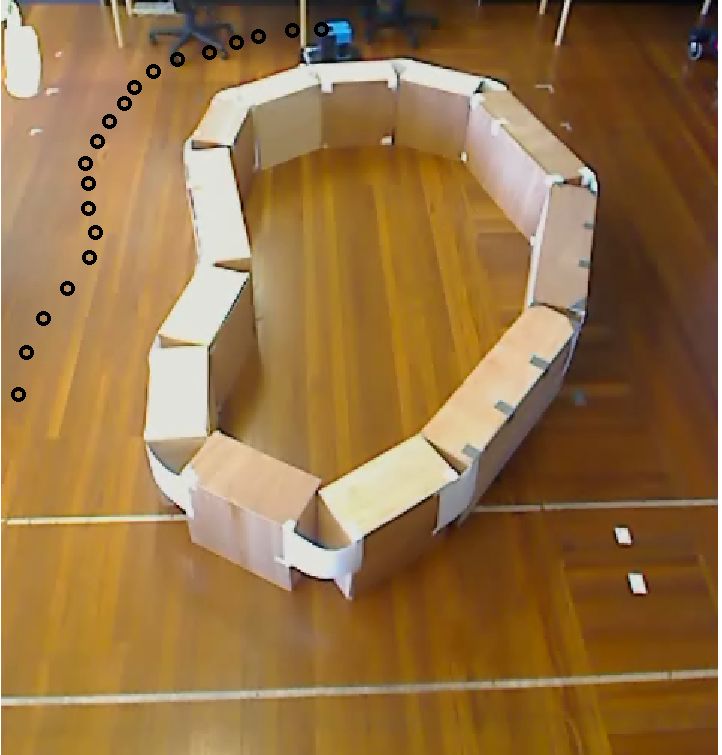}}}
  \subfigure[]{\scalebox{0.35}{\includegraphics[width=\columnwidth]{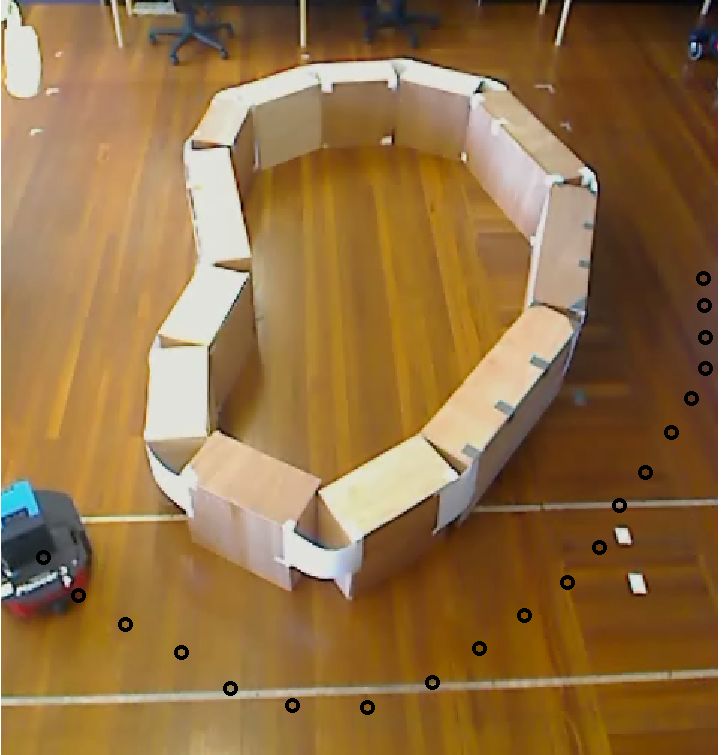}}}
  \subfigure[]{\scalebox{0.35}{\includegraphics[width=\columnwidth]{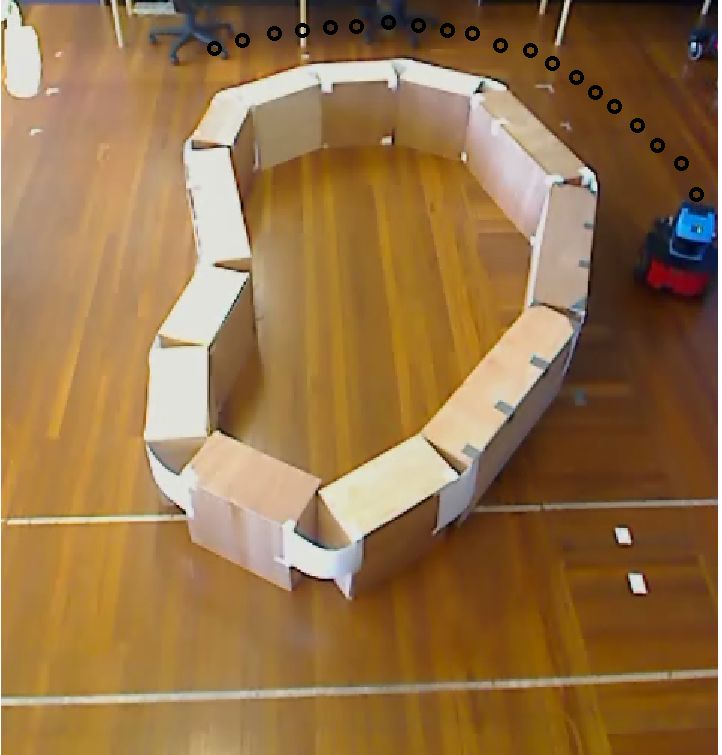}}}
  \subfigure[]{\scalebox{0.35}{\includegraphics[width=\columnwidth]{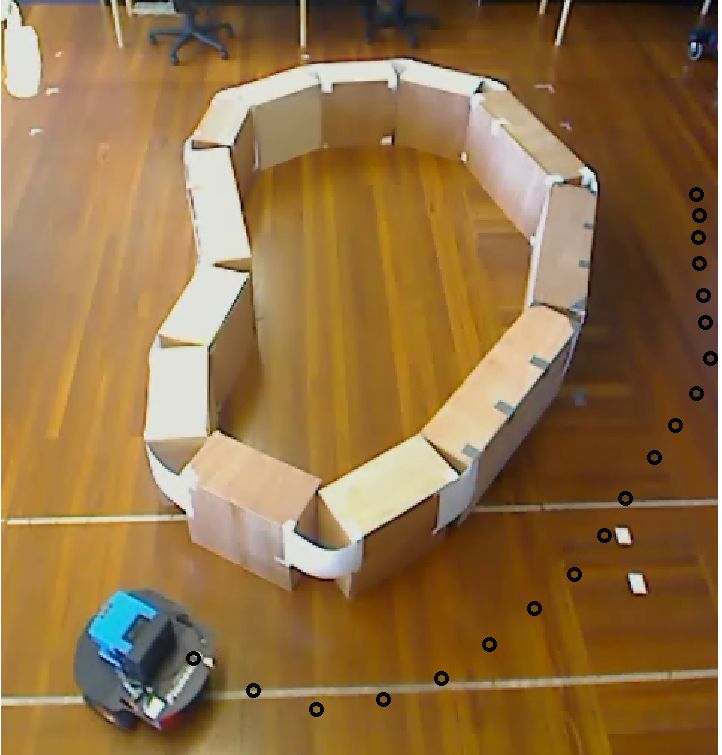}}}
  \caption{Sequence of images showing the experiment.}
  \label{fbf:fig:exp}
\end{figure}

\begin{figure}[ht]
  \centering
  \includegraphics[width=10cm]{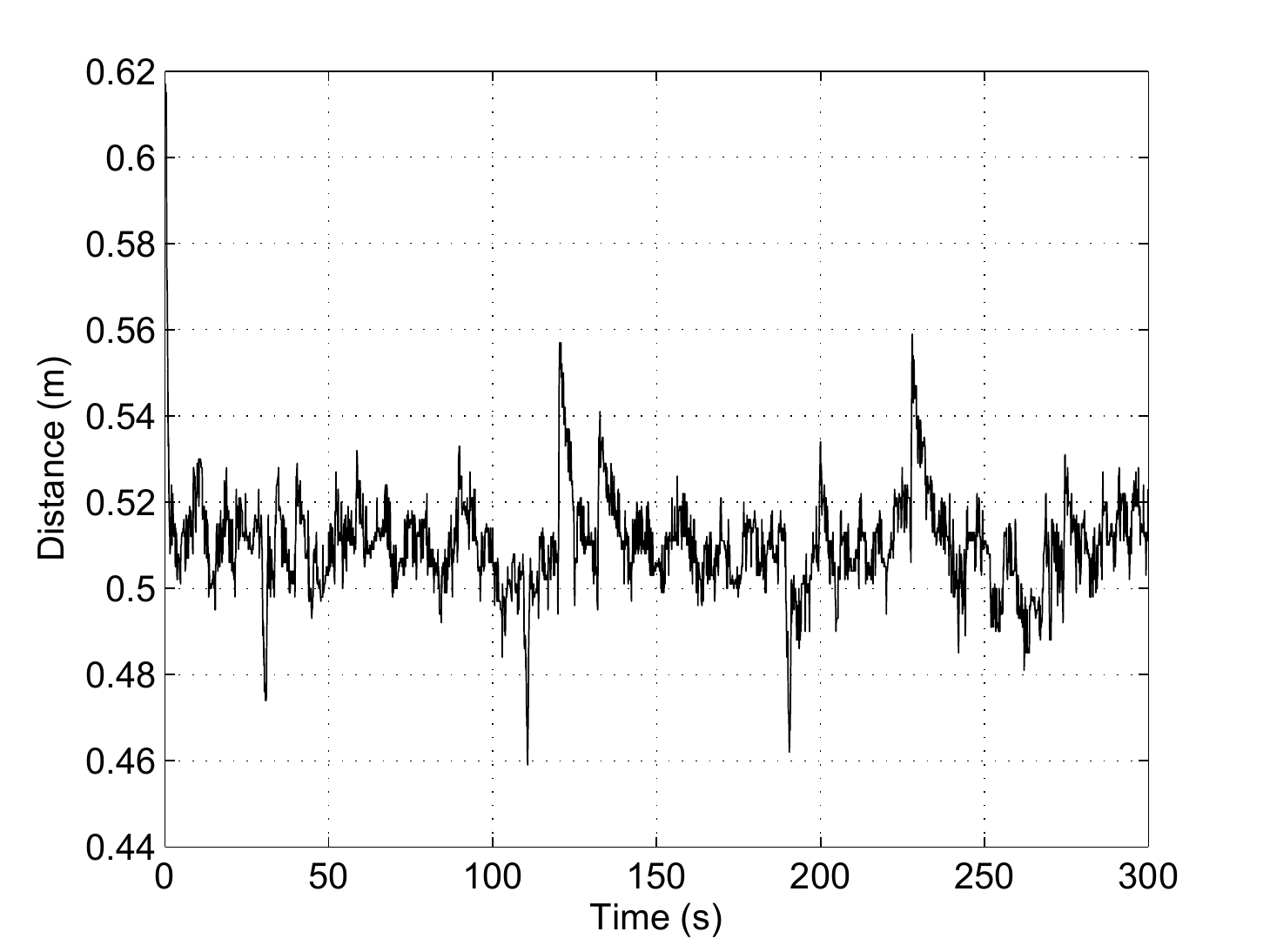}
  \caption{Distance measurement $d$ during the
    experiment.}
  \label{fbf:fig:disexp}
\end{figure}

\begin{figure}[ht]
  \centering
  \includegraphics[width=10cm]{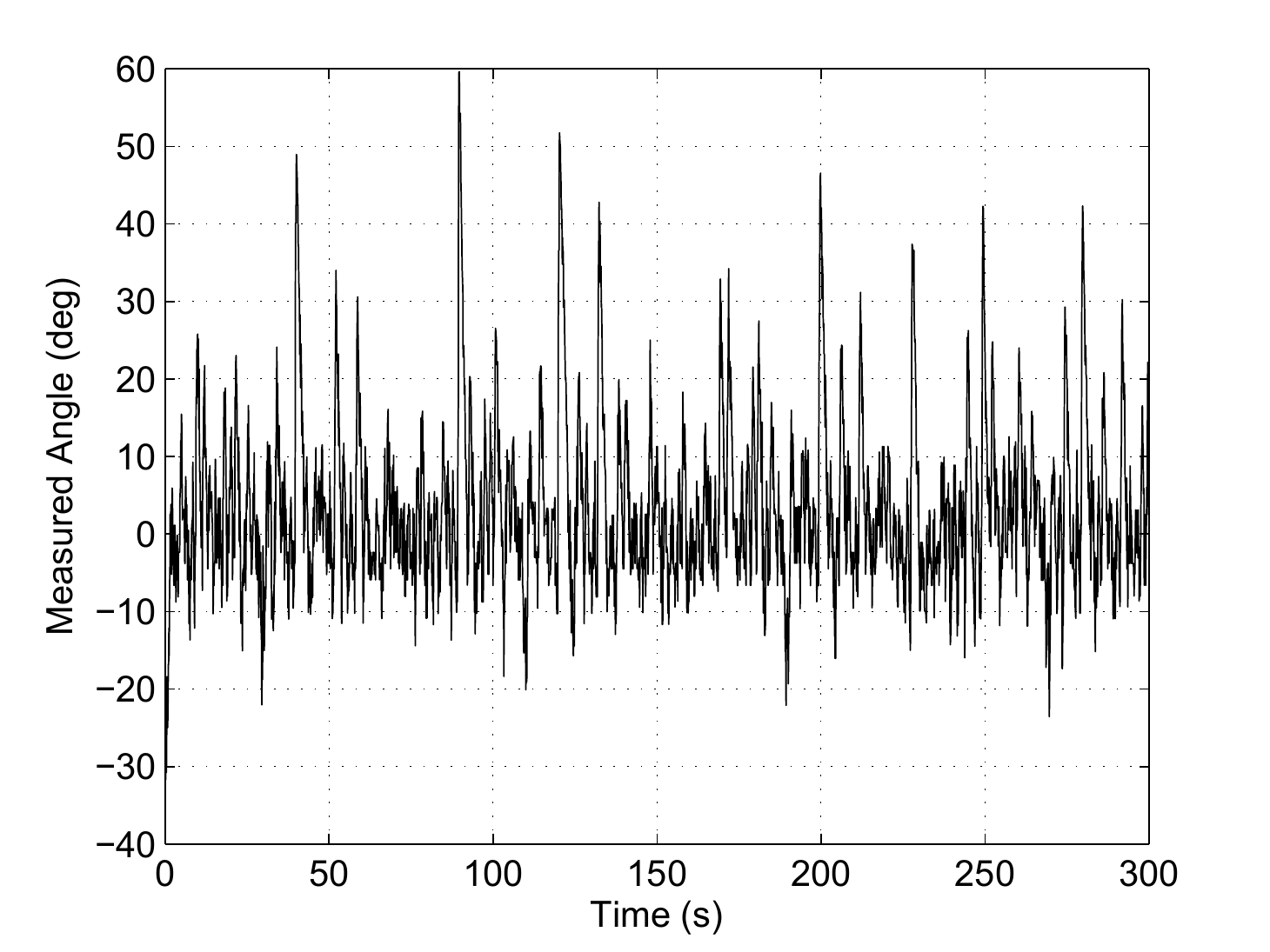}
  \caption{Estimate of the relative tangent angle
    $\varphi$ during the experiment.}
  \label{fbf:fig:angexp}
\end{figure}

The objective of the second experiment was to test the effect of a
thin obstacle for which the obstacle curvature assumptions adopted in
the theoretical part of the chapter are heavily violated at the end of
the obstacle. The results are shown in
Figs.~\ref{fbf:fig:expalt}--\ref{fbf:fig:angexpalt}.  In accordance
with explanatory remarks on Fig.~\ref{fbf:fig:bypassfrac}, passing the
obstacle end is accompanied with detection of a different and farther
obstacle, as can be seen in Fig.~\ref{fbf:fig:disexpalt}, which
displays the farthest distance detected within the employed narrow
beam of detection rays for illustration purposes. However due to the
reasons disclosed at the end of Sec.~\ref{fbf:sec.sim}, the robot
successfully copes with the thin obstacle and behaves as expected.

\begin{figure}[ht]
  \centering
  \subfigure[]{\scalebox{0.38}{\includegraphics[width=\columnwidth]{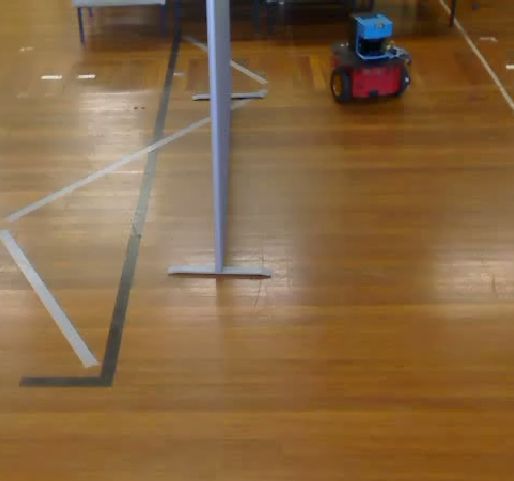}}}
  \subfigure[]{\scalebox{0.38}{\includegraphics[width=\columnwidth]{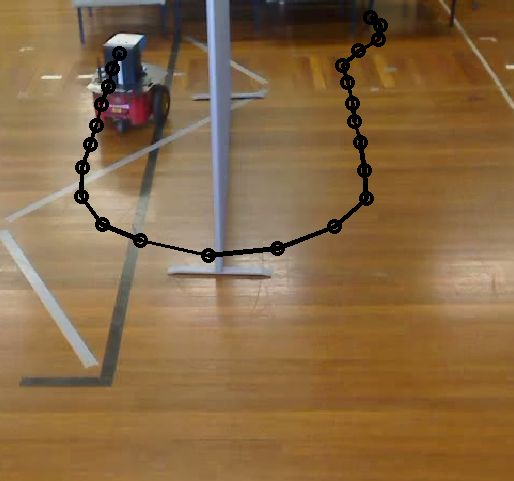}}}
  \caption{Sequence of images showing the experiment (thin obstacle).}
  \label{fbf:fig:expalt}
\end{figure}

\begin{figure}[ht]
  \centering
  \includegraphics[width=10cm]{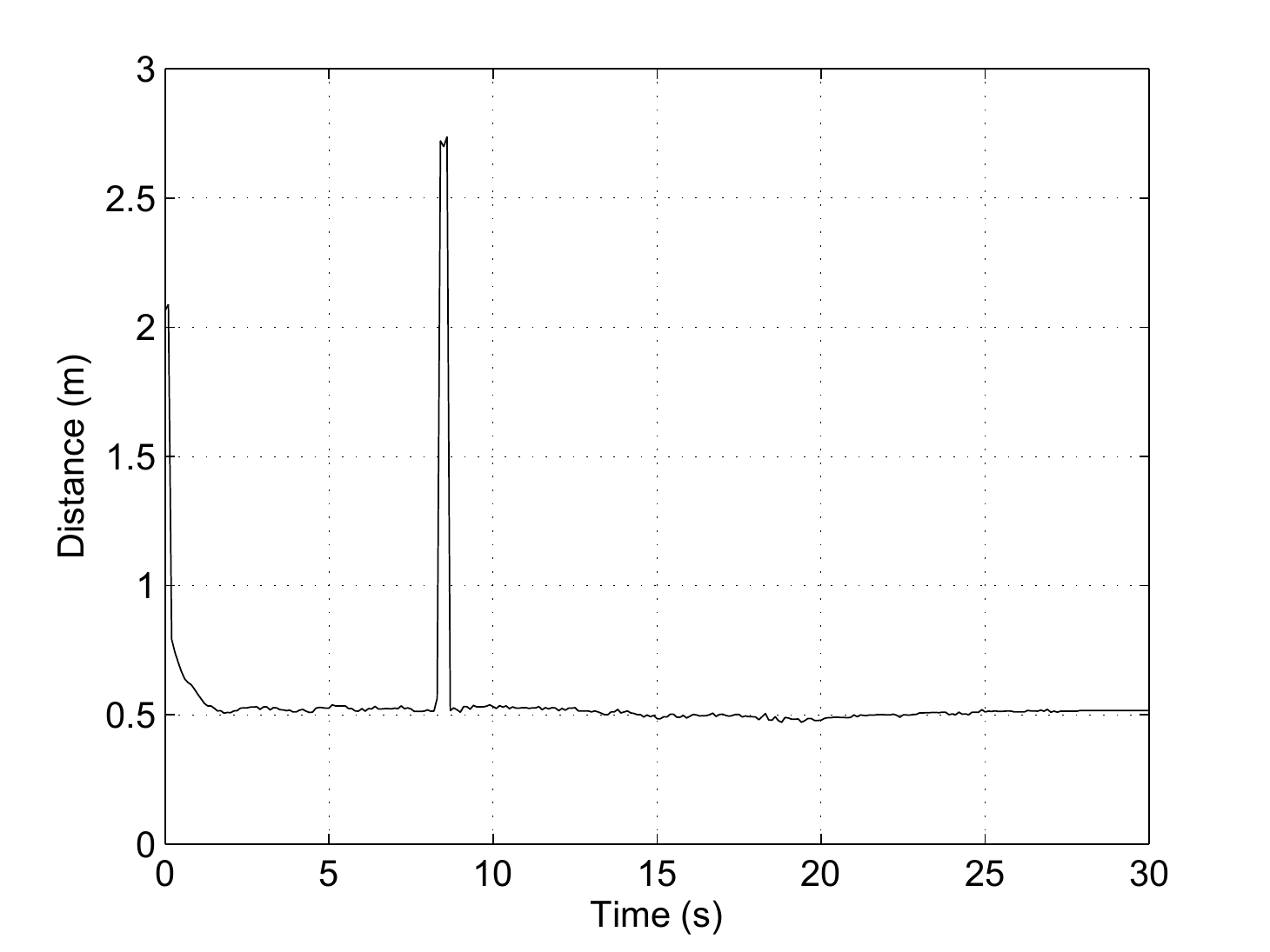}
  \caption{Distance measurement during the experiment
    with a thin obstacle.}
  \label{fbf:fig:disexpalt}
\end{figure}

\begin{figure}[ht]
  \centering
  \includegraphics[width=10cm]{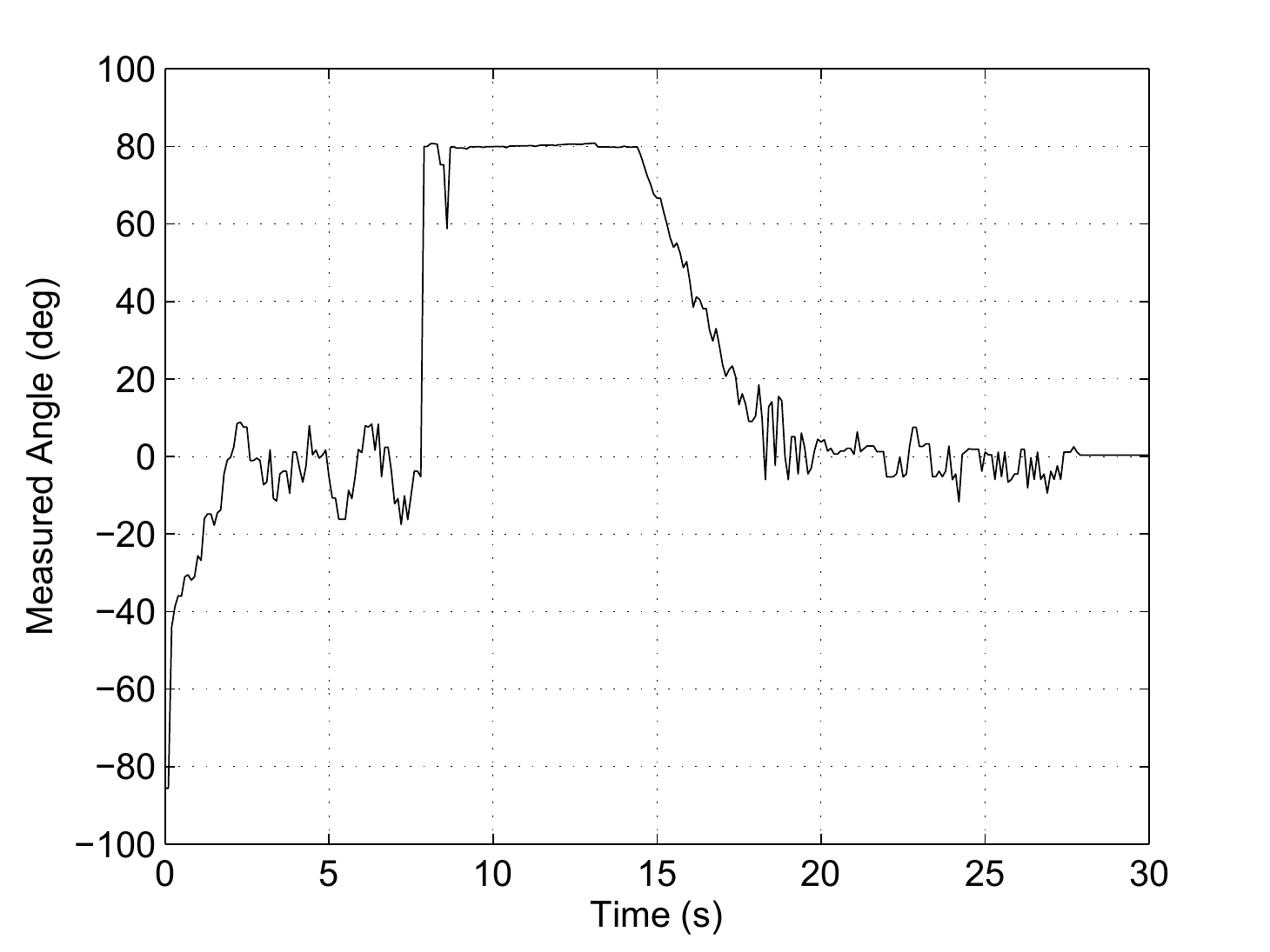}
  \caption{Estimate of the relative tangent angle
    $\varphi$ during the experiment with a thin obstacle.}
  \label{fbf:fig:angexpalt}
\end{figure}

Figs.~\ref{fbf:fig:expclut}--\ref{fbf:fig:angexpclut} present typical
results of experiments in more cluttered environments. As compared
with the previous experiments, the speed of the robot was reduced to
$0.05 ms^{-1}$ to make its minimal turning radius Eq.\eqref{fbf:Rmin} smaller
less than spacing between the obstacles. It can be seen that the robot
correctly navigates around the selected obstacle despite the presence
of the others.

\begin{figure}[ht]
  \centering
  \subfigure[]{\scalebox{0.38}{\includegraphics[width=\columnwidth]{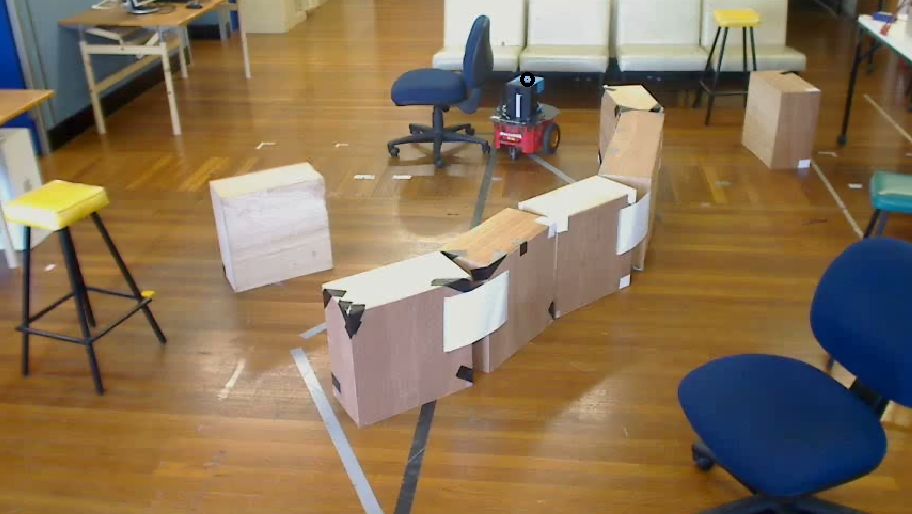}}}
  \subfigure[]{\scalebox{0.38}{\includegraphics[width=\columnwidth]{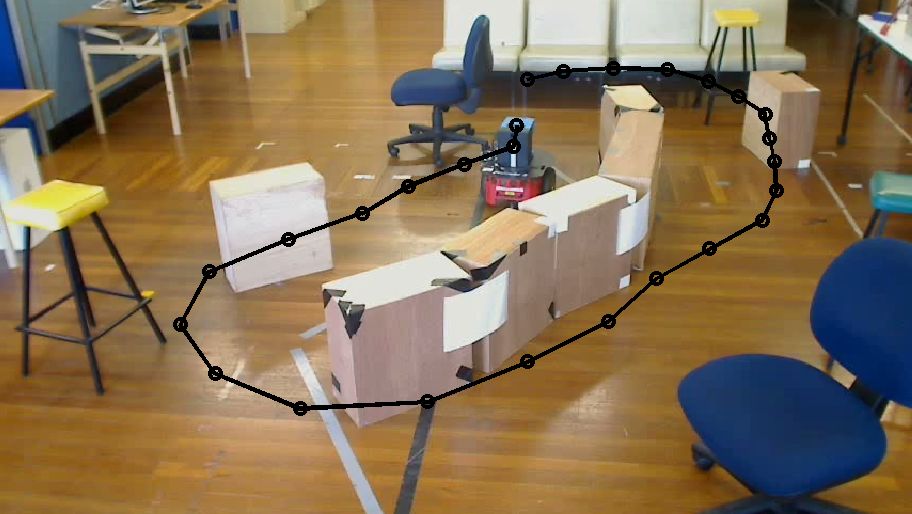}}}
  \caption{Sequence of images showing the experiment (cluttered environment).}
  \label{fbf:fig:expclut}
\end{figure}

\begin{figure}[ht]
  \centering
  \includegraphics[width=10cm]{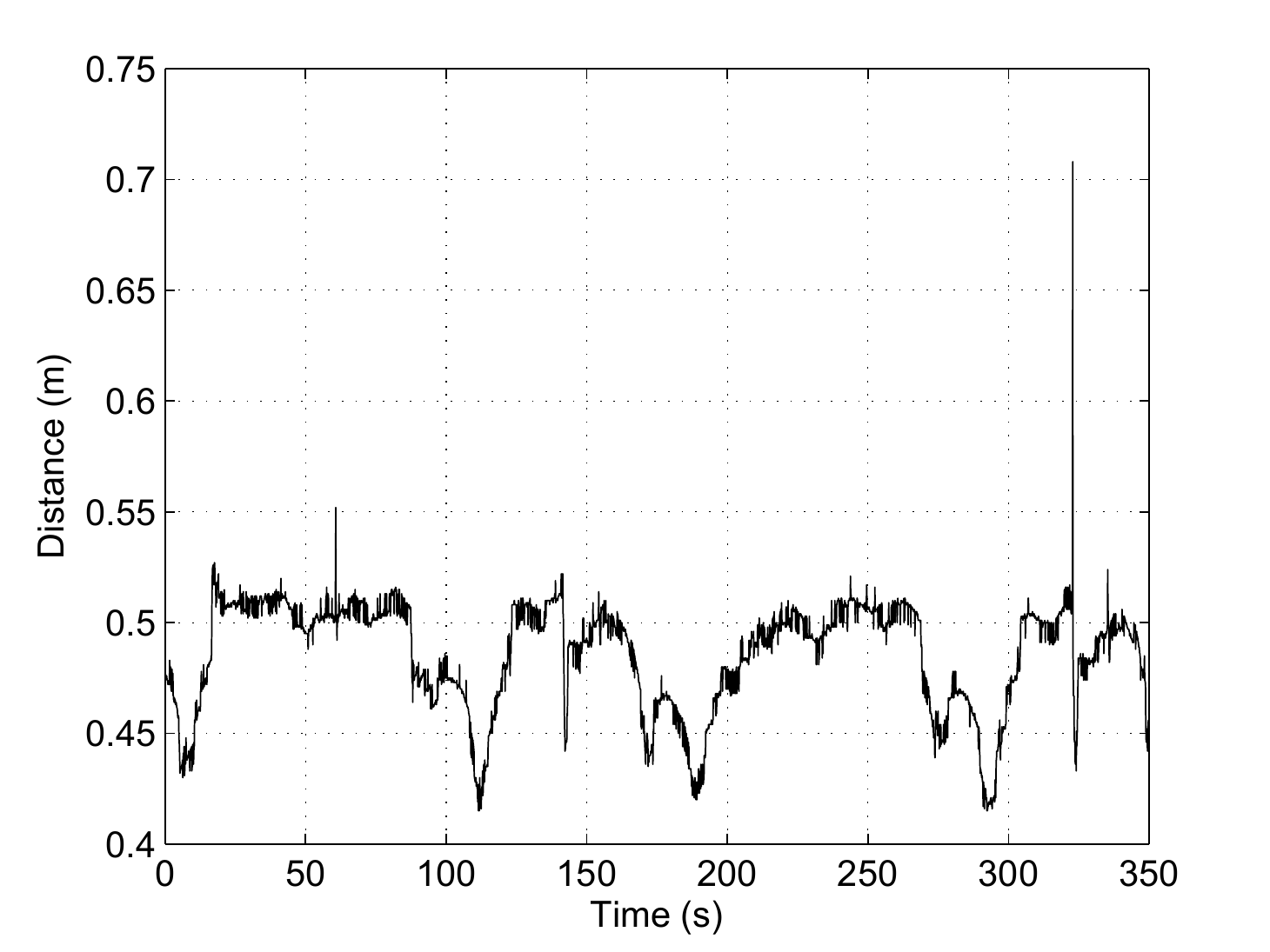}
  \caption{Distance measurement during the experiment
    in the cluttered environment.}
  \label{fbf:fig:disexpclut}
\end{figure}

\begin{figure}[ht]
  \centering
  \includegraphics[width=10cm]{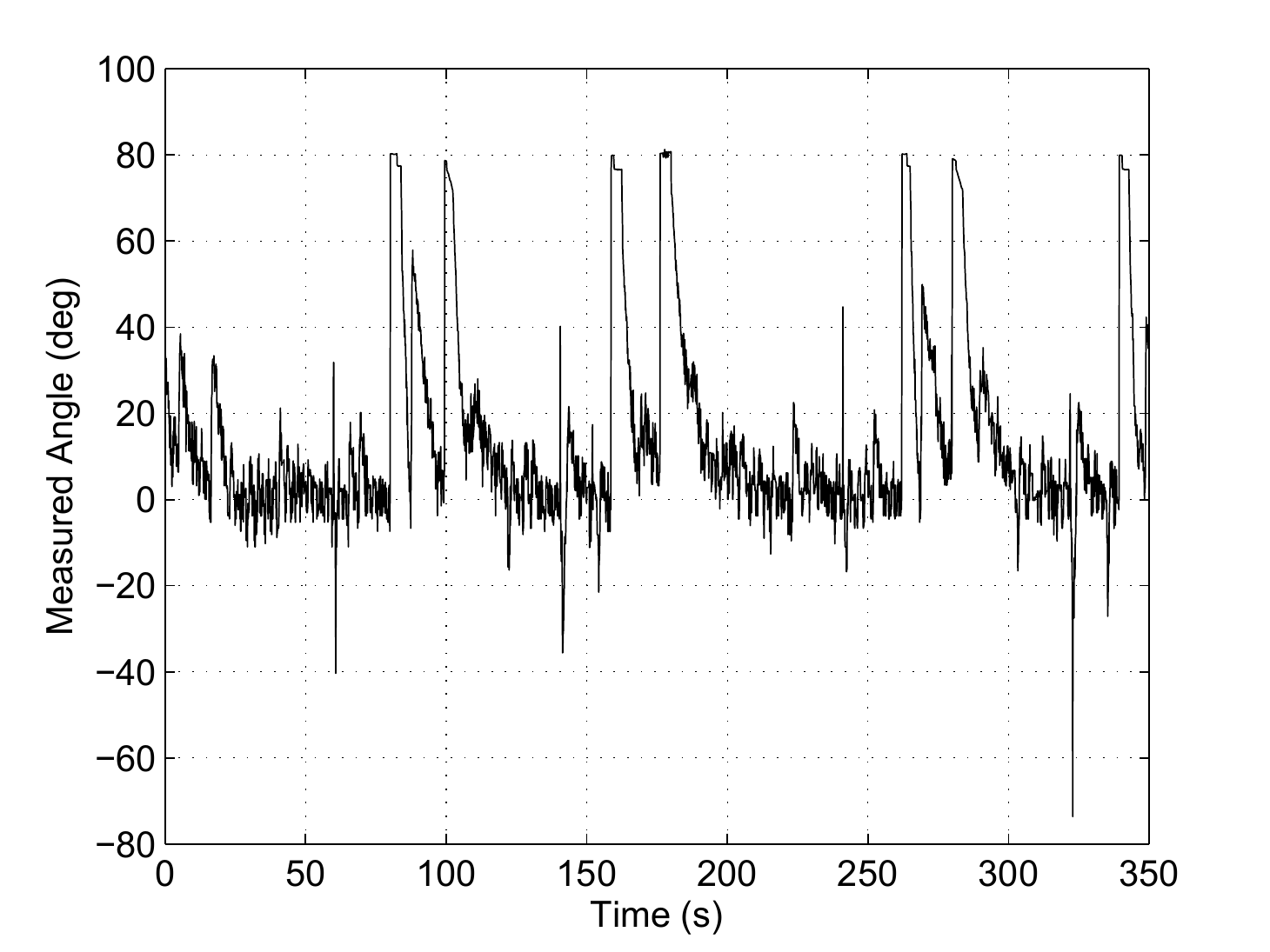}
  \caption{Estimate of the relative tangent angle
    $\varphi$ during the experiment in the cluttered environment.}
  \label{fbf:fig:angexpclut}
\end{figure}

\clearpage \section{Summary}
\label{fbf:sec.concl}
In this chapter, a novel approach for navigating an unmanned
unicycle-like vehicle along an obstacle boundary is proposed. It deals
with the situation where knowledge of the boundary is related to a
single detection ray directed perpendicularly to the vehicle
centerline. A sliding mode navigation law is proposed, which is able
to drive the vehicle at a fixed distance from this boundary. Computer simulations and
experimental results with a real wheeled robot confirm the viability
of the method.

\chapter{Extremum Seeking Navigation in a Scalar Field}
\label{chapt:ext}

In this chapter, a method for driving a vehicle to the maximal point
of a scalar environmental field is proposed. The main advantage of the
proposed method is that it requires estimation of neither the spatial
or temporal derivatives of the field value, yet convergence to the
maximal point can be proven under certain assumptions. Simulations and
experimental results are given to confirm the viability of the
proposed method.

\section{Introduction}

The chapter addresses the problem of driving a single robot to the
maximizer of an unknown scalar environmental field. For example, this may be
thermal, magnetic, electric, or optical field; or concentration of a
chemical, physical, or biological agent.  Some examples of missions where this
problem is of interest include environmental studies, geological
exploration, detecting and localizing the source of hazardous
chemical, vapor, radiation emission, or also the pursuit of
a moving target. In the last case, the objective is to approach the
target and follow it with a pre-specified margin.  Apart from source
seeking/localizing \cite{ZhArGhSiKr07,PoNe96,CoKr09}, this problem is
often referred to as gradient climbing/ascent
\cite{OgFiLe04,BiAr07,BaLe02,BuYoBrSi96}, which highlights the common
kinematic control paradigm -- try to align the velocity vector of the
robot with the field gradient.
\par
Recent surveys on extremum seeking control methods and
algorithms are available \cite{KoRu08,DoPeGu11}.  A good
deal of related research was concerned with gradient climbing based on
direct on-line gradient estimation. This approach is especially
beneficial for mobile sensor networks thanks to collaborative field
measurements in many locations and data exchange
\cite{PoNe96,OgFiLe04,BiAr07,BaLe02,GaPa04,MoCa10}.  However, even in
this scenario, data exchange degradation may require each robot in the team to operate autonomously
over considerable time. Similar algorithms can be basically used for a
single robot equipped with several sensors that are distant enough
from each other and thus provide the field values at several
essentially diverse locations. In any case, multiple vehicle/sensor
scenario means complicated and costly hardware.
\par
The lack of multiple sensor data can be compensated via exploring
multiple nearby locations by ``dithering'' the position of the single
sensor during special maneuvers, which may be excited by probing
high-frequency sinusoidal
\cite{BuYoBrSi96,ZhArGhSiKr07,CoKr09,CoSiGhKr09} or stochastic
\cite{LiKr10} inputs. A similar in spirit approach is extremum seeking
by means of many robots performing biased random walks \cite{MeHeAs08}
or by two robots with access to relative positions of each other and
rotational actuation \cite{ElBr12}.  These methods rely, either
implicitly or explicitly, on systematic side exploration maneuvers to
collect rich enough data. Another approach limits the field gradient
information to only the time-derivative of the measured field value
obtained by numerical differentiation
\cite{BaLe02,MaSaTe08,BarBail08,Matveev2011journ9}, and partly employs
switching controllers \cite{MaSaTe08,Matveev2011journ9}. These give
rise to concerns about amplification of the measurement noises and
chattering. 
\par
The common feature of the previous research is that it dealt with only
steady fields. However in real world, environmental fields are almost
never steady and often cannot be well approximated by steady fields,
whereas the theory of extremum seeking for dynamic fields lies in the
uncharted territory. As a particular case, this topic includes the
problem of navigation and guidance of a mobile robot towards an
unknowingly maneuvering target based on a single measurement that
decays as the sensor goes away from the target, like the strength of
the infrared, acoustic, or electromagnetic signal, or minus the
distance to the target. Such navigation is of interest in many areas
\cite{ADB04,GS04,Matveev2011journ4}; it carries a potential to reduce
the hardware complexity and cost and improve target pursuit
reliability. A solution for such problem in the very special case of
the unsteady field -- minus the distance to an unknown moving
Dubins-like target -- was proposed and justified in
\cite{Matveev2011journ4}.  However the results of
\cite{Matveev2011journ4} are not applicable to more general dynamic
fields.
\par
Contrary to the previous research, this chapter addresses the source
seeking problem for general dynamic fields. In this context, it
justifies a new kinematic control paradigm that offers to keep the
velocity orientation angle proportional to the discrepancy between the
field value and a given linear ascending function of time, as opposed
to conventionally trying to align the velocity vector with the
gradient.  This control law is free from evaluation of any
field-derivative data, uses only finite gains instead of switching
control, and demands only minor memory and computational robot's
capacities, being reactive in its nature.  Mathematically rigorous
justification of convergence and performance of this control law is
available in the case of a general dynamic field. In particular, it may be
shown that the closed-loop system is prone to monotonic,
non-oscillatory behavior during the transient to the source provided
that the controller parameters are properly tuned. Recommendations 
are available for the choice of these parameters under which the
robot inevitably reaches the desired vicinity of the moving field
maximizer in a finite time and remains there afterwards. The
applicability and performance of the control law are confirmed via
extensive computer simulations and experiments with a real wheeled
robot.
\par
For complex dynamical systems, kinematic control is often the first
step in controller design whose objective is to generate the velocity
reference signal.  The next step is to design a controller that tracks
this signal by means of forces and torques. This two-stage design
works well in many situations, including these experiments, and is
especially popular in the face of uncertainties in the dynamics
loops. 
\par
An algorithm similar to ours has been proposed which
uses the estimated time-derivative of the measurement \cite{BarBail08}. Also,
this approach only considered a steady harmonic field, the performance during the transient and the behavior
after reaching a vicinity of the maximizer were not addressed even for
general harmonic fields, and the convergence conditions were partly
implicit by giving no explicit bound on some entities that were
assumed sufficiently large. 
\par
All proofs of mathematical statements are omitted here; they are
available in the original manuscript \cite{MaHoSa13}.
\par
The body of this chapter is organized as follows. In Sec.~\ref{ext:sec1} the problem is formally defined,
 and in Sec.~\ref{ext:sec.ass} the
the main assumptions are described. The main results are outlined in Sec.~\ref{ext:sec.mr}. Simulations and
experiments are presented in Secs.~\ref{ext:sec.simtest} and \ref{ext:sec.exper}. Finally, brief conclusions are given
in Sec.~\ref{ext:sum}.

 \section{Problem Statement}
\label{ext:sec1}

A planar point-wise robot traveling in a two-dimensional
workspace is considered.  The robot is controlled by the time-varying linear
velocity $\vec{v}$ whose magnitude does not exceed a given constant
$v$. The workspace hosts an unknown scalar time-varying field
$D(t,\boldsymbol{r}) \in \br$, where $t$ is time, $\boldsymbol{r}:=
(x,y)^\trs$, and $x,y$ are the absolute Cartesian coordinates in the
plane $\br^2$. The objective is to drive the robot to the point
$\bldr^0(t)$ where $D(t,\bldr)$ attains its maximum over $\bldr$ and
then to keep it in a vicinity of $\bldr^0(t)$, thus displaying the
approximate location of $\bldr^0(t)$. The on-board control system has
access only to the field value $d(t):= D[t,\bldr(t)]$ at the robot
current location $\bldr(t) = [x(t),y(t)]^\trs$.  No data about the
derivatives of $D$ are available; in particular, the robot is aware of
neither the partial derivatives of $D(\cdot)$ nor the time-derivative
$\dot{d}$ of the measurement $d$.
\par
The kinematic model of the robot is as follows:

\begin{equation}
  \label{ext:1}
  \dot{\bldr} = \vec{v}, \qquad \bldr (0) = \bldr_{\text{in}}\qquad \|\vec{v}\| \leq v,
\end{equation}

The problem is to design a
controller that drives the robot into the desired vicinity
$V_\star(t)$
of the time-varying maximizer $\boldsymbol{r}^0(t)$ in a finite time
$t_0$ and then $t\geq t_0$ keeps the robot within $V_\star(t)$.
\par
In this chapter, the following control algorithm is examined:

\begin{equation}
  \label{ext:c.a}
  \vec{v}(t) = v \vec{e}\Big\{\mu\big[d(t) - v_\ast t \big] + \theta_0 \Big\}, \;
  \text{where} \;
  \vec{e}(\theta) := \left( \begin{smallmatrix} \cos \theta \\ \sin \theta \end{smallmatrix}\right)
\end{equation}

Here the variables $v_\ast , \mu >0, \theta_0 \in \br$ are tunable
parameters. Implementation of this algorithm requires access to an
absolute direction, which may be obtained for example from a compass.

\begin{figure}
  \subfigure[]{\includegraphics[width=7cm]{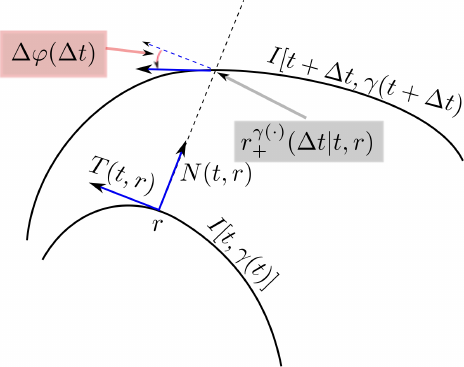}}
  \hspace{20pt}
  \subfigure[]{\includegraphics[width=7cm]{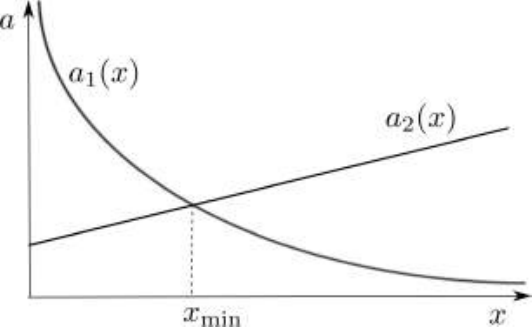}}
  \caption{(a) Two close isolines; (b) The graphs of $a_i(x)$.}
  \label{ext:fig.isl}
\end{figure}

To discuss this control law, the following notations and
quantities characterizing the moving field $D(\cdot)$ are needed:

\begin{itemize}
\item $\spr{\cdot}{\cdot}$ -- the standard inner product in the
  plane;
\item $\Phi_\beta = \left( \begin{smallmatrix} \cos \beta & - \sin
      \beta \\ \sin \beta & \cos \beta\end{smallmatrix}\right)$ --
  the matrix of counter-clockwise rotation through angle $\beta$;
\item $\nabla = \left( \begin{smallmatrix} \frac{\partial}{\partial x}
      \\ \frac{\partial}{\partial y} \end{smallmatrix}\right)$ -- the
  spatial gradient;
\item $D^{\prime\prime}$ -- the spatial Hessian, i.e., the matrix of
  the second derivatives with respect to $x$ and $y$;
\item $I(t,\gamma) := \{\bldr: D(t,\bldr) = \gamma \}$ -- the spatial
  isoline, i.e., the level curve of $D(t,\cdot)$ with the field level
  $\gamma$;
\item $[T,N]= [T(t,\bldr),N(t,\bldr)]$ -- the Frenet frame of the
  spatial isoline $I[t,\gamma]$ with $\gamma:= D(t,\bldr)$ at the
  point $\bldr$, i.e., $ N(t,\bldr) = \frac{\nabla
    D(t,\bldr)}{\|\nabla D (t,\bldr)\|} $ and the unit tangent vector
  $T(t,\bldr)$ is oriented so that when traveling on $I[t,\gamma]$ one
  has the domain of grater values $G_t^\gamma:=\{ \bldr^\prime :
  D(t,\bldr^\prime) > \gamma \}$ to the left;
\item $\varkappa$ -- the signed curvature of the spatial
  isoline;\footnote{This is positive and negative on the convexities
    and concavities, respectively, of the boundary of $G^\gamma_t$.}
\item $\bldr_+^{\gamma(\cdot)}(\Delta t| t,\bldr)$ -- the nearest (to
  $\bldr$) point of intersection between the ordinate axis of the
  Frenet frame and the displaced isoline $I[t+\Delta t,\gamma(t+\Delta
  t)]$, where the smooth function $\gamma(\cdot)$ is such that
  $\gamma(t)= D(t,\bldr)$; see Fig.~\ref{ext:fig.isl}(a);
\item $p^{\gamma(\cdot)}(\Delta t| t,\bldr)$ -- the ordinate of
  $\bldr_+^{\gamma(\cdot)}(\Delta t| t,\bldr)$;
\item $\lambda^{\gamma(\cdot)}(t,\bldr)$ -- the front velocity of the
  spatial isoline: $ \lambda^{\gamma(\cdot)}(t,\bldr):=\lim_{\Delta
    t\to 0}$ $\frac{p^{\gamma(\cdot)}(\Delta t| t,\bldr)}{\Delta t }; $
  if $\gamma(\cdot) \equiv \text{const}$ ($=D(t,\bldr)$), the upper
  index $^{\gamma(\cdot)}$ is dropped in the last three notations;
\item $\alpha(t,\bldr)$ -- the front acceleration of the spatial
  isoline $I[t,\gamma], \gamma:= D(t,\bldr)$:
  \begin{equation}
    \label{ext:alpha_def}
    \alpha(t,r):= \lim_{\Delta t\to 0}\frac{\lambda[t+\Delta t, \bldr_+(\Delta t)] - \lambda[t,\bldr]}{\Delta t }, \qquad \bldr_+(\Delta t) := \bldr_+(\Delta t| t,\bldr);
  \end{equation}
\item $\Delta \varphi(\Delta t|t,\bldr)$ -- the angular displacement
  of $T[t+\Delta t, \bldr_+(\Delta t) ]$ with respect to $T[t,\bldr]$;
  see Fig.~\ref{ext:fig.isl}(a);
\item $\omega(t,\bldr)$ -- the angular velocity of rotation of the
  spatial isoline $I[t,\gamma], \gamma:= D(t,\bldr)$, i.e., $
  \omega(t,\bldr):=\lim_{\Delta t\to 0}\frac{\Delta \varphi(\Delta
    t|t,\bldr)}{\Delta t } $;
\item $\rho(t,\bldr)$ -- the density of isolines at time $t$ at point
  $\bldr$: $ \rho(t,\bldr):=\lim_{\Delta \gamma \to 0}\frac{\Delta
    \gamma}{q(\Delta \gamma| t,\bldr)}, $ where $q(\Delta \gamma|
  t,\bldr)$ is the ordinate of the nearest (to $\bldr$) point of
  intersection between the ordinate axis of the Frenet frame
  $[T(t,\bldr),N(t,\bldr)]$ and the close isoline $I(t|
  t,\gamma+\Delta \gamma), \gamma:= D(t,\bldr)$; \footnote{This
    density characterizes the "number" of isolines within the unit
    distance from the basic one $I(t,\gamma)$, where the "number" is
    evaluated by the discrepancy in the values of $D(\cdot)$ observed
    on these isolines.}
\item $ v_\rho(t,\bldr)$ -- the evolutional (proportional) growth
  rate of the above density at time $t$ at point $\bldr$:
  \begin{equation}
    \label{ext:vrho.def}
    v_\rho(t,\bldr):= \frac{1}{\rho(t,\bldr)}\lim_{\Delta t\to 0}\frac{\rho[t+\Delta t, \bldr_+(\Delta t)]- \rho(t,\bldr)}{\Delta t};
  \end{equation}
\item $ \tau_\rho(t,\bldr)$ -- the tangential (proportional) growth
  rate of the isolines density at time $t$ at point $\bldr$:
  \begin{equation}
    \label{ext:beta_def}
    \tau_\rho (t,r):= \frac{1}{\rho(t,\bldr)} \lim_{\Delta s \to 0}\frac{\rho(t, r+ T \Delta s ) -  \rho(t,r)}{\Delta s} ;
  \end{equation}
\item $ n_\rho(t,\bldr)$ -- the normal (proportional) growth rate of
  the isolines density at time $t$ at point $\bldr$:
  \begin{equation}
    \label{ext:njrm_def}
    n_\rho (t,r):= \frac{1}{\rho(t,\bldr)}\lim_{\Delta s \to 0}\frac{ \rho(t, r+ N \Delta s ) -  \rho(t,r)}{\Delta s} .
  \end{equation}
\item $\omega_\nabla(t,\bldr)$ -- the angular velocity of the
  gradient $\nabla D$ rotation at time $t$ at point $\bldr$.
\end{itemize}

 \section{Main Assumptions}
\label{ext:sec.ass}

Local extrema are
allowed but only in a "vicinity" $Z_{\text{near}}$ of the maximizer
and at the outskirts $Z_{\text{far}}$, where the field is not assumed
even smooth. In the intermediary $Z_{\text{reg}}$, the field is smooth
and has no critical points and thus local extrema. To make the problem
tractable, $Z_{\text{near}}$ is assumed to lie within the desired
vicinity $V_\star$ of the maximizer, whereas the initial location is
in $Z_{\text{reg}} \cup Z_{\text{near}}$. The controller should keep
the robot in $Z_{\text{reg}}$ until reaching $V_\star$, thus avoiding
detrimental effects of local extrema. To reduce the amount of
technicalities, it is also assumed that the zones $Z_{\text{near}},
Z_{\text{far}}$, $Z_{\text{reg}}$, and $V_\star$ are separated by
isolines. This requirement can usually be met by properly reducing
$Z_{\text{reg}}$ and $V_\star$, if necessary. To encompass theoretical
distributions like $D(\bldr) = c/\|\bldr-\bldr^0(t)\|$ or $D(\bldr) =
-c \ln \|\bldr-\bldr^0(t)\|$ and their sums, the field is allowed to
be undefined at finitely many moving exceptional points.  This is
addressed in the following:

\begin{Assumption}
  \label{ext:ass.0}
  There exist smooth $\gamma_-(t) < \gamma_\star(t) < \gamma_+(t) \in
  \br$ and continuous functions $\bldr_1(t), \ldots, \bldr_k(t) \in
  \br^2$ of time $t$ such that the following statements hold:
  \begin{enumerate}[i)]
  \item On $Z_{\text{reg}}:= \{(t,\bldr) : \gamma_-(t) \leq D(t,\bldr)
    \leq \gamma_+(t) \}$, the distribution $D(\cdot)$ is identical to
    a $C^2$-smooth function defined on a larger and open set, and
    $\nabla D(t,\bldr) \neq 0$;
  \item At any time $t$, the spatial isoline $I[t, \gamma_-(t)]$ is a
    Jordan curve that encircles $I[t,\gamma_+(t)]$;
  \item At any time $t$, the points $\bldr_1(t), \ldots, \bldr_k(t)$
    lie inside $I[t,\gamma_+(t)]$;
  \item In the domain $Z_{\text{near}}:= \{(t,\bldr): \bldr \,\text{is
      inside}\, I[t,\gamma_+(t)] \, \text{and}\, \bldr \neq \bldr_j(t)
    \; \forall j\}$, the field $D(\cdot)$ takes values greater than
    $\gamma_+(t)$, converges to a finite or infinite limit $L_j(t)$ as
    $\bldr \to \bldr_j(t)$ for any $t$ and $j$, and is continuous in
    both $Z_{\text{near}}$ and its outer boundary $\{(t,\bldr): \bldr
    \in I[t,\gamma_+(t)] \}$;
  \item In the domain $Z_{\text{far}}:=\{(t,\bldr): \bldr \,\text{is
      outside}\, I[t,\gamma_-(t)] \}$, the field $D(\cdot)$ is
    everywhere defined, takes values lesser than $\gamma_-(t)$, and is
    continuous in both $Z_{\text{far}}$ and its boundary $\{(t,\bldr):
    \bldr \in I[t,\gamma_-(t)] \}$;
  \item The desired vicinity of the maximizer has the form $V_\star(t)
    = \{\bldr: D[t,\bldr(t)] \geq \gamma_\star(t)\}$;
  \item \label{ext:vii} For $t=0$, the initial location
    $\bldr_{\text{\rm in}}$ lies in the domain of $D(\cdot)$ and
    inside $I[0,\gamma_-(0)]$ (i.e., $D[0,\bldr_{\text{\rm in}}] >
    \gamma_-(0)$).
  \end{enumerate}
\end{Assumption}

The next assumption is typically fulfilled in real world, where
physical quantities take bounded values:

\begin{Assumption}
  \label{ext:ass.last}
  There exists constants $b_\omega^\nabla$, $b_{\aleph}$ for $\aleph =
  \rho, \lambda, \omega, \varkappa, v, \alpha, n, \tau$ and
  $\gamma^0_+, \ov{\gamma}$ such that
  \begin{multline}
    \label{ext:estim}
    \begin{array}{l}
      |\rho| \leq b_\rho, \;
      |\lambda| \leq b_\lambda, \; |\omega| \leq b_\omega,
      \; |\omega_\nabla| \leq b_\omega^\nabla,
      \\
      |\varkappa| \leq b_\varkappa, \; |v_\rho| \leq b_v, \; |\alpha| \leq b_\alpha, \; |\tau_\rho| \leq b_\tau, \; |n_\rho| \leq b_n
    \end{array}
    \; \forall (t,\bldr) \in Z_{\text{reg}}; \\ \gamma_+(t) \leq \gamma^0_+, |\dot{\gamma}_i(t)| \leq  \ov{\gamma} \quad \forall t, i=\pm,\star.
  \end{multline}
\end{Assumption}

The only assumption about the robot capacity with respect to the field
is that the mobility of the former exceeds that of the latter: in the
main operational zone $Z_{\text{reg}}$, the maximal speed of the robot
is greater than the front speed of the concerned isoline of the field:
$v > |\lambda(t,\bldr)|$. Moreover, if the level $\gamma_-(\cdot)$ or
$\gamma_\star(\cdot)$ is not constant, the robot is capable to remain
inside the moving isoline $I[t,\gamma_-(t)]$ and inside
$I[t,\gamma_\star(t)]$. 

All afore-mentioned strict inequalities are protected from degradation
as time progresses. To quanify this the notation {\it $f
  \ntriangleright g$ in $Z$} is used to express that $\exists \ve >0:
f(t,\bldr) \geq g(t,\bldr) + \ve \; \forall (t,\bldr) \in Z$; the
relation $\ntriangleleft$ is defined likewise.

\begin{Assumption}
  \label{ext:ass.speed}
  The following inequalities hold:
  \begin{equation}
    \label{ext:esstt.grad}
    \rho(t,\bldr) \ntriangleright
 0 , \quad v \ntriangleright |\lambda(t,\bldr)| + \rho^{-1} \ov{\gamma}, \; \quad \text{\rm in}\; Z_{\text{reg}}; \qquad \gamma_-(t) \ntriangleleft \gamma_\star(t) \ntriangleleft \gamma_+(t) \quad \text{\rm in}\; [0,\infty).
  \end{equation}
\end{Assumption}

 \section{Summary of Main Results}
\label{ext:sec.mr}

The first theorem shows that the control
objective can always be achieved by means of the control law
Eq.\eqref{ext:c.a}:

\begin{Theorem}
  \label{ext:th.m0} Suppose that Assumptions~{\rm
    \ref{ext:ass.0}--\ref{ext:ass.speed}} hold. Then there exist
  parameters $v_\ast, \mu, \theta_0$ of the controller Eq.\eqref{ext:c.a}
  such that the following claim holds:
  \begin{enumerate}[{\rm (i)}]
  \item The controller Eq.\eqref{ext:c.a} brings the robot to the desired
    vicinity of a maximizer in a finite time $t_0$ and keeps it there
    afterwards: $\bldr(t) \in V_\star(t)$ for $t\geq t_0$.
  \end{enumerate}
  Moreover, for any compact domain $D \subset \interior \{ \bldr:
  (0,\bldr) \in Z_{\text{reg}} \}$, there exist common values of the
  parameters for which {\rm (i)} holds whenever the initial location
  $\bldr_{\text{\rm in}} \in D$.
\end{Theorem}

 The remainder of the section is devoted
to discussion of controller parameters tuning.
The parameters $\theta_0,v_\ast$ in Eq.\eqref{ext:c.a} are chosen prior to
$\mu$. Whereas $\theta_0$ is arbitrary, $v_\ast$ is chosen so that:

\begin{equation}
  \label{ext:v.astchoice}
  v \ntriangleright |\lambda| + \rho^{-1} v_\ast  \quad  \text{in} \quad Z_{\text{reg}}, \qquad v_\ast > \ov{\gamma},
\end{equation}

This is possible thanks to Eq.\eqref{ext:esstt.grad}. The choice of $\mu$ is
prefaced by picking constants $\Delta_\nabla >0$ and $\Delta_\gamma
>0$ such that

\begin{multline}
  \label{ext:rho.lower}
  \rho(t,\bldr) \geq \Delta_\nabla  \quad \forall (t,\bldr) \in Z_{\text{reg}}, \quad \gamma_-(t)+\Delta_\gamma \leq \gamma_\star(t) \leq \gamma_+(t) - \Delta_\gamma, \\ D[0,\bldr_{\text{\rm in}}] \geq \gamma_-(0) + \Delta_\gamma/2,
\end{multline}

This is possible by the same argument. The parameter $\mu$ is chosen
so that for some $k=1,2,\ldots$,

\begin{subequations}
\begin{eqnarray}
  \label{ext:mu0.choice}
  \mu(v_\ast - \ov{\gamma}) \ntriangleright -2 \omega - \varkappa v_T
  + 2 \frac{\tau_\rho \ov{\gamma}}{\rho}  + 2 \frac{v_\rho
    \ov{\gamma}}{v_T \rho} - \frac{\alpha}{v_T} +  \frac{n_\rho
    \ov{\gamma}^2}{v_T \rho^2},\\ \nonumber \text{where} \; v_T := \pm \sqrt{v^2 -\left[ \lambda  + \rho^{-1} \ov{\gamma} \right]^2}, \quad \text{in} \quad Z_{\text{reg}};
  \\
  \label{ext:mu1.choice}
  \mu \ntriangleright \frac{\omega + \tau_\rho (\lambda - v)}{\rho ( v - \lambda) - v_\ast}
  \quad \text{in} \quad Z_{\text{reg}};
  \\
  \label{ext:mu2.choice}
  \mu (v_\ast - \ov{\gamma}) >
  \left\{
    \begin{array}{l}
      a_1(k):=2 b^\nabla_\omega \left[ 1+ \frac{1}{k}\right] + 2 v \sqrt{b_\varkappa^2+b_\tau^2} \left[ 2+ \frac{1}{k}\right],
      \\
      a_2(k):= 2b_\rho \frac{(2k+1)v+ 2 \pi (k+1) \left( b_\lambda + \Delta_\nabla^{-1} \ov{\gamma}\right)}{\Delta_\gamma}
    \end{array}
  \right. .
\end{eqnarray}
\end{subequations}

The following theorem ensures the correct behaviour is exhibited by the robot:

\begin{Theorem}
  \label{ext:th.m1} Suppose that Assumptions~{\rm
    \ref{ext:ass.0}--\ref{ext:ass.speed}} hold and the controller
  parameters satisfy Eq.\eqref{ext:v.astchoice} and
  Eq.\eqref{ext:mu0.choice}--Eq.\eqref{ext:mu2.choice}. Then {\rm (i)} from
  Theorem~{\rm \ref{ext:th.m0}} is true.
\end{Theorem}

 \section{Simulations}
\label{ext:sec.simtest}

Simulations were carried out with the point-wise robot Eq.\eqref{ext:1}
driven by the control law Eq.\eqref{ext:c.a}.  The numerical values of the
parameters used for simulations are shown in Table~\ref{ext:fig:param}
(where $u_d$ is the unit of measurement of $d=D(t,\bldr)$ and the
controller parameters were chosen according to recommendations from
Theorem~\ref{ext:th.m1}). The control was updated with a sampling time
of $0.1 s$.

\begin{table}[ht]
  \centering
  \begin{tabular}{| l | c |}
    \hline
    $v$ & $1 m/s$  \\
    \hline
    $v_{*}$ & $0.299 u_d/s$ \\
    \hline
  \end{tabular}
\hspace{10pt}
 \begin{tabular}{| l | c |}
    \hline
    $\mu$ & $0.8 rad/u_d$  \\
    \hline
    $\theta_0$  & $1.5 rad$\\
    \hline
  \end{tabular}
  \caption{Simulation parameters for extremum-seeking controller.}
  \label{ext:fig:param}
\end{table}

Fig.~\ref{ext:fig:lin} displays the results of tests in a linear field
with the orientation angle of $0.5$ rad and the ascension rate of $0.3
m^{-1}$; the grey intensity is proportional to the field value. The
steady state angular error may be computed to be $\approx 0.082
rad$. Fig.~\ref{ext:fig:lin} shows that the heading converges to
approximately $0.582 rad$, thus displaying a good match.

\begin{figure}[ht]
  \centering
  \subfigure[]{\includegraphics[width=0.75\columnwidth]{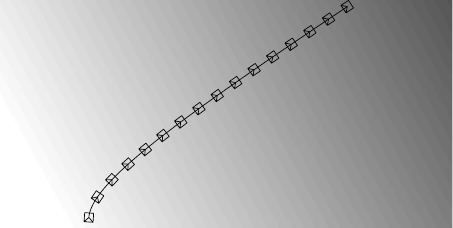}}
  \subfigure[]{\includegraphics[width=10cm]{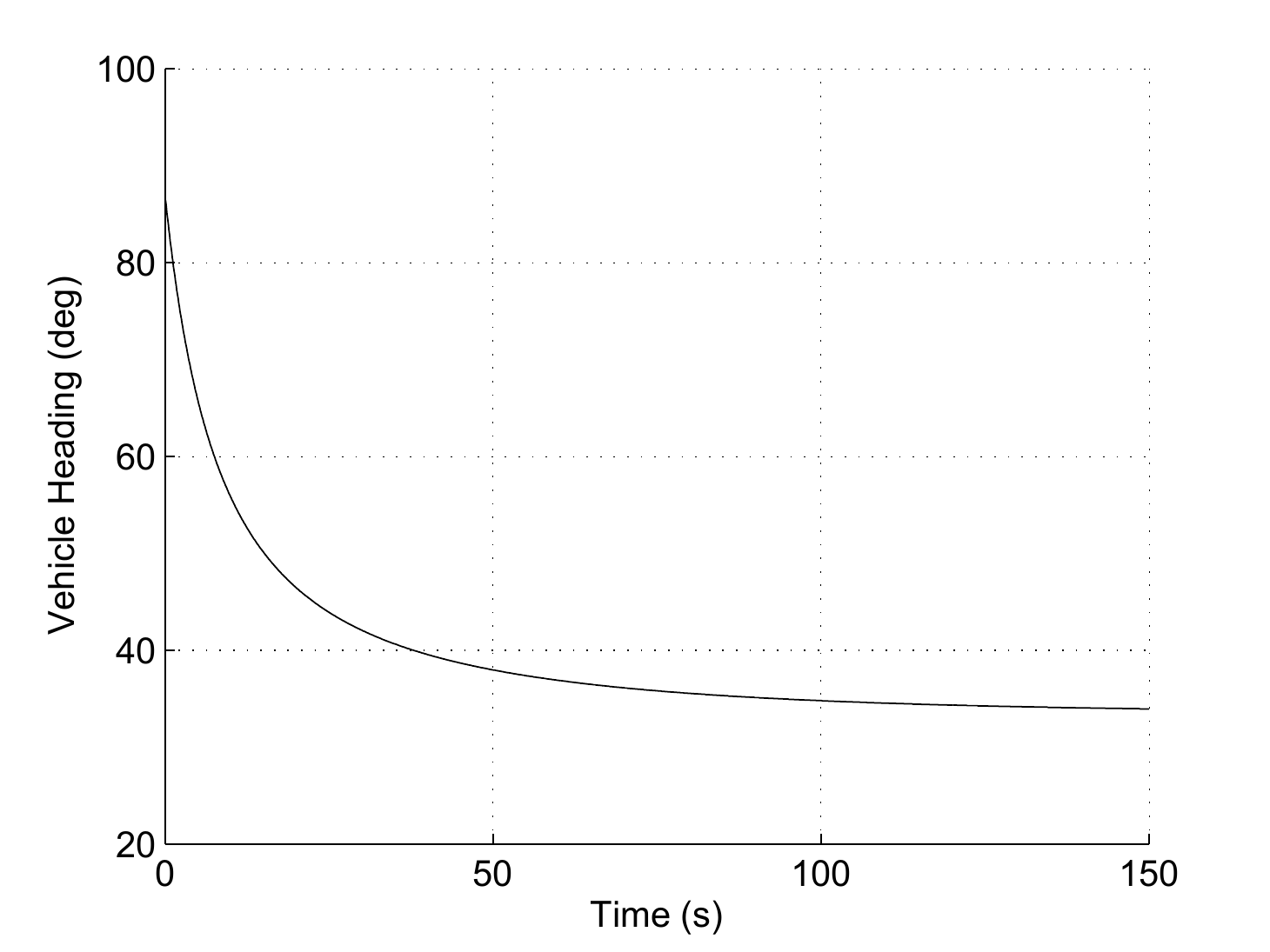}}
  \caption{Behavior in a linear field; (a) Path; (b) Robot's
    orientation.}
  \label{ext:fig:lin}
\end{figure}

Fig.~\ref{ext:fig:mov} presents the results of tests in an unsteady
field with a moving source.  The source $\bldr^0(t)$ moves to the
right at the speed of $0.3 ms^{-1}$, $D(t,\bldr)$ is a radial field
corrupted by two plain waves and is given by:

$$
D(t,\bldr) = -0.8 \cdot \lVert \bldr - \bldr^0(t)\rVert +
5\cdot\left[\sin(0.05 \cdot x) + \sin (0.05 \cdot y) \right].
$$

As can be seen, the robot converges to the source (whose path is
depicted by the solid black line) and then wheels around it in a close
proximity, thus displaying its current location and highlighting its
displacement. `Wheeling' commences when the robot achieves the desired
vicinity of the source and is unavoidable since the robot's speed
exceeds that of the source.

\begin{figure}[ht]
  \centering
  \subfigure[]{\includegraphics[width=0.75\columnwidth]{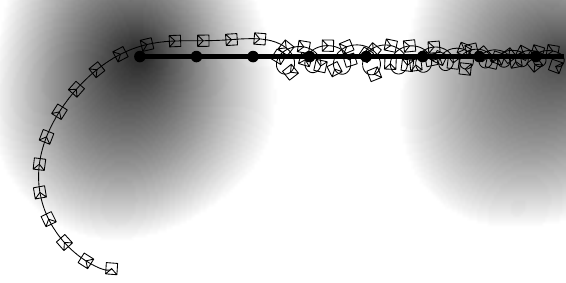}}
  \subfigure[]{\includegraphics[width=10cm]{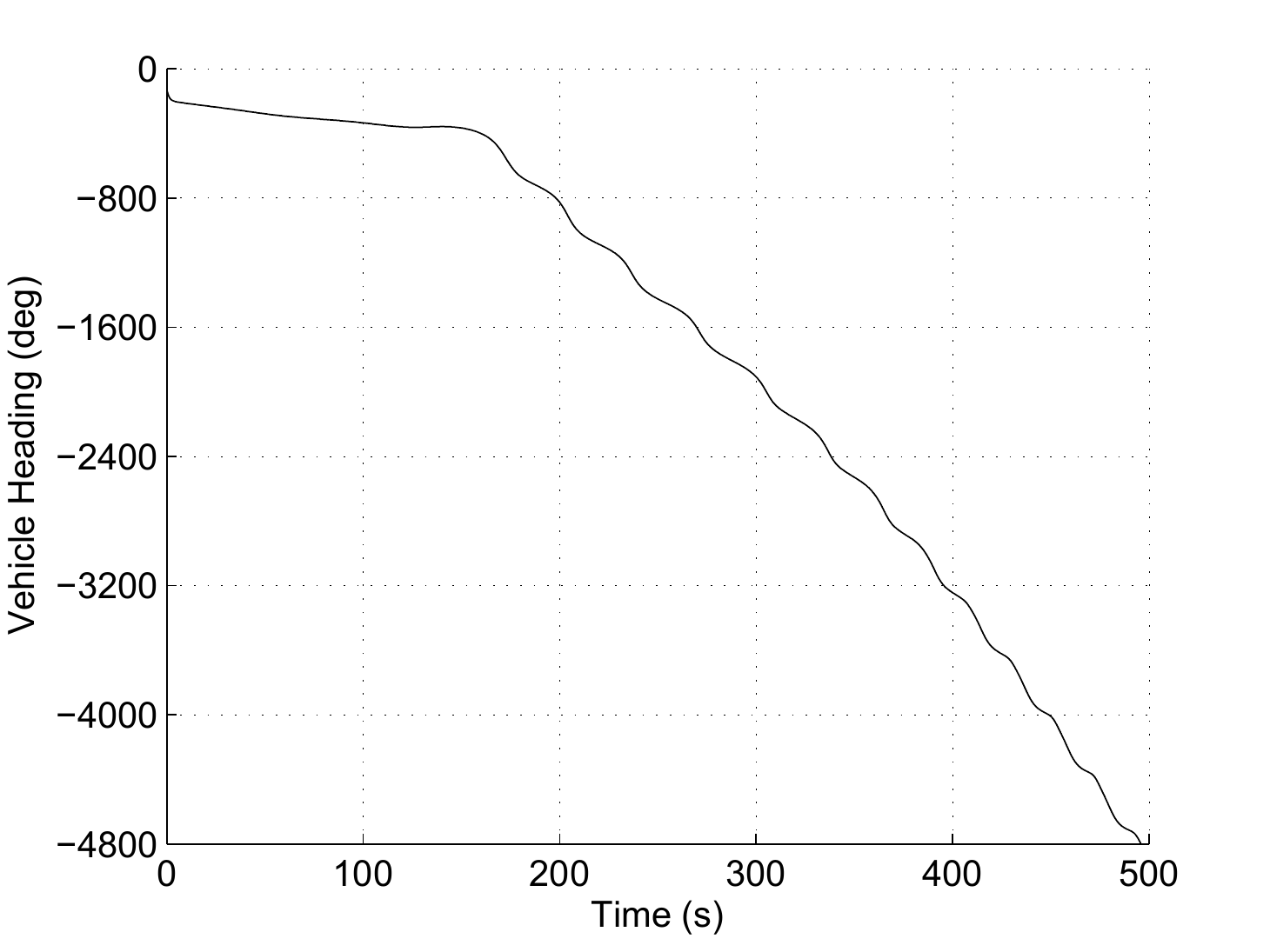}}
  \caption{Seeking a moving source; (a) Path; (b) Robot's orientation.}
  \label{ext:fig:mov}
\end{figure}

In Fig.~\ref{ext:fig:noise}, the same simulation setup was used,
except measurement noise and kinematic constraints were added.  The
robot's heading was not allowed to change faster than $0.5 rad
s^{-1}$, which in fact transforms Eq.\eqref{ext:1} into the non-holonomic
Dubins-car model since the robot's speed is constant by
Eq.\eqref{ext:c.a}. The field readings were corrupted by a random additive
noise uniformly distributed over the interval $[-2.5, 2.5]$. It may be
seen that nearly the same behavior is observed.
\begin{figure}[ht]
  \centering
  \subfigure[]{\includegraphics[width=0.75\columnwidth]{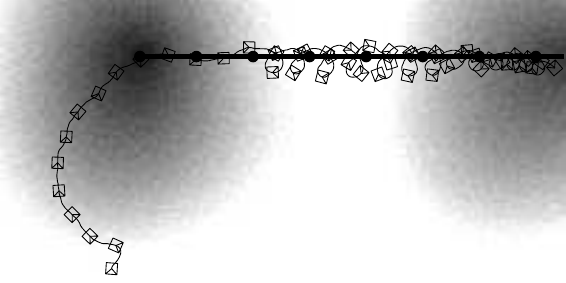}}
  \subfigure[]{\includegraphics[width=10cm]{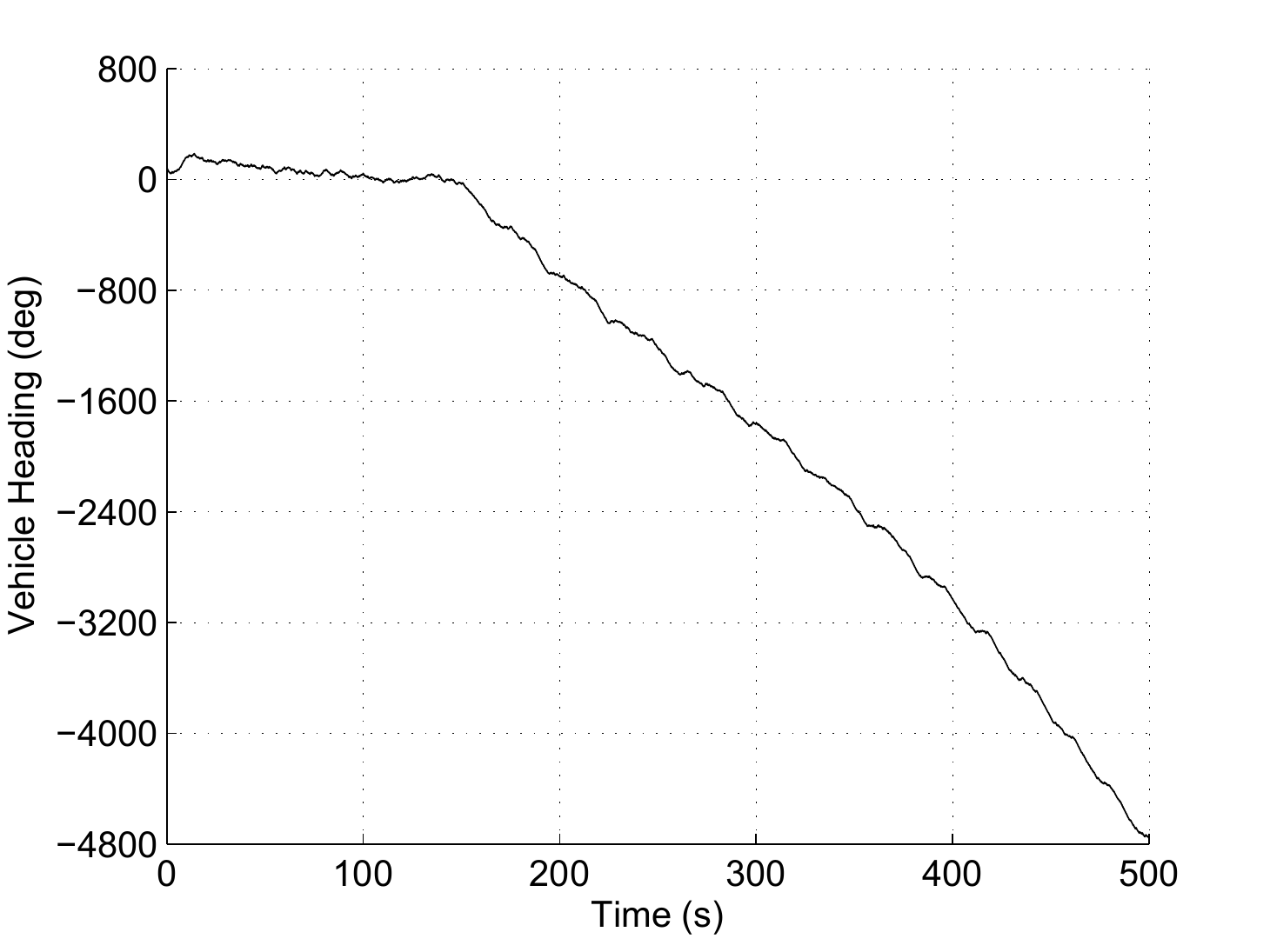}}
  \caption{Seeking a moving source under measurement noise and
    kinematic constraints; (a) Path; (b) Robot's orientation.}
  \label{ext:fig:noise}
\end{figure}

In Fig.~\ref{ext:fig:lmax}, the plain waves were enhanced to produce a
large number of local maxima and the source was stopped: the updated
field distribution is given by:
$$
D(t,\bldr) = -0.8 \cdot \lVert \bldr - \bldr^0\rVert +
10\cdot\left[\sin(1.0 \cdot x) + \sin (1.0 \cdot y) \right].
$$
Fig.~\ref{ext:fig:lmax} shows that despite the local maxima, the
vehicle stills converges to the source of the distribution. This is in
a sharp contrast to the pure gradient ascent method for which every
local maximum constitutes a trap. A presumable reason for this
property of the proposed algorithm is that it needs special tuning to
be trapped by a small enough vicinity of a maximizer. Whenever it is
tuned for a specific vicinity "size", it demonstrates the capability
to go through and escape from smaller traps. A detailed analysis of
this interesting and promising property is a subject of future
research.

\begin{figure}[ht]
  \centering
  \subfigure[]{\includegraphics[width=0.75\columnwidth]{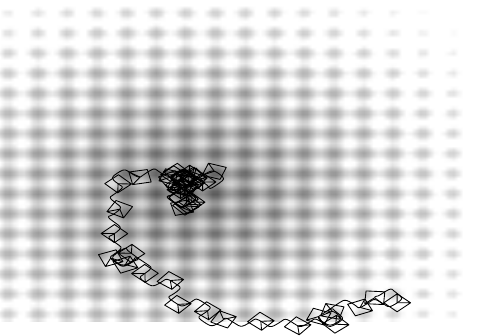}}
  \subfigure[]{\includegraphics[width=10cm]{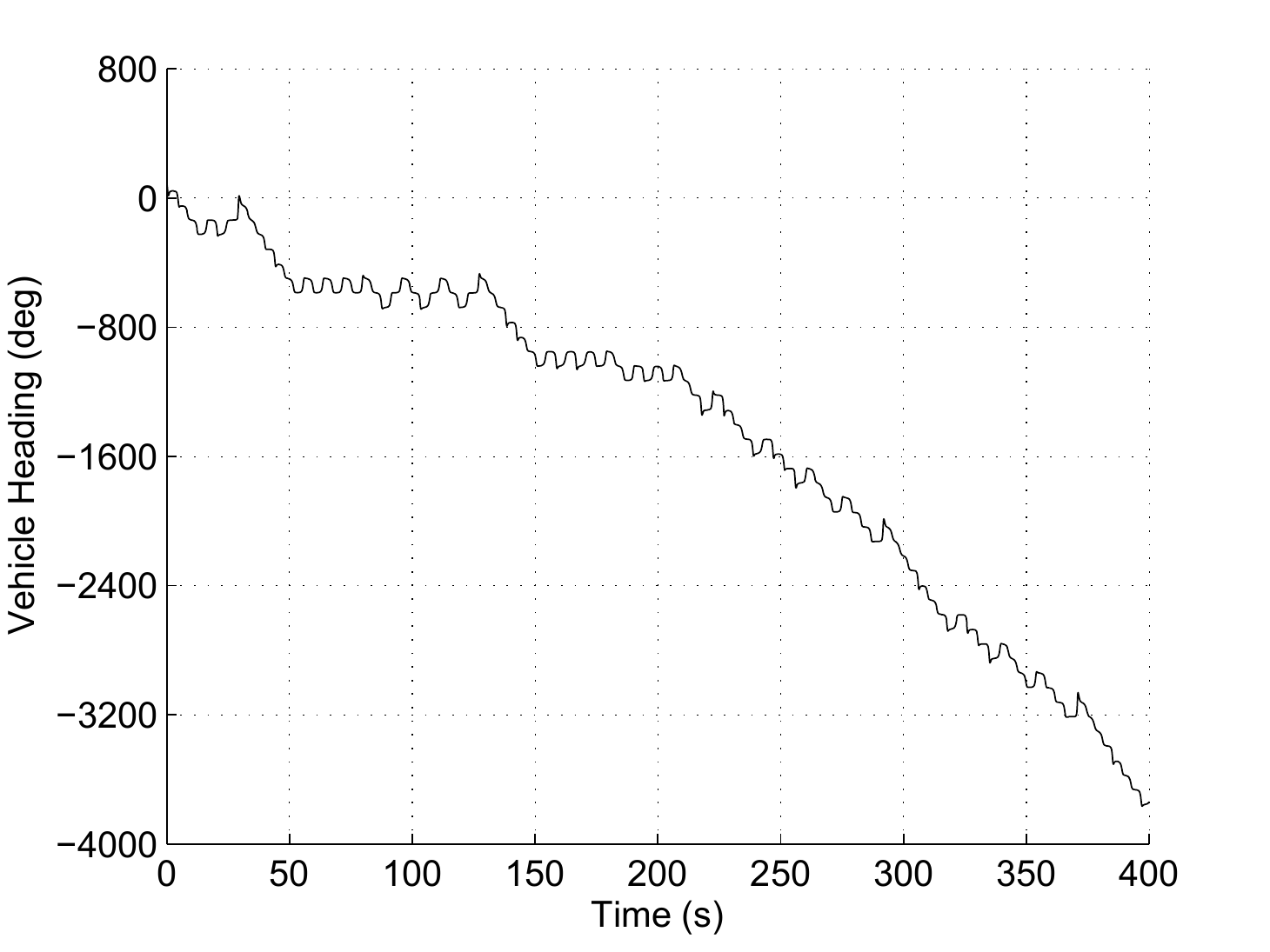}}
  \caption{Seeking a source in the presence of multiple local maxima;
    (a) Path; (b) Robot's orientation.}
  \label{ext:fig:lmax}
\end{figure}

Simulations were also carried out for a realistic model of a time
varying field caused by a constant-rate emanation of a certain
substance (such as heat or gas) from a moving point-wise source and its
subsequent diffusion in an isotropic two-dimensional medium. In many
cases, this process is described by the heat equation $\partial
D/\partial t = \rho \Delta D + \delta[\bldr - \bldr^0(t)]$. Here
$\Delta$ is the spatial Laplacian, $\rho = 16000 m^2s^{-1}$ is the
diffusion rate, $\delta$ is the spatial Dirac delta-function,
$\bldr^0(t)$ is the source location, and the emanation rate was set to
unity.  The field distribution was calculated prior to navigation
tests by the finite difference method. In doing so, the time and space
steps were $0.001 s$ and $4 m$, respectively; the results were stored
with the sampling rate $1 s$.  During the navigation test, the
distribution value was obtained through trilinear interpolation over
spatial and temporal variables.  To initialize the distribution, the
source stayed still for the first $100 s$ and only then commenced
motion. The vehicle turning rate was bounded by $0.5 rad s^{-1}$.
The results of these simulations are shown in
Fig.~\ref{ext:fig:heat}. It may be seen that the robot solves the
source seeking task.

\begin{figure}[ht]
  \centering
  \subfigure[]{\includegraphics[width=0.75\columnwidth]{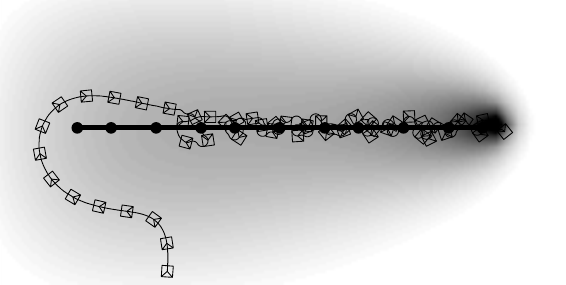}}
  \subfigure[]{\includegraphics[width=10cm]{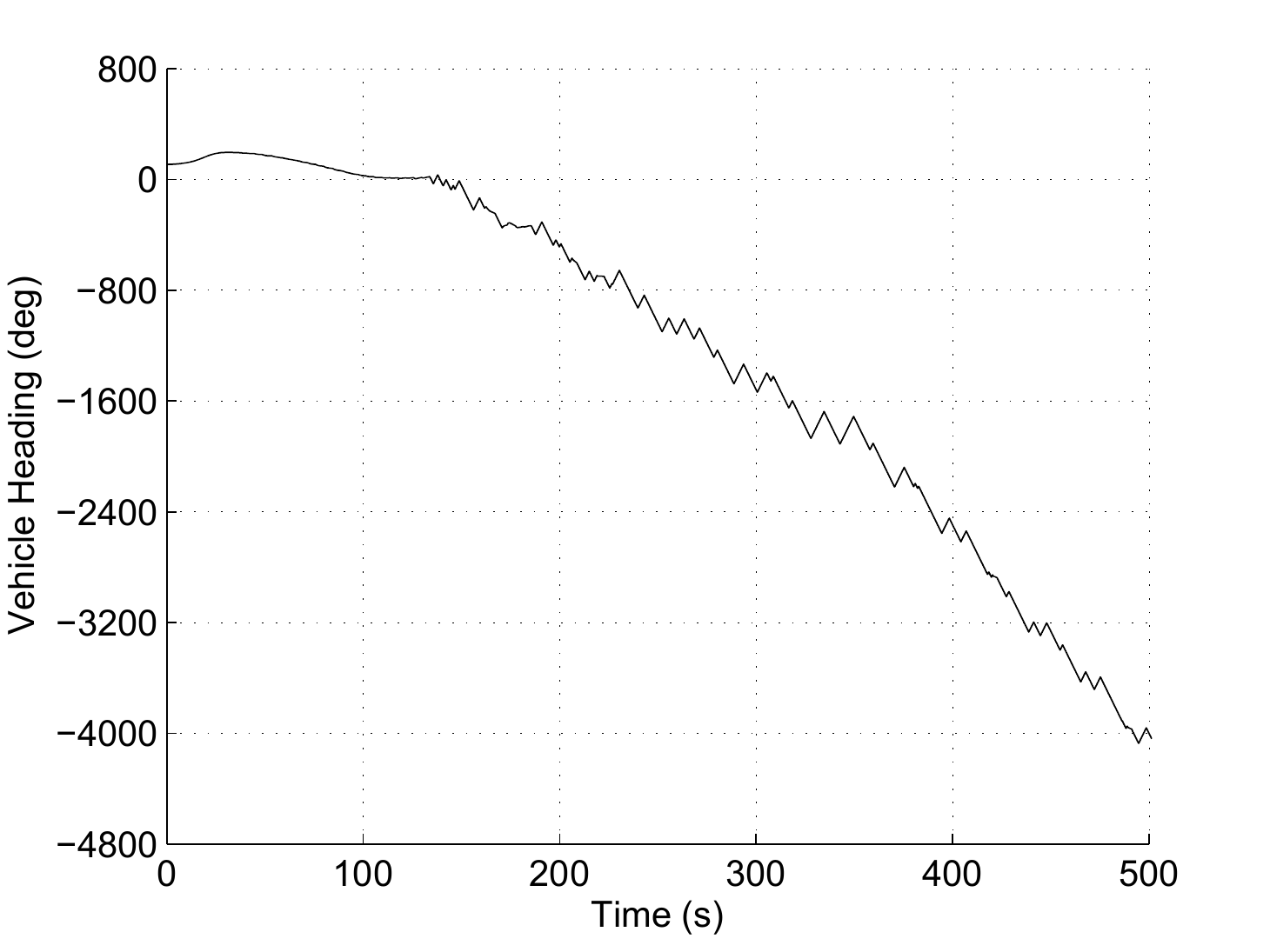}}
  \caption{Seeking a moving diffusion source; (a) Path; (b) Robot's
    orientation.}
  \label{ext:fig:heat}
\end{figure}

\begin{figure}[ht]
  \centering
  \subfigure[]{\includegraphics[width=0.75\columnwidth]{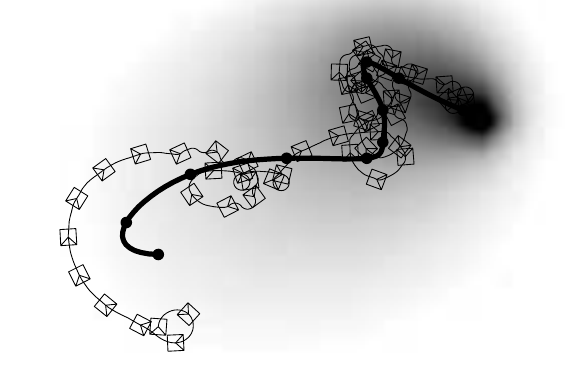}}
  \subfigure[]{\includegraphics[width=10cm]{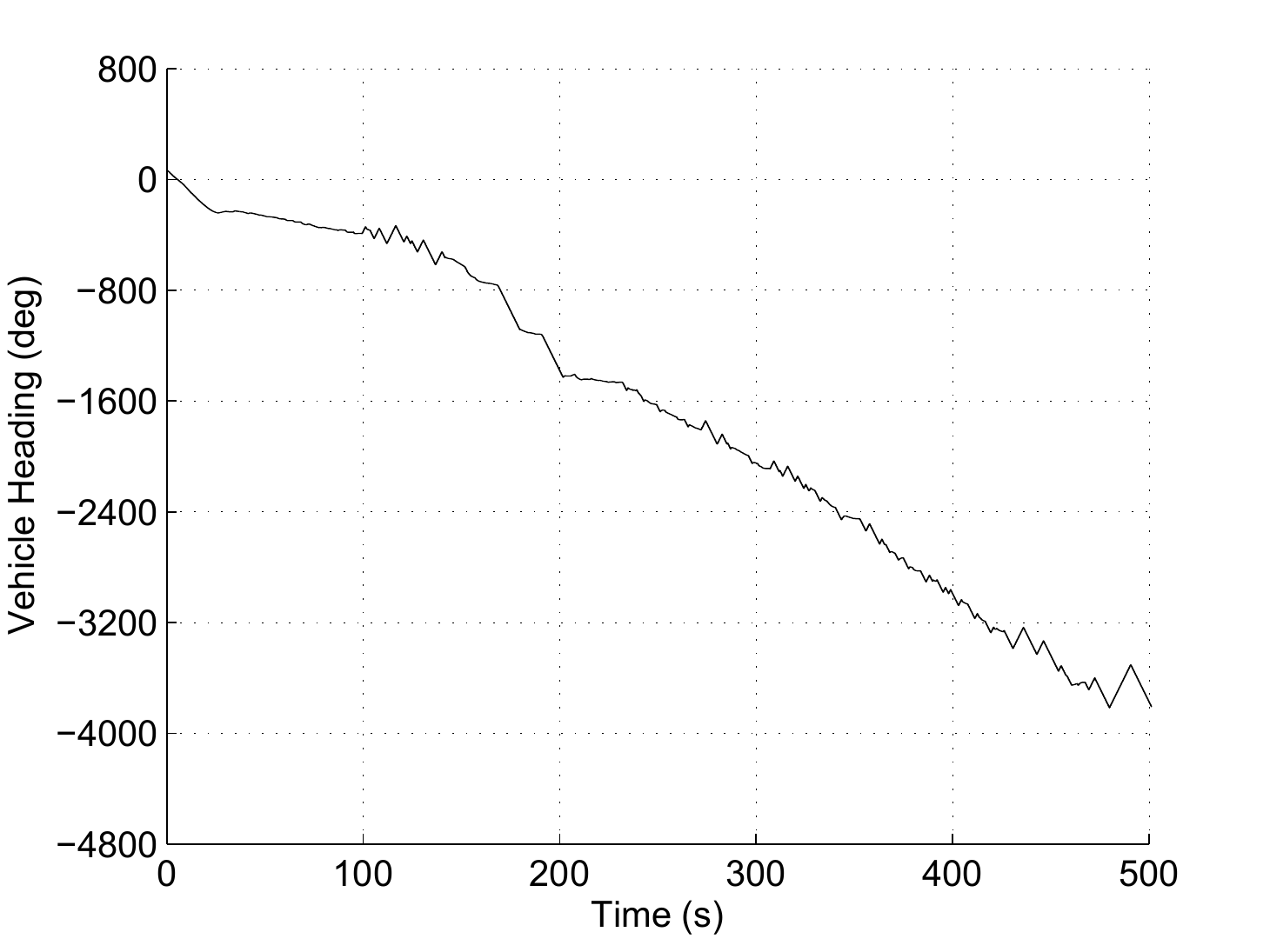}}
  \caption{Seeking an irregularly moving diffusion source; (a) Path;
    (b) Robot's orientation. }
  \label{ext:fig:heatmov}
\end{figure}

\clearpage \section{Experiments}
\label{ext:sec.exper}

Experiments were carried out with an Activ-Media Pioneer 3-DX
differential drive wheeled robot using its on-board PC and the Advanced
Robot Interface for Applications (ARIA 2.7.2), which is a C++ library
providing an interface to the robot's angular and translational
velocity set-points.

The origin of the reference frame was co-located with the center of
the robot in its initial position, its ordinate axis is directed
towards the viewer in Figs.~\ref{ext:fig.test1}, \ref{ext:fig.test4},
and \ref{ext:fig.test7}. Three scenarios were tested:
\begin{itemize}
\item A steady virtual point-wise source was located at the point with
  coordinates $(0, 3.5) m$ in Fig.~\ref{ext:fig.test1}. The experiment
  was run for $175 s$.
\item A virtual point-wise source moved from the point with the
  coordinates $(0, 3.5) m$ with a constant translational velocity $(0,
  -0.02) ms^{-1}$ outwards the viewer in Fig.~\ref{ext:fig.test4}. The
  experiment was run for $175 s$ so that the final position of the
  source coincided with the origin of the reference frame.
\item A virtual point-wise source moved at the constant speed $0.02
  ms^{-1}$ from the same point along the piece-wise linear path
  through the points $(0, 3.5) m$, $(0.6, 2.5) m$, $(-0.6, 1.0) m$ and
  $(0, 0) m$ in Fig.~\ref{ext:fig.test7}. The experiment was run for
  $213 s$.
\end{itemize}
The path of the source is displayed in Fig.~\ref{ext:fig.test4} by a
long black tape\footnote{A short perpendicular segment was added for
  calibration purposes and is not a part of the path.}, and in
Fig.~\ref{ext:fig.test7} by a long gray tape. The examined field was
minis the distance to the source, which was accessed via odometry,
whereas the source motion was virtual and emulated by computer.  The
orientation of the robot $\theta_{act}(t)$ was also captured via
odometry.
\par
The parameters used in the experiments are shown in
Table~\ref{ext:fig:paramexp}.

\begin{table}[ht]
  \centering
  \begin{tabular}{| l | c |}
    \hline
    $v$ & $0.15 m/s$  \\
    \hline
    $v_\ast$ & $0.12 u_d/s$ \\
    \hline
  \end{tabular}
\hspace{10pt}
  \begin{tabular}{| l | c |}
    \hline
$\mu$ & $4.0 rad/u_d$  \\
   \hline
    $\theta_0$  & $0 rad$\\
    \hline
  \end{tabular}
  \caption{Experimental parameters for extremum-seeking controller.}
  \label{ext:fig:paramexp}
\end{table}

The control law was updated at the rate of $0.1 s$.  The control law
Eq.\eqref{ext:c.a} was used as a generator of the velocity reference
signal. Tracking of the generated velocity profile was accomplished by
means of simple PD controllers.

Typical experimental results are presented on
Figs.~\ref{ext:fig.test1}--\ref{ext:fig.test9}. In these experiments,
like in the others, the robot successfully reached the source and then
tracked it until the end of the experiment. In
Figs.~\ref{ext:fig.test2}, \ref{ext:fig.test5}, and
\ref{ext:fig.test8}, the final maximum steady state distance to the
source is upper bounded by $0.5 m$. This approximately equals the
radius of the desired margin $R_\star$ employed in computation of the
controller parameters, which was chosen with regard to the turning
capacity of the robot at the selected speed.

\begin{figure}[ht]
  \centering
  \subfigure[]{\scalebox{0.25}{\includegraphics{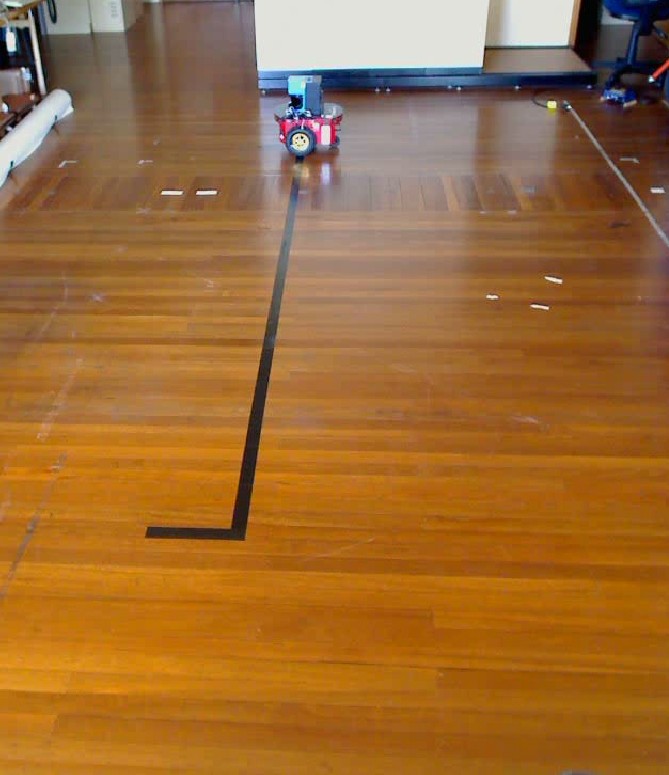}}}
  \subfigure[]{\scalebox{0.25}{\includegraphics{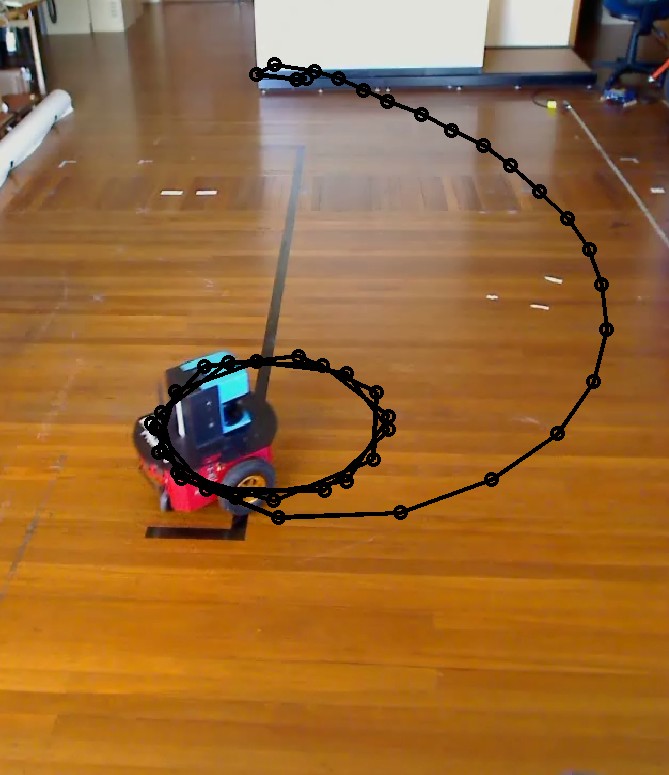}}}
  \subfigure[]{\scalebox{0.25}{\includegraphics{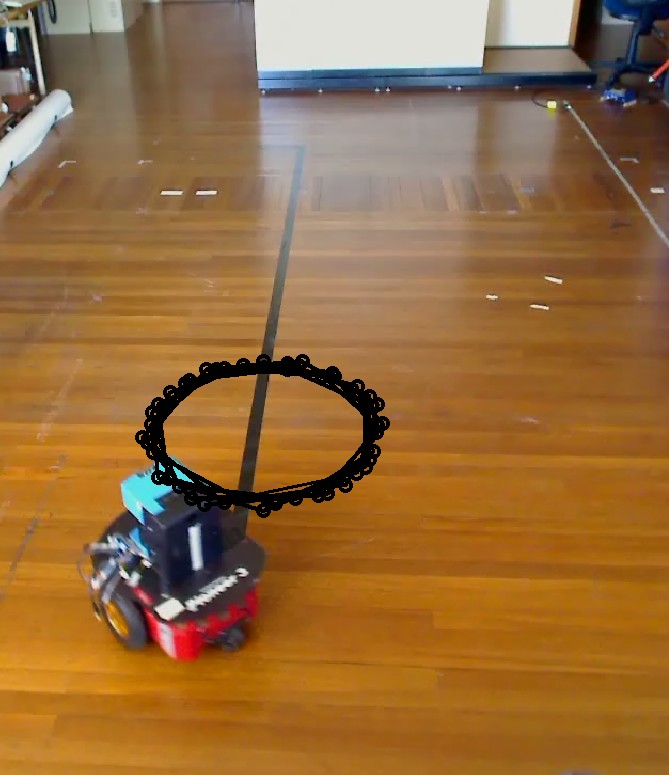}}}
  \caption{Sequence of images showing the experiment.} \label{ext:fig.test1}
\end{figure}

\begin{figure}[ht]
  \centering \includegraphics[width=10cm]{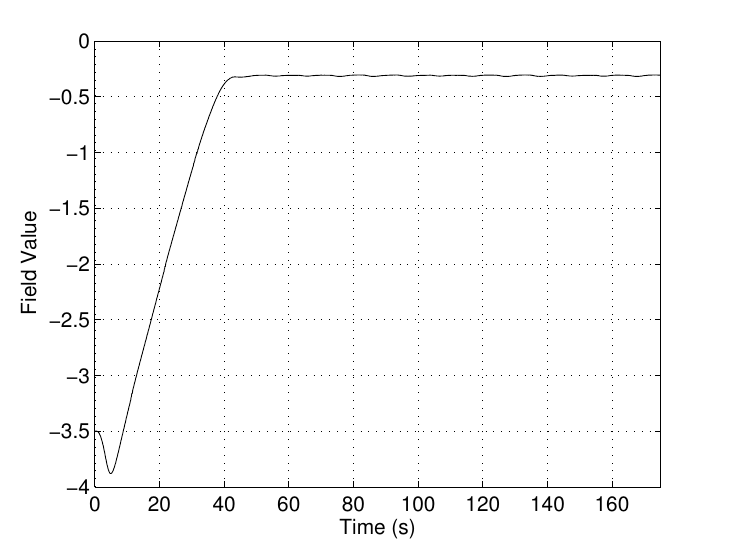}
  \caption{Evolution of the field value over the experiment. Field value is in metres.}
  \label{ext:fig.test2}
\end{figure}

\begin{figure}[ht]
  \centering
 \includegraphics[width=10cm]{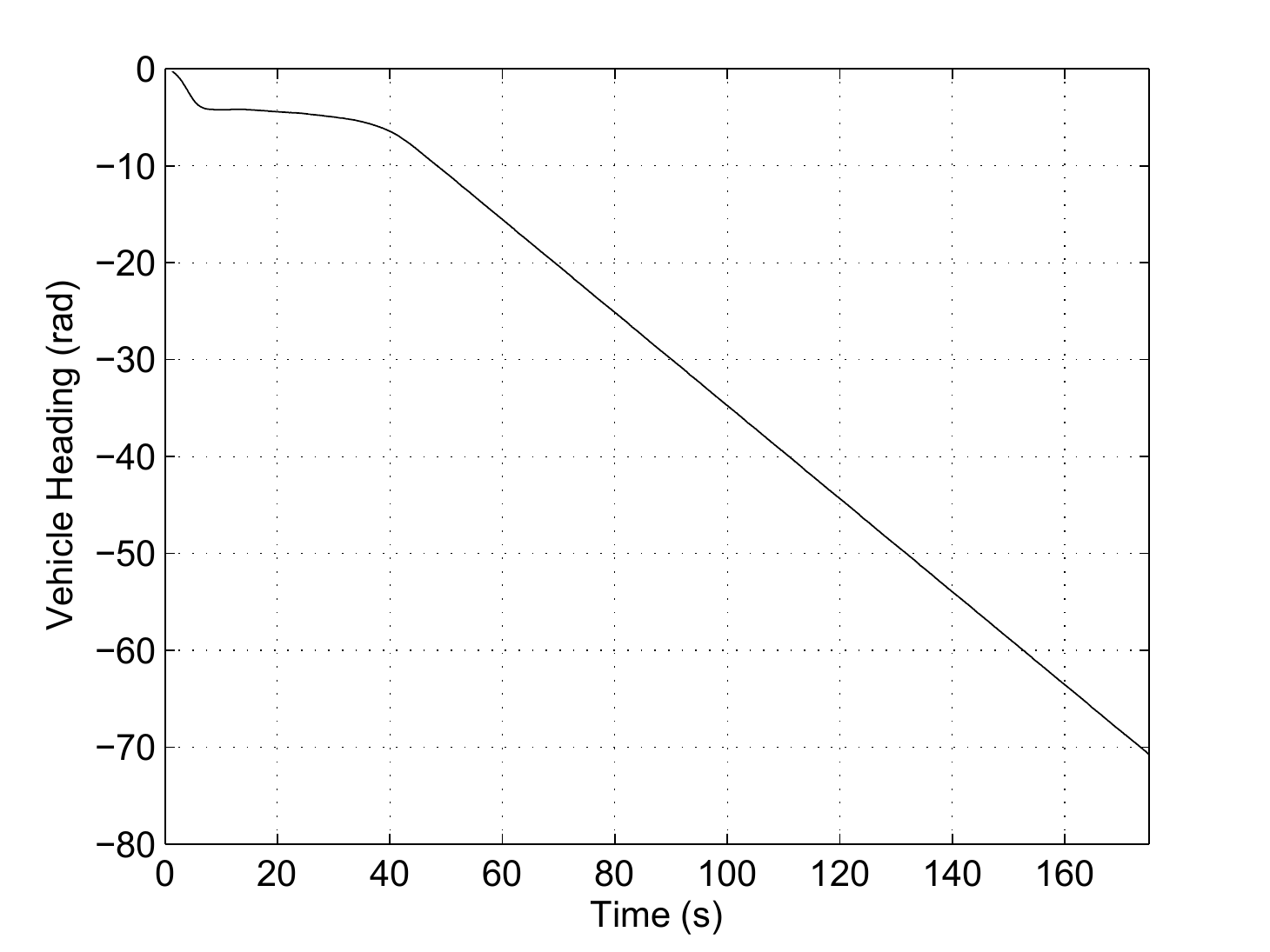}
  \caption{Evolution of the vehicle orientation over the
    experiment.}
  \label{ext:fig.test3}
\end{figure}

\begin{figure}[ht]
  \centering
  \subfigure[]{\scalebox{0.25}{\includegraphics{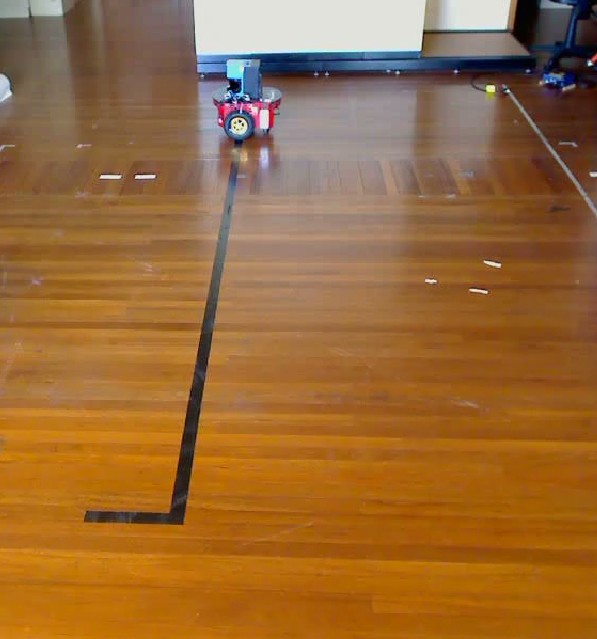}}}
  \subfigure[]{\scalebox{0.25}{\includegraphics{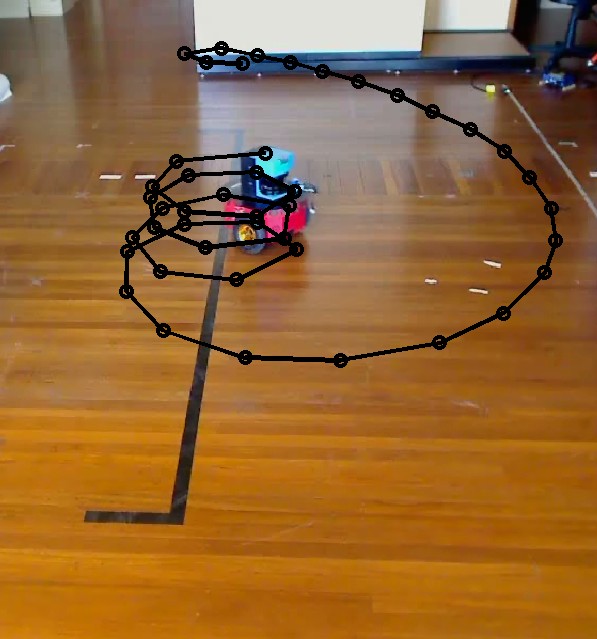}}}
  \subfigure[]{\scalebox{0.25}{\includegraphics{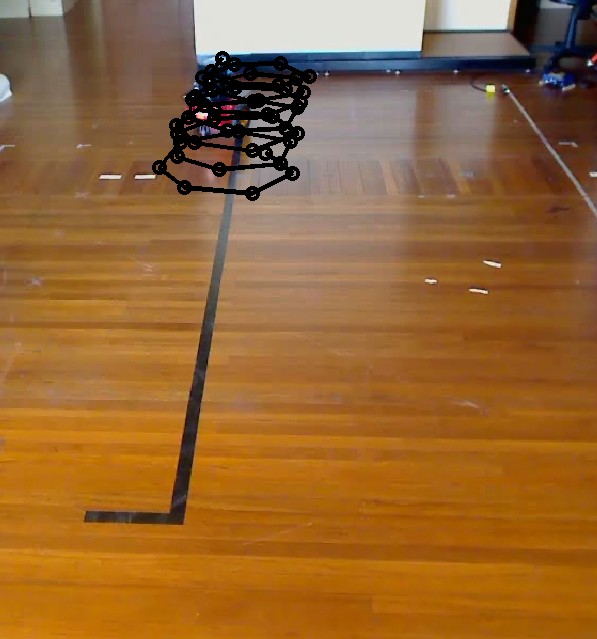}}}
  \caption{Sequence of images showing the experiment (moving field).} \label{ext:fig.test4}
\end{figure}

\begin{figure}[ht]
  \centering \includegraphics[width=10cm]{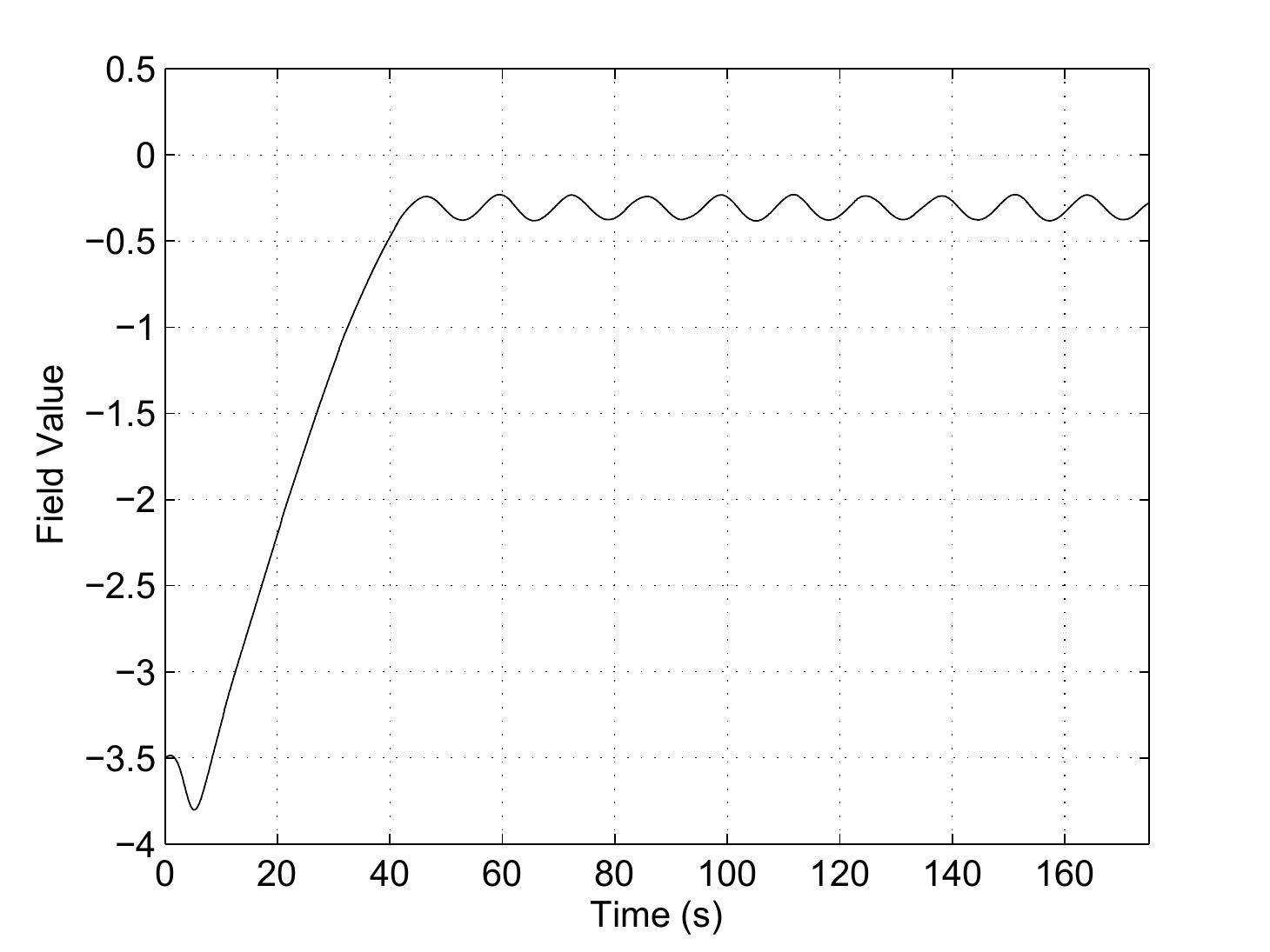}
  \caption{Evolution of the field value over the experiment. Field value is in metres.}
  \label{ext:fig.test5}
\end{figure}

\begin{figure}[ht]
  \centering \includegraphics[width=10cm]{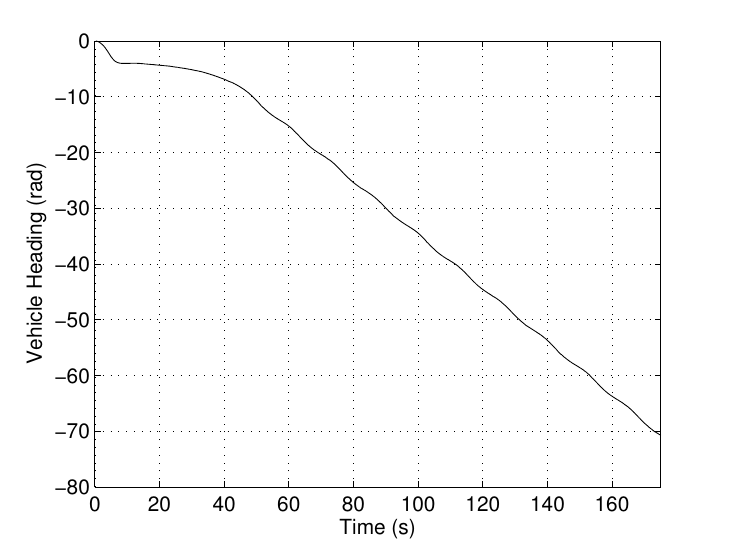}
  \caption{Evolution of the vehicle orientation over the
    experiment.}
  \label{ext:fig.test6}
\end{figure}

\begin{figure}[ht]
  \centering
  \subfigure[]{\scalebox{0.25}{\includegraphics{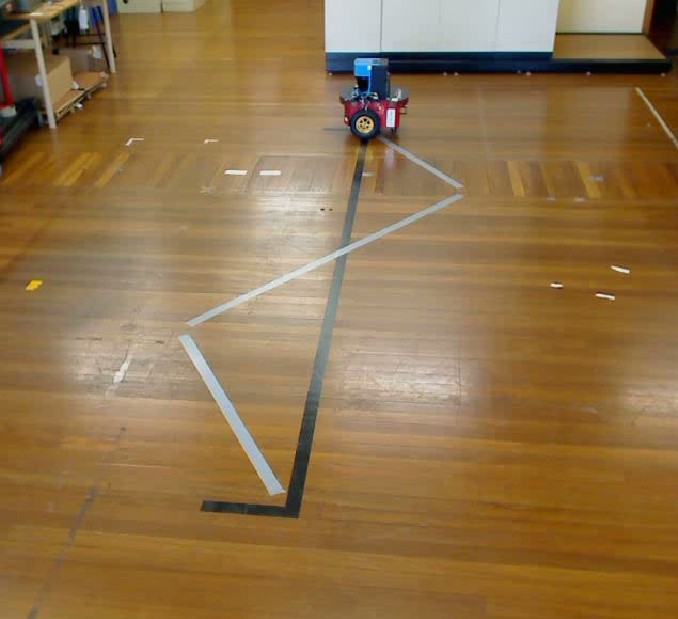}}}
  \subfigure[]{\scalebox{0.25}{\includegraphics{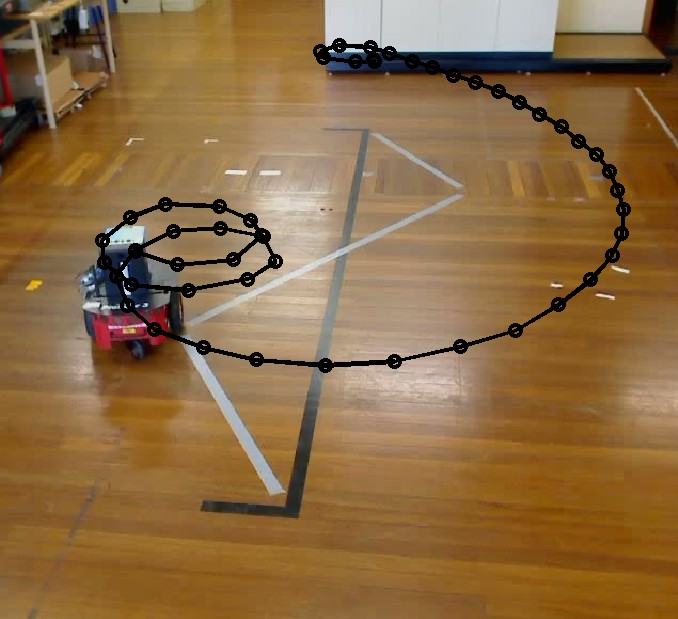}}}
  \subfigure[]{\scalebox{0.25}{\includegraphics{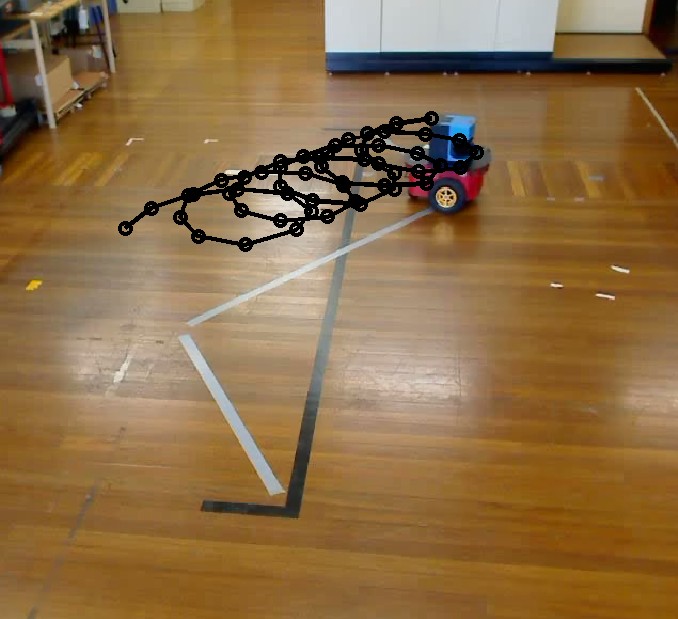}}}
  \subfigure[]{\scalebox{0.25}{\includegraphics{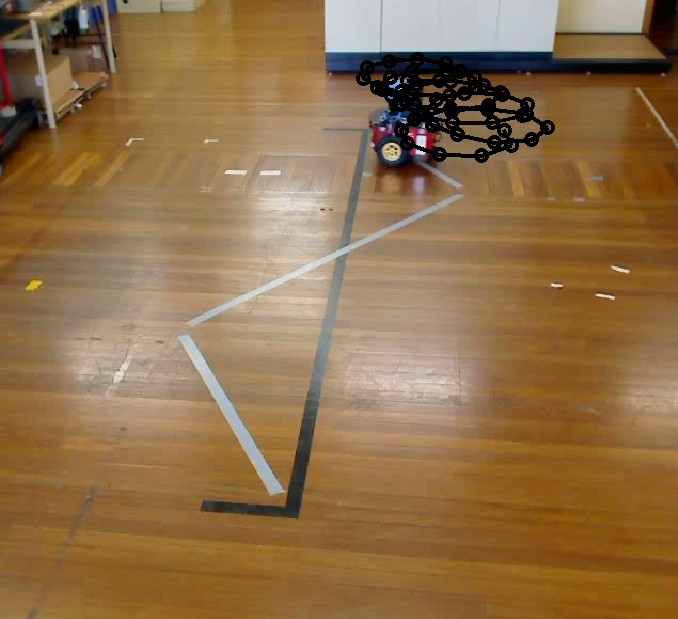}}}
  \caption{Sequence of images showing the experiment (irregularly moving field).} \label{ext:fig.test7}
\end{figure}

\begin{figure}[ht]
  \centering \includegraphics[width=10cm]{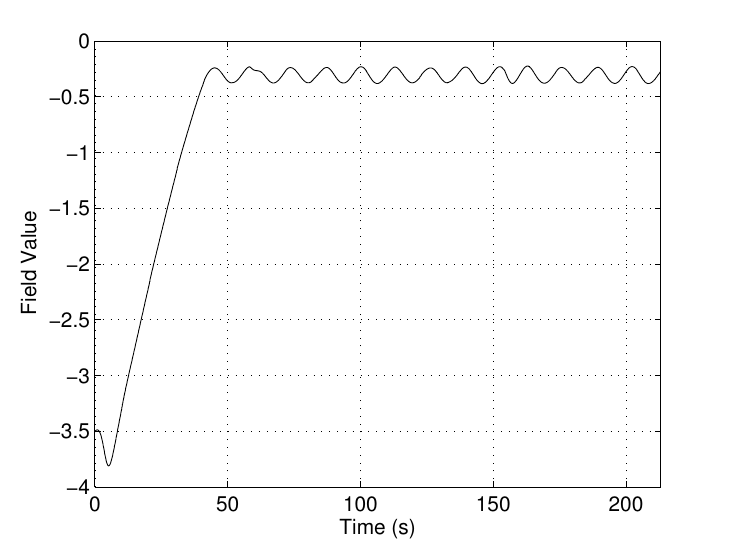}
  \caption{Evolution of the field value over the experiment. Field value is in metres.}
  \label{ext:fig.test8}
\end{figure}

\begin{figure}[ht]
  \centering
  \includegraphics[width=10cm]{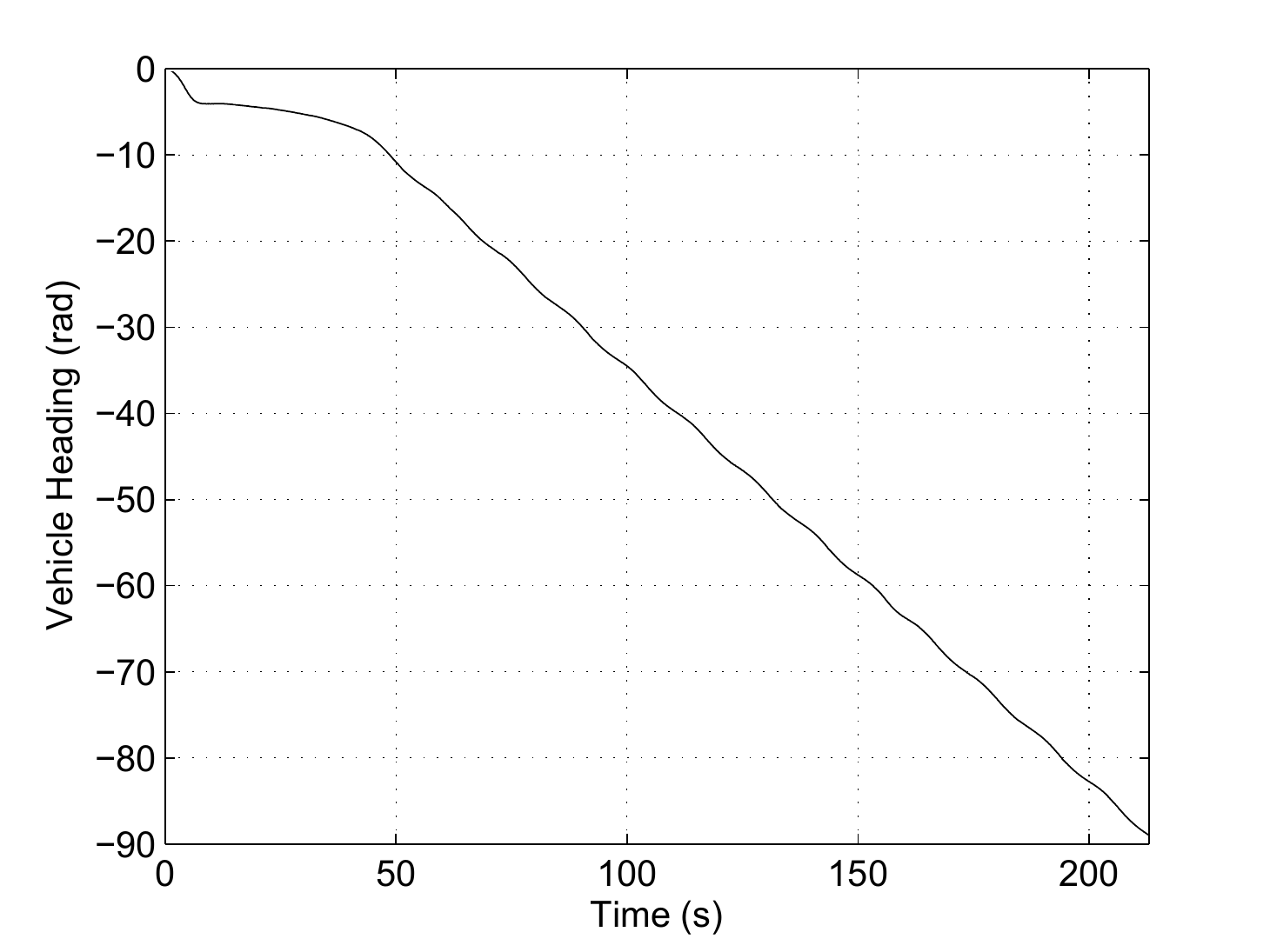}
  \caption{Evolution of the vehicle orientation over the
    experiment.}
  \label{ext:fig.test9}
\end{figure}

\clearpage \section{Summary}
\label{ext:sum}

A single kinematically controlled mobile robot traveling in a planar region supporting an unknown
and unsteady field distribution is considered, where only the current value of the field is known to the robot.
A reactive navigation strategy is presented that drives the robot towards the time-varying location where the field
distribution attains its spatial maximum and then keeps the robot in the pre-specified vicinity of that point. The
applicability and performance of the proposed guidance approach are confirmed by extensive simulation tests and
experiments with a real wheeled robot.
\chapter{Tracking the Level Set of a Scalar Field}
\label{chapt:lst}

This chapter introduces a method for tracking a specified level set or
isoline of an environmental scalar field. To achieve this,
a sliding mode control law is designed which can provably achieve the
desired behaviour, given set of assumptions about the field. Simulations and
experimental results are given to confirm the viability of the
proposed method.

\section{Introduction}

Recent environmental disasters
have highlighted the need for effective tools capable of timely
detection, exploration, and monitoring of environmental boundaries.
Some of these distributed phenomena can be observed as a whole using
large scale sensors. However, there are
many scenarios where such observation is troublesome, the quantity
can only be realistically measured using local sensors in a point-wise fashion at a
particular location (see e.g. \cite{CKBMLM06,HaMaGo02,LiDu03}). 
 Observation by means of a distributed network of static
sensors typically requires high deployment density and considerable
computational and communication loads to provide good accuracy
\cite{LiChGuZh02,NoMi03}.  More effective use of sensors is achieved
in mobile networks, where each sensor explores many locations. This is
especially beneficial whenever the interest is focused not so much on
the entire distributed phenomenon as on some dependent structure of a
lower dimension, like the boundary of an oil spill or radioactively
contaminated area. To derive this benefit, the mobile sensor should be
equipped with a motion control system by which it can detect and track
this structure.
\par
Recently, such problems have gained much interest in control
community. Among their basic setups, there is tracking of
environmental level sets: the robot should reach and then track the
curve where an unknown scalar field assumes a specific value and which
is thus the boundary of the area with greater values. In doing so, the
control law should use only point-wise field measurements along the
robot path.
\par
Many works in this area assume access to the field gradient or even
higher derivatives data such as the curvature of the isoline (see
e.g. \cite{MaBe03,SrRaKu08,ZhLe10,HsLoKu07}).
Some examples include gradient-based contour estimation by mobile
sensor networks \cite{MaBe03,BertKeMar04,SrRaKu08}, centralized
control laws originating from the `snake' algorithms in image
segmentation \cite{MaBe03,BertKeMar04}, cooperative distribution of
the sensors over the estimated contour with minimal latency
\cite{SrRaKu08}, artificial potential approach based on direct access
to the gradient \cite{HsLoKu07}, collaborative estimation of the
gradient and Hessian of the noise-corrupted field as the basis for
driving the center of a rigid formation of multiple sensors along a
level curve \cite{ZhLe10}.  However derivative-dependent data is often
unavailable, whereas its estimation requires access to the field
values at several nearby locations.  Even in the multiple sensor
scenario, such access may be degraded by limitations on communication
and may require ineffective concentration of sensors into a compact
cluster.
\par
A single mobile sensor with access to only point field values is the main
target for gradient-free approaches (see e.g.
\cite{BaRe03,KeBeMa04,CKBMLM06,Anders07}).
Control via switches between two steering angles depending on whether
the current field value is above or below the threshold have been proposed
\cite{ZhBer07,JAHB09}.  Similar approach with a larger set of
alternatives has been applied to an underwater vehicle equipped with a
profile sonar \cite{BaRe03}.  These methods typically result in a
zigzagging behavior and rely in effect on systematic side exploration
maneuvers to collect rich enough data.  A method to control an
unmanned aerial vehicle has been proposed based on segmentation of the infrared local
images of a forest fire \cite{CLBMM05}. These works
are based, more or less, on heuristics and provide no rigorous and
completed justification of the proposed control laws. A linear PD
controller fed by the current field value has been proposed 
for steering a unicycle-like vehicle along a level
curve of a field given by a radial harmonic function, and a local
convergence result was established for a vehicle with unlimited
control range \cite{BarBail07}. A sliding mode control method for tracking
environmental level sets without gradient estimation was offered in
\cite{MaTeSa12}.
\par
The characteristic feature of the previous research is that it dealt
with only steady fields, which means that the reference level curve
does not move or deform. However in real world, environmental fields
are almost never steady and often cannot be well approximated by
steady fields, whereas the theory of tracking the level sets for
dynamic fields lies in uncharted territory. As a particular case,
this topic includes navigation and guidance of a mobile robot towards
an unknowingly maneuvering target and further escorting it with a
pre-specified margin on the basis of a single measurement that decays as
the sensor goes away from the target, like the strength of the
infrared, acoustic, or electromagnetic signal, or minus the distance
to the target. Such navigation is of interest in many areas
\cite{ADB04,GS04,Matveev2011journ4}; it carries a potential to reduce
the hardware complexity and cost and improve target pursuit
reliability.  The mathematically
rigorous analysis of a navigation law for such problem was offered in
\cite{Matveev2011journ4} in the very special case of the unsteady
field -- the distance to an unknowingly moving Dubins-like target.
However the results of \cite{Matveev2011journ4} are not applicable to
more general dynamic fields.
\par
The navigation strategy considered in this chapter develops some ideas
set forth in \cite{MaTeSa12}; in particular, spatial gradient estimates and
systematic exploration maneuvers are not employed. However in
\cite{MaTeSa12}, only the case of static fields was examined. In this
work, it is shown that those ideas remain viable for much more general
scenarios with dynamic fields. Conditions
necessary for a Dubins-like vehicle to be capable of tracking the
moving and deforming level set of a dynamic field are established. It is shown that
whenever slight and partly unavoidable enhancements of these necessary
conditions hold, the problem can be solved by the proposed controller.
This is done by means of a mathematically rigorous non-local
convergence result, which contains recommendations on the controller
parameters tuning. 
\par
All proofs of mathematical statements are omitted here; they are
available in the original manuscript \cite{MaHoAnSa_sb}.
\par
The body of this chapter is organized as follows. In Sec.~\ref{lst:sec1} the problem is formally defined, 
 in Sec.~\ref{lst:sec.ass} the
the main assumptions are described. The main results are outlined in Sec.~\ref{lst:sec.mr}. Simulations and
experiments are presented in Secs.~\ref{lst:sec.simtest} and \ref{lst:sec.exper}. Finally, brief conclusions are given
in Sec.~\ref{lst:sec.soncl}.

 \section{Problem Statement}
\label{lst:sec1}

A planar mobile robot is considered, which travels with a constant speed
$v$ and is controlled by the time-varying angular velocity $u$ limited
by a given constant $\ov{u}$. The robot's workspace hosts an unknown
and time-varying scalar field $D(t,\boldsymbol{r}) \in \br$.  Here
$\boldsymbol{r}:= (x,y)^\trs$ is the vector of the absolute Cartesian
coordinates $x,y$ in the plane $\br^2$ and $t$ is time.  The objective
is to steer the robot to the level curve $D(t,\bldr) = d_0$ where the
distribution assumes a given value $d_0$ and to ensure that the robot
remains on this curve afterwards, circulating along it at the given
absolute speed $v$. The on-board control system has access to the
distribution value $d(t):= D(t,x,y)$ at the vehicle current location
$x=x(t), y =y(t)$ and is capable to access the rate $\dot{d}(t)$ at
which this measurement evolves over time $t$.  However, neither the
partial derivative $D^\prime_t$, nor $D^\prime_x$, nor $D^\prime_y$ is
accessible.
\par
The kinematic model of the robot is as follows:

\begin{equation}
  \label{lst:1}
  \begin{array}{l}
    \dot{x} = v \cos \theta,
    \\
    \dot{y} = v \sin \theta,
  \end{array},
  \quad
  \dot{\theta} = u , \quad |u|\leq \overline{u}, \quad
  \begin{array}{l}
    x(0) = x_{\text{in}}
    \\
    y(0) = y_{\text{in}}
  \end{array}, \quad \theta(0) = \theta_{\text{in}},
\end{equation}

Here $\theta$ gives the robot orientation.
It is required to design a controller that ensures the convergence
$D[t,x(t),y(t)] \to d_0$ as $t \to \infty$.
In this chapter, the following navigation law is examined:

\begin{equation}
  \label{lst:c.a}
  u(t)=-\sgn\{\dot{d}(t)+\chi[d(t)-d_{0} ] \} \bar{u}, \qquad d(t) = D[t,x(t),y(t)],
\end{equation}

Here $\chi(\cdot)$ is
a linear function with saturation:

\begin{equation}
  \label{lst:chi}
  \chi(p):=
  \begin{cases}
    \gamma p & \text{if}\; |p|\leq \delta \\
    \sgn(p)\mu & \text{otherwise}
  \end{cases}, \qquad \mu :=\gamma \delta .
\end{equation}

The gain coefficient $\gamma >0$ and the saturation threshold $\delta
>0$ are design parameters.
\par
For the discontinuous control law Eq.\eqref{lst:c.a}, the desired dynamics
\cite{Utkin1992book1} is given by $\dot{d}(t) = - \chi[d(t)-d_0]$.

 \section{Main Assumptions}
\label{lst:sec.ass}

The notation describing the field is taken from Chapt.~\ref{chapt:ext}.
The conditions employed in this chapter are described by the following:

\begin{Proposition}
  \label{lst:lem.23}
  Suppose that the robot moves so that it remains on the required
  isoline $D[t,r(t)] \equiv d_0$ and in a vicinity of this isoline,
  the function $D(\cdot,\cdot)$ is twice continuously differentiable
  and $\nabla D(\cdot,\cdot) \neq 0$. Then at any time the front speed
  of the isoline at the robot location does not exceed the speed of
  the robot
  \begin{equation}
    \label{lst:speed}
    \left| \lambda[t,r(t)] \right| \leq v,
  \end{equation}
  the robot's velocity
  \begin{equation}
    \label{lst:e}
    \vec{v} = v \vec{e}, \qquad \vec{e}:=\left( \begin{array}{c}
        \cos \theta
        \\
        \sin \theta
      \end{array}
    \right)
  \end{equation}
  has the form
  \begin{equation}
    \label{lst:vec.c}
    \vec{v} = \lambda N \pm  T \sqrt{v^2-\lambda^2},
  \end{equation}
  and the following inequality is true
  \begin{equation}
    \label{lst:accel}
    \left| \pm 2 \omega  + \frac{\alpha}{\sqrt{v^2-\lambda^2}} + \kappa \sqrt{v^2-\lambda^2}\right| \leq \overline{u}.
  \end{equation}
  In $\pm$, the sign $+$ is taken if the robot travels along the
  isoline in the positive direction (i.e., so that the domain $\{\bldr
  : D(t,\bldr) > d_0\}$ is to the left), and $-$ is taken otherwise.
\end{Proposition}

For the control objective to be attainable, Eq.\eqref{lst:speed} and
Eq.\eqref{lst:accel} should hold at the current location of the robot at any
time.  Since this location is not known in advance, it is reasonable
to extend this requirement on all points on the isoline.  However,
some form of controllability is also required to make the objective of
driving the robot to the pre-specified isoline realistic. For example,
the robot should be capable of moving from a given isoline
$I(t,d_\ast)$ to the area of larger field values $\{\bldr: D(t,\bldr)
> d_\ast\}$, as well as to that of smaller ones. Since remaining on
the isoline implies $\dot{d}=0$, this property means that the sign of
the second derivative $\ddot{d}$ can be made both positive and
negative by respective choices of feasible controls $u \in [- \ov{u},
\ov{u}]$. 
\par
It is assumed that such controllability holds in the entire zone of the
robot's maneuver $\mathcal{M}$, which is characterized by the extreme
values $d_- \leq d_+$ taken by the field in this zone:

\begin{equation}
  \label{lst:m}
  \mathcal{M} := \{(t,\bldr): d_- \leq D(t,\bldr) \leq d_+\}
\end{equation}

It is also assumed this contains the required isoline $d_- \leq d_0 \leq d_+$.  Finally,
it is assumed that the above strict inequalities do not degrade as time
progresses or location $\bldr$ goes to infinity. As a result, the following 
assumption is arrived at:

\begin{Assumption}
  \label{lst:ass.upr}
  The field $D(\cdot,\cdot)$ is twice continuously differentiable in
  the domain Eq.\eqref{lst:m} and there exist constants $\Delta_\lambda >0$
  and $\Delta_u >0$ such that the following enhanced analogs of
  Eq.\eqref{lst:speed} and Eq.\eqref{lst:accel} hold:
  \begin{equation}
    \label{lst:deltalambda}
    \left| \lambda \right| \leq v - \Delta_\lambda, \quad \left| \pm 2 \omega  + \frac{\alpha}{\sqrt{v^2-\lambda^2}} + \kappa \sqrt{v^2-\lambda^2}\right| \leq \overline{u} - \Delta_u \qquad \forall (t,r) \in \mathcal{M},
  \end{equation}
  where the second inequality is true with the both signs in $\pm$.
\end{Assumption}

Since the trajectory of the robot is smooth, it is natural to exclude
isoline singularities:

\begin{Assumption}
  \label{lst:ass.grad}
  The field has no spatial singularities $\nabla D \neq 0$ in the
  domain Eq.\eqref{lst:m}, and this property does not degrade as time
  progresses: there exists $b_\rho >0$ such that
  $\rho(t,\bldr)=\|\nabla D (t,\bldr)\| \geq b_\rho^{-1} \; \forall
  (t,\bldr) \in \mathcal{M}$.
\end{Assumption}

The next assumption is typically fulfilled in real world, where
physical quantities take bounded values:

\begin{Assumption}
  \label{lst:ass.last}
  There exist constants $b_\lambda, b_\tau, b_n, b_\varkappa, b_v,
  b_\alpha$ such that the following inequalities hold:
  \begin{equation}
    \label{lst:estim}
    |\lambda| \leq b_\lambda, \quad |\tau_\rho| \leq b_\tau, \quad |n_\rho| \leq b_n,
    \quad |\varkappa| \leq b_\varkappa, \quad |v_\rho| \leq b_v, \quad |\alpha| \leq b_\alpha \qquad \forall (t,r) \in \mathcal{M}.
  \end{equation}
\end{Assumption}

The last assumption is partly underlaid by the fact that under the
control law Eq.\eqref{lst:c.a} with properly tuned parameters, the vehicle
initially moves with $u\equiv \pm \ov{u}$ over an {\em initial circle}
$C_\pm^{\text{in}}$, which is defined as the related path starting
with the given initial data from Eq.\eqref{lst:1}.  It is required that these
circles lie in the operational zone Eq.\eqref{lst:m}, along with the
encircled disc's $D_\pm^{\text{in}}$ (also called {\it initial}).
Furthermore, there exists some initial time interval during which the
average angular speed of the spatial gradient rotation is less than
the maximal turning rate of the robot. This leads to the
following:

\begin{Assumption}
  \label{lst:ass.disk}
  There exists a natural $k$ such that during the time interval
  $\left[ 0,T_k\right],$ $T_k$ $:=$ $\frac{2 \pi k}{\ov{u}}$, (a) the gradient
  $\nabla D(t,\bldr_{\text{in}}), \bldr_{\text{in}}:= (x_{\text{in}},
  y_{\text{in}})^\trs$ rotates through an angle that does not exceed
  $2 \pi (k-1)$ and (b) the both initial disc's lie in the domain
  Eq.\eqref{lst:m}, i.e., $[0,T_k] \times D_\pm \subset \mathcal{M}$.
\end{Assumption}

 \section{Summary of Main Results}
\label{lst:sec.mr}

The main theoretical result may now be stated:

\begin{Theorem}
  \label{lst:th.main}
  Suppose that Assumptions~{\rm \ref{lst:ass.upr}}--{\rm
    \ref{lst:ass.disk}} hold and the parameters $\gamma,
  \mu=\gamma\delta$ of the controller Eq.\eqref{lst:c.a} satisfy the
  following inequalities, where $\mu_\ast$ $:=$ $b_\rho \mu$ and
  $\sigma(\mu_\ast)$ $:=$ $\sqrt{\Delta_\lambda^2 - 2 v \mu_\ast -
    \mu^2_\ast}$:
  \begin{equation}
    \label{lst:limit}
    0< \mu_\ast  <  \sqrt{v^2+\Delta_\lambda^2}-v, \quad
    \left(3 b_\tau + \frac{b_\varkappa+2b_v + \gamma +b_n \mu_\ast}{\sigma(\mu_\ast)} + \frac{b_\alpha}{\sigma(\mu_\ast)^3}\right) \mu_\ast < \Delta_u .
  \end{equation}
  Then the robot driven by the navigation law Eq.\eqref{lst:c.a} achieves the
  control objective $d(t)\xrightarrow{t \to \infty} d_{0}$ and moves
  in the domain Eq.\eqref{lst:m}.
\end{Theorem}

 \section{Simulations}
\label{lst:sec.simtest} 

Simulations were carried out with the Dubins-like robot Eq.\eqref{lst:1}
driven by the control law Eq.\eqref{lst:c.a}.  The numerical values of the
parameters used for simulations are shown in
Table~\ref{lst:table:param}, where $u_d$ is the unit of measurement of
$d=D(t,\bldr)$ and the controller parameters were chosen based on
recommendations from Theorem~\ref{lst:th.main}. The control was
updated with the sampling time of $0.1 s$.

\begin{table}[ht]
  \centering
  \begin{tabular}{| l | c |}
    \hline
    $v$ & $1 m/sec$  \\
    \hline
    $\mu$ & $0.5 u_d/sec$  \\
    \hline
    $\gamma$ & $0.04 sec^{-1}$ \\
    \hline
  \end{tabular}
  \caption{Simulation parameters for level set tracking controller.}
  \label{lst:table:param}
\end{table}

Figs.~\ref{lst:fig:lin} and \ref{lst:fig:linang} present the results
of tests in the moving radial field
$$
D(t,\bldr) = 70-0.8 \cdot \lVert \bldr - \bldr^0(t)\rVert ,
$$
where the source $\bldr^0(t)$ moves to the right at the speed of $0.3
m/sec$.  Fig.~\ref{lst:fig:linang} demonstrates successful convergence
to the desired field value $d_0=30$. Fig.~\ref{lst:fig:lin} shows the
related path of the robot; in Figs.~\ref{lst:fig:lin},
\ref{lst:fig:noise} and \ref{lst:fig:heat}, the position of the field
source is depicted by the solid black line. Since the isoline of the
unsteady radial field at hand is a circle undergoing a constant
velocity displacement, the robot that moves with a constant speed and
does not leave this isoline should trace a cycloid (either curtate or
prolate). As can be seen in Fig.~\ref{lst:fig:lin}, it does trace a
prolate cycloid-like path, up to the transient.

 \begin{figure}[ht]
   \centering
   \subfigure[]{\includegraphics[width=0.75\columnwidth]{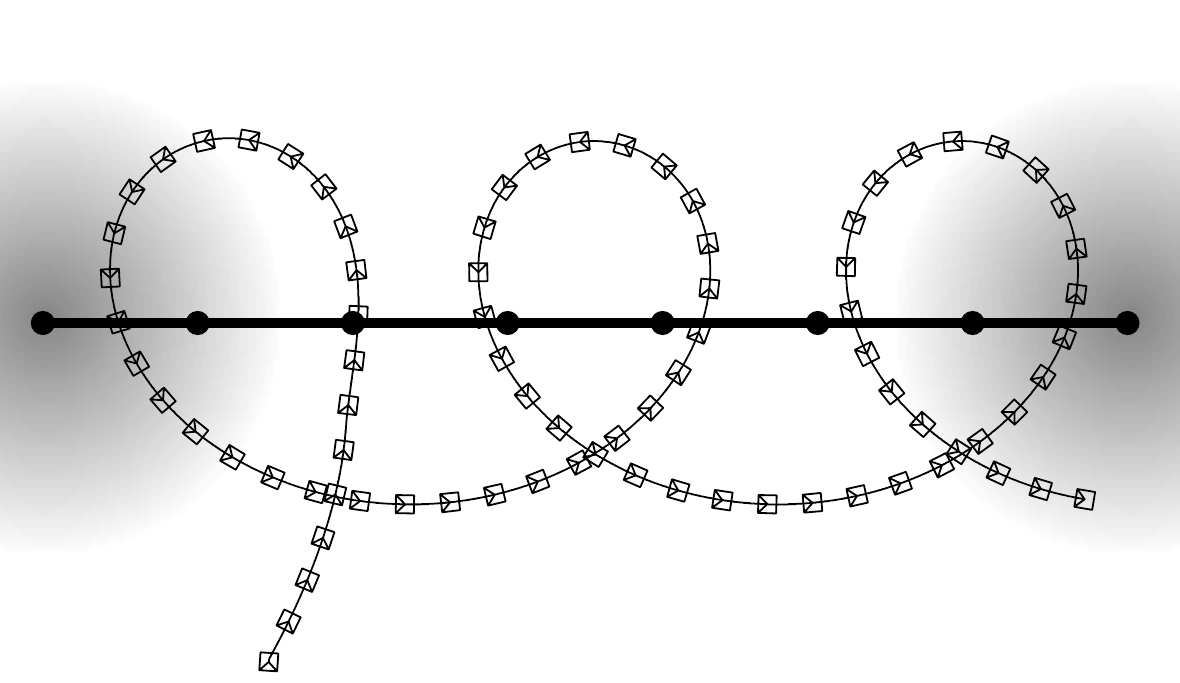}\label{lst:fig:lin}}
   \subfigure[]{\includegraphics[width=10cm]{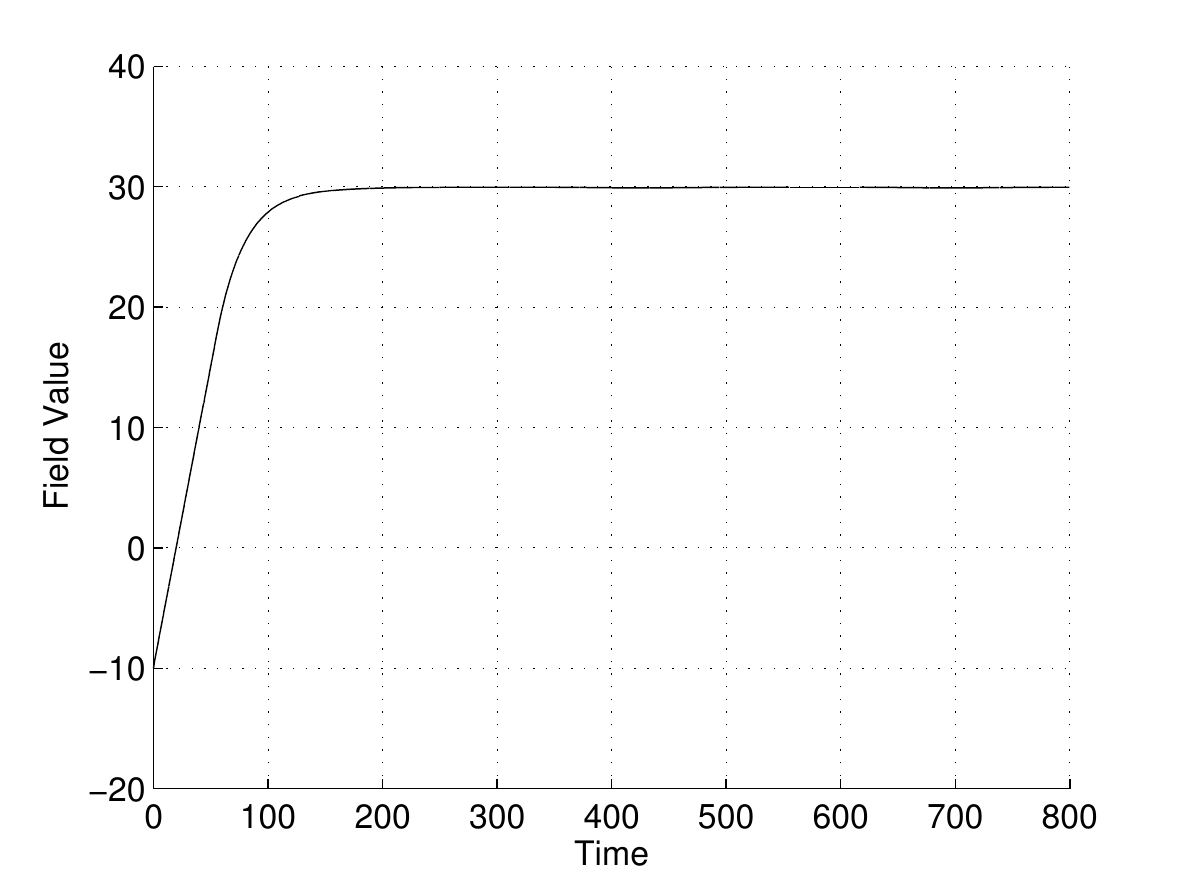}\label{lst:fig:linang}}
   \caption{Simulations with an unsteady
     radial field; (a) Path; (b) Robot's orientation. Time is in seconds, field value is in arbitrary units.}
 \end{figure}

 In Figs.~\ref{lst:fig:noise} and \ref{lst:fig:noiseang}, the effect
 of the measurement noise was examined in a similar simulation setup.
 Specifically, the measurements $d(t)$ and $\dot{d}(t)$ were
 individually corrupted by random additive noises uniformly
 distributed over the intervals $[-2.5 u_d, 2.5 u_d]$ and $[-2.0
 u_d/sec, 2.0 u_d/sec]$, respectively. The scalar field was corrupted
 by sinusoidal plane waves:
$$
D(t,\bldr) = 70-0.8 \cdot \lVert \bldr - \bldr^0(t)\rVert +
5\cdot\left[\sin(0.05 \cdot x) + \sin (0.05 \cdot y) \right].
$$

Figs.~\ref{lst:fig:noiseang} and \ref{lst:fig:noise} show that the
control objective is still achieved with a good exactness, though the
path becomes less regular.

 \begin{figure}[ht]
   \centering
   \subfigure[]{\includegraphics[width=0.75\columnwidth]{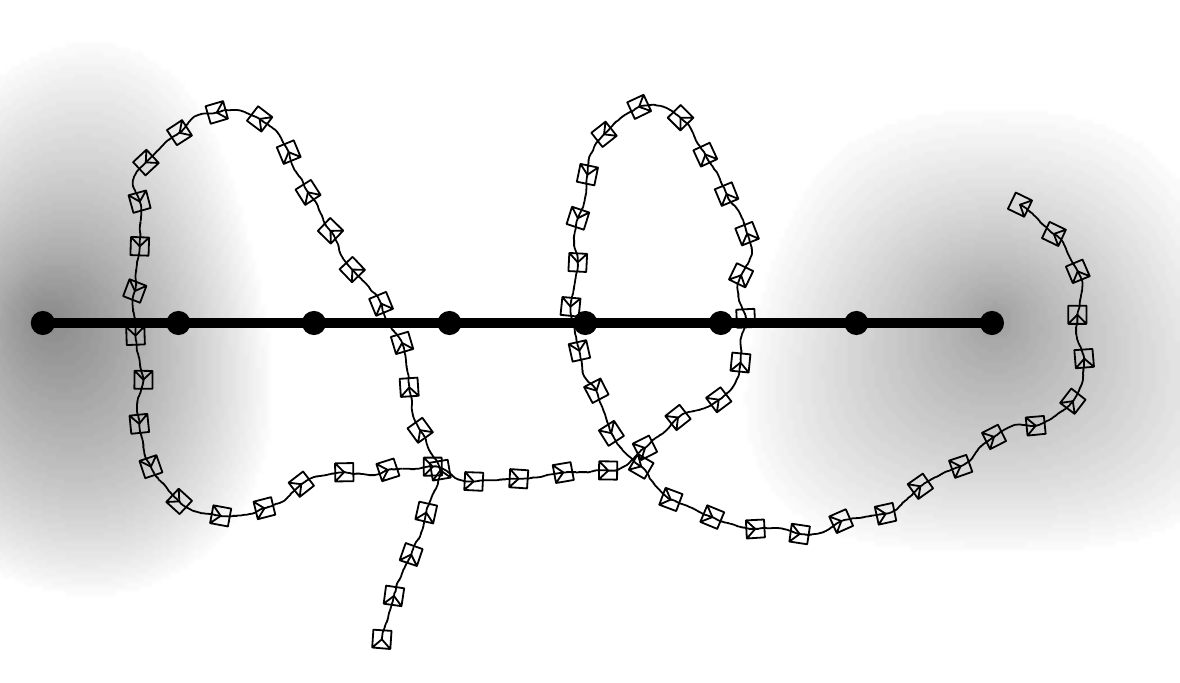}\label{lst:fig:noise}}
   \subfigure[]{\includegraphics[width=10cm]{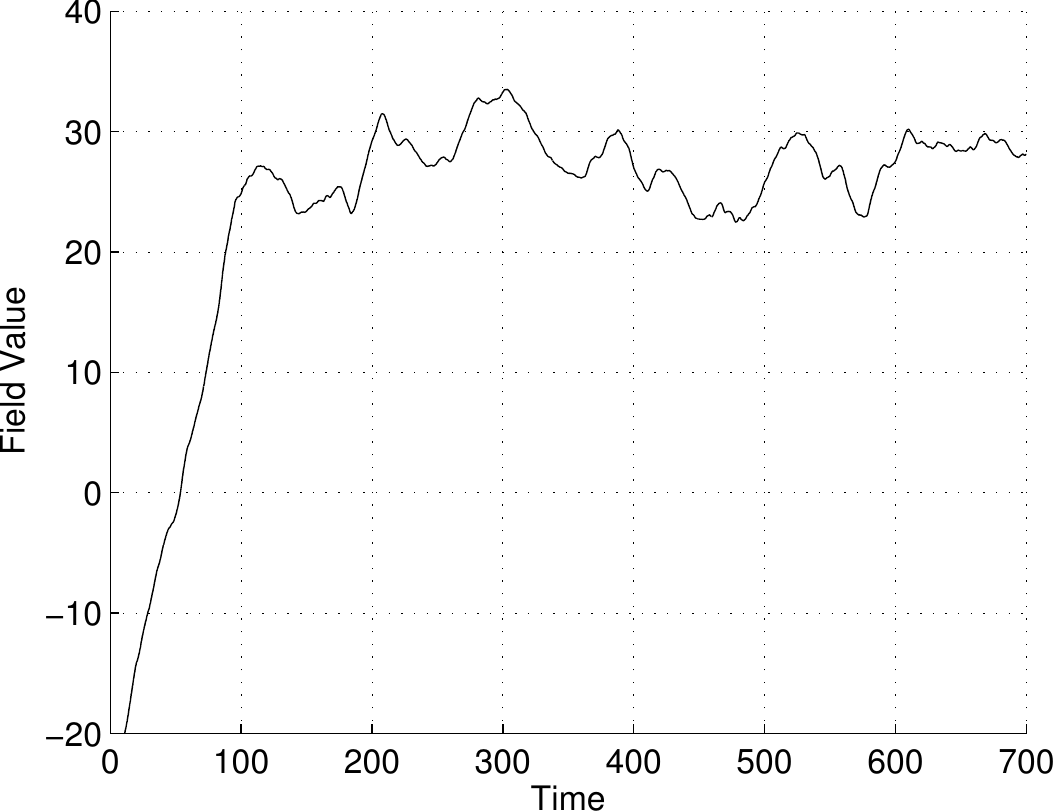}\label{lst:fig:noiseang}}
   \caption{Simulations with measurement noise; (a) Path; (b) Robot's orientation. Time is in seconds, field value is in arbitrary units.}
 \end{figure}

 Simulations were also carried out for a more realistic model of a
 time-varying field caused by a constant-rate emanation of a certain
 substance, which was taken to be the same as that in Chapt.~\ref{chapt:ext}.
 The results of these simulations are shown
 Figs.~\ref{lst:fig:heatang} and \ref{lst:fig:heat}, where the source
 of diffusion undergoes irregular motion depicted by the solid black
 curve in Fig.~\ref{lst:fig:heat}. The small impulse in the field
 value at $\approx 600 s$ is a result of breaking the limitations
 revealed by Proposition~\ref{lst:lem.23}. As can be seen, the control
 objective is still achieved with a good accuracy.

  \begin{figure}[ht]
   \centering
   \subfigure[]{\includegraphics[width=0.75\columnwidth]{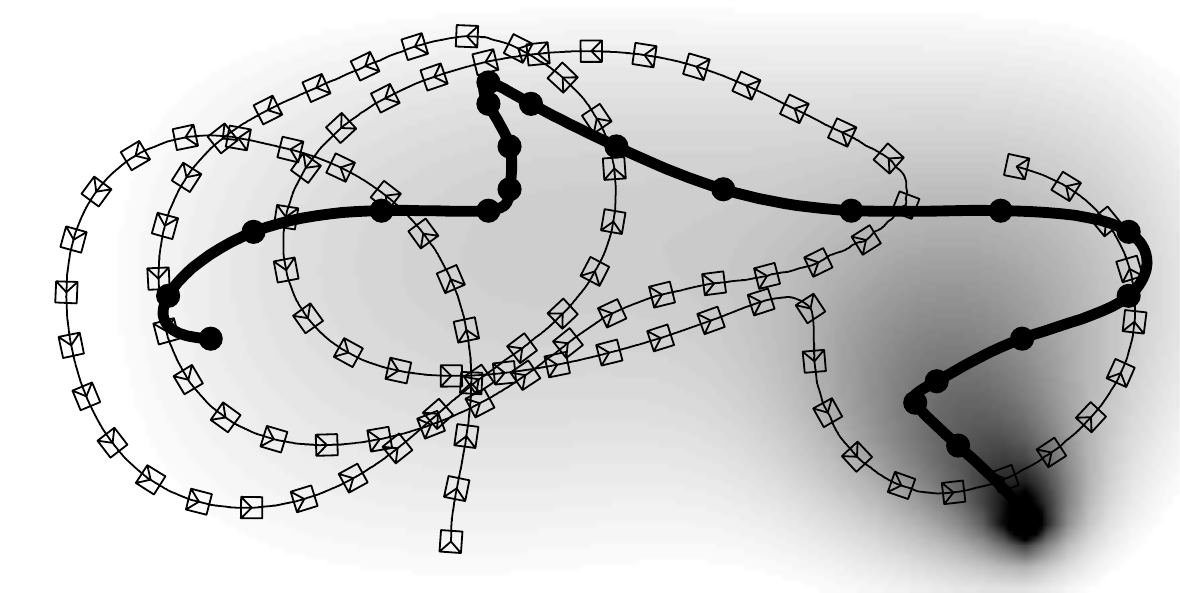}\label{lst:fig:heat}}
   \subfigure[]{\includegraphics[width=10cm]{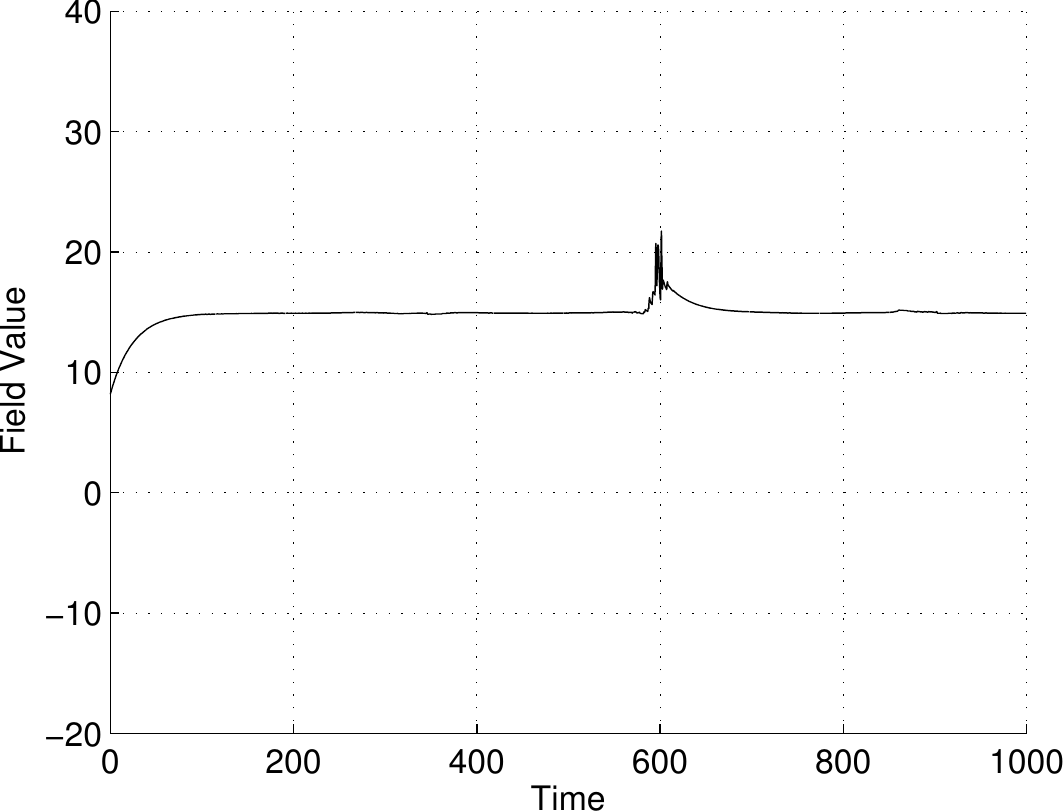}\label{lst:fig:heatang}}
   \caption{Simulations with a heat source; (a) Path; (b) Robot's orientation. Time is in seconds, field value is in arbitrary units.}
 \end{figure}

 \clearpage \section{Experiments}
 \label{lst:sec.exper}

 Experiments were carried out with an Activ-Media Pioneer 3-DX wheeled
 robot using its on-board PC and the Advanced Robot Interface for
 Applications (ARIA 2.7.2), which is a C++ library providing an
 interface to the robot's angular and translational velocity
 set-points. 

 The origin of the reference frame was co-located with the center of
 the robot in its initial position, its ordinate axis is directed
 towards the viewer in Figs.~\ref{lst:fig.test1} and
 Eq.\eqref{lst:fig.test3}. In two groups of experiments, a virtual point-wise
 source of the scalar field moved, respectively, as follows:

 \begin{enumerate}
 \item The source moved from the point with the coordinates $(0, 3.5)
   m$ with a constant translational velocity $(0, -0.02) ms^{-1}$
   outwards the viewer in Fig.~\ref{lst:fig.test1}. The experiment was
   run for $175 s$ so that the final position of the source coincided
   with the origin of the reference frame.
 \item The source moved at the constant speed $0.02 ms^{-1}$ from the
   same initial position along a piece-wise linear path through the
   points $(0, 3.5) m$, $(0.6, 2.5) m$, $(-0.6, 1.0) m$ and $(0, 0)
   m$; the motion was still outwards the viewer in
   Fig.~\ref{lst:fig.test3}. The experiment was run for $213 s$.
 \end{enumerate}

 The path of the source is displayed in Fig.~\ref{lst:fig.test1} by a
 long black tape, and in Fig.~\ref{lst:fig.test3} by long gray
 tape.\footnote{A short perpendicular black segment was added for
   calibration purposes and is not a part of the path.} The examined
 field was minus the distance to the moving source, which was accessed
 via odometry, whereas the source motion was virtual and emulated by
 computer.
 \par
 The parameters used in the experiments are shown in
 Table~\ref{lst:fig:paramexp}.

 \begin{table}[ht]
   \centering
   \begin{tabular}{| l | c |}
     \hline
     $v$ & $0.15 m/s$  \\
     \hline
     $\gamma$ & $0.2 s^{-1}$ \\   
     \hline
   \end{tabular}
\hspace{10pt}
\begin{tabular}{| l | c |}
     \hline
$\mu$ & $0.12 m/s$  \\
\hline
 $d_0$  & $- 0.8 m$\\
\hline   
\end{tabular}
   \caption{Experimental parameters for level set tracking controller.}
   \label{lst:fig:paramexp}
 \end{table}

 In the control law Eq.\eqref{lst:c.a}, the signum function $\sgn(\cdot)$ was
 replaced by a linear function with saturation $\text{\bf sat}\,(S) :=
 \sgn(S) \cdot \text{min} \{10\cdot|S|, \bar{u}\}$ with $\bar{u}=1.0$.
 The control law was updated at the rate of $0.1 s$.
 \par
 Typical results of experiments from the first and second groups are
 presented on Figs.~\ref{lst:fig.test1}, \ref{lst:fig.test2} and
 \ref{lst:fig.test3}, \ref{lst:fig.test4}, respectively. The desired
 field value is indicated by the thick black line in
 Figs.~\ref{lst:fig.test2} and \ref{lst:fig.test4}. In these
 experiments, like in the others, the robot successfully arrives at
 the desired field value and then maintains it until the end of the
 experiment. By Figs.~\ref{lst:fig.test2} and \ref{lst:fig.test4}, the
 steady state tracking error approximately equals $0.1 m$.

\begin{figure}[ht]
  \centering
  \subfigure[]{\scalebox{0.25}{\includegraphics{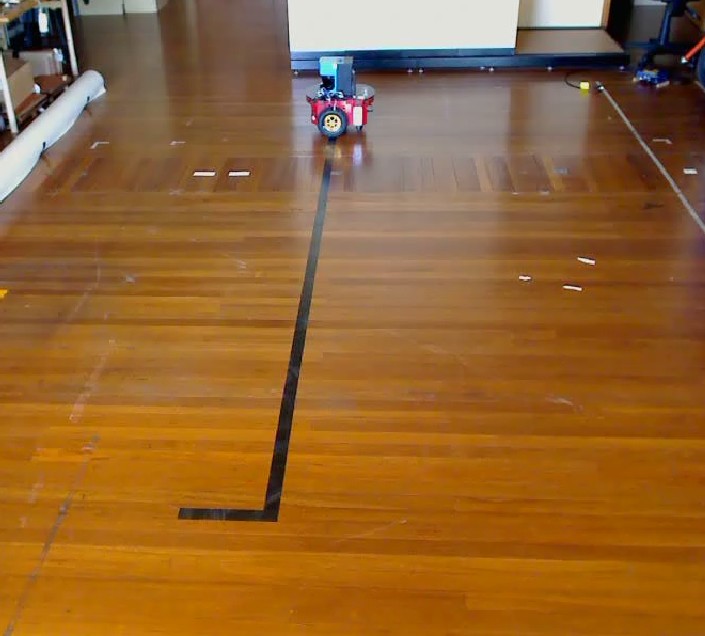}}}
  \subfigure[]{\scalebox{0.25}{\includegraphics{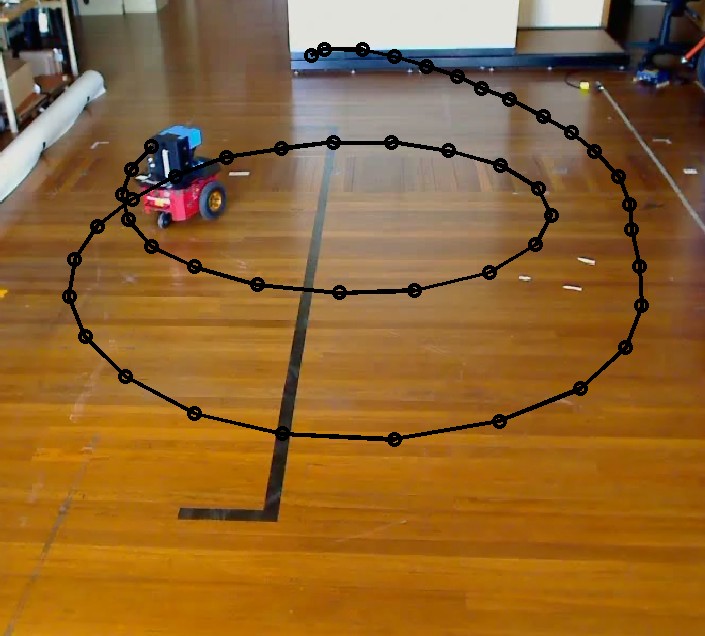}}}
  \subfigure[]{\scalebox{0.25}{\includegraphics{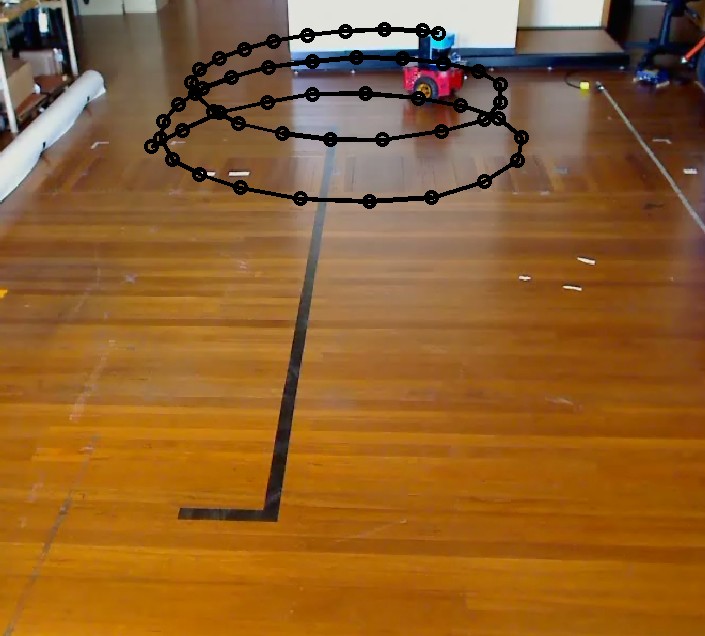}}}
  \caption{Sequence of images showing the experiment.} \label{lst:fig.test1}
\end{figure}

\begin{figure}[ht]
  \centering \includegraphics[width=10cm]{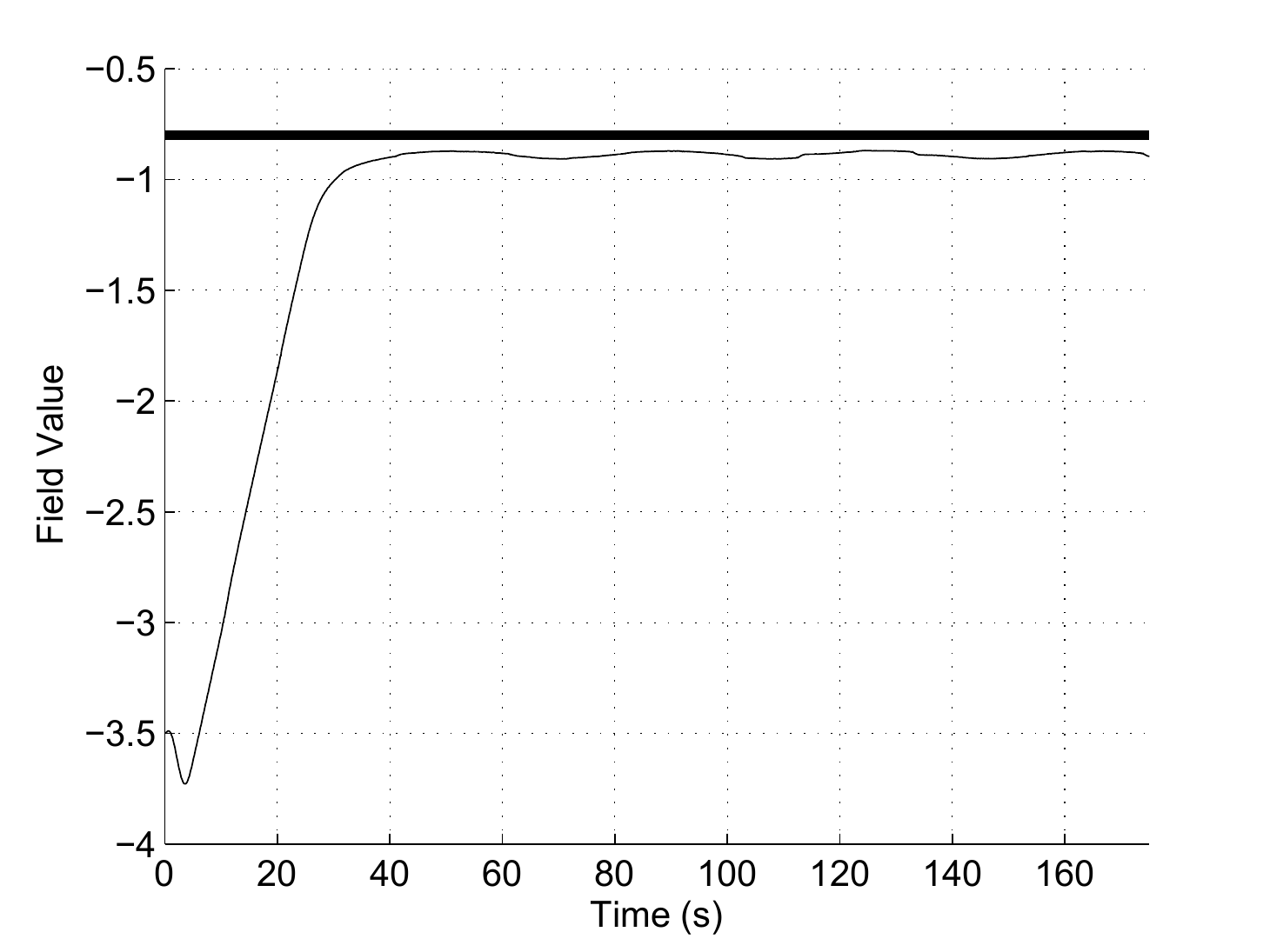}
  \caption{Evolution of the field value over the experiment. Field value is in metres.}
  \label{lst:fig.test2}
\end{figure}

\begin{figure}[ht]
  \centering
  \subfigure[]{\scalebox{0.25}{\includegraphics{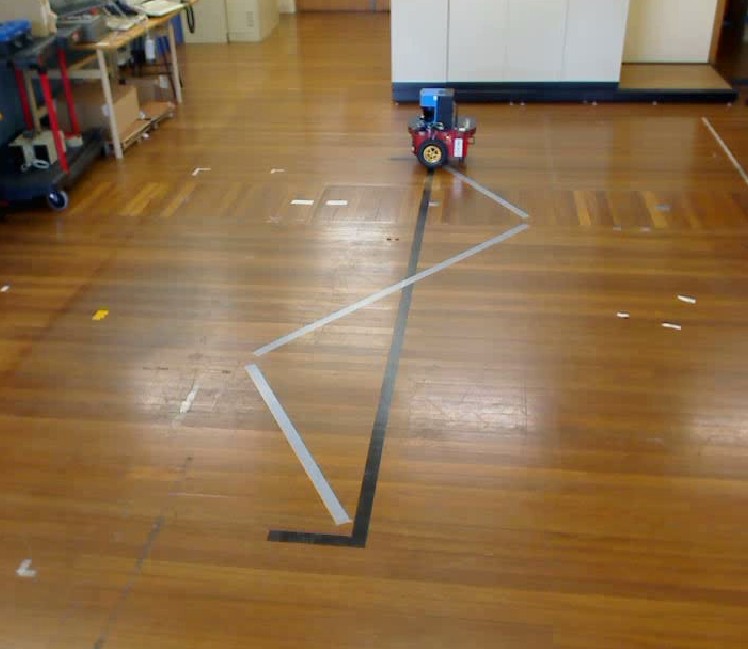}}}
  \subfigure[]{\scalebox{0.25}{\includegraphics{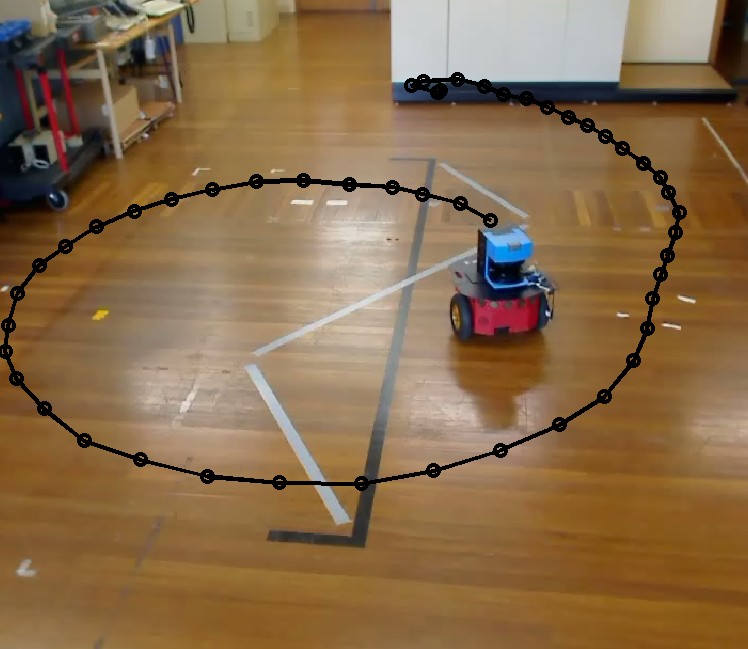}}}
  \subfigure[]{\scalebox{0.25}{\includegraphics{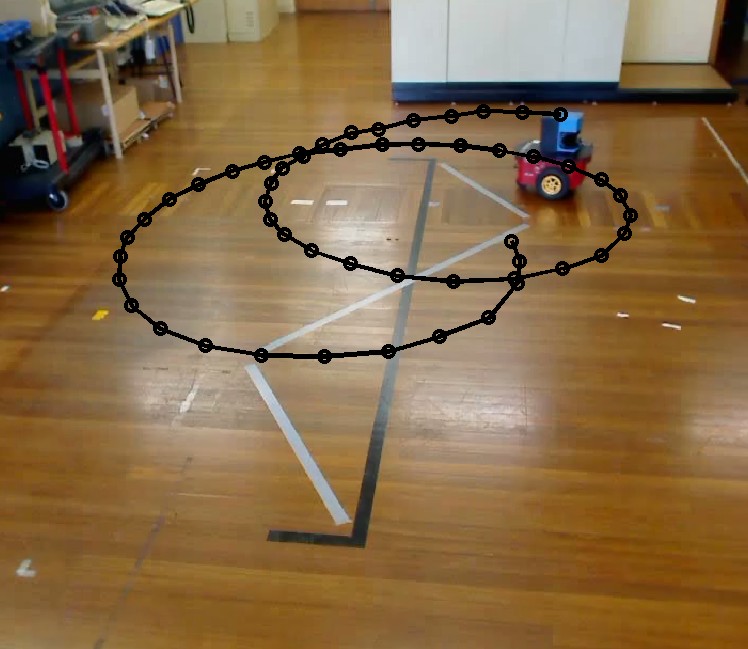}}}
  \subfigure[]{\scalebox{0.25}{\includegraphics{levelset/imagelevz3.jpg}}}
  \caption{Sequence of images showing the experiment.} \label{lst:fig.test3}
\end{figure}

\begin{figure}[ht]
  \centering \includegraphics[width=10cm]{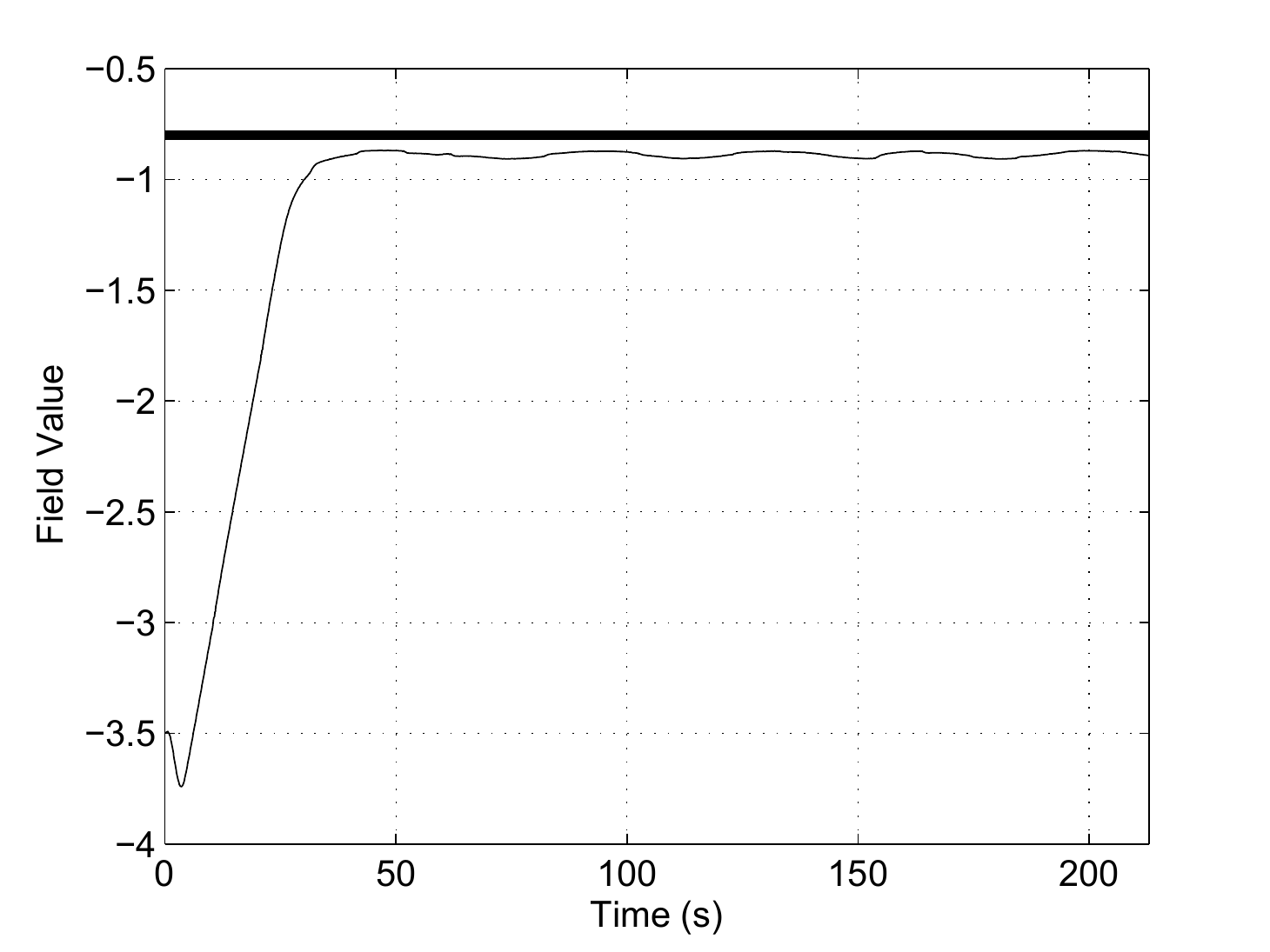}
  \caption{Evolution of the field value over the experiment. Field value is in metres.}
  \label{lst:fig.test4}
\end{figure}

\clearpage \section{Summary}
\label{lst:sec.soncl}
The chapter presented a sliding-mode control law that drives a single
non-holonomic Dubins-car like vehicle to the moving and deforming
level curve of an unknown and time-varying scalar field distribution
and ensures its ultimate circulation along this curve. The vehicle
travels at a constant speed and is controlled by the turning radius;
the sensor data are limited to the distribution value at the vehicle
current location. This proposed control law does not employ gradient
estimation and is non-demanding with respect to both motion and
computation. Its
performance was confirmed by computer simulations and experiments with
real robots.

\chapter{Decentralized Target Capturing Formation Control}
\label{chapt:tcp}

In this chapter a decentralized formation control problem is discussed, where
the aim is to drive a group of vehicles into an equidistant formation
along a circle, which is centred on a given target point. Sensor information is limited to scalar 
measurements about the target and surrounding vehicles. Simulations and
experimental results are given to confirm the viability of the
proposed method.

\section{Introduction}

A typical objective in multiagent systems research is to drive a team of vehicles to a
pre-specified formation. Formation control is one of the basic
technologies that enables multiple vehicles to cover large operational
areas and accomplish complex tasks \cite{SePaLe06}.  Although a
rigorous stability analysis of such systems is generally an extremely
hard task, partly due to time-varying topology of the information
flows \cite{FM02}, some theoretical results have been obtained.
\par
Recently, a substantial body of research was devoted to designs of
control strategies that drive multiple vehicles to a capturing
formation around a target object (see e.g.
\cite{MBF04,SiGh07,MBF06,SePaLe06,KiSu07,HiNa08,ShFiAn08,HiNa08a,KaBeAb08,HiNa09,KiHaHo10,KoHo11,Yamaguchi03,ShKoTaPo10,TsHaSaKi10,CeMaGaGi08,GYL10,LYL10}}.
This is motivated by various practical applications such as security
and rescue operations, explorations in hazardous terrestrial or marine
environments, deployments of mobile sensor networks, and
patrolling missions using multiple unmanned vehicles \cite{MBF04,MBF06,SePaLe06}. 
In fact, the maneuver itself
consists in enclosing and grasping. The former is to enclose the
target object while approaching it, whereas the latter is to grasp the
object by driving the vehicles to an optimal configuration around it,
which is typically the equal spacing formation.  Basically, the both
objectives should be achieved in a decentralized fashion.
\par
A decentralized capturing kinematic control law for multiple
point-wise planar vehicles was proposed in \cite{KoHo11}. This law
follows the gradient decent approach and employs only local data at
any agent (about the target and angularly closest neighbors). A
distributed motion coordination strategy for multiple Hilare-type
robots in cooperative hunting operations has been proposed
 and presents arguments supporting local stability of the
closed-loop system \cite{Yamaguchi03}.
\par
Conditions have been established under which simple
decentralized linear control laws can drive a team of identical linear
planar fully actuated agents into a given target-centric formation 
\cite{TsHaSaKi10}. However, the target is assumed
to obey a known and favorable for the capturing task escaping rule,
which is hardly realistic for most target-capturing scenario.
\par
Several algorithms are inspired by the
cyclic pursuit animal behavior \cite{MBF04,SiGh07,MBF06}. Though the pure cyclic pursuit
typically results in rendezvous for linear models, it was shown that
in the nonlinear case (planar unicycle-like models), identical
multiple vehicles can assemble a locally stable circular formation
under certain circumstances \cite{MBF04,MBF06}. For non-identical
unicycles and control gains, necessary conditions for the existence of a
circular equilibrium have been obtained \cite{SiGh07}. However, the
resultant circle parameters (the center and radius) are not under
control and are not related to a specific target, and not much was
established about the global convergence of the proposed algorithms.
By modifying the cyclic pursuit strategy, a
distributed cooperative control scheme for capturing a target in 3D
space by a team of velocity-controlled identical point-mass
vehicles has been proposed \cite{KiSu07}. This result was generalized for the case
of under-actuated identical stable vehicles with common dynamics
described by a linear MIMO model \cite{KiHaHo10}.
\par
A serious limitation of all pursuit-based approaches is that they
assume a fixed ring-like information-flow graph.\footnote{The pursuer
  is able to identify its prey agent and constantly sees it or
  communicates with it.} More general and time-varying
information-flow topologies have been examined, where a
control strategy can proposed that drives a set of identical
under-actuated planar unicycles into a cyclic formation \cite{SePaLe06}. However, the
center of this formation is out of control and is not concerned with a
particular target. Algorithms for enclosing a given moving target by a
set of fully actuated planar kinematically controlled and identical
but having an identity unicycles have been proposed for the case of general fixed communication
topology with some special properties \cite{HiNa08,HiNa08a}. This result was extended on a
time-varying topology and kinematically controlled point-mass vehicles
in 3D space \cite{HiNa09} and planar dynamically controlled ones with
inaccurate target information \cite{ShKoTaPo10} or information
exchange uncertainty \cite{SaMa10}. The information
graph may also be assumed independent of the control law, and the required properties are treated as
granted \cite{SePaLe06,HiNa08,HiNa08a,HiNa09,ShKoTaPo10}. However, in many cases this graph is determined by the
current locations of the vehicles, along with their limited visibility
or communication ranges, and so its properties essentially depend on
the control laws driving the vehicles.
\par
A scenario with position-dependent information graph has been considered for
identical and in\-dis\-tin\-guish\-able kinematically
controlled point-mass vehicles \cite{GYL10}. Each of them has access to the
positions of the target and two neighbors, which are the angularly
closest predecessor and follower in the circular order around the
moving target, irrespective of their metric distances from the
sensor. This does not match the capabilities of most sensors, for
which the visibility of an object may essentially depend on this
distance, whereas the angular discrepancy with respect to a third
party (target) position does not matter.
\par
All aforementioned works assume no limitations on both the control
range and the distance sensible by the vehicle, which can hardly be
qualified as realistic, and neglects the important issue of collision
avoidance. Realistic models of sensors were examined in
\cite{CeMaGaGi08}, where the visibility region of the vehicle is the
union of a disc sector and a full disc of a smaller radius, which are
associated with a narrow aperture long-distance sensor and an
omnidirectional short distance one. For a steady target
and a team of identical under-actuated unicycles with unlimited
control ranges and full observation, including the mutual orientation
of the vehicles, this approach proposes decentralized control laws
and shows that they transform the desired uniform capturing
configuration into a locally stable equilibrium of the closed-loop
system. However, no rigorous results on global stability or collision
avoidance are provided.
\par
The issue of collision avoidance has also been addressed
\cite{KaBeAb08,ShFiAn08,LYL10}. A
control strategy has been proposed which drives a team of fully actuated unconstrained
unicycles with all-to-all communication capability to a prescribed
distance to a steady target, while avoiding collisions with each
other \cite{KaBeAb08}. However uniform distribution of the vehicles around the target
is not ensured. A decentralized control rule was
offered and shown to drive kinematically controlled and labeled
point-mass unit-speed planar robots into an uniform circular formation
around a steady beacon while avoiding collisions with each other and
obstacles \cite{ShFiAn08}. However, this rule handles the team with only three agents,
which is a severe limitation. Another approach addresses both limited
control range and collision avoidance \cite{LYL10}. This considers a
team of identical and unlabeled unicycle-like planar vehicles and is
basically focused on the capturing maneuver within a vicinity of a
steady beacon, where the vehicles are assumed to acquire all-to-all
visibility. The vehicles are kinematically
controlled and fully actuated, with the ability to both stop,
immediately reverse the direction of the motion, and track an
arbitrarily contorted path. The proposed control strategy
asymptotically steers the vehicles at the desired distance to the
beacon and uniformly distributes them around it, with the speed of
formation rotation about the beacon being unprescribed. Collision
avoidance is guaranteed only if initially the vehicles lie on distinct
rays issued from the target.
\par
Thus the combination of realistic limited both control and sensor
ranges was not addressed in the literature; moreover for realistic
sensor models \cite{CeMaGaGi08}, no rigorous non-local convergence
results were offered. Such results taking into account the issue of
collision avoidance were established only for the team of three agents
\cite{ShFiAn08} or under the restricted assumption about the
all-to-all communication \cite{LYL10}. Furthermore, all
above references assume that both the line-of-sight angle (bearing)
and the relative distance (range) of visible objects are available to
the controller. However, the problem of range-only based
navigation is critical for many areas such as wireless networks,
unmanned vehicles and surveillance services \cite{ADB04,GS04}.  Many
sensors typical for these areas, like sonars or range-only radars,
provide only the relative distance between the sensor and the object. 
Range-only based navigation has a potential to reduce the
hardware complexity and cost. However this benefit is undermined by
the lack of suitable control design techniques.
\par
In this chapter, a decentralized control strategy is described which drives a
team of unlabeled and kinematically controlled planar non-holonomic
Dubins-like vehicles into an equi-spaced circular formation at a given
distance from a steady target. Unlike all papers in the area except
for \cite{LYL10}, the velocity range is limited, however unlike
\cite{LYL10}, the speed is not only upper but also lower bounded by a
given constant so that the vehicle cannot stop or immediately reverse
the direction of the motion. The speed lower bound provides an
additional challenge since it causes restriction of the paths along
which the vehicle can travel to the curves with upper limited
curvatures. Another distinction from \cite{LYL10} is that the vehicles
are not identical and the angular velocity of the formation rotation
about the target is pre-specified.
\par
The crucial difference of the approach described in this chapter from the entire previous research
in the area is that first, the proposed navigation algorithm is based
on range-only measurements; and second, the distance to a companion
vehicle is accessible only when it lies within a given disc sector
centered at the sensor (see Fig.~\ref{tcp:rob_sens}(a)), which would hold
for rigidly mounted narrow aperture distance sensors. It follows that
dangerous sideways convergence of two vehicles may be undetectable by
each of them, and this provides an additional challenge in design of the
control law -- it should exclude convergence in `invisible'
directions. It is assumed that the distance to the
target is constantly available to any vehicle, which would hold if
the target is endowed with a beacon facility or its visibility
essentially exceeds that of the vehicles.
\par
A sliding mode control strategy inspired by some ideas from
\cite{Teimoori2010journ1a} is proposed and shown to solve the
problem. It may be proven via rigorous mathematical analysis that under
some minor and partly unavoidable technical assumptions, the control
objective is achieved without fail. It is possible to show that despite of range-only
measurements, not only the required angular velocity but also the
direction of the formation rotation about the target are
ensured.
\par
The applicability of the proposed controller is confirmed by computer
simulations and experiments with real robots. In the extensive
literature on the multi-vehicle coordination, real-world tests have
been scarcely concerned up to now, especially with respect to the
particular problem treated in this chapter. Accordingly, there is still
considerable potential research that should be done on developing the
practical aspects of these implementations. Such tests are
particularly important since un-modeled system dynamics are a
classical reason for control malfunction.
\par
All proofs of mathematical statements are omitted here; they are
available in the original manuscript \cite{matveev2013capt}.
\par
The body of this chapter is organized as follows. In Sec.~\ref{tcp:sec.pst} 
the problem is formally defined, and in Sec.~\ref{tcp:sec.assr} 
the main analytical results are outlined. Simulations and
experiments are presented in Secs.~\ref{tcp:sec.simul} and \ref{tcp:sec.exp}. Finally, brief conclusions are given
in Sec.~\ref{tcp:sum}.

 \section{Problem Statement}
\label{tcp:sec.pst} 

In this chapter, a team of $N$ autonomous unicycle-like
robots are considered, enumerated by $i \in [0:N-1]$. They travel in the plane with
time-varying turning radii and longitudinal speeds limited by known
constants; the speed is bounded from not only above but also below by
positive constants. There also is a point-wise target or beacon $\bt$
in the plane. The objective is to drive all robots to the circle of
the pre-specified radius $R$ centered at the target, to uniformly
distribute them along this circle, and to maintain this formation
afterwards. The number $N$ is known to
every robot, and the formation should rotate with the prescribed
angular velocity $\omega$.
\par
\begin{figure}[ht]
  \centering \subfigure[]{
    \scalebox{0.3}{\includegraphics{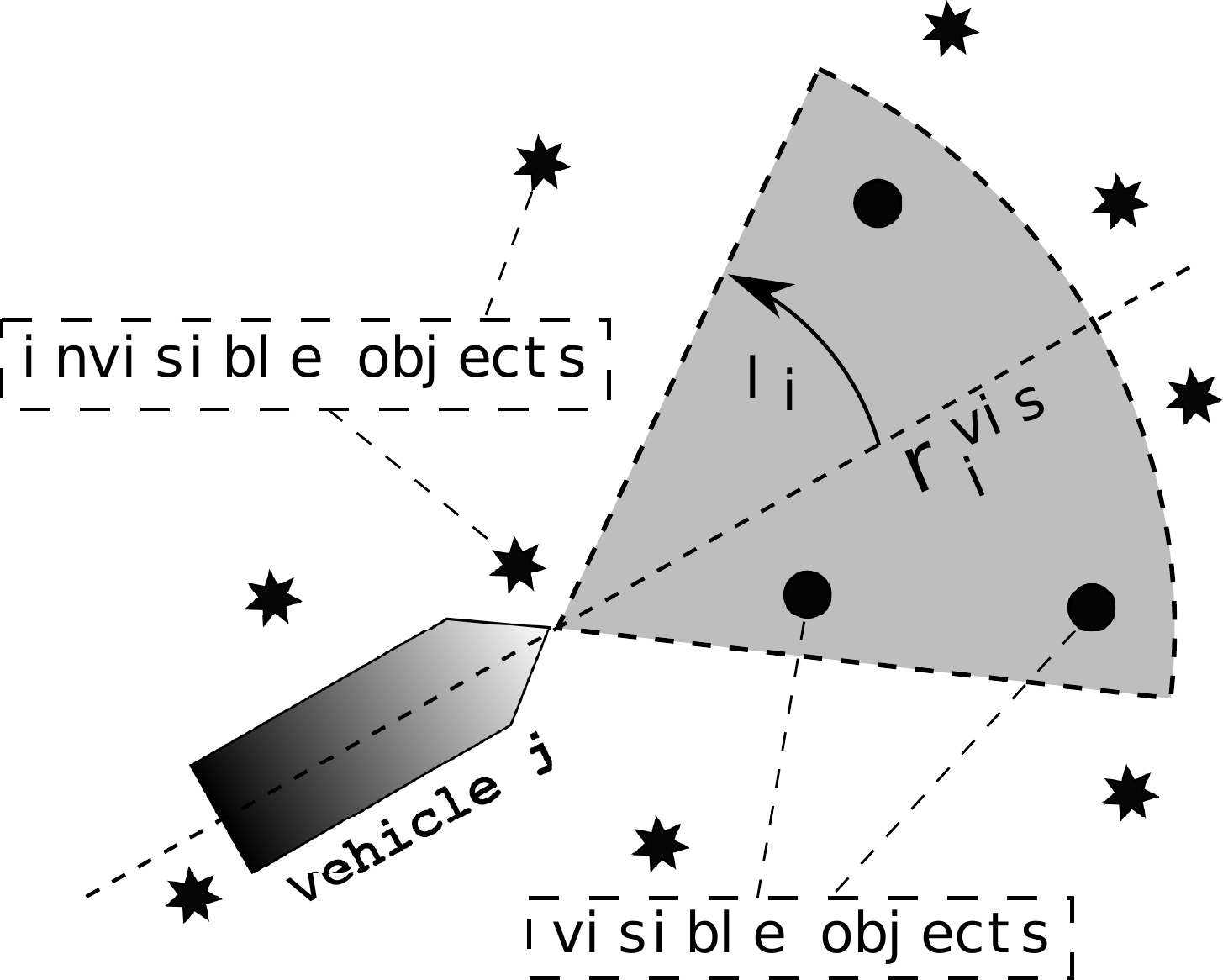}}}
  \hspace{20pt}
  \subfigure[]{\scalebox{0.3}{\includegraphics{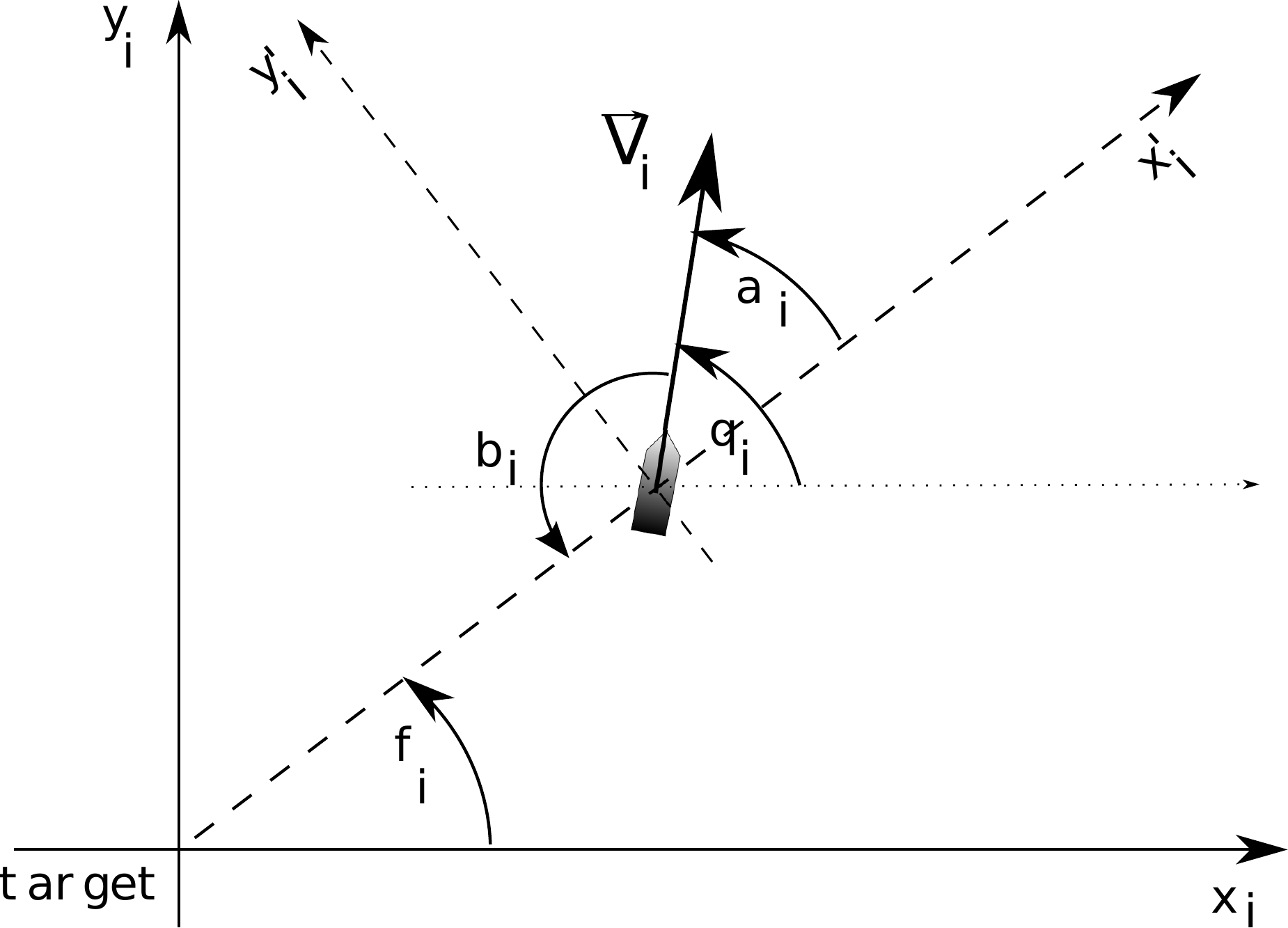}} }
  \subfigure[]{\scalebox{0.3}{\includegraphics{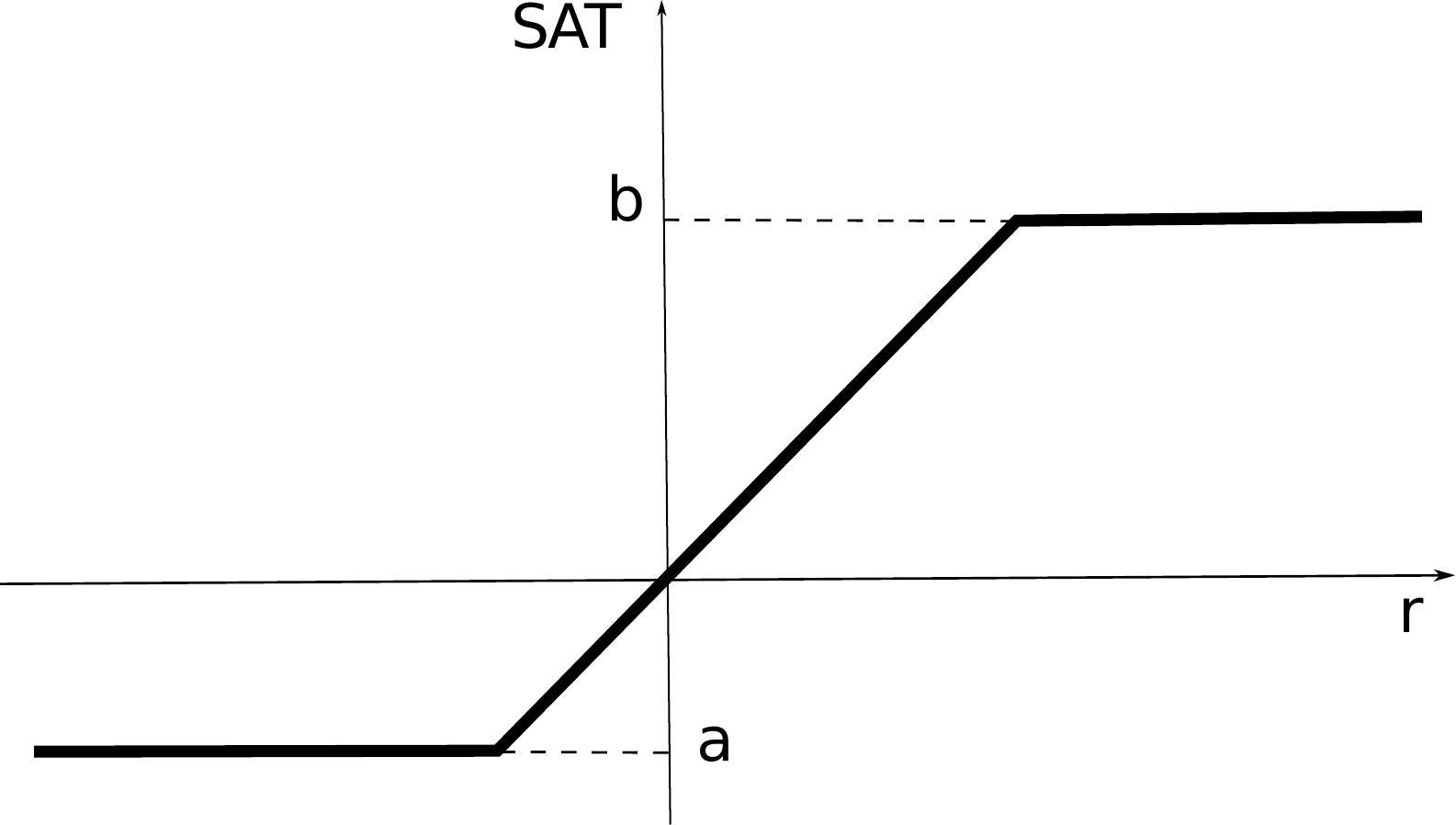}}}
  \caption{(a) The sensor capability of the robot; (b) Coordinate
    frames and variables; (c) The saturation function.}
  \label{tcp:rob_sens}
\end{figure}

To accomplish the mission, robot $i$ has access to the distance $d_{i
  \to \bt}$ from its center-point to the target. It also has access to
the distance $d_i$ to the nearest companion robot among those from its
visibility region. This is the sector between the arc and two radii of
the circle of the given radius $r^{\text{vis}}_i>0$ centered at the
robot, both radii are at the given angular distance $\lambda_i \in (0,
\pi/2)$ from the forward centerline ray of the robot; see
Fig.~\ref{tcp:rob_sens}(a). Whenever there are no robots in this region,
$d_i:= \infty$.

The following unicycle-like robot model is employed:

\begin{equation}
  \label{tcp:rob.eq}
  \begin{array}{r c l}
    \dot{x_i}& = & v_i \cos {\theta_i}
    \\
    \dot{y_i}& = & v_i \sin {\theta_i}
  \end{array}, \quad
  \begin{array}{l}
    \dot{\theta_i}=u_i,
    \\
    |u_i|\leq \overline{u}_i, \underline{v}\,_i \leq v_i\leq \overline{v}_i.
  \end{array}
\end{equation}

Here $x_i,y_i$ are the Cartesian coordinates of the $i$th robot in
the world frame centered at the target and $\theta_i$ gives the
orientation of the robot (see Fig.~\ref{tcp:rob_sens}(b)). The angular
$u_i$ and longitudinal velocities $v_i$ are the controls; the bounds
$\overline{u}_i >0, 0< \underline{v}\,_i < \overline{v}_i$ are
given. For the problem to be realistic, the required angular velocity
$\omega$ of the formation rotation should be such that
$\underline{v}\,_i \leq |\omega| R \leq \overline{v}_i, |\omega| \leq
\ov{u}_i$. This is slightly enhanced by assuming that:

\begin{equation}
  \label{tcp:bound.omega} \underline{v}\,_i < |\omega| R < \overline{v}_i,
  \qquad |\omega| < \ov{u}_i \qquad \forall i.
\end{equation}

To simplify the subsequent formulas, it is also assumed that the formation
should rotate about the target counter-clockwise $\omega >0$.
\par
In this chapter, the following simple switching
control strategy is employed:

\begin{equation}
  \label{tcp:control_unic} \left[
    \begin{array}{c}
      v_i
      \\
      u_i
    \end{array}
  \right] := \left[
    \begin{array}{c}
      \Psi_i(d_i)
      \\
      \ov{u}_i \sgn \left\{ \dot{d}_{i \to \bt} + v_i \sat^{b_i}_{-b_i}\left[\big(m_i [d_{i
            \to \bt} - R]\big)\right] \right\}
    \end{array}
  \right].
\end{equation}

Here the parameters $m_i,b_i >0$ and the functions $\Psi_i(d) \in
[\unnn{v}_i, \ov{v}_i]$ are chosen by the designer of the controller
subject to requirements to be disclosed further.

    The desired angular velocity $\omega$ of the formation will be
    taken into account by the requirements to $\Psi_i(\cdot)$. Note
    that the control law Eq.\eqref{tcp:control_unic} produces only feasible
    controls $v_i \in \big[ \underline{v}_i, \overline{v}_i \big],
    |u_i| \leq \ov{u}_i$.

     \section{Summary of Main Results}
    \label{tcp:sec.assr}

	In this section the control law
    Eq.\eqref{tcp:control_unic} is demonstrated to achieve the control
    objective. This holds if apart from Eq.\eqref{tcp:bound.omega}, another
    natural assumption is satisfied and the controller parameters are
    properly tuned:

    \begin{Assumption}
      \label{tcp:as_4} Under the uniform distribution over the required
      circle, every robot has at least one companion in its visibility
      region (see Fig.~{\rm \ref{tcp:fig.domains}(a)}):
      \begin{equation}
        \label{tcp:visibility.condr}
        \lambda_i > \Delta_\ast := \frac{\pi}{N}, \qquad r^{\text{vis}}_i > 2R \sin
        \Delta_\ast .
      \end{equation}
    \end{Assumption}

    To specify the choice of the speed gain functions $\Psi_i(\cdot)$
    in Eq.\eqref{tcp:control_unic}, the following is introduced:
	
    \begin{Definition}
      The function $g(\cdot)$ defined on $E \subset$ $[-\infty,
      +\infty]$ is said to be {\em uniformly H\"{o}lder continuous} is
      there exists $c \in [0,+\infty)$ and $\alpha \in (0,1]$ such
      that $|g(t^{\prime\prime}) - $ $g(t^\prime)|$ $\leq$ $
      c|t^{\prime\prime}-t^\prime|^\alpha$ for all $t^{\prime\prime},
      t^\prime \in E \cap \br$.
    \end{Definition}
	
    The choice of $\Psi_i(\cdot)$ proceeds from $\omega$ and two
    auxiliary velocity parameters $v_{i,\ast}, v^\ast_i$ chosen so
    that:
	
    \begin{equation}
      \label{tcp:range.v} \unnn{v}_i \leq v_{i,\ast} < \omega R < v^\ast_i <
      \min \big\{ \ov{u}_i R ; \ov{v}_i \big\}.
    \end{equation}
	
    Such choice is feasible thanks to Eq.\eqref{tcp:bound.omega}. As
    $v_{i,\ast}, v^\ast_i$ are picked, the choices of the controller
    parameters $\Psi_i(\cdot)$ and ($m_i, b_i$) for the speed and
    steering gains, respectively, become independent of each other.

    \begin{Requirement}
      \label{tcp:req.psi} The maps $\Psi_i(\cdot):[0,+\infty] \to
      [v_{i,\ast}, v^\ast_i]$ are uniformly H\"{o}lder continuous and
      such that
      \begin{enumerate}[i)]
      \item $\Psi_i(d)=R\omega \;\forall d \in [0, d^0];
        \displaystyle{\inf_{d \geq d^0+\varepsilon}} \Psi_i(d) > R
        \omega\;\forall \varepsilon >0$, where $d^0:=2 R \sin\frac
        {\pi}{N}$ is the distance between any two neighboring robots
        under their uniform distribution over the desired circle;
      \item whenever $i \neq j$, the equation $\Psi_i(d)= \Psi_j(d)$
        has no roots on $(d^0, 2R]$.
      \end{enumerate}
    \end{Requirement}

    An example of the functions satisfying this requirement is as
    follows:
	
    \begin{equation}
      \label{tcp:control1_angle} \Psi_i(d_i) := \sat^{v_i^\ast}_{v_{i,\ast}} \left[ R \omega
        + k_i \left[ d_i - 2 R \sin\frac {\pi}{N} \right]_+ \right].
    \end{equation}

    Here $[s]_+:=\max\{x,0\}$, whereas $k_i >0$ are tunable parameters
    and $k_i \neq k_j \; \forall i \neq j, \; k_i < k_j \Rightarrow
    v_i^\ast < v_j^\ast.$
    \par
    The proposed controller in fact employs the reduced speed range
    $[v_{i,\ast}, v^\ast_i] \subset [\unnn{v}_i, \ov{v}_i]$. As will be
    shown, the closer this range to $\omega R$, the larger the
    estimated convergence domain.
    \par
    The choice of the steering gain parameters is subjected to the
    following limitations:

    \begin{gather}
      \label{tcp:vst}
      0 < b_i < 1, \qquad 0 < m_i < \frac{\sqrt{1-b_i^2}}{b_i} \left[
        \frac{\ov{u}_i}{v^\ast_i} - \frac{1}{R}\right] \qquad i \in
      [0:N-1],
    \end{gather}

    Here the right-hand side of the last inequality is positive due
    to the last inequality from Eq.\eqref{tcp:range.v}. Such choice is always
    possible -- chosen $b_i$, all small enough $m_i$ are feasible.
    \par
    To state the main results of the chapter, the rotating
    normally oriented Cartesian frame RCF$_i$ is introduced, which is centered at the target
    with the abscissa axis directed towards robot $i$ (see
    Fig.~\ref{tcp:rob_sens}(b)). Let $\varphi_i$ stand for the angle from
    the abscissa axis of the world frame to that of RCF$_i$, with the
    anticlockwise angles being positive. The correct behaviour of the system is
	defined as follows:

    \begin{figure}[ht]
      \centering
      \subfigure[]{\scalebox{0.15}{\rotatebox{60}{\includegraphics{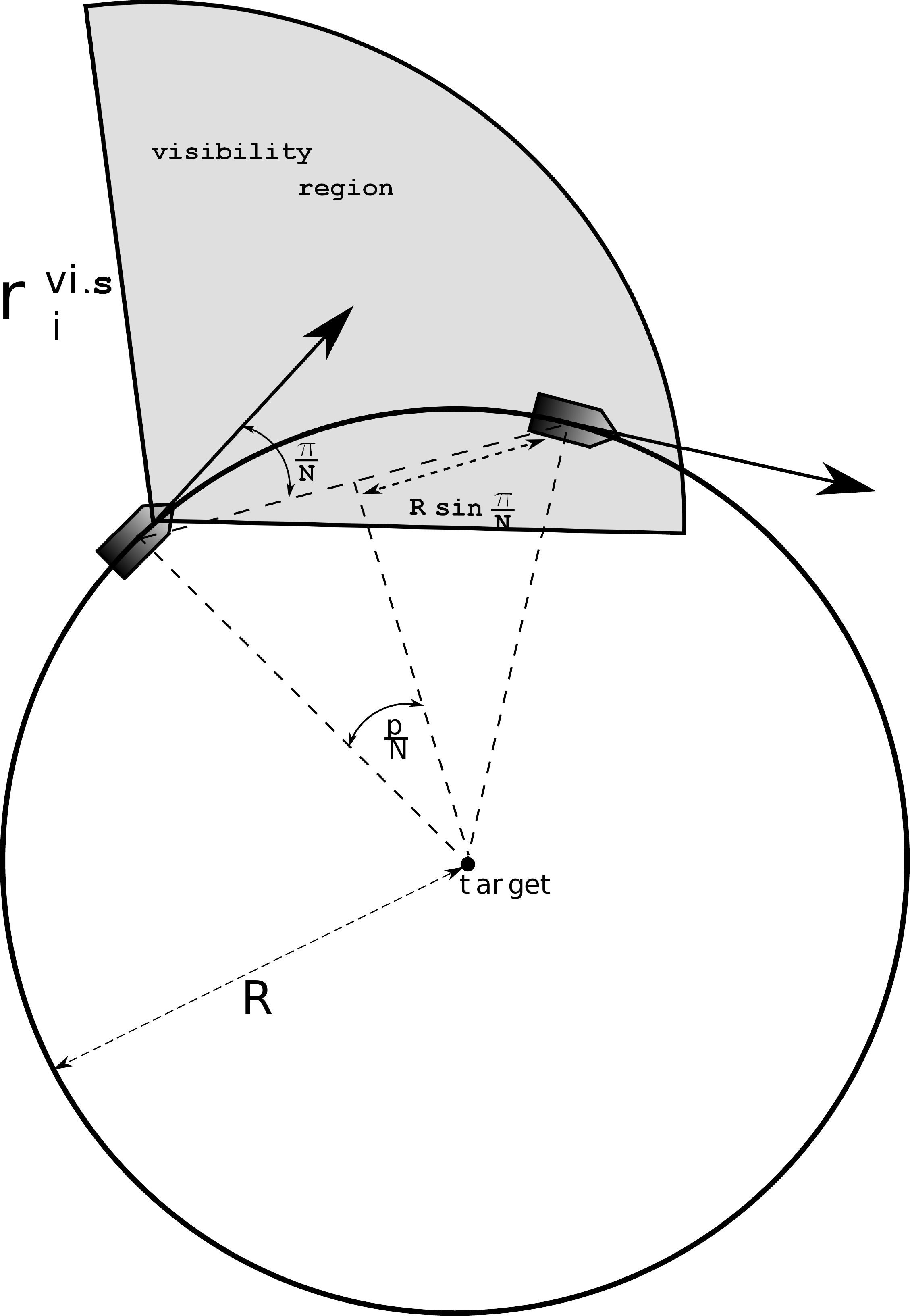}}}}
      \hspace{20pt}
      \subfigure[]{\scalebox{0.5}{\includegraphics{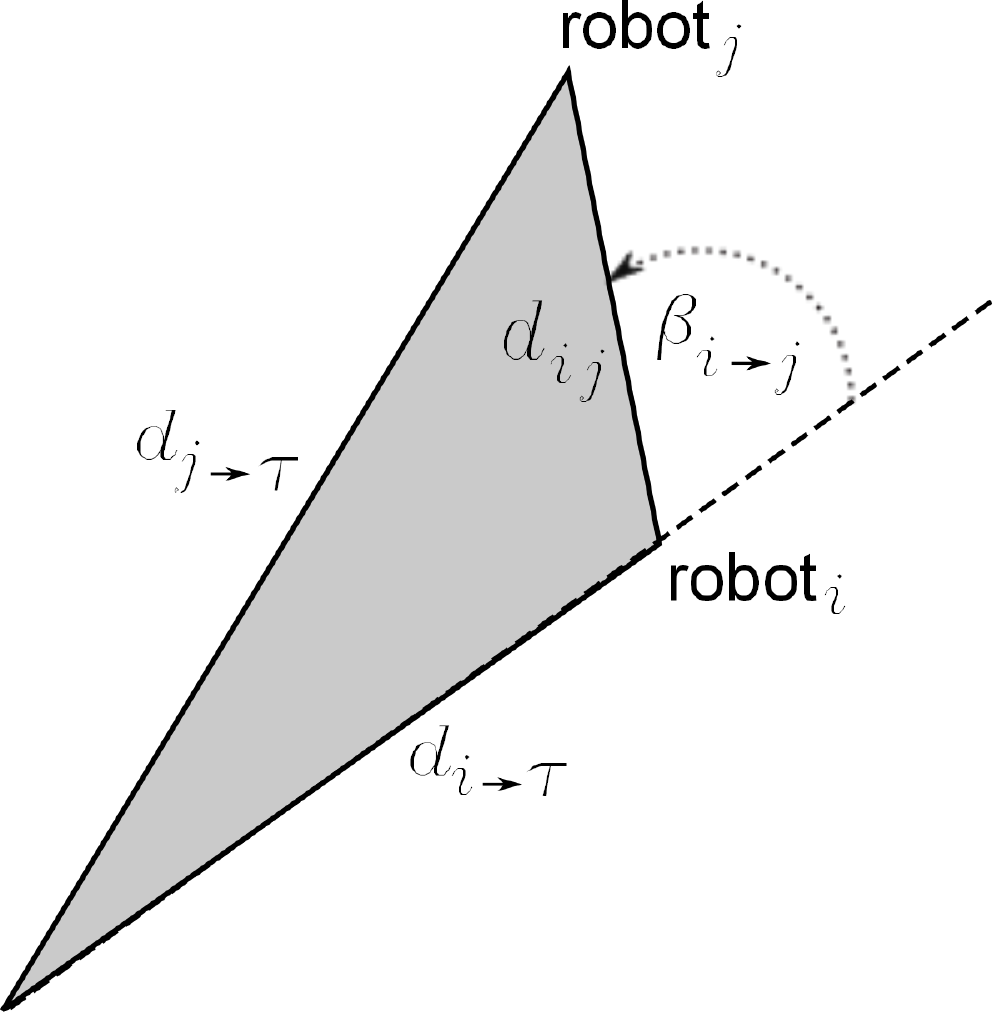}}}
      \caption{(a) Visibility region under uniform distribution; (b)
        Relative polar coordinates.} \label{tcp:fig.domains}
    \end{figure}

    \begin{Definition}
      The cumulative trajectory of the team of $N$ robots is said to
      be {\em target capturing} if
      \begin{equation}
        \label{tcp:objective}  d_{i \to \bt} (t) \to R,\qquad
        \dot{\varphi_{i}}(t) \to \omega \ \quad \text{\rm as} \quad t \to
        \infty \quad \forall i \in [0:N-1],
      \end{equation}
      vehicle-to-target collisions are excluded $d_{i \to \bt}(t) > 0
      \; \forall t,i$, and the robots can be enumerated so that
      \begin{equation}
        \label{tcp:phi.conv} \varphi_{i\oplus 1}(t) - \varphi_i(t) \to 2
        \Delta_\ast = \frac{2 \pi}{N}  \quad \text{\rm as} \quad t \to
        \infty.
      \end{equation}
    \end{Definition}

    Here and throughout $\oplus$ denotes the addition modulus $N$. It is
    underscored that this enumeration should be static -- it does not
    alter as time progresses.
    \par
    The following theorem is the first main result of the chapter:

    \begin{Theorem}
      \label{tcp:th.main} Suppose that Assumption~{\rm \ref{tcp:as_4}} and
      Eq.\eqref{tcp:bound.omega} hold and the controller parameters satisfy
      Requirement~{\rm \ref{tcp:req.psi}} and Eq.\eqref{tcp:vst}. Then the
      control law Eq.\eqref{tcp:control_unic} gives rise to a target
      capturing trajectory whenever the robots are initially far
      enough from the target:
      \begin{equation}
        \label{tcp:conv.domm}
        d_{i \to \bt}(0) > 4 r^\ast_i - 2 r_{i,\ast} +
        \frac{1}{\frac{1}{r^\ast_i} - \frac{m_ib_i}{\sqrt{1-b_i^2}}} \quad
        \forall i \in [0:N-1], \quad \text{\rm where} \quad r^\ast_i:=
        \frac{v^\ast_i}{\ov{u}_i}, r_{i,\ast}:= \frac{v_{i,\ast}}{\ov{u}_i}.
      \end{equation}
    \end{Theorem}

Theorem~\ref{tcp:th.main} ignores the issue of possible collisions between
the vehicles. It is tacitly taken for granted that the collisions are
resolved by an extra controller. Moreover, this controller does not
essentially corrupt the trajectories implemented under the proposed
control law. However, under a mild additional assumption, the
proposed control strategy automatically prevents inter-vehicle
distance converging to zero:

\begin{Proposition}
  \label{tcp:prop.main}
  Suppose that the assumptions of Theorem~{\rm \ref{tcp:th.main}} are
  true, the gain coefficients $m_i, b_i$ in Eq.\eqref{tcp:control_unic} are
  taken common for all robots, and the initial distance $d_{ij}(0)$
  between any pair of robots $i \neq j$ is large enough
  \begin{equation}
    \label{tcp:rob_x0}
    d_{ij}(0) >  3\pi(\ov{v}_i+\ov{v}_j) \max_{\nu=1,\ldots,N} \ov{u}_\nu^{\,-1}.
  \end{equation}
  Then the inter-vehicle collisions are excluded: $d_{ij}(t) >0 \;
  \forall t,i$.
\end{Proposition}

 \section{Simulations}
\label{tcp:sec.simul}

\begin{figure}[ht]
  \centering
  \subfigure[]{\includegraphics[width=0.45\columnwidth]{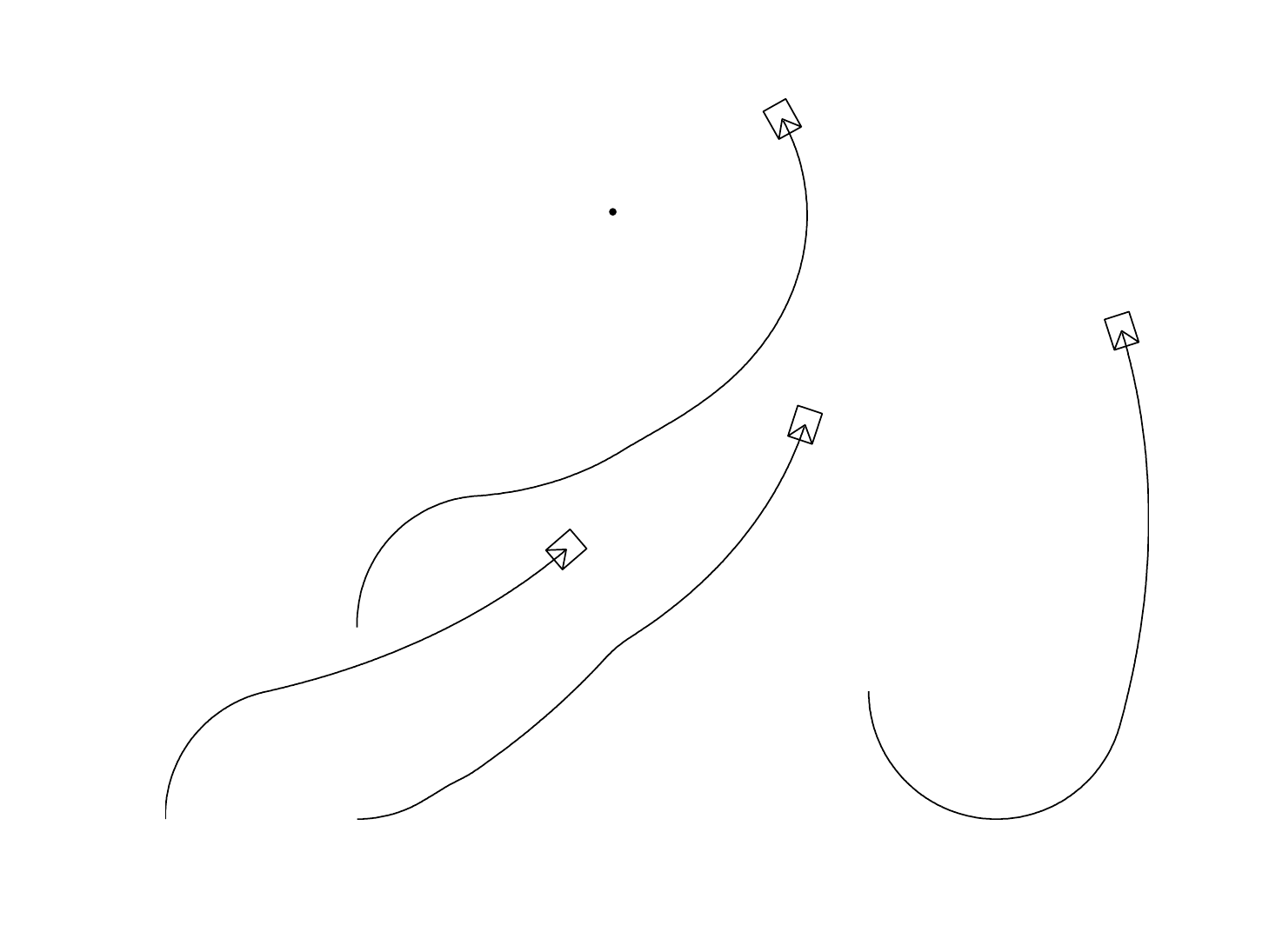}}
  \subfigure[]{\includegraphics[width=0.45\columnwidth]{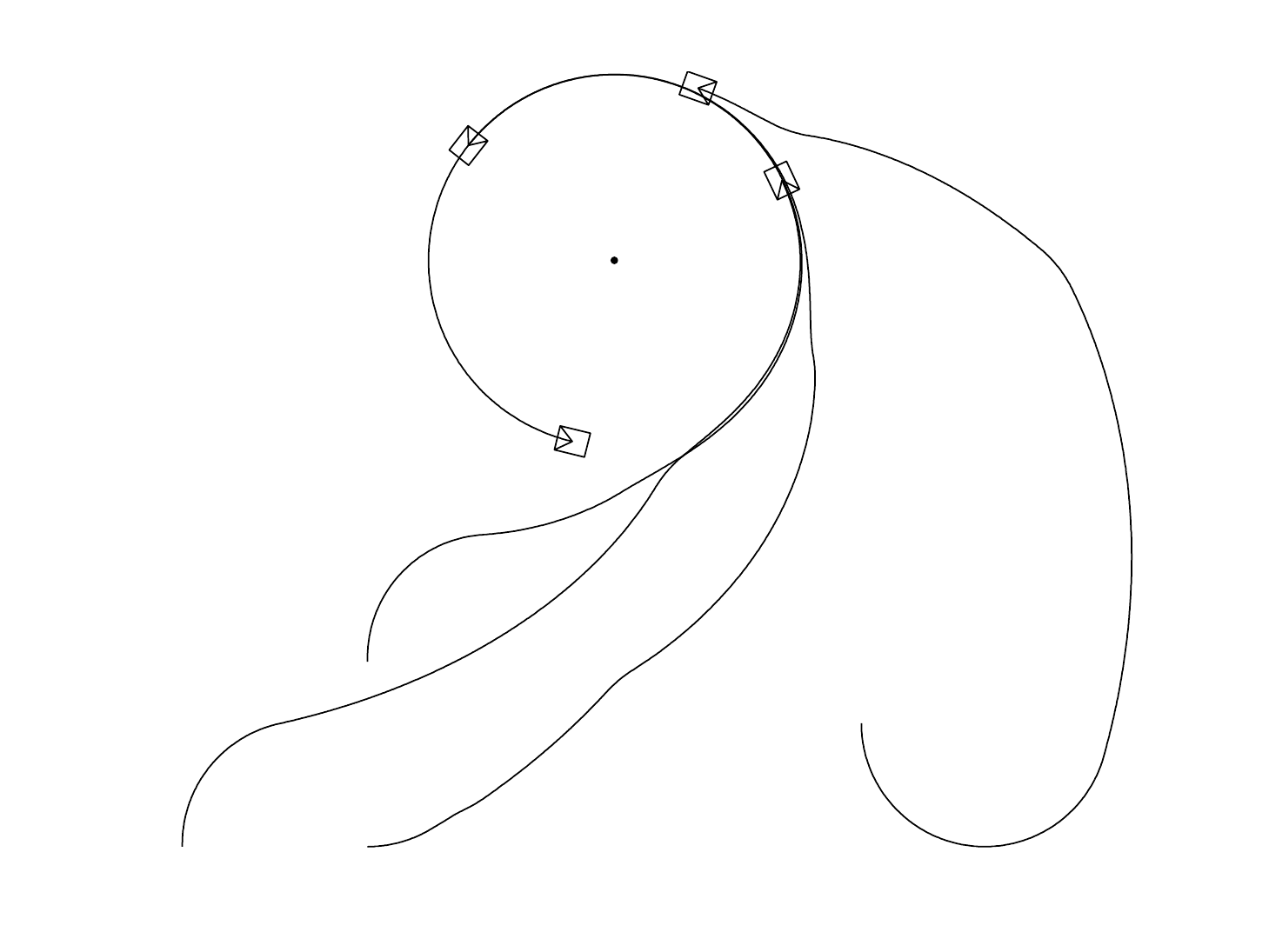}}
  \subfigure[]{\includegraphics[width=0.45\columnwidth]{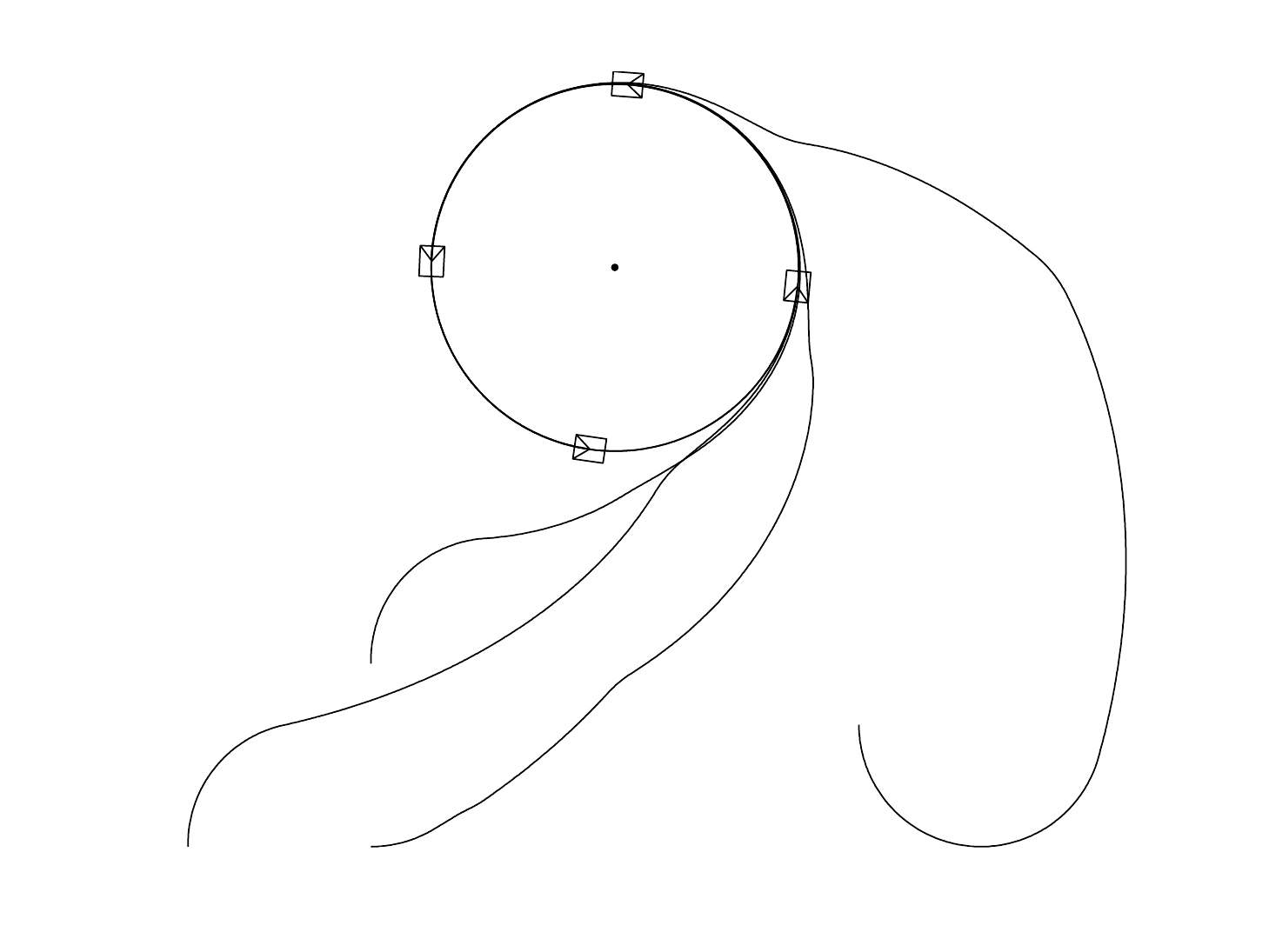}}
  \caption{Simulations with four vehicles converging to equal spacing around a target.}
  \label{tcp:fig:rig}
\end{figure}

To verify the correct operation of the proposed control law, the
system was simulated in a range of scenarios. The control was updated
with a sampling period of 0.01 s.
The vehicle and controller parameters are shown in
Table~\ref{tcp:fig:param}. The closed loop trajectories for each vehicle
are shown in Fig.~\ref{tcp:fig:rig}, which displays the expected
behavior of the vehicles: they converge to the desired circle around
the target, while equalizing the inter-vehicle distances.

\begin{table}[ht]
  \centering
  \begin{tabular}{| l | c |}
    \hline
    $\omega$ & $\frac{\pi}{8} rad s^{-1}$ \\
    \hline
    $R$ & $3 m$ \\
    \hline
    $\overline{u}_i$ & $1.5  rad s^{-1}$  \\
    \hline
  \end{tabular}
  \hspace{10pt}
  \begin{tabular}{| l | c |}
    \hline
    $k_i$ & $10 s^{-1}$ \\
    \hline
    $v_i^*$ & $3 ms^{-1}$  \\
    \hline
    $v_{i,*}$ & $2 ms^{-1}$  \\
    \hline
  \end{tabular}
  \hspace{10pt}
  \begin{tabular}{| l | c |}
    \hline
    $m_i$ & $2 m^{-1}$ \\
    \hline
    $b_i$ & $0.5$  \\
    \hline
  \end{tabular}
  \hspace{10pt}
  \begin{tabular}{| l | c |}
    \hline
    $\lambda_i$ & $\frac{\pi}{2} rad$  \\
    \hline
    $r_i^{vis}$ & $10 m$  \\
    \hline

  \end{tabular}
  \caption{Simulation parameters for target-capturing controller.}
  \label{tcp:fig:param}
\end{table}

\begin{figure}[ht]
  \centering
  \includegraphics[width=10cm]{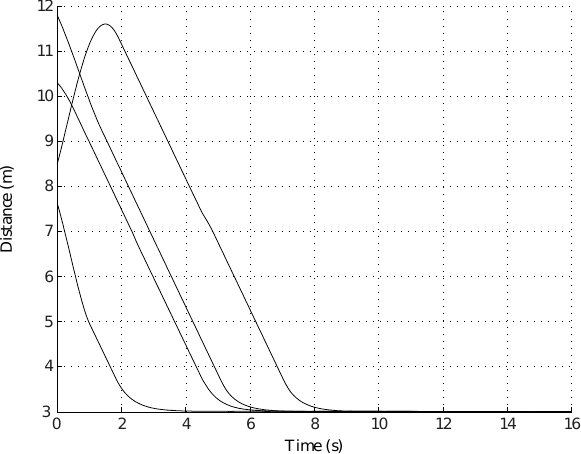}
  \caption{Distance to the target for each of the four vehicles in
    Fig.~\ref{tcp:fig:rig}.}
  \label{tcp:fig:dis}
\end{figure}

\begin{figure}[ht]
  \centering
  \includegraphics[width=10cm]{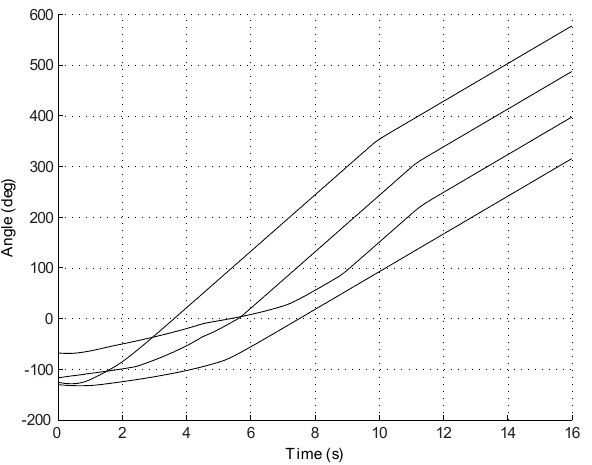}
  \caption{Angle subtended from the target for each of the four
    vehicles in Fig.~\ref{tcp:fig:rig}.}
  \label{tcp:fig:phi}
\end{figure}

This is confirmed by Figures \ref{tcp:fig:dis} and \ref{tcp:fig:phi}. The
first of them shows that the distance from each vehicle to the target
converges to the prescribed value $3 m$. The initial orientation of
the vehicle determines whether this convergence is monotonic or not -
the vehicle initially oriented away from the target unexpectedly
experiences a temporary increase in the target distance, while the
others do not. Fig.~\ref{tcp:fig:phi} shows that eventually the angle
subtended around the target
increases at nearly the desired rate $\omega=\frac{\pi}{8} rad
s^{-1}$, with nearly equal spacing between vehicles.

\begin{figure}[ht]
  \centering
  \subfigure[]{\includegraphics[width=0.6\columnwidth]{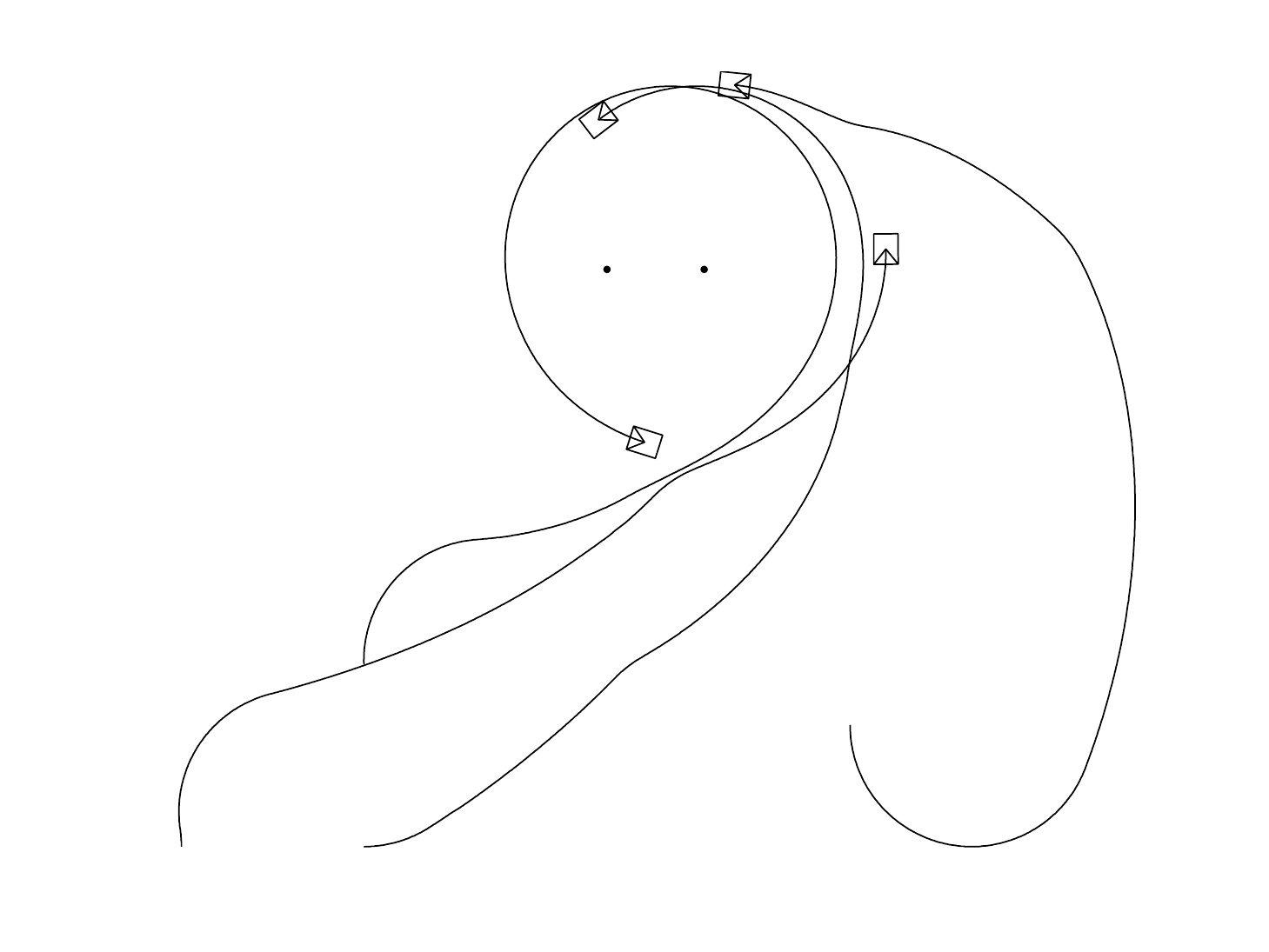}}
  \subfigure[]{\includegraphics[width=0.6\columnwidth]{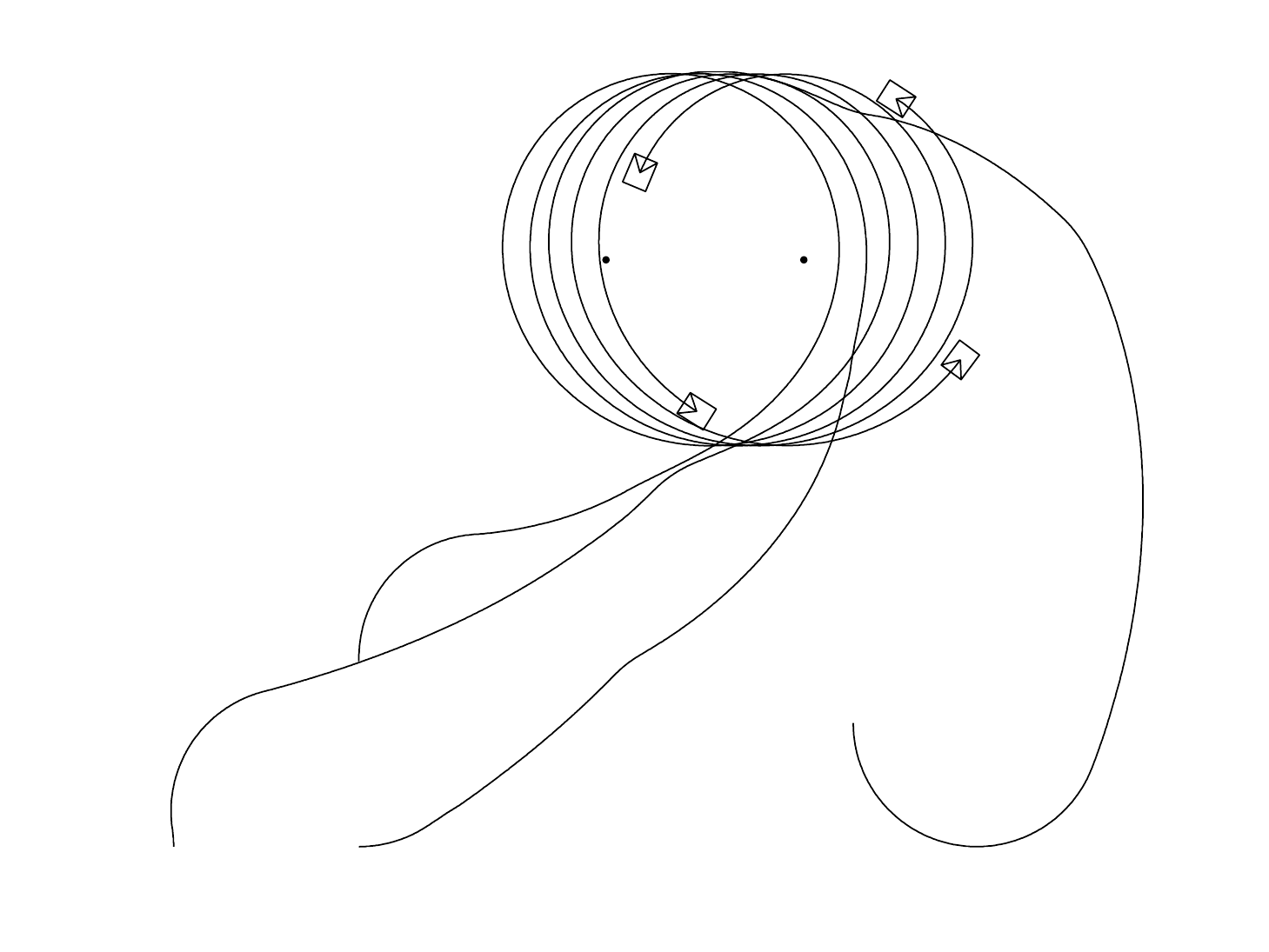}}
  \caption{Simulations with enclosing and grasping of a moving target.}
  \label{tcp:fig:mov}
\end{figure}

The system was also tested in a situation where the target
moves. Though this case is not covered by Theorem~\ref{tcp:th.main},
Fig.~\ref{tcp:fig:mov} demonstrates that the control law still ensures
the desired behavior. Since the target is slowly moving to the left,
the final paths of the vehicles unexpectedly resemble cycloids.

The choice Eq.\eqref{tcp:control1_angle} of the function $\Psi_i(d_i)$ in the
control law Eq.\eqref{tcp:control_unic} assumes that the team size $N$ is
known. If this size is unknown, the control law can be modified so
that this information is not required:

\begin{equation}
  \label{tcp:law.flexible}
  \Psi_i(d_i) := \sat^{v_i^\ast}_{v_{i,\ast}} \left[  k d_i \right]
\end{equation}

\begin{figure}[ht]
  \centering
  \subfigure[]{\includegraphics[width=0.6\columnwidth]{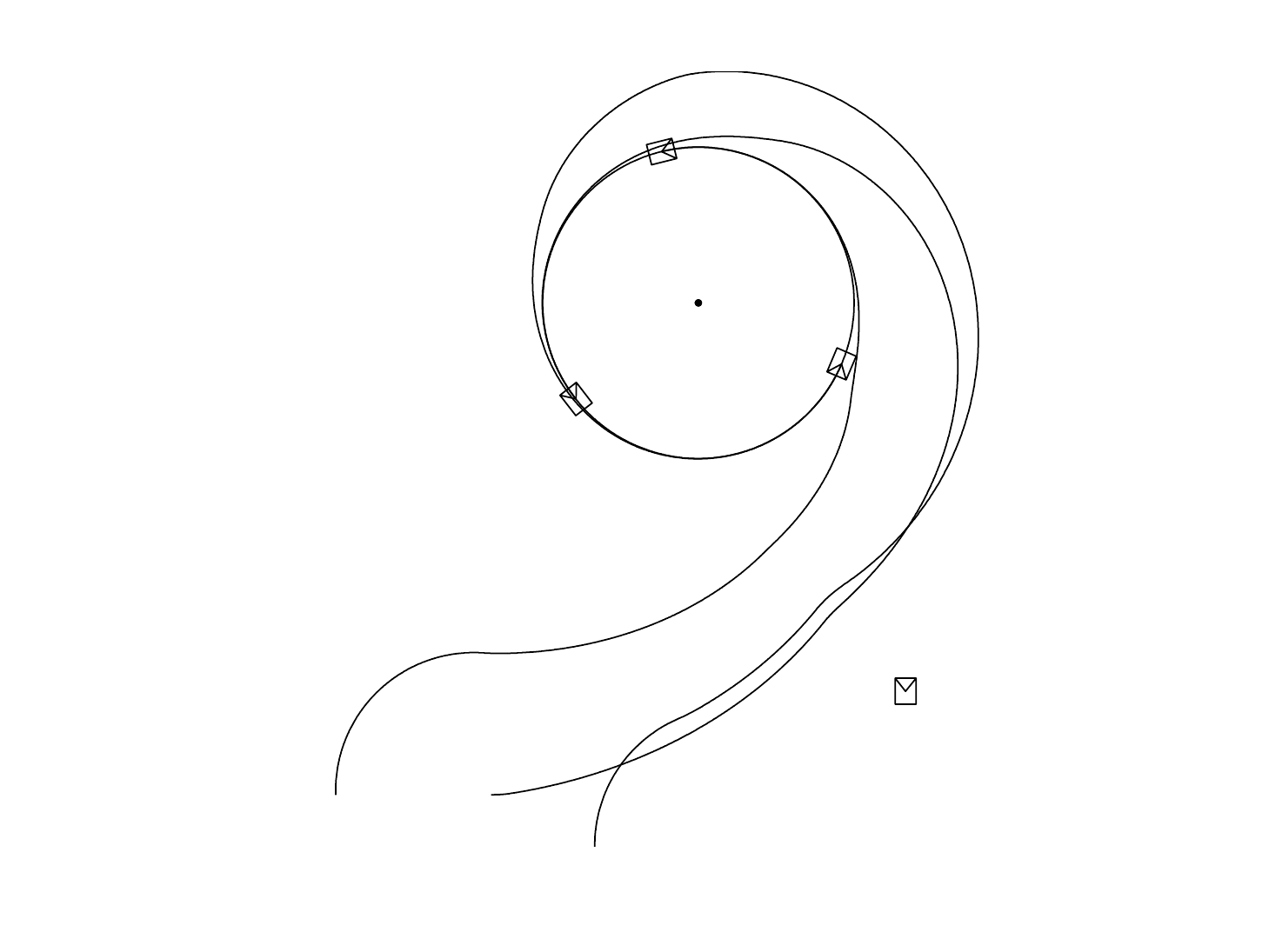}}
  \subfigure[]{\includegraphics[width=0.6\columnwidth]{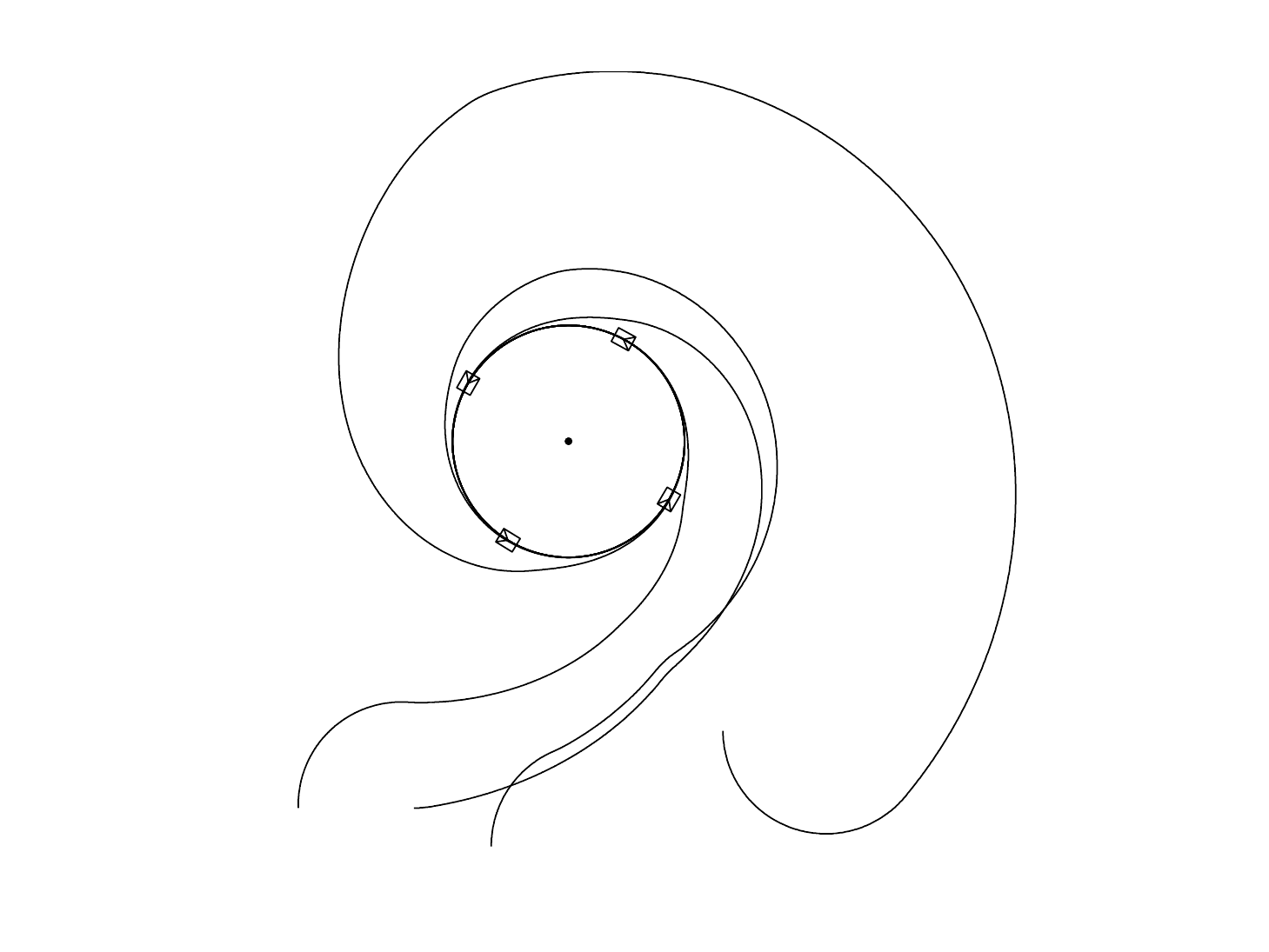}}
  \caption{Simulations with addition of vehicles during the target capturing maneuver.}
  \label{tcp:fig:flex}
\end{figure}

Contrary to Eq.\eqref{tcp:control1_angle}, the gain $k$ is common for all
vehicles here. Theoretical analysis of this control law, including
recommendations on the choice of the controller parameters, is outside
of the scope of this discussion and is a topic of ongoing research. For this
preliminary simulation testing of this law, the
parameters used in the test are given in
Table~\ref{tcp:fig:paramflex}. Simulations showed that the system still
behaves in a similar manner as before: the vehicles converge to the
desired circle centered at the target, while equalizing the
inter-vehicle spacing. A typical simulation result is shown in Figure
\ref{tcp:fig:flex}, where initially three vehicles converge to an uniform
formation, which is proceeded by an extra vehicle joining the team, causing
rearrangement of the uniform formation.

\begin{table}[ht]
  \centering
  \begin{tabular}{| l | c |}
    \hline
    $k$ & $0.25 s$ \\
    \hline
    $v_i^*$ & $4 ms^{-1}$  \\
    \hline
    $v_{i,*}$ & $1 ms^{-1}$  \\
    \hline
  \end{tabular}
  \hspace{10pt}
  \begin{tabular}{| l | c |}
    \hline
    $m_i$ & $1 m^{-1}$ \\
    \hline
    $b_i$ & $0.125$  \\
    \hline
  \end{tabular}
  \caption{Simulation parameters for target-capturing with unknown team size.}
  \label{tcp:fig:paramflex}
\end{table}

\clearpage \section{Experiments}
\label{tcp:sec.exp}

Experiments were carried out with three pioneer P3-DX mobile robots to
show real-time applicability of the proposed control system. A SICK
LMS-200 laser range-finding device was used to detect the vehicles and
target. This device has a nominal accuracy of $15mm$ along each
detection ray, however the accuracy may be significantly worse in real
world circumstances.
The range sensor was rotated 45 deg to the left to provide a better
view of the target. The algorithm for determining the position of the
target and companion vehicles is as follows:

\begin{itemize}
\item The detected points were segmented into clusters by identifying
  consecutive detection points with an Euclidean separation smaller
  than $0.2 m$. The resultant sequence of points was filtered based on
  the distance between the edges of the cluster.  In the experiment,
  the cutoff was set so that the cluster diameter was between $0.1 m$
  and $0.4 m$.
\item The distance from the vehicle to the cluster was taken to be the
  distance along the detection ray bisecting the edges of the
  cluster. The clusters were filtered based on this distance and a
  cutoff was set so that this measurement was under $3 m$ in the
  experiment. The closest cluster was taken to be the target, and the
  next closest cluster to the right of the target was taken to be the
  next vehicle\footnote{This heuristic was sufficient to emulate the
    examined control law in the particular scenario involved in this
    experiment, where the vehicles starting locations and orientations
    played an essential role (the vehicles were initially arranged
    approximately evenly around the target with moderate offsets from
    the equilibrium positions and so that during convergence the
    nearest object was always the target). However, this heuristic may
    be insufficient for other scenarios.
  }.
\end{itemize}

At each control update, the speed and turning rate of the low level
wheel controllers were set in accordance with the output of the
navigation algorithm. The control was updated with a period of $0.2
s$. The values of the parameters used in the experiments are listed in
Fig. \ref{tcp:fig:paramtest}.

\begin{figure}[ht]
    \centering
    \includegraphics[width=7cm]{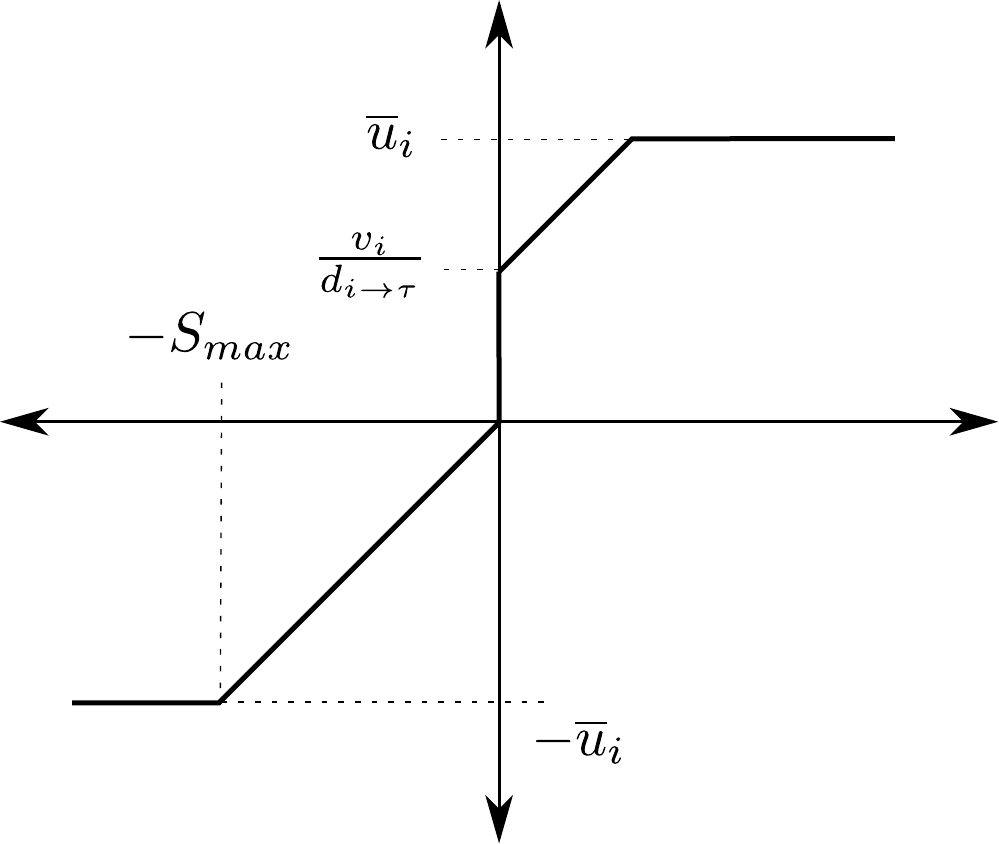}
    \caption{Saturation function used to reduce chattering during the
      experiment.}
    \label{tcp:fig:sat}
  \end{figure}

It is common to implement chattering reduction in sliding mode control
systems by smooth approximation of the signum function, and this can be
achieved using a linear function with saturation. However this may cause
static error of the turning rate, thus entailing a tracking error. As
an alternative, the steering control (commissioned to navigate the
robot to the desired distance to the target point) was approximated by
replacing $\overline{u}_i\sgn$ in the second line from
Eq.\eqref{tcp:control_unic} by the function from Figure~\ref{tcp:fig:sat}, where
$S_{\max}$ is an upper estimate of the absolute value of the function
in the curly brackets from Eq.\eqref{tcp:control_unic} in the domain where
this function is negative and the concerned $\sat$ is not
saturated. In this experiment,  $S_{\max}$ was taken to be $1.0 m
s^{-1}$. According to Figure~\ref{tcp:fig:sat}, this uniquely determines
the slope of the function on the negative ray of the abscissa axis; on
the positive ray, the slope is the same. Such choice does not
influence the equilibrium turning rate, thus eliminating the above
static error issue.

  \begin{table}[ht]
\centering
    \begin{tabular}{| l | c |}
      \hline
      $N$ & 3
      \\
      \hline
      $\overline{u}_i$ & $1.5  rad s^{-1}$  \\
      \hline
      $\lambda_i$ & $\frac{\pi}{2}$  \\
      \hline
    \end{tabular}
    \hspace{10pt}
    \begin{tabular}{| l | c |}
      \hline
      $r_i^{vis}$ & $10 m$  \\
      \hline
      $v_i^*$ & $0.1 ms^{-1}$  \\
      \hline
      $v_{i,*}$ & $0.2 ms^{-1}$  \\
      \hline
    \end{tabular}
    \hspace{10pt}
    \begin{tabular}{| l | c |}
      \hline
      $\omega$ & $\frac{\pi}{25} rad s^{-1}$ \\
      \hline
      $R$ & $1 m$ \\
      \hline
    \end{tabular}
    \hspace{10pt}
    \begin{tabular}{| l | c |}
      \hline
      $k_i$ & $0.06 s^{-1}$ \\
      \hline
      $m_i$ & $0.2 m^{-1}$ \\
      \hline
      $b_i$ & $0.05$  \\
      \hline
    \end{tabular}
    \caption{Experimental parameters for target-capturing controller.}
    \label{tcp:fig:paramtest}
  \end{table}

 \begin{figure}[ht]
   \centering
     \subfigure[]{\includegraphics[width=0.4\columnwidth]{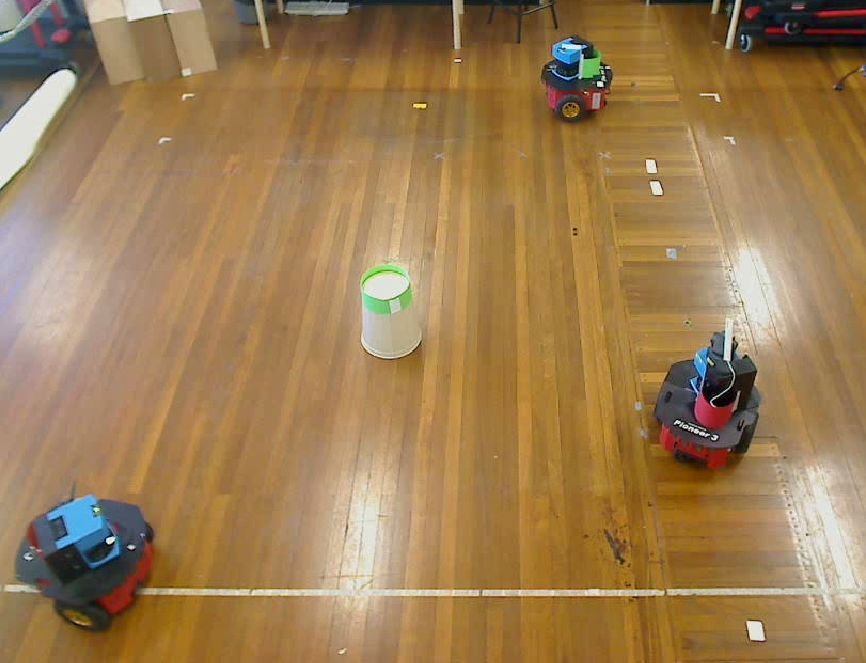}}
	 \subfigure[]{\includegraphics[width=0.4\columnwidth]{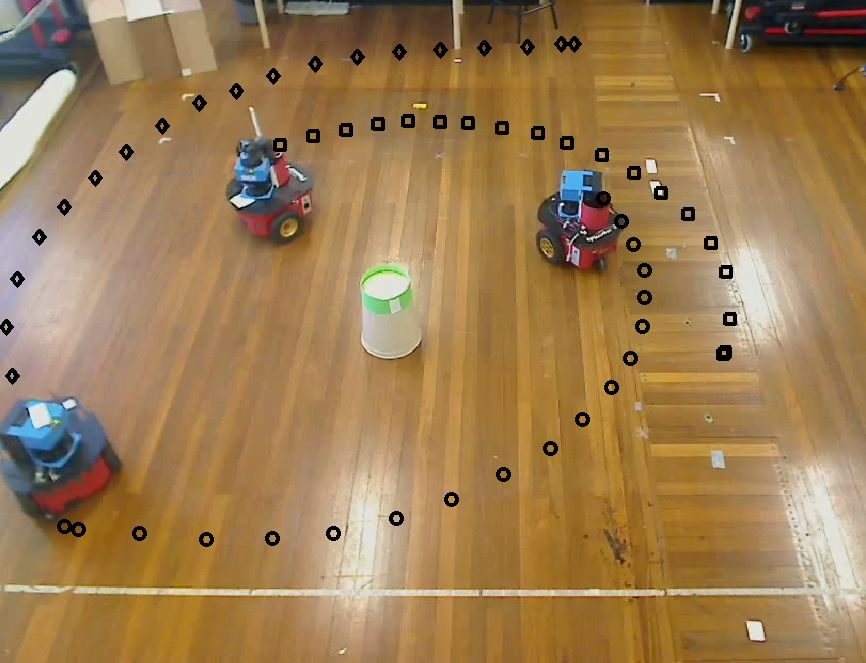}}
     \subfigure[]{\includegraphics[width=0.4\columnwidth]{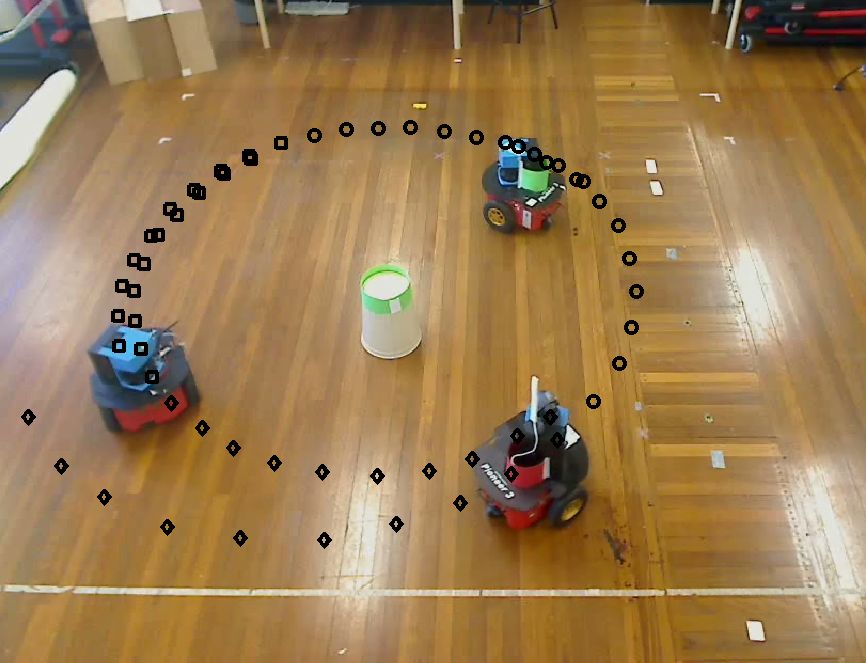}}
	 \subfigure[]{\includegraphics[width=0.4\columnwidth]{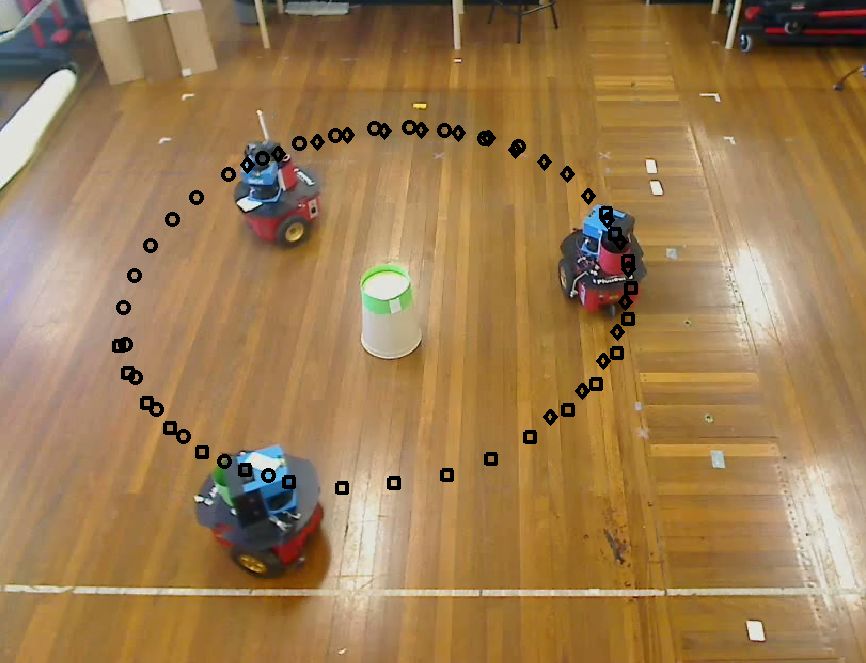}}
     \subfigure[]{\includegraphics[width=0.4\columnwidth]{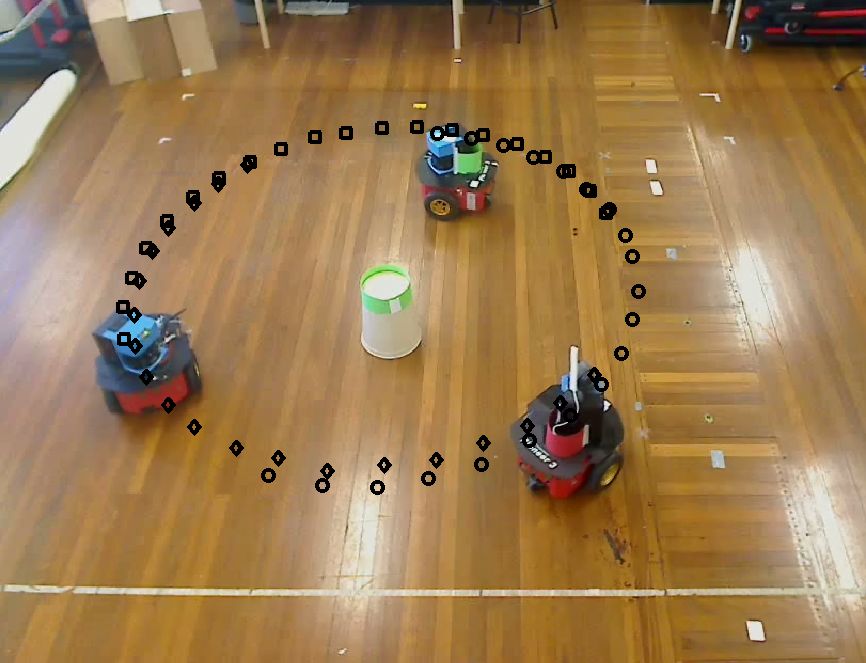}} 
	 \subfigure[]{\includegraphics[width=0.4\columnwidth]{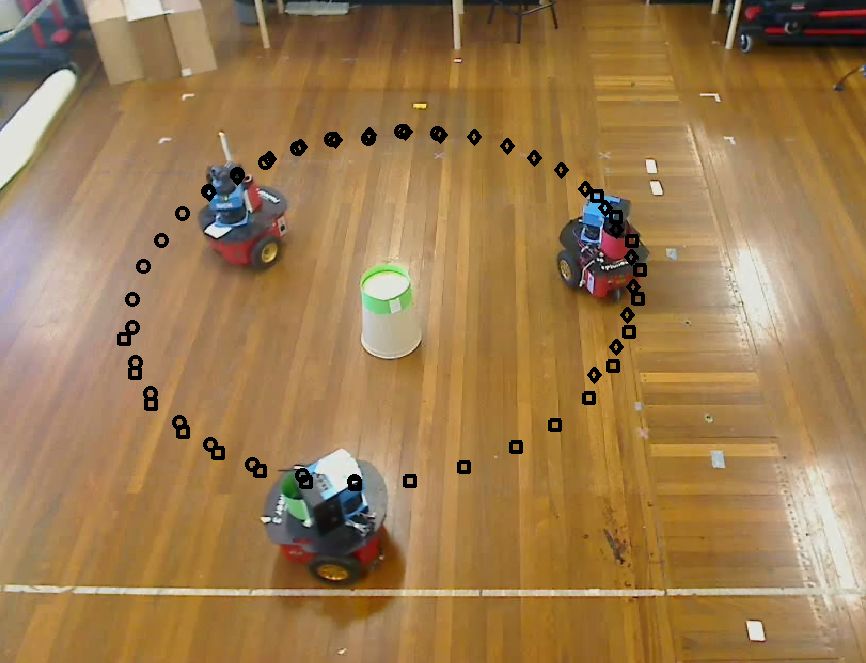}}
   \caption{Sequence of images showing the experiment.}
   \label{tcp:fig:pics}
 \end{figure}

 \begin{figure}[ht]
   \centering
   \includegraphics[width=10cm]{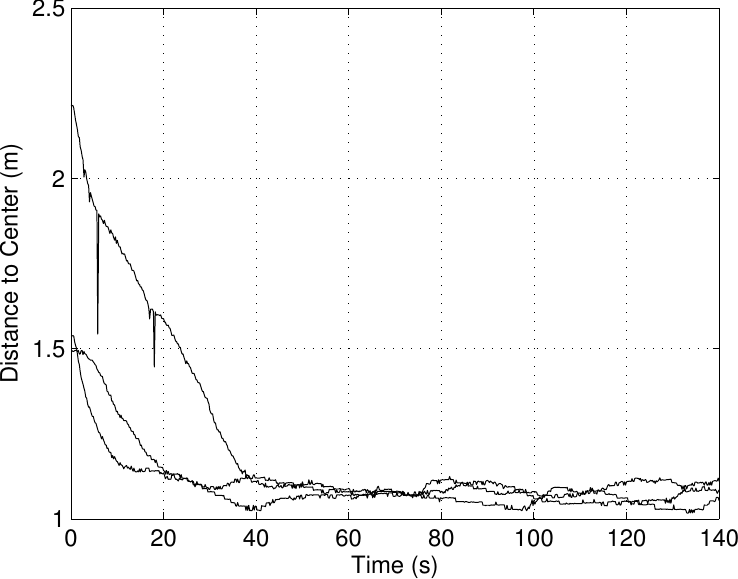}
   \caption{Distance to the center during the
     experiment.}
   \label{tcp:fig:dcen}
 \end{figure}

 \begin{figure}[ht]
   \centering
   \includegraphics[width=10cm]{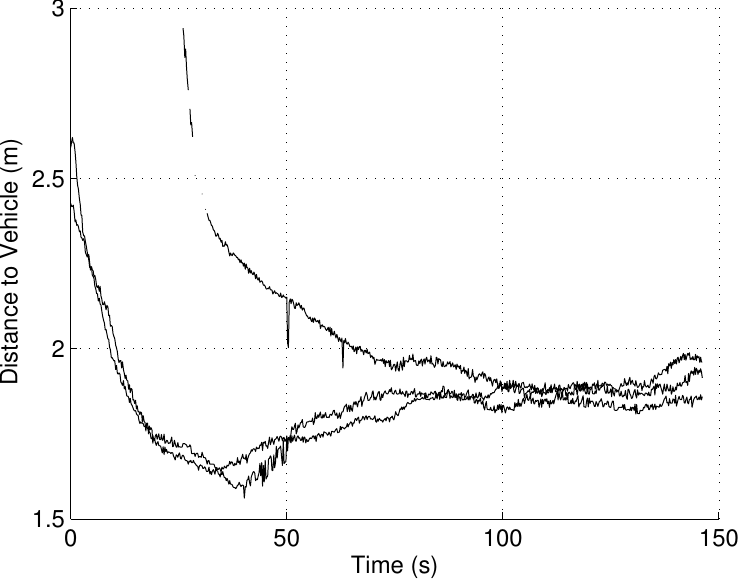}
   \caption{Distance to the nearest visible robot.}
   \label{tcp:fig:drob}
 \end{figure}

 A sequence of images obtained during the experiment is shown in
 Fig.~\ref{tcp:fig:pics}, where the positions of the robots were
 manually marked in the video frames. Figures \ref{tcp:fig:dcen} and
 \ref{tcp:fig:drob} demonstrate convergence to the required formation and
 confirm real-life applicability of the proposed control law.

 \clearpage \section{Summary}
\label{tcp:sum}

 In this chapter the problem of capturing a target using a team of
 decentralized non-holonomic vehicles was investigated, where the
 objective is to drive all vehicles to the circle of the prescribed
 radius centered at the target with uniform separation. If every
 vehicle has access to the distance to the target and the distances to
 the companions from the given disc sector centered at this vehicle,
 it may be shown that the objective is achieved. The performance of the
 control law is illustrated by computer simulations and experiments
 with real robots.

\chapter{Conclusions}
\label{chap:conclusions}

This report is concerned with provable collision avoidance of multiple autonomous
vehicles operating in unknown environments. The main contributions of this work
are a novel MPC-based strategy applicable to multiple vehicle systems
(see Chapt.~\ref{chap:multiple}), along with an approach
 to navigation problems where information about the
obstacle is described by a ray-based sensor model (see
Chapt.~\ref{chap:convsingle}).
\par
The decentralized coordination of multiple vehicle systems is a relatively
recent area of research, and the approach described in this report provides
some advantages over other proposed solutions to this problem. In particular it allows
vehicles to simultaneously update trajectories, use a sign-board based
communication model, and avoid artificial limitations to the magnitude
of trajectory alterations. Through it currently only works with simplified planning algorithms of
the type proposed in Chapt.~\ref{chap:singlevehicle}, it seems to be a useful
approach for practical implementations in vehicles.
\par
Concurrent use of both a ray-based sensor model, and a realistic vehicle kinematic
model to achieve boundary following has been previously proposed, but is limited to the use of a
single detection ray.
The application of MPC to boundary following problems proposed in this report is
original, and has many desirable properties such as adaptable speed and offset
distance. Additionally, while assumptions about the obstacle are unavoidable, these are 
likely more general than previously proposed approaches.
It seems that combined sensor and MPC-based navigation is an approach which could be
useful for many applications in the future.
\par
An additional contribution of this work is a proposed extension to allow
deadlock avoidance in Chapt.~\ref{chap:dead}. While only applicable to only two vehicles, the proposed idea
 is novel since all previously proposed deadlock
avoidance systems either rely on intractable centralized coordination,
discrete graph-based state abstractions, or an obstacle free workspace. 
Showing local deadlock avoidance for even two vehicles seems to be a useful property for a navigation law.
\par
The
local reactive navigation strategy proposed in Chapt.~\ref{chap:singlevehicle} 
is based on a combination of simplified
planning algorithms and robust MPC, and some of the details involved make it original.
The method is applicable to both holonomic and
unicycle vehicle models with bounded acceleration and disturbance.
Many simulations and real world tests throughout this report confirm the
viability of the proposed methods.
While Chapts.~\ref{chap:singlevehicle} to \ref{chap:dead} are all based on the 
same basic control structure, solving a generalised problem with multiple vehicles
and limited sensor data was not attempted, and this would be the logical next step to
extend this work.
\par
In most cases analytical justification for showing the correct behavior of the
proposed methods was offered. This type of justification is becoming increasing
important in robotics research. Properties such as provable collision avoidance
under disturbance are generally only easily provable using robust MPC, however
these often require significant computation and communication capabilities, making real
time use more difficult. This is especially apparent in the development of micro UAV vehicles,
which have limited resources available. It is hoped the approaches described in this report
provide a suitable trade-off between tractability, optimality and
robustness.
\par
In addition to collision avoidance, task achievement is a desirable attribute
to prove. In Chapt.~\ref{chap:convsingle} complete transversal of
the obstacle and finite completion time is shown. If high level navigation system functions 
can be delegated to low level controllers without significant complication,
it may increase the overall durability of the system (it seems likely that high level 
decision making system are slower and more likely to fail in the real world). 
\par
Lastly, Chapts.~\ref{chapt:rbf} to \ref{chapt:tcp}, provide simulations 
and real world testing performed to validate a variety of other navigation 
systems. The contributions are listed as follows:

\begin{itemize}
\item Chapt.~\ref{chapt:rbf} solves the problem of reactively avoiding obstacles 
and provably converging to a target using very limited scalar measurements. The advantage of the proposed method is that it can be analytically shown to have the correct behavior, despite the extremely limited sensor information assumed to be available.

\item Chapt.~\ref{chapt:tf} provides a method for reactively avoiding obstacles 
and provably converging to a target using a tangent sensor. The advantage of the proposed method is that it explicitly allows for the kinematics of the vehicle.

\item Chapt.~\ref{chapt:pf} achieves a method for path following with side slip allowance suitable for an agricultural vehiclele. The advantage of the proposed method is that it explicitly allows for steering angle limits, and it has been shown to have good comparative performance in certain situations.

\item Chapt.~\ref{chapt:fbf} contributes a method which allows the boundary of an obstacle to be followed using only a rigid range sensor. The advantage of the proposed method is that it provides a single contiguous controller, and is analytically correct at transitions from concave to convex boundary segments.

\item Chapt.~\ref{chapt:ext} provides a novel method for seeking the maximal point of a scalar environmental field. The advantage of the proposed method is that it des not require any type of derivative estimation, and it may be analytically proven to be correct in the case of time--varying environmental fields.

\item Chapt.~\ref{chapt:lst} achieves a new approach of
tracking level sets of an environmental field. The advantage of the proposed method is that it may be analytically proven to be correct in the case of time--varying environmental fields.

\item Chapt.~\ref{chapt:tcp} achieves a novel method of
 decentralized formation control for a group of robots. The advantage of the proposed method is that it only requires local sensor information, allows for vehicle kinematics and does not require communication between vehicles.

\end{itemize}

\subsection*{Future Work}

The work in this report opens up a number of future research problems as follows:

\begin{itemize}
  \item \textit{Deadlock avoidance}. In Chapt.~\ref{chap:dead}, the deadlock avoidance problem
  for a pair of vehicles was explored. While an approach applicable to an arbitrary number 
  of vehicles would extremely non-trivial, it remains an open question as to whether decentralized MPC
  approaches can be extended to cover this case. Also, a version which do not require the use of a navigation
  function could be useful, through this may rely on Bug type behavior and would complicate analysis.
  
  \item \textit{Broader robustness analysis}. Noise and model disturbance was not considered in 
  Chapts.~\ref{chap:convsingle} or \ref{chap:dead}. Through it may complicate the analysis, 
  employing the full robust MPC technique from Chapt.~\ref{chap:singlevehicle} should be possible.
  
  \item \textit{Navigation system unification}. Currently the navigation system in Chapts.~\ref{chap:multiple} and \ref{chap:dead} has not been 
  unified with the navigation system in Chapt.~\ref{chap:convsingle}. A system which coordinates multiple vehicles while
  taking sensor constraints into account may be possible, however it
  may complicate subsequent analysis.

  \item \textit{Trajectory planning heuristic}. When designing incomplete trajectory planning systems, the design of the set
  of possible trajectories is largely heuristic, and could be tuned for better
  performance. A wider range of possible trajectory
  shapes has previously been considered \cite{Blanco2008journ}, through there is no evidence this is the best possible 
  choice. A possible avenue of research would be to apply some global optimization method
  to determine the best possible set of trajectories to consider for a particular class of scenarios.
 One possible advantage is that the trajectory planner may favor vehicle positions which provide a 
 marginally more informative view of the obstacle.
  
  \item \textit{Extensions to 3D}. There is no reason why the navigation framework presented here cannot be extended
  to a three dimensional workspace. In these cases the tractability of the 
  trajectory planner becomes a much more critical issue, so simplified planning approaches  
 would be well suited to this problem.
  
  \item \textit{Other navigation tasks}. A highly active area of research currently is formation control
  of a group of vehicles \cite{Weihua2011journ}. The multiplexed MPC framework was recently extended to achieve collision avoidance,
  successfully achieving decentralized, robust formation control \cite{Weihua2011journ}. It may be possible to extend the approach proposed in 
  Chapt.~\ref{chap:singlevehicle} to formation control problems, while retaining the relevant advantages.
  
  \item \textit{Potential field methods for acceleration bounded vehicles}. As mentioned in Chapt.~\ref{chap:lit}, the problem 
  of achieving collision avoidance using an APF method for acceleration bounded vehicles has not been fully 
  solved. Through it seems unlikely an APF method would achieve the same closed loop performance as MPC, 
  solution to this problem may be useful to provide more objective comparisons between the methods.

\end{itemize}

\addcontentsline{toc}{chapter}{Bibliography} \bibliography{bib}
\chapter{Simulations with a Realistic Helicopter Model}
\label{apdx:helpexp}

This appendix outlines preliminary simulations that were carried out with a realistic helicopter model.

\clearpage \section{Helicopter Model}

In order to further test the proposed navigation approach, the
basic navigation approach was also tested
against a realistic helicopter model
\cite{Garratt2006conf2, garratt_biologically_2007}. The helicopter model is based on
an autonomous 8.2 kg helicopter with a main rotor diameter of
$1.52m$, based on a Hirobo Eagle helicopter with conventional main
rotor and tail rotor configuration. The Eagle, described further in
\cite{Ahmed2008conf4} is electrically propelled and instrumented with
differential GPS, Inertial Measurement Unit and an on-board
autopilot.
\par
The simulation combines a linearized rotor aerodynamic model with non-linear
rigid body equations. A non-linear thrust and rotor inflow model is also
incorporated, involving an iterative scheme to simultaneously solve for the
induced down wash velocity through the rotor and the corresponding thrust.  This
approach produces a minimum complexity model which captures the essential
dynamics of the helicopter.
Balancing simulation fidelity against
practicality, a simulation has been created that is capable of simulating the
following effects:
 
\begin{itemize}
\item Exact non-linear rigid body equations of motion;
\item Wind gusts and turbulence;
\item First order main rotor flapping dynamics;
\item Hover, rear-wards, sideways and forward flight;
\item Dynamic effects of the Bell-Hiller stabilizer bar;
\item Fuselage and tail-plane aerodynamic forces;
\item Approximate servo dynamics; and
\item Sensor lags, filtering, offsets and noise.
\end{itemize}

The main rotor forces and moments are controlled by the collective and cyclic
pitch channels. The collective pitch control varies the average blade incidence of all
of the blades. Increasing the collective pitch control results in an increased angle of
attack of each blade and a subsequent increase in main rotor thrust. Decreasing the
collective pitch has the opposite effect. The vertical motion of the helicopter is thus
controlled by varying the collective pitch.

In order to achieve pitching and rolling moments, the orientation or tilt of the
rotor disk is changed by applying pitch that varies cyclically, once per revolution.
Increasing the blade pitch on one side of the rotor disk and decreasing it on the
opposite side causes the path of the blade tips, known as the Tip Path Plane (TPP)
to be tilted. Since the thrust vector acts essentially perpendicular to the TPP,
this can be used to change the trim of the helicopter.

\begin{figure}[ht]
    \centering
        \includegraphics[width=0.7\columnwidth]{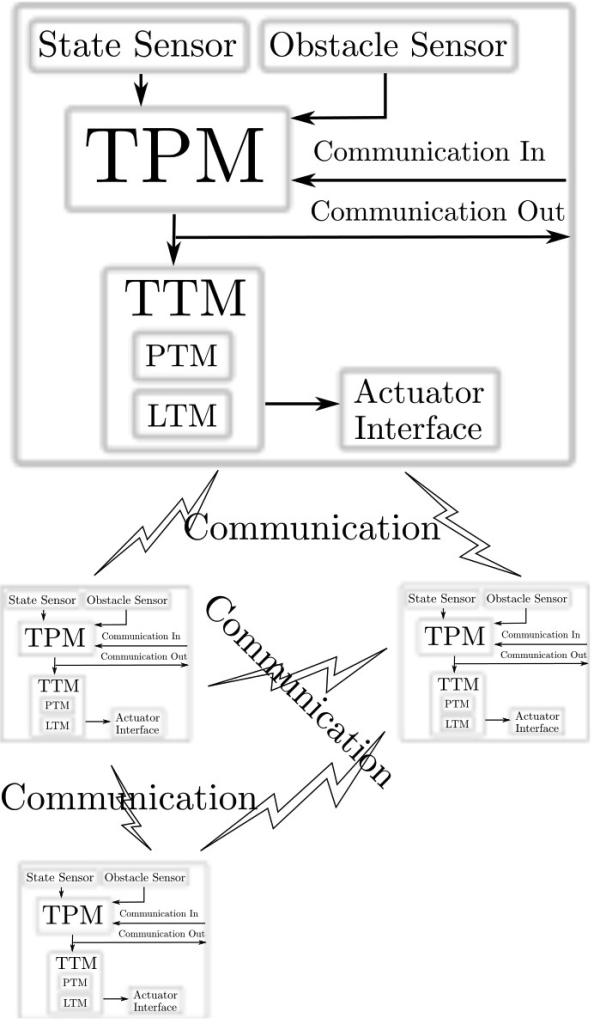}
    \caption{Block view of the helicopter model under test.}
    \label{F:main_sim}
\end{figure}

To simplify debugging, the simulation is divided into a number of blocks and
subsystems. On the highest level, the simulation consists of an aerodynamics
subsystem, sensor subsystem, sensor fusion block and controller subsystem as
shown in Fig.~\ref{F:main_sim}. Most of the computational blocks have been
implemented as C code S-functions.  A detailed discussion on the principles
behind this model can be found in \cite{garratt_biologically_2007}.

Initial adjustments of the simulation were made to match the trim control
settings of the simulation to the collective pitch, aileron and elevator
settings observed from flight test. The simulation was validated against actual
flight test data using frequency response techniques based on chirp and doublet
waveforms.

A commonly used scheme for controlling a helicopter is used for this part of the
controller and consists of an attitude feedback inner loop implemented as a PD
controller, combined with a PI based velocity outer loop.
This scheme effectively deals with the problem of the helicopter being
under-actuated. In the inner loop, cyclic pitch is used to control the
helicopter pitch and roll attitude. For the outer loop controller, the desired
pitch and roll attitude is set in response to the velocity error using PI
feedback.

The gains were tuned systematically using trial and error to converge on an
acceptable solution. In the first instance, the attitude control loop was tuned independently
by turning off the outer loop and stopping the integration of velocity
and position in the dynamics block. Once satisfactory stability was demonstrated,
the outer loop was re-activated and the velocity gains were then tuned.

The main inputs to this system are the desired velocities in the longitudinal
and lateral directions, and the desired heading. Height and heading were
controlled using separate PID control of collective pitch and tail rotor pitch
respectively. A constant height above ground was maintained by the collective
pitch PID controller.

\section{Testing}

To interface the proposed navigation law with the low level controllers in the helicopter model, the desired
velocity in the longitudinal and lateral directions were set to $[v_{long},$
$v_{lat}]$ $=$ $R(-\theta(k))\cdot v^*(1|k)$, where $R(\cdot)$ is the $2\times 2$
rotation matrix converting the coordinates from the world frame to the relative
helicopter reference frame. The desired heading was set to $\theta(1|k)$. 

For all these simulations, the controller refresh rate of 5 Hz was used. The parameters used for control can be found in Table~\ref{apxa:param}.

\begin{table}
\centering
\begin{tabular}{| l | c |}
\hline
 $u_{max}$ & $1.9 ms^{-2}$ \\
\hline
 $u_{p}$ & $1.5 ms^{-2}$  \\
\hline
 $u_{sp}$ & $1.0 ms^{-2}$  \\
\hline
\end{tabular} \hspace{10pt}
\begin{tabular}{| l | c |}
\hline
 $d_{tar}$ & $5 m$  \\
\hline
 $v_{max}$ & $4.5 ms^{-1}$  \\
\hline
 $\gamma_0$ & $10$  \\
\hline
\end{tabular} \\

\caption{Simulation parameters for collision avoidance with a realistic helicopter vehicle model.}
\label{apxa:param}
\end{table}

\begin{figure}[ht]
    \centering
        \includegraphics[width=0.5\columnwidth]{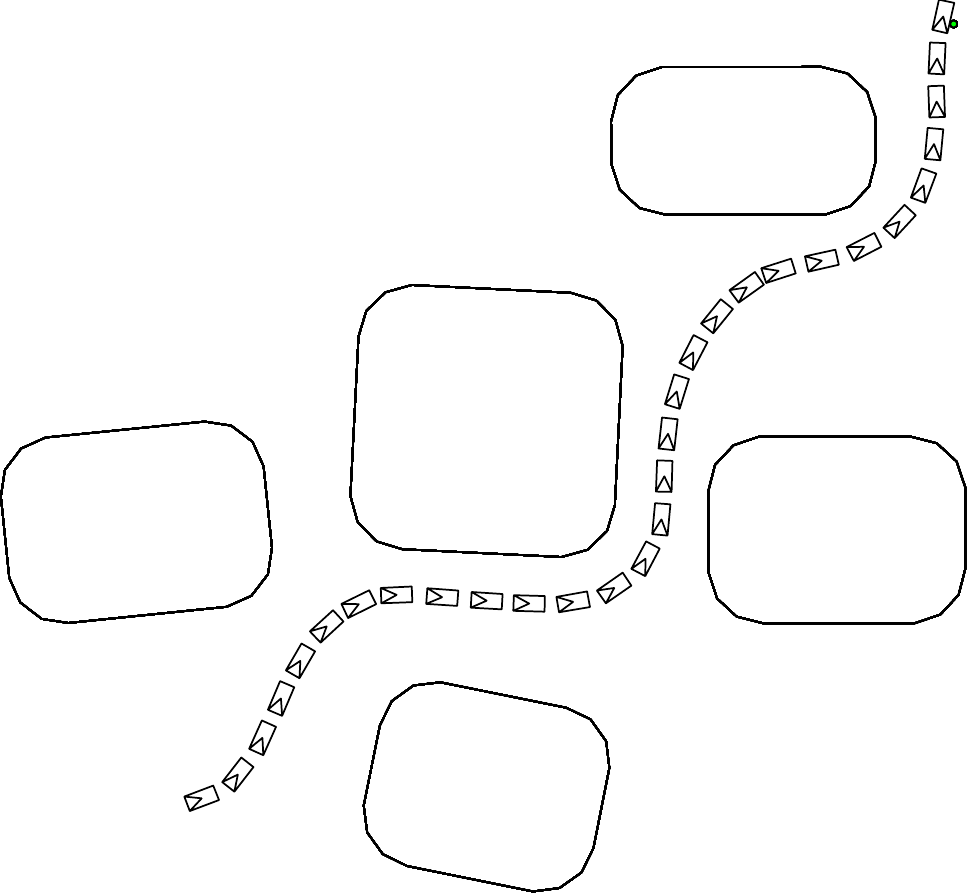}
    \caption{Simulations with a realistic helicopter model.}
    \label{fig:hfmpc}
\end{figure}

Fig.~\ref{fig:hfmpc} indicates the simulated helicopter was successfully able to
navigate a cluttered environment using the proposed method.
 This introductory experiment gives promising results for using the proposed method
 to navigate real world helicopters.

\end{document}